\documentclass[11pt,leqno]{amsart}
\usepackage{amsmath}
\usepackage{epsfig}
\usepackage{amssymb,latexsym}
\usepackage{texdraw}

\newtheorem{thm}{\bf Theorem}[section]
\newtheorem{defi}[thm]{\bf Definition}
\newtheorem{prop}[thm]{\bf Proposition}
\newtheorem{cor}[thm]{\bf Corollary}
\newtheorem{lem}[thm]{\bf Lemma}
\newtheorem{rem}[thm]{\bf Remark}
\newtheorem{ex}[thm]{\bf Example}
\numberwithin{equation}{section}

\DeclareMathOperator{\wt}{wt}
\newcommand{\nc}{\newcommand}
\newcommand{\pf}{{\bfseries Proof. }}

\nc{\op}{\oplus} \nc{\pv}{P^{\vee}}

\newcommand{\A}{\mathbb{A}}

\newcommand{\Q}{\mathbb{Q}}
\newcommand{\Z}{\mathbb{Z}}
\newcommand{\La}{\Lambda}
\newcommand{\la}{\lambda}
\newcommand{\ali}{\alpha_i}
\newcommand{\g}{\mathfrak{g}}
\newcommand{\ot}{\otimes}
\nc{\B}{\mathbf{B}} \nc{\V}{\mathbf{V}}

\nc{\nbinom}[2]{\genfrac{}{}{0pt}{1}{{#1}}{{#2}}}
\nc{\qbinom}[2]{\left[\genfrac{}{}{0pt}{1}{{#1}}{{#2}}\right]}
\nc{\fit}{\tilde{f}_i} \nc{\eit}{\tilde{e}_i} \nc{\Y}{\mathbf{Y}}

\nc{\ra}{\rightarrow} \nc{\vep}{\varepsilon} \nc{\vphi}{\varphi}
\nc{\h}{\mathfrak{h}} \nc{\oP}{\overline{P}}

\newbox{\tmppic}
\newbox{\tmpdraw}
\newbox{\tmpfig}
\newbox{\tmpa}
\newbox{\tmpb}

\begin{document}

\title[Fock space representation of quantum affine algebras]
      {Fock space representations of quantum affine algebras
      and generalized Lascoux-Leclerc-Thibon algorithm}
\author[S.-J. Kang and J.-H. Kwon]{Seok-Jin Kang$^{*}$ and Jae-Hoon Kwon$^{\dagger}$}
\address{School of Mathematics\\ Korea Institute for Advanced Study \\
         207-43 Cheongryangri-dong Dongdaemun-gu \\
         Seoul 130-012, Korea}
\thanks{$^{*}$This research was supported by KOSEF Grant
\# 98-0701-01-5-L and the Young Scientist Award, Korean Academy of
Science and Technology. \\
$^{\dagger}$This research was supported by KOSEF Grant \#
98-0701-01-05-L} \email{sjkang@kias.re.kr, jhkwon@kias.re.kr}

\maketitle

\begin{abstract}
We construct the Fock space representations of classical quantum
affine algebras using combinatorics of Young walls. We also show
that the crystal graphs of the Fock space representations can be
realized as the abstract crystal consisting of proper Young walls.
Finally, we give a generalized version of Lascoux-Leclerc-Thibon
algorithm for computing the global bases of the basic
representations of classical quantum affine algebras.
\end{abstract}


\section*{Introduction}

The {\it crystal basis theory} developed by Kashiwara (\cite
{Kas90,Kas91}) provides us with a very powerful combinatorial
method of studying the structure of integrable modules over
quantum groups. Let $M$ be an integrable module over a quantum
group $U_q({\mathfrak g})$ in the category ${\mathcal{O}}_{int}$
and let $\A_0$ denote the subring of $\mathbb{Q}(q)$ consisting of
regular functions at $q=0$. A {\it crystal basis} of $M$ is a pair
$(L, B)$, where $L$ is an $\A_0$-lattice of $M$ and $B$ is a
$\mathbb{Q}$-basis of $L/qL$ satisfying certain properties
involving {\it Kashiwara operators}. Thus a crystal basis can be
understood as a basis of $M$ at $q=0$ and the set $B$ is given a
structure of colored oriented graph, called the {\it crystal
graph}, that reflects the combinatorial structure of $M$.

It is known that every $U_q({\mathfrak g})$-module $M$ in the
category ${\mathcal{O}}_{int}$ is a direct sum of irreducible
highest weight modules with dominant integral highest weights and
that the crystal bases are preserved under this decomposition
(see, for example, \cite{HK2002,Kas91}). Hence it is a very
natural problem to find a concrete realization of the crystal
graph $B(\la)$ of an irreducible highest weight $U_q({\mathfrak
g})$-module $V(\la)$ with dominant integral highest weight $\la$.

Moreover, one can {\it globalize} the main idea of crystal basis
theory. Let $V(\la)$ be an irreducible highest weight
$U_q({\mathfrak g})$-module with dominant integral highest weight
$\la$, and let $(L(\la), B(\la))$ be the crystal basis of
$V(\la)$. Consider the $\mathbb{Q}$-algebra automorphism of
$U_q({\mathfrak g})$ defined by
$$\overline {q} = q^{-1}, \quad
\overline {e_i} = e_i, \quad \overline {f_i} = f_i, \quad
\overline {q^h} = q^{-h} \quad \text{for} \ \  i\in I, h\in
P^{\vee}.$$ Then we get a $\mathbb{Q}$-linear automorphism of
$V(\la)$ given by
$$P v_{\la} \longmapsto \overline{P} v_{\la}
\quad \text{for} \ \ P \in U_q({\mathfrak g}),$$ where $v_{\la}$
denotes the highest weight vector of $V(\la)$. In \cite{Kas91} and
\cite{Lus90}, Kashiwara and Lusztig independently showed that
there exists a unique {\it global basis} (or {\it canonical
basis}) $G(\la)=\{\,G(b)\,|\, b \in B(\la) \,\}$ of $V(\la)$
satisfying the properties:
\begin{equation*}
\overline {G(b)} =G(b), \quad G(b) \equiv b \ \ \text{mod} \,
qL(\la) \quad \text{for all} \ b\in B(\la).
\end{equation*}
Therefore one naturally gets interested in the following problem:
Given a realization of the crystal graph $B(\la)$ of $V(\la)$, can
we find an effective algorithm for constructing the global basis
$G(\la)$?

\vskip 2mm In \cite{Ka2000}, Kang introduced the notion of {\it
Young walls} as a new combinatorial scheme for realizing the
crystal bases for quantum affine algebras. The Young walls consist
of colored blocks with various shapes that are built on a given
{\it ground-state wall}, and they can be viewed as generalizations
of colored Young diagrams. More precisely, let $U_q({\mathfrak
g})$ be a classical quantum affine algebra of type $A_n^{(1)}$
$(n\ge 1)$, $A_{2n-1}^{(2)}$ $(n\ge 3)$, $D_n^{(1)}$ $(n\ge 4)$,
$A_{2n}^{(2)}$ $(n\ge 1)$, $D_{n+1}^{(2)}$ $(n\ge 2)$ and
$B_n^{(1)}$ $(n\ge 3)$, and let $\Lambda$ be a dominant integral
weight of level 1. Then the description of ground-state wall
$Y_{\La}$, the rules and patterns for building Young walls, and
the action of Kashiwara operators are given explicitly in terms of
combinatorics of Young walls. In this way, the set
$\mathcal{Z}(\La)$ of {\it proper Young walls} is given a
structure of {\it abstract crystal} for the quantum affine algebra
$U_q({\mathfrak g})$, and the crystal $B(\La)$ of the basic
representation $V(\La)$ is realized as the abstract crystal
${\mathcal Y}(\La)$ consisting of {\it reduced proper Young walls}
(see \cite{Ka2000} for more details).

The goal of this paper is to find an effective algorithm for
computing the global basis element $G(Y)$ in $G(\La)$ for each
reduced proper Young wall $Y$ in $\mathcal{Y}(\La)$. For this
purpose, we first construct the {\it Fock space representations}
of quantum affine algebras in a purely combinatorial way. We take
the Fock space ${\mathcal F}(\La)$ to be the
$\mathbb{Q}(q)$-vector space with a basis consisting of proper
Young walls, and define the $U_q({\mathfrak g})$-action on
${\mathcal F}(\La)$ in terms of combinatorics of Young walls. Then
the Fock space ${\mathcal F}(\La)$ becomes an integrable
$U_q({\mathfrak g})$-module in the category ${\mathcal{O}}_{int}$
(see Section 5).

We then show that the crystal of ${\mathcal F}(\La)$ is isomorphic
to the crystal ${\mathcal Z}(\La)$ consisting of proper Young
walls (Theorem \ref{crystal basis for Fock space}). As a
corollary, we get an explicit decomposition of the Fock space
${\mathcal F}(\La)$ into a direct sum of irreducible highest
weight $U_q({\mathfrak g})$-modules by locating the maximal
vectors in the crystal ${\mathcal Z}(\La)$ (Corollary
\ref{decomposition}).

In \cite{KMPY}, Kashiwara, Miwa, Petersen and Yung gave a more
abstract construction of the Fock space representations of quantum
affine algebras. For a dominant integral weight $\La$ of level $l
\ge 1$, the Fock space ${\mathcal F}(\La)$ was realized as the
inductive limit of $q$-deformed wedge spaces arising from a level
$l$ perfect representation. We expect one can construct the higher
level Fock space representations of quantum affine algebras using
combinatorics of Young walls.

Finally, we give a generalized version of Lascoux-Leclerc-Thibon
algorithm (\cite{LLT}) for constructing the global bases of basic
representations of classical quantum affine algebras of type
$A_n^{(1)}$ $(n\ge 1)$, $A_{2n-1}^{(2)}$ $(n\ge 3)$, $D_n^{(1)}$
$(n\ge 4)$, $A_{2n}^{(2)}$ $(n\ge 1)$, $D_{n+1}^{(2)}$ $(n\ge 2)$
and $B_n^{(1)}$ $(n\ge 3)$. More precisely, for each reduced
proper Young wall $Y$ in $\mathcal{Y}(\La)$, we obtain an
effective algorithm for computing the global basis element $G(Y)$
in $G(\Lambda)$ that can be expressed as a $\mathbb{Z}[q]$-linear
combination of proper Young walls (Theorem \ref{LLT}):
$$G(Y)= \sum_{Z \in {\mathcal Z}(\La)} G_{Y, Z}(q) Z \quad
\text{for some} \quad G_{Y, Z}(q) \in \mathbb{Z}[q].$$ By
construction, the matrix coefficients $G_{Y, Z}(q)$ satisfy
certain unitriangular conditions. We expect that there exist some
interesting algebraic structures such that the irreducible modules
at some specializations are parametrized by reduced proper Young
walls and that the decomposition matrices are determined by the
polynomials $G_{Y,Z}(q)$ giving the global basis elements
(cf.\cite{Ar,BK,Gro}).

\vskip 1cm


\section{Quantum Groups}

Let $I$ be a finite index set. A square matrix $A=(a_{ij})_{i,j
\in I}$ is called a {\it generalized Cartan matrix} if it
satisfies\,: (i) $a_{ii}=2$ for all $i\in I$, \ \ (ii) $a_{ij} \in
\Z_{\leq 0}$ for all $i,j\in I$, \ \ (iii) $a_{ij}=0$ implies
$a_{ji}=0$. In this paper, we assume that $A$ is {\it
symmetrizable}\,; i.e., there is a diagonal matrix $D=\text{diag}
(s_i \in \Z_{>0} | \, i\in I)$ with positive integral entries such
that $DA$ is symmetric.

Consider the free abelian group
\begin{equation}
P^{\vee} = \left(\bigoplus_{i\in I} \Z h_i \right) \oplus \left(
\bigoplus_{j=1}^{\text{corank}\, A} \Z d_j \right),
\end{equation}
and let ${\mathfrak h}= \Q \otimes_{\Z} P^{\vee}$. The free
abelian group $P^{\vee}$ is called the {\it dual weight lattice}
and the $\Q$-vector space ${\mathfrak h}$ is called the {\it
Cartan subalgebra}.

The {\it weight lattice} and the set of {\it simple coroots} are
defined to be
\begin{equation}
P = \{\, \la \in \mathfrak{h}^* \,|\, \ \la(P^{\vee}) \subset \Z
\,\}, \qquad \Pi^{\vee} = \{\, h_i \,|\, \ i\in I \,\}.
\end{equation}
We denote by $\Pi = \{\, \alpha_i \,|\, \ i\in I \,\}$ the set of
{\it simple roots}, which is a linearly independent subset of
${\mathfrak h}^*$ satisfying
\begin{equation}
\alpha_i(h_j)=a_{ji} \quad \text{for all} \ \ i, j\in I.
\end{equation}

\begin{defi}
{\rm The quintuple $(A, P^{\vee}, P, \Pi^{\vee}, \Pi)$ defined
above is called a {\it Cartan datum} associated with $A$. }
\end{defi}

\vskip 3mm

We denote by $P^{+}=\{\,\la \in P \,|\, \ \la(h_i) \ge 0 \ \
\text{for all} \ i\in I \,\}$ the set of {\it dominant integral
weights}. The free abelian group $Q=\bigoplus_{i\in I} \Z
\alpha_i$ is called the {\it root lattice}. We set $Q_{+} =
\sum_{i\in I} \Z_{\ge 0} \alpha_i$ and $Q_{-} = -Q_{+}$. There is
a partial ordering on ${\mathfrak h}^*$ defined by $\la \ge \mu$
if and only if $\la - \mu \in Q_{+}$. Since the generalized Cartan
matrix $A$ is symmetrizable, there is a nondegenerate symmetric
bilinear form $(\ | \ )$ on ${\mathfrak h}^*$ satisfying
\begin{equation*}
s_i = \dfrac{(\alpha_i|\alpha_i)}{2} \quad \text{and} \quad
\dfrac{2(\alpha_i | \alpha_j)}{(\alpha_i | \alpha_i)} = a_{ij}
\quad \text{for all} \ \ i,j \in I.
\end{equation*}

For an indeterminate $q$, set $q_i = q^{s_i}$ and define
\begin{equation*}
[n]_i = \dfrac{q_i^{n} - q_i^{-n}} {q_i - q_i^{-1}}, \quad [n]_i!
= \prod_{k=1}^n [k]_i, \quad \left[\begin{matrix} m \\ n
\end{matrix}\right]_{i}
 = \dfrac{[m]_i!}{[m-n]_i! [n]_i!}.
\end{equation*}

\begin{defi}\label{defining relations}
{\rm The {\it quantum group} $U_q({\mathfrak g})$ associated with
a Cartan datum $(A, P^{\vee}, P, \Pi^{\vee}, \Pi)$ is the
associative algebra over $\Q(q)$ with 1 generated by the symbols
$e_i$, $f_i$ $(i\in I)$ and $q^h$ $(h\in P^{\vee})$ subject to the
following defining relations\,:
\begin{equation}\label{defining rels}
\begin{aligned}
\ & q^0 = 1, \ \ q^h q^{h'} = q^{h+ h'} \quad (h, h'\in P^{\vee}),\\
\ & q^h e_i q^{-h} = q^{\ali(h)} e_i, \quad
q^h f_i q^{-h} = q^{-\ali(h)} f_i \quad (h\in P^{\vee}, \ i\in I), \\
\ & e_i f_j - f_j e_i = \delta_{ij} \frac{K_i - K_i^{-1}}{q_i -
q_i^{-1}},
\quad \text{where} \quad K_i =  q^{s_i h_i}, \\
\ & \sum_{k=0}^{1-a_{ij}} (-1)^k {\begin{bmatrix} 1-a_{ij} \\ k
\end{bmatrix}}_{i}
e_i^{1-a_{ij}-k} e_j e_i^{k} = 0 \qquad (i \neq j), \\
\ & \sum_{k=0}^{1-a_{ij}} (-1)^k {\begin{bmatrix} 1-a_{ij} \\ k
\end{bmatrix}}_{i} f_i^{1-a_{ij}-k} f_j f_i^{k} = 0 \qquad (i \neq
j).
\end{aligned}
\end{equation}
}
\end{defi}

The quantum group $U_q(\g)$ has a {\it Hopf algebra} structure
with the comultiplication $\Delta$, counit $\varepsilon$, and
antipode $S$ defined by
\begin{equation}
\begin{aligned}\text{}
& \Delta(q^h) = q^h \ot q^h, \\
& \Delta(e_i) = e_i \ot K_i^{-1} + 1 \ot e_i,
\quad \Delta(f_i) = f_i \ot 1 + K_i \ot f_i, \\
& \varepsilon(q^h)=1, \quad
\varepsilon(e_i) = \varepsilon(f_i) = 0, \\
& S(q^h) = q^{-h}, \quad S(e_i) = - e_i K_i, \quad S(f_i) = -
K_i^{-1} f_i
\end{aligned}
\end{equation}
for $h\in P^{\vee}$ and $i\in I$.

\vskip 3mm Let $U^+$ (resp. $U^-$) be the subalgebra of
$U_q(\mathfrak g)$ generated by the elements $e_i$ (resp. $f_i$)
for $i\in I$, and let $U^0$ be the subalgebra of $U_q({\mathfrak
g})$ generated by $q^h$ $(h\in P^{\vee})$. Then we have the {\it
triangular decomposition}\,:
\begin{equation*}
U_q({\mathfrak g}) \cong U^- \otimes U^0 \otimes U^+.
\end{equation*}

A $U_q(\g)$-module $V$ is called a {\it weight module} if it
admits a {\it weight space decomposition} $V = \bigoplus_{\mu \in
P} V_{\mu}$, where $V_{\mu} = \{ v \in V | \ q^h v = q^{\mu(h)} v
\ \ \text{for all} \ h\in P^{\vee} \}$. If $\dim_{\Q(q)} V_{\mu} <
\infty$ for all $\mu \in P$, we define the {\it character} of $V$
by
$$\text{ch} V = \sum_{\mu \in P}
\left(\dim_{\Q(q)} V_{\mu}\right) e^{\mu},$$ where $e^{\mu}$ are
basis elements of the group algebra $\Q(q)[P]$ with the
multiplication given by $e^{\mu} e^{\nu} = e^{\mu+\nu}$ for all
$\mu, \nu \in P$.

A weight module $V$ over $U_q(\g)$ is called a {\it highest weight
module with highest weight $\la$} ($\la\in P$) if there exists a
non-zero vector $v_{\la} \in V$ (called the {\it highest weight
vector}) such that \ (i) $e_i \, v_{\la} = 0$ for all $i\in I$, \
\ (ii) $q^{h} v_{\la} = q^{\la(h)} v_{\la}$ for all $h\in
P^{\vee}$, \ \ (iii) $V= U_q({\mathfrak g}) v_{\la}$.

For example, let $J(\la)$ denote the left ideal of $U_q({\mathfrak
g})$ generated by $e_i$, $q^h - q^{\la(h)} 1$ for $i\in I$, $h\in
P^{\vee}$, and set $M(\la) = U_q({\mathfrak g}) / J(\la)$. Then,
via left multiplication, $M(\la)$ becomes a highest weight
$U_q({\mathfrak g})$-module with highest weight $\la$, called the
{\it Verma module}, and it satisfies the following properties\,:

\begin{prop} {\rm (cf. \cite{HK2002})} \hfill

{\rm (a)} $M(\la)$ is a free $U^-$-module of rank 1.

{\rm (b)} Every highest weight $U_q({\mathfrak g})$-module with
highest weight $\la$ is a homomorphic image of $M(\la)$.

{\rm (c)} $M(\la)$ contains a unique maximal submodule $R(\la)$.

\end{prop}

\noindent The unique irreducible quotient $V(\la) = M(\la) /
R(\la)$ is called the {\it irreducible highest weight
$U_q({\mathfrak g})$-module with highest weight $\la$}.

\begin{defi} \label{defi:O^q_{int}}
{\rm The {\it category ${\mathcal{O}}_{int}$} consists of
$U_q(\g)$-modules $M$ satisfying the following properties\,:

\hskip 2mm (i) $M$ is a weight module,

\hskip 2mm (ii) there exist finitely many $\la_1, \cdots, \la_s
\in P$ such that
$$\text{wt}(M) := \{ \mu \in P | \
M_{\mu} \neq 0 \} \subset \bigcup_{j=1}^s \left(\la_j -
Q_{+}\right),$$

\hskip 3mm (iii) $e_i$ and $f_i$ $(i\in I)$ are locally nilpotent
on $M$. }
\end{defi}

The basic properties of the category ${\mathcal{O}}_{int}$ are
given in the following proposition.

\begin{prop} \label{prop:O^q_{int}} {\rm (cf. \cite{HK2002})} \hfill

{\rm (a)} For each $i\in I$, let $U_{(i)}$ be the subalgebra of
$U_q(\g)$ generated by $e_i$, $f_i$, $K_i^{\pm 1}$, which is
isomorphic to the quantum group $U_q(\frak{sl}_2)$. Then every
$U_q(\g)$-module $M$ in the category ${\mathcal O}_{int}$ is a
direct sum of finite dimensional irreducible $U_{(i)}$-submodules.

{\rm (b)} The category ${\mathcal O}_{int}$ is semisimple.
Moreover, every irreducible object in the category ${\mathcal
O}_{int}$ has the form $V(\la)$ with $\la \in P^{+}$.\qed

\end{prop}

\vskip 1cm


\section{Crystal Bases}

In this section, we briefly review the {\it crystal basis theory}
for quantum groups developed by Kashiwara (\cite{Kas90, Kas91}).
We will also use the following notation for {\it divided
powers}\,:
$$e_i^{(n)} = e_i^{n} / [n]_i!, \qquad
f_i^{(n)} = f_i^{n} / [n]_i!.$$

Fix an index $i\in I$ and let $M=\bigoplus_{\la \in P} M_{\la}$ be
a $U_q(\g)$-module in the category ${\mathcal O}_{int}$. By the
representation theory of $U_q(\frak{sl}_2)$, every element $v\in
M_\la$ can  be written uniquely as
\begin{equation*}
v = \sum_{k\geq0} f_i^{(k)} v_k,
\end{equation*}
where $k \ge -\la(h_i)$ and $v_k \in \ker e_i \cap M_{\la+k\ali}$.
We define the endomorphisms $\eit$ and $\fit$ on $M$, called the
{\it Kashiwara operators},  by
\begin{equation}
\eit v = \sum_{k\geq1} f_i^{(k-1)} v_k, \qquad \fit v =
\sum_{k\geq0} f_i^{(k+1)} v_k.
\end{equation}

Let $\A_{0} = \{ f/g \in \Q(q)  \, | \, f, g \in \Q[q], \,
g(0)\neq 0 \}$ be the subring of $\Q(q)$ consisting of the
rational functions in $q$ that are regular at $q=0$.

\begin{defi} \label{defi:crystal basis}
{\rm A {\it crystal basis} of $M$ is a pair $(L, B)$, where

\hskip 2mm (i) $L$ is a free $\A_{0}$-submodule $M$ such that $M
\cong \Q(q) \otimes_{\A_{0}} L$,

\hskip 2mm (ii) $B$ is a basis of the $\Q$-vector space $L/qL$,

\hskip 2mm (iii) $L=\bigoplus_{\lambda \in P} L_{\la}$, where
$L_{\la} = L \cap M_{\la}$,

\hskip 2mm (iv) $B=\bigsqcup_{\la \in P} B_{\la}$, where
$B_{\la}=B \cap \left(L_{\la} / q L_{\la} \right)$,

\hskip 2mm (v) $\eit L \subset L$, $\fit L \subset L$ for all
$i\in I$,

\hskip 2mm (vi) $\eit B\subset B\cup \{0\}$, \ $\fit B\subset
B\cup \{0\}$ for all $i\in I$,

\hskip 2mm (vii) for $b, b'\in B$, $\fit b = b'$ if and only if $b
= \eit b'$. }
\end{defi}

\noindent The set $B$ is given a colored oriented graph structure
with the arrows defined by
$$b \stackrel{i} \longrightarrow b' \quad \text{if and only if}
\quad \fit b = b'.$$ The graph associated with $B$ is called the
{\it crystal graph} of $M$ and it reflects the combinatorial
structure of $M$. For instance, we have
\begin{equation*}
\dim_{\Q(q)} M_{\lambda} =\# B_{\lambda} \quad \text{for all} \ \
\lambda \in P.
\end{equation*}

\vskip 3mm

Let $(L,B)$ be a crystal basis of a $U_q(\g)$-module $M$ in the
category ${\mathcal{O}}_{int}$. For each $b\in B$ and $i\in I$, we
define
\begin{equation*}
\varepsilon_i(b)=\max \{\,k\ge 0 \,|\, \, \eit^k b \in B \,\},
\qquad \varphi_i(b) = \max \{\, k\ge 0 \,|\, \, \fit^k b \in B
\,\}.
\end{equation*}

\noindent Then the set $B$ satisfies the following properties.

\vskip 3mm

\begin{prop} {\rm (\cite{Kas91, Kas93, Kas94})}

{\rm (a)} For all $i\in I$ and $b\in B$, we have
\begin{equation*}
\begin{aligned}\text{}
& \varphi_i(b) = \varepsilon_i(b) + \langle h_i, \wt(b) \rangle, \\
& \wt(\eit b)=\wt(b) + \alpha_i, \\
& \wt(\fit b)=\wt(b) - \alpha_i.
\end{aligned}
\end{equation*}

{\rm (b)} If $\eit b \in B$, then
$$\varepsilon_i(\eit b) = \varepsilon_i(b) - 1, \quad
\varphi_i(\eit b) = \varphi_i(b) + 1.$$

{\rm (c)} If $\fit b \in B$, then
$$\varepsilon_i(\fit b) = \varepsilon_i(b) + 1, \quad
\varphi_i(\fit b) = \varphi_i(b) - 1.$$

\end{prop}

\vskip 3mm Moreover, the crystal bases have extremely simple
behavior with respect to taking the tensor product.

\vskip 3mm

\begin{prop} {\rm (\cite{Kas90, Kas91})}
\label{prop:tensor product}  \hfill

Let $M_j$ $(j=1, 2)$ be a  $U_q(\g)$-module in the category
${\mathcal O}_{int}$ and let $(L_j, B_j)$ be its crystal basis.
Set
\begin{equation*}
L = L_1 \otimes_{\A_{0}} L_2, \quad B= B_1 \times B_2.
\end{equation*}
Then $(L,B)$ is a crystal basis of $M_1 \otimes_{\Q(q)} M_2$ with
the Kashiwara operators on $B$ given by
\begin{equation*}
\begin{aligned}\text{}
\eit(b_1\ot b_2)
  &=
  \begin{cases}
  \eit b_1 \ot b_2 & \text{if $\vphi_i(b_1) \geq \vep_i(b_2)$,}\\
  b_1\ot\eit b_2 & \text{if $\vphi_i(b_1) < \vep_i(b_2)$,}
  \end{cases}\\
\fit(b_1\ot b_2)
  &=
  \begin{cases}
  \fit b_1 \ot b_2 & \text{if $\vphi_i(b_1) > \vep_i(b_2)$,}\\
  b_1\ot\fit b_2 & \text{if $\vphi_i(b_1) \leq \vep_i(b_2)$.}
  \end{cases}
\end{aligned}
\end{equation*}\qed
\end{prop}

\vskip 3mm

The set $B_1\times B_2$ will be denoted by $B_1\otimes B_2$. The
{\it tensor product rule} in Proposition \ref{prop:tensor product}
gives a very convenient combinatorial description of the action of
Kashiwara operators on the multi-fold tensor product of crystal
basis. Let $M_j$ be a $U_q(\g)$-module in the category ${\mathcal
O}_{int}$ with a crystal basis $(L_j, B_j)$ \ $(j=1, \cdots, N)$.
Fix an index $i\in I$ and consider a vector $b=b_1 \ot \cdots \ot
b_N \in B_1 \ot \cdots \ot B_N$. To each $b_j \in B_j$ \, $(j=1,
\cdots , N)$, we assign a sequence of $-$'s and $+$'s with as many
$-$'s as $\vep_i(b_j)$ followed by as many $+$'s as
$\vphi_i(b_j)$\,:
\begin{equation*}
\begin{aligned}\text{}
b \, &  =  \, b_1 \otimes b_2 \otimes \cdots \otimes b_N \\
& \longmapsto
 (\underbrace{-, \cdots, -}_{\vep_i(b_1)},
\underbrace{+, \cdots, +}_{\vphi_i(b_1)}, \cdots \cdots,
\underbrace{-, \cdots, -}_{\vep_i(b_N)}, \underbrace{+, \cdots,
+}_{\vphi_i(b_N)}).
\end{aligned}
\end{equation*}
In this sequence, we cancel out all the $(+,-)$-pairs to obtain a
sequence of $-$'s followed by $+$'s\,:
\begin{equation} \label{eq:signature}
\text{$i$-sgn}(b) = (-, -, \cdots, -, +, +, \cdots, +).
\end{equation}
The sequence $\text{$i$-sgn}(b)$ is called the {\it $i$-signature}
of $b$.

Now the tensor product rule tells that $\eit$ acts on $b_j$
corresponding to the rightmost $-$ in the $i$-signature of $b$ and
$\fit$ acts on $b_k$ corresponding to the leftmost $+$ in the
$i$-signature of $b$\,:
\begin{equation}\label{eq:multifold tensor product}
\begin{aligned}\text{}
& \eit b = b_1 \otimes \cdots \ot \eit b_j \ot \cdots \ot b_N, \\
& \fit b = b_1 \otimes \cdots \ot \fit b_k \ot \cdots \ot b_N.
\end{aligned}
\end{equation}
We define $\eit b=0$ (resp. $\fit b=0$) if there is no $-$ (resp.
$+$) in the $i$-signature of $b$.

\vskip 3mm

By extracting the properties of crystal graphs, we define the
notion of abstract {\it crystals} as follows (\cite{Kas93,
Kas94}).

\begin{defi}
{\rm Let $(A, P^{\vee}, P, \Pi^{\vee}, \Pi)$ be a Cartan datum and
let $U_q(\g)$ be the corresponding quantum group.

A {\it crystal associated with} $(A, P^{\vee}, P, \Pi^{\vee},
\Pi)$ (or a {\it $U_q(\g)$-crystal}) is a set $B$ together with
the maps $\wt : B \ra P$, $\vep_i : B \ra \Z\cup\{-\infty\}$,
$\vphi_i : B \ra \Z\cup\{-\infty\}$, $\eit : B \ra B\cup\{0\}$,
and $\fit : B \ra B\cup\{0\}$ satisfying the following
conditions\,:

\hskip 3mm (i) for all $i\in I$, $b\in B$, we have
\begin{equation*}
\begin{aligned}\text{}
& \varphi_i(b) = \varepsilon_i(b) + \langle h_i, \wt(b) \rangle, \\
& \wt(\eit b)=\wt(b) + \alpha_i, \\
& \wt(\fit b)=\wt(b) - \alpha_i,
\end{aligned}
\end{equation*}

\hskip 3mm (ii) if $\eit b \in B$, then
$$\varepsilon_i(\eit b) = \varepsilon_i(b) - 1, \quad
\varphi_i(\eit b) = \varphi_i(b) + 1,$$

\hskip 3mm (iii) if $\fit b \in B$, then
$$\varepsilon_i(\fit b) = \varepsilon_i(b) + 1, \quad
\varphi_i(\fit b) = \varphi_i(b) - 1,$$

\hskip 3mm (iv) $\fit b = b'$ if and only if $b = \eit b'$ for all
$i\in I$, $b, b' \in B$,

\vskip 3mm \hskip 3mm (v) if $\vep_i(b) = - \infty$, then $\eit b
= \fit b = 0$. }
\end{defi}

For instance, if $(L,B)$ is a crystal basis of a $U_q(\g)$-module
in the category ${\mathcal O}_{int}$, then $B$ is a
$U_q(\g)$-crystal.

\vskip 3mm

\begin{defi}
{\rm Let $B_1$ and $B_2$ be crystals. A {\it morphism of crystals}
(or a {\it crystal morphism}) $\psi:B_1 \ra B_2$ is a map
$\psi:B_1 \cup \{0\} \ra B_2 \cup \{0\}$ satisfying the
conditions\,:

\hskip 3mm (i) $\psi(0)=0$,

\hskip 3mm (ii) if $b\in B_1$ and $\psi(b) \in B_2$, then
$$\wt(\psi(b))= \wt(b), \ \  \vep_i(\psi(b))=\vep_i(b),
\ \ \vphi_i(\psi(b)) = \vphi_i(b),$$

\hskip 3mm (iii) if $b, b'\in B_1$, $\psi(b), \psi(b') \in B_2$
and $\fit b =b'$, then $\fit \psi(b) = \psi(b')$. }
\end{defi}


\begin{defi}
{\rm The {\it tensor product} $B_1 \otimes B_2$ of the crystals
$B_1$ and $B_2$  is defined to be the set $B_1 \times B_2$ whose
crystal structure is given by
\begin{equation}
\begin{aligned}\text{}
\wt(b_1\ot b_2) &= \wt(b_1) + \wt(b_2),\\
\vep_i(b_1\ot b_2)
  &= \max(\vep_i(b_1), \vep_i(b_2) - \langle h_i, \wt(b_1) \rangle),\\
\vphi_i(b_1\ot b_2)
  &= \max(\vphi_i(b_2), \vep_i(b_1) + \langle h_i, \wt(b_2) \rangle),\\
\eit(b_1\ot b_2)
  &=
  \begin{cases}
  \eit b_1 \ot b_2 & \text{if $\vphi_i(b_1) \geq \vep_i(b_2)$,}\\
  b_1\ot\eit b_2 & \text{if $\vphi_i(b_1) < \vep_i(b_2)$,}
  \end{cases}\\
\fit(b_1\ot b_2)
  &=
  \begin{cases}
  \fit b_1 \ot b_2 & \text{if $\vphi_i(b_1) > \vep_i(b_2)$,}\\
  b_1\ot\fit b_2 & \text{if $\vphi_i(b_1) \leq \vep_i(b_2)$.}
  \end{cases}
\end{aligned}
\end{equation}
}
\end{defi}
\noindent Here, we denote $b_1 \otimes b_2 = (b_1, b_2)$ and use
the convention that $b_1 \ot 0 = 0 \ot b_2 =0$. With the above
definitions, we can check that the category of crystals becomes a
{\it tensor category}.

\vskip 3mm The existence and uniqueness theorem for crystal bases
is given as follows.

\begin{thm} {\rm (\cite{Kas91})} \label{thm:existence1}
Let $V(\la)$ be the irreducible highest weight $U_q(\g)$-module
with highest weight $\la \in P^{+}$ and highest weight vector
$v_{\la}$. Let $L(\la)$ be the free $\A_{0}$-submodule of $V(\la)$
spanned by the vectors of the form $\tilde f_{i_1} \cdots \tilde
f_{i_r} u_{\la}$ $(i_k \in I, r\in \Z_{\ge 0})$ and set
\begin{equation*}
B(\la) = \{ \tilde f_{i_1} \cdots \tilde f_{i_r} u_{\la} + q
L(\la) \in  L(\la) / q L(\la) \} \setminus \{0\}.
\end{equation*}
Then $(L(\la), B(\la))$ is a crystal basis of $V(\la)$ and every
crystal basis of $V(\la)$ is isomorphic to $(L(\la), B(\la))$.\qed
\end{thm}
One can globalize the main idea of crystal basis theory. Consider
the $\mathbb{Q}$-algebra automorphism of $U_q(\frak{g})$ defined
by
\begin{equation}
\overline {q} = q^{-1}, \quad \overline {e_i} = e_i, \quad
\overline {f_i} = f_i, \quad \overline {q^h} = q^{-h} \quad
\text{for} \ \  i\in I, h\in P^{\vee}.
\end{equation}
This induces a $\mathbb{Q}$-linear automorphism of $V(\la)$ given
by
\begin{equation}
P v_{\la} \longmapsto \overline{P} v_{\la} \quad \text{for} \ \ P
\in U_q({\mathfrak g}),
\end{equation}
where $v_{\la}$ denotes the highest weight vector of $V(\la)$. Let
$\A = \Q[q, q^{-1}]$ and define $V(\la)_{\A} = U_{\A}^{-}
v_{\la}$, where $U_{\A}^{-}$ is the $\A$-subalgebra of $U_q(\g)$
generated by $f_i^{(n)}$ $(i\in I, n \in \Z_{\ge 0})$. Then we
have
\begin{thm}{\rm (\cite{Kas91,Lus90})} \label{thm:existence2}
There exists a unique $\A$-basis $G(\la) =\{\,G(b)\in
V(\la)_{\A}\cap L(\la)  \,|\, \, b\in B(\la) \,\}$ of
$V(\la)_{\A}$ such that
\begin{equation*}
\overline {G(b)} =G(b), \quad G(b) \equiv b \mod{qL(\la)} \quad
\text{for all} \ b\in B(\la).
\end{equation*}\qed
\end{thm}

\vskip 3mm

The basis $G(\la)$ of $V(\la)$ given in Theorem
\ref{thm:existence2} is called the {\it global basis} or the {\it
canonical basis} of $V(\la)$ associated with the crystal basis
$(L(\la),B(\la))$.

\vskip 1cm


\section{Quantum Affine Algebras}

Let $I=\{0,1,\cdots, n \}$ be an index set and let
$A=(a_{ij})_{i,j \in I}$ be a generalized Cartan matrix of {\it
affine type}. We denote by
\begin{equation*}
P^{\vee} = \Z h_0 \op \Z h_1 \op \cdots \op \Z h_n \op \Z d
\end{equation*}
the {\it dual weight lattice} and $\Pi^{\vee}=\{h_i | \ i \in I
\}$ the {\it simple coroots}. The {\it simple roots} $\alpha_i$
and the {\it fundamental weights} $\La_i$ are given by
\begin{equation*}
\begin{aligned}\text{}
\alpha_i (h_j) & = a_{ji}, \qquad \alpha_i (d) = \delta_{0, i}, \\
\Lambda_i(h_j) & = \delta_{ij}, \qquad \Lambda_i(d)=0 \qquad (i,j
\in I).
\end{aligned}
\end{equation*}
We define the {\it affine weight lattice} to be
\begin{equation*}
P=\{ \lambda \in {\mathfrak h}^* | \ \lambda (P^{\vee}) \subset \Z
\}.
\end{equation*}

\vskip 3mm

The quintuple $(A, \Pi, \Pi^{\vee}, P, P^{\vee})$ is called an
{\it affine Cartan datum}. To each affine Cartan datum, we can
associate an infinite dimensional Lie algebra  $\mathfrak{g}$
called the {\it affine Kac-Moody algebra} (\cite{Kac90}). The
center of the affine Kac-Moody algebra $\mathfrak{g}$ is
1-dimensional and is generated by the {\it canonical central
element}
\begin{equation*}
c=c_0 h_0 + c_1 h_1 + \cdots + c_n h_n.
\end{equation*}
Moreover, the imaginary roots of ${\mathfrak g}$ are nonzero
integral multiples of the {\it null root}
\begin{equation*}
\delta = d_0 \alpha_0 + d_1 \alpha_1 + \cdots + d_n \alpha_n.
\end{equation*}
Here, $c_i$ and $d_i$ $(i\in I)$ are the non-negative integers
given in \cite{Kac90}.

Using the fundamental weights and the null root, the affine weight
lattice can be written as
\begin{equation*}
P=\Z \Lambda_0 \op \Z \Lambda_1 \op \cdots \op \Z \Lambda_n \op \Z
\delta.
\end{equation*}
Set
\begin{equation*}
P^+=\{ \lambda \in P \, | \, \lambda(h_i) \in \Z_{\ge 0} \ \ \
\text{for all} \ \ i\in I \}.
\end{equation*}
The elements of $P$ (resp. $P^{+}$) are called the {\it affine
weights} (resp. {\it affine dominant integral weights}). The {\it
level} of an affine dominant integral weight $\lambda \in P^+$ is
defined to be the nonnegative integer $\lambda(c)$.

\vskip 3mm

\begin{defi}
{\rm The {\it quantum affine algebra} $U_q({\mathfrak g})$ is the
quantum group associated with the affine Cartan datum $(A, \Pi,
\Pi^{\vee}, P, P^{\vee})$.}
\end{defi}

\vskip 3mm The subalgebra of $U_q({\mathfrak g})$ generated by
$e_i$, $f_i$, $K_i^{\pm1}$ $(i\in I)$ is denoted by
$U'_q({\mathfrak g})$, and is also called the {\it quantum affine
algebra}.




In this paper, we will focus on the quantum affine algebras of
type $A_n^{(1)}$ ($n\geq 1$), $A_{2n-1}^{(2)}$ ($n\geq 3$),
$D_n^{(1)}$ ($n\geq 4$), $A_{2n}^{(2)}$ ($n\geq 1$),
$D_{n+1}^{(2)}$ ($n\geq 2$) and $B_n^{(1)}$ ($n\geq 3$).

\vskip 1cm


\section{Combinatorics of Young Walls}

In \cite{Ka2000}, Kang introduced a new family of combinatorial
objects called the {\it Young walls} which can be viewed as
generalizations of colored Young diagrams, and gave a realization
of crystals $B(\La)$ for the basic representations of quantum
affine algebras in terms of {\it reduced proper Young walls}. In
this section, we briefly explain the combinatorics of Young walls.

The Young walls are built of colored blocks of three different
shapes. They are called the blocks of {\it type {\rm I}, type {\rm
II}}, and {\it type {\rm III}}, respectively. For each type of
quantum affine algebras, we use different sets of colored blocks.

\vskip 5mm

\renewcommand{\arraystretch}{2}
\begin{center}
\begin{tabular}{c|c|c|c|c|c}
Type & Shape & Width & Thickness & Height & Volume \\
\hline I & \raisebox{-0.4\height}{
\begin{texdraw}
\drawdim em \setunitscale 0.1 \linewd 0.5 \move(-10 0)\lvec(0
0)\lvec(0 10)\lvec(-10 10)\lvec(-10 0) \move(0 0)\lvec(5 5)\lvec(5
15)\lvec(-5 15)\lvec(-10 10) \move(0 10)\lvec(5 15)
\end{texdraw}
} $=$ \raisebox{-0.4\height}{
\begin{texdraw}
\drawdim em \setunitscale 0.1 \linewd 0.5 \move(-10 0)\lvec(0
0)\lvec(0 10)\lvec(-10 10)\lvec(-10 0)
\end{texdraw}}
 & 1 & 1 & 1 & 1 \\
II & \raisebox{-0.4\height}{
\begin{texdraw}
\drawdim em \setunitscale 0.1 \linewd 0.5 \move(0 0)\lvec(10
0)\lvec(10 5)\lvec(0 5)\lvec(0 0) \move(10 0)\lvec(15 5)\lvec(15
10)\lvec(5 10)\lvec(0 5) \move(10 5)\lvec(15 10)
\end{texdraw}
} $=$ \raisebox{-0.4\height}{
\begin{texdraw}
\drawdim em \setunitscale 0.1 \linewd 0.5 \textref h:C v:C \move(0
0)\lvec(10 0)\lvec(10 5)\lvec(0 5)\lvec(0 0)
\end{texdraw}
}
& 1 & 1 & $\frac{1}{2}$ & $\frac{1}{2}$ \\
III &

\raisebox{-0.4\height}{
\begin{texdraw}
\drawdim em \setunitscale 0.1 \linewd 0.5 \move(-10 0)\lvec(0
0)\lvec(0 10)\lvec(-10 10)\lvec(-10 0) \move(0 0)\lvec(2.5
2.5)\lvec(2.5 12.5)\lvec(-7.5 12.5)\lvec(-10 10) \move(0
10)\lvec(2.5 12.5) \lpatt(0.3 1) \move(0 0)\lvec(-2.5
-2.5)\lvec(-12.5 -2.5)\lvec(-10 0)
\end{texdraw}
} $=$ \raisebox{-0.4\height}{
\begin{texdraw}
\drawdim em \setunitscale 0.1 \linewd 0.5 \move(-10 0)\lvec(0
0)\lvec(0 10)\lvec(-10 0) \lpatt(0.3 1)\move(-10 0)\lvec(-10
10)\lvec(0 10)
\end{texdraw}}
, \hskip 2mm \raisebox{-0.4\height}{
\begin{texdraw}
\drawdim em \setunitscale 0.1 \linewd 0.5 \move(-10 0)\lvec(0
0)\lvec(0 10)\lvec(-10 10)\lvec(-10 0) \move(0 0)\lvec(2.5
2.5)\lvec(2.5 12.5)\lvec(-7.5 12.5)\lvec(-10 10) \move(0
10)\lvec(2.5 12.5) \lpatt(0.3 1) \move(2.5 2.5)\lvec(5 5)\lvec(2.5
5)
\end{texdraw}
} $=$ \raisebox{-0.4\height}{
\begin{texdraw}
\drawdim em \setunitscale 0.1 \linewd 0.5 \move(-10 0)\lvec(-10
10)\lvec(0 10)\lvec(-10 0)\lpatt(0.3 1)\move(-10 0)\lvec(0
0)\lvec(0 10)
\end{texdraw}
}

& 1 & $\frac{1}{2}$ & 1 & $\frac{1}{2}$ \\

\end{tabular}
\end{center}
\renewcommand{\arraystretch}{1}

\vskip 5mm

For each dominant integral weight $\La$ of level 1, we fix a frame
$Y_{\La}$ called the {\it ground-state wall of weight $\Lambda$},
and on this frame, we build a wall of thickness less than or equal
to one unit which extends infinitely to the left. The rules for
building the walls are as follows:

\begin{itemize}
\item[(1)] The colored blocks should be stacked in columns. No block
can be placed on top of a column of half-unit thickness.

\item[(2)] Except for the right-most column, there should be no free
space to the right of any block.

\item[(3)] The colored blocks should be stacked in a specified
pattern which is determined by the type of the quantum affine
algebra $U_q(\frak{g})$ and the level 1 dominant integral weight
$\Lambda$.
\end{itemize}

\noindent The coloring of blocks, the description of ground-state
walls and the patterns for building the walls are given in
\cite{Ka2000} (see also Appendix).

\vskip 5mm

\begin{ex}
{\rm If $\frak{g}=B^{(1)}_3$ and $\Lambda=\Lambda_0$, we use the
colored blocks

\vskip 5mm \hskip 4cm
\raisebox{-0.33\height}[0.69\height][0.35\height]
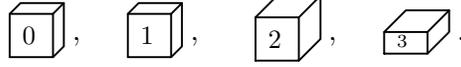
{%
\begin{texdraw}
\fontsize{10}{10}\drawdim em \setunitscale 0.15 \linewd 0.4
\move(0 0)\lvec(10 0)\lvec(10 10)\lvec(0 10)\lvec(0 0) \move(10
0)\lvec(12.5 2.5)\lvec(12.5 12.5)\lvec(2.5 12.5)\lvec(0 10)
\move(10 10)\lvec(12.5 12.5) \htext(3 3){$0$}
\end{texdraw}%
} , \hskip 5mm
\raisebox{-0.33\height}[0.69\height][0.35\height]{%
\begin{texdraw}
\fontsize{10}{10}\drawdim em \setunitscale 0.15 \linewd 0.4
\move(0 0)\lvec(10 0)\lvec(10 10)\lvec(0 10)\lvec(0 0) \move(10
0)\lvec(12.5 2.5)\lvec(12.5 12.5)\lvec(2.5 12.5)\lvec(0 10)
\move(10 10)\lvec(12.5 12.5) \htext(3 3){$1$}
\end{texdraw}%
} , \hskip 5mm \raisebox{-0.33\height}[0.69\height][0.35\height]{
\begin{texdraw}
\fontsize{10}{10}\drawdim em \setunitscale 0.15 \linewd 0.4
\move(0 0)\lvec(10 0)\lvec(10 10)\lvec(0 10)\lvec(0 0) \move(10
0)\lvec(15 5)\lvec(15 15)\lvec(5 15)\lvec(0 10) \move(10
10)\lvec(15 15) \htext(3 3){$2$}
\end{texdraw}
}, \hskip 5mm
\raisebox{-0.33\height}[0.69\height][0.35\height]{%
\begin{texdraw}
\fontsize{10}{10}\drawdim em \setunitscale 0.15 \linewd 0.4
\move(0 0)\lvec(10 0)\lvec(10 5)\lvec(0 5)\lvec(0 0) \move(10
0)\lvec(15 5)\lvec(15 10)\lvec(5 10)\lvec(0 5) \move(10 5)\lvec(15
10) \htext(3 1.5){\tiny $3$}
\end{texdraw}%
} .

\vskip 5mm

The walls are built on the ground-state wall

$$Y_{\Lambda_0} = \cdots \raisebox{-0.3\height}{\begin{texdraw} \drawdim
em \setunitscale 0.15 \linewd 0.4 \move(-42 0)\lvec(0 0)\lvec(2.5
2.5)\lvec(2.5 12.5)\lvec(-39.5 12.5) \move(-42 10)\lvec(0 10)
\move(0 0)\lvec(0 10)\lvec(2.5 12.5) \move(-10 0)\lvec(-10
10)\lvec(-7.5 12.5) \move(-20 0)\lvec(-20 10)\lvec(-17.5 12.5)
\move(-30 0)\lvec(-30 10)\lvec(-27.5 12.5) \move(-40 0)\lvec(-40
10)\lvec(-37.5 12.5) \move(0 0)\lvec(-2.5 -2.5)\lvec(-44.5 -2.5)
\move(-10 0)\rlvec(-2.5 -2.5) \move(-20 0)\rlvec(-2.5 -2.5)
\move(-30 0)\rlvec(-2.5 -2.5) \move(-40 0)\rlvec(-2.5 -2.5)
\htext(-7 3){$_1$} \htext(-17 3){$_0$} \htext(-27 3){$_1$}
\htext(-37 3){$_0$}
\end{texdraw}} = \cdots \raisebox{-0.3\height}{\begin{texdraw} \drawdim
em \setunitscale 0.15 \linewd 0.4 \move(-5 0)\lvec(40 0)\lvec(40
10)\lvec(-5 10)\move(0 0)\rlvec(0 10)\move(10 0)\rlvec(0
10)\lvec(0 0)\move(20 0)\rlvec(0 10)\lvec(10 0)\move(30 0)\rlvec(0
10)\lvec(20 0)\move(40 10)\lvec(30 0) \move(1 0.5) \bsegment
\htext(5 1){{\tiny$0$}}\htext(15 1){{\tiny$1$}}\htext(25
1){{\tiny$0$}}\htext(35 1){{\tiny$1$}} \esegment
\end{texdraw}}$$

\vskip 5mm

\noindent following the pattern given below.

\vskip 5mm

\begin{center}
\begin{texdraw}\drawdim
em \setunitscale 0.15 \linewd 0.4

\move(-10 0)\lvec(40 0)\move(-10 10)\lvec(40 10)\move(-10
20)\lvec(40 20)\move(-10 30)\lvec(40 30)\move(-10 40)\lvec(40
40)\move(-10 50)\lvec(40 50)\move(0 0)\lvec(0 55)\move(10
0)\lvec(10 55)\move(20 0)\lvec(20 55)\move(30 0)\lvec(30
55)\move(40 0)\lvec(40 55)\move(-10 25)\lvec(40 25)

\move(10 10)\lvec(0 0)\lvec(10 0)\lvec(10 10)\lfill f:0.8

\move(20 10)\lvec(10 0)\lvec(20 0)\lvec(20 10)\lfill f:0.8

\move(30 10)\lvec(20 0)\lvec(30 0)\lvec(30 10)\lfill f:0.8

\move(40 10)\lvec(30 0)\lvec(40 0)\lvec(40 10)\lfill f:0.8

\move(10 50)\lvec(0 40)

\move(20 50)\lvec(10 40)

\move(30 50)\lvec(20 40)

\move(40 50)\lvec(30 40)

\htext(2 6){\tiny$0$} \htext(16 2){\tiny $0$} \htext(22 6){\tiny
$0$}\htext(36 2){$_0$}

\htext(6 2){\tiny $1$} \htext(12 6){\tiny $1$} \htext(26 2){\tiny
$1$}\htext(32 6){$_1$}

\htext(2 46){\tiny $0$} \htext(16 42){\tiny $0$} \htext(22
46){\tiny $0$}\htext(36 42){$_0$}

\htext(6 42){\tiny $1$} \htext(12 46){\tiny $1$} \htext(26
42){\tiny $1$}\htext(32 46){$_1$}

\move(0.5 0.5) \bsegment \htext(3 13){$_2$}\htext(13 13){$_2$}

\htext(23 13){$_2$}\htext(33 13){$_2$}

\htext(3 33){$_2$}\htext(13 33){$_2$}

\htext(23 33){$_2$}\htext(33 33){$_2$}

\htext(3 20.5){$_3$}\htext(13 20.5){$_3$}

\htext(23 20.5){$_3$}\htext(33 20.5){$_3$}

\htext(3 25.5){$_3$}\htext(13 25.5){$_3$}

\htext(23 25.5){$_3$}\htext(33 25.5){$_3$} \esegment
\end{texdraw}
\end{center}
}\end{ex}

\vskip 5mm

A wall $Y$ built on the ground-state wall $Y_{\Lambda}$ following
the rules given above is called a {\it Young wall on
$Y_{\Lambda}$}, for the heights of its columns are weakly
decreasing as we proceed from right to left. We often write
$Y=(y_k)_{k=0}^{\infty}=(\cdots,y_2,y_1,y_0)$ as an infinite
sequence of its columns.

\vskip 5mm

\begin{defi}{\rm \mbox{}

\begin{itemize}
\item[(1)] A column of a Young wall is called a {\it full
column} if its height is a multiple of the unit length and its top
is of unit thickness.

\item[(2)] For quantum affine algebras of type $A_{2n-1}^{(2)}$
$(n\ge 3)$, $D_n^{(1)}$ $(n\ge 4)$, $A_{2n}^{(2)}$ $(n\ge 1)$,
$D_{n+1}^{(2)}$ $(n\ge 2)$ and $B_n^{(1)}$ $(n\ge 3)$, a Young
wall is said to be {\it proper} if none of the full columns have
the same heights.

\item[(3)] For quantum affine algebras of type $A_n^{(1)}$ $(n\geq 1)$, every
Young wall is defined to be {\it proper}.
\end{itemize}
}
\end{defi}

\vskip 5mm

Let ${\mathcal Z}(\La)$ denote the set of all proper Young walls
on $Y_{\Lambda}$. Then ${\mathcal Z}(\La)$ can be given a
$U_q(\frak{g})$-crystal structure as follows.

\vskip 3mm

\begin{defi}\label {admrmv}
{\rm Let $Y=(y_k)_{k=0}^{\infty}$ be a proper Young wall on
$Y_{\Lambda}$.

\begin{enumerate}
\item[(1)] A block of color $i$ (in short, an $i$-block)
in $Y$ is called a {\it removable $i$-block} if $Y$ remains a
proper Young wall after removing the block. A column in $Y$ is
said to be {\it $i$-removable} if the column has a removable
$i$-block.

\item[(2)] A place in $Y$ is called an {\it admissible $i$-slot} if one
may add an $i$-block to obtain another proper Young wall. A column
in $Y$ is said to be {\it $i$-admissible} if the column has an
admissible $i$-slot.
\end{enumerate}
}
\end{defi}

\vskip 5mm

\begin{ex}
{\rm In the following figure, we consider a proper Young wall for
$\frak{g}=B_3^{(1)}$ built on the ground-state wall
$Y_{\Lambda_0}$ and indicate all the removable blocks and
admissible slots.

\vskip 1cm

\hskip 5mm\raisebox{-0.4\height}
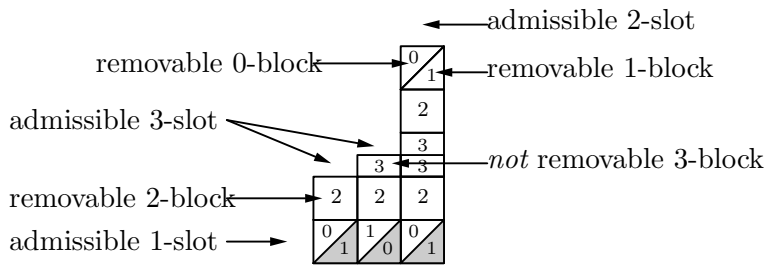
{
\begin{texdraw}
\drawdim em \setunitscale 0.15 \linewd 0.4 \arrowheadtype t:F
\arrowheadsize l:4 w:2

\move(20 0)\lvec(30 0)\lvec(30 10)\lvec(20 10)\lvec(20 0)\htext(22
6){\tiny $0$}

\move(30 0)\lvec(40 0)\lvec(40 10)\lvec(30 10)\lvec(30 0)\htext(32
6){\tiny $1$}

\move(40 0)\lvec(50 0)\lvec(50 10)\lvec(40 10)\lvec(40 0)\htext(42
6){\tiny $0$}
\move(20 0)\lvec(30 10)\lvec(30 0)\lvec(20 0)\lfill f:0.8
\htext(26 2){\tiny $1$}

\move(30 0)\lvec(40 10)\lvec(40 0)\lvec(30 0)\lfill f:0.8
\htext(36 2){\tiny $0$}

\move(40 0)\lvec(50 10)\lvec(50 0)\lvec(40 0)\lfill f:0.8
\htext(46 2){\tiny $1$}
\move(20 10)\lvec(30 10)\lvec(30 20)\lvec(20 20)\lvec(20
10)\htext(24 14){$_2$}

\move(30 10)\lvec(40 10)\lvec(40 20)\lvec(30 20)\lvec(30
10)\htext(34 14){$_2$}

\move(40 10)\lvec(50 10)\lvec(50 20)\lvec(40 20)\lvec(40
10)\htext(44 14){$_2$}
\move(40 20)\lvec(50 20)\lvec(50 30)\lvec(40 30)\lvec(40
20)\htext(44 26){\tiny $3$}

\move(30 20)\lvec(40 20)\lvec(40 25)\lvec(30 25)\lvec(30 20)
\htext(34 21){\tiny $3$}

\move(40 20)\lvec(50 20)\lvec(50 25)\lvec(40 25)\lvec(40 20)
\htext(44 21){\tiny $3$}
\move(40 30)\lvec(50 30)\lvec(50 40)\lvec(40 40)\lvec(40
30)\htext(44 34){$_2$}
\move(40 40)\lvec(50 50)\lvec(40 50)\lvec(40 40) \htext(42
46){\tiny $0$}

\move(40 40)\lvec(50 50)\lvec(50 40)\lvec(40 40) \htext(46
42){\tiny $1$}

\move(0 5) \avec(14 5)

\htext(-50 3){admissible $1$-slot}

\move(0 15) \avec(22 15)

\htext(-50 13){removable $2$-block}

\move(0 33) \avec(24 23)

\move(0 33) \avec(34 27)

\htext(-50 31){admissible $3$-slot}

\move(20 46) \avec(42 46)

\htext(-30 44){removable $0$-block}

\move(60 44) \avec(48 44)

\htext(60 43){removable $1$-block}

\move(60 55) \avec(45 55)

\htext(60 54){admissible $2$-slot}

\move(60 23) \avec(38 23)

\htext(60 22){{\em not} removable $3$-block}
\end{texdraw}}

}
\end{ex}

\vskip 5mm

Fix an index $i\in I$ and let $Y=(y_k)_{k=0}^{\infty}$ be a proper
Young wall on $Y_{\Lambda}$. To each column $y_k$ of $Y$, we
assign its {\it $i$-signature} as follows\,:

\begin{enumerate}
\item we assign $-\,-$ if the column $y_k$ is twice $i$-removable
      (the $i$-block will be of type II);
\item we assign $-$ if the column is once $i$-removable,
      but not $i$-admissible (the $i$-block may be of type I, II, III);
\item we assign $-\,+$ if the column is once $i$-removable and
      once $i$-admissible (the $i$-block will be of type II);
\item we assign $+$ if the column is once $i$-admissible, but not
      $i$-removable (the $i$-block may be of type I, II, III);
\item we assign $+\,+$ if the column is twice $i$-admissible
      (the $i$-block will be of type II).
\end{enumerate}

Then we get a finite sequence of $+$'s and $-$'s for $Y$. From
this sequence, we cancel out every $(+,-)$-pair to obtain a finite
sequence of $-$'s followed by $+$'s, reading from left to right.
This sequence is called the {\it $i$-signature} of $Y$.

\vskip 3mm

We now define the {\it abstract Kashiwara operators}
$\tilde{E}_i$, $\tilde{F}_i$ $(i\in I)$ on ${\mathcal Z}(\Lambda)$
as follows.
\begin{enumerate}
\item We define $\tilde{E}_i Y$ to be the proper Young wall
obtained from $Y$ by removing the $i$-block corresponding to the
rightmost $-$ in the $i$-signature of $Y$. We define $\tilde{E}_i
Y = 0$ if there exists no $-$ in the $i$-signature of $Y$.

\item We define $\tilde{F}_i Y$ to be the proper Young wall obtained from $Y$
by adding an $i$-block to the column corresponding to the leftmost
$+$ in the $i$-signature of $Y$. We define $\tilde{F}_i Y = 0$ if
there exists no $+$ in the $i$-signature of $Y$.
\end{enumerate}
\vskip 5mm

\begin{ex}
{\rm (a) Let $\g = B_3^{(1)}$ and let

\vskip 3mm

\begin{center}
$Y=$\raisebox{-0.5\height}{
\begin{texdraw}
\drawdim em \setunitscale 0.13 \linewd 0.5

\move(10 0)\lvec(20 0)\lvec(20 10)\lvec(10 10)\lvec(10 0)\htext(12
6){\tiny $1$}

\move(20 0)\lvec(30 0)\lvec(30 10)\lvec(20 10)\lvec(20 0)\htext(22
6){\tiny $0$}

\move(30 0)\lvec(40 0)\lvec(40 10)\lvec(30 10)\lvec(30 0)\htext(32
6){\tiny $1$}

\move(40 0)\lvec(50 0)\lvec(50 10)\lvec(40 10)\lvec(40 0)\htext(42
6){\tiny $0$}

\move(50 0)\lvec(60 0)\lvec(60 10)\lvec(50 10)\lvec(50 0)\htext(52
6){\tiny $1$}

\move(60 0)\lvec(70 0)\lvec(70 10)\lvec(60 10)\lvec(60 0)\htext(62
6){\tiny $0$}

\move(70 0)\lvec(80 0)\lvec(80 10)\lvec(70 10)\lvec(70 0)\htext(72
6){\tiny $1$}

\move(80 0)\lvec(90 0)\lvec(90 10)\lvec(80 10)\lvec(80 0)\htext(82
6){\tiny $0$}

%

\move(10 0)\lvec(20 10)\lvec(20 0)\lvec(10 0)\lfill f:0.8
\htext(16 2){\tiny $0$}

\move(20 0)\lvec(30 10)\lvec(30 0)\lvec(20 0)\lfill f:0.8
\htext(26 2){\tiny $1$}

\move(30 0)\lvec(40 10)\lvec(40 0)\lvec(30 0)\lfill f:0.8
\htext(36 2){\tiny $0$}

\move(40 0)\lvec(50 10)\lvec(50 0)\lvec(40 0)\lfill f:0.8
\htext(46 2){\tiny $1$}

\move(50 0)\lvec(60 10)\lvec(60 0)\lvec(50 0)\lfill f:0.8
\htext(56 2){\tiny $0$}

\move(60 0)\lvec(70 10)\lvec(70 0)\lvec(60 0)\lfill f:0.8
\htext(66 2){\tiny $1$}

\move(70 0)\lvec(80 10)\lvec(80 0)\lvec(70 0)\lfill f:0.8
\htext(76 2){\tiny $0$}

\move(80 0)\lvec(90 10)\lvec(90 0)\lvec(80 0)\lfill f:0.8
\htext(86 2){\tiny $1$}

\move(10 10)\lvec(20 10)\lvec(20 20)\lvec(10 20)\lvec(10
10)\htext(13 13){$_2$}

\move(20 10)\lvec(30 10)\lvec(30 20)\lvec(20 20)\lvec(20
10)\htext(23.5 13){$_2$}

\move(30 10)\lvec(40 10)\lvec(40 20)\lvec(30 20)\lvec(30
10)\htext(33.5 13){$_2$}

\move(40 10)\lvec(50 10)\lvec(50 20)\lvec(40 20)\lvec(40
10)\htext(43.5 13){$_2$}

\move(50 10)\lvec(60 10)\lvec(60 20)\lvec(50 20)\lvec(50
10)\htext(53.5 13){$_2$}

\move(60 10)\lvec(70 10)\lvec(70 20)\lvec(60 20)\lvec(60
10)\htext(63.5 13){$_2$}

\move(70 10)\lvec(80 10)\lvec(80 20)\lvec(70 20)\lvec(70
10)\htext(73.5 13){$_2$}

\move(80 10)\lvec(90 10)\lvec(90 20)\lvec(80 20)\lvec(80
10)\htext(83.5 13){$_2$}
\move(40 20)\lvec(50 20)\lvec(50 30)\lvec(40 30)\lvec(40
20)\htext(43.5 26){\tiny $3$}

\move(50 20)\lvec(60 20)\lvec(60 30)\lvec(50 30)\lvec(50
20)\htext(53.5 26){\tiny $3$}

\move(60 20)\lvec(70 20)\lvec(70 30)\lvec(60 30)\lvec(60
20)\htext(63.5 26){\tiny $3$}

\move(70 20)\lvec(80 20)\lvec(80 30)\lvec(70 30)\lvec(70
20)\htext(73.5 26){\tiny $3$}

\move(80 20)\lvec(90 20)\lvec(90 30)\lvec(80 30)\lvec(80
20)\htext(83.5 26){\tiny $3$}

\move(20 20)\lvec(30 20)\lvec(30 25)\lvec(20 25)\lvec(20 20)
\htext(23.5 21){\tiny $3$}

\move(30 20)\lvec(40 20)\lvec(40 25)\lvec(30 25)\lvec(30 20)
\htext(33.5 21){\tiny $3$}

\move(40 20)\lvec(50 20)\lvec(50 25)\lvec(40 25)\lvec(40 20)
\htext(43.5 21){\tiny $3$}

\move(50 20)\lvec(60 20)\lvec(60 25)\lvec(50 25)\lvec(50
20)\htext(53.5 21){\tiny $3$}

\move(60 20)\lvec(70 20)\lvec(70 25)\lvec(60 25)\lvec(60
20)\htext(63.5 21){\tiny $3$}

\move(70 20)\lvec(80 20)\lvec(80 25)\lvec(70 25)\lvec(70
20)\htext(73.5 21){\tiny $3$}

\move(80 20)\lvec(90 20)\lvec(90 25)\lvec(80 25)\lvec(80
20)\htext(83.5 21){\tiny $3$}

\move(50 30)\lvec(60 30)\lvec(60 40)\lvec(50 40)\lvec(50
30)\htext(53.5 33){$_2$}

\move(60 30)\lvec(70 30)\lvec(70 40)\lvec(60 40)\lvec(60
30)\htext(63.5 33){$_2$}

\move(70 30)\lvec(80 30)\lvec(80 40)\lvec(70 40)\lvec(70
30)\htext(73.5 33){$_2$}

\move(80 30)\lvec(90 30)\lvec(90 40)\lvec(80 40)\lvec(80
30)\htext(83.5 33){$_2$}
%


\move(50 40)\lvec(60 40)\lvec(60 50)\lvec(50 50) \lvec(50
40)\htext(52 46){\tiny $1$}

\move(60 40)\lvec(70 40)\lvec(70 50)\lvec(60 50)\lvec(60
40)\htext(62 46){\tiny $0$}

\move(70 40)\lvec(80 40)\lvec(80 50)\lvec(70 50)\lvec(70
40)\htext(72 46){\tiny $1$}

\move(80 40)\lvec(90 40)\lvec(90 50)\lvec(80 50)\lvec(80
40)\htext(82 46){\tiny $0$}

\move(50 40)\lvec(60 50)\lvec(60 40)\lvec(50 40) \htext(56
42){\tiny $0$}

\move(60 40)\lvec(70 50)\lvec(70 40)\lvec(60 40) \htext(66
42){\tiny $1$}

\move(70 40)\lvec(80 50)\lvec(80 40)\lvec(70 40) \htext(76
42){\tiny $0$}

\move(80 40)\lvec(90 50)\lvec(90 40)\lvec(80 40) \htext(86
42){\tiny $1$}
\move(50 50)\lvec(60 50)\lvec(60 60)\lvec(50 60)\lvec(50
50)\htext(53.5 53){$_2$}

\move(60 50)\lvec(70 50)\lvec(70 60)\lvec(60 60)\lvec(60
50)\htext(63.5 53){$_2$}

\move(70 50)\lvec(80 50)\lvec(80 60)\lvec(70 60)\lvec(70
50)\htext(73.5 53){$_2$}

\move(80 50)\lvec(90 50)\lvec(90 60)\lvec(80 60)\lvec(80
50)\htext(83.5 53){$_2$}
\move(60 60)\lvec(70 60)\lvec(70 65)\lvec(60 65)\lvec(60
60)\htext(63.5 61){\tiny $3$}

\move(70 60)\lvec(80 60)\lvec(80 65)\lvec(70 65)\lvec(70
60)\htext(73.5 61){\tiny $3$}

\move(80 60)\lvec(90 60)\lvec(90 65)\lvec(80 65)\lvec(80
60)\htext(83.5 61){\tiny $3$}
\move(60 60)\lvec(70 60)\lvec(70 70)\lvec(60 70)\lvec(60
60)\htext(63.5 66){\tiny $3$}

\move(70 60)\lvec(80 60)\lvec(80 70)\lvec(70 70)\lvec(70
60)\htext(73.5 66){\tiny $3$}

\move(80 60)\lvec(90 60)\lvec(90 70)\lvec(80 70)\lvec(80
60)\htext(83.5 66){\tiny $3$}
\move(60 70)\lvec(70 70)\lvec(70 80)\lvec(60 80)\lvec(60
70)\htext(63.5 73){$_2$}

\move(70 70)\lvec(80 70)\lvec(80 80)\lvec(70 80)\lvec(70
70)\htext(73.5 73){$_2$}

\move(80 70)\lvec(90 70)\lvec(90 80)\lvec(80 80)\lvec(80
70)\htext(83.5 73){$_2$}
\move(70 80)\lvec(80 90)\lvec(80 80)\lvec(70 80)\htext(76
82){\tiny $0$}

\move(80 80)\lvec(90 90)\lvec(90 80)\lvec(80 80)\htext(86
82){\tiny $1$}

\htext(-25 -9){$_{\cdot}$}\htext(-15 -9){$_{\cdot}$}\htext(-5
-9){$_{\cdot}$}

\htext(5 -9){$_{\cdot}$}\htext(12 -10){$_+$}\htext(25
-9){$_{\cdot}$}\htext(35 -9){$_{\cdot}$}\htext(41 -10){$_-$}

\htext(50 -10){$_{++}$}\htext(65 -9){$_{\cdot}$}\htext(75
-9){$_{\cdot}$}\htext(85 -9){$_{\cdot}$}
\end{texdraw}}\quad .
\end{center}

\vskip 5mm

If $i=3$, we first get the sequence $(\cdots, +, \cdot, \cdot, -,
++, \cdot, \cdot, \cdot)$. After cancelling out $(+,-)$-pairs, we
obtain the $3$-signature $(+,+)$ of $Y$.

Therefore, we have $\tilde{E}_3Y=0$ and

\vskip 3mm

\begin{center}
$\tilde{F}_3 Y=$\raisebox{-0.5\height}{
\begin{texdraw}
\drawdim em \setunitscale 0.13 \linewd 0.5

\move(10 0)\lvec(20 0)\lvec(20 10)\lvec(10 10)\lvec(10 0)\htext(12
6){\tiny $1$}

\move(20 0)\lvec(30 0)\lvec(30 10)\lvec(20 10)\lvec(20 0)\htext(22
6){\tiny $0$}

\move(30 0)\lvec(40 0)\lvec(40 10)\lvec(30 10)\lvec(30 0)\htext(32
6){\tiny $1$}

\move(40 0)\lvec(50 0)\lvec(50 10)\lvec(40 10)\lvec(40 0)\htext(42
6){\tiny $0$}

\move(50 0)\lvec(60 0)\lvec(60 10)\lvec(50 10)\lvec(50 0)\htext(52
6){\tiny $1$}

\move(60 0)\lvec(70 0)\lvec(70 10)\lvec(60 10)\lvec(60 0)\htext(62
6){\tiny $0$}

\move(70 0)\lvec(80 0)\lvec(80 10)\lvec(70 10)\lvec(70 0)\htext(72
6){\tiny $1$}

\move(80 0)\lvec(90 0)\lvec(90 10)\lvec(80 10)\lvec(80 0)\htext(82
6){\tiny $0$}

%

\move(10 0)\lvec(20 10)\lvec(20 0)\lvec(10 0)\lfill f:0.8
\htext(16 2){\tiny $0$}

\move(20 0)\lvec(30 10)\lvec(30 0)\lvec(20 0)\lfill f:0.8
\htext(26 2){\tiny $1$}

\move(30 0)\lvec(40 10)\lvec(40 0)\lvec(30 0)\lfill f:0.8
\htext(36 2){\tiny $0$}

\move(40 0)\lvec(50 10)\lvec(50 0)\lvec(40 0)\lfill f:0.8
\htext(46 2){\tiny $1$}

\move(50 0)\lvec(60 10)\lvec(60 0)\lvec(50 0)\lfill f:0.8
\htext(56 2){\tiny $0$}

\move(60 0)\lvec(70 10)\lvec(70 0)\lvec(60 0)\lfill f:0.8
\htext(66 2){\tiny $1$}

\move(70 0)\lvec(80 10)\lvec(80 0)\lvec(70 0)\lfill f:0.8
\htext(76 2){\tiny $0$}

\move(80 0)\lvec(90 10)\lvec(90 0)\lvec(80 0)\lfill f:0.8
\htext(86 2){\tiny $1$}

\move(10 10)\lvec(20 10)\lvec(20 20)\lvec(10 20)\lvec(10
10)\htext(13 13){$_2$}

\move(20 10)\lvec(30 10)\lvec(30 20)\lvec(20 20)\lvec(20
10)\htext(23.5 13){$_2$}

\move(30 10)\lvec(40 10)\lvec(40 20)\lvec(30 20)\lvec(30
10)\htext(33.5 13){$_2$}

\move(40 10)\lvec(50 10)\lvec(50 20)\lvec(40 20)\lvec(40
10)\htext(43.5 13){$_2$}

\move(50 10)\lvec(60 10)\lvec(60 20)\lvec(50 20)\lvec(50
10)\htext(53.5 13){$_2$}

\move(60 10)\lvec(70 10)\lvec(70 20)\lvec(60 20)\lvec(60
10)\htext(63.5 13){$_2$}

\move(70 10)\lvec(80 10)\lvec(80 20)\lvec(70 20)\lvec(70
10)\htext(73.5 13){$_2$}

\move(80 10)\lvec(90 10)\lvec(90 20)\lvec(80 20)\lvec(80
10)\htext(83.5 13){$_2$}
\move(40 20)\lvec(50 20)\lvec(50 30)\lvec(40 30)\lvec(40
20)\htext(43.5 26){\tiny $3$}

\move(50 20)\lvec(60 20)\lvec(60 30)\lvec(50 30)\lvec(50
20)\htext(53.5 26){\tiny $3$}

\move(60 20)\lvec(70 20)\lvec(70 30)\lvec(60 30)\lvec(60
20)\htext(63.5 26){\tiny $3$}

\move(70 20)\lvec(80 20)\lvec(80 30)\lvec(70 30)\lvec(70
20)\htext(73.5 26){\tiny $3$}

\move(80 20)\lvec(90 20)\lvec(90 30)\lvec(80 30)\lvec(80
20)\htext(83.5 26){\tiny $3$}

\move(20 20)\lvec(30 20)\lvec(30 25)\lvec(20 25)\lvec(20 20)
\htext(23.5 21){\tiny $3$}

\move(30 20)\lvec(40 20)\lvec(40 25)\lvec(30 25)\lvec(30 20)
\htext(33.5 21){\tiny $3$}

\move(40 20)\lvec(50 20)\lvec(50 25)\lvec(40 25)\lvec(40 20)
\htext(43.5 21){\tiny $3$}

\move(50 20)\lvec(60 20)\lvec(60 25)\lvec(50 25)\lvec(50
20)\htext(53.5 21){\tiny $3$}

\move(60 20)\lvec(70 20)\lvec(70 25)\lvec(60 25)\lvec(60
20)\htext(63.5 21){\tiny $3$}

\move(70 20)\lvec(80 20)\lvec(80 25)\lvec(70 25)\lvec(70
20)\htext(73.5 21){\tiny $3$}

\move(80 20)\lvec(90 20)\lvec(90 25)\lvec(80 25)\lvec(80
20)\htext(83.5 21){\tiny $3$}

\move(50 30)\lvec(60 30)\lvec(60 40)\lvec(50 40)\lvec(50
30)\htext(53.5 33){$_2$}

\move(60 30)\lvec(70 30)\lvec(70 40)\lvec(60 40)\lvec(60
30)\htext(63.5 33){$_2$}

\move(70 30)\lvec(80 30)\lvec(80 40)\lvec(70 40)\lvec(70
30)\htext(73.5 33){$_2$}

\move(80 30)\lvec(90 30)\lvec(90 40)\lvec(80 40)\lvec(80
30)\htext(83.5 33){$_2$}
%


\move(50 40)\lvec(60 40)\lvec(60 50)\lvec(50 50) \lvec(50
40)\htext(52 46){\tiny $1$}

\move(60 40)\lvec(70 40)\lvec(70 50)\lvec(60 50)\lvec(60
40)\htext(62 46){\tiny $0$}

\move(70 40)\lvec(80 40)\lvec(80 50)\lvec(70 50)\lvec(70
40)\htext(72 46){\tiny $1$}

\move(80 40)\lvec(90 40)\lvec(90 50)\lvec(80 50)\lvec(80
40)\htext(82 46){\tiny $0$}

\move(50 40)\lvec(60 50)\lvec(60 40)\lvec(50 40) \htext(56
42){\tiny $0$}

\move(60 40)\lvec(70 50)\lvec(70 40)\lvec(60 40) \htext(66
42){\tiny $1$}

\move(70 40)\lvec(80 50)\lvec(80 40)\lvec(70 40) \htext(76
42){\tiny $0$}

\move(80 40)\lvec(90 50)\lvec(90 40)\lvec(80 40) \htext(86
42){\tiny $1$}
\move(50 50)\lvec(60 50)\lvec(60 60)\lvec(50 60)\lvec(50
50)\htext(53.5 53){$_2$}

\move(60 50)\lvec(70 50)\lvec(70 60)\lvec(60 60)\lvec(60
50)\htext(63.5 53){$_2$}

\move(70 50)\lvec(80 50)\lvec(80 60)\lvec(70 60)\lvec(70
50)\htext(73.5 53){$_2$}

\move(80 50)\lvec(90 50)\lvec(90 60)\lvec(80 60)\lvec(80
50)\htext(83.5 53){$_2$}
\move(55 65)\lcir r:6

\move(60 60)\lvec(70 60)\lvec(70 65)\lvec(60 65)\lvec(60
60)\htext(63.5 61){\tiny $3$}

\move(50 60)\lvec(60 60)\lvec(60 65)\lvec(50 65)\lvec(50
60)\htext(53.5 61){\tiny $3$}

\move(70 60)\lvec(80 60)\lvec(80 65)\lvec(70 65)\lvec(70
60)\htext(73.5 61){\tiny $3$}

\move(80 60)\lvec(90 60)\lvec(90 65)\lvec(80 65)\lvec(80
60)\htext(83.5 61){\tiny $3$}
\move(60 60)\lvec(70 60)\lvec(70 70)\lvec(60 70)\lvec(60
60)\htext(63.5 66){\tiny $3$}

\move(70 60)\lvec(80 60)\lvec(80 70)\lvec(70 70)\lvec(70
60)\htext(73.5 66){\tiny $3$}

\move(80 60)\lvec(90 60)\lvec(90 70)\lvec(80 70)\lvec(80
60)\htext(83.5 66){\tiny $3$}
\move(60 70)\lvec(70 70)\lvec(70 80)\lvec(60 80)\lvec(60
70)\htext(63.5 73){$_2$}

\move(70 70)\lvec(80 70)\lvec(80 80)\lvec(70 80)\lvec(70
70)\htext(73.5 73){$_2$}

\move(80 70)\lvec(90 70)\lvec(90 80)\lvec(80 80)\lvec(80
70)\htext(83.5 73){$_2$}
\move(70 80)\lvec(80 90)\lvec(80 80)\lvec(70 80)\htext(76
82){\tiny $0$}

\move(80 80)\lvec(90 90)\lvec(90 80)\lvec(80 80)\htext(86
82){\tiny $1$}

\end{texdraw}}\quad .
\end{center}

\vskip 5mm

(b) Let $\g=A_{4}^{(2)}$ and let \vskip 5mm

\begin{center}
$Y=$ \raisebox{-0.5\height}{
\begin{texdraw}
\drawdim em \setunitscale 0.13 \linewd 0.5

\move(-10 0)\lvec(0 0)\lvec(0 10)\lvec(-10 10)\lvec(-10
0)\htext(-7 1){\tiny $0$}

\move(0 0)\lvec(10 0)\lvec(10 10)\lvec(0 10)\lvec(0 0)\htext(3
1){\tiny $0$}

\move(10 0)\lvec(20 0)\lvec(20 10)\lvec(10 10)\lvec(10 0)\htext(13
1){\tiny $0$}

\move(20 0)\lvec(30 0)\lvec(30 10)\lvec(20 10)\lvec(20 0)\htext(23
1){\tiny $0$}

\move(30 0)\lvec(40 0)\lvec(40 10)\lvec(30 10)\lvec(30 0)\htext(33
1){\tiny $0$}
\move(-10 0)\lvec(0 0)\lvec(0 5)\lvec(-10 5)\lvec(-10 0)\lfill
f:0.8 \htext(-7 6){\tiny $0$}

\move(0 0)\lvec(10 0)\lvec(10 5)\lvec(0 5)\lvec(0 0)\lfill f:0.8
\htext(3 6){\tiny $0$}

\move(10 0)\lvec(20 0)\lvec(20 5)\lvec(10 5)\lvec(10 0)\lfill
f:0.8 \htext(13 6){\tiny $0$}

\move(20 0)\lvec(30 0)\lvec(30 5)\lvec(20 5)\lvec(20 0)\lfill
f:0.8 \htext(23 6){\tiny $0$}

\move(30 0)\lvec(40 0)\lvec(40 5)\lvec(30 5)\lvec(30 0)\lfill
f:0.8 \htext(33 6){\tiny $0$}
\move(0 10)\lvec(10 10)\lvec(10 20)\lvec(0 20)\lvec(0 10)\htext(3
13){$_1$}

\move(10 10)\lvec(20 10)\lvec(20 20)\lvec(10 20)\lvec(10
10)\htext(13 13){$_1$}

\move(20 10)\lvec(30 10)\lvec(30 20)\lvec(20 20)\lvec(20
10)\htext(23 13){$_1$}

\move(30 10)\lvec(40 10)\lvec(40 20)\lvec(30 20)\lvec(30
10)\htext(33 13){$_1$}
\move(0 20)\lvec(10 20)\lvec(10 30)\lvec(0 30)\lvec(0 20)\htext(3
23){$_2$}

\move(10 20)\lvec(20 20)\lvec(20 30)\lvec(10 30)\lvec(10
20)\htext(13 23){$_2$}

\move(20 20)\lvec(30 20)\lvec(30 30)\lvec(20 30)\lvec(20
20)\htext(23 23){$_2$}

\move(30 20)\lvec(40 20)\lvec(40 30)\lvec(30 30)\lvec(30
20)\htext(33 23){$_2$}
\move(0 30)\lvec(10 30)\lvec(10 40)\lvec(0 40)\lvec(0 30)\htext(3
33){$_1$}

\move(10 30)\lvec(20 30)\lvec(20 40)\lvec(10 40)\lvec(10
30)\htext(13 33){$_1$}

\move(20 30)\lvec(30 30)\lvec(30 40)\lvec(20 40)\lvec(20
30)\htext(23 33){$_1$}

\move(30 30)\lvec(40 30)\lvec(40 40)\lvec(30 40)\lvec(30
30)\htext(33 33){$_1$}
\move(0 40)\lvec(10 40)\lvec(10 45)\lvec(0 45)\lvec(0 40)\htext(3
41){\tiny $0$}

\move(10 40)\lvec(20 40)\lvec(20 45)\lvec(10 45)\lvec(10
40)\htext(13 41){\tiny $0$}

\move(20 40)\lvec(30 40)\lvec(30 45)\lvec(20 45)\lvec(20
40)\htext(23 41){\tiny $0$}

\move(30 40)\lvec(40 40)\lvec(40 45)\lvec(30 45)\lvec(30
40)\htext(33 41){\tiny $0$}
\move(30 45)\lvec(40 45)\lvec(40 50)\lvec(30 50)\lvec(30
45)\htext(33 46){\tiny $0$}
\move(30 50)\lvec(40 50)\lvec(40 60)\lvec(30 60)\lvec(30
50)\htext(33 53){$_1$}
\move(30 60)\lvec(40 60)\lvec(40 70)\lvec(30 70)\lvec(30
60)\htext(33 63){$_2$}
\move(30 70)\lvec(40 70)\lvec(40 80)\lvec(30 80)\lvec(30
70)\htext(33 73){$_1$}
\move(30 80)\lvec(40 80)\lvec(40 85)\lvec(30 85)\lvec(30
80)\htext(33 81){\tiny $0$}
\move(2 0)\bsegment \htext(-40 -10){$_{\cdot}$}\htext(-30
-10){$_{\cdot}$}\htext(-20 -10){$_{\cdot}$}\htext(-10
-10){$_-$}\htext(0 -10){$_-$}\htext(12 -10){$_{\cdot}$}\htext(20
-10){$_+$}\htext(30 -10){$_{-+}$}\esegment
\end{texdraw}}\quad .
\end{center}
\vskip 5mm

If $i=0$, we fist get the sequence $(\cdots, -,-,\cdot, +, -+)$.
After cancelling out $(+,-)$-pairs, we obtain the $0$-signature
$(-,-,+)$ of $Y$.

Therefore, we have \vskip 5mm

\begin{center}
$\tilde{E}_0 Y=$\raisebox{-0.5\height}{
\begin{texdraw}
\drawdim em \setunitscale 0.13 \linewd 0.5

\move(-10 0)\lvec(0 0)\lvec(0 10)\lvec(-10 10)\lvec(-10
0)\htext(-7 1){\tiny $0$}

\move(0 0)\lvec(10 0)\lvec(10 10)\lvec(0 10)\lvec(0 0)\htext(3
1){\tiny $0$}

\move(10 0)\lvec(20 0)\lvec(20 10)\lvec(10 10)\lvec(10 0)\htext(13
1){\tiny $0$}

\move(20 0)\lvec(30 0)\lvec(30 10)\lvec(20 10)\lvec(20 0)\htext(23
1){\tiny $0$}

\move(30 0)\lvec(40 0)\lvec(40 10)\lvec(30 10)\lvec(30 0)\htext(33
1){\tiny $0$}
\move(-10 0)\lvec(0 0)\lvec(0 5)\lvec(-10 5)\lvec(-10 0)\lfill
f:0.8 \htext(-7 6){\tiny $0$}

\move(0 0)\lvec(10 0)\lvec(10 5)\lvec(0 5)\lvec(0 0)\lfill f:0.8
\htext(3 6){\tiny $0$}

\move(10 0)\lvec(20 0)\lvec(20 5)\lvec(10 5)\lvec(10 0)\lfill
f:0.8 \htext(13 6){\tiny $0$}

\move(20 0)\lvec(30 0)\lvec(30 5)\lvec(20 5)\lvec(20 0)\lfill
f:0.8 \htext(23 6){\tiny $0$}

\move(30 0)\lvec(40 0)\lvec(40 5)\lvec(30 5)\lvec(30 0)\lfill
f:0.8 \htext(33 6){\tiny $0$}
\move(0 10)\lvec(10 10)\lvec(10 20)\lvec(0 20)\lvec(0 10)\htext(3
13){$_1$}

\move(10 10)\lvec(20 10)\lvec(20 20)\lvec(10 20)\lvec(10
10)\htext(13 13){$_1$}

\move(20 10)\lvec(30 10)\lvec(30 20)\lvec(20 20)\lvec(20
10)\htext(23 13){$_1$}

\move(30 10)\lvec(40 10)\lvec(40 20)\lvec(30 20)\lvec(30
10)\htext(33 13){$_1$}
\move(0 20)\lvec(10 20)\lvec(10 30)\lvec(0 30)\lvec(0 20)\htext(3
23){$_2$}

\move(10 20)\lvec(20 20)\lvec(20 30)\lvec(10 30)\lvec(10
20)\htext(13 23){$_2$}

\move(20 20)\lvec(30 20)\lvec(30 30)\lvec(20 30)\lvec(20
20)\htext(23 23){$_2$}

\move(30 20)\lvec(40 20)\lvec(40 30)\lvec(30 30)\lvec(30
20)\htext(33 23){$_2$}
\move(0 30)\lvec(10 30)\lvec(10 40)\lvec(0 40)\lvec(0 30)\htext(3
33){$_1$}

\move(10 30)\lvec(20 30)\lvec(20 40)\lvec(10 40)\lvec(10
30)\htext(13 33){$_1$}

\move(20 30)\lvec(30 30)\lvec(30 40)\lvec(20 40)\lvec(20
30)\htext(23 33){$_1$}

\move(30 30)\lvec(40 30)\lvec(40 40)\lvec(30 40)\lvec(30
30)\htext(33 33){$_1$}
\move(5 45)\lcir r:6

\move(10 40)\lvec(20 40)\lvec(20 45)\lvec(10 45)\lvec(10
40)\htext(13 41){\tiny $0$}

\move(20 40)\lvec(30 40)\lvec(30 45)\lvec(20 45)\lvec(20
40)\htext(23 41){\tiny $0$}

\move(30 40)\lvec(40 40)\lvec(40 45)\lvec(30 45)\lvec(30
40)\htext(33 41){\tiny $0$}
\move(30 45)\lvec(40 45)\lvec(40 50)\lvec(30 50)\lvec(30
45)\htext(33 46){\tiny $0$}
\move(30 50)\lvec(40 50)\lvec(40 60)\lvec(30 60)\lvec(30
50)\htext(33 53){$_1$}
\move(30 60)\lvec(40 60)\lvec(40 70)\lvec(30 70)\lvec(30
60)\htext(33 63){$_2$}
\move(30 70)\lvec(40 70)\lvec(40 80)\lvec(30 80)\lvec(30
70)\htext(33 73){$_1$}
\move(30 80)\lvec(40 80)\lvec(40 85)\lvec(30 85)\lvec(30
80)\htext(33 81){\tiny $0$}
\end{texdraw}}\quad , \quad
$\tilde{F}_0 Y=$ \raisebox{-0.5\height}{
\begin{texdraw}
\drawdim em \setunitscale 0.13 \linewd 0.5

\move(-10 0)\lvec(0 0)\lvec(0 10)\lvec(-10 10)\lvec(-10
0)\htext(-7 1){\tiny $0$}

\move(0 0)\lvec(10 0)\lvec(10 10)\lvec(0 10)\lvec(0 0)\htext(3
1){\tiny $0$}

\move(10 0)\lvec(20 0)\lvec(20 10)\lvec(10 10)\lvec(10 0)\htext(13
1){\tiny $0$}

\move(20 0)\lvec(30 0)\lvec(30 10)\lvec(20 10)\lvec(20 0)\htext(23
1){\tiny $0$}

\move(30 0)\lvec(40 0)\lvec(40 10)\lvec(30 10)\lvec(30 0)\htext(33
1){\tiny $0$}
\move(-10 0)\lvec(0 0)\lvec(0 5)\lvec(-10 5)\lvec(-10 0)\lfill
f:0.8 \htext(-7 6){\tiny $0$}

\move(0 0)\lvec(10 0)\lvec(10 5)\lvec(0 5)\lvec(0 0)\lfill f:0.8
\htext(3 6){\tiny $0$}

\move(10 0)\lvec(20 0)\lvec(20 5)\lvec(10 5)\lvec(10 0)\lfill
f:0.8 \htext(13 6){\tiny $0$}

\move(20 0)\lvec(30 0)\lvec(30 5)\lvec(20 5)\lvec(20 0)\lfill
f:0.8 \htext(23 6){\tiny $0$}

\move(30 0)\lvec(40 0)\lvec(40 5)\lvec(30 5)\lvec(30 0)\lfill
f:0.8 \htext(33 6){\tiny $0$}
\move(0 10)\lvec(10 10)\lvec(10 20)\lvec(0 20)\lvec(0 10)\htext(3
13){$_1$}

\move(10 10)\lvec(20 10)\lvec(20 20)\lvec(10 20)\lvec(10
10)\htext(13 13){$_1$}

\move(20 10)\lvec(30 10)\lvec(30 20)\lvec(20 20)\lvec(20
10)\htext(23 13){$_1$}

\move(30 10)\lvec(40 10)\lvec(40 20)\lvec(30 20)\lvec(30
10)\htext(33 13){$_1$}
\move(0 20)\lvec(10 20)\lvec(10 30)\lvec(0 30)\lvec(0 20)\htext(3
23){$_2$}

\move(10 20)\lvec(20 20)\lvec(20 30)\lvec(10 30)\lvec(10
20)\htext(13 23){$_2$}

\move(20 20)\lvec(30 20)\lvec(30 30)\lvec(20 30)\lvec(20
20)\htext(23 23){$_2$}

\move(30 20)\lvec(40 20)\lvec(40 30)\lvec(30 30)\lvec(30
20)\htext(33 23){$_2$}
\move(0 30)\lvec(10 30)\lvec(10 40)\lvec(0 40)\lvec(0 30)\htext(3
33){$_1$}

\move(10 30)\lvec(20 30)\lvec(20 40)\lvec(10 40)\lvec(10
30)\htext(13 33){$_1$}

\move(20 30)\lvec(30 30)\lvec(30 40)\lvec(20 40)\lvec(20
30)\htext(23 33){$_1$}

\move(30 30)\lvec(40 30)\lvec(40 40)\lvec(30 40)\lvec(30
30)\htext(33 33){$_1$}
\move(0 40)\lvec(10 40)\lvec(10 45)\lvec(0 45)\lvec(0 40)\htext(3
41){\tiny $0$}

\move(10 40)\lvec(20 40)\lvec(20 45)\lvec(10 45)\lvec(10
40)\htext(13 41){\tiny $0$}

\move(20 40)\lvec(30 40)\lvec(30 45)\lvec(20 45)\lvec(20
40)\htext(23 41){\tiny $0$}

\move(30 40)\lvec(40 40)\lvec(40 45)\lvec(30 45)\lvec(30
40)\htext(33 41){\tiny $0$}
\move(30 45)\lvec(40 45)\lvec(40 50)\lvec(30 50)\lvec(30
45)\htext(33 46){\tiny $0$}
\move(30 50)\lvec(40 50)\lvec(40 60)\lvec(30 60)\lvec(30
50)\htext(33 53){$_1$}
\move(30 60)\lvec(40 60)\lvec(40 70)\lvec(30 70)\lvec(30
60)\htext(33 63){$_2$}
\move(30 70)\lvec(40 70)\lvec(40 80)\lvec(30 80)\lvec(30
70)\htext(33 73){$_1$}
\move(35 85)\lcir r:6

\move(30 80)\lvec(40 80)\lvec(40 85)\lvec(30 85)\lvec(30
80)\htext(33 81){\tiny $0$}
\move(30 85)\lvec(40 85)\lvec(40 90)\lvec(30 90)\lvec(30
85)\htext(33 86){\tiny $0$}
\end{texdraw}}\quad .
\end{center}

}
\end{ex}

\vskip 5mm

Next, we define the maps
$$ \wt : {\mathcal Z}(\La) \longrightarrow P, \quad
\varepsilon_i : {\mathcal Z}(\La) \longrightarrow \Z, \quad
\varphi_i : {\mathcal Z}(\La) \longrightarrow \Z$$ by
\begin{equation}
\begin{aligned}\text{}
{\rm wt}(Y)&=\Lambda-\sum_{i\in I}k_i\alpha_i, \\
\varepsilon_i(Y)&=\text{the number of $-$'s in the $i$-signature of $Y$}, \\
\varphi_i(Y)&=\text{the number of $+$'s in the $i$-signature of
$Y$},
\end{aligned}
\end{equation}
where $k_i$ is the number of $i$-blocks in $Y$ that have been
added to $Y_{\Lambda}$.

\vskip 3mm

Then we have :

\begin{thm}\label {Kang1}{\rm (\cite{Ka2000})}
The set ${\mathcal Z}(\Lambda)$ together with the maps $\wt:
{\mathcal Z}(\Lambda) \longrightarrow  P$, $\tilde{E}_i,
\tilde{F}_i : {\mathcal Z}(\Lambda) \longrightarrow {\mathcal
Z}(\Lambda) \cup \{0\}$, and $\varepsilon_i, \varphi_i : {\mathcal
Z}(\Lambda) \longrightarrow \Z$ {\rm ($i\in I$)} becomes a
$U_q(\g)$-crystal.  \qed
\end{thm}

\vskip 3mm

Let $\delta=d_0\alpha_0+\cdots +d_n\alpha_n$ be the null root of
$U_q(\frak{g})$, and set
\begin{equation}
\begin{aligned}\text{}
a_i = \begin{cases} d_i \quad & \text{if} \ \ \g \neq D_{n+1}^{(2)}, \\
2 d_i \quad & \text{if} \ \ \g = D_{n+1}^{(2)}.
\end{cases}
\end{aligned}
\end{equation}

\vskip 3mm
\begin{defi}{\rm \mbox{}

\begin{itemize}
\item[(1)] The part of a column in a proper Young wall
is called a {\it $\delta$-column} if it contains $a_{0}$-many
$0$-blocks, $a_{1}$-many $1$-blocks, $\cdots$, $a_{n}$-many
$n$-blocks in some cyclic order.

\item[(2)] A $\delta$-column in a proper Young wall is
called {\it removable} if it can be removed from the top to yield
another proper Young wall.

\item[(3)] A proper Young wall is said to be {\it reduced}
if none of its columns contain a removable $\delta$-column.
\end{itemize}
}
\end{defi}

\vskip 5mm

\begin{ex}{\rm
(a) The following are $\delta$-columns for $\g = B_3^{(1)}$.

\vskip 5mm

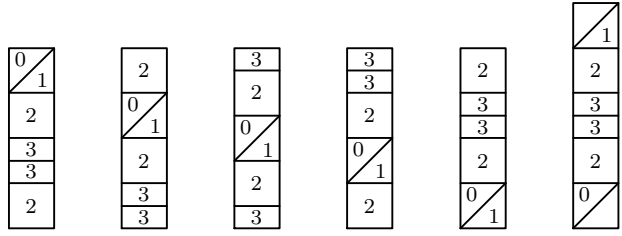
\begin{center}
\begin{texdraw}
\textref h:C v:C \fontsize{8}{8}\selectfont \drawdim mm \linewd
0.25 \setunitscale 6  \move(0 0) \bsegment \move(0 0)\lvec(1
0)\lvec(1 4)\lvec(0 4)\lvec(0 0) \move(0 1)\lvec(1 1) \move(0
1.5)\lvec(1 1.5) \move(0 2)\lvec(1 2) \move(0 3)\lvec(1 3) \move(0
3)\lvec(1 4) \htext(0.5 0.5){$2$} \htext(0.5 1.25){$3$} \htext(0.5
1.75){$3$} \htext(0.5 2.5){$2$} \htext(0.25 3.75){$0$} \htext(0.75
3.27){$1$} \esegment \move(2.5 0) \bsegment \move(0 0)\lvec(1
0)\lvec(1 4)\lvec(0 4)\lvec(0 0) \move(0 0.5)\lvec(1 0.5) \move(0
1)\lvec(1 1) \move(0 2)\lvec(1 2) \move(0 3)\lvec(1 3) \move(0
2)\lvec(1 3) \htext(0.5 3.5){$2$} \htext(0.5 0.25){$3$} \htext(0.5
0.75){$3$} \htext(0.5 1.5){$2$} \htext(0.25 2.75){$0$} \htext(0.75
2.27){$1$} \esegment \move(5 0) \bsegment \move(0 0)\lvec(1
0)\lvec(1 4)\lvec(0 4)\lvec(0 0) \move(0 0.5)\lvec(1 0.5) \move(0
1.5)\lvec(1 1.5) \move(0 2.5)\lvec(1 2.5) \move(0 3.5)\lvec(1 3.5)
\move(0 1.5)\lvec(1 2.5) \htext(0.5 1){$2$} \htext(0.5 0.25){$3$}
\htext(0.5 3.75){$3$} \htext(0.5 3){$2$} \htext(0.75 1.72){$1$}
\htext(0.25 2.25){$0$} \esegment \move(7.5 0) \bsegment \move(0
0)\lvec(1 0)\lvec(1 4)\lvec(0 4)\lvec(0 0) \move(0 3.5)\lvec(1
3.5) \move(0 1)\lvec(1 1) \move(0 2)\lvec(1 2) \move(0 3)\lvec(1
3) \move(0 1)\lvec(1 2) \htext(0.5 2.5){$2$} \htext(0.5 3.25){$3$}
\htext(0.5 3.75){$3$} \htext(0.5 0.5){$2$} \htext(0.25 1.75){$0$}
\htext(0.75 1.27){$1$} \esegment \move(10 0) \bsegment \move(0
0)\lvec(1 0)\lvec(1 4)\lvec(0 4)\lvec(0 0) \move(0 2.5)\lvec(1
2.5) \move(0 1)\lvec(1 1) \move(0 2)\lvec(1 2) \move(0 3)\lvec(1
3) \move(0 0)\lvec(1 1) \htext(0.5 3.5){$2$} \htext(0.5 2.25){$3$}
\htext(0.5 2.75){$3$} \htext(0.5 1.5){$2$} \htext(0.25 0.75){$0$}
\htext(0.75 0.27){$1$} \esegment \move(12.5 0) \bsegment \move(0
0)\lvec(1 0)\lvec(1 4)\lvec(0 4)\lvec(0 0) \move(0 2.5)\lvec(1
2.5) \move(0 1)\lvec(1 1) \move(0 2)\lvec(1 2) \move(0 3)\lvec(1
3) \move(0 0)\lvec(1 1) \htext(0.5 3.5){$2$} \htext(0.5 2.25){$3$}
\htext(0.5 2.75){$3$} \htext(0.5 1.5){$2$} \htext(0.25 0.75){$0$}
\htext(0.75 4.27){$1$} \move(1 4)\lvec(1 5)\lvec(0 5)\lvec(0
4)\lvec(1 5) \esegment
\end{texdraw}
\end{center}

(b) Consider the following proper Young walls for $\g= B_3^{(1)}$.
The first one is reduced, but the second one is not. Note that the
second Young wall contains a removable $\delta$-column.

\vskip 5mm

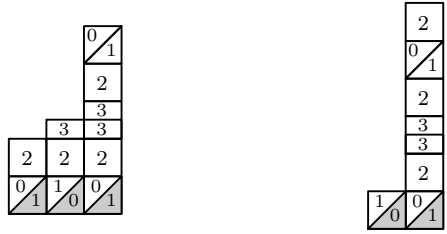
\begin{center}
\raisebox{-0.4\height}{
\begin{texdraw}
\drawdim em \setunitscale 0.13 \linewd 0.5

\move(20 0)\lvec(30 0)\lvec(30 10)\lvec(20 10)\lvec(20 0)\htext(22
6){\tiny $0$}

\move(30 0)\lvec(40 0)\lvec(40 10)\lvec(30 10)\lvec(30 0)\htext(32
6){\tiny $1$}

\move(40 0)\lvec(50 0)\lvec(50 10)\lvec(40 10)\lvec(40 0)\htext(42
6){\tiny $0$}
\move(20 0)\lvec(30 10)\lvec(30 0)\lvec(20 0)\lfill f:0.8
\htext(26 2){\tiny $1$}

\move(30 0)\lvec(40 10)\lvec(40 0)\lvec(30 0)\lfill f:0.8
\htext(36 2){\tiny $0$}

\move(40 0)\lvec(50 10)\lvec(50 0)\lvec(40 0)\lfill f:0.8
\htext(46 2){\tiny $1$}
\move(20 10)\lvec(30 10)\lvec(30 20)\lvec(20 20)\lvec(20
10)\htext(23.5 13){$_2$}

\move(30 10)\lvec(40 10)\lvec(40 20)\lvec(30 20)\lvec(30
10)\htext(33.5 13){$_2$}

\move(40 10)\lvec(50 10)\lvec(50 20)\lvec(40 20)\lvec(40
10)\htext(43.5 13){$_2$}
\move(40 20)\lvec(50 20)\lvec(50 30)\lvec(40 30)\lvec(40
20)\htext(43.5 26){\tiny $3$}

\move(30 20)\lvec(40 20)\lvec(40 25)\lvec(30 25)\lvec(30 20)
\htext(33.5 21){\tiny $3$}

\move(40 20)\lvec(50 20)\lvec(50 25)\lvec(40 25)\lvec(40 20)
\htext(43.5 21){\tiny $3$}
\move(40 30)\lvec(50 30)\lvec(50 40)\lvec(40 40)\lvec(40
30)\htext(43.5 33){$_2$}
\move(40 40)\lvec(50 50)\lvec(40 50)\lvec(40 40) \htext(41
46){\tiny $0$}

\move(40 40)\lvec(50 50)\lvec(50 40)\lvec(40 40) \htext(46
42){\tiny $1$}
\end{texdraw}}
\hskip 3cm \raisebox{-0.4\height}{
\begin{texdraw}
\drawdim em \setunitscale 0.13 \linewd 0.5

\move(30 0)\lvec(40 0)\lvec(40 10)\lvec(30 10)\lvec(30 0)\htext(32
6){\tiny $1$}

\move(40 0)\lvec(50 0)\lvec(50 10)\lvec(40 10)\lvec(40 0)\htext(42
6){\tiny $0$}

\move(30 0)\lvec(40 10)\lvec(40 0)\lvec(30 0)\lfill f:0.8
\htext(36 2){\tiny $0$}

\move(40 0)\lvec(50 10)\lvec(50 0)\lvec(40 0)\lfill f:0.8
\htext(46 2){\tiny $1$}
\move(40 10)\lvec(50 10)\lvec(50 20)\lvec(40 20)\lvec(40
10)\htext(43.5 13){$_2$}
\move(40 20)\lvec(50 20)\lvec(50 30)\lvec(40 30)\lvec(40
20)\htext(43.5 26){\tiny $3$}

\move(40 20)\lvec(50 20)\lvec(50 25)\lvec(40 25)\lvec(40 20)
\htext(43.5 21){\tiny $3$}
\move(40 30)\lvec(50 30)\lvec(50 40)\lvec(40 40)\lvec(40
30)\htext(43.5 33){$_2$}

\move(40 50)\lvec(50 50)\lvec(50 60)\lvec(40 60)\lvec(40
50)\htext(43.5 53){$_2$}
\move(40 40)\lvec(50 50)\lvec(40 50)\lvec(40 40) \htext(41
46){\tiny $0$}

\move(40 40)\lvec(50 50)\lvec(50 40)\lvec(40 40) \htext(46
42){\tiny $1$}
\end{texdraw}}
\end{center}
}
\end{ex}
\vskip 5mm

Let $\Delta$ be the volume of the $\delta$-column. We list the
value of $\Delta$ for each quantum affine algebra $U_q(\frak{g})$
in the following table.\vskip 3mm

\begin{center}
\begin{tabular}{c|c}
$U_q(\frak{g})$ & $\Delta$ \\ \hline
$A^{(1)}_n$ & $n$ \\
$A^{(2)}_{2n-1}$ & $2n-2$ \\
$D^{(1)}_n$ & $2n-4$ \\
$A^{(2)}_{2n}$ & $2n$ \\
$D^{(2)}_{n+1}$ & $2n$ \\
$B^{(1)}_n$ & $2n-2$ \\ 
\end{tabular}
\end{center}\vskip 3mm

Note that $\Delta$ is not necessarily equal to the Coxeter number
or to the dual Coxeter number for $\frak{g}$ (cf. \cite{Kac90}).

Let ${\mathcal Y}(\Lambda)\subset{\mathcal Z}(\Lambda) $ be the
set of all reduced proper Young walls on $Y_{\Lambda}$. Then we
have:

\begin{thm}{\rm (\cite{Ka2000})}
For all $i\in I$ and $Y\in {\mathcal Y}(\La)$, we have
$$\tilde{E}_i Y \in {\mathcal Y} (\La) \cup \{0\},
\qquad \tilde{F}_i Y \in {\mathcal Y}(\La) \cup \{0\}.$$ Hence
${\mathcal Y}(\Lambda)$ is a $U_q(\frak{g})$-crystal. Moreover,
there exists an isomorphism of $U_q(\frak{g})$-crystals
$${\mathcal Y}(\La) \stackrel{\sim}{\longrightarrow} B(\Lambda)
\qquad \text{given by} \quad Y_{\La} \longmapsto u_{\La},$$ where
$B(\La)$ is the crystal of the basic representation $V(\La)$ of
$U_q(\g)$ and $u_{\La}$ is the highest weight vector in $B(\La)$.
\qed
\end{thm}
\newpage

\vskip 3mm

\begin{ex}
{\rm The crystal ${\mathcal Y}(\La_0)$ for $U_q(B_3^{(1)})$ is
given in Figure 1.}
\end{ex}

\begin{rem}{\rm
When $\zeta$ is a primitive $n$-th root of unity, the finite
dimensional irreducible representations of the (finite) Hecke
algebra ${\mathcal H}_{N}(\zeta)$ can be parametrized by {\it
$n$-reduced} colored Young diagrams. Observe that they are the
same as the reduced proper Young walls of type $A_{n-1}^{(1)}$. We
expect that for each type of classical quantum affine algebras and
level 1 dominant integral weights, there exist some interesting
algebraic structures whose irreducible representations (at some
specialization) are parametrized by reduced proper Young walls. In
\cite{BK}, Brundan and Kleshchev verified this idea by showing
that the irreducible representations of the Heck-Clifford
superalgebra ${\mathcal H}_N(\zeta)$ with $\zeta$ a primitive
$(2n+1)$-th root of unity are parametrized by the set of reduced
proper Young walls of type $A_{2n}^{(2)}$ with $N$ blocks.
}\end{rem}

\newpage

\begin{center}
\begin{texdraw}
\textref h:C v:C \fontsize{8}{8}\selectfont \drawdim mm
\setunitscale 5 \arrowheadtype t:F \arrowheadsize l:.4 w:.2
\nc{\dtri}{ \bsegment \move(-1 0)\lvec(0 1)\lvec(0 0)\lvec(-1
0)\ifill f:0.7 \esegment }
\htext(-1.4 24.5){$Y_{\Lambda_0}$} \move(-1 20.5) \bsegment
\move(0 0)\dtri \move(0 1)\lvec(-1 0)\lvec(-1 1)\lvec(0 1)\lvec(0
0)\lvec(-1 0) \htext(-0.3 0.3){$1$} \htext(-0.7 0.75){$0$}
\esegment \move(-1 16) \bsegment \move(0 0)\dtri \move(-1
1)\lvec(0 1)\lvec(-1 0)\lvec(-1 2)\lvec(0 2)\lvec(0 0)\lvec(-1 0)
\htext(-0.3 0.3){$1$} \htext(-0.7 0.75){$0$} \htext(-0.5 1.5){$2$}
\esegment \move(-3 11.5) \bsegment \move(0 0)\dtri \move(-1
0)\dtri \move(-1 1)\lvec(-2 0)\lvec(-2 1)\lvec(0 1)\lvec(-1
0)\lvec(-1 2)\lvec(0 2) \lvec(0 0)\lvec(-2 0) \htext(-0.3
0.3){$1$} \htext(-0.7 0.75){$0$} \htext(-1.3 0.3){$0$} \htext(-1.7
0.75){$1$} \htext(-0.5 1.5){$2$} \esegment \move(1 11.5) \bsegment
\move(0 0)\dtri \move(-1 1)\lvec(0 1)\lvec(-1 0)\lvec(-1
2.5)\lvec(0 2.5)\lvec(0 0)\lvec(-1 0) \move(-1 2)\lvec(0 2)
\htext(-0.3 0.3){$1$} \htext(-0.7 0.75){$0$} \htext(-0.5 1.5){$2$}
\htext(-0.5 2.25){$3$} \esegment \move(-3 6) \bsegment \move(0
0)\dtri \move(-1 0)\dtri \move(-1 1)\lvec(-2 0)\lvec(-2 1)\lvec(0
1)\lvec(-1 0)\lvec(-1 2.5)\lvec(0 2.5) \lvec(0 0)\lvec(-2 0)
\move(0 2)\lvec(-1 2) \htext(-0.3 0.3){$1$} \htext(-0.7 0.75){$0$}
\htext(-1.3 0.3){$0$} \htext(-1.7 0.75){$1$} \htext(-0.5 1.5){$2$}
\htext(-0.5 2.25){$3$} \esegment \move(1 6) \bsegment \move(0
0)\dtri \move(-1 1)\lvec(0 1)\lvec(-1 0)\lvec(-1 3)\lvec(0
3)\lvec(0 0)\lvec(-1 0) \move(-1 2)\lvec(0 2) \move(-1 2.5)\lvec(0
2.5) \htext(-0.3 0.3){$1$} \htext(-0.7 0.75){$0$} \htext(-0.5
1.5){$2$} \htext(-0.5 2.25){$3$} \htext(-0.5 2.75){$3$} \esegment
\move(-5 0) \bsegment \move(0 0)\dtri \move(-1 0)\dtri \move(-2
1)\lvec(0 1)\lvec(-1 0)\lvec(-1 2.5)\lvec(0 2.5)\lvec(0 0)
\lvec(-2 0)\lvec(-2 2)\lvec(0 2) \move(-2 0)\lvec(-1 1)
\htext(-0.3 0.3){$1$} \htext(-0.7 0.75){$0$} \htext(-1.3 0.3){$0$}
\htext(-1.7 0.75){$1$} \htext(-0.5 1.5){$2$} \htext(-1.5 1.5){$2$}
\htext(-0.5 2.25){$3$} \esegment \move(-1 0) \bsegment \move(0
0)\dtri \move(-1 0)\dtri \move(-1 1)\lvec(-2 0)\lvec(-2 1)\lvec(0
1)\lvec(-1 0)\lvec(-1 3)\lvec(0 3) \lvec(0 0)\lvec(-2 0) \move(0
2)\lvec(-1 2) \move(0 2.5)\lvec(-1 2.5) \htext(-0.3 0.3){$1$}
\htext(-0.7 0.75){$0$} \htext(-1.3 0.3){$0$} \htext(-1.7
0.75){$1$} \htext(-0.5 1.5){$2$} \htext(-0.5 2.25){$3$}
\htext(-0.5 2.75){$3$} \esegment \move(3 0) \bsegment \move(0
0)\dtri \move(-1 1)\lvec(0 1)\lvec(-1 0)\lvec(-1 4)\lvec(0
4)\lvec(0 0)\lvec(-1 0) \move(-1 2)\lvec(0 2) \move(-1 2.5)\lvec(0
2.5) \move(-1 3)\lvec(0 3) \htext(-0.3 0.3){$1$} \htext(-0.7
0.75){$0$} \htext(-0.5 1.5){$2$} \htext(-0.5 2.25){$3$}
\htext(-0.5 2.75){$3$} \htext(-0.5 3.5){$2$} \esegment \move(-6
-7) \bsegment \move(0 0)\dtri \move(-1 0)\dtri \move(-2 0)\dtri
\move(-2 1)\lvec(-3 0)\lvec(-3 1)\lvec(0 1)\lvec(-1 0)\lvec(-1
2.5) \lvec(0 2.5)\lvec(0 0)\lvec(-3 0) \move(-1 1)\lvec(-2
0)\lvec(-2 2)\lvec(0 2) \htext(-0.3 0.3){$1$} \htext(-0.7
0.75){$0$} \htext(-1.3 0.3){$0$} \htext(-1.7 0.75){$1$}
\htext(-2.3 0.3){$1$} \htext(-2.7 0.75){$0$} \htext(-0.5 1.5){$2$}
\htext(-1.5 1.5){$2$} \htext(-0.5 2.25){$3$} \esegment \move(-3
-7) \bsegment \move(0 0)\dtri \move(-1 0)\dtri \move(0 0)\lvec(-2
0)\lvec(-2 2.5)\lvec(0 2.5)\lvec(0 0) \move(-2 0)\lvec(-1 1)
\move(-1 0)\lvec(0 1)\lvec(-2 1) \move(-1 0)\lvec(-1 2.5) \move(0
2)\lvec(-2 2) \htext(-0.3 0.3){$1$} \htext(-0.7 0.75){$0$}
\htext(-1.3 0.3){$0$} \htext(-1.7 0.75){$1$} \htext(-0.5 1.5){$2$}
\htext(-1.5 1.5){$2$} \htext(-0.5 2.25){$3$} \htext(-1.5
2.25){$3$} \esegment \move(0 -7) \bsegment \move(0 0)\dtri
\move(-1 0)\dtri \move(-1 0)\lvec(-1 3)\lvec(0 3)\lvec(0
0)\lvec(-2 0)\lvec(-2 2)\lvec(0 2) \move(-2 1)\lvec(0 1)\lvec(-1
0) \move(-2 0)\lvec(-1 1) \move(0 2.5)\lvec(-1 2.5) \htext(-0.3
0.3){$1$} \htext(-0.7 0.75){$0$} \htext(-1.3 0.3){$0$} \htext(-1.7
0.75){$1$} \htext(-0.5 1.5){$2$} \htext(-1.5 1.5){$2$} \htext(-0.5
2.25){$3$} \htext(-0.5 2.75){$3$} \esegment \move(2 -7) \bsegment
\move(0 0)\dtri \move(-1 1)\lvec(0 1)\lvec(-1 0)\lvec(-1 5)\lvec(0
5)\lvec(0 0)\lvec(-1 0) \move(0 2)\lvec(-1 2) \move(0 2.5)\lvec(-1
2.5) \move(0 3)\lvec(-1 3) \move(0 4)\lvec(-1 4)\lvec(0 5)
\htext(-0.3 0.3){$1$} \htext(-0.7 0.75){$0$} \htext(-0.5 1.5){$2$}
\htext(-0.5 2.25){$3$} \htext(-0.5 2.75){$3$} \htext(-0.5
3.5){$2$} \htext(-0.7 4.75){$0$} \esegment \move(5 -7) \bsegment
\move(0 0)\dtri \move(-1 0)\dtri \move(-1 1)\lvec(-2 0)\lvec(-2
1)\lvec(0 1)\lvec(-1 0)\lvec(-1 4)\lvec(0 4) \lvec(0 0)\lvec(-2 0)
\move(0 2)\lvec(-1 2) \move(0 2.5)\lvec(-1 2.5) \move(0 3)\lvec(-1
3) \htext(-0.3 0.3){$1$} \htext(-0.7 0.75){$0$} \htext(-1.3
0.3){$0$} \htext(-1.7 0.75){$1$} \htext(-0.5 1.5){$2$} \htext(-0.5
2.25){$3$} \htext(-0.5 2.75){$3$} \htext(-0.5 3.5){$2$} \esegment
\move(-1.5 23.8)\avec(-1.5 21.7)\htext(-1.1 22.9){$0$} \move(-1.5
20.2)\avec(-1.5 18.2)\htext(-1.1 19.4){$2$} \move(-1.75
15.7)\avec(-3.4 13.7)\htext(-2.9 14.9){$1$} \move(-1.25
15.7)\avec(0.3 14.2)\htext(-0.1 15.2){$3$} \move(-3.5
11.2)\avec(-3.5 8.7)\htext(-3.9 10.1){$3$} \move(0.1
11.2)\avec(-3.1 8.7)\htext(-1.9 10.2){$1$} \move(0.5
11.2)\avec(0.5 9.2)\htext(0.9 10.3){$3$} \move(-4.2 5.7)\avec(-5.8
2.7)\htext(-5.4 4.4){$2$} \move(-3.8 5.7)\avec(-1.7
3.2)\htext(-2.4 4.7){$3$} \move(0.3 5.7)\avec(-1.3 3.2)\htext(-0.8
4.7){$1$} \move(0.7 5.7)\avec(2.3 4.2)\htext(1.8 5.2){$2$}
\move(-6.2 -0.3)\avec(-6.9 -4.3)\htext(-6.9 -2.1){$0$} \move(-5.8
-0.3)\avec(-4.1 -4.3)\htext(-4.6 -2.1){$3$} \move(-1.8
-0.3)\avec(-0.9 -3.8)\htext(-1.0 -1.9){$2$} \move(2.3
-0.3)\avec(1.6 -1.8)\htext(1.7 -0.8){$0$} \move(2.7 -0.3)\avec(4.4
-2.8)\htext(3.9 -1.3){$1$} \vtext(-7.5 -8){$\cdots$} \vtext(-4
-8){$\cdots$} \vtext(-1 -8){$\cdots$} \vtext(1.5 -8){$\cdots$}
\vtext(4 -8){$\cdots$} \move(-10 -8.5)\move(7 25)
\end{texdraw}

Figure 1
\end{center}

\vskip 1cm

\section{Fock Space Representation}

Let ${\mathcal F}(\Lambda)=\bigoplus_{Y\in {\mathcal
Z}(\Lambda)}\mathbb{Q}(q)Y$ be the $\mathbb{Q}(q)$-vector space
with a basis ${\mathcal Z}(\Lambda)$ consisting of proper Young
walls. The goal of this section is to define a
$U_q(\frak{g})$-module structure on ${\mathcal F}(\Lambda)$, the
{\it Fock space representation of $U_q(\frak{g})$}.

For this purpose, we introduce some terminologies. Let
$Y=(y_k)_{k=0}^{\infty}$ be a proper Young wall on $Y_{\Lambda}$,
and let $|y_k|$ denote the number of blocks in $y_k$ added to
$Y_{\Lambda}$. We define the {\it associated partition of $Y$} to
be $|Y|=(\cdots,|y_k|,\cdots,|y_1|,|y_0|)$. For proper Young walls
$Y=(y_k)_{k=0}^{\infty}$ and $Z=(z_k)_{k=0}^{\infty}$ in
${\mathcal Z}(\Lambda)$, we define $|Y|\unrhd|Z|$ if and only if
$\sum_{k=l}^{\infty}|y_k|\geq\sum_{k=l}^{\infty}|z_k|$ for all
$l\geq 0$. For example, if \vskip 5mm

\begin{center}
$Y=$ \raisebox{-0.4\height}{
\begin{texdraw}
\drawdim em \setunitscale 0.13 \linewd 0.5

\move(20 0)\lvec(30 0)\lvec(30 10)\lvec(20 10)\lvec(20 0)\htext(22
6){\tiny $0$}

\move(30 0)\lvec(40 0)\lvec(40 10)\lvec(30 10)\lvec(30 0)\htext(32
6){\tiny $1$}

\move(40 0)\lvec(50 0)\lvec(50 10)\lvec(40 10)\lvec(40 0)\htext(42
6){\tiny $0$}

\move(20 0)\lvec(30 10)\lvec(30 0)\lvec(20 0)\lfill f:0.8
\htext(26 2){\tiny $1$}

\move(30 0)\lvec(40 10)\lvec(40 0)\lvec(30 0)\lfill f:0.8
\htext(36 2){\tiny $0$}

\move(40 0)\lvec(50 10)\lvec(50 0)\lvec(40 0)\lfill f:0.8
\htext(46 2){\tiny $1$}
\move(20 10)\lvec(30 10)\lvec(30 20)\lvec(20 20)\lvec(20
10)\htext(23 13){$_2$}

\move(30 10)\lvec(40 10)\lvec(40 20)\lvec(30 20)\lvec(30
10)\htext(33 13){$_2$}

\move(40 10)\lvec(50 10)\lvec(50 20)\lvec(40 20)\lvec(40
10)\htext(43 13){$_2$}

\move(30 20)\lvec(40 20)\lvec(40 30)\lvec(30 30)\lvec(30
20)\htext(33 26){\tiny $3$}

\move(40 20)\lvec(50 20)\lvec(50 30)\lvec(40 30)\lvec(40
20)\htext(43 26){\tiny $3$}

\move(20 20)\lvec(30 20)\lvec(30 25)\lvec(20 25)\lvec(20 20)
\htext(23 21){\tiny $3$}

\move(30 20)\lvec(40 20)\lvec(40 25)\lvec(30 25)\lvec(30 20)
\htext(33 21){\tiny $3$}

\move(40 20)\lvec(50 20)\lvec(50 25)\lvec(40 25)\lvec(40 20)
\htext(43 21){\tiny $3$}

\move(30 30)\lvec(40 30)\lvec(40 40)\lvec(30 40)\lvec(30
30)\htext(33 33){$_2$}

\move(40 30)\lvec(50 30)\lvec(50 40)\lvec(40 40)\lvec(40
30)\htext(43 33){$_2$}

\move(40 40)\lvec(40 50)\lvec(50 50)\lvec(40 40) \htext(42
46){\tiny $0$}

\move(40 40)\lvec(50 50)\lvec(50 40)\lvec(40 40) \htext(46
42){\tiny $1$}

\end{texdraw}}\quad\quad and \quad
$Z=$ \raisebox{-0.4\height}{
\begin{texdraw}
\drawdim em \setunitscale 0.13 \linewd 0.5

\move(20 0)\lvec(30 0)\lvec(30 10)\lvec(20 10)\lvec(20 0)\htext(22
6){\tiny $0$}

\move(30 0)\lvec(40 0)\lvec(40 10)\lvec(30 10)\lvec(30 0)\htext(32
6){\tiny $1$}

\move(40 0)\lvec(50 0)\lvec(50 10)\lvec(40 10)\lvec(40 0)\htext(42
6){\tiny $0$}

\move(20 0)\lvec(30 10)\lvec(30 0)\lvec(20 0)\lfill f:0.8
\htext(26 2){\tiny $1$}

\move(30 0)\lvec(40 10)\lvec(40 0)\lvec(30 0)\lfill f:0.8
\htext(36 2){\tiny $0$}

\move(40 0)\lvec(50 10)\lvec(50 0)\lvec(40 0)\lfill f:0.8
\htext(46 2){\tiny $1$}
\move(40 50)\lvec(50 50)\lvec(50 60)\lvec(40 60)\lvec(40
50)\htext(43 53){$_2$}

\move(30 10)\lvec(40 10)\lvec(40 20)\lvec(30 20)\lvec(30
10)\htext(33 13){$_2$}

\move(40 10)\lvec(50 10)\lvec(50 20)\lvec(40 20)\lvec(40
10)\htext(43 13){$_2$}

\move(30 20)\lvec(40 20)\lvec(40 30)\lvec(30 30)\lvec(30
20)\htext(33 26){\tiny $3$}

\move(40 20)\lvec(50 20)\lvec(50 30)\lvec(40 30)\lvec(40
20)\htext(43 26){\tiny $3$}

\move(40 60)\lvec(50 60)\lvec(50 65)\lvec(40 65)\lvec(40 60)
\htext(43 61){\tiny $3$}

\move(30 20)\lvec(40 20)\lvec(40 25)\lvec(30 25)\lvec(30 20)
\htext(33 21){\tiny $3$}

\move(40 20)\lvec(50 20)\lvec(50 25)\lvec(40 25)\lvec(40 20)
\htext(43 21){\tiny $3$}

\move(30 30)\lvec(40 30)\lvec(40 40)\lvec(30 40)\lvec(30
30)\htext(33 33){$_2$}

\move(40 30)\lvec(50 30)\lvec(50 40)\lvec(40 40)\lvec(40
30)\htext(43 33){$_2$}

\move(40 40)\lvec(40 50)\lvec(50 50)\lvec(40 40) \htext(42
46){\tiny $0$}

\move(40 40)\lvec(50 50)\lvec(50 40)\lvec(40 40) \htext(46
42){\tiny $1$}

\end{texdraw}}\quad ,
\end{center}
then we have $|Y|=(3,5,7)\rhd |Z|=(1,5,9)$.

Note that it is not a partial ordering on ${\mathcal Z}(\Lambda)$
since there exist $Y\neq Z$ in ${\mathcal Z}(\Lambda)$ such that
$|Y|=|Z|$. However, it induces a partial ordering on the set of
associated partitions. (The readers may want to compare it with
the usual dominance ordering. See, for example, \cite{Mac})

Let $S=\{\,k\,|\,s\leq k < t\,\}$ for some $0\leq s <
t\leq\infty$. Then $S$ is a finite or an infinite interval in
$\mathbb{Z}_{\geq 0}$. We define  the {\it {\rm (}$S$-{\rm )}part
of $Y$} to be $Y_S=(y_k)_{k\in S}$. For example, if
$S=\{\,k\in\mathbb{Z}_{\geq 0 }\,|\,k\geq s\,\}$, the
$Y_S=(y_k)_{k=s}^{\infty}$ is itself a proper Young wall in
${\mathcal Z}(\Lambda')$ for some level 1 dominant integral weight
$\Lambda'$. On the other hand, if $S=\{\,k\,|\,s\leq k < t<\infty
\,\}$ is a finite interval, $Y_S$ is no more a proper Young wall,
but a finite collection of columns in $Y$. By restricting our
attentions to the columns in $Y_S$, we may define the notions of
admissible $i$-slots, removable $i$-blocks, the $i$-signature of
$Y_S$, $\varepsilon_i(Y_S)$, and $\varphi_i(Y_S)$ (however, {\rm
wt} can be defined for proper Young walls only).
\begin{ex}{\rm When $\frak{g}=A^{(2)}_4$ and $\Lambda=\Lambda_0$, consider
\vskip 5mm

\begin{center}
$Y=(y_k)_{k=0}^{\infty}=$\hskip
3mm\raisebox{-0.5\height}{\begin{texdraw}\drawdim em \setunitscale
0.13 \linewd 0.5

\move(0 0)\lvec(10 0)\lvec(10 10)\lvec(0 10)\lvec(0 0)\htext(3
1){\tiny $0$}

\move(10 0)\lvec(20 0)\lvec(20 10)\lvec(10 10)\lvec(10 0)\htext(13
1){\tiny $0$}

\move(20 0)\lvec(30 0)\lvec(30 10)\lvec(20 10)\lvec(20 0)\htext(23
1){\tiny $0$}

\move(30 0)\lvec(40 0)\lvec(40 10)\lvec(30 10)\lvec(30 0)\htext(33
1){\tiny $0$}

\move(40 0)\lvec(50 0)\lvec(50 10)\lvec(40 10)\lvec(40 0)\htext(43
1){\tiny $0$}

\move(50 0)\lvec(60 0)\lvec(60 10)\lvec(50 10)\lvec(50 0)\htext(53
1){\tiny $0$}
%

\move(0 0)\lvec(10 0)\lvec(10 5)\lvec(0 5)\lvec(0 0)\lfill f:0.8
\htext(3 6){\tiny $0$}

\move(10 0)\lvec(20 0)\lvec(20 5)\lvec(10 5)\lvec(10 0)\lfill
f:0.8 \htext(13 6){\tiny $0$}

\move(20 0)\lvec(30 0)\lvec(30 5)\lvec(20 5)\lvec(20 0)\lfill
f:0.8 \htext(23 6){\tiny $0$}

\move(30 0)\lvec(40 0)\lvec(40 5)\lvec(30 5)\lvec(30 0)\lfill
f:0.8 \htext(33 6){\tiny $0$}

\move(40 0)\lvec(50 0)\lvec(50 5)\lvec(40 5)\lvec(40 0)\lfill
f:0.8 \htext(43 6){\tiny $0$}

\move(50 0)\lvec(60 0)\lvec(60 5)\lvec(50 5)\lvec(50 0)\lfill
f:0.8 \htext(53 6){\tiny $0$}

\move(0 10)\lvec(10 10)\lvec(10 20)\lvec(0 20)\lvec(0 10)\htext(3
13){$_1$}

\move(10 10)\lvec(20 10)\lvec(20 20)\lvec(10 20)\lvec(10
10)\htext(13 13){$_1$}

\move(20 10)\lvec(30 10)\lvec(30 20)\lvec(20 20)\lvec(20
10)\htext(23 13){$_1$}

\move(30 10)\lvec(40 10)\lvec(40 20)\lvec(30 20)\lvec(30
10)\htext(33 13){$_1$}

\move(40 10)\lvec(50 10)\lvec(50 20)\lvec(40 20)\lvec(40
10)\htext(43 13){$_1$}

\move(50 10)\lvec(60 10)\lvec(60 20)\lvec(50 20)\lvec(50
10)\htext(53 13){$_1$}
\move(0 20)\lvec(10 20)\lvec(10 30)\lvec(0 30)\lvec(0 20)\htext(3
23){$_2$}

\move(10 20)\lvec(20 20)\lvec(20 30)\lvec(10 30)\lvec(10
20)\htext(13 23){$_2$}

\move(20 20)\lvec(30 20)\lvec(30 30)\lvec(20 30)\lvec(20
20)\htext(23 23){$_2$}

\move(30 20)\lvec(40 20)\lvec(40 30)\lvec(30 30)\lvec(30
20)\htext(33 23){$_2$}

\move(40 20)\lvec(50 20)\lvec(50 30)\lvec(40 30)\lvec(40
20)\htext(43 23){$_2$}

\move(50 20)\lvec(60 20)\lvec(60 30)\lvec(50 30)\lvec(50
20)\htext(53 23){$_2$}

\move(0 30)\lvec(10 30)\lvec(10 40)\lvec(0 40)\lvec(0 30)\htext(3
33){$_1$}

\move(10 30)\lvec(20 30)\lvec(20 40)\lvec(10 40)\lvec(10
30)\htext(13 33){$_1$}

\move(20 30)\lvec(30 30)\lvec(30 40)\lvec(20 40)\lvec(20
30)\htext(23 33){$_1$}

\move(30 30)\lvec(40 30)\lvec(40 40)\lvec(30 40)\lvec(30
30)\htext(33 33){$_1$}

\move(40 30)\lvec(50 30)\lvec(50 40)\lvec(40 40)\lvec(40
30)\htext(43 33){$_1$}

\move(50 30)\lvec(60 30)\lvec(60 40)\lvec(50 40)\lvec(50
30)\htext(53 33){$_1$}
\move(10 40)\lvec(20 40)\lvec(20 45)\lvec(10 45)\lvec(10
40)\htext(13 41){\tiny $0$}

\move(20 40)\lvec(30 40)\lvec(30 45)\lvec(20 45)\lvec(20
40)\htext(23 41){\tiny $0$}

\move(30 40)\lvec(40 40)\lvec(40 45)\lvec(30 45)\lvec(30
40)\htext(33 41){\tiny $0$}

\move(40 40)\lvec(50 40)\lvec(50 45)\lvec(40 45)\lvec(40
40)\htext(43 41){\tiny $0$}

\move(50 40)\lvec(60 40)\lvec(60 45)\lvec(50 45)\lvec(50
40)\htext(53 41){\tiny $0$}

\move(50 45)\lvec(60 45)\lvec(60 50)\lvec(50 50)\lvec(50
45)\htext(53 46){\tiny $0$}

\end{texdraw}}\quad .
\end{center}\vskip 5mm

Note that the $0$-signature of $Y$ is $+$ and  hence
$\varepsilon_0(Y)=0$, $\varphi_0(Y)=1$. If $S=\{\,0,1,2,3,4\,\}$,
then we have \vskip 5mm

\begin{center}
$Y_S=$\hskip 3mm\raisebox{-0.5\height}{\begin{texdraw}\drawdim em
\setunitscale 0.13 \linewd 0.5

\move(10 0)\lvec(20 0)\lvec(20 10)\lvec(10 10)\lvec(10 0)\htext(13
1){\tiny $0$}

\move(20 0)\lvec(30 0)\lvec(30 10)\lvec(20 10)\lvec(20 0)\htext(23
1){\tiny $0$}

\move(30 0)\lvec(40 0)\lvec(40 10)\lvec(30 10)\lvec(30 0)\htext(33
1){\tiny $0$}

\move(40 0)\lvec(50 0)\lvec(50 10)\lvec(40 10)\lvec(40 0)\htext(43
1){\tiny $0$}

\move(50 0)\lvec(60 0)\lvec(60 10)\lvec(50 10)\lvec(50 0)\htext(53
1){\tiny $0$}
%

\move(10 0)\lvec(20 0)\lvec(20 5)\lvec(10 5)\lvec(10 0)\lfill
f:0.8 \htext(13 6){\tiny $0$}

\move(20 0)\lvec(30 0)\lvec(30 5)\lvec(20 5)\lvec(20 0)\lfill
f:0.8 \htext(23 6){\tiny $0$}

\move(30 0)\lvec(40 0)\lvec(40 5)\lvec(30 5)\lvec(30 0)\lfill
f:0.8 \htext(33 6){\tiny $0$}

\move(40 0)\lvec(50 0)\lvec(50 5)\lvec(40 5)\lvec(40 0)\lfill
f:0.8 \htext(43 6){\tiny $0$}

\move(50 0)\lvec(60 0)\lvec(60 5)\lvec(50 5)\lvec(50 0)\lfill
f:0.8 \htext(53 6){\tiny $0$}

\move(10 10)\lvec(20 10)\lvec(20 20)\lvec(10 20)\lvec(10
10)\htext(13 13){$_1$}

\move(20 10)\lvec(30 10)\lvec(30 20)\lvec(20 20)\lvec(20
10)\htext(23 13){$_1$}

\move(30 10)\lvec(40 10)\lvec(40 20)\lvec(30 20)\lvec(30
10)\htext(33 13){$_1$}

\move(40 10)\lvec(50 10)\lvec(50 20)\lvec(40 20)\lvec(40
10)\htext(43 13){$_1$}

\move(50 10)\lvec(60 10)\lvec(60 20)\lvec(50 20)\lvec(50
10)\htext(53 13){$_1$}
\move(10 20)\lvec(20 20)\lvec(20 30)\lvec(10 30)\lvec(10
20)\htext(13 23){$_2$}

\move(20 20)\lvec(30 20)\lvec(30 30)\lvec(20 30)\lvec(20
20)\htext(23 23){$_2$}

\move(30 20)\lvec(40 20)\lvec(40 30)\lvec(30 30)\lvec(30
20)\htext(33 23){$_2$}

\move(40 20)\lvec(50 20)\lvec(50 30)\lvec(40 30)\lvec(40
20)\htext(43 23){$_2$}

\move(50 20)\lvec(60 20)\lvec(60 30)\lvec(50 30)\lvec(50
20)\htext(53 23){$_2$}

\move(10 30)\lvec(20 30)\lvec(20 40)\lvec(10 40)\lvec(10
30)\htext(13 33){$_1$}

\move(20 30)\lvec(30 30)\lvec(30 40)\lvec(20 40)\lvec(20
30)\htext(23 33){$_1$}

\move(30 30)\lvec(40 30)\lvec(40 40)\lvec(30 40)\lvec(30
30)\htext(33 33){$_1$}

\move(40 30)\lvec(50 30)\lvec(50 40)\lvec(40 40)\lvec(40
30)\htext(43 33){$_1$}

\move(50 30)\lvec(60 30)\lvec(60 40)\lvec(50 40)\lvec(50
30)\htext(53 33){$_1$}
\move(10 40)\lvec(20 40)\lvec(20 45)\lvec(10 45)\lvec(10
40)\htext(13 41){\tiny $0$}

\move(20 40)\lvec(30 40)\lvec(30 45)\lvec(20 45)\lvec(20
40)\htext(23 41){\tiny $0$}

\move(30 40)\lvec(40 40)\lvec(40 45)\lvec(30 45)\lvec(30
40)\htext(33 41){\tiny $0$}

\move(40 40)\lvec(50 40)\lvec(50 45)\lvec(40 45)\lvec(40
40)\htext(43 41){\tiny $0$}

\move(50 40)\lvec(60 40)\lvec(60 45)\lvec(50 45)\lvec(50
40)\htext(53 41){\tiny $0$}

\move(50 45)\lvec(60 45)\lvec(60 50)\lvec(50 50)\lvec(50
45)\htext(53 46){\tiny $0$}

\end{texdraw}}\quad ,
\end{center}\vskip 5mm

\noindent $\varepsilon_0(Y_S)=2$ and $\varphi_0(Y_S)=0$ since the
$0$-signature of $Y_S=--$. }
\end{ex}

We will now define the action of $U_q(\frak{g})$ on
$\mathcal{Z}(\Lambda)$. Since the action of $q^h$ ($h\in
P^{\vee}$) is easily defined by
\begin{equation}
q^h Y=q^{{\rm wt}(Y)(h)}Y \quad \text{for $Y\in
\mathcal{Z}(\Lambda)$,}
\end{equation}
we will focus on the actions of $e_i$ and $f_i$ ($i\in I$).

\vskip 2mm \noindent {\bf Case 1.} Suppose that the $i$-blocks are
of type I.\vskip 2mm

Let $b$ be a removable $i$-block in $y_k$ of $Y$. We define
$Y_R(b)=(y_{k-1},\cdots,y_{1},y_0)$ to be the part of $Y$
consisting of the columns lying at the right of $b$, and set
\begin{equation}
R_i(b;Y)=\varphi_i(Y_R(b))-\varepsilon_i(Y_R(b))
\end{equation}
If $k=0$, we understand $Y_R(b)=\emptyset$ and $R_i(b;Y)=0$. We
denote by $Y\nearrow b$ the Young wall obtained by removing $b$
from $Y$. Then we define
\begin{equation}
e_i\,Y=\sum_{b}q_i^{-R_i(b;Y)}(Y\nearrow b),
\end{equation}
\noindent where $b$ runs over all removable $i$-blocks in $Y$.

On the other hand, if $b$ is an admissible $i$-slot in $y_k$ of
$Y$, then we define $Y_L(b)=(\cdots,y_{k+2},y_{k+1})$ to be the
Young wall consisting of the columns in $Y$ lying at the left of
$b$, and set
\begin{equation}
L_i(b;Y)=\varphi_i(Y_L(b))-\varepsilon_i(Y_L(b)).
\end{equation}
Here, the wall $Y_L(b)$ is a proper Young wall on $Y_{\Lambda'}$
for some level 1 dominant integral weight $\Lambda'$. We denote by
$Y\swarrow b$ the Young wall obtained by adding an $i$-block at
$b$. Then we define
\begin{equation}
f_i\,Y=\sum_{b}q_i^{L_i(b;Y)}(Y\swarrow b),
\end{equation}
where $b$ runs over all admissible $i$-slots in $Y$.
\begin{ex}{\rm \mbox{}

(a) If $\frak{g}=A^{(1)}_2$, $\Lambda=\Lambda_0$ and $i=0$, then
$q_0=q$ and we have\vskip 5mm

$e_0$ \raisebox{-0.5\height}{
\begin{texdraw}\drawdim
em \setunitscale 0.13 \linewd 0.5

\move(10 0)\lvec(20 0)\lvec(20 10)\lvec(10 10)\lvec(10 0)\htext(13
3){$_0$}

\move(20 0)\lvec(30 0)\lvec(30 10)\lvec(20 10)\lvec(20 0)\htext(23
3){$_1$}

\move(30 0)\lvec(40 0)\lvec(40 10)\lvec(30 10)\lvec(30 0)\htext(33
3){$_2$}

\move(40 0)\lvec(50 0)\lvec(50 10)\lvec(40 10)\lvec(40 0)\htext(43
3){$_0$}

\move(20 10)\lvec(30 10)\lvec(30 20)\lvec(20 20)\lvec(20
10)\htext(23 13){$_2$}

\move(30 10)\lvec(40 10)\lvec(40 20)\lvec(30 20)\lvec(30
10)\htext(33 13){$_0$}

\move(40 10)\lvec(50 10)\lvec(50 20)\lvec(40 20)\lvec(40
10)\htext(43 13){$_1$}

\move(20 20)\lvec(30 20)\lvec(30 30)\lvec(20 30)\lvec(20
20)\htext(23 23){$_0$}

\move(30 20)\lvec(40 20)\lvec(40 30)\lvec(30 30)\lvec(30
20)\htext(33 23){$_1$}

\move(40 20)\lvec(50 20)\lvec(50 30)\lvec(40 30)\lvec(40
20)\htext(43 23){$_2$}

\move(30 30)\lvec(40 30)\lvec(40 40)\lvec(30 40)\lvec(30
30)\htext(33 33){$_2$}

\move(40 30)\lvec(50 30)\lvec(50 40)\lvec(40 40)\lvec(40
30)\htext(43 33){$_0$}

\move(30 40)\lvec(40 40)\lvec(40 50)\lvec(30 50)\lvec(30
40)\htext(33 43){$_0$}

\move(40 40)\lvec(50 40)\lvec(50 50)\lvec(40 50)\lvec(40
40)\htext(43 43){$_1$}

\move(40 50)\lvec(50 50)\lvec(50 60)\lvec(40 60)\lvec(40
50)\htext(43 53){$_2$}
\end{texdraw}}\vskip 5mm

$=$ \hskip 5mm $q$ \raisebox{-0.5\height}{
\begin{texdraw}\drawdim
em \setunitscale 0.13 \linewd 0.5

\move(20 0)\lvec(30 0)\lvec(30 10)\lvec(20 10)\lvec(20 0)\htext(23
3){$_1$}

\move(30 0)\lvec(40 0)\lvec(40 10)\lvec(30 10)\lvec(30 0)\htext(33
3){$_2$}

\move(40 0)\lvec(50 0)\lvec(50 10)\lvec(40 10)\lvec(40 0)\htext(43
3){$_0$}

\move(20 10)\lvec(30 10)\lvec(30 20)\lvec(20 20)\lvec(20
10)\htext(23 13){$_2$}

\move(30 10)\lvec(40 10)\lvec(40 20)\lvec(30 20)\lvec(30
10)\htext(33 13){$_0$}

\move(40 10)\lvec(50 10)\lvec(50 20)\lvec(40 20)\lvec(40
10)\htext(43 13){$_1$}

\move(20 20)\lvec(30 20)\lvec(30 30)\lvec(20 30)\lvec(20
20)\htext(23 23){$_0$}

\move(30 20)\lvec(40 20)\lvec(40 30)\lvec(30 30)\lvec(30
20)\htext(33 23){$_1$}

\move(40 20)\lvec(50 20)\lvec(50 30)\lvec(40 30)\lvec(40
20)\htext(43 23){$_2$}

\move(30 30)\lvec(40 30)\lvec(40 40)\lvec(30 40)\lvec(30
30)\htext(33 33){$_2$}

\move(40 30)\lvec(50 30)\lvec(50 40)\lvec(40 40)\lvec(40
30)\htext(43 33){$_0$}

\move(30 40)\lvec(40 40)\lvec(40 50)\lvec(30 50)\lvec(30
40)\htext(33 43){$_0$}

\move(40 40)\lvec(50 40)\lvec(50 50)\lvec(40 50)\lvec(40
40)\htext(43 43){$_1$}

\move(40 50)\lvec(50 50)\lvec(50 60)\lvec(40 60)\lvec(40
50)\htext(43 53){$_2$}
\end{texdraw}}\hskip 5mm $+$ \hskip 5mm \raisebox{-0.5\height}{
\begin{texdraw}\drawdim
em \setunitscale 0.13 \linewd 0.5

\move(10 0)\lvec(20 0)\lvec(20 10)\lvec(10 10)\lvec(10 0)\htext(13
3){$_0$}

\move(20 0)\lvec(30 0)\lvec(30 10)\lvec(20 10)\lvec(20 0)\htext(23
3){$_1$}

\move(30 0)\lvec(40 0)\lvec(40 10)\lvec(30 10)\lvec(30 0)\htext(33
3){$_2$}

\move(40 0)\lvec(50 0)\lvec(50 10)\lvec(40 10)\lvec(40 0)\htext(43
3){$_0$}

\move(20 10)\lvec(30 10)\lvec(30 20)\lvec(20 20)\lvec(20
10)\htext(23 13){$_2$}

\move(30 10)\lvec(40 10)\lvec(40 20)\lvec(30 20)\lvec(30
10)\htext(33 13){$_0$}

\move(40 10)\lvec(50 10)\lvec(50 20)\lvec(40 20)\lvec(40
10)\htext(43 13){$_1$}

\move(30 20)\lvec(40 20)\lvec(40 30)\lvec(30 30)\lvec(30
20)\htext(33 23){$_1$}

\move(40 20)\lvec(50 20)\lvec(50 30)\lvec(40 30)\lvec(40
20)\htext(43 23){$_2$}

\move(30 30)\lvec(40 30)\lvec(40 40)\lvec(30 40)\lvec(30
30)\htext(33 33){$_2$}

\move(40 30)\lvec(50 30)\lvec(50 40)\lvec(40 40)\lvec(40
30)\htext(43 33){$_0$}

\move(30 40)\lvec(40 40)\lvec(40 50)\lvec(30 50)\lvec(30
40)\htext(33 43){$_0$}

\move(40 40)\lvec(50 40)\lvec(50 50)\lvec(40 50)\lvec(40
40)\htext(43 43){$_1$}

\move(40 50)\lvec(50 50)\lvec(50 60)\lvec(40 60)\lvec(40
50)\htext(43 53){$_2$}
\end{texdraw}}
\hskip 5mm $+$ \hskip 5mm $q^{-1}$\raisebox{-0.5\height}{
\begin{texdraw}\drawdim
em \setunitscale 0.13 \linewd 0.5

\move(10 0)\lvec(20 0)\lvec(20 10)\lvec(10 10)\lvec(10 0)\htext(13
3){$_0$}

\move(20 0)\lvec(30 0)\lvec(30 10)\lvec(20 10)\lvec(20 0)\htext(23
3){$_1$}

\move(30 0)\lvec(40 0)\lvec(40 10)\lvec(30 10)\lvec(30 0)\htext(33
3){$_2$}

\move(40 0)\lvec(50 0)\lvec(50 10)\lvec(40 10)\lvec(40 0)\htext(43
3){$_0$}

\move(20 10)\lvec(30 10)\lvec(30 20)\lvec(20 20)\lvec(20
10)\htext(23 13){$_2$}

\move(30 10)\lvec(40 10)\lvec(40 20)\lvec(30 20)\lvec(30
10)\htext(33 13){$_0$}

\move(40 10)\lvec(50 10)\lvec(50 20)\lvec(40 20)\lvec(40
10)\htext(43 13){$_1$}

\move(20 20)\lvec(30 20)\lvec(30 30)\lvec(20 30)\lvec(20
20)\htext(23 23){$_0$}

\move(30 20)\lvec(40 20)\lvec(40 30)\lvec(30 30)\lvec(30
20)\htext(33 23){$_1$}

\move(40 20)\lvec(50 20)\lvec(50 30)\lvec(40 30)\lvec(40
20)\htext(43 23){$_2$}

\move(30 30)\lvec(40 30)\lvec(40 40)\lvec(30 40)\lvec(30
30)\htext(33 33){$_2$}

\move(40 30)\lvec(50 30)\lvec(50 40)\lvec(40 40)\lvec(40
30)\htext(43 33){$_0$}

\move(40 40)\lvec(50 40)\lvec(50 50)\lvec(40 50)\lvec(40
40)\htext(43 43){$_1$}

\move(40 50)\lvec(50 50)\lvec(50 60)\lvec(40 60)\lvec(40
50)\htext(43 53){$_2$}
\end{texdraw}}\quad ,
\vskip 10mm

$f_0$ \raisebox{-0.5\height}{
\begin{texdraw}\drawdim
em \setunitscale 0.13 \linewd 0.5

\move(10 0)\lvec(20 0)\lvec(20 10)\lvec(10 10)\lvec(10 0)\htext(13
3){$_0$}

\move(20 0)\lvec(30 0)\lvec(30 10)\lvec(20 10)\lvec(20 0)\htext(23
3){$_1$}

\move(30 0)\lvec(40 0)\lvec(40 10)\lvec(30 10)\lvec(30 0)\htext(33
3){$_2$}

\move(40 0)\lvec(50 0)\lvec(50 10)\lvec(40 10)\lvec(40 0)\htext(43
3){$_0$}

\move(20 10)\lvec(30 10)\lvec(30 20)\lvec(20 20)\lvec(20
10)\htext(23 13){$_2$}

\move(30 10)\lvec(40 10)\lvec(40 20)\lvec(30 20)\lvec(30
10)\htext(33 13){$_0$}

\move(40 10)\lvec(50 10)\lvec(50 20)\lvec(40 20)\lvec(40
10)\htext(43 13){$_1$}

\move(30 20)\lvec(40 20)\lvec(40 30)\lvec(30 30)\lvec(30
20)\htext(33 23){$_1$}

\move(40 20)\lvec(50 20)\lvec(50 30)\lvec(40 30)\lvec(40
20)\htext(43 23){$_2$}

\move(30 30)\lvec(40 30)\lvec(40 40)\lvec(30 40)\lvec(30
30)\htext(33 33){$_2$}

\move(40 30)\lvec(50 30)\lvec(50 40)\lvec(40 40)\lvec(40
30)\htext(43 33){$_0$}

\move(40 40)\lvec(50 40)\lvec(50 50)\lvec(40 50)\lvec(40
40)\htext(43 43){$_1$}

\move(40 50)\lvec(50 50)\lvec(50 60)\lvec(40 60)\lvec(40
50)\htext(43 53){$_2$}
\end{texdraw}}\vskip 5mm

$=$ \hskip 5mm $q^{-1}$ \raisebox{-0.5\height}{
\begin{texdraw}\drawdim
em \setunitscale 0.13 \linewd 0.5

\move(10 0)\lvec(20 0)\lvec(20 10)\lvec(10 10)\lvec(10 0)\htext(13
3){$_0$}

\move(20 0)\lvec(30 0)\lvec(30 10)\lvec(20 10)\lvec(20 0)\htext(23
3){$_1$}

\move(30 0)\lvec(40 0)\lvec(40 10)\lvec(30 10)\lvec(30 0)\htext(33
3){$_2$}

\move(40 0)\lvec(50 0)\lvec(50 10)\lvec(40 10)\lvec(40 0)\htext(43
3){$_0$}

\move(20 10)\lvec(30 10)\lvec(30 20)\lvec(20 20)\lvec(20
10)\htext(23 13){$_2$}

\move(30 10)\lvec(40 10)\lvec(40 20)\lvec(30 20)\lvec(30
10)\htext(33 13){$_0$}

\move(40 10)\lvec(50 10)\lvec(50 20)\lvec(40 20)\lvec(40
10)\htext(43 13){$_1$}

\move(20 20)\lvec(30 20)\lvec(30 30)\lvec(20 30)\lvec(20
20)\htext(23 23){$_0$}

\move(30 20)\lvec(40 20)\lvec(40 30)\lvec(30 30)\lvec(30
20)\htext(33 23){$_1$}

\move(40 20)\lvec(50 20)\lvec(50 30)\lvec(40 30)\lvec(40
20)\htext(43 23){$_2$}

\move(30 30)\lvec(40 30)\lvec(40 40)\lvec(30 40)\lvec(30
30)\htext(33 33){$_2$}

\move(40 30)\lvec(50 30)\lvec(50 40)\lvec(40 40)\lvec(40
30)\htext(43 33){$_0$}

\move(40 40)\lvec(50 40)\lvec(50 50)\lvec(40 50)\lvec(40
40)\htext(43 43){$_1$}

\move(40 50)\lvec(50 50)\lvec(50 60)\lvec(40 60)\lvec(40
50)\htext(43 53){$_2$}
\end{texdraw}}\hskip 5mm $+$ \hskip 5mm \raisebox{-0.5\height}{
\begin{texdraw}\drawdim
em \setunitscale 0.13 \linewd 0.5

\move(10 0)\lvec(20 0)\lvec(20 10)\lvec(10 10)\lvec(10 0)\htext(13
3){$_0$}

\move(20 0)\lvec(30 0)\lvec(30 10)\lvec(20 10)\lvec(20 0)\htext(23
3){$_1$}

\move(30 0)\lvec(40 0)\lvec(40 10)\lvec(30 10)\lvec(30 0)\htext(33
3){$_2$}

\move(40 0)\lvec(50 0)\lvec(50 10)\lvec(40 10)\lvec(40 0)\htext(43
3){$_0$}

\move(20 10)\lvec(30 10)\lvec(30 20)\lvec(20 20)\lvec(20
10)\htext(23 13){$_2$}

\move(30 10)\lvec(40 10)\lvec(40 20)\lvec(30 20)\lvec(30
10)\htext(33 13){$_0$}

\move(40 10)\lvec(50 10)\lvec(50 20)\lvec(40 20)\lvec(40
10)\htext(43 13){$_1$}

\move(30 20)\lvec(40 20)\lvec(40 30)\lvec(30 30)\lvec(30
20)\htext(33 23){$_1$}

\move(40 20)\lvec(50 20)\lvec(50 30)\lvec(40 30)\lvec(40
20)\htext(43 23){$_2$}

\move(30 30)\lvec(40 30)\lvec(40 40)\lvec(30 40)\lvec(30
30)\htext(33 33){$_2$}

\move(40 30)\lvec(50 30)\lvec(50 40)\lvec(40 40)\lvec(40
30)\htext(43 33){$_0$}

\move(30 40)\lvec(40 40)\lvec(40 50)\lvec(30 50)\lvec(30
40)\htext(33 43){$_0$}

\move(40 40)\lvec(50 40)\lvec(50 50)\lvec(40 50)\lvec(40
40)\htext(43 43){$_1$}

\move(40 50)\lvec(50 50)\lvec(50 60)\lvec(40 60)\lvec(40
50)\htext(43 53){$_2$}
\end{texdraw}}\hskip 5mm
$+$ \hskip 5mm $q$\raisebox{-0.5\height}{
\begin{texdraw}\drawdim
em \setunitscale 0.13 \linewd 0.5

\move(10 0)\lvec(20 0)\lvec(20 10)\lvec(10 10)\lvec(10 0)\htext(13
3){$_0$}

\move(20 0)\lvec(30 0)\lvec(30 10)\lvec(20 10)\lvec(20 0)\htext(23
3){$_1$}

\move(30 0)\lvec(40 0)\lvec(40 10)\lvec(30 10)\lvec(30 0)\htext(33
3){$_2$}

\move(40 0)\lvec(50 0)\lvec(50 10)\lvec(40 10)\lvec(40 0)\htext(43
3){$_0$}

\move(20 10)\lvec(30 10)\lvec(30 20)\lvec(20 20)\lvec(20
10)\htext(23 13){$_2$}

\move(30 10)\lvec(40 10)\lvec(40 20)\lvec(30 20)\lvec(30
10)\htext(33 13){$_0$}

\move(40 10)\lvec(50 10)\lvec(50 20)\lvec(40 20)\lvec(40
10)\htext(43 13){$_1$}

\move(30 20)\lvec(40 20)\lvec(40 30)\lvec(30 30)\lvec(30
20)\htext(33 23){$_1$}

\move(40 20)\lvec(50 20)\lvec(50 30)\lvec(40 30)\lvec(40
20)\htext(43 23){$_2$}

\move(30 30)\lvec(40 30)\lvec(40 40)\lvec(30 40)\lvec(30
30)\htext(33 33){$_2$}

\move(40 30)\lvec(50 30)\lvec(50 40)\lvec(40 40)\lvec(40
30)\htext(43 33){$_0$}

\move(40 40)\lvec(50 40)\lvec(50 50)\lvec(40 50)\lvec(40
40)\htext(43 43){$_1$}

\move(40 50)\lvec(50 50)\lvec(50 60)\lvec(40 60)\lvec(40
50)\htext(43 53){$_2$}

\move(40 60)\lvec(50 60)\lvec(50 70)\lvec(40 70)\lvec(40
60)\htext(43 63){$_0$}

\end{texdraw}}\quad .
\vskip 5mm

(b) If $\frak{g}=A^{(2)}_5$, $\Lambda=\Lambda_0$ and $i=2$, then
$q_2=q$ and we have \vskip 5mm

$e_2$ \raisebox{-0.5\height}{\begin{texdraw}\drawdim em
\setunitscale 0.13 \linewd 0.5

\move(10 0)\lvec(20 0)\lvec(20 10)\lvec(10 10)\lvec(10 0)

\move(20 0)\lvec(30 0)\lvec(30 10)\lvec(20 10)\lvec(20 0)

\move(30 0)\lvec(40 0)\lvec(40 10)\lvec(30 10)\lvec(30 0)

\move(40 0)\lvec(50 0)\lvec(50 10)\lvec(40 10)\lvec(40 0)

\move(10 10)\lvec(20 10)\lvec(20 20)\lvec(10 20)\lvec(10
10)\htext(13 13){$_2$}

\move(20 10)\lvec(30 10)\lvec(30 20)\lvec(20 20)\lvec(20
10)\htext(23 13){$_2$}

\move(30 10)\lvec(40 10)\lvec(40 20)\lvec(30 20)\lvec(30
10)\htext(33 13){$_2$}

\move(40 10)\lvec(50 10)\lvec(50 20)\lvec(40 20)\lvec(40
10)\htext(43 13){$_2$}

\move(20 20)\lvec(30 20)\lvec(30 30)\lvec(20 30)\lvec(20
20)\htext(23 23){$_3$}

\move(30 20)\lvec(40 20)\lvec(40 30)\lvec(30 30)\lvec(30
20)\htext(33 23){$_3$}

\move(40 20)\lvec(50 20)\lvec(50 30)\lvec(40 30)\lvec(40
20)\htext(43 23){$_3$}

\move(20 30)\lvec(30 30)\lvec(30 40)\lvec(20 40)\lvec(20
30)\htext(23 33){$_2$}

\move(30 30)\lvec(40 30)\lvec(40 40)\lvec(30 40)\lvec(30
30)\htext(33 33){$_2$}

\move(40 30)\lvec(50 30)\lvec(50 40)\lvec(40 40)\lvec(40
30)\htext(43 33){$_2$}

\move(30 40)\lvec(40 40)\lvec(40 50)\lvec(30 50)\lvec(30 40)

\move(40 40)\lvec(50 40)\lvec(50 50)\lvec(40 50)\lvec(40 40)

\move(30 50)\lvec(40 50)\lvec(40 60)\lvec(30 60)\lvec(30
50)\htext(33 53){$_2$}

\move(40 50)\lvec(50 50)\lvec(50 60)\lvec(40 60)\lvec(40
50)\htext(43 53){$_2$}

\move(40 60)\lvec(50 60)\lvec(50 70)\lvec(40 70)\lvec(40
60)\htext(43 63){$_3$}

\move(10 0)\lvec(20 10)\lvec(20 0)\lvec(10 0)\lfill f:0.8
\htext(11 5){\tiny $1$}\htext(16 1){\tiny $0$}

\move(20 0)\lvec(30 10)\lvec(30 0)\lvec(20 0)\lfill f:0.8
\htext(22 5){\tiny $0$}\htext(26 2){\tiny $1$}

\move(30 0)\lvec(40 10)\lvec(40 0)\lvec(30 0)\lfill f:0.8
\htext(32 5){\tiny $1$}\htext(36 2){\tiny $0$}

\move(40 0)\lvec(50 10)\lvec(50 0)\lvec(40 0)\lfill f:0.8
\htext(42 5){\tiny $0$}\htext(46 2){\tiny $1$}

\move(30 40)\lvec(40 50)\lvec(40 40)\lvec(30 40)\htext(32
45){\tiny $1$}\htext(36 42){\tiny $0$}

\move(40 40)\lvec(50 50)\lvec(50 40)\lvec(40 40)\htext(42
45){\tiny $0$}\htext(46 42){\tiny $1$}
\end{texdraw}}\vskip 5mm

$=$\hskip 5mm $q$
 \raisebox{-0.5\height}{\begin{texdraw}\drawdim em
\setunitscale 0.13 \linewd 0.5

\move(10 0)\lvec(20 0)\lvec(20 10)\lvec(10 10)\lvec(10 0)

\move(20 0)\lvec(30 0)\lvec(30 10)\lvec(20 10)\lvec(20 0)

\move(30 0)\lvec(40 0)\lvec(40 10)\lvec(30 10)\lvec(30 0)

\move(40 0)\lvec(50 0)\lvec(50 10)\lvec(40 10)\lvec(40 0)


\move(20 10)\lvec(30 10)\lvec(30 20)\lvec(20 20)\lvec(20
10)\htext(23 13){$_2$}

\move(30 10)\lvec(40 10)\lvec(40 20)\lvec(30 20)\lvec(30
10)\htext(33 13){$_2$}

\move(40 10)\lvec(50 10)\lvec(50 20)\lvec(40 20)\lvec(40
10)\htext(43 13){$_2$}

\move(20 20)\lvec(30 20)\lvec(30 30)\lvec(20 30)\lvec(20
20)\htext(23 23){$_3$}

\move(30 20)\lvec(40 20)\lvec(40 30)\lvec(30 30)\lvec(30
20)\htext(33 23){$_3$}

\move(40 20)\lvec(50 20)\lvec(50 30)\lvec(40 30)\lvec(40
20)\htext(43 23){$_3$}

\move(20 30)\lvec(30 30)\lvec(30 40)\lvec(20 40)\lvec(20
30)\htext(23 33){$_2$}

\move(30 30)\lvec(40 30)\lvec(40 40)\lvec(30 40)\lvec(30
30)\htext(33 33){$_2$}

\move(40 30)\lvec(50 30)\lvec(50 40)\lvec(40 40)\lvec(40
30)\htext(43 33){$_2$}

\move(30 40)\lvec(40 40)\lvec(40 50)\lvec(30 50)\lvec(30 40)

\move(40 40)\lvec(50 40)\lvec(50 50)\lvec(40 50)\lvec(40 40)

\move(30 50)\lvec(40 50)\lvec(40 60)\lvec(30 60)\lvec(30
50)\htext(33 53){$_2$}

\move(40 50)\lvec(50 50)\lvec(50 60)\lvec(40 60)\lvec(40
50)\htext(43 53){$_2$}

\move(40 60)\lvec(50 60)\lvec(50 70)\lvec(40 70)\lvec(40
60)\htext(43 63){$_3$}

\move(10 0)\lvec(20 10)\lvec(20 0)\lvec(10 0)\lfill f:0.8
\htext(11 5){\tiny $1$}\htext(16 1){\tiny $0$}

\move(20 0)\lvec(30 10)\lvec(30 0)\lvec(20 0)\lfill f:0.8
\htext(22 5){\tiny $0$}\htext(26 2){\tiny $1$}

\move(30 0)\lvec(40 10)\lvec(40 0)\lvec(30 0)\lfill f:0.8
\htext(32 5){\tiny $1$}\htext(36 2){\tiny $0$}

\move(40 0)\lvec(50 10)\lvec(50 0)\lvec(40 0)\lfill f:0.8
\htext(42 5){\tiny $0$}\htext(46 2){\tiny $1$}

\move(30 40)\lvec(40 50)\lvec(40 40)\lvec(30 40)\htext(32
45){\tiny $1$}\htext(36 42){\tiny $0$}

\move(40 40)\lvec(50 50)\lvec(50 40)\lvec(40 40)\htext(42
45){\tiny $0$}\htext(46 42){\tiny $1$}
\end{texdraw}}\hskip 5mm $+$\hskip 5mm
\raisebox{-0.5\height}{\begin{texdraw}\drawdim em \setunitscale
0.13 \linewd 0.5

\move(10 0)\lvec(20 0)\lvec(20 10)\lvec(10 10)\lvec(10 0)

\move(20 0)\lvec(30 0)\lvec(30 10)\lvec(20 10)\lvec(20 0)

\move(30 0)\lvec(40 0)\lvec(40 10)\lvec(30 10)\lvec(30 0)

\move(40 0)\lvec(50 0)\lvec(50 10)\lvec(40 10)\lvec(40 0)

\move(10 10)\lvec(20 10)\lvec(20 20)\lvec(10 20)\lvec(10
10)\htext(13 13){$_2$}

\move(20 10)\lvec(30 10)\lvec(30 20)\lvec(20 20)\lvec(20
10)\htext(23 13){$_2$}

\move(30 10)\lvec(40 10)\lvec(40 20)\lvec(30 20)\lvec(30
10)\htext(33 13){$_2$}

\move(40 10)\lvec(50 10)\lvec(50 20)\lvec(40 20)\lvec(40
10)\htext(43 13){$_2$}

\move(20 20)\lvec(30 20)\lvec(30 30)\lvec(20 30)\lvec(20
20)\htext(23 23){$_3$}

\move(30 20)\lvec(40 20)\lvec(40 30)\lvec(30 30)\lvec(30
20)\htext(33 23){$_3$}

\move(40 20)\lvec(50 20)\lvec(50 30)\lvec(40 30)\lvec(40
20)\htext(43 23){$_3$}

\move(30 30)\lvec(40 30)\lvec(40 40)\lvec(30 40)\lvec(30
30)\htext(33 33){$_2$}

\move(40 30)\lvec(50 30)\lvec(50 40)\lvec(40 40)\lvec(40
30)\htext(43 33){$_2$}

\move(30 40)\lvec(40 40)\lvec(40 50)\lvec(30 50)\lvec(30 40)

\move(40 40)\lvec(50 40)\lvec(50 50)\lvec(40 50)\lvec(40 40)

\move(30 50)\lvec(40 50)\lvec(40 60)\lvec(30 60)\lvec(30
50)\htext(33 53){$_2$}

\move(40 50)\lvec(50 50)\lvec(50 60)\lvec(40 60)\lvec(40
50)\htext(43 53){$_2$}

\move(40 60)\lvec(50 60)\lvec(50 70)\lvec(40 70)\lvec(40
60)\htext(43 63){$_3$}

\move(10 0)\lvec(20 10)\lvec(20 0)\lvec(10 0)\lfill f:0.8
\htext(11 5){\tiny $1$}\htext(16 1){\tiny $0$}

\move(20 0)\lvec(30 10)\lvec(30 0)\lvec(20 0)\lfill f:0.8
\htext(22 5){\tiny $0$}\htext(26 2){\tiny $1$}

\move(30 0)\lvec(40 10)\lvec(40 0)\lvec(30 0)\lfill f:0.8
\htext(32 5){\tiny $1$}\htext(36 2){\tiny $0$}

\move(40 0)\lvec(50 10)\lvec(50 0)\lvec(40 0)\lfill f:0.8
\htext(42 5){\tiny $0$}\htext(46 2){\tiny $1$}

\move(30 40)\lvec(40 50)\lvec(40 40)\lvec(30 40)\htext(32
45){\tiny $1$}\htext(36 42){\tiny $0$}

\move(40 40)\lvec(50 50)\lvec(50 40)\lvec(40 40)\htext(42
45){\tiny $0$}\htext(46 42){\tiny $1$}
\end{texdraw}}\hskip 5mm
$+$\hskip 5mm $q^{-1}$ \raisebox{-0.5\height}{
\begin{texdraw}\drawdim
em \setunitscale 0.13 \linewd 0.5

\move(10 0)\lvec(20 0)\lvec(20 10)\lvec(10 10)\lvec(10 0)

\move(20 0)\lvec(30 0)\lvec(30 10)\lvec(20 10)\lvec(20 0)

\move(30 0)\lvec(40 0)\lvec(40 10)\lvec(30 10)\lvec(30 0)

\move(40 0)\lvec(50 0)\lvec(50 10)\lvec(40 10)\lvec(40 0)

\move(10 10)\lvec(20 10)\lvec(20 20)\lvec(10 20)\lvec(10
10)\htext(13 13){$_2$}

\move(20 10)\lvec(30 10)\lvec(30 20)\lvec(20 20)\lvec(20
10)\htext(23 13){$_2$}

\move(30 10)\lvec(40 10)\lvec(40 20)\lvec(30 20)\lvec(30
10)\htext(33 13){$_2$}

\move(40 10)\lvec(50 10)\lvec(50 20)\lvec(40 20)\lvec(40
10)\htext(43 13){$_2$}

\move(20 20)\lvec(30 20)\lvec(30 30)\lvec(20 30)\lvec(20
20)\htext(23 23){$_3$}

\move(30 20)\lvec(40 20)\lvec(40 30)\lvec(30 30)\lvec(30
20)\htext(33 23){$_3$}

\move(40 20)\lvec(50 20)\lvec(50 30)\lvec(40 30)\lvec(40
20)\htext(43 23){$_3$}

\move(20 30)\lvec(30 30)\lvec(30 40)\lvec(20 40)\lvec(20
30)\htext(23 33){$_2$}

\move(30 30)\lvec(40 30)\lvec(40 40)\lvec(30 40)\lvec(30
30)\htext(33 33){$_2$}

\move(40 30)\lvec(50 30)\lvec(50 40)\lvec(40 40)\lvec(40
30)\htext(43 33){$_2$}

\move(30 40)\lvec(40 40)\lvec(40 50)\lvec(30 50)\lvec(30 40)

\move(40 40)\lvec(50 40)\lvec(50 50)\lvec(40 50)\lvec(40 40)


\move(40 50)\lvec(50 50)\lvec(50 60)\lvec(40 60)\lvec(40
50)\htext(43 53){$_2$}

\move(40 60)\lvec(50 60)\lvec(50 70)\lvec(40 70)\lvec(40
60)\htext(43 63){$_3$}

\move(10 0)\lvec(20 10)\lvec(20 0)\lvec(10 0)\lfill f:0.8
\htext(11 5){\tiny $1$}\htext(16 1){\tiny $0$}

\move(20 0)\lvec(30 10)\lvec(30 0)\lvec(20 0)\lfill f:0.8
\htext(22 5){\tiny $0$}\htext(26 2){\tiny $1$}

\move(30 0)\lvec(40 10)\lvec(40 0)\lvec(30 0)\lfill f:0.8
\htext(32 5){\tiny $1$}\htext(36 2){\tiny $0$}

\move(40 0)\lvec(50 10)\lvec(50 0)\lvec(40 0)\lfill f:0.8
\htext(42 5){\tiny $0$}\htext(46 2){\tiny $1$}

\move(30 40)\lvec(40 50)\lvec(40 40)\lvec(30 40)\htext(32
45){\tiny $1$}\htext(36 42){\tiny $0$}

\move(40 40)\lvec(50 50)\lvec(50 40)\lvec(40 40)\htext(42
45){\tiny $0$}\htext(46 42){\tiny $1$}
\end{texdraw}}\quad ,
\vskip 10mm

$f_2$ \raisebox{-0.5\height}{\begin{texdraw}\drawdim em
\setunitscale 0.13 \linewd 0.5

\move(10 0)\lvec(20 0)\lvec(20 10)\lvec(10 10)\lvec(10 0)

\move(20 0)\lvec(30 0)\lvec(30 10)\lvec(20 10)\lvec(20 0)

\move(30 0)\lvec(40 0)\lvec(40 10)\lvec(30 10)\lvec(30 0)

\move(40 0)\lvec(50 0)\lvec(50 10)\lvec(40 10)\lvec(40 0)

\move(10 10)\lvec(20 10)\lvec(20 20)\lvec(10 20)\lvec(10
10)\htext(13 13){$_2$}

\move(20 10)\lvec(30 10)\lvec(30 20)\lvec(20 20)\lvec(20
10)\htext(23 13){$_2$}

\move(30 10)\lvec(40 10)\lvec(40 20)\lvec(30 20)\lvec(30
10)\htext(33 13){$_2$}

\move(40 10)\lvec(50 10)\lvec(50 20)\lvec(40 20)\lvec(40
10)\htext(43 13){$_2$}

\move(20 20)\lvec(30 20)\lvec(30 30)\lvec(20 30)\lvec(20
20)\htext(23 23){$_3$}

\move(30 20)\lvec(40 20)\lvec(40 30)\lvec(30 30)\lvec(30
20)\htext(33 23){$_3$}

\move(40 20)\lvec(50 20)\lvec(50 30)\lvec(40 30)\lvec(40
20)\htext(43 23){$_3$}

\move(30 30)\lvec(40 30)\lvec(40 40)\lvec(30 40)\lvec(30
30)\htext(33 33){$_2$}

\move(40 30)\lvec(50 30)\lvec(50 40)\lvec(40 40)\lvec(40
30)\htext(43 33){$_2$}

\move(30 40)\lvec(40 40)\lvec(40 50)\lvec(30 50)\lvec(30 40)

\move(40 40)\lvec(50 40)\lvec(50 50)\lvec(40 50)\lvec(40 40)

\move(40 50)\lvec(50 50)\lvec(50 60)\lvec(40 60)\lvec(40
50)\htext(43 53){$_2$}

\move(40 60)\lvec(50 60)\lvec(50 70)\lvec(40 70)\lvec(40
60)\htext(43 63){$_3$}

\move(10 0)\lvec(20 10)\lvec(20 0)\lvec(10 0)\lfill f:0.8
\htext(11 5){\tiny $1$}\htext(16 1){\tiny $0$}

\move(20 0)\lvec(30 10)\lvec(30 0)\lvec(20 0)\lfill f:0.8
\htext(22 5){\tiny $0$}\htext(26 2){\tiny $1$}

\move(30 0)\lvec(40 10)\lvec(40 0)\lvec(30 0)\lfill f:0.8
\htext(32 5){\tiny $1$}\htext(36 2){\tiny $0$}

\move(40 0)\lvec(50 10)\lvec(50 0)\lvec(40 0)\lfill f:0.8
\htext(42 5){\tiny $0$}\htext(46 2){\tiny $1$}

\move(30 40)\lvec(40 50)\lvec(40 40)\lvec(30 40)\htext(32
45){\tiny $1$}\htext(36 42){\tiny $0$}

\move(40 40)\lvec(50 50)\lvec(50 40)\lvec(40 40)\htext(42
45){\tiny $0$}\htext(46 42){\tiny $1$}
\end{texdraw}}\vskip 5mm

$=$\hskip 5mm $q^{-1}$
\raisebox{-0.5\height}{\begin{texdraw}\drawdim em \setunitscale
0.13 \linewd 0.5

\move(10 0)\lvec(20 0)\lvec(20 10)\lvec(10 10)\lvec(10 0)

\move(20 0)\lvec(30 0)\lvec(30 10)\lvec(20 10)\lvec(20 0)

\move(30 0)\lvec(40 0)\lvec(40 10)\lvec(30 10)\lvec(30 0)

\move(40 0)\lvec(50 0)\lvec(50 10)\lvec(40 10)\lvec(40 0)

\move(10 10)\lvec(20 10)\lvec(20 20)\lvec(10 20)\lvec(10
10)\htext(13 13){$_2$}

\move(20 10)\lvec(30 10)\lvec(30 20)\lvec(20 20)\lvec(20
10)\htext(23 13){$_2$}

\move(30 10)\lvec(40 10)\lvec(40 20)\lvec(30 20)\lvec(30
10)\htext(33 13){$_2$}

\move(40 10)\lvec(50 10)\lvec(50 20)\lvec(40 20)\lvec(40
10)\htext(43 13){$_2$}

\move(20 20)\lvec(30 20)\lvec(30 30)\lvec(20 30)\lvec(20
20)\htext(23 23){$_3$}

\move(30 20)\lvec(40 20)\lvec(40 30)\lvec(30 30)\lvec(30
20)\htext(33 23){$_3$}

\move(40 20)\lvec(50 20)\lvec(50 30)\lvec(40 30)\lvec(40
20)\htext(43 23){$_3$}

\move(20 30)\lvec(30 30)\lvec(30 40)\lvec(20 40)\lvec(20
30)\htext(23 33){$_2$}

\move(30 30)\lvec(40 30)\lvec(40 40)\lvec(30 40)\lvec(30
30)\htext(33 33){$_2$}

\move(40 30)\lvec(50 30)\lvec(50 40)\lvec(40 40)\lvec(40
30)\htext(43 33){$_2$}

\move(30 40)\lvec(40 40)\lvec(40 50)\lvec(30 50)\lvec(30 40)

\move(40 40)\lvec(50 40)\lvec(50 50)\lvec(40 50)\lvec(40 40)

\move(40 50)\lvec(50 50)\lvec(50 60)\lvec(40 60)\lvec(40
50)\htext(43 53){$_2$}

\move(40 60)\lvec(50 60)\lvec(50 70)\lvec(40 70)\lvec(40
60)\htext(43 63){$_3$}

\move(10 0)\lvec(20 10)\lvec(20 0)\lvec(10 0)\lfill f:0.8
\htext(11 5){\tiny $1$}\htext(16 1){\tiny $0$}

\move(20 0)\lvec(30 10)\lvec(30 0)\lvec(20 0)\lfill f:0.8
\htext(22 5){\tiny $0$}\htext(26 2){\tiny $1$}

\move(30 0)\lvec(40 10)\lvec(40 0)\lvec(30 0)\lfill f:0.8
\htext(32 5){\tiny $1$}\htext(36 2){\tiny $0$}

\move(40 0)\lvec(50 10)\lvec(50 0)\lvec(40 0)\lfill f:0.8
\htext(42 5){\tiny $0$}\htext(46 2){\tiny $1$}

\move(30 40)\lvec(40 50)\lvec(40 40)\lvec(30 40)\htext(32
45){\tiny $1$}\htext(36 42){\tiny $0$}

\move(40 40)\lvec(50 50)\lvec(50 40)\lvec(40 40)\htext(42
45){\tiny $0$}\htext(46 42){\tiny $1$}
\end{texdraw}}\hskip 5mm $+$\hskip 5mm
\raisebox{-0.5\height}{\begin{texdraw}\drawdim em \setunitscale
0.13 \linewd 0.5

\move(10 0)\lvec(20 0)\lvec(20 10)\lvec(10 10)\lvec(10 0)

\move(20 0)\lvec(30 0)\lvec(30 10)\lvec(20 10)\lvec(20 0)

\move(30 0)\lvec(40 0)\lvec(40 10)\lvec(30 10)\lvec(30 0)

\move(40 0)\lvec(50 0)\lvec(50 10)\lvec(40 10)\lvec(40 0)

\move(10 10)\lvec(20 10)\lvec(20 20)\lvec(10 20)\lvec(10
10)\htext(13 13){$_2$}

\move(20 10)\lvec(30 10)\lvec(30 20)\lvec(20 20)\lvec(20
10)\htext(23 13){$_2$}

\move(30 10)\lvec(40 10)\lvec(40 20)\lvec(30 20)\lvec(30
10)\htext(33 13){$_2$}

\move(40 10)\lvec(50 10)\lvec(50 20)\lvec(40 20)\lvec(40
10)\htext(43 13){$_2$}

\move(20 20)\lvec(30 20)\lvec(30 30)\lvec(20 30)\lvec(20
20)\htext(23 23){$_3$}

\move(30 20)\lvec(40 20)\lvec(40 30)\lvec(30 30)\lvec(30
20)\htext(33 23){$_3$}

\move(40 20)\lvec(50 20)\lvec(50 30)\lvec(40 30)\lvec(40
20)\htext(43 23){$_3$}

\move(30 30)\lvec(40 30)\lvec(40 40)\lvec(30 40)\lvec(30
30)\htext(33 33){$_2$}

\move(40 30)\lvec(50 30)\lvec(50 40)\lvec(40 40)\lvec(40
30)\htext(43 33){$_2$}

\move(30 40)\lvec(40 40)\lvec(40 50)\lvec(30 50)\lvec(30 40)

\move(40 40)\lvec(50 40)\lvec(50 50)\lvec(40 50)\lvec(40 40)

\move(30 50)\lvec(40 50)\lvec(40 60)\lvec(30 60)\lvec(30
50)\htext(33 53){$_2$}

\move(40 50)\lvec(50 50)\lvec(50 60)\lvec(40 60)\lvec(40
50)\htext(43 53){$_2$}

\move(40 60)\lvec(50 60)\lvec(50 70)\lvec(40 70)\lvec(40
60)\htext(43 63){$_3$}

\move(10 0)\lvec(20 10)\lvec(20 0)\lvec(10 0)\lfill f:0.8
\htext(11 5){\tiny $1$}\htext(16 1){\tiny $0$}

\move(20 0)\lvec(30 10)\lvec(30 0)\lvec(20 0)\lfill f:0.8
\htext(22 5){\tiny $0$}\htext(26 2){\tiny $1$}

\move(30 0)\lvec(40 10)\lvec(40 0)\lvec(30 0)\lfill f:0.8
\htext(32 5){\tiny $1$}\htext(36 2){\tiny $0$}

\move(40 0)\lvec(50 10)\lvec(50 0)\lvec(40 0)\lfill f:0.8
\htext(42 5){\tiny $0$}\htext(46 2){\tiny $1$}

\move(30 40)\lvec(40 50)\lvec(40 40)\lvec(30 40)\htext(32
45){\tiny $1$}\htext(36 42){\tiny $0$}

\move(40 40)\lvec(50 50)\lvec(50 40)\lvec(40 40)\htext(42
45){\tiny $0$}\htext(46 42){\tiny $1$}
\end{texdraw}}
\hskip 5mm $+$\hskip 5mm $q$ \raisebox{-0.5\height}{
\begin{texdraw}\drawdim
em \setunitscale 0.13 \linewd 0.5

\move(10 0)\lvec(20 0)\lvec(20 10)\lvec(10 10)\lvec(10 0)

\move(20 0)\lvec(30 0)\lvec(30 10)\lvec(20 10)\lvec(20 0)

\move(30 0)\lvec(40 0)\lvec(40 10)\lvec(30 10)\lvec(30 0)

\move(40 0)\lvec(50 0)\lvec(50 10)\lvec(40 10)\lvec(40 0)

\move(10 10)\lvec(20 10)\lvec(20 20)\lvec(10 20)\lvec(10
10)\htext(13 13){$_2$}

\move(20 10)\lvec(30 10)\lvec(30 20)\lvec(20 20)\lvec(20
10)\htext(23 13){$_2$}

\move(30 10)\lvec(40 10)\lvec(40 20)\lvec(30 20)\lvec(30
10)\htext(33 13){$_2$}

\move(40 10)\lvec(50 10)\lvec(50 20)\lvec(40 20)\lvec(40
10)\htext(43 13){$_2$}

\move(20 20)\lvec(30 20)\lvec(30 30)\lvec(20 30)\lvec(20
20)\htext(23 23){$_3$}

\move(30 20)\lvec(40 20)\lvec(40 30)\lvec(30 30)\lvec(30
20)\htext(33 23){$_3$}

\move(40 20)\lvec(50 20)\lvec(50 30)\lvec(40 30)\lvec(40
20)\htext(43 23){$_3$}

\move(30 30)\lvec(40 30)\lvec(40 40)\lvec(30 40)\lvec(30
30)\htext(33 33){$_2$}

\move(40 30)\lvec(50 30)\lvec(50 40)\lvec(40 40)\lvec(40
30)\htext(43 33){$_2$}

\move(30 40)\lvec(40 40)\lvec(40 50)\lvec(30 50)\lvec(30 40)

\move(40 40)\lvec(50 40)\lvec(50 50)\lvec(40 50)\lvec(40 40)


\move(40 50)\lvec(50 50)\lvec(50 60)\lvec(40 60)\lvec(40
50)\htext(43 53){$_2$}

\move(40 60)\lvec(50 60)\lvec(50 70)\lvec(40 70)\lvec(40
60)\htext(43 63){$_3$}

\move(40 70)\lvec(50 70)\lvec(50 80)\lvec(40 80)\lvec(40
70)\htext(43 73){$_2$}

\move(10 0)\lvec(20 10)\lvec(20 0)\lvec(10 0)\lfill f:0.8
\htext(11 5){\tiny $1$}\htext(16 1){\tiny $0$}

\move(20 0)\lvec(30 10)\lvec(30 0)\lvec(20 0)\lfill f:0.8
\htext(22 5){\tiny $0$}\htext(26 2){\tiny $1$}

\move(30 0)\lvec(40 10)\lvec(40 0)\lvec(30 0)\lfill f:0.8
\htext(32 5){\tiny $1$}\htext(36 2){\tiny $0$}

\move(40 0)\lvec(50 10)\lvec(50 0)\lvec(40 0)\lfill f:0.8
\htext(42 5){\tiny $0$}\htext(46 2){\tiny $1$}

\move(30 40)\lvec(40 50)\lvec(40 40)\lvec(30 40)\htext(32
45){\tiny $1$}\htext(36 42){\tiny $0$}

\move(40 40)\lvec(50 50)\lvec(50 40)\lvec(40 40)\htext(42
45){\tiny $0$}\htext(46 42){\tiny $1$}
\end{texdraw}}\quad .
\vskip 5mm }
\end{ex}

\vskip 2mm \noindent {\bf Case 2.} Suppose that the $i$-blocks are
of type II.\vskip 2mm

In this case, we have $q_i=q$. Let $b$ be a removable $i$-block in
$y_k$ of $Y$. If the $i$-signature of $y_k$ is $--$, or if the
$i$-signature of $y_k$ is $-$ and there is another $i$-block
beneath $b$, we define $Y\nearrow b$ to be the Young wall obtained
by removing the block $b$ from $Y$. If the $i$-signature of $y_k$
is $-+$, or if the $i$-signature of $y_k$ is $-$ and there is no
$i$-block beneath $b$, then we define $Y\nearrow
b=q^{-1}(1-(-q^2)^{l(b)+1})Z$, where $Z$ is the Young wall
obtained by removing the block $b$ from $Y$ and $l(b)$ is the
number of $y_l$'s with $l<k$ such that $|y_l|=|y_k|$. That is, if

\vskip 5mm \hskip 4cm $Y=$\raisebox{-0.5\height}{
\begin{texdraw}\fontsize{9}{9}
\drawdim em \setunitscale 0.1 \linewd 0.5

\move(-20 -20)\lvec(-20 0)\lvec(0 0)\lvec(0 10)\lvec(50
10)\lvec(50 30)\lvec(60 30)\lvec(60 40)\lvec(80 40)\lvec(80
-20)\lvec(-20 -20)

\move(50 10)\lvec(50 -20)

\htext(57 10){$Y_R(b)$}

\move(0 10) \lvec(0 15) \lvec(10 15)\lvec(10 10)\lvec(0 10)\lfill
f:0.8

\move(10 10)\lvec(10 15)\lvec(50 15)\lvec(50 10)\lvec(10 10)
\lfill f:0.8

\move(20 10) \lvec(20 15) \move(40 10) \lvec(40 15) \htext(25
11){$\cdots$}

\move(10 0) \arrowheadtype t:F \arrowheadsize l:4 w:2 \avec(3 12)

\htext(12 -1){$_b$}

\move(10 15)\clvec(12 20)(23 20)(25 20)

\move(50 15)\clvec(48 20)(37 20)(35 20)

\htext(26 18){$_{l(b)}$}
\end{texdraw}}\quad ,\vskip 5mm

\hskip 5mm then \hskip 5mm $Y\nearrow b=$\raisebox{-0.5\height}{
\begin{texdraw}\fontsize{9}{9}
\drawdim em \setunitscale 0.1 \linewd 0.5

\move(-20 -20)\lvec(-20 0)\lvec(0 0)\lvec(0 10)\lvec(50
10)\lvec(50 30)\lvec(60 30)\lvec(60 40)\lvec(80 40)\lvec(80
-20)\lvec(-20 -20)

\move(10 10)\lvec(10 15)\lvec(50 15)\lvec(50 10)\lvec(10 10)
\lfill f:0.8

\move(20 10) \lvec(20 15) \move(40 10) \lvec(40 15) \htext(25
11){$\cdots$}

\htext(-90 3){$\dfrac{(1-(-q^2)^{l(b)+1})}{q}\times$}
\end{texdraw}}\quad .
\vskip 5mm

In either case, we define $Y_R(b)=(y_{l},\cdots,y_{0})$, where $l$
is the integer such that $|y_k|=|y_{k-1}|=\cdots=|y_{l+1}|<|y_l|$,
and set
\begin{equation}
R_i(b;Y)=\varphi_i(Y_R(b))-\varepsilon_i(Y_R(b)).
\end{equation}
If $k=0$, we understand $Y_R(b)=\emptyset$ and $R_i(b;Y)=0$. Then
we define
\begin{equation}
e_i\,Y=\sum_{b}q_i^{-R_i(b;Y)}(Y\nearrow b),
\end{equation}
\vskip -2mm \noindent where $b$ runs over all removable $i$-blocks
in $Y$.

On the other hand, suppose that $b$ is an admissible $i$-slot in
$y_k$ of $Y$. If the $i$-signature of $y_k$ is $++$, or if the
$i$-signature of $y_k$ is $+$ and there is no $i$-block beneath
$b$, then we define $Y\swarrow b$ to be the Young wall obtained by
adding an $i$-block at $b$. If the $i$-signature of $y_k$ is $-+$,
or if the $i$-signature of $y_k$ is $+$ and there is another
$i$-block beneath $b$, then we define $Y\swarrow
b=q^{-1}(1-(-q^2)^{l(b)+1})Z$, where $Z$ is the Young wall
obtained by adding an $i$-block at $b$ and $l(b)$ is the number of
$y_l$'s with $l>k$ such that $|y_l|=|y_k|$. That is, if \vskip 5mm

\hskip 4cm $Y=$\raisebox{-0.5\height}{
\begin{texdraw}\fontsize{9}{9}
\drawdim em \setunitscale 0.1 \linewd 0.5

\move(-30 -20)\lvec(-30 0)\lvec(-10 0)\lvec(-10 10)\lvec(50
10)\lvec(50 30)\lvec(65 30)\lvec(65 -20)\lvec(-30 -20)

\move(0 10)\lvec(0 -20)\htext(-24 -14){$Y_L(b)$}

\move(40 15) \lvec(40 20)\lvec(50 20)

\move(0 10) \lvec(0 15) \lvec(10 15)\lvec(10 10)\lvec(0 10)\lfill
f:0.8

\move(10 10)\lvec(10 15)\lvec(50 15)\lvec(50 10)\lvec(10 10)
\lfill f:0.8

\move(30 10) \lvec(30 15) \move(40 10) \lvec(40 15) \htext(15
11){$\cdots$}

\move(55 0) \arrowheadtype t:F \arrowheadsize l:4 w:2 \avec(43 17)

\htext(60 -3){$_b$}

\move(0 15)\clvec(2 20)(13 20)(15 20)

\move(40 15)\clvec(38 20)(27 20)(25 20)

\htext(16 18){$_{l(b)}$}
\end{texdraw}}\quad ,
\vskip 0.5cm

\hskip 5mm then \hskip 0.5cm $Y\swarrow b=$\raisebox{-0.5\height}{
\begin{texdraw}\fontsize{9}{9}
\drawdim em \setunitscale 0.1 \linewd 0.5

\move(-30 -20)\lvec(-30 0)\lvec(-10 0)\lvec(-10 10)\lvec(50
10)\lvec(50 30)\lvec(65 30)\lvec(65 -20)\lvec(-30 -20)

\move(40 15) \lvec(40 20)\lvec(50 20)\lvec(50 15)\lfill f:0.8

\move(0 10) \lvec(0 15) \lvec(10 15)\lvec(10 10)\lvec(0 10)\lfill
f:0.8

\move(10 10)\lvec(10 15)\lvec(50 15)\lvec(50 10)\lvec(10 10)
\lfill f:0.8

\move(30 10) \lvec(30 15) \move(40 10) \lvec(40 15) \htext(15
11){$\cdots$}

\htext(-100 0){$\dfrac{(1-(-q^2)^{l(b)+1})}{q}\times$}
\end{texdraw}}\quad .\vskip 5mm

In either case, we define $Y_L(b)=(\cdots,y_{l+2},y_{l+1})$, where
$l$ is the integer such that
$|y_{l+1}|<|y_l|=|y_{l-1}|=\cdots=|y_k|$, and set
\begin{equation}
L_i(b;Y)=\varphi_i(Y_L(b))-\varepsilon_i(Y_L(b)).
\end{equation}
Then we define
\begin{equation}
f_i\,Y=\sum_{b}q_i^{L_i(b;Y)}(Y\swarrow b),
\end{equation}
\vskip -2mm \noindent where $b$ runs over all admissible $i$-slots
in $Y$.
\begin{ex}{\rm
If $\frak{g}=A^{(2)}_4$, $\Lambda=\Lambda_0$ and $i=0$, then we
have\vskip 5mm

$e_0$\raisebox{-0.5\height}{
\begin{texdraw}
\drawdim em \setunitscale 0.13 \linewd 0.5

\move(-10 0)\lvec(0 0)\lvec(0 10)\lvec(-10 10)\lvec(-10
0)\htext(-7 1){\tiny $0$}

\move(0 0)\lvec(10 0)\lvec(10 10)\lvec(0 10)\lvec(0 0)\htext(3
1){\tiny $0$}

\move(10 0)\lvec(20 0)\lvec(20 10)\lvec(10 10)\lvec(10 0)\htext(13
1){\tiny $0$}

\move(20 0)\lvec(30 0)\lvec(30 10)\lvec(20 10)\lvec(20 0)\htext(23
1){\tiny $0$}

\move(30 0)\lvec(40 0)\lvec(40 10)\lvec(30 10)\lvec(30 0)\htext(33
1){\tiny $0$}

\move(40 0)\lvec(50 0)\lvec(50 10)\lvec(40 10)\lvec(40 0)\htext(43
1){\tiny $0$}

\move(50 0)\lvec(60 0)\lvec(60 10)\lvec(50 10)\lvec(50 0)\htext(53
1){\tiny $0$}
\move(-10 0)\lvec(0 0)\lvec(0 5)\lvec(-10 5)\lvec(-10 0)\lfill
f:0.8 \htext(-7 6){\tiny $0$}

\move(0 0)\lvec(10 0)\lvec(10 5)\lvec(0 5)\lvec(0 0)\lfill f:0.8
\htext(3 6){\tiny $0$}

\move(10 0)\lvec(20 0)\lvec(20 5)\lvec(10 5)\lvec(10 0)\lfill
f:0.8 \htext(13 6){\tiny $0$}

\move(20 0)\lvec(30 0)\lvec(30 5)\lvec(20 5)\lvec(20 0)\lfill
f:0.8 \htext(23 6){\tiny $0$}

\move(30 0)\lvec(40 0)\lvec(40 5)\lvec(30 5)\lvec(30 0)\lfill
f:0.8 \htext(33 6){\tiny $0$}

\move(40 0)\lvec(50 0)\lvec(50 5)\lvec(40 5)\lvec(40 0)\lfill
f:0.8 \htext(43 6){\tiny $0$}

\move(50 0)\lvec(60 0)\lvec(60 5)\lvec(50 5)\lvec(50 0)\lfill
f:0.8 \htext(53 6){\tiny $0$}

\move(0 10)\lvec(10 10)\lvec(10 20)\lvec(0 20)\lvec(0 10)\htext(3
13){$_1$}

\move(10 10)\lvec(20 10)\lvec(20 20)\lvec(10 20)\lvec(10
10)\htext(13 13){$_1$}

\move(20 10)\lvec(30 10)\lvec(30 20)\lvec(20 20)\lvec(20
10)\htext(23 13){$_1$}

\move(30 10)\lvec(40 10)\lvec(40 20)\lvec(30 20)\lvec(30
10)\htext(33 13){$_1$}

\move(40 10)\lvec(50 10)\lvec(50 20)\lvec(40 20)\lvec(40
10)\htext(43 13){$_1$}

\move(50 10)\lvec(60 10)\lvec(60 20)\lvec(50 20)\lvec(50
10)\htext(53 13){$_1$}
\move(0 20)\lvec(10 20)\lvec(10 30)\lvec(0 30)\lvec(0 20)\htext(3
23){$_2$}

\move(10 20)\lvec(20 20)\lvec(20 30)\lvec(10 30)\lvec(10
20)\htext(13 23){$_2$}

\move(20 20)\lvec(30 20)\lvec(30 30)\lvec(20 30)\lvec(20
20)\htext(23 23){$_2$}

\move(30 20)\lvec(40 20)\lvec(40 30)\lvec(30 30)\lvec(30
20)\htext(33 23){$_2$}

\move(40 20)\lvec(50 20)\lvec(50 30)\lvec(40 30)\lvec(40
20)\htext(43 23){$_2$}

\move(50 20)\lvec(60 20)\lvec(60 30)\lvec(50 30)\lvec(50
20)\htext(53 23){$_2$}

\move(0 30)\lvec(10 30)\lvec(10 40)\lvec(0 40)\lvec(0 30)\htext(3
33){$_1$}

\move(10 30)\lvec(20 30)\lvec(20 40)\lvec(10 40)\lvec(10
30)\htext(13 33){$_1$}

\move(20 30)\lvec(30 30)\lvec(30 40)\lvec(20 40)\lvec(20
30)\htext(23 33){$_1$}

\move(30 30)\lvec(40 30)\lvec(40 40)\lvec(30 40)\lvec(30
30)\htext(33 33){$_1$}

\move(40 30)\lvec(50 30)\lvec(50 40)\lvec(40 40)\lvec(40
30)\htext(43 33){$_1$}

\move(50 30)\lvec(60 30)\lvec(60 40)\lvec(50 40)\lvec(50
30)\htext(53 33){$_1$}
\move(0 40)\lvec(10 40)\lvec(10 45)\lvec(0 45)\lvec(0 40)\htext(3
41){\tiny $0$}

\move(10 40)\lvec(20 40)\lvec(20 45)\lvec(10 45)\lvec(10
40)\htext(13 41){\tiny $0$}

\move(20 40)\lvec(30 40)\lvec(30 45)\lvec(20 45)\lvec(20
40)\htext(23 41){\tiny $0$}

\move(30 40)\lvec(40 40)\lvec(40 45)\lvec(30 45)\lvec(30
40)\htext(33 41){\tiny $0$}

\move(40 40)\lvec(50 40)\lvec(50 45)\lvec(40 45)\lvec(40
40)\htext(43 41){\tiny $0$}

\move(50 40)\lvec(60 40)\lvec(60 45)\lvec(50 45)\lvec(50
40)\htext(53 41){\tiny $0$}

\move(50 45)\lvec(60 45)\lvec(60 50)\lvec(50 50)\lvec(50
45)\htext(53 46){\tiny $0$}

\end{texdraw}}\hskip 3mm $=$\hskip 3mm $q^{2}$
\raisebox{-0.5\height}{
\begin{texdraw}
\drawdim em \setunitscale 0.13 \linewd 0.5

\move(0 0)\lvec(10 0)\lvec(10 10)\lvec(0 10)\lvec(0 0)\htext(3
1){\tiny $0$}

\move(10 0)\lvec(20 0)\lvec(20 10)\lvec(10 10)\lvec(10 0)\htext(13
1){\tiny $0$}

\move(20 0)\lvec(30 0)\lvec(30 10)\lvec(20 10)\lvec(20 0)\htext(23
1){\tiny $0$}

\move(30 0)\lvec(40 0)\lvec(40 10)\lvec(30 10)\lvec(30 0)\htext(33
1){\tiny $0$}

\move(40 0)\lvec(50 0)\lvec(50 10)\lvec(40 10)\lvec(40 0)\htext(43
1){\tiny $0$}

\move(50 0)\lvec(60 0)\lvec(60 10)\lvec(50 10)\lvec(50 0)\htext(53
1){\tiny $0$}
\move(0 0)\lvec(10 0)\lvec(10 5)\lvec(0 5)\lvec(0 0)\lfill f:0.8
\htext(3 6){\tiny $0$}

\move(10 0)\lvec(20 0)\lvec(20 5)\lvec(10 5)\lvec(10 0)\lfill
f:0.8 \htext(13 6){\tiny $0$}

\move(20 0)\lvec(30 0)\lvec(30 5)\lvec(20 5)\lvec(20 0)\lfill
f:0.8 \htext(23 6){\tiny $0$}

\move(30 0)\lvec(40 0)\lvec(40 5)\lvec(30 5)\lvec(30 0)\lfill
f:0.8 \htext(33 6){\tiny $0$}

\move(40 0)\lvec(50 0)\lvec(50 5)\lvec(40 5)\lvec(40 0)\lfill
f:0.8 \htext(43 6){\tiny $0$}

\move(50 0)\lvec(60 0)\lvec(60 5)\lvec(50 5)\lvec(50 0)\lfill
f:0.8 \htext(53 6){\tiny $0$}

\move(0 10)\lvec(10 10)\lvec(10 20)\lvec(0 20)\lvec(0 10)\htext(3
13){$_1$}

\move(10 10)\lvec(20 10)\lvec(20 20)\lvec(10 20)\lvec(10
10)\htext(13 13){$_1$}

\move(20 10)\lvec(30 10)\lvec(30 20)\lvec(20 20)\lvec(20
10)\htext(23 13){$_1$}

\move(30 10)\lvec(40 10)\lvec(40 20)\lvec(30 20)\lvec(30
10)\htext(33 13){$_1$}

\move(40 10)\lvec(50 10)\lvec(50 20)\lvec(40 20)\lvec(40
10)\htext(43 13){$_1$}

\move(50 10)\lvec(60 10)\lvec(60 20)\lvec(50 20)\lvec(50
10)\htext(53 13){$_1$}
\move(0 20)\lvec(10 20)\lvec(10 30)\lvec(0 30)\lvec(0 20)\htext(3
23){$_2$}

\move(10 20)\lvec(20 20)\lvec(20 30)\lvec(10 30)\lvec(10
20)\htext(13 23){$_2$}

\move(20 20)\lvec(30 20)\lvec(30 30)\lvec(20 30)\lvec(20
20)\htext(23 23){$_2$}

\move(30 20)\lvec(40 20)\lvec(40 30)\lvec(30 30)\lvec(30
20)\htext(33 23){$_2$}

\move(40 20)\lvec(50 20)\lvec(50 30)\lvec(40 30)\lvec(40
20)\htext(43 23){$_2$}

\move(50 20)\lvec(60 20)\lvec(60 30)\lvec(50 30)\lvec(50
20)\htext(53 23){$_2$}

\move(0 30)\lvec(10 30)\lvec(10 40)\lvec(0 40)\lvec(0 30)\htext(3
33){$_1$}

\move(10 30)\lvec(20 30)\lvec(20 40)\lvec(10 40)\lvec(10
30)\htext(13 33){$_1$}

\move(20 30)\lvec(30 30)\lvec(30 40)\lvec(20 40)\lvec(20
30)\htext(23 33){$_1$}

\move(30 30)\lvec(40 30)\lvec(40 40)\lvec(30 40)\lvec(30
30)\htext(33 33){$_1$}

\move(40 30)\lvec(50 30)\lvec(50 40)\lvec(40 40)\lvec(40
30)\htext(43 33){$_1$}

\move(50 30)\lvec(60 30)\lvec(60 40)\lvec(50 40)\lvec(50
30)\htext(53 33){$_1$}
\move(0 40)\lvec(10 40)\lvec(10 45)\lvec(0 45)\lvec(0 40)\htext(3
41){\tiny $0$}

\move(10 40)\lvec(20 40)\lvec(20 45)\lvec(10 45)\lvec(10
40)\htext(13 41){\tiny $0$}

\move(20 40)\lvec(30 40)\lvec(30 45)\lvec(20 45)\lvec(20
40)\htext(23 41){\tiny $0$}

\move(30 40)\lvec(40 40)\lvec(40 45)\lvec(30 45)\lvec(30
40)\htext(33 41){\tiny $0$}

\move(40 40)\lvec(50 40)\lvec(50 45)\lvec(40 45)\lvec(40
40)\htext(43 41){\tiny $0$}

\move(50 40)\lvec(60 40)\lvec(60 45)\lvec(50 45)\lvec(50
40)\htext(53 41){\tiny $0$}

\move(50 45)\lvec(60 45)\lvec(60 50)\lvec(50 50)\lvec(50
45)\htext(53 46){\tiny $0$}

\end{texdraw}}\vskip 5mm

$+$ \hskip 3mm $q(1+q^{10})$ \raisebox{-0.5\height}{
\begin{texdraw}
\drawdim em \setunitscale 0.13 \linewd 0.5

\move(-10 0)\lvec(0 0)\lvec(0 10)\lvec(-10 10)\lvec(-10
0)\htext(-7 1){\tiny $0$}

\move(0 0)\lvec(10 0)\lvec(10 10)\lvec(0 10)\lvec(0 0)\htext(3
1){\tiny $0$}

\move(10 0)\lvec(20 0)\lvec(20 10)\lvec(10 10)\lvec(10 0)\htext(13
1){\tiny $0$}

\move(20 0)\lvec(30 0)\lvec(30 10)\lvec(20 10)\lvec(20 0)\htext(23
1){\tiny $0$}

\move(30 0)\lvec(40 0)\lvec(40 10)\lvec(30 10)\lvec(30 0)\htext(33
1){\tiny $0$}

\move(40 0)\lvec(50 0)\lvec(50 10)\lvec(40 10)\lvec(40 0)\htext(43
1){\tiny $0$}

\move(50 0)\lvec(60 0)\lvec(60 10)\lvec(50 10)\lvec(50 0)\htext(53
1){\tiny $0$}
\move(-10 0)\lvec(0 0)\lvec(0 5)\lvec(-10 5)\lvec(-10 0)\lfill
f:0.8 \htext(-7 6){\tiny $0$}

\move(0 0)\lvec(10 0)\lvec(10 5)\lvec(0 5)\lvec(0 0)\lfill f:0.8
\htext(3 6){\tiny $0$}

\move(10 0)\lvec(20 0)\lvec(20 5)\lvec(10 5)\lvec(10 0)\lfill
f:0.8 \htext(13 6){\tiny $0$}

\move(20 0)\lvec(30 0)\lvec(30 5)\lvec(20 5)\lvec(20 0)\lfill
f:0.8 \htext(23 6){\tiny $0$}

\move(30 0)\lvec(40 0)\lvec(40 5)\lvec(30 5)\lvec(30 0)\lfill
f:0.8 \htext(33 6){\tiny $0$}

\move(40 0)\lvec(50 0)\lvec(50 5)\lvec(40 5)\lvec(40 0)\lfill
f:0.8 \htext(43 6){\tiny $0$}

\move(50 0)\lvec(60 0)\lvec(60 5)\lvec(50 5)\lvec(50 0)\lfill
f:0.8 \htext(53 6){\tiny $0$}

\move(0 10)\lvec(10 10)\lvec(10 20)\lvec(0 20)\lvec(0 10)\htext(3
13){$_1$}

\move(10 10)\lvec(20 10)\lvec(20 20)\lvec(10 20)\lvec(10
10)\htext(13 13){$_1$}

\move(20 10)\lvec(30 10)\lvec(30 20)\lvec(20 20)\lvec(20
10)\htext(23 13){$_1$}

\move(30 10)\lvec(40 10)\lvec(40 20)\lvec(30 20)\lvec(30
10)\htext(33 13){$_1$}

\move(40 10)\lvec(50 10)\lvec(50 20)\lvec(40 20)\lvec(40
10)\htext(43 13){$_1$}

\move(50 10)\lvec(60 10)\lvec(60 20)\lvec(50 20)\lvec(50
10)\htext(53 13){$_1$}
\move(0 20)\lvec(10 20)\lvec(10 30)\lvec(0 30)\lvec(0 20)\htext(3
23){$_2$}

\move(10 20)\lvec(20 20)\lvec(20 30)\lvec(10 30)\lvec(10
20)\htext(13 23){$_2$}

\move(20 20)\lvec(30 20)\lvec(30 30)\lvec(20 30)\lvec(20
20)\htext(23 23){$_2$}

\move(30 20)\lvec(40 20)\lvec(40 30)\lvec(30 30)\lvec(30
20)\htext(33 23){$_2$}

\move(40 20)\lvec(50 20)\lvec(50 30)\lvec(40 30)\lvec(40
20)\htext(43 23){$_2$}

\move(50 20)\lvec(60 20)\lvec(60 30)\lvec(50 30)\lvec(50
20)\htext(53 23){$_2$}

\move(0 30)\lvec(10 30)\lvec(10 40)\lvec(0 40)\lvec(0 30)\htext(3
33){$_1$}

\move(10 30)\lvec(20 30)\lvec(20 40)\lvec(10 40)\lvec(10
30)\htext(13 33){$_1$}

\move(20 30)\lvec(30 30)\lvec(30 40)\lvec(20 40)\lvec(20
30)\htext(23 33){$_1$}

\move(30 30)\lvec(40 30)\lvec(40 40)\lvec(30 40)\lvec(30
30)\htext(33 33){$_1$}

\move(40 30)\lvec(50 30)\lvec(50 40)\lvec(40 40)\lvec(40
30)\htext(43 33){$_1$}

\move(50 30)\lvec(60 30)\lvec(60 40)\lvec(50 40)\lvec(50
30)\htext(53 33){$_1$}
\move(10 40)\lvec(20 40)\lvec(20 45)\lvec(10 45)\lvec(10
40)\htext(13 41){\tiny $0$}

\move(20 40)\lvec(30 40)\lvec(30 45)\lvec(20 45)\lvec(20
40)\htext(23 41){\tiny $0$}

\move(30 40)\lvec(40 40)\lvec(40 45)\lvec(30 45)\lvec(30
40)\htext(33 41){\tiny $0$}

\move(40 40)\lvec(50 40)\lvec(50 45)\lvec(40 45)\lvec(40
40)\htext(43 41){\tiny $0$}

\move(50 40)\lvec(60 40)\lvec(60 45)\lvec(50 45)\lvec(50
40)\htext(53 41){\tiny $0$}

\move(50 45)\lvec(60 45)\lvec(60 50)\lvec(50 50)\lvec(50
45)\htext(53 46){\tiny $0$}

\end{texdraw}}\hskip 5mm $+$ \hskip 5mm
\raisebox{-0.5\height}{
\begin{texdraw}
\drawdim em \setunitscale 0.13 \linewd 0.5

\move(-10 0)\lvec(0 0)\lvec(0 10)\lvec(-10 10)\lvec(-10
0)\htext(-7 1){\tiny $0$}

\move(0 0)\lvec(10 0)\lvec(10 10)\lvec(0 10)\lvec(0 0)\htext(3
1){\tiny $0$}

\move(10 0)\lvec(20 0)\lvec(20 10)\lvec(10 10)\lvec(10 0)\htext(13
1){\tiny $0$}

\move(20 0)\lvec(30 0)\lvec(30 10)\lvec(20 10)\lvec(20 0)\htext(23
1){\tiny $0$}

\move(30 0)\lvec(40 0)\lvec(40 10)\lvec(30 10)\lvec(30 0)\htext(33
1){\tiny $0$}

\move(40 0)\lvec(50 0)\lvec(50 10)\lvec(40 10)\lvec(40 0)\htext(43
1){\tiny $0$}

\move(50 0)\lvec(60 0)\lvec(60 10)\lvec(50 10)\lvec(50 0)\htext(53
1){\tiny $0$}
\move(-10 0)\lvec(0 0)\lvec(0 5)\lvec(-10 5)\lvec(-10 0)\lfill
f:0.8 \htext(-7 6){\tiny $0$}

\move(0 0)\lvec(10 0)\lvec(10 5)\lvec(0 5)\lvec(0 0)\lfill f:0.8
\htext(3 6){\tiny $0$}

\move(10 0)\lvec(20 0)\lvec(20 5)\lvec(10 5)\lvec(10 0)\lfill
f:0.8 \htext(13 6){\tiny $0$}

\move(20 0)\lvec(30 0)\lvec(30 5)\lvec(20 5)\lvec(20 0)\lfill
f:0.8 \htext(23 6){\tiny $0$}

\move(30 0)\lvec(40 0)\lvec(40 5)\lvec(30 5)\lvec(30 0)\lfill
f:0.8 \htext(33 6){\tiny $0$}

\move(40 0)\lvec(50 0)\lvec(50 5)\lvec(40 5)\lvec(40 0)\lfill
f:0.8 \htext(43 6){\tiny $0$}

\move(50 0)\lvec(60 0)\lvec(60 5)\lvec(50 5)\lvec(50 0)\lfill
f:0.8 \htext(53 6){\tiny $0$}

\move(0 10)\lvec(10 10)\lvec(10 20)\lvec(0 20)\lvec(0 10)\htext(3
13){$_1$}

\move(10 10)\lvec(20 10)\lvec(20 20)\lvec(10 20)\lvec(10
10)\htext(13 13){$_1$}

\move(20 10)\lvec(30 10)\lvec(30 20)\lvec(20 20)\lvec(20
10)\htext(23 13){$_1$}

\move(30 10)\lvec(40 10)\lvec(40 20)\lvec(30 20)\lvec(30
10)\htext(33 13){$_1$}

\move(40 10)\lvec(50 10)\lvec(50 20)\lvec(40 20)\lvec(40
10)\htext(43 13){$_1$}

\move(50 10)\lvec(60 10)\lvec(60 20)\lvec(50 20)\lvec(50
10)\htext(53 13){$_1$}
\move(0 20)\lvec(10 20)\lvec(10 30)\lvec(0 30)\lvec(0 20)\htext(3
23){$_2$}

\move(10 20)\lvec(20 20)\lvec(20 30)\lvec(10 30)\lvec(10
20)\htext(13 23){$_2$}

\move(20 20)\lvec(30 20)\lvec(30 30)\lvec(20 30)\lvec(20
20)\htext(23 23){$_2$}

\move(30 20)\lvec(40 20)\lvec(40 30)\lvec(30 30)\lvec(30
20)\htext(33 23){$_2$}

\move(40 20)\lvec(50 20)\lvec(50 30)\lvec(40 30)\lvec(40
20)\htext(43 23){$_2$}

\move(50 20)\lvec(60 20)\lvec(60 30)\lvec(50 30)\lvec(50
20)\htext(53 23){$_2$}

\move(0 30)\lvec(10 30)\lvec(10 40)\lvec(0 40)\lvec(0 30)\htext(3
33){$_1$}

\move(10 30)\lvec(20 30)\lvec(20 40)\lvec(10 40)\lvec(10
30)\htext(13 33){$_1$}

\move(20 30)\lvec(30 30)\lvec(30 40)\lvec(20 40)\lvec(20
30)\htext(23 33){$_1$}

\move(30 30)\lvec(40 30)\lvec(40 40)\lvec(30 40)\lvec(30
30)\htext(33 33){$_1$}

\move(40 30)\lvec(50 30)\lvec(50 40)\lvec(40 40)\lvec(40
30)\htext(43 33){$_1$}

\move(50 30)\lvec(60 30)\lvec(60 40)\lvec(50 40)\lvec(50
30)\htext(53 33){$_1$}
\move(0 40)\lvec(10 40)\lvec(10 45)\lvec(0 45)\lvec(0 40)\htext(3
41){\tiny $0$}

\move(10 40)\lvec(20 40)\lvec(20 45)\lvec(10 45)\lvec(10
40)\htext(13 41){\tiny $0$}

\move(20 40)\lvec(30 40)\lvec(30 45)\lvec(20 45)\lvec(20
40)\htext(23 41){\tiny $0$}

\move(30 40)\lvec(40 40)\lvec(40 45)\lvec(30 45)\lvec(30
40)\htext(33 41){\tiny $0$}

\move(40 40)\lvec(50 40)\lvec(50 45)\lvec(40 45)\lvec(40
40)\htext(43 41){\tiny $0$}

\move(50 40)\lvec(60 40)\lvec(60 45)\lvec(50 45)\lvec(50
40)\htext(53 41){\tiny $0$}

\end{texdraw}}\quad ,\vskip 10mm

$f_0$\raisebox{-0.5\height}{
\begin{texdraw}
\drawdim em \setunitscale 0.13 \linewd 0.5

\move(-10 0)\lvec(0 0)\lvec(0 10)\lvec(-10 10)\lvec(-10
0)\htext(-7 1){\tiny $0$}

\move(0 0)\lvec(10 0)\lvec(10 10)\lvec(0 10)\lvec(0 0)\htext(3
1){\tiny $0$}

\move(10 0)\lvec(20 0)\lvec(20 10)\lvec(10 10)\lvec(10 0)\htext(13
1){\tiny $0$}

\move(20 0)\lvec(30 0)\lvec(30 10)\lvec(20 10)\lvec(20 0)\htext(23
1){\tiny $0$}

\move(30 0)\lvec(40 0)\lvec(40 10)\lvec(30 10)\lvec(30 0)\htext(33
1){\tiny $0$}

\move(40 0)\lvec(50 0)\lvec(50 10)\lvec(40 10)\lvec(40 0)\htext(43
1){\tiny $0$}

\move(50 0)\lvec(60 0)\lvec(60 10)\lvec(50 10)\lvec(50 0)\htext(53
1){\tiny $0$}
\move(-10 0)\lvec(0 0)\lvec(0 5)\lvec(-10 5)\lvec(-10 0)\lfill
f:0.8 \htext(-7 6){\tiny $0$}

\move(0 0)\lvec(10 0)\lvec(10 5)\lvec(0 5)\lvec(0 0)\lfill f:0.8
\htext(3 6){\tiny $0$}

\move(10 0)\lvec(20 0)\lvec(20 5)\lvec(10 5)\lvec(10 0)\lfill
f:0.8 \htext(13 6){\tiny $0$}

\move(20 0)\lvec(30 0)\lvec(30 5)\lvec(20 5)\lvec(20 0)\lfill
f:0.8 \htext(23 6){\tiny $0$}

\move(30 0)\lvec(40 0)\lvec(40 5)\lvec(30 5)\lvec(30 0)\lfill
f:0.8 \htext(33 6){\tiny $0$}

\move(40 0)\lvec(50 0)\lvec(50 5)\lvec(40 5)\lvec(40 0)\lfill
f:0.8 \htext(43 6){\tiny $0$}

\move(50 0)\lvec(60 0)\lvec(60 5)\lvec(50 5)\lvec(50 0)\lfill
f:0.8 \htext(53 6){\tiny $0$}

\move(0 10)\lvec(10 10)\lvec(10 20)\lvec(0 20)\lvec(0 10)\htext(3
13){$_1$}

\move(10 10)\lvec(20 10)\lvec(20 20)\lvec(10 20)\lvec(10
10)\htext(13 13){$_1$}

\move(20 10)\lvec(30 10)\lvec(30 20)\lvec(20 20)\lvec(20
10)\htext(23 13){$_1$}

\move(30 10)\lvec(40 10)\lvec(40 20)\lvec(30 20)\lvec(30
10)\htext(33 13){$_1$}

\move(40 10)\lvec(50 10)\lvec(50 20)\lvec(40 20)\lvec(40
10)\htext(43 13){$_1$}

\move(50 10)\lvec(60 10)\lvec(60 20)\lvec(50 20)\lvec(50
10)\htext(53 13){$_1$}
\move(0 20)\lvec(10 20)\lvec(10 30)\lvec(0 30)\lvec(0 20)\htext(3
23){$_2$}

\move(10 20)\lvec(20 20)\lvec(20 30)\lvec(10 30)\lvec(10
20)\htext(13 23){$_2$}

\move(20 20)\lvec(30 20)\lvec(30 30)\lvec(20 30)\lvec(20
20)\htext(23 23){$_2$}

\move(30 20)\lvec(40 20)\lvec(40 30)\lvec(30 30)\lvec(30
20)\htext(33 23){$_2$}

\move(40 20)\lvec(50 20)\lvec(50 30)\lvec(40 30)\lvec(40
20)\htext(43 23){$_2$}

\move(50 20)\lvec(60 20)\lvec(60 30)\lvec(50 30)\lvec(50
20)\htext(53 23){$_2$}

\move(0 30)\lvec(10 30)\lvec(10 40)\lvec(0 40)\lvec(0 30)\htext(3
33){$_1$}

\move(10 30)\lvec(20 30)\lvec(20 40)\lvec(10 40)\lvec(10
30)\htext(13 33){$_1$}

\move(20 30)\lvec(30 30)\lvec(30 40)\lvec(20 40)\lvec(20
30)\htext(23 33){$_1$}

\move(30 30)\lvec(40 30)\lvec(40 40)\lvec(30 40)\lvec(30
30)\htext(33 33){$_1$}

\move(40 30)\lvec(50 30)\lvec(50 40)\lvec(40 40)\lvec(40
30)\htext(43 33){$_1$}

\move(50 30)\lvec(60 30)\lvec(60 40)\lvec(50 40)\lvec(50
30)\htext(53 33){$_1$}
\move(10 40)\lvec(20 40)\lvec(20 45)\lvec(10 45)\lvec(10
40)\htext(13 41){\tiny $0$}

\move(20 40)\lvec(30 40)\lvec(30 45)\lvec(20 45)\lvec(20
40)\htext(23 41){\tiny $0$}

\move(30 40)\lvec(40 40)\lvec(40 45)\lvec(30 45)\lvec(30
40)\htext(33 41){\tiny $0$}

\move(40 40)\lvec(50 40)\lvec(50 45)\lvec(40 45)\lvec(40
40)\htext(43 41){\tiny $0$}

\move(50 40)\lvec(60 40)\lvec(60 45)\lvec(50 45)\lvec(50
40)\htext(53 41){\tiny $0$}

\end{texdraw}}\vskip 5mm

$=$\hskip 3mm $q^{-1}$ \raisebox{-0.5\height}{
\begin{texdraw}
\drawdim em \setunitscale 0.13 \linewd 0.5

\move(-10 0)\lvec(0 0)\lvec(0 10)\lvec(-10 10)\lvec(-10
0)\htext(-7 1){\tiny $0$}

\move(0 0)\lvec(10 0)\lvec(10 10)\lvec(0 10)\lvec(0 0)\htext(3
1){\tiny $0$}

\move(10 0)\lvec(20 0)\lvec(20 10)\lvec(10 10)\lvec(10 0)\htext(13
1){\tiny $0$}

\move(20 0)\lvec(30 0)\lvec(30 10)\lvec(20 10)\lvec(20 0)\htext(23
1){\tiny $0$}

\move(30 0)\lvec(40 0)\lvec(40 10)\lvec(30 10)\lvec(30 0)\htext(33
1){\tiny $0$}

\move(40 0)\lvec(50 0)\lvec(50 10)\lvec(40 10)\lvec(40 0)\htext(43
1){\tiny $0$}

\move(50 0)\lvec(60 0)\lvec(60 10)\lvec(50 10)\lvec(50 0)\htext(53
1){\tiny $0$}
\move(-10 0)\lvec(0 0)\lvec(0 5)\lvec(-10 5)\lvec(-10 0)\lfill
f:0.8 \htext(-7 6){\tiny $0$}

\move(0 0)\lvec(10 0)\lvec(10 5)\lvec(0 5)\lvec(0 0)\lfill f:0.8
\htext(3 6){\tiny $0$}

\move(10 0)\lvec(20 0)\lvec(20 5)\lvec(10 5)\lvec(10 0)\lfill
f:0.8 \htext(13 6){\tiny $0$}

\move(20 0)\lvec(30 0)\lvec(30 5)\lvec(20 5)\lvec(20 0)\lfill
f:0.8 \htext(23 6){\tiny $0$}

\move(30 0)\lvec(40 0)\lvec(40 5)\lvec(30 5)\lvec(30 0)\lfill
f:0.8 \htext(33 6){\tiny $0$}

\move(40 0)\lvec(50 0)\lvec(50 5)\lvec(40 5)\lvec(40 0)\lfill
f:0.8 \htext(43 6){\tiny $0$}

\move(50 0)\lvec(60 0)\lvec(60 5)\lvec(50 5)\lvec(50 0)\lfill
f:0.8 \htext(53 6){\tiny $0$}

\move(0 10)\lvec(10 10)\lvec(10 20)\lvec(0 20)\lvec(0 10)\htext(3
13){$_1$}

\move(10 10)\lvec(20 10)\lvec(20 20)\lvec(10 20)\lvec(10
10)\htext(13 13){$_1$}

\move(20 10)\lvec(30 10)\lvec(30 20)\lvec(20 20)\lvec(20
10)\htext(23 13){$_1$}

\move(30 10)\lvec(40 10)\lvec(40 20)\lvec(30 20)\lvec(30
10)\htext(33 13){$_1$}

\move(40 10)\lvec(50 10)\lvec(50 20)\lvec(40 20)\lvec(40
10)\htext(43 13){$_1$}

\move(50 10)\lvec(60 10)\lvec(60 20)\lvec(50 20)\lvec(50
10)\htext(53 13){$_1$}
\move(0 20)\lvec(10 20)\lvec(10 30)\lvec(0 30)\lvec(0 20)\htext(3
23){$_2$}

\move(10 20)\lvec(20 20)\lvec(20 30)\lvec(10 30)\lvec(10
20)\htext(13 23){$_2$}

\move(20 20)\lvec(30 20)\lvec(30 30)\lvec(20 30)\lvec(20
20)\htext(23 23){$_2$}

\move(30 20)\lvec(40 20)\lvec(40 30)\lvec(30 30)\lvec(30
20)\htext(33 23){$_2$}

\move(40 20)\lvec(50 20)\lvec(50 30)\lvec(40 30)\lvec(40
20)\htext(43 23){$_2$}

\move(50 20)\lvec(60 20)\lvec(60 30)\lvec(50 30)\lvec(50
20)\htext(53 23){$_2$}

\move(0 30)\lvec(10 30)\lvec(10 40)\lvec(0 40)\lvec(0 30)\htext(3
33){$_1$}

\move(10 30)\lvec(20 30)\lvec(20 40)\lvec(10 40)\lvec(10
30)\htext(13 33){$_1$}

\move(20 30)\lvec(30 30)\lvec(30 40)\lvec(20 40)\lvec(20
30)\htext(23 33){$_1$}

\move(30 30)\lvec(40 30)\lvec(40 40)\lvec(30 40)\lvec(30
30)\htext(33 33){$_1$}

\move(40 30)\lvec(50 30)\lvec(50 40)\lvec(40 40)\lvec(40
30)\htext(43 33){$_1$}

\move(50 30)\lvec(60 30)\lvec(60 40)\lvec(50 40)\lvec(50
30)\htext(53 33){$_1$}
\move(0 40)\lvec(10 40)\lvec(10 45)\lvec(0 45)\lvec(0 40)\htext(3
41){\tiny $0$}

\move(10 40)\lvec(20 40)\lvec(20 45)\lvec(10 45)\lvec(10
40)\htext(13 41){\tiny $0$}

\move(20 40)\lvec(30 40)\lvec(30 45)\lvec(20 45)\lvec(20
40)\htext(23 41){\tiny $0$}

\move(30 40)\lvec(40 40)\lvec(40 45)\lvec(30 45)\lvec(30
40)\htext(33 41){\tiny $0$}

\move(40 40)\lvec(50 40)\lvec(50 45)\lvec(40 45)\lvec(40
40)\htext(43 41){\tiny $0$}

\move(50 40)\lvec(60 40)\lvec(60 45)\lvec(50 45)\lvec(50
40)\htext(53 41){\tiny $0$}

\end{texdraw}}\hskip 5mm
$+$ \hskip 3mm $(1+q^{10})$ \raisebox{-0.5\height}{
\begin{texdraw}
\drawdim em \setunitscale 0.13 \linewd 0.5

\move(-10 0)\lvec(0 0)\lvec(0 10)\lvec(-10 10)\lvec(-10
0)\htext(-7 1){\tiny $0$}

\move(0 0)\lvec(10 0)\lvec(10 10)\lvec(0 10)\lvec(0 0)\htext(3
1){\tiny $0$}

\move(10 0)\lvec(20 0)\lvec(20 10)\lvec(10 10)\lvec(10 0)\htext(13
1){\tiny $0$}

\move(20 0)\lvec(30 0)\lvec(30 10)\lvec(20 10)\lvec(20 0)\htext(23
1){\tiny $0$}

\move(30 0)\lvec(40 0)\lvec(40 10)\lvec(30 10)\lvec(30 0)\htext(33
1){\tiny $0$}

\move(40 0)\lvec(50 0)\lvec(50 10)\lvec(40 10)\lvec(40 0)\htext(43
1){\tiny $0$}

\move(50 0)\lvec(60 0)\lvec(60 10)\lvec(50 10)\lvec(50 0)\htext(53
1){\tiny $0$}
\move(-10 0)\lvec(0 0)\lvec(0 5)\lvec(-10 5)\lvec(-10 0)\lfill
f:0.8 \htext(-7 6){\tiny $0$}

\move(0 0)\lvec(10 0)\lvec(10 5)\lvec(0 5)\lvec(0 0)\lfill f:0.8
\htext(3 6){\tiny $0$}

\move(10 0)\lvec(20 0)\lvec(20 5)\lvec(10 5)\lvec(10 0)\lfill
f:0.8 \htext(13 6){\tiny $0$}

\move(20 0)\lvec(30 0)\lvec(30 5)\lvec(20 5)\lvec(20 0)\lfill
f:0.8 \htext(23 6){\tiny $0$}

\move(30 0)\lvec(40 0)\lvec(40 5)\lvec(30 5)\lvec(30 0)\lfill
f:0.8 \htext(33 6){\tiny $0$}

\move(40 0)\lvec(50 0)\lvec(50 5)\lvec(40 5)\lvec(40 0)\lfill
f:0.8 \htext(43 6){\tiny $0$}

\move(50 0)\lvec(60 0)\lvec(60 5)\lvec(50 5)\lvec(50 0)\lfill
f:0.8 \htext(53 6){\tiny $0$}

\move(0 10)\lvec(10 10)\lvec(10 20)\lvec(0 20)\lvec(0 10)\htext(3
13){$_1$}

\move(10 10)\lvec(20 10)\lvec(20 20)\lvec(10 20)\lvec(10
10)\htext(13 13){$_1$}

\move(20 10)\lvec(30 10)\lvec(30 20)\lvec(20 20)\lvec(20
10)\htext(23 13){$_1$}

\move(30 10)\lvec(40 10)\lvec(40 20)\lvec(30 20)\lvec(30
10)\htext(33 13){$_1$}

\move(40 10)\lvec(50 10)\lvec(50 20)\lvec(40 20)\lvec(40
10)\htext(43 13){$_1$}

\move(50 10)\lvec(60 10)\lvec(60 20)\lvec(50 20)\lvec(50
10)\htext(53 13){$_1$}
\move(0 20)\lvec(10 20)\lvec(10 30)\lvec(0 30)\lvec(0 20)\htext(3
23){$_2$}

\move(10 20)\lvec(20 20)\lvec(20 30)\lvec(10 30)\lvec(10
20)\htext(13 23){$_2$}

\move(20 20)\lvec(30 20)\lvec(30 30)\lvec(20 30)\lvec(20
20)\htext(23 23){$_2$}

\move(30 20)\lvec(40 20)\lvec(40 30)\lvec(30 30)\lvec(30
20)\htext(33 23){$_2$}

\move(40 20)\lvec(50 20)\lvec(50 30)\lvec(40 30)\lvec(40
20)\htext(43 23){$_2$}

\move(50 20)\lvec(60 20)\lvec(60 30)\lvec(50 30)\lvec(50
20)\htext(53 23){$_2$}

\move(0 30)\lvec(10 30)\lvec(10 40)\lvec(0 40)\lvec(0 30)\htext(3
33){$_1$}

\move(10 30)\lvec(20 30)\lvec(20 40)\lvec(10 40)\lvec(10
30)\htext(13 33){$_1$}

\move(20 30)\lvec(30 30)\lvec(30 40)\lvec(20 40)\lvec(20
30)\htext(23 33){$_1$}

\move(30 30)\lvec(40 30)\lvec(40 40)\lvec(30 40)\lvec(30
30)\htext(33 33){$_1$}

\move(40 30)\lvec(50 30)\lvec(50 40)\lvec(40 40)\lvec(40
30)\htext(43 33){$_1$}

\move(50 30)\lvec(60 30)\lvec(60 40)\lvec(50 40)\lvec(50
30)\htext(53 33){$_1$}
\move(10 40)\lvec(20 40)\lvec(20 45)\lvec(10 45)\lvec(10
40)\htext(13 41){\tiny $0$}

\move(20 40)\lvec(30 40)\lvec(30 45)\lvec(20 45)\lvec(20
40)\htext(23 41){\tiny $0$}

\move(30 40)\lvec(40 40)\lvec(40 45)\lvec(30 45)\lvec(30
40)\htext(33 41){\tiny $0$}

\move(40 40)\lvec(50 40)\lvec(50 45)\lvec(40 45)\lvec(40
40)\htext(43 41){\tiny $0$}

\move(50 40)\lvec(60 40)\lvec(60 45)\lvec(50 45)\lvec(50
40)\htext(53 41){\tiny $0$}

\move(50 45)\lvec(60 45)\lvec(60 50)\lvec(50 50)\lvec(50
45)\htext(53 46){\tiny $0$}

\end{texdraw}}\quad .\vskip 5mm

}\end{ex}

\vskip 2mm \noindent {\bf Case 3.} Suppose that the $i$-blocks are
of type III.\vskip 2mm

If $b$ is a removable $i$-block in $y_k$ of $Y$, then we define
$Y\nearrow b$ to be the Young wall obtained by removing the block
$b$ from $Y$. We also consider the following $i$-block $b$ in
$y_k$ of $Y$, which we call a {\it virtually removable
$i$-block}.\vskip 5mm

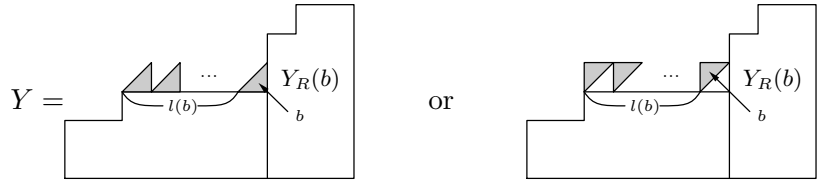
\begin{center}
$Y=$\raisebox{-0.4\height}{\begin{texdraw}\fontsize{9}{9} \drawdim
em \setunitscale 0.1 \linewd 0.5

\move(-10 -20)\lvec(-10 0)\lvec(10 0)\lvec(10 10)\lvec(60
10)\lvec(60 30)\lvec(70 30)\lvec(70 40)\lvec(90 40)\lvec(90
-20)\lvec(-10 -20)

\move(60 10)\lvec(60 -20)

\move(10 10)\lvec(20 20)\lvec(20 10)\lvec(10 10)\lfill f:0.8

\move(20 10)\lvec(30 20)\lvec(30 10)\lvec(20 10)\lfill f:0.8

\move(50 10)\lvec(60 20)\lvec(60 10)\lvec(50 10)\lfill f:0.8

\htext(65 10){$Y_R(b)$}

\htext(37 14){$_{\cdots}$}

\move(67 3) \arrowheadtype t:F \arrowheadsize l:4 w:2 \avec(57 13)

\htext(70 0){$_b$}

\move(24 5)\clvec(16 5)(13 5)(10 10)

\move(36 5)\clvec(44 5)(47 5)(50 10)

\htext(26 2){$_{l(b)}$}
\end{texdraw}}\hskip 1cm or
\hskip .7cm\raisebox{-0.4\height}{
\begin{texdraw}\fontsize{9}{9}
\drawdim em \setunitscale 0.1 \linewd 0.5

\move(-10 -20)\lvec(-10 0)\lvec(10 0)\lvec(10 10)\lvec(60
10)\lvec(60 30)\lvec(70 30)\lvec(70 40)\lvec(90 40)\lvec(90
-20)\lvec(-10 -20)

\move(60 10)\lvec(60 -20)

\htext(65 10){$Y_R(b)$}

\move(10 10)\lvec(20 20)\lvec(10 20)\lvec(10 10)\lfill f:0.8

\move(20 10)\lvec(30 20)\lvec(20 20)\lvec(20 10)\lfill f:0.8

\move(50 10)\lvec(60 20)\lvec(50 20)\lvec(50 10)\lfill f:0.8


\htext(37 14){$_{\cdots}$}

\move(67 3) \arrowheadtype t:F \arrowheadsize l:4 w:2 \avec(53 17)

\htext(70 0){$_b$}

\move(24 5)\clvec(16 5)(13 5)(10 10)

\move(36 5)\clvec(44 5)(47 5)(50 10)

\htext(26 2){$_{l(b)}$}
\end{texdraw}}
\end{center}\vskip 5mm

\noindent In this case, we define $Y\nearrow b$ to be \vskip 5mm

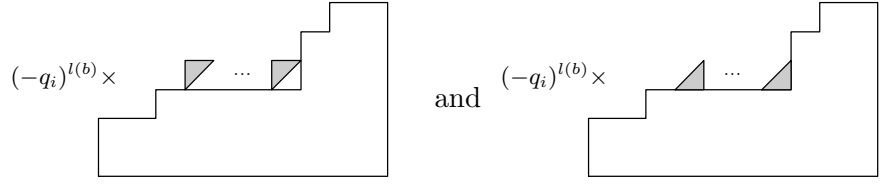
\begin{center}
\raisebox{-0.4\height}{
\begin{texdraw}\fontsize{9}{9}
\drawdim em \setunitscale 0.1 \linewd 0.5

\htext(-40 10){${(-q_i)^{l(b)}\times}$}

\move(-10 -20)\lvec(-10 0)\lvec(10 0)\lvec(10 10)\lvec(60
10)\lvec(60 30)\lvec(70 30)\lvec(70 40)\lvec(90 40)\lvec(90
-20)\lvec(-10 -20)

\move(20 10)\lvec(30 20)\lvec(20 20)\lvec(20 10)\lfill f:0.8

\move(50 10)\lvec(60 20)\lvec(50 20)\lvec(50 10)\lfill f:0.8


\htext(37 14){$_{\cdots}$}
\end{texdraw}} \hskip 5mm and
\raisebox{-0.4\height}{
\begin{texdraw}\fontsize{9}{9}
\drawdim em \setunitscale 0.1 \linewd 0.5

\htext(-40 10){${(-q_i)^{l(b)}\times}$}

\move(-10 -20)\lvec(-10 0)\lvec(10 0)\lvec(10 10)\lvec(60
10)\lvec(60 30)\lvec(70 30)\lvec(70 40)\lvec(90 40)\lvec(90
-20)\lvec(-10 -20)

\move(20 10)\lvec(30 20)\lvec(30 10)\lvec(20 10)\lfill f:0.8

\move(50 10)\lvec(60 20)\lvec(60 10)\lvec(50 10)\lfill f:0.8


\htext(37 14){$_{\cdots}$}

\end{texdraw}}
\end{center}\vskip 5mm

\noindent respectively, where $l(b)\geq 1$ is given in the above
figure. Note that, unlike Case 1 and Case 2, we need to {\it
shift} the blocks from left to right and from back to front (resp.
from front to back). In either case, we define
$Y_R(b)=(y_{k-1},\cdots,y_{0})$, and set
\begin{equation}
R_i(b;Y)=\varphi_i(Y_R(b))-\varepsilon_i(Y_R(b)).
\end{equation}
If $k=0$, we understand $Y_R(b)=\emptyset$ and $R_i(b;Y)=0$. Then
we define
\begin{equation}
e_i\,Y=\sum_{b}q_i^{-R_i(b;Y)}(Y\nearrow b),
\end{equation}
\vskip -2mm \noindent where $b$ runs over all removable and
virtually removable $i$-blocks in $Y$.

On the other hand, if $b$ is an admissible $i$-slot in $y_k$ of
$Y$, then we define $Y\swarrow b$ to be the Young wall obtained by
adding an $i$-block at $b$. We also consider the following
$i$-slot $b$ in $y_k$ of $Y$, which we call a {\it virtually
admissible $i$-slot}:\vskip 5mm

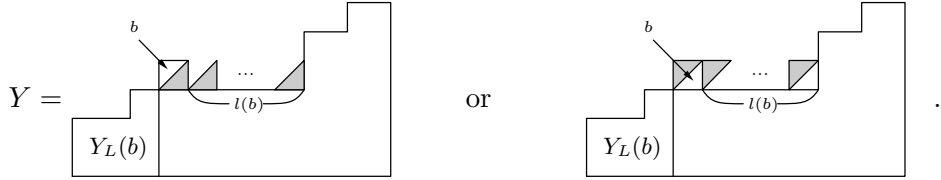
\begin{center}
$Y=$\raisebox{-0.4\height}{
\begin{texdraw}\fontsize{9}{9}
\drawdim em \setunitscale 0.1 \linewd 0.5

\move(-20 -20)\lvec(-20 0)\lvec(0 0)\lvec(0 10)\lvec(60
10)\lvec(60 30)\lvec(75 30)\lvec(75 40)\lvec(90 40)\lvec(90
-20)\lvec(-20 -20)

\move(10 10)\lvec(10 -20)\htext(-14 -14){$Y_L(b)$}

\move(10 10)\lvec(10 20)\lvec(20 20)

\move(10 10)\lvec(20 20)\lvec(20 10)\lvec(10 10)\lfill f:0.8

\move(20 10)\lvec(30 20)\lvec(30 10)\lvec(20 10)\lfill f:0.8

\move(50 10)\lvec(60 20)\lvec(60 10)\lvec(50 10)\lfill f:0.8


\htext(37 14){$_{\cdots}$}

\move(3 27) \arrowheadtype t:F \arrowheadsize l:4 w:2 \avec(13 17)

\htext(0 30){$_b$}

\move(34 5)\clvec(26 5)(23 5)(20 10)

\move(46 5)\clvec(54 5)(57 5)(60 10)

\htext(36 2){$_{l(b)}$}
\end{texdraw}}\hskip 1cm or
\hskip 1cm\raisebox{-0.4\height}{
\begin{texdraw}\fontsize{9}{9}
\drawdim em \setunitscale 0.1 \linewd 0.5

\move(-20 -20)\lvec(-20 0)\lvec(0 0)\lvec(0 10)\lvec(60
10)\lvec(60 30)\lvec(75 30)\lvec(75 40)\lvec(90 40)\lvec(90
-20)\lvec(-20 -20)

\move(10 10)\lvec(10 -20)\htext(-14 -14){$Y_L(b)$}

\move(10 10)\lvec(20 20)\lvec(10 20)\lvec(10 10)\lfill f:0.8

\move(20 10)\lvec(30 20)\lvec(20 20)\lvec(20 10)\lfill f:0.8

\move(50 10)\lvec(60 20)\lvec(50 20)\lvec(50 10)\lfill f:0.8

\htext(37 14){$_{\cdots}$}

\move(3 27) \arrowheadtype t:F \arrowheadsize l:4 w:2 \avec(17 13)

\htext(0 30){$_b$}

\move(34 5)\clvec(26 5)(23 5)(20 10)

\move(46 5)\clvec(54 5)(57 5)(60 10)

\htext(36 2){$_{l(b)}$}
\end{texdraw}}\quad .
\end{center}\vskip 5mm

\noindent In this case, we define $Y\swarrow b$ to be \vskip 5mm

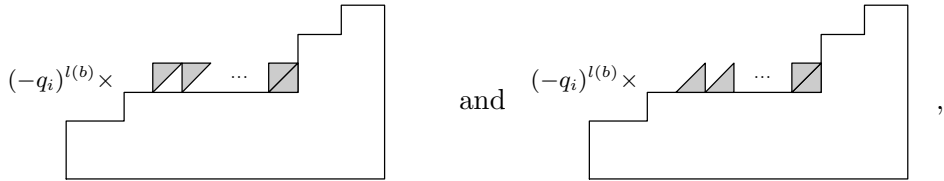
\begin{center}
\raisebox{-0.4\height}{
\begin{texdraw}\fontsize{9}{9}
\drawdim em \setunitscale 0.1 \linewd 0.5

\htext(-40 10){${(-q_i)^{l(b)}\times}$}

\move(-20 -20)\lvec(-20 0)\lvec(0 0)\lvec(0 10)\lvec(60
10)\lvec(60 30)\lvec(75 30)\lvec(75 40)\lvec(90 40)\lvec(90
-20)\lvec(-20 -20)

\move(10 10)\lvec(20 20)\lvec(10 20)\lvec(10 10)\lfill f:0.8

\move(20 10)\lvec(30 20)\lvec(20 20)\lvec(20 10)\lfill f:0.8

\move(50 10)\lvec(60 20)\lvec(50 20)\lvec(50 10)\lfill f:0.8

\move(50 10)\lvec(60 20)\lvec(60 10)\lvec(50 10)\lfill f:0.8


\htext(37 14){$_{\cdots}$}

\end{texdraw}}\hskip 10mm and
\hskip 2mm \raisebox{-0.4\height}{\begin{texdraw}\fontsize{9}{9}
\drawdim em \setunitscale 0.1 \linewd 0.5

\htext(-40 10){${(-q_i)^{l(b)}\times}$}

\move(-20 -20)\lvec(-20 0)\lvec(0 0)\lvec(0 10)\lvec(60
10)\lvec(60 30)\lvec(75 30)\lvec(75 40)\lvec(90 40)\lvec(90
-20)\lvec(-20 -20)

\move(10 10)\lvec(20 20)\lvec(20 10)\lvec(10 10)\lfill f:0.8

\move(20 10)\lvec(30 20)\lvec(30 10)\lvec(20 10)\lfill f:0.8

\move(50 10)\lvec(60 20)\lvec(60 10)\lvec(50 10)\lfill f:0.8

\move(50 10)\lvec(60 20)\lvec(50 20)\lvec(50 10)\lfill f:0.8


\htext(37 14){$_{\cdots}$}

\end{texdraw}}\quad ,
\end{center}\vskip 5mm

\noindent respectively, where $l(b)\geq 1$ is given in the above
figure. Here again, one can observe that we need to {\it shift}
the blocks from left to right and from back to front (resp. from
front to back). In either case, we define
$Y_L(b)=(\cdots,y_{k+2},y_{k+1})$ and set
\begin{equation}
L_i(b;Y)=\varphi_i(Y_L(b))-\varepsilon_i(Y_L(b)).
\end{equation}
Then we define
\begin{equation}
f_i\,Y=\sum_{b}q_i^{L_i(b;Y)}(Y\swarrow b),
\end{equation}
\noindent where $b$ runs over all admissible and virtually
admissible $i$-slots in $Y$.
\begin{ex}{\rm
If $\frak{g}=B^{(1)}_3$, $\Lambda=\Lambda_0$ and $i=0$, then
$q_0=q^2$ and we have \vskip 5mm

$e_0$\raisebox{-0.5\height}{
\begin{texdraw}
\drawdim em \setunitscale 0.13 \linewd 0.5


\move(10 0)\lvec(20 0)\lvec(20 10)\lvec(10 10)\lvec(10 0)\htext(12
6){\tiny $1$}

\move(20 0)\lvec(30 0)\lvec(30 10)\lvec(20 10)\lvec(20 0)\htext(22
6){\tiny $0$}

\move(30 0)\lvec(40 0)\lvec(40 10)\lvec(30 10)\lvec(30 0)\htext(32
6){\tiny $1$}

\move(40 0)\lvec(50 0)\lvec(50 10)\lvec(40 10)\lvec(40 0)\htext(42
6){\tiny $0$}
%

\move(10 0)\lvec(20 10)\lvec(20 0)\lvec(10 0)\lfill f:0.8
\htext(16 2){\tiny $0$}

\move(20 0)\lvec(30 10)\lvec(30 0)\lvec(20 0)\lfill f:0.8
\htext(26 2){\tiny $1$}

\move(30 0)\lvec(40 10)\lvec(40 0)\lvec(30 0)\lfill f:0.8
\htext(36 2){\tiny $0$}

\move(40 0)\lvec(50 10)\lvec(50 0)\lvec(40 0)\lfill f:0.8
\htext(46 2){\tiny $1$}
\move(10 10)\lvec(20 10)\lvec(20 20)\lvec(10 20)\lvec(10
10)\htext(13 13){$_2$}

\move(20 10)\lvec(30 10)\lvec(30 20)\lvec(20 20)\lvec(20
10)\htext(23 13){$_2$}

\move(30 10)\lvec(40 10)\lvec(40 20)\lvec(30 20)\lvec(30
10)\htext(33 13){$_2$}

\move(40 10)\lvec(50 10)\lvec(50 20)\lvec(40 20)\lvec(40
10)\htext(43 13){$_2$}
\move(10 20)\lvec(20 20)\lvec(20 30)\lvec(10 30)\lvec(10
20)\htext(13 26){\tiny $3$}

\move(20 20)\lvec(30 20)\lvec(30 30)\lvec(20 30)\lvec(20
20)\htext(23 26){\tiny $3$}

\move(30 20)\lvec(40 20)\lvec(40 30)\lvec(30 30)\lvec(30
20)\htext(33 26){\tiny $3$}

\move(40 20)\lvec(50 20)\lvec(50 30)\lvec(40 30)\lvec(40
20)\htext(43 26){\tiny $3$}

\move(10 20)\lvec(20 20)\lvec(20 25)\lvec(10 25)\lvec(10 20)
\htext(13 21){\tiny $3$}

\move(20 20)\lvec(30 20)\lvec(30 25)\lvec(20 25)\lvec(20 20)
\htext(23 21){\tiny $3$}

\move(30 20)\lvec(40 20)\lvec(40 25)\lvec(30 25)\lvec(30 20)
\htext(33 21){\tiny $3$}

\move(40 20)\lvec(50 20)\lvec(50 25)\lvec(40 25)\lvec(40 20)
\htext(43 21){\tiny $3$}
\move(10 30)\lvec(20 30)\lvec(20 40)\lvec(10 40)\lvec(10
30)\htext(13 33){$_2$}

\move(20 30)\lvec(30 30)\lvec(30 40)\lvec(20 40)\lvec(20
30)\htext(23 33){$_2$}

\move(30 30)\lvec(40 30)\lvec(40 40)\lvec(30 40)\lvec(30
30)\htext(33 33){$_2$}

\move(40 30)\lvec(50 30)\lvec(50 40)\lvec(40 40)\lvec(40
30)\htext(43 33){$_2$}
\move(10 40)\lvec(20 50)\lvec(20 40)\lvec(10 40) \htext(16
42){\tiny $0$}

\move(20 40)\lvec(30 50)\lvec(30 40)\lvec(20 40) \htext(26
42){\tiny $1$}

\move(30 40)\lvec(40 50)\lvec(40 40)\lvec(30 40) \htext(36
42){\tiny $0$}

\move(40 40)\lvec(50 50)\lvec(50 40)\lvec(40 40) \htext(46
42){\tiny $1$}

\move(40 40)\lvec(50 50)\lvec(40 50)\lvec(40 40) \htext(42
46){\tiny $0$}
\end{texdraw}}\vskip 5mm

$=$ \hskip 3mm $q^2$ \raisebox{-0.5\height}{
\begin{texdraw}
\drawdim em \setunitscale 0.13 \linewd 0.5


\move(10 0)\lvec(20 0)\lvec(20 10)\lvec(10 10)\lvec(10 0)\htext(12
6){\tiny $1$}

\move(20 0)\lvec(30 0)\lvec(30 10)\lvec(20 10)\lvec(20 0)\htext(22
6){\tiny $0$}

\move(30 0)\lvec(40 0)\lvec(40 10)\lvec(30 10)\lvec(30 0)\htext(32
6){\tiny $1$}

\move(40 0)\lvec(50 0)\lvec(50 10)\lvec(40 10)\lvec(40 0)\htext(42
6){\tiny $0$}
%

\move(10 0)\lvec(20 10)\lvec(20 0)\lvec(10 0)\lfill f:0.8
\htext(16 2){\tiny $0$}

\move(20 0)\lvec(30 10)\lvec(30 0)\lvec(20 0)\lfill f:0.8
\htext(26 2){\tiny $1$}

\move(30 0)\lvec(40 10)\lvec(40 0)\lvec(30 0)\lfill f:0.8
\htext(36 2){\tiny $0$}

\move(40 0)\lvec(50 10)\lvec(50 0)\lvec(40 0)\lfill f:0.8
\htext(46 2){\tiny $1$}
\move(10 10)\lvec(20 10)\lvec(20 20)\lvec(10 20)\lvec(10
10)\htext(13 13){$_2$}

\move(20 10)\lvec(30 10)\lvec(30 20)\lvec(20 20)\lvec(20
10)\htext(23 13){$_2$}

\move(30 10)\lvec(40 10)\lvec(40 20)\lvec(30 20)\lvec(30
10)\htext(33 13){$_2$}

\move(40 10)\lvec(50 10)\lvec(50 20)\lvec(40 20)\lvec(40
10)\htext(43 13){$_2$}
\move(10 20)\lvec(20 20)\lvec(20 30)\lvec(10 30)\lvec(10
20)\htext(13 26){\tiny $3$}

\move(20 20)\lvec(30 20)\lvec(30 30)\lvec(20 30)\lvec(20
20)\htext(23 26){\tiny $3$}

\move(30 20)\lvec(40 20)\lvec(40 30)\lvec(30 30)\lvec(30
20)\htext(33 26){\tiny $3$}

\move(40 20)\lvec(50 20)\lvec(50 30)\lvec(40 30)\lvec(40
20)\htext(43 26){\tiny $3$}

\move(10 20)\lvec(20 20)\lvec(20 25)\lvec(10 25)\lvec(10 20)
\htext(13 21){\tiny $3$}

\move(20 20)\lvec(30 20)\lvec(30 25)\lvec(20 25)\lvec(20 20)
\htext(23 21){\tiny $3$}

\move(30 20)\lvec(40 20)\lvec(40 25)\lvec(30 25)\lvec(30 20)
\htext(33 21){\tiny $3$}

\move(40 20)\lvec(50 20)\lvec(50 25)\lvec(40 25)\lvec(40 20)
\htext(43 21){\tiny $3$}
\move(10 30)\lvec(20 30)\lvec(20 40)\lvec(10 40)\lvec(10
30)\htext(13 33){$_2$}

\move(20 30)\lvec(30 30)\lvec(30 40)\lvec(20 40)\lvec(20
30)\htext(23 33){$_2$}

\move(30 30)\lvec(40 30)\lvec(40 40)\lvec(30 40)\lvec(30
30)\htext(33 33){$_2$}

\move(40 30)\lvec(50 30)\lvec(50 40)\lvec(40 40)\lvec(40
30)\htext(43 33){$_2$}
%

\move(20 40)\lvec(30 50)\lvec(30 40)\lvec(20 40) \htext(26
42){\tiny $1$}

\move(30 40)\lvec(40 50)\lvec(40 40)\lvec(30 40) \htext(36
42){\tiny $0$}

\move(40 40)\lvec(50 50)\lvec(50 40)\lvec(40 40) \htext(46
42){\tiny $1$}

\move(40 40)\lvec(50 50)\lvec(40 50)\lvec(40 40) \htext(42
46){\tiny $0$}
\end{texdraw}}\hskip 5mm $+$ \hskip 5mm  $q^6$  \raisebox{-0.5\height}{
\begin{texdraw}
\drawdim em \setunitscale 0.13 \linewd 0.5


\move(10 0)\lvec(20 0)\lvec(20 10)\lvec(10 10)\lvec(10 0)\htext(12
6){\tiny $1$}

\move(20 0)\lvec(30 0)\lvec(30 10)\lvec(20 10)\lvec(20 0)\htext(22
6){\tiny $0$}

\move(30 0)\lvec(40 0)\lvec(40 10)\lvec(30 10)\lvec(30 0)\htext(32
6){\tiny $1$}

\move(40 0)\lvec(50 0)\lvec(50 10)\lvec(40 10)\lvec(40 0)\htext(42
6){\tiny $0$}
%

\move(10 0)\lvec(20 10)\lvec(20 0)\lvec(10 0)\lfill f:0.8
\htext(16 2){\tiny $0$}

\move(20 0)\lvec(30 10)\lvec(30 0)\lvec(20 0)\lfill f:0.8
\htext(26 2){\tiny $1$}

\move(30 0)\lvec(40 10)\lvec(40 0)\lvec(30 0)\lfill f:0.8
\htext(36 2){\tiny $0$}

\move(40 0)\lvec(50 10)\lvec(50 0)\lvec(40 0)\lfill f:0.8
\htext(46 2){\tiny $1$}
\move(10 10)\lvec(20 10)\lvec(20 20)\lvec(10 20)\lvec(10
10)\htext(13 13){$_2$}

\move(20 10)\lvec(30 10)\lvec(30 20)\lvec(20 20)\lvec(20
10)\htext(23 13){$_2$}

\move(30 10)\lvec(40 10)\lvec(40 20)\lvec(30 20)\lvec(30
10)\htext(33 13){$_2$}

\move(40 10)\lvec(50 10)\lvec(50 20)\lvec(40 20)\lvec(40
10)\htext(43 13){$_2$}
\move(10 20)\lvec(20 20)\lvec(20 30)\lvec(10 30)\lvec(10
20)\htext(13 26){\tiny $3$}

\move(20 20)\lvec(30 20)\lvec(30 30)\lvec(20 30)\lvec(20
20)\htext(23 26){\tiny $3$}

\move(30 20)\lvec(40 20)\lvec(40 30)\lvec(30 30)\lvec(30
20)\htext(33 26){\tiny $3$}

\move(40 20)\lvec(50 20)\lvec(50 30)\lvec(40 30)\lvec(40
20)\htext(43 26){\tiny $3$}

\move(10 20)\lvec(20 20)\lvec(20 25)\lvec(10 25)\lvec(10 20)
\htext(13 21){\tiny $3$}

\move(20 20)\lvec(30 20)\lvec(30 25)\lvec(20 25)\lvec(20 20)
\htext(23 21){\tiny $3$}

\move(30 20)\lvec(40 20)\lvec(40 25)\lvec(30 25)\lvec(30 20)
\htext(33 21){\tiny $3$}

\move(40 20)\lvec(50 20)\lvec(50 25)\lvec(40 25)\lvec(40 20)
\htext(43 21){\tiny $3$}
\move(10 30)\lvec(20 30)\lvec(20 40)\lvec(10 40)\lvec(10
30)\htext(13 33){$_2$}

\move(20 30)\lvec(30 30)\lvec(30 40)\lvec(20 40)\lvec(20
30)\htext(23 33){$_2$}

\move(30 30)\lvec(40 30)\lvec(40 40)\lvec(30 40)\lvec(30
30)\htext(33 33){$_2$}

\move(40 30)\lvec(50 30)\lvec(50 40)\lvec(40 40)\lvec(40
30)\htext(43 33){$_2$}
%

\move(20 40)\lvec(30 50)\lvec(20 50)\lvec(20 40) \htext(22
46){\tiny $0$}

\move(30 40)\lvec(40 50)\lvec(30 50)\lvec(30 40) \htext(32
46){\tiny $1$}

\move(40 40)\lvec(50 50)\lvec(40 50)\lvec(40 40) \htext(42
46){\tiny $0$}

\move(40 40)\lvec(50 50)\lvec(50 40)\lvec(40 40) \htext(46
42){\tiny $1$}
\end{texdraw}}
\hskip 5mm $+$ \hskip 5mm \raisebox{-0.5\height}{
\begin{texdraw}
\drawdim em \setunitscale 0.13 \linewd 0.5


\move(10 0)\lvec(20 0)\lvec(20 10)\lvec(10 10)\lvec(10 0)\htext(12
6){\tiny $1$}

\move(20 0)\lvec(30 0)\lvec(30 10)\lvec(20 10)\lvec(20 0)\htext(22
6){\tiny $0$}

\move(30 0)\lvec(40 0)\lvec(40 10)\lvec(30 10)\lvec(30 0)\htext(32
6){\tiny $1$}

\move(40 0)\lvec(50 0)\lvec(50 10)\lvec(40 10)\lvec(40 0)\htext(42
6){\tiny $0$}
%

\move(10 0)\lvec(20 10)\lvec(20 0)\lvec(10 0)\lfill f:0.8
\htext(16 2){\tiny $0$}

\move(20 0)\lvec(30 10)\lvec(30 0)\lvec(20 0)\lfill f:0.8
\htext(26 2){\tiny $1$}

\move(30 0)\lvec(40 10)\lvec(40 0)\lvec(30 0)\lfill f:0.8
\htext(36 2){\tiny $0$}

\move(40 0)\lvec(50 10)\lvec(50 0)\lvec(40 0)\lfill f:0.8
\htext(46 2){\tiny $1$}
\move(10 10)\lvec(20 10)\lvec(20 20)\lvec(10 20)\lvec(10
10)\htext(13 13){$_2$}

\move(20 10)\lvec(30 10)\lvec(30 20)\lvec(20 20)\lvec(20
10)\htext(23 13){$_2$}

\move(30 10)\lvec(40 10)\lvec(40 20)\lvec(30 20)\lvec(30
10)\htext(33 13){$_2$}

\move(40 10)\lvec(50 10)\lvec(50 20)\lvec(40 20)\lvec(40
10)\htext(43 13){$_2$}
\move(10 20)\lvec(20 20)\lvec(20 30)\lvec(10 30)\lvec(10
20)\htext(13 26){\tiny $3$}

\move(20 20)\lvec(30 20)\lvec(30 30)\lvec(20 30)\lvec(20
20)\htext(23 26){\tiny $3$}

\move(30 20)\lvec(40 20)\lvec(40 30)\lvec(30 30)\lvec(30
20)\htext(33 26){\tiny $3$}

\move(40 20)\lvec(50 20)\lvec(50 30)\lvec(40 30)\lvec(40
20)\htext(43 26){\tiny $3$}

\move(10 20)\lvec(20 20)\lvec(20 25)\lvec(10 25)\lvec(10 20)
\htext(13 21){\tiny $3$}

\move(20 20)\lvec(30 20)\lvec(30 25)\lvec(20 25)\lvec(20 20)
\htext(23 21){\tiny $3$}

\move(30 20)\lvec(40 20)\lvec(40 25)\lvec(30 25)\lvec(30 20)
\htext(33 21){\tiny $3$}

\move(40 20)\lvec(50 20)\lvec(50 25)\lvec(40 25)\lvec(40 20)
\htext(43 21){\tiny $3$}
\move(10 30)\lvec(20 30)\lvec(20 40)\lvec(10 40)\lvec(10
30)\htext(13 33){$_2$}

\move(20 30)\lvec(30 30)\lvec(30 40)\lvec(20 40)\lvec(20
30)\htext(23 33){$_2$}

\move(30 30)\lvec(40 30)\lvec(40 40)\lvec(30 40)\lvec(30
30)\htext(33 33){$_2$}

\move(40 30)\lvec(50 30)\lvec(50 40)\lvec(40 40)\lvec(40
30)\htext(43 33){$_2$}
\move(10 40)\lvec(20 50)\lvec(20 40)\lvec(10 40) \htext(16
42){\tiny $0$}

\move(20 40)\lvec(30 50)\lvec(30 40)\lvec(20 40) \htext(26
42){\tiny $1$}

\move(30 40)\lvec(40 50)\lvec(40 40)\lvec(30 40) \htext(36
42){\tiny $0$}

\move(40 40)\lvec(50 50)\lvec(50 40)\lvec(40 40) \htext(46
42){\tiny $1$}

\end{texdraw}}\quad ,

\vskip 1cm $f_0$\raisebox{-0.5\height}{
\begin{texdraw}
\drawdim em \setunitscale 0.13 \linewd 0.5


\move(10 0)\lvec(20 0)\lvec(20 10)\lvec(10 10)\lvec(10 0)\htext(12
6){\tiny $1$}

\move(20 0)\lvec(30 0)\lvec(30 10)\lvec(20 10)\lvec(20 0)\htext(22
6){\tiny $0$}

\move(30 0)\lvec(40 0)\lvec(40 10)\lvec(30 10)\lvec(30 0)\htext(32
6){\tiny $1$}

\move(40 0)\lvec(50 0)\lvec(50 10)\lvec(40 10)\lvec(40 0)\htext(42
6){\tiny $0$}
%

\move(10 0)\lvec(20 10)\lvec(20 0)\lvec(10 0)\lfill f:0.8
\htext(16 2){\tiny $0$}

\move(20 0)\lvec(30 10)\lvec(30 0)\lvec(20 0)\lfill f:0.8
\htext(26 2){\tiny $1$}

\move(30 0)\lvec(40 10)\lvec(40 0)\lvec(30 0)\lfill f:0.8
\htext(36 2){\tiny $0$}

\move(40 0)\lvec(50 10)\lvec(50 0)\lvec(40 0)\lfill f:0.8
\htext(46 2){\tiny $1$}
\move(10 10)\lvec(20 10)\lvec(20 20)\lvec(10 20)\lvec(10
10)\htext(13 13){$_2$}

\move(20 10)\lvec(30 10)\lvec(30 20)\lvec(20 20)\lvec(20
10)\htext(23 13){$_2$}

\move(30 10)\lvec(40 10)\lvec(40 20)\lvec(30 20)\lvec(30
10)\htext(33 13){$_2$}

\move(40 10)\lvec(50 10)\lvec(50 20)\lvec(40 20)\lvec(40
10)\htext(43 13){$_2$}
\move(10 20)\lvec(20 20)\lvec(20 30)\lvec(10 30)\lvec(10
20)\htext(13 26){\tiny $3$}

\move(20 20)\lvec(30 20)\lvec(30 30)\lvec(20 30)\lvec(20
20)\htext(23 26){\tiny $3$}

\move(30 20)\lvec(40 20)\lvec(40 30)\lvec(30 30)\lvec(30
20)\htext(33 26){\tiny $3$}

\move(40 20)\lvec(50 20)\lvec(50 30)\lvec(40 30)\lvec(40
20)\htext(43 26){\tiny $3$}

\move(10 20)\lvec(20 20)\lvec(20 25)\lvec(10 25)\lvec(10 20)
\htext(13 21){\tiny $3$}

\move(20 20)\lvec(30 20)\lvec(30 25)\lvec(20 25)\lvec(20 20)
\htext(23 21){\tiny $3$}

\move(30 20)\lvec(40 20)\lvec(40 25)\lvec(30 25)\lvec(30 20)
\htext(33 21){\tiny $3$}

\move(40 20)\lvec(50 20)\lvec(50 25)\lvec(40 25)\lvec(40 20)
\htext(43 21){\tiny $3$}
\move(10 30)\lvec(20 30)\lvec(20 40)\lvec(10 40)\lvec(10
30)\htext(13 33){$_2$}

\move(20 30)\lvec(30 30)\lvec(30 40)\lvec(20 40)\lvec(20
30)\htext(23 33){$_2$}

\move(30 30)\lvec(40 30)\lvec(40 40)\lvec(30 40)\lvec(30
30)\htext(33 33){$_2$}

\move(40 30)\lvec(50 30)\lvec(50 40)\lvec(40 40)\lvec(40
30)\htext(43 33){$_2$}
%

\move(20 40)\lvec(30 50)\lvec(30 40)\lvec(20 40) \htext(26
42){\tiny $1$}

\move(30 40)\lvec(40 50)\lvec(40 40)\lvec(30 40) \htext(36
42){\tiny $0$}

\move(40 40)\lvec(50 50)\lvec(50 40)\lvec(40 40) \htext(46
42){\tiny $1$}

\end{texdraw}}\hskip 3mm $=$ \hskip 3mm
\raisebox{-0.5\height}{
\begin{texdraw}
\drawdim em \setunitscale 0.13 \linewd 0.5

\move(0 0)\lvec(10 0)\lvec(10 10)\lvec(0 10)\lvec(0 0)\htext(2
6){\tiny $0$}

\move(10 0)\lvec(20 0)\lvec(20 10)\lvec(10 10)\lvec(10 0)\htext(12
6){\tiny $1$}

\move(20 0)\lvec(30 0)\lvec(30 10)\lvec(20 10)\lvec(20 0)\htext(22
6){\tiny $0$}

\move(30 0)\lvec(40 0)\lvec(40 10)\lvec(30 10)\lvec(30 0)\htext(32
6){\tiny $1$}

\move(40 0)\lvec(50 0)\lvec(50 10)\lvec(40 10)\lvec(40 0)\htext(42
6){\tiny $0$}
\move(0 0)\lvec(10 10)\lvec(10 0)\lvec(0 0)\lfill f:0.8 \htext(6
2){\tiny $1$}

\move(10 0)\lvec(20 10)\lvec(20 0)\lvec(10 0)\lfill f:0.8
\htext(16 2){\tiny $0$}

\move(20 0)\lvec(30 10)\lvec(30 0)\lvec(20 0)\lfill f:0.8
\htext(26 2){\tiny $1$}

\move(30 0)\lvec(40 10)\lvec(40 0)\lvec(30 0)\lfill f:0.8
\htext(36 2){\tiny $0$}

\move(40 0)\lvec(50 10)\lvec(50 0)\lvec(40 0)\lfill f:0.8
\htext(46 2){\tiny $1$}
\move(10 10)\lvec(20 10)\lvec(20 20)\lvec(10 20)\lvec(10
10)\htext(13 13){$_2$}

\move(20 10)\lvec(30 10)\lvec(30 20)\lvec(20 20)\lvec(20
10)\htext(23 13){$_2$}

\move(30 10)\lvec(40 10)\lvec(40 20)\lvec(30 20)\lvec(30
10)\htext(33 13){$_2$}

\move(40 10)\lvec(50 10)\lvec(50 20)\lvec(40 20)\lvec(40
10)\htext(43 13){$_2$}
\move(10 20)\lvec(20 20)\lvec(20 30)\lvec(10 30)\lvec(10
20)\htext(13 26){\tiny $3$}

\move(20 20)\lvec(30 20)\lvec(30 30)\lvec(20 30)\lvec(20
20)\htext(23 26){\tiny $3$}

\move(30 20)\lvec(40 20)\lvec(40 30)\lvec(30 30)\lvec(30
20)\htext(33 26){\tiny $3$}

\move(40 20)\lvec(50 20)\lvec(50 30)\lvec(40 30)\lvec(40
20)\htext(43 26){\tiny $3$}

\move(10 20)\lvec(20 20)\lvec(20 25)\lvec(10 25)\lvec(10 20)
\htext(13 21){\tiny $3$}

\move(20 20)\lvec(30 20)\lvec(30 25)\lvec(20 25)\lvec(20 20)
\htext(23 21){\tiny $3$}

\move(30 20)\lvec(40 20)\lvec(40 25)\lvec(30 25)\lvec(30 20)
\htext(33 21){\tiny $3$}

\move(40 20)\lvec(50 20)\lvec(50 25)\lvec(40 25)\lvec(40 20)
\htext(43 21){\tiny $3$}
\move(10 30)\lvec(20 30)\lvec(20 40)\lvec(10 40)\lvec(10
30)\htext(13 33){$_2$}

\move(20 30)\lvec(30 30)\lvec(30 40)\lvec(20 40)\lvec(20
30)\htext(23 33){$_2$}

\move(30 30)\lvec(40 30)\lvec(40 40)\lvec(30 40)\lvec(30
30)\htext(33 33){$_2$}

\move(40 30)\lvec(50 30)\lvec(50 40)\lvec(40 40)\lvec(40
30)\htext(43 33){$_2$}
%

\move(20 40)\lvec(30 50)\lvec(30 40)\lvec(20 40) \htext(26
42){\tiny $1$}

\move(30 40)\lvec(40 50)\lvec(40 40)\lvec(30 40) \htext(36
42){\tiny $0$}

\move(40 40)\lvec(50 50)\lvec(50 40)\lvec(40 40) \htext(46
42){\tiny $1$}

\end{texdraw}}\hskip 3mm $+$ \hskip 3mm $q^2$
\raisebox{-0.5\height}{
\begin{texdraw}
\drawdim em \setunitscale 0.13 \linewd 0.5


\move(10 0)\lvec(20 0)\lvec(20 10)\lvec(10 10)\lvec(10 0)\htext(12
6){\tiny $1$}

\move(20 0)\lvec(30 0)\lvec(30 10)\lvec(20 10)\lvec(20 0)\htext(22
6){\tiny $0$}

\move(30 0)\lvec(40 0)\lvec(40 10)\lvec(30 10)\lvec(30 0)\htext(32
6){\tiny $1$}

\move(40 0)\lvec(50 0)\lvec(50 10)\lvec(40 10)\lvec(40 0)\htext(42
6){\tiny $0$}
%

\move(10 0)\lvec(20 10)\lvec(20 0)\lvec(10 0)\lfill f:0.8
\htext(16 2){\tiny $0$}

\move(20 0)\lvec(30 10)\lvec(30 0)\lvec(20 0)\lfill f:0.8
\htext(26 2){\tiny $1$}

\move(30 0)\lvec(40 10)\lvec(40 0)\lvec(30 0)\lfill f:0.8
\htext(36 2){\tiny $0$}

\move(40 0)\lvec(50 10)\lvec(50 0)\lvec(40 0)\lfill f:0.8
\htext(46 2){\tiny $1$}
\move(10 10)\lvec(20 10)\lvec(20 20)\lvec(10 20)\lvec(10
10)\htext(13 13){$_2$}

\move(20 10)\lvec(30 10)\lvec(30 20)\lvec(20 20)\lvec(20
10)\htext(23 13){$_2$}

\move(30 10)\lvec(40 10)\lvec(40 20)\lvec(30 20)\lvec(30
10)\htext(33 13){$_2$}

\move(40 10)\lvec(50 10)\lvec(50 20)\lvec(40 20)\lvec(40
10)\htext(43 13){$_2$}
\move(10 20)\lvec(20 20)\lvec(20 30)\lvec(10 30)\lvec(10
20)\htext(13 26){\tiny $3$}

\move(20 20)\lvec(30 20)\lvec(30 30)\lvec(20 30)\lvec(20
20)\htext(23 26){\tiny $3$}

\move(30 20)\lvec(40 20)\lvec(40 30)\lvec(30 30)\lvec(30
20)\htext(33 26){\tiny $3$}

\move(40 20)\lvec(50 20)\lvec(50 30)\lvec(40 30)\lvec(40
20)\htext(43 26){\tiny $3$}

\move(10 20)\lvec(20 20)\lvec(20 25)\lvec(10 25)\lvec(10 20)
\htext(13 21){\tiny $3$}

\move(20 20)\lvec(30 20)\lvec(30 25)\lvec(20 25)\lvec(20 20)
\htext(23 21){\tiny $3$}

\move(30 20)\lvec(40 20)\lvec(40 25)\lvec(30 25)\lvec(30 20)
\htext(33 21){\tiny $3$}

\move(40 20)\lvec(50 20)\lvec(50 25)\lvec(40 25)\lvec(40 20)
\htext(43 21){\tiny $3$}
\move(10 30)\lvec(20 30)\lvec(20 40)\lvec(10 40)\lvec(10
30)\htext(13 33){$_2$}

\move(20 30)\lvec(30 30)\lvec(30 40)\lvec(20 40)\lvec(20
30)\htext(23 33){$_2$}

\move(30 30)\lvec(40 30)\lvec(40 40)\lvec(30 40)\lvec(30
30)\htext(33 33){$_2$}

\move(40 30)\lvec(50 30)\lvec(50 40)\lvec(40 40)\lvec(40
30)\htext(43 33){$_2$}
\move(10 40)\lvec(20 50)\lvec(20 40)\lvec(10 40) \htext(16
42){\tiny $0$}

\move(20 40)\lvec(30 50)\lvec(30 40)\lvec(20 40) \htext(26
42){\tiny $1$}

\move(30 40)\lvec(40 50)\lvec(40 40)\lvec(30 40) \htext(36
42){\tiny $0$}

\move(40 40)\lvec(50 50)\lvec(50 40)\lvec(40 40) \htext(46
42){\tiny $1$}

\end{texdraw}}\vskip 5mm

\hskip 4cm $+$ \hskip 3mm $q^8$ \raisebox{-0.5\height}{
\begin{texdraw}
\drawdim em \setunitscale 0.13 \linewd 0.5


\move(10 0)\lvec(20 0)\lvec(20 10)\lvec(10 10)\lvec(10 0)\htext(12
6){\tiny $1$}

\move(20 0)\lvec(30 0)\lvec(30 10)\lvec(20 10)\lvec(20 0)\htext(22
6){\tiny $0$}

\move(30 0)\lvec(40 0)\lvec(40 10)\lvec(30 10)\lvec(30 0)\htext(32
6){\tiny $1$}

\move(40 0)\lvec(50 0)\lvec(50 10)\lvec(40 10)\lvec(40 0)\htext(42
6){\tiny $0$}
%

\move(10 0)\lvec(20 10)\lvec(20 0)\lvec(10 0)\lfill f:0.8
\htext(16 2){\tiny $0$}

\move(20 0)\lvec(30 10)\lvec(30 0)\lvec(20 0)\lfill f:0.8
\htext(26 2){\tiny $1$}

\move(30 0)\lvec(40 10)\lvec(40 0)\lvec(30 0)\lfill f:0.8
\htext(36 2){\tiny $0$}

\move(40 0)\lvec(50 10)\lvec(50 0)\lvec(40 0)\lfill f:0.8
\htext(46 2){\tiny $1$}
\move(10 10)\lvec(20 10)\lvec(20 20)\lvec(10 20)\lvec(10
10)\htext(13 13){$_2$}

\move(20 10)\lvec(30 10)\lvec(30 20)\lvec(20 20)\lvec(20
10)\htext(23 13){$_2$}

\move(30 10)\lvec(40 10)\lvec(40 20)\lvec(30 20)\lvec(30
10)\htext(33 13){$_2$}

\move(40 10)\lvec(50 10)\lvec(50 20)\lvec(40 20)\lvec(40
10)\htext(43 13){$_2$}
\move(10 20)\lvec(20 20)\lvec(20 30)\lvec(10 30)\lvec(10
20)\htext(13 26){\tiny $3$}

\move(20 20)\lvec(30 20)\lvec(30 30)\lvec(20 30)\lvec(20
20)\htext(23 26){\tiny $3$}

\move(30 20)\lvec(40 20)\lvec(40 30)\lvec(30 30)\lvec(30
20)\htext(33 26){\tiny $3$}

\move(40 20)\lvec(50 20)\lvec(50 30)\lvec(40 30)\lvec(40
20)\htext(43 26){\tiny $3$}

\move(10 20)\lvec(20 20)\lvec(20 25)\lvec(10 25)\lvec(10 20)
\htext(13 21){\tiny $3$}

\move(20 20)\lvec(30 20)\lvec(30 25)\lvec(20 25)\lvec(20 20)
\htext(23 21){\tiny $3$}

\move(30 20)\lvec(40 20)\lvec(40 25)\lvec(30 25)\lvec(30 20)
\htext(33 21){\tiny $3$}

\move(40 20)\lvec(50 20)\lvec(50 25)\lvec(40 25)\lvec(40 20)
\htext(43 21){\tiny $3$}
\move(10 30)\lvec(20 30)\lvec(20 40)\lvec(10 40)\lvec(10
30)\htext(13 33){$_2$}

\move(20 30)\lvec(30 30)\lvec(30 40)\lvec(20 40)\lvec(20
30)\htext(23 33){$_2$}

\move(30 30)\lvec(40 30)\lvec(40 40)\lvec(30 40)\lvec(30
30)\htext(33 33){$_2$}

\move(40 30)\lvec(50 30)\lvec(50 40)\lvec(40 40)\lvec(40
30)\htext(43 33){$_2$}
%

\move(20 40)\lvec(30 50)\lvec(20 50)\lvec(20 40) \htext(22
46){\tiny $0$}

\move(30 40)\lvec(40 50)\lvec(30 50)\lvec(30 40) \htext(32
46){\tiny $1$}

\move(40 40)\lvec(50 50)\lvec(40 50)\lvec(40 40) \htext(42
46){\tiny $0$}

\move(40 40)\lvec(50 50)\lvec(50 40)\lvec(40 40) \htext(46
42){\tiny $1$}
\end{texdraw}}\hskip 5mm
$+$\hskip 3mm $q^4$ \raisebox{-0.5\height}{
\begin{texdraw}
\drawdim em \setunitscale 0.13 \linewd 0.5


\move(10 0)\lvec(20 0)\lvec(20 10)\lvec(10 10)\lvec(10 0)\htext(12
6){\tiny $1$}

\move(20 0)\lvec(30 0)\lvec(30 10)\lvec(20 10)\lvec(20 0)\htext(22
6){\tiny $0$}

\move(30 0)\lvec(40 0)\lvec(40 10)\lvec(30 10)\lvec(30 0)\htext(32
6){\tiny $1$}

\move(40 0)\lvec(50 0)\lvec(50 10)\lvec(40 10)\lvec(40 0)\htext(42
6){\tiny $0$}
%

\move(10 0)\lvec(20 10)\lvec(20 0)\lvec(10 0)\lfill f:0.8
\htext(16 2){\tiny $0$}

\move(20 0)\lvec(30 10)\lvec(30 0)\lvec(20 0)\lfill f:0.8
\htext(26 2){\tiny $1$}

\move(30 0)\lvec(40 10)\lvec(40 0)\lvec(30 0)\lfill f:0.8
\htext(36 2){\tiny $0$}

\move(40 0)\lvec(50 10)\lvec(50 0)\lvec(40 0)\lfill f:0.8
\htext(46 2){\tiny $1$}
\move(10 10)\lvec(20 10)\lvec(20 20)\lvec(10 20)\lvec(10
10)\htext(13 13){$_2$}

\move(20 10)\lvec(30 10)\lvec(30 20)\lvec(20 20)\lvec(20
10)\htext(23 13){$_2$}

\move(30 10)\lvec(40 10)\lvec(40 20)\lvec(30 20)\lvec(30
10)\htext(33 13){$_2$}

\move(40 10)\lvec(50 10)\lvec(50 20)\lvec(40 20)\lvec(40
10)\htext(43 13){$_2$}
\move(10 20)\lvec(20 20)\lvec(20 30)\lvec(10 30)\lvec(10
20)\htext(13 26){\tiny $3$}

\move(20 20)\lvec(30 20)\lvec(30 30)\lvec(20 30)\lvec(20
20)\htext(23 26){\tiny $3$}

\move(30 20)\lvec(40 20)\lvec(40 30)\lvec(30 30)\lvec(30
20)\htext(33 26){\tiny $3$}

\move(40 20)\lvec(50 20)\lvec(50 30)\lvec(40 30)\lvec(40
20)\htext(43 26){\tiny $3$}

\move(10 20)\lvec(20 20)\lvec(20 25)\lvec(10 25)\lvec(10 20)
\htext(13 21){\tiny $3$}

\move(20 20)\lvec(30 20)\lvec(30 25)\lvec(20 25)\lvec(20 20)
\htext(23 21){\tiny $3$}

\move(30 20)\lvec(40 20)\lvec(40 25)\lvec(30 25)\lvec(30 20)
\htext(33 21){\tiny $3$}

\move(40 20)\lvec(50 20)\lvec(50 25)\lvec(40 25)\lvec(40 20)
\htext(43 21){\tiny $3$}
\move(10 30)\lvec(20 30)\lvec(20 40)\lvec(10 40)\lvec(10
30)\htext(13 33){$_2$}

\move(20 30)\lvec(30 30)\lvec(30 40)\lvec(20 40)\lvec(20
30)\htext(23 33){$_2$}

\move(30 30)\lvec(40 30)\lvec(40 40)\lvec(30 40)\lvec(30
30)\htext(33 33){$_2$}

\move(40 30)\lvec(50 30)\lvec(50 40)\lvec(40 40)\lvec(40
30)\htext(43 33){$_2$}
%

\move(20 40)\lvec(30 50)\lvec(30 40)\lvec(20 40) \htext(26
42){\tiny $1$}

\move(30 40)\lvec(40 50)\lvec(40 40)\lvec(30 40) \htext(36
42){\tiny $0$}

\move(40 40)\lvec(50 50)\lvec(50 40)\lvec(40 40) \htext(46
42){\tiny $1$}

\move(40 40)\lvec(50 50)\lvec(40 50)\lvec(40 40) \htext(42
46){\tiny $0$}
\end{texdraw}}\quad .\vskip 5mm
}
\end{ex}

With these actions, we have
\begin{thm}\label{Fock space}
The Fock space ${\mathcal F}(\Lambda)$ is a $U_q(\frak{g})$-module
in the category ${\mathcal O}_{int}$.
\end{thm}
To prove this theorem, we need to verify that all the defining
relations in \eqref{defining rels} hold in $\mathcal{F}(\Lambda)$.
First, it is straightforward to verify that the following
relations hold
\begin{equation}\label {rel1}
\begin{split}
q^hq^{h'}Y&=q^{h+h'}Y, \\
q^h e_i q^{-h}Y&=q^{\alpha_i(h)}e_i Y, \\
q^h f_i q^{-h}Y&=q^{-\alpha_i(h)}f_i Y,
\end{split}
\end{equation}
for $Y\in {\mathcal Z}(\Lambda)$, $i\in I$ and $h,h'\in P^{\vee}$.
Also, it is clear that $e_i$ and $f_i$ ($i\in I$) act locally
nilpotently on ${\mathcal F}(\Lambda)$. Therefore, by Proposition
B.1 in \cite{KMPY}, we have only to show that
\begin{equation}\label {rel2}
\begin{split}
(e_if_j-f_je_i)Y=\delta_{ij}\frac{K_i-K_i^{-1}}{q_i-q_i^{-1}}Y
\end{split}
\end{equation}
for $Y\in {\mathcal Z}(\Lambda)$ and $i,j\in I$.

The rest of this subsection will be devoted to proving the
relation \eqref{rel2}. We first investigate the {\it local
behavior} of $U_q(\frak{g})$-action on ${\mathcal F}(\Lambda)$.

Fix $i\in I$ and let $Y$ be a proper Young wall. We will decompose
$Y$ into a sequence $Y=(Y_0,Y_1,\cdots,Y_r)$ of parts, reading
from left to right, which are called the {\it $i$-component} of
$Y$.\vskip 2mm

{\bf Case 1.} Suppose that $i$-blocks are of type I. Observe that
the admissible $i$-slots and removable $i$-blocks in $Y$ appear in
one of the following situations:\vskip 5mm

\noindent ( $\rm I_{\infty}$) \hskip 1cm \raisebox{-0.6\height}
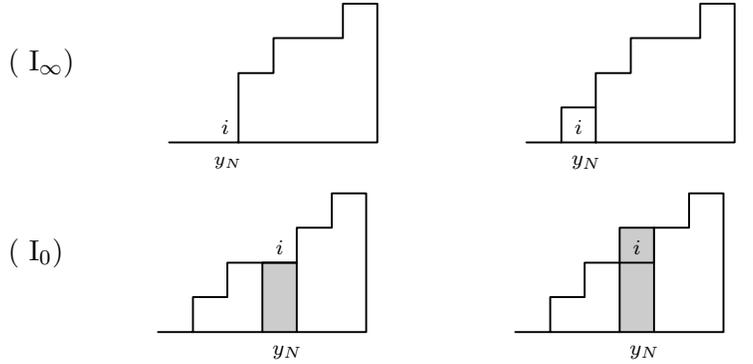
{
\begin{texdraw} \drawdim em \setunitscale 0.12 \linewd 0.5
\move(0 0)\lvec(60 0)\lvec(60 40)\lvec(50 40)\lvec(50 30)\lvec(30
30)\lvec(30 20)\lvec(20 20)\lvec(20 0)\htext(15
2.5){$_i$}\htext(13 -7.5){\small $_{y_N}$}
\end{texdraw}\hskip 2cm
\begin{texdraw} \drawdim em \setunitscale 0.12 \linewd 0.5
\move(0 0)\lvec(60 0)\lvec(60 40)\lvec(50 40)\lvec(50 30)\lvec(30
30)\lvec(30 20)\lvec(20 20)\lvec(20 0)\move(10 0)\lvec(10
10)\lvec(20 10)\htext(14 2.5){$_i$}\htext(13 -7.5){$_{y_N}$}
\end{texdraw}}
\vskip 3mm

\noindent ( $\rm I_0$)  \hskip 1cm \raisebox{-0.6\height}{
\begin{texdraw} \drawdim em \setunitscale 0.12 \linewd 0.5
\move(30 0)\lvec(30 20)\lvec(40 20)\lvec(40 0)\lfill f:0.8 \move(0
0)\lvec(60 0)\lvec(60 40)\lvec(50 40)\lvec(50 30)\lvec(40
30)\lvec(40 20)\lvec(20 20)\lvec(20 10)\lvec(10 10)\lvec(10
0)\htext(34 22.5){$_i$}\htext(33 -7.5){$_{y_N}$}
\end{texdraw}\hskip 2cm
\begin{texdraw} \drawdim em \setunitscale 0.12 \linewd 0.5
\move(30 0)\lvec(30 30)\lvec(40 30)\lvec(40 0)\lfill f:0.8 \move(0
0)\lvec(60 0)\lvec(60 40)\lvec(50 40)\lvec(50 30)\lvec(40
30)\lvec(40 20)\lvec(20 20)\lvec(20 10)\lvec(10 10)\lvec(10
0)\htext(34 22.5){$_i$}\htext(33 -7.5){$_{y_N}$}
\end{texdraw}}\quad .\vskip 5mm

In case $\rm I_{\infty}$, let $y_N$ be the column containing the
admissible $i$-slot in the ground-state wall or the removable
$i$-block that touches the ground-state wall as indicated in the
above figure. We denote by $Y_0$ the part of $Y$ consisting of
$y_N$ and the blocks in the ground-state wall lying in the left of
$y_N$, and call it an {\it $i$-component of type $\rm
I_{\infty}$}.

In case $\rm I_0$, the column which is $i$-admissible or
$i$-removable is called an {\it $i$-component of type $\rm I_0$}.
If $y_N$ is the left-most $i$-component of type $\rm I_0$ and
there is no $i$-component of type $\rm I_{\infty}$, we denote by
$Y_0$ the part of $Y$ consisting of the blocks lying in the left
of $y_N$. In this case, we call $Y_0$ a {\it trivial
$i$-component}. The parts of $Y$ lying between two $i$-components
of type $\rm I_{\infty}$ or $\rm I_0$ will also be called the {\it
trivial $i$-components}.

In this way, we obtain a unique decomposition
$Y=(Y_0,Y_1,\cdots,Y_r)$ of $Y$, where each of $Y_k$ is an
$i$-component of type $\rm I_{\infty}$ or $\rm I_0$, or a trivial
$i$-component.

\begin{ex}\label{comp1}{\rm
Let $\frak{g}=A^{(2)}_5$, $\Lambda=\Lambda_0$ and $i=2$. If
\begin{center}
$Y=$\hskip 3mm\raisebox{-0.5\height}{\begin{texdraw}\drawdim em
\setunitscale 0.13 \linewd 0.5

\move(10 0)\lvec(20 0)\lvec(20 10)\lvec(10 10)\lvec(10 0)

\move(20 0)\lvec(30 0)\lvec(30 10)\lvec(20 10)\lvec(20 0)

\move(30 0)\lvec(40 0)\lvec(40 10)\lvec(30 10)\lvec(30 0)

\move(40 0)\lvec(50 0)\lvec(50 10)\lvec(40 10)\lvec(40 0)
\move(10 10)\lvec(20 10)\lvec(20 20)\lvec(10 20)\lvec(10
10)\htext(13 13){$_2$}

\move(20 10)\lvec(30 10)\lvec(30 20)\lvec(20 20)\lvec(20
10)\htext(23 13){$_2$}

\move(30 10)\lvec(40 10)\lvec(40 20)\lvec(30 20)\lvec(30
10)\htext(33 13){$_2$}

\move(40 10)\lvec(50 10)\lvec(50 20)\lvec(40 20)\lvec(40
10)\htext(43 13){$_2$}
\move(20 20)\lvec(30 20)\lvec(30 30)\lvec(20 30)\lvec(20
20)\htext(23 23){$_3$}

\move(30 20)\lvec(40 20)\lvec(40 30)\lvec(30 30)\lvec(30
20)\htext(33 23){$_3$}

\move(40 20)\lvec(50 20)\lvec(50 30)\lvec(40 30)\lvec(40
20)\htext(43 23){$_3$}
\move(20 30)\lvec(30 30)\lvec(30 40)\lvec(20 40)\lvec(20
30)\htext(23 33){$_2$}

\move(30 30)\lvec(40 30)\lvec(40 40)\lvec(30 40)\lvec(30
30)\htext(33 33){$_2$}

\move(40 30)\lvec(50 30)\lvec(50 40)\lvec(40 40)\lvec(40
30)\htext(43 33){$_2$}
\move(20 40)\lvec(30 40)\lvec(30 50)\lvec(20 50)\lvec(20 40)

\move(30 40)\lvec(40 40)\lvec(40 50)\lvec(30 50)\lvec(30 40)

\move(40 40)\lvec(50 40)\lvec(50 50)\lvec(40 50)\lvec(40 40)
\move(30 50)\lvec(40 50)\lvec(40 60)\lvec(30 60)\lvec(30
50)\htext(33 53){$_2$}

\move(40 50)\lvec(50 50)\lvec(50 60)\lvec(40 60)\lvec(40
50)\htext(43 53){$_2$}
\move(40 60)\lvec(50 60)\lvec(50 70)\lvec(40 70)\lvec(40
60)\htext(43 63){$_3$}
\move(10 0)\lvec(20 10)\lvec(20 0)\lvec(10 0)\lfill f:0.8
\htext(11 5){\tiny $1$}\htext(16 1){\tiny $0$}

\move(20 0)\lvec(30 10)\lvec(30 0)\lvec(20 0)\lfill f:0.8
\htext(22 5){\tiny $0$}\htext(26 2){\tiny $1$}

\move(30 0)\lvec(40 10)\lvec(40 0)\lvec(30 0)\lfill f:0.8
\htext(32 5){\tiny $1$}\htext(36 2){\tiny $0$}

\move(40 0)\lvec(50 10)\lvec(50 0)\lvec(40 0)\lfill f:0.8
\htext(42 5){\tiny $0$}\htext(46 2){\tiny $1$}
\move(20 40)\lvec(30 50)\lvec(30 40)\lvec(20 40)\htext(22
45){\tiny $0$}\htext(26 42){\tiny $1$}

\move(30 40)\lvec(40 50)\lvec(40 40)\lvec(30 40)\htext(32
45){\tiny $1$}\htext(36 42){\tiny $0$}

\move(40 40)\lvec(50 50)\lvec(50 40)\lvec(40 40)\htext(42
45){\tiny $0$}\htext(46 42){\tiny $1$}
\end{texdraw}}\quad .
\end{center}\vskip 5mm

then we have $Y=(Y_0,Y_1,Y_2,Y_3)$, where\vskip 5mm

\hskip 1cm$Y_0=$ \raisebox{-0.4\height}{\begin{texdraw}\drawdim em
\setunitscale 0.13 \linewd 0.5

\htext(-.2 3){$\cdots$} \move(10 0)\lvec(20 10)\lvec(20 0)\lvec(10
0)\lfill f:0.8 \htext(16 1){\tiny $0$}

\move(20 0)\lvec(30 10)\lvec(30 0)\lvec(20 0)\lfill f:0.8
\htext(26 2){\tiny $1$}

\move(30 0)\lvec(40 10)\lvec(40 0)\lvec(30 0)\lfill f:0.8
\htext(36 2){\tiny $0$}

\move(40 0)\lvec(50 10)\lvec(50 0)\lvec(40 0)\lfill f:0.8
\htext(46 2){\tiny $1$}
\end{texdraw}} : trivial $i$-component, \vskip 5mm

\hskip 1cm $Y_1=$\hskip 2mm
\raisebox{-0.3\height}{\begin{texdraw}\drawdim em \setunitscale
0.13 \linewd 0.5

\move(10 10)\lvec(20 10)\lvec(20 20)\lvec(10 20)\lvec(10
10)\htext(13 13){$_2$}

\move(10 0)\lvec(20 10)\lvec(20 0)\lvec(10 0)\lvec(10 10)\lfill
f:0.8 \htext(11 5){\tiny $1$}\htext(16 1){\tiny $0$}
\end{texdraw}} : $i$-component of type $\rm I_0$,\vskip 5mm

\hskip 1cm$Y_2=$ \raisebox{-0.4\height}{\begin{texdraw}\drawdim em
\setunitscale 0.13 \linewd 0.5

\move(20 0)\lvec(30 10)\lvec(30 0)\lvec(20 0)\lvec(20 10)\lfill
f:0.8 \htext(22 5){\tiny $0$}\htext(26 2){\tiny $1$}

\move(30 0)\lvec(40 10)\lvec(40 0)\lvec(30 0)\lfill f:0.8
\htext(32 5){\tiny $1$}\htext(36 2){\tiny $0$}

\move(20 10)\lvec(30 10)\lvec(30 20)\lvec(20 20)\lvec(20
10)\htext(23 13){$_2$}

\move(30 10)\lvec(40 10)\lvec(40 20)\lvec(30 20)\lvec(30
10)\htext(33 13){$_2$}

\move(20 20)\lvec(30 20)\lvec(30 30)\lvec(20 30)\lvec(20
20)\htext(23 23){$_3$}

\move(30 20)\lvec(40 20)\lvec(40 30)\lvec(30 30)\lvec(30
20)\htext(33 23){$_3$}

\move(20 30)\lvec(30 30)\lvec(30 40)\lvec(20 40)\lvec(20
30)\htext(23 33){$_2$}

\move(30 30)\lvec(40 30)\lvec(40 40)\lvec(30 40)\lvec(30
30)\htext(33 33){$_2$}

\move(20 40)\lvec(30 50)\lvec(30 40)\lvec(20 40)\lvec(20
50)\lvec(30 50)\htext(22 45){\tiny $0$}\htext(26 42){\tiny $1$}

\move(30 40)\lvec(40 50)\lvec(40 40)\lvec(30 40)\htext(32
45){\tiny $1$}\htext(36 42){\tiny $0$}

\move(30 50)\lvec(40 50)\lvec(40 60)\lvec(30 60)\lvec(30
50)\htext(33 53){$_2$}

\end{texdraw}} : trivial $i$-component and\vskip 5mm

\hskip 1cm$Y_3=$ \raisebox{-0.4\height}{\begin{texdraw}\drawdim em
\setunitscale 0.13 \linewd 0.5

\move(40 0)\lvec(50 0)\lvec(50 10)\lvec(40 10)\lvec(40 0)

\move(40 10)\lvec(50 10)\lvec(50 20)\lvec(40 20)\lvec(40
10)\htext(43 13){$_2$}

\move(40 20)\lvec(50 20)\lvec(50 30)\lvec(40 30)\lvec(40
20)\htext(43 23){$_3$}

\move(40 30)\lvec(50 30)\lvec(50 40)\lvec(40 40)\lvec(40
30)\htext(43 33){$_2$}

\move(40 40)\lvec(50 40)\lvec(50 50)\lvec(40 50)\lvec(40 40)

\move(40 50)\lvec(50 50)\lvec(50 60)\lvec(40 60)\lvec(40
50)\htext(43 53){$_2$}
\move(40 60)\lvec(50 60)\lvec(50 70)\lvec(40 70)\lvec(40
60)\htext(43 63){$_3$}
\move(40 0)\lvec(50 10)\lvec(50 0)\lvec(40 0)\lfill f:0.8
\htext(41 5){\tiny $0$}\htext(46 2){\tiny $1$}

\move(40 40)\lvec(50 50)\lvec(50 40)\lvec(40 40)\htext(41
45){\tiny $0$}\htext(46 42){\tiny $1$}
\end{texdraw}} : $i$-component of type $\rm I_0$.\vskip 5mm
}
\end{ex}

Let $Y$ be a proper Young wall with the decomposition
$Y=(Y_0,Y_1,\cdots,Y_r)$ into $i$-components. To each $Y_k$, we
associate a $U_{(i)}$-module $V_k$ as follows. Then we will view
$Y$ as $Y_0\otimes Y_1 \otimes \cdots \otimes Y_r$ inside
$V_0\otimes V_1 \otimes \cdots \otimes V_r$.

If $Y_k$ is a trivial $i$-component, then we associate the trivial
representation $V_k=U=\mathbb{Q}(q)u$, and we identify $Y_k$ with
$u$. If $Y_k$ is an $i$-component of type $\rm I_{\infty}$ or $\rm
I_0$, then we associate the 2-dimensional representation
$V_k=V=\mathbb{Q}(q)v_0 \oplus \mathbb{Q}(q)v_1$, where the
$U_{(i)}$-action is given by
\begin{equation}\label{I+}
\begin{split}
K_iv_0&=q_iv_0,  \quad K_iv_1=q_i^{-1}v_1, \\
e_iv_0&=0, \quad e_iv_1=v_0,\\
f_iv_0&=v_1, \quad f_iv_1=0.
\end{split}
\end{equation}
We identify the $i$-component $Y_k$ with a basis element of $V$ as
follows :\vskip 5mm

($\rm I_{\infty}$) \hskip 1cm $v_0 \leftrightarrow
$\raisebox{-0.6\height}{
\begin{texdraw}
\drawdim em \setunitscale 0.12 \linewd 0.5 \move(0 0)\lvec(30
0)\lvec(30 15)\htext(25 2.5){$_i$}
\end{texdraw}}\hskip 1cm
$v_1 \leftrightarrow $\raisebox{-0.4\height}{
\begin{texdraw}
\drawdim em \setunitscale 0.12 \linewd 0.5 \move(0 0)\lvec(30
0)\lvec(30 15)\htext(24 2.5){$_i$}\move(20 0)\lvec(20 10)\lvec(30
10)
\end{texdraw}}\vskip 5mm

\vskip 3mm

($\rm I_0$) \hskip 1cm $v_0 \leftrightarrow
$\raisebox{-0.6\height}{
\begin{texdraw}
\drawdim em \setunitscale 0.12 \linewd 0.5 \move(0 -10)\lvec(0
10)\lvec(10 10)\lvec(10 -10)\htext(4 12.5){$_i$}
\end{texdraw}}\hskip 2cm
$v_1 \leftrightarrow $\raisebox{-0.4\height}{
\begin{texdraw}
\drawdim em \setunitscale 0.12 \linewd 0.5 \move(0 -10)\lvec(0
10)\lvec(10 10)\lvec(10 -10) \move(0 10)\lvec(0 20)\lvec(10
20)\lvec(10 10)\lvec(0 10) \htext(4 13){$_i$}
\end{texdraw}}\quad .\vskip 5mm

Note that $V$ is isomorphic to the 2-dimensional irreducible
$U_{(i)}$-module $V(1)$ with the crystal basis $(L,B)$, where
\begin{equation*}
L=\mathbb{A}_0v_0 \oplus \mathbb{A}_0v_1, \ \
B=\{\overline{v_0},\overline{v_1} \},
\end{equation*}
and the crystal graph is given by
$$\overline{v_0}\stackrel{i}{\longrightarrow}\overline{v_1}.$$

\begin{ex}{\rm In Example \ref{comp1}, we have \vskip 5mm

$Y=$  \raisebox{-0.4\height}{\begin{texdraw}\drawdim em
\setunitscale 0.13 \linewd 0.5

\htext(-.2 3){$\cdots$} \move(10 0)\lvec(20 10)\lvec(20 0)\lvec(10
0)\lfill f:0.8 \htext(16 2){\tiny $0$}

\move(20 0)\lvec(30 10)\lvec(30 0)\lvec(20 0)\lfill f:0.8
\htext(26 2){\tiny $1$}

\move(30 0)\lvec(40 10)\lvec(40 0)\lvec(30 0)\lfill f:0.8
\htext(36 2){\tiny $0$}

\move(40 0)\lvec(50 10)\lvec(50 0)\lvec(40 0)\lfill f:0.8
\htext(46 2){\tiny $1$}
\end{texdraw}}\hskip 2mm $\otimes$ \hskip 2mm
\raisebox{-0.3\height}{\begin{texdraw}\drawdim em \setunitscale
0.13 \linewd 0.5

\move(10 10)\lvec(20 10)\lvec(20 20)\lvec(10 20)\lvec(10
10)\htext(13 13){$_2$}

\move(10 0)\lvec(20 10)\lvec(20 0)\lvec(10 0)\lvec(10 10)\lfill
f:0.8 \htext(11 5){\tiny $1$}\htext(16 1){\tiny $0$}
\end{texdraw}}\hskip 2mm $\otimes$ \hskip 2mm
\raisebox{-0.4\height}{\begin{texdraw}\drawdim em \setunitscale
0.13 \linewd 0.5

\move(20 0)\lvec(30 10)\lvec(30 0)\lvec(20 0)\lvec(20 10)\lfill
f:0.8 \htext(22 5){\tiny $0$}\htext(26 2){\tiny $1$}

\move(30 0)\lvec(40 10)\lvec(40 0)\lvec(30 0)\lfill f:0.8
\htext(32 5){\tiny $1$}\htext(36 2){\tiny $0$}

\move(20 10)\lvec(30 10)\lvec(30 20)\lvec(20 20)\lvec(20
10)\htext(23 13){$_2$}

\move(30 10)\lvec(40 10)\lvec(40 20)\lvec(30 20)\lvec(30
10)\htext(33 13){$_2$}

\move(20 20)\lvec(30 20)\lvec(30 30)\lvec(20 30)\lvec(20
20)\htext(23 23){$_3$}

\move(30 20)\lvec(40 20)\lvec(40 30)\lvec(30 30)\lvec(30
20)\htext(33 23){$_3$}

\move(20 30)\lvec(30 30)\lvec(30 40)\lvec(20 40)\lvec(20
30)\htext(23 33){$_2$}

\move(30 30)\lvec(40 30)\lvec(40 40)\lvec(30 40)\lvec(30
30)\htext(33 33){$_2$}

\move(20 40)\lvec(30 50)\lvec(30 40)\lvec(20 40)\lvec(20
50)\lvec(30 50)\htext(22 45){\tiny $0$}\htext(26 42){\tiny $1$}

\move(30 40)\lvec(40 50)\lvec(40 40)\lvec(30 40)\htext(32
45){\tiny $1$}\htext(36 42){\tiny $0$}

\move(30 50)\lvec(40 50)\lvec(40 60)\lvec(30 60)\lvec(30
50)\htext(33 53){$_2$}

\end{texdraw}}\hskip 2mm $\otimes$\hskip 2mm
\raisebox{-0.4\height}{\begin{texdraw}\drawdim em \setunitscale
0.13 \linewd 0.5

\move(40 0)\lvec(50 0)\lvec(50 10)\lvec(40 10)\lvec(40 0)

\move(40 10)\lvec(50 10)\lvec(50 20)\lvec(40 20)\lvec(40
10)\htext(43 13){$_2$}

\move(40 20)\lvec(50 20)\lvec(50 30)\lvec(40 30)\lvec(40
20)\htext(43 23){$_3$}

\move(40 30)\lvec(50 30)\lvec(50 40)\lvec(40 40)\lvec(40
30)\htext(43 33){$_2$}

\move(40 40)\lvec(50 40)\lvec(50 50)\lvec(40 50)\lvec(40 40)

\move(40 50)\lvec(50 50)\lvec(50 60)\lvec(40 60)\lvec(40
50)\htext(43 53){$_2$}
\move(40 60)\lvec(50 60)\lvec(50 70)\lvec(40 70)\lvec(40
60)\htext(43 63){$_3$}
\move(40 0)\lvec(50 10)\lvec(50 0)\lvec(40 0)\lfill f:0.8
\htext(41 5){\tiny $0$}\htext(46 2){\tiny $1$}

\move(40 40)\lvec(50 50)\lvec(50 40)\lvec(40 40)\htext(41
45){\tiny $0$}\htext(46 42){\tiny $1$}
\end{texdraw}} \hskip 2mm $\in U\otimes V \otimes U \otimes V$.

}
\end{ex}

\vskip 3mm {\bf Case 2.} Suppose that the $i$-blocks are of type
II. In this case, the admissible $i$-slots and the removable
$i$-blocks in $Y$ appear in one of the following situations :
\vskip 5mm

($\rm II_{\infty}$) \hskip 1cm \raisebox{-0.6\height}
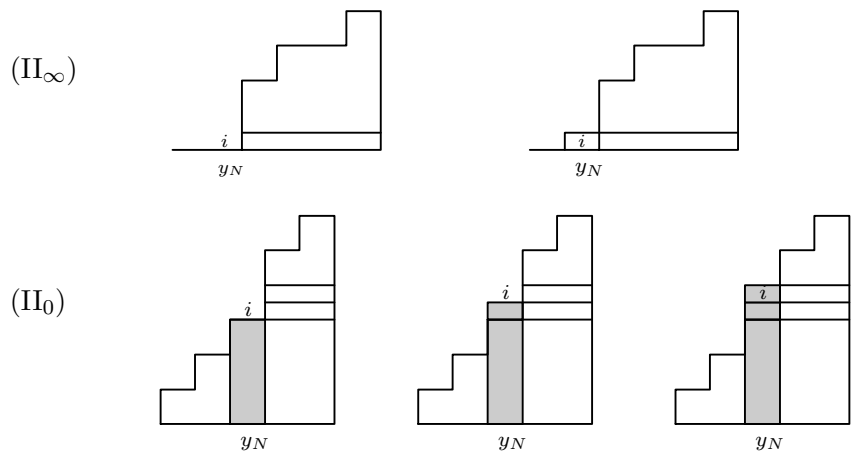
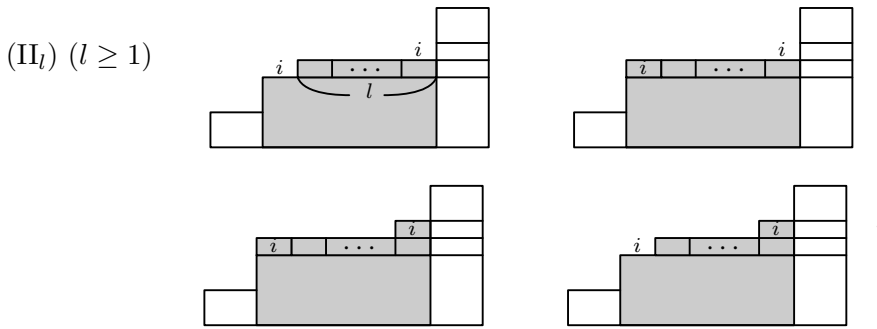
{
\begin{texdraw} \drawdim em \setunitscale 0.12 \linewd 0.5
\move(0 0)\lvec(60 0)\lvec(60 40)\lvec(50 40)\lvec(50 30)\lvec(30
30)\lvec(30 20)\lvec(20 20)\lvec(20 0)\htext(14 1){\small
$_i$}\htext(13 -7.5){\small $_{y_N}$}\move(20 5)\lvec(60 5)
\end{texdraw}\hskip 2cm
\begin{texdraw} \drawdim em \setunitscale 0.12 \linewd 0.5
\move(0 0)\lvec(60 0)\lvec(60 40)\lvec(50 40)\lvec(50 30)\lvec(30
30)\lvec(30 20)\lvec(20 20)\lvec(20 0)\move(10 0)\lvec(10
5)\lvec(60 5)\htext(14 1){\small $_i$}\htext(13 -7.5){$_{y_N}$}
\end{texdraw}}\vskip 5mm

($\rm II_0$) \hskip 1cm \raisebox{-0.6\height}{
\begin{texdraw}
\drawdim em \setunitscale 0.12 \linewd 0.5

\move(0 -20)\lvec(0 10)\lvec(10 10)\lvec(10 -20)\lvec(0 -20)\lfill
f:0.8 \move(-20 -20)\lvec(30 -20)\lvec(30 40)\lvec(20 40)\lvec(20
30)\lvec(10 30)\lvec(10 10)\lvec(0 10)\lvec(0 0)\lvec(-10
0)\lvec(-10 -10)\lvec(-20 -10)\lvec(-20 -20)

\move(10 10)\lvec(30 10)\move(10 15)\lvec(30 15)\move(10
20)\lvec(30 20)

\htext(3 -27.5){$_{y_N}$}\htext(4 11){$_i$}
\end{texdraw}}\hskip 10mm
\raisebox{-0.6\height}{
\begin{texdraw}
\drawdim em \setunitscale 0.12 \linewd 0.5

\move(0 -20)\lvec(0 10)\lvec(10 10)\lvec(10 -20)\lvec(0 -20)\lfill
f:0.8

\move(0 10)\lvec(0 15)\lvec(10 15)\lvec(10 10)\lvec(0 10)\lfill
f:0.8

\move(-20 -20)\lvec(30 -20)\lvec(30 40)\lvec(20 40)\lvec(20
30)\lvec(10 30)\lvec(10 10)\lvec(0 10)\lvec(0 0)\lvec(-10
0)\lvec(-10 -10)\lvec(-20 -10)\lvec(-20 -20)

\move(10 10)\lvec(30 10)\move(10 15)\lvec(30 15)\move(10
20)\lvec(30 20)

\htext(3 -27.5){$_{y_N}$}\htext(4 16){$_i$}
\end{texdraw}}\hskip 10mm
\raisebox{-0.6\height}{
\begin{texdraw}
\drawdim em \setunitscale 0.12 \linewd 0.5

\move(0 -20)\lvec(0 10)\lvec(10 10)\lvec(10 -20)\lvec(0 -20)\lfill
f:0.8

\move(0 10)\lvec(0 15)\lvec(10 15)\lvec(10 10)\lvec(0 10)\lfill
f:0.8

\move(0 15)\lvec(0 20)\lvec(10 20)\lvec(10 15)\lvec(0 15)\lfill
f:0.8

\move(-20 -20)\lvec(30 -20)\lvec(30 40)\lvec(20 40)\lvec(20
30)\lvec(10 30)\lvec(10 10)\lvec(0 10)\lvec(0 0)\lvec(-10
0)\lvec(-10 -10)\lvec(-20 -10)\lvec(-20 -20)

\move(10 10)\lvec(30 10)\move(10 15)\lvec(30 15)\move(10
20)\lvec(30 20)

\htext(3 -27.5){$_{y_N}$}\htext(4 16){$_i$}
\end{texdraw}}\vskip 5mm

(${\rm II}_l$) ($l\geq 1$) \hskip 5mm \raisebox{-0.6\height}{
\begin{texdraw}
\drawdim em \setunitscale 0.12 \linewd 0.5

\move(0 -10)\lvec(0 10)\lvec(50 10)\lvec(50 -10)\lvec(0 -10)\lfill
f:0.8

\move(10 10)\lvec(10 15)\lvec(50 15)\lvec(50 10)\lvec(10 10)
\lfill f:0.8

\move(20 10) \lvec(20 15) \move(40 10) \lvec(40 15) \htext(25
10){$\cdots$}

\move(50 10)\lvec(50 30)\lvec(65 30)\lvec(65 -10)\lvec(-15
-10)\lvec(-15 0)\lvec(0 0)

\move(50 10)\lvec(65 10)\move(50 15)\lvec(65 15)\move(50
20)\lvec(65 20)

\move(10 10)\clvec(12 5)(23 5)(25 5)

\move(50 10)\clvec(48 5)(37 5)(35 5)

\htext(30 4){$_{l}$}

\htext(4 11){$_i$}\htext(44 16){$_i$}
\end{texdraw}}\hskip 10mm
\raisebox{-0.6\height}{
\begin{texdraw}
\drawdim em \setunitscale 0.12 \linewd 0.5

\move(0 -10)\lvec(0 10)\lvec(50 10)\lvec(50 -10)\lvec(0 -10)\lfill
f:0.8

\move(0 10) \lvec(0 15) \lvec(10 15)\lvec(10 10)\lvec(0 10)\lfill
f:0.8

\move(10 10)\lvec(10 15)\lvec(50 15)\lvec(50 10)\lvec(10 10)
\lfill f:0.8

\move(20 10) \lvec(20 15) \move(40 10) \lvec(40 15) \htext(25
10){$\cdots$}

\move(50 10)\lvec(50 30)\lvec(65 30)\lvec(65 -10)\lvec(-15
-10)\lvec(-15 0)\lvec(0 0)

\move(50 10)\lvec(65 10)\move(50 15)\lvec(65 15)\move(50
20)\lvec(65 20)

\htext(4 11){$_i$}\htext(44 16){$_i$}
\end{texdraw}}\vskip 5mm

\hskip 2.5cm\raisebox{-0.7\height}{
\begin{texdraw}
\drawdim em \setunitscale 0.12 \linewd 0.5

\move(0 -10)\lvec(0 10)\lvec(50 10)\lvec(50 -10)\lvec(0 -10)\lfill
f:0.8

\move(0 10) \lvec(0 15) \lvec(10 15)\lvec(10 10)\lvec(0 10)\lfill
f:0.8

\move(40 15) \lvec(40 20) \lvec(50 20) \lvec(50 15)\lfill f:0.8

\move(10 10)\lvec(10 15)\lvec(50 15)\lvec(50 10)\lvec(10 10)
\lfill f:0.8

\move(20 10) \lvec(20 15) \move(40 10) \lvec(40 15) \htext(25
10){$\cdots$}

\move(50 10)\lvec(50 30)\lvec(65 30)\lvec(65 -10)\lvec(-15
-10)\lvec(-15 0)\lvec(0 0)

\move(50 10)\lvec(65 10)\move(50 15)\lvec(65 15)\move(50
20)\lvec(65 20)

\htext(4 11){$_i$}\htext(44 16){$_i$}
\end{texdraw}}\hskip 1cm
\raisebox{-0.7\height}{
\begin{texdraw}
\drawdim em \setunitscale 0.12 \linewd 0.5

\move(0 -10)\lvec(0 10)\lvec(50 10)\lvec(50 -10)\lvec(0 -10)\lfill
f:0.8

\move(40 15) \lvec(40 20) \lvec(50 20) \lvec(50 15)\lfill f:0.8

\move(10 10)\lvec(10 15)\lvec(50 15)\lvec(50 10)\lvec(10 10)
\lfill f:0.8

\move(20 10) \lvec(20 15) \move(40 10) \lvec(40 15) \htext(25
10){$\cdots$}

\move(50 10)\lvec(50 30)\lvec(65 30)\lvec(65 -10)\lvec(-15
-10)\lvec(-15 0)\lvec(0 0)

\move(50 10)\lvec(65 10)\move(50 15)\lvec(65 15)\move(50
20)\lvec(65 20)

\htext(4 11){$_i$}\htext(44 16){$_i$}
\end{texdraw}}\quad .\vskip 5mm

In case $\rm II_{\infty}$, let $y_N$ be the column containing the
admissible $i$-slot in the ground-state wall or the removable
$i$-block that touches the ground-state wall as indicated in the
above figure. We denote by $Y_0$ the part of $Y$ consisting of
$y_N$ and the blocks in the ground-state wall lying in the left of
$y_N$, and call it an {\it $i$-component of type $\rm
II_{\infty}$}.

In case $\rm II_0$, the column which is $i$-admissible or
$i$-removable is called an {\it $i$-component of type $\rm II_0$}.
In case ${\rm II}_l$, the whole shaded part containing an
admissible $i$-slot or a removable $i$-block will be called an
{\it $i$-component of type ${\rm II}_l$}.

Let $Y_1$ be the left-most $i$-component of type ${\rm II}_l$
($l\geq 0$). If there is no $i$-component of type $\rm
II_{\infty}$, then we denote by $Y_0$ the part of $Y$ consisting
of the blocks lying in the left of $Y_1$, and call it a {\it
trivial $i$-component}. The parts of $Y$ lying between two
$i$-components of type $\rm II_{\infty}$ or ${\rm II}_l$ ($l\geq
0$) will also be called the {\it trivial $i$-components}.

In this way, we obtain a unique decomposition
$Y=(Y_0,Y_1,\cdots,Y_r)$ of $Y$, where each $Y_k$ is an
$i$-component of type $\rm II_{\infty}$, ${\rm II}_l$ ($l\geq 0$),
or a trivial $i$-component.

\begin{ex}\label{comp2}{\rm
Let $\frak{g}=A^{(2)}_4$, $\Lambda=\Lambda_0$ and $i=0$. If \vskip
5mm

\begin{center}
$Y=$\raisebox{-0.5\height}{
\begin{texdraw}
\drawdim em \setunitscale 0.13 \linewd 0.5

\move(-10 0)\lvec(0 0)\lvec(0 10)\lvec(-10 10)\lvec(-10
0)\htext(-7 1){\tiny $0$}

\move(0 0)\lvec(10 0)\lvec(10 10)\lvec(0 10)\lvec(0 0)\htext(3
1){\tiny $0$}

\move(10 0)\lvec(20 0)\lvec(20 10)\lvec(10 10)\lvec(10 0)\htext(13
1){\tiny $0$}

\move(20 0)\lvec(30 0)\lvec(30 10)\lvec(20 10)\lvec(20 0)\htext(23
1){\tiny $0$}

\move(30 0)\lvec(40 0)\lvec(40 10)\lvec(30 10)\lvec(30 0)\htext(33
1){\tiny $0$}

\move(40 0)\lvec(50 0)\lvec(50 10)\lvec(40 10)\lvec(40 0)\htext(43
1){\tiny $0$}

\move(50 0)\lvec(60 0)\lvec(60 10)\lvec(50 10)\lvec(50 0)\htext(53
1){\tiny $0$}

\move(60 0)\lvec(70 0)\lvec(70 10)\lvec(60 10)\lvec(60 0)\htext(63
1){\tiny $0$}
\move(-10 0)\lvec(0 0)\lvec(0 5)\lvec(-10 5)\lvec(-10 0)\lfill
f:0.8 \htext(-7 6){\tiny $0$}

\move(0 0)\lvec(10 0)\lvec(10 5)\lvec(0 5)\lvec(0 0)\lfill f:0.8
\htext(3 6){\tiny $0$}

\move(10 0)\lvec(20 0)\lvec(20 5)\lvec(10 5)\lvec(10 0)\lfill
f:0.8 \htext(13 6){\tiny $0$}

\move(20 0)\lvec(30 0)\lvec(30 5)\lvec(20 5)\lvec(20 0)\lfill
f:0.8 \htext(23 6){\tiny $0$}

\move(30 0)\lvec(40 0)\lvec(40 5)\lvec(30 5)\lvec(30 0)\lfill
f:0.8 \htext(33 6){\tiny $0$}

\move(40 0)\lvec(50 0)\lvec(50 5)\lvec(40 5)\lvec(40 0)\lfill
f:0.8 \htext(43 6){\tiny $0$}

\move(50 0)\lvec(60 0)\lvec(60 5)\lvec(50 5)\lvec(50 0)\lfill
f:0.8 \htext(53 6){\tiny $0$}

\move(60 0)\lvec(70 0)\lvec(70 5)\lvec(60 5)\lvec(60 0)\lfill
f:0.8 \htext(63 6){\tiny $0$}
\move(0 10)\lvec(10 10)\lvec(10 20)\lvec(0 20)\lvec(0 10)\htext(3
13){$_1$}

\move(10 10)\lvec(20 10)\lvec(20 20)\lvec(10 20)\lvec(10
10)\htext(13 13){$_1$}

\move(20 10)\lvec(30 10)\lvec(30 20)\lvec(20 20)\lvec(20
10)\htext(23 13){$_1$}

\move(30 10)\lvec(40 10)\lvec(40 20)\lvec(30 20)\lvec(30
10)\htext(33 13){$_1$}

\move(40 10)\lvec(50 10)\lvec(50 20)\lvec(40 20)\lvec(40
10)\htext(43 13){$_1$}

\move(50 10)\lvec(60 10)\lvec(60 20)\lvec(50 20)\lvec(50
10)\htext(53 13){$_1$}

\move(60 10)\lvec(70 10)\lvec(70 20)\lvec(60 20)\lvec(60
10)\htext(63 13){$_1$}
\move(0 20)\lvec(10 20)\lvec(10 30)\lvec(0 30)\lvec(0 20)\htext(3
23){$_2$}

\move(10 20)\lvec(20 20)\lvec(20 30)\lvec(10 30)\lvec(10
20)\htext(13 23){$_2$}

\move(20 20)\lvec(30 20)\lvec(30 30)\lvec(20 30)\lvec(20
20)\htext(23 23){$_2$}

\move(30 20)\lvec(40 20)\lvec(40 30)\lvec(30 30)\lvec(30
20)\htext(33 23){$_2$}

\move(40 20)\lvec(50 20)\lvec(50 30)\lvec(40 30)\lvec(40
20)\htext(43 23){$_2$}

\move(50 20)\lvec(60 20)\lvec(60 30)\lvec(50 30)\lvec(50
20)\htext(53 23){$_2$}

\move(60 20)\lvec(70 20)\lvec(70 30)\lvec(60 30)\lvec(60
20)\htext(63 23){$_2$}

\move(0 30)\lvec(10 30)\lvec(10 40)\lvec(0 40)\lvec(0 30)\htext(3
33){$_1$}

\move(10 30)\lvec(20 30)\lvec(20 40)\lvec(10 40)\lvec(10
30)\htext(13 33){$_1$}

\move(20 30)\lvec(30 30)\lvec(30 40)\lvec(20 40)\lvec(20
30)\htext(23 33){$_1$}

\move(30 30)\lvec(40 30)\lvec(40 40)\lvec(30 40)\lvec(30
30)\htext(33 33){$_1$}

\move(40 30)\lvec(50 30)\lvec(50 40)\lvec(40 40)\lvec(40
30)\htext(43 33){$_1$}

\move(50 30)\lvec(60 30)\lvec(60 40)\lvec(50 40)\lvec(50
30)\htext(53 33){$_1$}

\move(60 30)\lvec(70 30)\lvec(70 40)\lvec(60 40)\lvec(60
30)\htext(63 33){$_1$}
\move(10 40)\lvec(20 40)\lvec(20 45)\lvec(10 45)\lvec(10
40)\htext(13 41){\tiny $0$}

\move(20 40)\lvec(30 40)\lvec(30 45)\lvec(20 45)\lvec(20
40)\htext(23 41){\tiny $0$}

\move(30 40)\lvec(40 40)\lvec(40 45)\lvec(30 45)\lvec(30
40)\htext(33 41){\tiny $0$}

\move(40 40)\lvec(50 40)\lvec(50 45)\lvec(40 45)\lvec(40
40)\htext(43 41){\tiny $0$}

\move(50 40)\lvec(60 40)\lvec(60 45)\lvec(50 45)\lvec(50
40)\htext(53 41){\tiny $0$}

\move(60 40)\lvec(70 40)\lvec(70 45)\lvec(60 45)\lvec(60
40)\htext(63 41){\tiny $0$}
\move(30 40)\lvec(40 40)\lvec(40 50)\lvec(30 50)\lvec(30
40)\htext(33 46){\tiny $0$}

\move(40 40)\lvec(50 40)\lvec(50 50)\lvec(40 50)\lvec(40
40)\htext(43 46){\tiny $0$}

\move(50 40)\lvec(60 40)\lvec(60 50)\lvec(50 50)\lvec(50
40)\htext(53 46){\tiny $0$}

\move(60 40)\lvec(70 40)\lvec(70 50)\lvec(60 50)\lvec(60
40)\htext(63 46){\tiny $0$}

\move(40 50)\lvec(50 50)\lvec(50 60)\lvec(40 60)\lvec(40
50)\htext(43 53){$_1$}

\move(50 50)\lvec(60 50)\lvec(60 60)\lvec(50 60)\lvec(50
50)\htext(53 53){$_1$}

\move(60 50)\lvec(70 50)\lvec(70 60)\lvec(60 60)\lvec(60
50)\htext(63 53){$_1$}
\move(50 60)\lvec(60 60)\lvec(60 70)\lvec(50 70)\lvec(50
60)\htext(53 63){$_2$}

\move(60 60)\lvec(70 60)\lvec(70 70)\lvec(60 70)\lvec(60
60)\htext(63 63){$_2$}
%

\move(60 70)\lvec(70 70)\lvec(70 80)\lvec(60 80)\lvec(60
70)\htext(63 73){$_1$}
%

\move(60 80)\lvec(70 80)\lvec(70 85)\lvec(60 85)\lvec(60
80)\htext(63 81){\tiny $0$}
\end{texdraw}}
\end{center}\vskip 5mm
then we have $Y=(Y_0,Y_1,Y_2,Y_3)$, where\vskip 5mm

\hskip 1cm$Y_0=$\raisebox{-0.5\height}{
\begin{texdraw}
\drawdim em \setunitscale 0.13 \linewd 0.5

\move(-10 0)\lvec(0 0)\lvec(0 5)\lvec(-10 5)\lvec(-10 0)\lfill
f:0.8 \htext(-7 1){\tiny $0$}

\move(0 0)\lvec(10 0)\lvec(10 5)\lvec(0 5)\lvec(0 0)\lfill f:0.8
\htext(3 1){\tiny $0$}

\move(10 0)\lvec(20 0)\lvec(20 5)\lvec(10 5)\lvec(10 0)\lfill
f:0.8 \htext(13 1){\tiny $0$}

\move(20 0)\lvec(30 0)\lvec(30 5)\lvec(20 5)\lvec(20 0)\lfill
f:0.8 \htext(23 1){\tiny $0$}

\move(30 0)\lvec(40 0)\lvec(40 5)\lvec(30 5)\lvec(30 0)\lfill
f:0.8 \htext(33 1){\tiny $0$}

\move(40 0)\lvec(50 0)\lvec(50 5)\lvec(40 5)\lvec(40 0)\lfill
f:0.8 \htext(43 1){\tiny $0$}

\move(40 5)\lvec(40 10)\lvec(50 10)\lvec(50 5)\htext(43 6){\tiny
$0$}

\htext(-23 1){$\cdots$}
\end{texdraw}} : $i$-component of type $\rm II_{\infty}$,\vskip 5mm

\hskip 1cm$Y_1=$\raisebox{-0.5\height}{
\begin{texdraw}
\drawdim em \setunitscale 0.13 \linewd 0.5

\move(10 0)\lvec(20 0)\lvec(20 10)\lvec(10 10)\lvec(10 0)\htext(13
1){\tiny $0$}

\move(20 0)\lvec(30 0)\lvec(30 10)\lvec(20 10)\lvec(20 0)\htext(23
1){\tiny $0$}

\move(30 0)\lvec(40 0)\lvec(40 10)\lvec(30 10)\lvec(30 0)\htext(33
1){\tiny $0$}

\move(40 0)\lvec(50 0)\lvec(50 10)\lvec(40 10)\lvec(40 0)\htext(43
1){\tiny $0$}

\move(10 0)\lvec(20 0)\lvec(20 5)\lvec(10 5)\lvec(10 0)\lfill
f:0.8 \htext(13 6){\tiny $0$}

\move(20 0)\lvec(30 0)\lvec(30 5)\lvec(20 5)\lvec(20 0)\lfill
f:0.8 \htext(23 6){\tiny $0$}

\move(30 0)\lvec(40 0)\lvec(40 5)\lvec(30 5)\lvec(30 0)\lfill
f:0.8 \htext(33 6){\tiny $0$}

\move(40 0)\lvec(50 0)\lvec(50 5)\lvec(40 5)\lvec(40 0)\lfill
f:0.8 \htext(43 6){\tiny $0$}

\move(10 10)\lvec(20 10)\lvec(20 20)\lvec(10 20)\lvec(10
10)\htext(13 13){$_1$}

\move(20 10)\lvec(30 10)\lvec(30 20)\lvec(20 20)\lvec(20
10)\htext(23 13){$_1$}

\move(30 10)\lvec(40 10)\lvec(40 20)\lvec(30 20)\lvec(30
10)\htext(33 13){$_1$}

\move(40 10)\lvec(50 10)\lvec(50 20)\lvec(40 20)\lvec(40
10)\htext(43 13){$_1$}

\move(10 20)\lvec(20 20)\lvec(20 30)\lvec(10 30)\lvec(10
20)\htext(13 23){$_2$}

\move(20 20)\lvec(30 20)\lvec(30 30)\lvec(20 30)\lvec(20
20)\htext(23 23){$_2$}

\move(30 20)\lvec(40 20)\lvec(40 30)\lvec(30 30)\lvec(30
20)\htext(33 23){$_2$}

\move(40 20)\lvec(50 20)\lvec(50 30)\lvec(40 30)\lvec(40
20)\htext(43 23){$_2$}

\move(10 30)\lvec(20 30)\lvec(20 40)\lvec(10 40)\lvec(10
30)\htext(13 33){$_1$}

\move(20 30)\lvec(30 30)\lvec(30 40)\lvec(20 40)\lvec(20
30)\htext(23 33){$_1$}

\move(30 30)\lvec(40 30)\lvec(40 40)\lvec(30 40)\lvec(30
30)\htext(33 33){$_1$}

\move(40 30)\lvec(50 30)\lvec(50 40)\lvec(40 40)\lvec(40
30)\htext(43 33){$_1$}

\move(20 40)\lvec(30 40)\lvec(30 45)\lvec(20 45)\lvec(20
40)\htext(23 41){\tiny $0$}

\move(30 40)\lvec(40 40)\lvec(40 45)\lvec(30 45)\lvec(30
40)\htext(33 41){\tiny $0$}

\move(40 40)\lvec(50 40)\lvec(50 45)\lvec(40 45)\lvec(40
40)\htext(43 41){\tiny $0$}

\move(40 40)\lvec(50 40)\lvec(50 50)\lvec(40 50)\lvec(40
40)\htext(43 46){\tiny $0$}
\end{texdraw}} : $i$-component of type $\rm II_3$,\vskip 5mm

\hskip 1cm$Y_2=$\raisebox{-0.5\height}{
\begin{texdraw}
\drawdim em \setunitscale 0.13 \linewd 0.5

\move(40 0)\lvec(50 0)\lvec(50 10)\lvec(40 10)\lvec(40 0)\htext(43
1){\tiny $0$}

\move(50 0)\lvec(60 0)\lvec(60 10)\lvec(50 10)\lvec(50 0)\htext(53
1){\tiny $0$}

\move(40 0)\lvec(50 0)\lvec(50 5)\lvec(40 5)\lvec(40 0)\lfill
f:0.8 \htext(43 6){\tiny $0$}

\move(50 0)\lvec(60 0)\lvec(60 5)\lvec(50 5)\lvec(50 0)\lfill
f:0.8 \htext(53 6){\tiny $0$}

\move(40 10)\lvec(50 10)\lvec(50 20)\lvec(40 20)\lvec(40
10)\htext(43 13){$_1$}

\move(50 10)\lvec(60 10)\lvec(60 20)\lvec(50 20)\lvec(50
10)\htext(53 13){$_1$}

\move(40 20)\lvec(50 20)\lvec(50 30)\lvec(40 30)\lvec(40
20)\htext(43 23){$_2$}

\move(50 20)\lvec(60 20)\lvec(60 30)\lvec(50 30)\lvec(50
20)\htext(53 23){$_2$}

\move(40 30)\lvec(50 30)\lvec(50 40)\lvec(40 40)\lvec(40
30)\htext(43 33){$_1$}

\move(50 30)\lvec(60 30)\lvec(60 40)\lvec(50 40)\lvec(50
30)\htext(53 33){$_1$}

\move(40 40)\lvec(50 40)\lvec(50 45)\lvec(40 45)\lvec(40
40)\htext(43 41){\tiny $0$}

\move(50 40)\lvec(60 40)\lvec(60 45)\lvec(50 45)\lvec(50
40)\htext(53 41){\tiny $0$}

\move(40 40)\lvec(50 40)\lvec(50 50)\lvec(40 50)\lvec(40
40)\htext(43 46){\tiny $0$}

\move(50 40)\lvec(60 40)\lvec(60 50)\lvec(50 50)\lvec(50
40)\htext(53 46){\tiny $0$}

\move(40 50)\lvec(50 50)\lvec(50 60)\lvec(40 60)\lvec(40
50)\htext(43 53){$_1$}

\move(50 50)\lvec(60 50)\lvec(60 60)\lvec(50 60)\lvec(50
50)\htext(53 53){$_1$}

\move(50 60)\lvec(60 60)\lvec(60 70)\lvec(50 70)\lvec(50
60)\htext(53 63){$_2$}

\end{texdraw}} : trivial $i$-component and \vskip 5mm

\hskip 1cm$Y_3=$\raisebox{-0.5\height}{
\begin{texdraw}
\drawdim em \setunitscale 0.13 \linewd 0.5

\move(60 0)\lvec(70 0)\lvec(70 10)\lvec(60 10)\lvec(60 0)\htext(63
1){\tiny $0$}

\move(60 0)\lvec(70 0)\lvec(70 5)\lvec(60 5)\lvec(60 0)\lfill
f:0.8 \htext(63 6){\tiny $0$}

\move(60 10)\lvec(70 10)\lvec(70 20)\lvec(60 20)\lvec(60
10)\htext(63 13){$_1$}

\move(60 20)\lvec(70 20)\lvec(70 30)\lvec(60 30)\lvec(60
20)\htext(63 23){$_2$}

\move(60 30)\lvec(70 30)\lvec(70 40)\lvec(60 40)\lvec(60
30)\htext(63 33){$_1$}

\move(60 40)\lvec(70 40)\lvec(70 45)\lvec(60 45)\lvec(60
40)\htext(63 41){\tiny $0$}

\move(60 40)\lvec(70 40)\lvec(70 50)\lvec(60 50)\lvec(60
40)\htext(63 46){\tiny $0$}
\move(60 50)\lvec(70 50)\lvec(70 60)\lvec(60 60)\lvec(60
50)\htext(63 53){$_1$}

\move(60 60)\lvec(70 60)\lvec(70 70)\lvec(60 70)\lvec(60
60)\htext(63 63){$_2$}

\move(60 70)\lvec(70 70)\lvec(70 80)\lvec(60 80)\lvec(60
70)\htext(63 73){$_1$}

\move(60 80)\lvec(70 80)\lvec(70 85)\lvec(60 85)\lvec(60
80)\htext(63 81){\tiny $0$}

\end{texdraw}} : $i$-component of type $\rm II_0$.

}
\end{ex}

Let $Y$ be a proper Young wall with the decomposition $Y=(Y_0,Y_1,
\cdots,Y_r)$ into $i$-components. To each $Y_k$, we associate a
$U_{(i)}$-module $V_k$ as follows. Then we will view $Y$ as
$Y_0\otimes Y_1 \otimes \cdots \otimes Y_r$ inside $V_0\otimes V_1
\otimes \cdots \otimes V_r$.

If $Y_k$ is a trivial $i$-component, then we associate the trivial
representation $V_k=U=\mathbb{Q}(q)u$, and we identify $Y_k$ with
$u$.

If $Y_k$ is an $i$-component of type $\rm II_{\infty}$, then we
must have $k=0$ and we associate the 2-dimensional representation
$V_0=V=\mathbb{Q}(q)v_0\oplus \mathbb{Q}(q)v_1$, where
$U_{(i)}$-module action is given by \eqref{I+}. We identify $Y_0$
with the basis element of $V$ as follows :\vskip 5mm

\begin{center}
$v_0 \leftrightarrow$ \raisebox{-0.5\height}{
\begin{texdraw}
\drawdim em \setunitscale 0.12 \linewd 0.5

\move(-10 0)\lvec(0 0)\lvec(0 5)\lvec(-10 5)\lvec(-10 0)\lfill
f:0.8 \htext(-7 1){\tiny $i$}

\move(0 0)\lvec(10 0)\lvec(10 5)\lvec(0 5)\lvec(0 0)\lfill f:0.8
\htext(3 1){\tiny $i$}

\move(10 0)\lvec(20 0)\lvec(20 5)\lvec(10 5)\lvec(10 0)\lfill
f:0.8 \htext(13 1){\tiny $i$}

\move(20 0)\lvec(30 0)\lvec(30 5)\lvec(20 5)\lvec(20 0)\lfill
f:0.8 \htext(23 1){\tiny $i$}

\move(30 0)\lvec(40 0)\lvec(40 5)\lvec(30 5)\lvec(30 0)\lfill
f:0.8 \htext(33 1){\tiny $i$}

\move(40 0)\lvec(50 0)\lvec(50 5)\lvec(40 5)\lvec(40 0)\lfill
f:0.8 \htext(43 1){\tiny $i$}

\htext(43 6){\tiny $i$}

\htext(-23 1){$\cdots$}
\end{texdraw}}\hskip 2cm
$v_1 \leftrightarrow$\raisebox{-0.5\height}{
\begin{texdraw}
\drawdim em \setunitscale 0.12 \linewd 0.5

\move(-10 0)\lvec(0 0)\lvec(0 5)\lvec(-10 5)\lvec(-10 0)\lfill
f:0.8 \htext(-7 1){\tiny $i$}

\move(0 0)\lvec(10 0)\lvec(10 5)\lvec(0 5)\lvec(0 0)\lfill f:0.8
\htext(3 1){\tiny $i$}

\move(10 0)\lvec(20 0)\lvec(20 5)\lvec(10 5)\lvec(10 0)\lfill
f:0.8 \htext(13 1){\tiny $i$}

\move(20 0)\lvec(30 0)\lvec(30 5)\lvec(20 5)\lvec(20 0)\lfill
f:0.8 \htext(23 1){\tiny $i$}

\move(30 0)\lvec(40 0)\lvec(40 5)\lvec(30 5)\lvec(30 0)\lfill
f:0.8 \htext(33 1){\tiny $i$}

\move(40 0)\lvec(50 0)\lvec(50 5)\lvec(40 5)\lvec(40 0)\lfill
f:0.8 \htext(43 1){\tiny $i$}

\move(40 5)\lvec(40 10)\lvec(50 10)\lvec(50 5)\htext(43 6){\tiny
$i$}

\htext(-23 1){$\cdots$}
\end{texdraw}}\quad.
\end{center}
\vskip 5mm

If $Y_k$ is an $i$-component of type $\rm II_0$, then we associate
the 3-dimensional representation $V_k=W_0=\mathbb{Q}(q)w_0\oplus
\mathbb{Q}(q)w_1\oplus \mathbb{Q}(q)w_2$, where the
$U_{(i)}$-module action is given by
\begin{equation}\label{II+}
\begin{split}
& K_iw_0=q_i^2w_0,\ K_iw_1=w_1, \ K_iw_2=q_i^{-2}w_2, \\
& e_iw_0=0, \ e_iw_1=(q_i+q_i^{-1})w_0, \ e_iw_2=w_1, \\
& f_iw_0=w_1, \ f_iw_1=(q_i+q_i^{-1})w_2, \ f_iw_2=0.
\end{split}
\end{equation}
We identify the $i$-component $Y_k$ with a basis element of $W$ as
follows :\vskip 5mm

\begin{center}
$w_0 \leftrightarrow$\raisebox{-0.6\height}{
\begin{texdraw}
\drawdim em \setunitscale 0.12 \linewd 0.5 \move(0 -10)\lvec(0
10)\lvec(10 10)\lvec(10 -10)

\end{texdraw}}\hskip 10mm
$w_1 \leftrightarrow$\raisebox{-0.5\height}{
\begin{texdraw}
\drawdim em \setunitscale 0.12 \linewd 0.5 \move(0 -10)\lvec(0
10)\lvec(10 10)\lvec(10 -10)

\move(0 10)\lvec(0 15)\lvec(10 15)\lvec(10 10)\lvec(0 10)

\htext(4 11){\tiny $i$}
\end{texdraw}}\hskip 10mm
$w_2 \leftrightarrow$\raisebox{-0.4\height}{
\begin{texdraw}
\drawdim em \setunitscale 0.12 \linewd 0.5 \move(0 -10)\lvec(0
10)\lvec(10 10)\lvec(10 -10)

\move(0 10)\lvec(0 15)\lvec(10 15)\lvec(10 10)\lvec(0 10)

\move(0 15)\lvec(0 20)\lvec(10 20)\lvec(10 15)\lvec(0 15)

\htext(4 11){\tiny $i$}\htext(4 16){\tiny $i$}
\end{texdraw}}\quad .
\end{center}\vskip 5mm

Note that $W$ is isomorphic to the 3-dimensional irreducible
$U_{(i)}$-module $V(2)$ with the crystal basis $(L,B)$, where
\begin{equation*}
L =\mathbb{A}_0w_0\oplus \mathbb{A}_0w_1\oplus\mathbb{A}_0w_2, \
B=\{\,\overline{w_k}\,|\,k=0,1,2\,\}.
\end{equation*}
and the crystal graph is given by
$$\overline{w_0}\stackrel{i}{\longrightarrow}\overline{w_1}
\stackrel{i}{\longrightarrow}\overline{w_2}.$$

If $Y_k$ is an $i$-component of type ${\rm II}_l$ ($l\geq 1$),
then we associate the 4-dimensional representation
\begin{equation*}
V_k=W_l=\mathbb{Q}(q)w_0\oplus\mathbb{Q}(q)w_1\oplus\mathbb{Q}(q)w_2\oplus\mathbb{Q}(q)u,
\end{equation*}
where the $U_{(i)}$-module action is given by
\begin{equation}
\begin{split}
& K_iw_0=q_i^2w_0, \ K_iw_1=w_1, \  K_iw_2=q_i^{-2}w_2 , \ K_iu=u, \\
& e_iw_0= 0, \ e_iw_1=q_i^{-1}(1-(-q_i^2)^{l+1})w_0, \\
& e_iw_2=w_1+q_i(1-(-q_i^2)^l)u, \ e_iu=w_0, \\
& f_iw_0=w_1+q_i(1-(-q_i^2)^l)u, \
f_iw_1=q_i^{-1}(1-(-q_i^2)^{l+1})w_2, \\
& f_iw_2=0, \ f_iu=w_2.
\end{split}
\end{equation}
We identify the $i$-component $Y_k$ with a basis element of $W_l$
as follows :\vskip 5mm

\begin{center}
$w_0 \leftrightarrow$\raisebox{-0.6\height}{
\begin{texdraw}
\drawdim em \setunitscale 0.1 \linewd 0.5 \move(0 0)\lvec(0
10)\lvec(50 10)\lvec(50 0)

\move(10 10)\lvec(10 15)\lvec(50 15)\lvec(50 10)\lvec(10 10)

\move(20 10) \lvec(20 15) \move(40 10) \lvec(40 15) \htext(25
10){$\cdots$}

\move(10 15)\clvec(12 20)(23 20)(25 20)

\move(50 15)\clvec(48 20)(37 20)(35 20)

\htext(30 20){$_{l}$}
\end{texdraw}}\hskip 10mm
$w_1\leftrightarrow$\raisebox{-0.7\height}{
\begin{texdraw}
\drawdim em \setunitscale 0.1 \linewd 0.5 \move(0 0)\lvec(0
10)\lvec(50 10)\lvec(50 0)

\move(0 10) \lvec(0 15) \lvec(10 15)\lvec(10 10)\lvec(0 10)

\move(10 10)\lvec(10 15)\lvec(50 15)\lvec(50 10)\lvec(10 10)

\move(20 10) \lvec(20 15) \move(40 10) \lvec(40 15) \htext(25
10){$\cdots$}

\end{texdraw}}\vskip 10mm

$w_2 \leftrightarrow$\raisebox{-0.7\height}{
\begin{texdraw}
\drawdim em \setunitscale 0.1 \linewd 0.5 \move(0 0)\lvec(0
10)\lvec(50 10)\lvec(50 0)

\move(0 10) \lvec(0 15) \lvec(10 15)\lvec(10 10)\lvec(0 10)

\move(40 15) \lvec(40 20) \lvec(50 20) \lvec(50 15)

\move(10 10)\lvec(10 15)\lvec(50 15)\lvec(50 10)\lvec(10 10)

\move(20 10) \lvec(20 15) \move(40 10) \lvec(40 15) \htext(25
10){$\cdots$}

\end{texdraw}}
\hskip 10mm $u \leftrightarrow$\raisebox{-0.7\height}{
\begin{texdraw}
\drawdim em \setunitscale 0.1 \linewd 0.5 \move(0 0)\lvec(0
10)\lvec(50 10)\lvec(50 0)

\move(40 15) \lvec(40 20) \lvec(50 20) \lvec(50 15)

\move(10 10)\lvec(10 15)\lvec(50 15)\lvec(50 10)\lvec(10 10)

\move(20 10) \lvec(20 15) \move(40 10) \lvec(40 15) \htext(25
10){$\cdots$}
\end{texdraw}}\quad .\vskip 5mm
\end{center}

The $U_{(i)}$-module $W_l$ is decomposed as $W_l\cong V(2)\oplus
V(0)$, where
\begin{equation*}
\begin{split}
V(2)&\cong \mathbb{Q}(q)w_0\oplus
\mathbb{Q}(q)(w_1+q_i(1-(-q_i^2)^{l})u)\oplus\mathbb{Q}(q)w_2,\\
V(0)&\cong \mathbb{Q}(q)(u-\frac{q_i}{(1-(-q_i^2)^{l+1})}w_1).
\end{split}
\end{equation*}
Hence, the crystal basis $(L,B)$ of $W_l$ is given by
\begin{equation*}
\begin{split}
L=&\mathbb{A}_0w_0\oplus
\mathbb{A}_0(w_1+q_i(1-(-q_i^2)^{l})u)\oplus\mathbb{A}_0w_2
\\ & \oplus
\mathbb{A}_0(u-\frac{q_i}{(1-(-q_i^2)^{l+1})}w_1), \\
B=&\{\,\overline{w_0},\overline{w_1},\overline{w_2},\overline{u}\,\}
\end{split}
\end{equation*}
with the crystal graph\vskip 5mm
\begin{center}
\raisebox{-0.5\height}{\begin{texdraw} \drawdim em \setunitscale
0.12 \linewd 0.3

\htext(0 0){$\overline{u}$} \htext(0
30){$\overline{w_0}$}\htext(30 30){$\overline{w_1}$}\htext(30
0){$\bar{w_2}$}

\move(15 33) \arrowheadtype t:V \arrowheadsize l:2 w:4 \avec(25
33)

\move(35 23) \arrowheadtype t:V \arrowheadsize l:2 w:4 \avec(35
12)

\htext(19 35){\tiny $i$} \htext(37 17.5){\tiny $i$}
\end{texdraw}}\quad .
\end{center}\vskip 5mm

\begin{ex}{\rm
In Example \ref{comp2}, we have \vskip 5mm

\begin{center}
$Y=$\hskip 2mm \raisebox{-0.5\height}{
\begin{texdraw}
\drawdim em \setunitscale 0.13 \linewd 0.5

\move(-10 0)\lvec(0 0)\lvec(0 5)\lvec(-10 5)\lvec(-10 0)\lfill
f:0.8 \htext(-7 1){\tiny $0$}

\move(0 0)\lvec(10 0)\lvec(10 5)\lvec(0 5)\lvec(0 0)\lfill f:0.8
\htext(3 1){\tiny $0$}

\move(10 0)\lvec(20 0)\lvec(20 5)\lvec(10 5)\lvec(10 0)\lfill
f:0.8 \htext(13 1){\tiny $0$}

\move(20 0)\lvec(30 0)\lvec(30 5)\lvec(20 5)\lvec(20 0)\lfill
f:0.8 \htext(23 1){\tiny $0$}

\move(30 0)\lvec(40 0)\lvec(40 5)\lvec(30 5)\lvec(30 0)\lfill
f:0.8 \htext(33 1){\tiny $0$}

\move(40 0)\lvec(50 0)\lvec(50 5)\lvec(40 5)\lvec(40 0)\lfill
f:0.8 \htext(43 1){\tiny $0$}

\move(40 5)\lvec(40 10)\lvec(50 10)\lvec(50 5)\htext(43 6){\tiny
$0$}

\htext(-23 1){$\cdots$}
\end{texdraw}}\hskip 2mm$\otimes$
\raisebox{-0.2\height}{
\begin{texdraw}
\drawdim em \setunitscale 0.13 \linewd 0.5

\move(10 0)\lvec(20 0)\lvec(20 10)\lvec(10 10)\lvec(10 0)\htext(13
1){\tiny $0$}

\move(20 0)\lvec(30 0)\lvec(30 10)\lvec(20 10)\lvec(20 0)\htext(23
1){\tiny $0$}

\move(30 0)\lvec(40 0)\lvec(40 10)\lvec(30 10)\lvec(30 0)\htext(33
1){\tiny $0$}

\move(40 0)\lvec(50 0)\lvec(50 10)\lvec(40 10)\lvec(40 0)\htext(43
1){\tiny $0$}

\move(10 0)\lvec(20 0)\lvec(20 5)\lvec(10 5)\lvec(10 0)\lfill
f:0.8 \htext(13 6){\tiny $0$}

\move(20 0)\lvec(30 0)\lvec(30 5)\lvec(20 5)\lvec(20 0)\lfill
f:0.8 \htext(23 6){\tiny $0$}

\move(30 0)\lvec(40 0)\lvec(40 5)\lvec(30 5)\lvec(30 0)\lfill
f:0.8 \htext(33 6){\tiny $0$}

\move(40 0)\lvec(50 0)\lvec(50 5)\lvec(40 5)\lvec(40 0)\lfill
f:0.8 \htext(43 6){\tiny $0$}

\move(10 10)\lvec(20 10)\lvec(20 20)\lvec(10 20)\lvec(10
10)\htext(13 13){$_1$}

\move(20 10)\lvec(30 10)\lvec(30 20)\lvec(20 20)\lvec(20
10)\htext(23 13){$_1$}

\move(30 10)\lvec(40 10)\lvec(40 20)\lvec(30 20)\lvec(30
10)\htext(33 13){$_1$}

\move(40 10)\lvec(50 10)\lvec(50 20)\lvec(40 20)\lvec(40
10)\htext(43 13){$_1$}

\move(10 20)\lvec(20 20)\lvec(20 30)\lvec(10 30)\lvec(10
20)\htext(13 23){$_2$}

\move(20 20)\lvec(30 20)\lvec(30 30)\lvec(20 30)\lvec(20
20)\htext(23 23){$_2$}

\move(30 20)\lvec(40 20)\lvec(40 30)\lvec(30 30)\lvec(30
20)\htext(33 23){$_2$}

\move(40 20)\lvec(50 20)\lvec(50 30)\lvec(40 30)\lvec(40
20)\htext(43 23){$_2$}

\move(10 30)\lvec(20 30)\lvec(20 40)\lvec(10 40)\lvec(10
30)\htext(13 33){$_1$}

\move(20 30)\lvec(30 30)\lvec(30 40)\lvec(20 40)\lvec(20
30)\htext(23 33){$_1$}

\move(30 30)\lvec(40 30)\lvec(40 40)\lvec(30 40)\lvec(30
30)\htext(33 33){$_1$}

\move(40 30)\lvec(50 30)\lvec(50 40)\lvec(40 40)\lvec(40
30)\htext(43 33){$_1$}

\move(20 40)\lvec(30 40)\lvec(30 45)\lvec(20 45)\lvec(20
40)\htext(23 41){\tiny $0$}

\move(30 40)\lvec(40 40)\lvec(40 45)\lvec(30 45)\lvec(30
40)\htext(33 41){\tiny $0$}

\move(40 40)\lvec(50 40)\lvec(50 45)\lvec(40 45)\lvec(40
40)\htext(43 41){\tiny $0$}

\move(40 40)\lvec(50 40)\lvec(50 50)\lvec(40 50)\lvec(40
40)\htext(43 46){\tiny $0$}
\end{texdraw}}\hskip 2mm$\otimes$
\raisebox{-0.2\height}{
\begin{texdraw}
\drawdim em \setunitscale 0.13 \linewd 0.5

\move(40 0)\lvec(50 0)\lvec(50 10)\lvec(40 10)\lvec(40 0)\htext(43
1){\tiny $0$}

\move(50 0)\lvec(60 0)\lvec(60 10)\lvec(50 10)\lvec(50 0)\htext(53
1){\tiny $0$}

\move(40 0)\lvec(50 0)\lvec(50 5)\lvec(40 5)\lvec(40 0)\lfill
f:0.8 \htext(43 6){\tiny $0$}

\move(50 0)\lvec(60 0)\lvec(60 5)\lvec(50 5)\lvec(50 0)\lfill
f:0.8 \htext(53 6){\tiny $0$}

\move(40 10)\lvec(50 10)\lvec(50 20)\lvec(40 20)\lvec(40
10)\htext(43 13){$_1$}

\move(50 10)\lvec(60 10)\lvec(60 20)\lvec(50 20)\lvec(50
10)\htext(53 13){$_1$}

\move(40 20)\lvec(50 20)\lvec(50 30)\lvec(40 30)\lvec(40
20)\htext(43 23){$_2$}

\move(50 20)\lvec(60 20)\lvec(60 30)\lvec(50 30)\lvec(50
20)\htext(53 23){$_2$}

\move(40 30)\lvec(50 30)\lvec(50 40)\lvec(40 40)\lvec(40
30)\htext(43 33){$_1$}

\move(50 30)\lvec(60 30)\lvec(60 40)\lvec(50 40)\lvec(50
30)\htext(53 33){$_1$}

\move(40 40)\lvec(50 40)\lvec(50 45)\lvec(40 45)\lvec(40
40)\htext(43 41){\tiny $0$}

\move(50 40)\lvec(60 40)\lvec(60 45)\lvec(50 45)\lvec(50
40)\htext(53 41){\tiny $0$}

\move(40 40)\lvec(50 40)\lvec(50 50)\lvec(40 50)\lvec(40
40)\htext(43 46){\tiny $0$}

\move(50 40)\lvec(60 40)\lvec(60 50)\lvec(50 50)\lvec(50
40)\htext(53 46){\tiny $0$}

\move(40 50)\lvec(50 50)\lvec(50 60)\lvec(40 60)\lvec(40
50)\htext(43 53){$_1$}

\move(50 50)\lvec(60 50)\lvec(60 60)\lvec(50 60)\lvec(50
50)\htext(53 53){$_1$}

\move(50 60)\lvec(60 60)\lvec(60 70)\lvec(50 70)\lvec(50
60)\htext(53 63){$_2$}

\end{texdraw}}\hskip 2mm$\otimes$
\raisebox{-0.2\height}{
\begin{texdraw}
\drawdim em \setunitscale 0.13 \linewd 0.5

\move(60 0)\lvec(70 0)\lvec(70 10)\lvec(60 10)\lvec(60 0)\htext(63
1){\tiny $0$}

\move(60 0)\lvec(70 0)\lvec(70 5)\lvec(60 5)\lvec(60 0)\lfill
f:0.8 \htext(63 6){\tiny $0$}

\move(60 10)\lvec(70 10)\lvec(70 20)\lvec(60 20)\lvec(60
10)\htext(63 13){$_1$}

\move(60 20)\lvec(70 20)\lvec(70 30)\lvec(60 30)\lvec(60
20)\htext(63 23){$_2$}

\move(60 30)\lvec(70 30)\lvec(70 40)\lvec(60 40)\lvec(60
30)\htext(63 33){$_1$}

\move(60 40)\lvec(70 40)\lvec(70 45)\lvec(60 45)\lvec(60
40)\htext(63 41){\tiny $0$}

\move(60 40)\lvec(70 40)\lvec(70 50)\lvec(60 50)\lvec(60
40)\htext(63 46){\tiny $0$}
\move(60 50)\lvec(70 50)\lvec(70 60)\lvec(60 60)\lvec(60
50)\htext(63 53){$_1$}

\move(60 60)\lvec(70 60)\lvec(70 70)\lvec(60 70)\lvec(60
60)\htext(63 63){$_2$}

\move(60 70)\lvec(70 70)\lvec(70 80)\lvec(60 80)\lvec(60
70)\htext(63 73){$_1$}

\move(60 80)\lvec(70 80)\lvec(70 85)\lvec(60 85)\lvec(60
80)\htext(63 81){\tiny $0$}

\end{texdraw}}\vskip 5mm

\hskip -2cm$\in V\otimes W_3\otimes U\otimes W_0$.
\end{center}\vskip 5mm

}
\end{ex}

{\bf Case 3.} $\Box={\rm III}$

Suppose that the $i$-blocks are of type III. For convenience, we
denote by $*$ the color $i$, and denote by $\cdot$ the color of
the blocks placed on the opposite side of the $i$-blocks. In this
case, the (virtually) admissible $i$-slots and the (virtually)
removable $i$-blocks in $Y$ appear in one of the following
situations: \vskip 5mm

($\rm III_{\infty}$)\hskip 1cm
\raisebox{-0.4\height}{\begin{texdraw}\drawdim em \setunitscale
0.1 \linewd 0.5

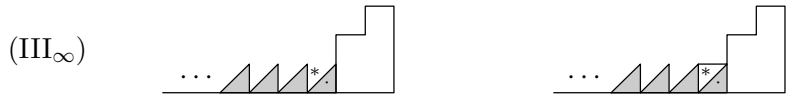
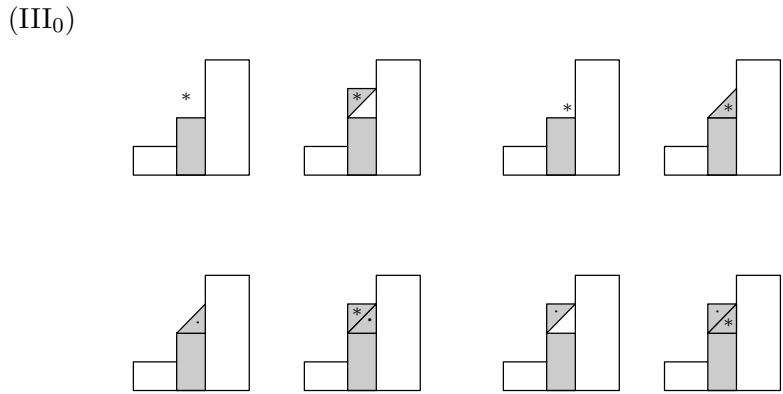
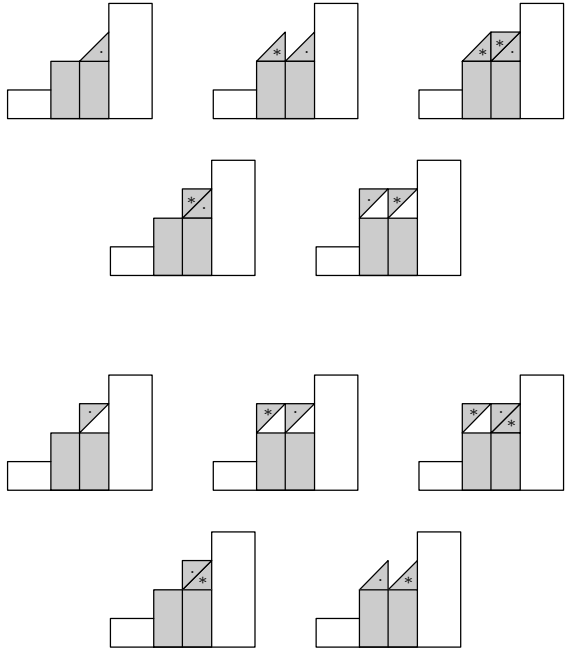
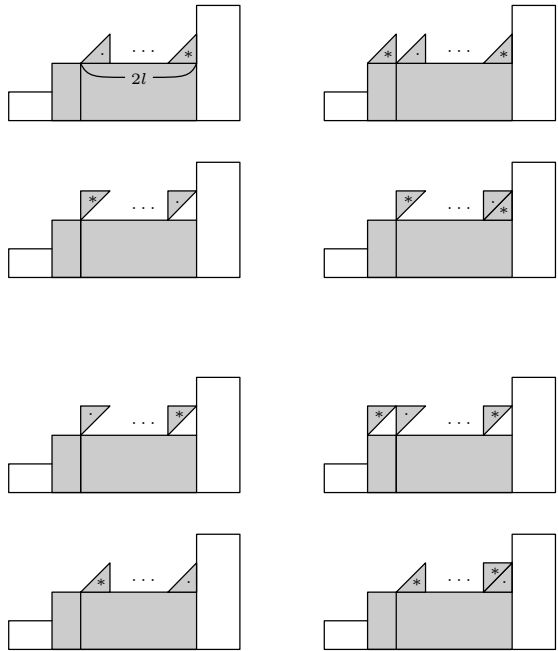
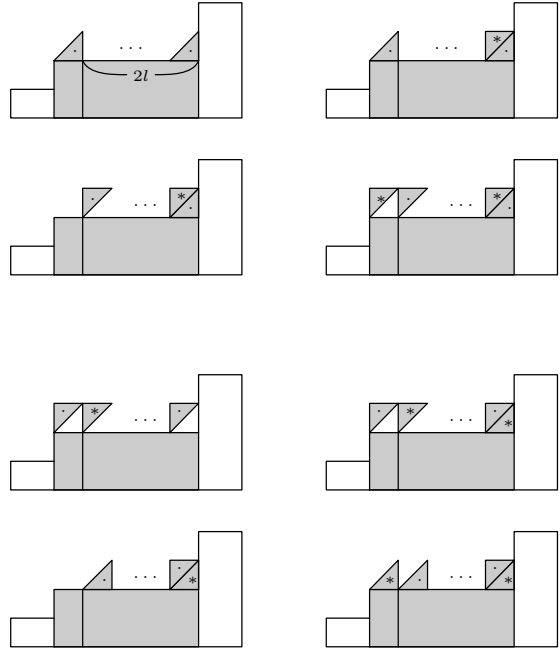
\htext(-4 3){$\cdots$} \move(10 0)\lvec(20 10)\lvec(20 0)\lvec(10
0)\lfill f:0.8

\move(20 0)\lvec(30 10)\lvec(30 0)\lvec(20 0)\lfill f:0.8

\move(30 0)\lvec(40 10)\lvec(40 0)\lvec(30 0)\lfill f:0.8

\move(40 0)\lvec(50 10)\lvec(50 0)\lvec(40 0)\lfill f:0.8

\htext(46 2){\tiny $\cdot$}\htext(41 6){\tiny $*$}

\move(-10 0)\lvec(70 0)\lvec(70 30)\lvec(60 30)\lvec(60
20)\lvec(50 20)\lvec(50 0)
\end{texdraw}}
\hskip 2cm \raisebox{-0.4\height}{\begin{texdraw}\drawdim em
\setunitscale 0.1 \linewd 0.5

\htext(-5 3){$\cdots$} \move(10 0)\lvec(20 10)\lvec(20 0)\lvec(10
0)\lfill f:0.8

\move(20 0)\lvec(30 10)\lvec(30 0)\lvec(20 0)\lfill f:0.8

\move(30 0)\lvec(40 10)\lvec(40 0)\lvec(30 0)\lfill f:0.8

\move(40 0)\lvec(50 10)\lvec(50 0)\lvec(40 0)\lfill f:0.8

\htext(46 2){\tiny $\cdot$}\htext(41 6){\tiny $*$}

\move(-10 0)\lvec(70 0)\lvec(70 30)\lvec(60 30)\lvec(60
20)\lvec(50 20)\lvec(50 10)\lvec(40 10)
\end{texdraw}}\vskip 10mm

($\rm III_0$)
\begin{center}

\raisebox{-0.6\height}{
\begin{texdraw}
\drawdim em \setunitscale 0.1 \linewd 0.5

\move(0 -10)\lvec(0 10)\lvec(10 10)\lvec(10 -10)\lvec(0 -10)\lfill
f:0.8

\move(0 0)\lvec(-15 0)\lvec(-15 -10)\lvec(25 -10)\lvec(25
30)\lvec(10 30)\lvec(10 10)

\htext(1 16){\tiny $*$}
\end{texdraw}}
\hskip 5mm \raisebox{-0.6\height}{
\begin{texdraw}
\drawdim em \setunitscale 0.1 \linewd 0.5

\move(0 -10)\lvec(0 10)\lvec(10 10)\lvec(10 -10)\lvec(0 -10)\lfill
f:0.8

\move(0 0)\lvec(-15 0)\lvec(-15 -10)\lvec(25 -10)\lvec(25
30)\lvec(10 30)\lvec(10 10)

\move(0 10)\lvec(0 20)\lvec(10 20)\lvec(0 10)\lfill f:0.8

\htext(1 16){\tiny $*$}
\end{texdraw}}\hskip 10mm \raisebox{-0.6\height}{
\begin{texdraw}
\drawdim em \setunitscale 0.1 \linewd 0.5

\move(0 -10)\lvec(0 10)\lvec(10 10)\lvec(10 -10)\lvec(0 -10)\lfill
f:0.8

\move(0 0)\lvec(-15 0)\lvec(-15 -10)\lvec(25 -10)\lvec(25
30)\lvec(10 30)\lvec(10 10)

\htext(5 12){\tiny $*$}
\end{texdraw}}\hskip 5mm \raisebox{-0.6\height}{
\begin{texdraw}
\drawdim em \setunitscale 0.1 \linewd 0.5

\move(0 -10)\lvec(0 10)\lvec(10 10)\lvec(10 -10)\lvec(0 -10)\lfill
f:0.8

\move(0 0)\lvec(-15 0)\lvec(-15 -10)\lvec(25 -10)\lvec(25
30)\lvec(10 30)\lvec(10 10)

\move(0 10)\lvec(10 20)\lvec(10 10)\lvec(0 10)\lfill f:0.8

\htext(5 12){\tiny $*$}
\end{texdraw}}
\end{center}\vskip 10mm

\begin{center}
\raisebox{-0.6\height}{
\begin{texdraw}
\drawdim em \setunitscale 0.1 \linewd 0.5

\move(0 -10)\lvec(0 10)\lvec(10 10)\lvec(10 -10)\lvec(0 -10)\lfill
f:0.8

\move(0 0)\lvec(-15 0)\lvec(-15 -10)\lvec(25 -10)\lvec(25
30)\lvec(10 30)\lvec(10 10)

\move(0 10)\lvec(10 20)\lvec(10 10)\lvec(0 10)\lfill f:0.8

\htext(6 12){\tiny $\cdot$}
\end{texdraw}}
\hskip 5mm \raisebox{-0.6\height}{
\begin{texdraw}
\drawdim em \setunitscale 0.1 \linewd 0.5

\move(0 -10)\lvec(0 10)\lvec(10 10)\lvec(10 -10)\lvec(0 -10)\lfill
f:0.8

\move(0 0)\lvec(-15 0)\lvec(-15 -10)\lvec(25 -10)\lvec(25
30)\lvec(10 30)\lvec(10 10)

\move(0 10)\lvec(0 20)\lvec(10 20)\lvec(0 10)\lfill f:0.8

\move(0 10)\lvec(10 20)\lvec(10 10)\lvec(0 10)\lfill f:0.8

\htext(1 16){\tiny $*$}\htext(6 12){$\cdot$}
\end{texdraw}}\hskip 10mm \raisebox{-0.6\height}{
\begin{texdraw}
\drawdim em \setunitscale 0.1 \linewd 0.5

\move(0 -10)\lvec(0 10)\lvec(10 10)\lvec(10 -10)\lvec(0 -10)\lfill
f:0.8

\move(0 0)\lvec(-15 0)\lvec(-15 -10)\lvec(25 -10)\lvec(25
30)\lvec(10 30)\lvec(10 10)

\move(0 10)\lvec(0 20)\lvec(10 20)\lvec(0 10)\lfill f:0.8

\htext(2 16){\tiny $\cdot$}
\end{texdraw}}\hskip 5mm \raisebox{-0.6\height}{
\begin{texdraw}
\drawdim em \setunitscale 0.1 \linewd 0.5

\move(0 -10)\lvec(0 10)\lvec(10 10)\lvec(10 -10)\lvec(0 -10)\lfill
f:0.8

\move(0 0)\lvec(-15 0)\lvec(-15 -10)\lvec(25 -10)\lvec(25
30)\lvec(10 30)\lvec(10 10)

\move(0 10)\lvec(10 20)\lvec(10 10)\lvec(0 10)\lfill f:0.8

\move(0 10)\lvec(0 20)\lvec(10 20)\lvec(0 10)\lfill f:0.8

\htext(2 16){\tiny $\cdot$}

\htext(5 12){\tiny $*$}
\end{texdraw}}
\end{center}\vskip 10mm

($\rm III_1$)

\begin{center}
\raisebox{-0.4\height}{
\begin{texdraw}
\fontsize{6}{6}\drawdim em \setunitscale 0.1 \linewd 0.5

\move(0 -10)\lvec(0 10)\lvec(20 10)\lvec(20 -10)\lvec(0 -10)\lfill
f:0.8

\move(10 -10)\lvec(10 10)

\move(10 10)\lvec(20 20)\lvec(20 10)\lvec(10 10)\lfill f:0.8

\htext(16 12){$\cdot$}

\move(0 0)\lvec(-15 0)\lvec(-15 -10)\lvec(35 -10)\lvec(35
30)\lvec(20 30)\lvec(20 20)

\end{texdraw}}\hskip 7mm
\raisebox{-0.4\height}{
\begin{texdraw}
\fontsize{6}{6}\drawdim em \setunitscale 0.1 \linewd 0.5

\move(10 10)\lvec(20 20)\lvec(20 10)\lvec(10 10)\lfill f:0.8

\htext(16 12){$\cdot$}

\move(0 10)\lvec(10 20)\lvec(10 10)\lvec(0 10)\lfill f:0.8

\htext(5 12){$*$}

\move(0 -10)\lvec(0 10)\lvec(20 10)\lvec(20 -10)\lvec(0 -10)\lfill
f:0.8

\move(10 -10)\lvec(10 10)

\move(0 0)\lvec(-15 0)\lvec(-15 -10)\lvec(35 -10)\lvec(35
30)\lvec(20 30)\lvec(20 20)

\end{texdraw}}\hskip 7mm
\raisebox{-0.4\height}{
\begin{texdraw}
\fontsize{6}{6}\drawdim em \setunitscale 0.1 \linewd 0.5

\move(10 10)\lvec(20 20)\lvec(20 10)\lvec(10 10)\lfill f:0.8

\htext(16 12){$\cdot$}

\move(0 10)\lvec(10 20)\lvec(10 10)\lvec(0 10)\lfill f:0.8

\htext(5 12){$*$}

\move(10 10)\lvec(20 20)\lvec(10 20)\lvec(10 10)\lfill f:0.8

\htext(11 15){$*$}

\move(0 -10)\lvec(0 10)\lvec(20 10)\lvec(20 -10)\lvec(0 -10)\lfill
f:0.8

\move(10 -10)\lvec(10 10)

\move(0 0)\lvec(-15 0)\lvec(-15 -10)\lvec(35 -10)\lvec(35
30)\lvec(20 30)\lvec(20 20)
\end{texdraw}}\vskip 5mm

\raisebox{-0.4\height}{
\begin{texdraw}
\fontsize{6}{6}\drawdim em \setunitscale 0.1 \linewd 0.5

\move(10 10)\lvec(20 20)\lvec(20 10)\lvec(10 10)\lfill f:0.8

\htext(16 12){$\cdot$}

\move(10 10)\lvec(20 20)\lvec(10 20)\lvec(10 10)\lfill f:0.8

\htext(11 15){$*$}

\move(0 -10)\lvec(0 10)\lvec(20 10)\lvec(20 -10)\lvec(0 -10)\lfill
f:0.8

\move(10 -10)\lvec(10 10)

\move(0 0)\lvec(-15 0)\lvec(-15 -10)\lvec(35 -10)\lvec(35
30)\lvec(20 30)\lvec(20 20)
\end{texdraw}}\hskip 7mm
\raisebox{-0.4\height}{
\begin{texdraw}
\fontsize{6}{6}\drawdim em \setunitscale 0.1 \linewd 0.5

\move(0 10)\lvec(10 20)\lvec(0 20)\lvec(0 10)\lfill f:0.8

\htext(2 15){$\cdot$}

\move(10 10)\lvec(20 20)\lvec(10 20)\lvec(10 10)\lfill f:0.8

\htext(11 15){$*$}

\move(0 -10)\lvec(0 10)\lvec(20 10)\lvec(20 -10)\lvec(0 -10)\lfill
f:0.8

\move(10 -10)\lvec(10 10)

\move(0 0)\lvec(-15 0)\lvec(-15 -10)\lvec(35 -10)\lvec(35
30)\lvec(20 30)\lvec(20 10)

\end{texdraw}}
\end{center}\vskip 10mm

\begin{center}
\raisebox{-0.4\height}{
\begin{texdraw}
\fontsize{6}{6}\drawdim em \setunitscale 0.1 \linewd 0.5

\move(0 -10)\lvec(0 10)\lvec(20 10)\lvec(20 -10)\lvec(0 -10)\lfill
f:0.8

\move(10 -10)\lvec(10 10)

\move(10 10)\lvec(20 20)\lvec(10 20)\lvec(10 10)\lfill f:0.8

\htext(12 16){$\cdot$}

\move(0 0)\lvec(-15 0)\lvec(-15 -10)\lvec(35 -10)\lvec(35
30)\lvec(20 30)\lvec(20 10)

\end{texdraw}}\hskip 7mm
\raisebox{-0.4\height}{
\begin{texdraw}
\fontsize{6}{6}\drawdim em \setunitscale 0.1 \linewd 0.5

\move(0 10)\lvec(10 20)\lvec(0 20)\lvec(0 10)\lfill f:0.8

\htext(2 16){$*$}

\move(10 10)\lvec(20 20)\lvec(10 20)\lvec(10 10)\lfill f:0.8

\htext(12 16){$\cdot$}

\move(0 -10)\lvec(0 10)\lvec(20 10)\lvec(20 -10)\lvec(0 -10)\lfill
f:0.8

\move(10 -10)\lvec(10 10)

\move(0 0)\lvec(-15 0)\lvec(-15 -10)\lvec(35 -10)\lvec(35
30)\lvec(20 30)\lvec(20 10)

\end{texdraw}}\hskip 7mm
\raisebox{-0.4\height}{
\begin{texdraw}
\fontsize{6}{6}\drawdim em \setunitscale 0.1 \linewd 0.5

\move(10 10)\lvec(20 20)\lvec(20 10)\lvec(10 10)\lfill f:0.8

\htext(15 12){$*$}

\move(0 10)\lvec(10 20)\lvec(0 20)\lvec(0 10)\lfill f:0.8

\htext(2 16){$*$}

\move(10 10)\lvec(20 20)\lvec(10 20)\lvec(10 10)\lfill f:0.8

\htext(12 16){$\cdot$}

\move(0 -10)\lvec(0 10)\lvec(20 10)\lvec(20 -10)\lvec(0 -10)\lfill
f:0.8

\move(10 -10)\lvec(10 10)

\move(0 0)\lvec(-15 0)\lvec(-15 -10)\lvec(35 -10)\lvec(35
30)\lvec(20 30)\lvec(20 20)
\end{texdraw}}\vskip 5mm

\raisebox{-0.4\height}{
\begin{texdraw}
\fontsize{6}{6}\drawdim em \setunitscale 0.1 \linewd 0.5

\move(10 10)\lvec(20 20)\lvec(20 10)\lvec(10 10)\lfill f:0.8

\htext(15 12){$*$}

\move(10 10)\lvec(20 20)\lvec(10 20)\lvec(10 10)\lfill f:0.8

\htext(12 15){$\cdot$}

\move(0 -10)\lvec(0 10)\lvec(20 10)\lvec(20 -10)\lvec(0 -10)\lfill
f:0.8

\move(10 -10)\lvec(10 10)

\move(0 0)\lvec(-15 0)\lvec(-15 -10)\lvec(35 -10)\lvec(35
30)\lvec(20 30)\lvec(20 20)
\end{texdraw}}\hskip 7mm
\raisebox{-0.4\height}{
\begin{texdraw}
\fontsize{6}{6}\drawdim em \setunitscale 0.1 \linewd 0.5

\move(10 10)\lvec(20 20)\lvec(20 10)\lvec(10 10)\lfill f:0.8

\htext(15 12){$*$}

\move(0 10)\lvec(10 20)\lvec(10 10)\lvec(0 10)\lfill f:0.8

\htext(6 12){$\cdot$}

\move(0 -10)\lvec(0 10)\lvec(20 10)\lvec(20 -10)\lvec(0 -10)\lfill
f:0.8

\move(10 -10)\lvec(10 10)

\move(0 0)\lvec(-15 0)\lvec(-15 -10)\lvec(35 -10)\lvec(35
30)\lvec(20 30)\lvec(20 20)

\end{texdraw}}\hskip 7mm
\end{center}\vskip 5mm

(${\rm III}_{2l}^+$) ($l\geq 1$)
\begin{center}
\raisebox{-0.4\height}{
\begin{texdraw}
\fontsize{6}{6}\drawdim em \setunitscale 0.1 \linewd 0.5

\move(0 -10)\lvec(0 10)\lvec(50 10)\lvec(50 -10)\lvec(0 -10)\lfill
f:0.8

\move(10 10)\lvec(10 -10)

\move(10 10)\lvec(20 20)\lvec(20 10)\lvec(10 10)\lfill f:0.8

\htext(16 12){$\cdot$}

\move(40 10)\lvec(50 20)\lvec(50 10)\lvec(40 10)\lfill f:0.8

\htext(45 12){$*$}

\htext(27 13){$\cdots$}

\move(10 10)\clvec(12 5)(15 5)(25 5)

\move(50 10)\clvec(48 5)(45 5)(35 5)

\htext(27 3){${2l}$}

\move(0 0)\lvec(-15 0)\lvec(-15 -10)\lvec(65 -10)\lvec(65
30)\lvec(50 30)\lvec(50 10)
\end{texdraw}}\hskip 10mm
\raisebox{-0.4\height}{
\begin{texdraw}
\fontsize{6}{6}\drawdim em \setunitscale 0.1 \linewd 0.5

\move(0 -10)\lvec(0 10)\lvec(50 10)\lvec(50 -10)\lvec(0 -10)\lfill
f:0.8

\move(10 10)\lvec(10 -10)

\move(10 10)\lvec(20 20)\lvec(20 10)\lvec(10 10)\lfill f:0.8

\htext(16 12){$\cdot$}

\move(40 10)\lvec(50 20)\lvec(50 10)\lvec(40 10)\lfill f:0.8

\htext(45 12){$*$}

\htext(27 13){$\cdots$}

\move(0 10)\lvec(10 20)\lvec(10 10)\lvec(0 10)\lfill f:0.8

\htext(5 12){$*$}

\move(0 0)\lvec(-15 0)\lvec(-15 -10)\lvec(65 -10)\lvec(65
30)\lvec(50 30)\lvec(50 10)
\end{texdraw}}\vskip 5mm

\raisebox{-0.4\height}{
\begin{texdraw}
\fontsize{6}{6}\drawdim em \setunitscale 0.1 \linewd 0.5

\move(0 -10)\lvec(0 10)\lvec(50 10)\lvec(50 -10)\lvec(0 -10)\lfill
f:0.8

\move(10 10)\lvec(10 -10)

\move(10 10)\lvec(10 -10)

\move(10 10)\lvec(20 20)\lvec(10 20)\lvec(10 10)\lfill f:0.8

\htext(12 16){$*$}

\htext(27 13){$\cdots$}

\move(40 10)\lvec(50 20)\lvec(40 20)\lvec(40 10)\lfill f:0.8

\htext(42 15){$\cdot$}

\move(0 0)\lvec(-15 0)\lvec(-15 -10)\lvec(65 -10)\lvec(65
30)\lvec(50 30)\lvec(50 10)
\end{texdraw}}\hskip 10mm
\raisebox{-0.4\height}{
\begin{texdraw}
\fontsize{6}{6}\drawdim em \setunitscale 0.1 \linewd 0.5

\move(0 -10)\lvec(0 10)\lvec(50 10)\lvec(50 -10)\lvec(0 -10)\lfill
f:0.8

\move(10 10)\lvec(10 -10)

\move(10 10)\lvec(20 20)\lvec(10 20)\lvec(10 10)\lfill f:0.8

\htext(12 16){$*$}

\move(40 10)\lvec(50 20)\lvec(50 10)\lvec(40 10)\lfill f:0.8

\htext(45 12){$*$}

\htext(27 13){$\cdots$}

\move(40 10)\lvec(50 20)\lvec(40 20)\lvec(40 10)\lfill f:0.8

\htext(42 15){$\cdot$}

\move(0 0)\lvec(-15 0)\lvec(-15 -10)\lvec(65 -10)\lvec(65
30)\lvec(50 30)\lvec(50 10)
\end{texdraw}}
\end{center}\vskip 10mm

\begin{center}
\raisebox{-0.4\height}{
\begin{texdraw}
\fontsize{6}{6}\drawdim em \setunitscale 0.1 \linewd 0.5

\move(0 -10)\lvec(0 10)\lvec(50 10)\lvec(50 -10)\lvec(0 -10)\lfill
f:0.8

\move(10 10)\lvec(10 -10)

\move(10 10)\lvec(10 -10)

\move(10 10)\lvec(20 20)\lvec(10 20)\lvec(10 10)\lfill f:0.8

\htext(12 16){$\cdot$}

\htext(27 13){$\cdots$}

\move(40 10)\lvec(50 20)\lvec(40 20)\lvec(40 10)\lfill f:0.8

\htext(42 16){$*$}

\move(0 0)\lvec(-15 0)\lvec(-15 -10)\lvec(65 -10)\lvec(65
30)\lvec(50 30)\lvec(50 10)

\end{texdraw}}\hskip 10mm
\raisebox{-0.4\height}{
\begin{texdraw}
\fontsize{6}{6}\drawdim em \setunitscale 0.1 \linewd 0.5

\move(0 -10)\lvec(0 10)\lvec(50 10)\lvec(50 -10)\lvec(0 -10)\lfill
f:0.8

\move(10 10)\lvec(10 -10)

\move(10 10)\lvec(10 -10)

\move(10 10)\lvec(20 20)\lvec(10 20)\lvec(10 10)\lfill f:0.8

\htext(12 16){$\cdot$}

\htext(27 13){$\cdots$}

\move(40 10)\lvec(50 20)\lvec(40 20)\lvec(40 10)\lfill f:0.8

\htext(42 16){$*$}

\move(0 10)\lvec(10 20)\lvec(0 20)\lvec(0 10)\lfill f:0.8

\htext(2 16){$*$}

\move(0 0)\lvec(-15 0)\lvec(-15 -10)\lvec(65 -10)\lvec(65
30)\lvec(50 30)\lvec(50 10)
\end{texdraw}}\vskip 5mm

\raisebox{-0.4\height}{
\begin{texdraw}
\fontsize{6}{6}\drawdim em \setunitscale 0.1 \linewd 0.5

\move(0 -10)\lvec(0 10)\lvec(50 10)\lvec(50 -10)\lvec(0 -10)\lfill
f:0.8

\move(10 10)\lvec(10 -10)

\move(10 10)\lvec(20 20)\lvec(20 10)\lvec(10 10)\lfill f:0.8

\htext(15 12){$*$}

\move(40 10)\lvec(50 20)\lvec(50 10)\lvec(40 10)\lfill f:0.8

\htext(46 12){$\cdot$}

\htext(27 13){$\cdots$}

\move(0 0)\lvec(-15 0)\lvec(-15 -10)\lvec(65 -10)\lvec(65
30)\lvec(50 30)\lvec(50 10)
\end{texdraw}}\hskip 10mm\raisebox{-0.4\height}{
\begin{texdraw}
\fontsize{6}{6}\drawdim em \setunitscale 0.1 \linewd 0.5

\move(0 -10)\lvec(0 10)\lvec(50 10)\lvec(50 -10)\lvec(0 -10)\lfill
f:0.8

\move(10 10)\lvec(10 -10)

\move(10 10)\lvec(20 20)\lvec(20 10)\lvec(10 10)\lfill f:0.8

\htext(15 12){$*$}

\move(40 10)\lvec(50 20)\lvec(50 10)\lvec(40 10)\lfill f:0.8

\htext(46 12){$\cdot$}

\move(40 10)\lvec(50 20)\lvec(40 20)\lvec(40 10)\lfill f:0.8

\htext(42 16){$*$}

\htext(27 13){$\cdots$}

\move(0 0)\lvec(-15 0)\lvec(-15 -10)\lvec(65 -10)\lvec(65
30)\lvec(50 30)\lvec(50 10)
\end{texdraw}}
\end{center}\vskip 5mm

(${\rm III}_{2l}^-$)  ($l\geq 1$)

\begin{center}
\raisebox{-0.4\height}{
\begin{texdraw}
\fontsize{6}{6}\drawdim em \setunitscale 0.1 \linewd 0.5

\move(0 -10)\lvec(0 10)\lvec(50 10)\lvec(50 -10)\lvec(0 -10)\lfill
f:0.8

\move(10 10)\lvec(10 -10)

\move(0 10)\lvec(10 20)\lvec(10 10)\lvec(0 10)\lfill f:0.8

\htext(6 12){$\cdot$}

\move(40 10)\lvec(50 20)\lvec(50 10)\lvec(40 10)\lfill f:0.8

\htext(46 12){$\cdot$}

\htext(22 13){$\cdots$}

\move(10 10)\clvec(12 5)(20 5)(25 5)

\move(50 10)\clvec(48 5)(40 5)(35 5)

\htext(27 3){$2l$}

\move(0 0)\lvec(-15 0)\lvec(-15 -10)\lvec(65 -10)\lvec(65
30)\lvec(50 30)\lvec(50 10)
\end{texdraw}}\hskip 10mm
\raisebox{-0.4\height}{
\begin{texdraw}
\fontsize{6}{6}\drawdim em \setunitscale 0.1 \linewd 0.5

\move(0 -10)\lvec(0 10)\lvec(50 10)\lvec(50 -10)\lvec(0 -10)\lfill
f:0.8

\move(10 10)\lvec(10 -10)

\move(0 10)\lvec(10 20)\lvec(10 10)\lvec(0 10)\lfill f:0.8

\htext(6 12){$\cdot$}

\move(40 10)\lvec(50 20)\lvec(50 10)\lvec(40 10)\lfill f:0.8

\htext(46 12){$\cdot$}

\htext(22 13){$\cdots$}

\move(40 10)\lvec(50 20)\lvec(40 20)\lvec(40 10)\lfill f:0.8

\htext(42 16){$*$}

\move(0 0)\lvec(-15 0)\lvec(-15 -10)\lvec(65 -10)\lvec(65
30)\lvec(50 30)\lvec(50 10)
\end{texdraw}}\vskip 5mm

\raisebox{-0.4\height}{
\begin{texdraw}
\fontsize{6}{6}\drawdim em \setunitscale 0.1 \linewd 0.5

\move(0 -10)\lvec(0 10)\lvec(50 10)\lvec(50 -10)\lvec(0 -10)\lfill
f:0.8

\move(10 10)\lvec(10 -10)

\move(10 10)\lvec(20 20)\lvec(10 20)\lvec(10 10)\lfill f:0.8

\htext(12 15){$\cdot$}

\move(40 10)\lvec(50 20)\lvec(50 10)\lvec(40 10)\lfill f:0.8

\htext(46 12){$\cdot$}

\htext(27 13){$\cdots$}

\move(40 10)\lvec(50 20)\lvec(40 20)\lvec(40 10)\lfill f:0.8

\htext(42 16){$*$}

\move(0 0)\lvec(-15 0)\lvec(-15 -10)\lvec(65 -10)\lvec(65
30)\lvec(50 30)\lvec(50 10)
\end{texdraw}}\hskip 10mm
\raisebox{-0.4\height}{
\begin{texdraw}
\fontsize{6}{6}\drawdim em \setunitscale 0.1 \linewd 0.5

\move(0 -10)\lvec(0 10)\lvec(50 10)\lvec(50 -10)\lvec(0 -10)\lfill
f:0.8

\move(10 10)\lvec(10 -10)

\move(0 10)\lvec(10 20)\lvec(0 20)\lvec(0 10)\lfill f:0.8

\htext(2 15){$*$}

\move(10 10)\lvec(20 20)\lvec(10 20)\lvec(10 10)\lfill f:0.8

\htext(12 15){$\cdot$}

\move(40 10)\lvec(50 20)\lvec(50 10)\lvec(40 10)\lfill f:0.8

\htext(47 12){$\cdot$}

\htext(27 13){$\cdots$}

\move(40 10)\lvec(50 20)\lvec(40 20)\lvec(40 10)\lfill f:0.8

\htext(42 16){$*$}

\move(0 0)\lvec(-15 0)\lvec(-15 -10)\lvec(65 -10)\lvec(65
30)\lvec(50 30)\lvec(50 10)
\end{texdraw}}
\end{center}\vskip 10mm

\begin{center}
\raisebox{-0.4\height}{
\begin{texdraw}
\fontsize{6}{6}\drawdim em \setunitscale 0.1 \linewd 0.5

\move(0 -10)\lvec(0 10)\lvec(50 10)\lvec(50 -10)\lvec(0 -10)\lfill
f:0.8

\move(10 10)\lvec(10 -10)

\move(0 10)\lvec(10 20)\lvec(0 20)\lvec(0 10)\lfill f:0.8

\htext(2 16){$\cdot$}

\move(10 10)\lvec(20 20)\lvec(10 20)\lvec(10 10)\lfill f:0.8

\htext(12 16){$*$}

\htext(27 13){$\cdots$}

\move(40 10)\lvec(50 20)\lvec(40 20)\lvec(40 10)\lfill f:0.8

\htext(42 16){$\cdot$}

\move(0 0)\lvec(-15 0)\lvec(-15 -10)\lvec(65 -10)\lvec(65
30)\lvec(50 30)\lvec(50 10)
\end{texdraw}}\hskip 10mm
\raisebox{-0.4\height}{
\begin{texdraw}
\fontsize{6}{6}\drawdim em \setunitscale 0.1 \linewd 0.5

\move(0 -10)\lvec(0 10)\lvec(50 10)\lvec(50 -10)\lvec(0 -10)\lfill
f:0.8

\move(10 10)\lvec(10 -10)

\move(0 10)\lvec(10 20)\lvec(0 20)\lvec(0 10)\lfill f:0.8

\htext(2 16){$\cdot$}

\move(10 10)\lvec(20 20)\lvec(10 20)\lvec(10 10)\lfill f:0.8

\htext(12 16){$*$}

\move(40 10)\lvec(50 20)\lvec(50 10)\lvec(40 10)\lfill f:0.8

\htext(46 12){$*$}

\htext(27 13){$\cdots$}

\move(40 10)\lvec(50 20)\lvec(40 20)\lvec(40 10)\lfill f:0.8

\htext(42 16){$\cdot$}

\move(0 0)\lvec(-15 0)\lvec(-15 -10)\lvec(65 -10)\lvec(65
30)\lvec(50 30)\lvec(50 10)
\end{texdraw}}\vskip 5mm

\raisebox{-0.4\height}{
\begin{texdraw}
\fontsize{6}{6}\drawdim em \setunitscale 0.1 \linewd 0.5

\move(0 -10)\lvec(0 10)\lvec(50 10)\lvec(50 -10)\lvec(0 -10)\lfill
f:0.8

\move(10 10)\lvec(10 -10)

\move(10 10)\lvec(20 20)\lvec(20 10)\lvec(10 10)\lfill f:0.8

\htext(16 12){$\cdot$}

\move(40 10)\lvec(50 20)\lvec(50 10)\lvec(40 10)\lfill f:0.8

\htext(46 12){$*$}

\move(40 10)\lvec(50 20)\lvec(40 20)\lvec(40 10)\lfill f:0.8

\htext(42 16){$\cdot$}

\htext(27 13){$\cdots$}

\move(0 0)\lvec(-15 0)\lvec(-15 -10)\lvec(65 -10)\lvec(65
30)\lvec(50 30)\lvec(50 10)
\end{texdraw}}\hskip 10mm
\raisebox{-0.4\height}{
\begin{texdraw}
\fontsize{6}{6}\drawdim em \setunitscale 0.1 \linewd 0.5

\move(0 -10)\lvec(0 10)\lvec(50 10)\lvec(50 -10)\lvec(0 -10)\lfill
f:0.8

\move(10 10)\lvec(10 -10)

\move(10 10)\lvec(20 20)\lvec(20 10)\lvec(10 10)\lfill f:0.8

\htext(16 12){$\cdot$}

\move(40 10)\lvec(50 20)\lvec(50 10)\lvec(40 10)\lfill f:0.8

\htext(46 12){$*$}

\htext(27 13){$\cdots$}

\move(0 10)\lvec(10 20)\lvec(10 10)\lvec(0 10)\lfill f:0.8

\htext(5 12){$*$}

\move(40 10)\lvec(50 20)\lvec(40 20)\lvec(40 10)\lfill f:0.8

\htext(42 16){$\cdot$}

\move(0 0)\lvec(-15 0)\lvec(-15 -10)\lvec(65 -10)\lvec(65
30)\lvec(50 30)\lvec(50 10)
\end{texdraw}}
\end{center}\vskip 5mm

(${\rm III}_{2l+1}$)  ($l\geq 1$)

\begin{center}
\raisebox{-0.4\height}{
\begin{texdraw}
\fontsize{6}{6}\drawdim em \setunitscale 0.1 \linewd 0.5

\move(0 -10)\lvec(0 10)\lvec(50 10)\lvec(50 -10)\lvec(0 -10)\lfill
f:0.8

\move(10 10)\lvec(10 -10)

\move(10 10)\lvec(20 20)\lvec(20 10)\lvec(10 10)\lfill f:0.8

\htext(16 12){$\cdot$}

\move(40 10)\lvec(50 20)\lvec(50 10)\lvec(40 10)\lfill f:0.8

\htext(46 12){$\cdot$}

\htext(27 13){$\cdots$}

\move(10 10)\clvec(12 5)(15 5)(20 5)

\move(50 10)\clvec(48 5)(45 5)(40 5)

\htext(21 3){$2l+1$}

\move(0 0)\lvec(-15 0)\lvec(-15 -10)\lvec(65 -10)\lvec(65
30)\lvec(50 30)\lvec(50 10)
\end{texdraw}}\hskip 10mm
\raisebox{-0.4\height}{
\begin{texdraw}
\fontsize{6}{6}\drawdim em \setunitscale 0.1 \linewd 0.5

\move(0 -10)\lvec(0 10)\lvec(50 10)\lvec(50 -10)\lvec(0 -10)\lfill
f:0.8

\move(10 10)\lvec(10 -10)

\move(10 10)\lvec(20 20)\lvec(20 10)\lvec(10 10)\lfill f:0.8

\htext(16 12){$\cdot$}

\move(40 10)\lvec(50 20)\lvec(50 10)\lvec(40 10)\lfill f:0.8

\htext(46 12){$\cdot$}

\htext(27 13){$\cdots$}

\move(0 10)\lvec(10 20)\lvec(10 10)\lvec(0 10)\lfill f:0.8

\htext(5 12){$*$}

\move(0 0)\lvec(-15 0)\lvec(-15 -10)\lvec(65 -10)\lvec(65
30)\lvec(50 30)\lvec(50 10)
\end{texdraw}}\vskip 5mm

\raisebox{-0.4\height}{
\begin{texdraw}
\fontsize{6}{6}\drawdim em \setunitscale 0.1 \linewd 0.5

\move(0 -10)\lvec(0 10)\lvec(50 10)\lvec(50 -10)\lvec(0 -10)\lfill
f:0.8

\move(10 10)\lvec(10 -10)

\move(10 10)\lvec(20 20)\lvec(20 10)\lvec(10 10)\lfill f:0.8

\htext(16 12){$\cdot$}

\move(40 10)\lvec(50 20)\lvec(50 10)\lvec(40 10)\lfill f:0.8

\htext(46 12){$\cdot$}

\htext(27 13){$\cdots$}

\move(0 10)\lvec(10 20)\lvec(10 10)\lvec(0 10)\lfill f:0.8

\htext(5 12){$*$}

\move(40 10)\lvec(50 20)\lvec(40 20)\lvec(40 10)\lfill f:0.8

\htext(42 16){$*$}

\move(0 0)\lvec(-15 0)\lvec(-15 -10)\lvec(65 -10)\lvec(65
30)\lvec(50 30)\lvec(50 10)
\end{texdraw}}\hskip 10mm
\raisebox{-0.4\height}{
\begin{texdraw}
\fontsize{6}{6}\drawdim em \setunitscale 0.1 \linewd 0.5

\move(0 -10)\lvec(0 10)\lvec(50 10)\lvec(50 -10)\lvec(0 -10)\lfill
f:0.8

\move(10 10)\lvec(10 -10)

\move(10 10)\lvec(20 20)\lvec(20 10)\lvec(10 10)\lfill f:0.8

\htext(16 12){$\cdot$}

\move(40 10)\lvec(50 20)\lvec(50 10)\lvec(40 10)\lfill f:0.8

\htext(46 12){$\cdot$}

\htext(27 13){$\cdots$}

\move(40 10)\lvec(50 20)\lvec(40 20)\lvec(40 10)\lfill f:0.8

\htext(42 16){$*$}

\move(0 0)\lvec(-15 0)\lvec(-15 -10)\lvec(65 -10)\lvec(65
30)\lvec(50 30)\lvec(50 10)
\end{texdraw}}\vskip 5mm

\raisebox{-0.4\height}{
\begin{texdraw}
\fontsize{6}{6}\drawdim em \setunitscale 0.1 \linewd 0.5

\move(0 -10)\lvec(0 10)\lvec(50 10)\lvec(50 -10)\lvec(0 -10)\lfill
f:0.8

\move(10 10)\lvec(10 -10)

\move(0 10)\lvec(10 20)\lvec(0 20)\lvec(0 10)\lfill f:0.8

\htext(2 16){$\cdot$}

\htext(20 13){$\cdots$}

\move(40 10)\lvec(50 20)\lvec(40 20)\lvec(40 10)\lfill f:0.8

\htext(42 16){$*$}

\move(0 0)\lvec(-15 0)\lvec(-15 -10)\lvec(65 -10)\lvec(65
30)\lvec(50 30)\lvec(50 10)
\end{texdraw}}\hskip 1cm
\raisebox{-0.4\height}{
\begin{texdraw}
\fontsize{6}{6}\drawdim em \setunitscale 0.1 \linewd 0.5

\move(0 -10)\lvec(0 10)\lvec(50 10)\lvec(50 -10)\lvec(0 -10)\lfill
f:0.8

\move(10 10)\lvec(10 -10)

\move(10 10)\lvec(20 20)\lvec(10 20)\lvec(10 10)\lfill f:0.8

\htext(12 16){$*$}

\move(40 10)\lvec(50 20)\lvec(50 10)\lvec(40 10)\lfill f:0.8

\htext(46 12){$\cdot$}

\htext(27 13){$\cdots$}

\move(40 10)\lvec(50 20)\lvec(40 20)\lvec(40 10)\lfill f:0.8

\htext(42 16){$*$}

\move(0 0)\lvec(-15 0)\lvec(-15 -10)\lvec(65 -10)\lvec(65
30)\lvec(50 30)\lvec(50 10)
\end{texdraw}}
\end{center}\vskip 10mm

\begin{center}
\raisebox{-0.4\height}{
\begin{texdraw}
\fontsize{6}{6}\drawdim em \setunitscale 0.1 \linewd 0.5

\move(0 -10)\lvec(0 10)\lvec(50 10)\lvec(50 -10)\lvec(0 -10)\lfill
f:0.8

\move(10 10)\lvec(10 -10)

\move(10 10)\lvec(10 -10)

\move(10 10)\lvec(20 20)\lvec(10 20)\lvec(10 10)\lfill f:0.8

\htext(12 16){$\cdot$}

\htext(27 13){$\cdots$}

\move(40 10)\lvec(50 20)\lvec(40 20)\lvec(40 10)\lfill f:0.8

\htext(42 16){$\cdot$}

\move(0 0)\lvec(-15 0)\lvec(-15 -10)\lvec(65 -10)\lvec(65
30)\lvec(50 30)\lvec(50 10)

\end{texdraw}}\hskip 10mm
\raisebox{-0.4\height}{
\begin{texdraw}
\fontsize{6}{6}\drawdim em \setunitscale 0.1 \linewd 0.5

\move(0 -10)\lvec(0 10)\lvec(50 10)\lvec(50 -10)\lvec(0 -10)\lfill
f:0.8

\move(10 10)\lvec(10 -10)

\move(0 10)\lvec(10 20)\lvec(0 20)\lvec(0 10)\lfill f:0.8

\htext(2 16){$*$}

\move(10 10)\lvec(20 20)\lvec(10 20)\lvec(10 10)\lfill f:0.8

\htext(12 16){$\cdot$}

\htext(27 13){$\cdots$}

\move(40 10)\lvec(50 20)\lvec(40 20)\lvec(40 10)\lfill f:0.8

\htext(42 16){$\cdot$}

\move(0 0)\lvec(-15 0)\lvec(-15 -10)\lvec(65 -10)\lvec(65
30)\lvec(50 30)\lvec(50 10)
\end{texdraw}}\vskip 5mm

\raisebox{-0.4\height}{
\begin{texdraw}
\fontsize{6}{6}\drawdim em \setunitscale 0.1 \linewd 0.5

\move(0 -10)\lvec(0 10)\lvec(50 10)\lvec(50 -10)\lvec(0 -10)\lfill
f:0.8

\move(10 10)\lvec(10 -10)

\move(0 10)\lvec(10 20)\lvec(0 20)\lvec(0 10)\lfill f:0.8

\htext(2 16){$*$}

\move(10 10)\lvec(20 20)\lvec(10 20)\lvec(10 10)\lfill f:0.8

\htext(12 16){$\cdot$}

\move(40 10)\lvec(50 20)\lvec(50 10)\lvec(40 10)\lfill f:0.8

\htext(46 12){$*$}

\htext(27 13){$\cdots$}

\move(40 10)\lvec(50 20)\lvec(40 20)\lvec(40 10)\lfill f:0.8

\htext(42 16){$\cdot$}

\move(0 0)\lvec(-15 0)\lvec(-15 -10)\lvec(65 -10)\lvec(65
30)\lvec(50 30)\lvec(50 10)
\end{texdraw}}\hskip 10mm
\raisebox{-0.4\height}{
\begin{texdraw}
\fontsize{6}{6}\drawdim em \setunitscale 0.1 \linewd 0.5

\move(0 -10)\lvec(0 10)\lvec(50 10)\lvec(50 -10)\lvec(0 -10)\lfill
f:0.8

\move(10 10)\lvec(10 -10)

\move(10 10)\lvec(20 20)\lvec(10 20)\lvec(10 10)\lfill f:0.8

\htext(12 16){$\cdot$}

\move(40 10)\lvec(50 20)\lvec(50 10)\lvec(40 10)\lfill f:0.8

\htext(46 12){$*$}

\htext(27 13){$\cdots$}

\move(40 10)\lvec(50 20)\lvec(40 20)\lvec(40 10)\lfill f:0.8

\htext(42 16){$\cdot$}

\move(0 0)\lvec(-15 0)\lvec(-15 -10)\lvec(65 -10)\lvec(65
30)\lvec(50 30)\lvec(50 10)
\end{texdraw}}\vskip 5mm

\raisebox{-0.4\height}{
\begin{texdraw}
\fontsize{6}{6}\drawdim em \setunitscale 0.1 \linewd 0.5

\move(0 -10)\lvec(0 10)\lvec(50 10)\lvec(50 -10)\lvec(0 -10)\lfill
f:0.8

\move(10 10)\lvec(10 -10)

\move(10 10)\lvec(20 20)\lvec(20 10)\lvec(10 10)\lfill f:0.8

\htext(16 12){$*$}

\move(40 10)\lvec(50 20)\lvec(50 10)\lvec(40 10)\lfill f:0.8

\htext(46 12){$*$}

\htext(27 13){$\cdots$}

\move(0 10)\lvec(10 20)\lvec(10 10)\lvec(0 10)\lfill f:0.8

\htext(6 12){$\cdot$}

\move(0 0)\lvec(-15 0)\lvec(-15 -10)\lvec(65 -10)\lvec(65
30)\lvec(50 30)\lvec(50 10)
\end{texdraw}}\hskip 10mm
\raisebox{-0.4\height}{
\begin{texdraw}
\fontsize{6}{6}\drawdim em \setunitscale 0.1 \linewd 0.5

\move(0 -10)\lvec(0 10)\lvec(50 10)\lvec(50 -10)\lvec(0 -10)\lfill
f:0.8

\move(10 10)\lvec(10 -10)

\move(10 10)\lvec(20 20)\lvec(20 10)\lvec(10 10)\lfill f:0.8

\htext(16 12){$*$}

\move(40 10)\lvec(50 20)\lvec(50 10)\lvec(40 10)\lfill f:0.8

\htext(46 12){$*$}

\htext(27 13){$\cdots$}

\move(40 10)\lvec(50 20)\lvec(40 20)\lvec(40 10)\lfill f:0.8

\htext(42 16){$\cdot$}

\move(0 0)\lvec(-15 0)\lvec(-15 -10)\lvec(65 -10)\lvec(65
30)\lvec(50 30)\lvec(50 10)
\end{texdraw}}\quad .
\end{center}\vskip 5mm

In case $\rm III_{\infty}$, let $y_N$ be the column containing the
admissible $i$-slot in the ground-state wall or the removable
$i$-block that touches the ground-state wall as indicated in the
above figure. We denote by $Y_0$ the part of $Y$ consisting of
$y_N$ and the blocks in the ground-state wall lying in the left of
$y_N$, and call it an {\it $i$-component of type $\rm
III_{\infty}$}.

In case ${\rm III}_{0}$, ${\rm III}_{1}$, ${\rm III}^{\pm}_{2l}$
and ${\rm III}_{2l+1}$ ($l\geq 1$), the whole shaded part
containing (virtually) admissible $i$-slots and (virtually)
removable $i$-blocks will be called an {\it $i$-component of type}
${\rm III}_{0}$, ${\rm III}_{1}$, ${\rm III}^{\pm}_{2l}$ and ${\rm
III}_{2l+1}$ ($l\geq 1$), respectively.

Let $Y_1$ be the left-most $i$-component of type ${\rm III}_{0}$,
${\rm III}_{1}$, ${\rm III}^{\pm}_{2l}$ and ${\rm III}_{2l+1}$. If
there is no $i$-component of type $\rm III_{\infty}$, we denote by
$Y_0$ the part of $Y$ consisting of blocks lying in the left of
$Y_1$, and call it a {\it trivial $i$-component}. The parts of $Y$
lying between two $i$-components of type ${\rm III}_{\infty}$,
${\rm III}_{0}$, ${\rm III}_{1}$, ${\rm III}^{\pm}_{2l}$ and ${\rm
III}_{2l+1}$ will also be called the {\it trivial $i$-components}.

In this way, we obtain a unique decomposition
$Y=(Y_0,Y_1,\cdots,Y_r)$ of $Y$, where each $Y_k$ is an
$i$-component of type ${\rm III}_{\infty}$, ${\rm III}_{0}$, ${\rm
III}_{1}$, ${\rm III}^{\pm}_{2l}$, ${\rm III}_{2l+1}$ or a trivial
$i$-component.

\begin{ex}\label{comp3}{\rm
Let $\frak{g}=B^{(1)}_3$, $\Lambda=\Lambda_0$ and $i=0$. If \vskip
5mm

\begin{center}
$Y=$\raisebox{-0.5\height}{
\begin{texdraw}
\fontsize{8}{8}\drawdim em \setunitscale 0.13 \linewd 0.5


\move(10 0)\lvec(20 0)\lvec(20 10)\lvec(10 10)\lvec(10 0)\htext(12
6){\tiny $1$}

\move(20 0)\lvec(30 0)\lvec(30 10)\lvec(20 10)\lvec(20 0)\htext(22
6){\tiny$0$}

\move(30 0)\lvec(40 0)\lvec(40 10)\lvec(30 10)\lvec(30 0)\htext(32
6){\tiny $1$}

\move(40 0)\lvec(50 0)\lvec(50 10)\lvec(40 10)\lvec(40 0)\htext(42
6){\tiny $0$}

\move(50 0)\lvec(60 0)\lvec(60 10)\lvec(50 10)\lvec(50 0)\htext(52
6){\tiny $1$}

\move(60 0)\lvec(70 0)\lvec(70 10)\lvec(60 10)\lvec(60 0)\htext(62
6){\tiny $0$}

\move(70 0)\lvec(80 0)\lvec(80 10)\lvec(70 10)\lvec(70 0)\htext(72
6){\tiny $1$}

\move(80 0)\lvec(90 0)\lvec(90 10)\lvec(80 10)\lvec(80 0)\htext(82
6){\tiny $0$}

%

\move(10 0)\lvec(20 10)\lvec(20 0)\lvec(10 0)\lfill f:0.8
\htext(16 2){\tiny $0$}

\move(20 0)\lvec(30 10)\lvec(30 0)\lvec(20 0)\lfill f:0.8
\htext(26 2){\tiny $1$}

\move(30 0)\lvec(40 10)\lvec(40 0)\lvec(30 0)\lfill f:0.8
\htext(36 2){\tiny $0$}

\move(40 0)\lvec(50 10)\lvec(50 0)\lvec(40 0)\lfill f:0.8
\htext(46 2){\tiny $1$}

\move(50 0)\lvec(60 10)\lvec(60 0)\lvec(50 0)\lfill f:0.8
\htext(56 2){\tiny $0$}

\move(60 0)\lvec(70 10)\lvec(70 0)\lvec(60 0)\lfill f:0.8
\htext(66 2){\tiny $1$}

\move(70 0)\lvec(80 10)\lvec(80 0)\lvec(70 0)\lfill f:0.8
\htext(76 2){\tiny $0$}

\move(80 0)\lvec(90 10)\lvec(90 0)\lvec(80 0)\lfill f:0.8
\htext(86 2){\tiny $1$}

\move(10 10)\lvec(20 10)\lvec(20 20)\lvec(10 20)\lvec(10
10)\htext(13 13){$2$}

\move(20 10)\lvec(30 10)\lvec(30 20)\lvec(20 20)\lvec(20
10)\htext(23 13){$2$}

\move(30 10)\lvec(40 10)\lvec(40 20)\lvec(30 20)\lvec(30
10)\htext(33 13){$2$}

\move(40 10)\lvec(50 10)\lvec(50 20)\lvec(40 20)\lvec(40
10)\htext(43 13){$2$}

\move(50 10)\lvec(60 10)\lvec(60 20)\lvec(50 20)\lvec(50
10)\htext(53 13){$2$}

\move(60 10)\lvec(70 10)\lvec(70 20)\lvec(60 20)\lvec(60
10)\htext(63 13){$2$}

\move(70 10)\lvec(80 10)\lvec(80 20)\lvec(70 20)\lvec(70
10)\htext(73 13){$2$}

\move(80 10)\lvec(90 10)\lvec(90 20)\lvec(80 20)\lvec(80
10)\htext(83 13){$2$}
\move(10 20)\lvec(20 20)\lvec(20 30)\lvec(10 30)\lvec(10
20)\htext(13 26){\tiny $3$}

\move(20 20)\lvec(30 20)\lvec(30 30)\lvec(20 30)\lvec(20
20)\htext(23 26){\tiny $3$}

\move(30 20)\lvec(40 20)\lvec(40 30)\lvec(30 30)\lvec(30
20)\htext(33 26){\tiny $3$}

\move(40 20)\lvec(50 20)\lvec(50 30)\lvec(40 30)\lvec(40
20)\htext(43 26){\tiny $3$}

\move(50 20)\lvec(60 20)\lvec(60 30)\lvec(50 30)\lvec(50
20)\htext(53 26){\tiny $3$}

\move(60 20)\lvec(70 20)\lvec(70 30)\lvec(60 30)\lvec(60
20)\htext(63 26){\tiny $3$}

\move(70 20)\lvec(80 20)\lvec(80 30)\lvec(70 30)\lvec(70
20)\htext(73 26){\tiny $3$}

\move(80 20)\lvec(90 20)\lvec(90 30)\lvec(80 30)\lvec(80
20)\htext(83 26){\tiny $3$}

\move(10 20)\lvec(20 20)\lvec(20 25)\lvec(10 25)\lvec(10 20)
\htext(13 21){\tiny $3$}

\move(20 20)\lvec(30 20)\lvec(30 25)\lvec(20 25)\lvec(20 20)
\htext(23 21){\tiny $3$}

\move(30 20)\lvec(40 20)\lvec(40 25)\lvec(30 25)\lvec(30 20)
\htext(33 21){\tiny $3$}

\move(40 20)\lvec(50 20)\lvec(50 25)\lvec(40 25)\lvec(40 20)
\htext(43 21){\tiny $3$}

\move(50 20)\lvec(60 20)\lvec(60 25)\lvec(50 25)\lvec(50
20)\htext(53 21){\tiny $3$}

\move(60 20)\lvec(70 20)\lvec(70 25)\lvec(60 25)\lvec(60
20)\htext(63 21){\tiny $3$}

\move(70 20)\lvec(80 20)\lvec(80 25)\lvec(70 25)\lvec(70
20)\htext(73 21){\tiny $3$}

\move(80 20)\lvec(90 20)\lvec(90 25)\lvec(80 25)\lvec(80
20)\htext(83 21){\tiny $3$}
\move(10 30)\lvec(20 30)\lvec(20 40)\lvec(10 40)\lvec(10
30)\htext(13 33){$2$}

\move(20 30)\lvec(30 30)\lvec(30 40)\lvec(20 40)\lvec(20
30)\htext(23 33){$2$}

\move(30 30)\lvec(40 30)\lvec(40 40)\lvec(30 40)\lvec(30
30)\htext(33 33){$2$}

\move(40 30)\lvec(50 30)\lvec(50 40)\lvec(40 40)\lvec(40
30)\htext(43 33){$2$}

\move(50 30)\lvec(60 30)\lvec(60 40)\lvec(50 40)\lvec(50
30)\htext(53 33){$2$}

\move(60 30)\lvec(70 30)\lvec(70 40)\lvec(60 40)\lvec(60
30)\htext(63 33){$2$}

\move(70 30)\lvec(80 30)\lvec(80 40)\lvec(70 40)\lvec(70
30)\htext(73 33){$2$}

\move(80 30)\lvec(90 30)\lvec(90 40)\lvec(80 40)\lvec(80
30)\htext(83 33){$2$}
%


\move(20 40)\lvec(30 50)\lvec(20 50)\lvec(20 40) \htext(22
46){\tiny $0$}

\move(30 40)\lvec(40 50)\lvec(30 50)\lvec(30 40) \htext(32
46){\tiny $1$}

\move(40 40)\lvec(50 50)\lvec(40 50)\lvec(40 40) \htext(42
46){\tiny $0$}

\move(40 40)\lvec(50 50)\lvec(50 40)\lvec(40 40) \htext(46
42){\tiny $1$}

\move(50 40)\lvec(60 40)\lvec(60 50)\lvec(50 50) \lvec(50
40)\htext(52 46){\tiny $1$}

\move(60 40)\lvec(70 40)\lvec(70 50)\lvec(60 50)\lvec(60
40)\htext(62 46){\tiny $0$}

\move(70 40)\lvec(80 40)\lvec(80 50)\lvec(70 50)\lvec(70
40)\htext(72 46){\tiny $1$}

\move(80 40)\lvec(90 40)\lvec(90 50)\lvec(80 50)\lvec(80
40)\htext(82 46){\tiny $0$}

\move(50 40)\lvec(60 50)\lvec(60 40)\lvec(50 40) \htext(56
42){\tiny $0$}

\move(60 40)\lvec(70 50)\lvec(70 40)\lvec(60 40) \htext(66
42){\tiny $1$}

\move(70 40)\lvec(80 50)\lvec(80 40)\lvec(70 40) \htext(76
42){\tiny $0$}

\move(80 40)\lvec(90 50)\lvec(90 40)\lvec(80 40) \htext(86
42){\tiny $1$}
\move(50 50)\lvec(60 50)\lvec(60 60)\lvec(50 60)\lvec(50
50)\htext(53 53){$2$}

\move(60 50)\lvec(70 50)\lvec(70 60)\lvec(60 60)\lvec(60
50)\htext(63 53){$2$}

\move(70 50)\lvec(80 50)\lvec(80 60)\lvec(70 60)\lvec(70
50)\htext(73 53){$2$}

\move(80 50)\lvec(90 50)\lvec(90 60)\lvec(80 60)\lvec(80
50)\htext(83 53){$2$}
\move(60 60)\lvec(70 60)\lvec(70 65)\lvec(60 65)\lvec(60
60)\htext(63 61){\tiny $3$}

\move(70 60)\lvec(80 60)\lvec(80 65)\lvec(70 65)\lvec(70
60)\htext(73 61){\tiny $3$}

\move(80 60)\lvec(90 60)\lvec(90 65)\lvec(80 65)\lvec(80
60)\htext(83 61){\tiny $3$}
\move(60 60)\lvec(70 60)\lvec(70 70)\lvec(60 70)\lvec(60
60)\htext(63 66){\tiny $3$}

\move(70 60)\lvec(80 60)\lvec(80 70)\lvec(70 70)\lvec(70
60)\htext(73 66){\tiny $3$}

\move(80 60)\lvec(90 60)\lvec(90 70)\lvec(80 70)\lvec(80
60)\htext(83 66){\tiny $3$}
\move(60 70)\lvec(70 70)\lvec(70 80)\lvec(60 80)\lvec(60
70)\htext(63 73){$2$}

\move(70 70)\lvec(80 70)\lvec(80 80)\lvec(70 80)\lvec(70
70)\htext(73 73){$2$}

\move(80 70)\lvec(90 70)\lvec(90 80)\lvec(80 80)\lvec(80
70)\htext(83 73){$2$}
\move(60 80)\lvec(70 90)\lvec(70 80)\lvec(60 80)\htext(66
82){\tiny $1$}

\move(70 80)\lvec(80 90)\lvec(80 80)\lvec(70 80)\htext(76
82){\tiny $0$}

\move(80 80)\lvec(90 90)\lvec(90 80)\lvec(80 80)\htext(86
82){\tiny $1$}
\end{texdraw}}\quad ,
\end{center}\vskip 5mm

\noindent then we have $Y=(Y_0,Y_1,Y_2,Y_3)$, where \vskip 5mm

\hskip 1cm$Y_0=$ \raisebox{-0.4\height}{\begin{texdraw}\drawdim em
\setunitscale 0.13 \linewd 0.5

\htext(-.2 3){$\cdots$} \move(10 0)\lvec(20 10)\lvec(20 0)\lvec(10
0)\lfill f:0.8 \htext(16 2){\tiny $0$}

\move(20 0)\lvec(30 10)\lvec(30 0)\lvec(20 0)\lfill f:0.8
\htext(26 2){\tiny $1$}

\move(30 0)\lvec(40 10)\lvec(40 0)\lvec(30 0)\lfill f:0.8
\htext(36 2){\tiny $0$}

\move(40 0)\lvec(50 10)\lvec(50 0)\lvec(40 0)\lfill f:0.8
\htext(46 2){\tiny $1$}
\end{texdraw}} : $i$-component of type $\rm III_{\infty}$, \vskip 5mm

\hskip 10mm$Y_1=$ \raisebox{-0.4\height}{\begin{texdraw}
\fontsize{8}{8}\drawdim em \setunitscale 0.13 \linewd 0.5

\move(10 0)\lvec(20 0)\lvec(20 10)\lvec(10 10)\lvec(10 0)\htext(12
6){\tiny $1$}

\move(20 0)\lvec(30 0)\lvec(30 10)\lvec(20 10)\lvec(20 0)\htext(22
6){\tiny $0$}

\move(30 0)\lvec(40 0)\lvec(40 10)\lvec(30 10)\lvec(30 0)\htext(32
6){\tiny $1$}

\move(40 0)\lvec(50 0)\lvec(50 10)\lvec(40 10)\lvec(40 0)\htext(42
6){\tiny $0$}

\move(10 0)\lvec(20 10)\lvec(20 0)\lvec(10 0)\lfill f:0.8
\htext(16 2){\tiny $0$}

\move(20 0)\lvec(30 10)\lvec(30 0)\lvec(20 0)\lfill f:0.8
\htext(26 2){\tiny $1$}

\move(30 0)\lvec(40 10)\lvec(40 0)\lvec(30 0)\lfill f:0.8
\htext(36 2){\tiny $0$}

\move(40 0)\lvec(50 10)\lvec(50 0)\lvec(40 0)\lfill f:0.8
\htext(46 2){\tiny $1$}

\move(10 10)\lvec(20 10)\lvec(20 20)\lvec(10 20)\lvec(10
10)\htext(13 13){$2$}

\move(20 10)\lvec(30 10)\lvec(30 20)\lvec(20 20)\lvec(20
10)\htext(23 13){$2$}

\move(30 10)\lvec(40 10)\lvec(40 20)\lvec(30 20)\lvec(30
10)\htext(33 13){$2$}

\move(40 10)\lvec(50 10)\lvec(50 20)\lvec(40 20)\lvec(40
10)\htext(43 13){$2$}

\move(10 20)\lvec(20 20)\lvec(20 30)\lvec(10 30)\lvec(10
20)\htext(13 26){\tiny $3$}

\move(20 20)\lvec(30 20)\lvec(30 30)\lvec(20 30)\lvec(20
20)\htext(23 26){\tiny $3$}

\move(30 20)\lvec(40 20)\lvec(40 30)\lvec(30 30)\lvec(30
20)\htext(33 26){\tiny $3$}

\move(40 20)\lvec(50 20)\lvec(50 30)\lvec(40 30)\lvec(40
20)\htext(43 26){\tiny $3$}

\move(10 20)\lvec(20 20)\lvec(20 25)\lvec(10 25)\lvec(10 20)
\htext(13 21){\tiny $3$}

\move(20 20)\lvec(30 20)\lvec(30 25)\lvec(20 25)\lvec(20 20)
\htext(23 21){\tiny $3$}

\move(30 20)\lvec(40 20)\lvec(40 25)\lvec(30 25)\lvec(30 20)
\htext(33 21){\tiny $3$}

\move(40 20)\lvec(50 20)\lvec(50 25)\lvec(40 25)\lvec(40 20)
\htext(43 21){\tiny $3$}

\move(10 30)\lvec(20 30)\lvec(20 40)\lvec(10 40)\lvec(10
30)\htext(13 33){$2$}

\move(20 30)\lvec(30 30)\lvec(30 40)\lvec(20 40)\lvec(20
30)\htext(23 33){$2$}

\move(30 30)\lvec(40 30)\lvec(40 40)\lvec(30 40)\lvec(30
30)\htext(33 33){$2$}

\move(40 30)\lvec(50 30)\lvec(50 40)\lvec(40 40)\lvec(40
30)\htext(43 33){$2$}

\move(20 40)\lvec(30 50)\lvec(20 50)\lvec(20 40) \htext(22
46){\tiny $0$}

\move(30 40)\lvec(40 50)\lvec(30 50)\lvec(30 40) \htext(32
46){\tiny $1$}

\move(40 40)\lvec(50 50)\lvec(40 50)\lvec(40 40) \htext(42
46){\tiny $0$}

\move(40 40)\lvec(50 50)\lvec(50 40)\lvec(40 40) \htext(46
42){\tiny $1$}
\end{texdraw}} : $i$-component of type $\rm III_3$,\vskip 5mm

\hskip 10mm$Y_2=$ \raisebox{-0.4\height}{\begin{texdraw}
\fontsize{8}{8}\drawdim em \setunitscale 0.13 \linewd 0.5

\move(50 0)\lvec(60 0)\lvec(60 10)\lvec(50 10)\lvec(50 0)\htext(52
6){\tiny $1$}

\move(50 0)\lvec(60 10)\lvec(60 0)\lvec(50 0)\lfill f:0.8
\htext(56 2){\tiny $0$}

\move(50 10)\lvec(60 10)\lvec(60 20)\lvec(50 20)\lvec(50
10)\htext(53 13){$2$}

\move(50 20)\lvec(60 20)\lvec(60 30)\lvec(50 30)\lvec(50
20)\htext(53 26){\tiny $3$}

\move(50 20)\lvec(60 20)\lvec(60 25)\lvec(50 25)\lvec(50
20)\htext(53 21){\tiny $3$}

\move(50 30)\lvec(60 30)\lvec(60 40)\lvec(50 40)\lvec(50
30)\htext(53 33){$2$}

\move(50 40)\lvec(60 40)\lvec(60 50)\lvec(50 50) \lvec(50
40)\htext(52 46){\tiny $1$}

\move(50 40)\lvec(60 50)\lvec(60 40)\lvec(50 40) \htext(56
42){\tiny $0$}

\move(50 50)\lvec(60 50)\lvec(60 60)\lvec(50 60)\lvec(50
50)\htext(53 53){$2$}

\end{texdraw}} : trivial $i$-component and\vskip 5mm

\hskip 10mm$Y_3=$ \raisebox{-0.4\height}{\begin{texdraw}
\fontsize{8}{8}\drawdim em \setunitscale 0.13 \linewd 0.5

\move(60 0)\lvec(70 0)\lvec(70 10)\lvec(60 10)\lvec(60 0)\htext(62
6){\tiny $0$}

\move(70 0)\lvec(80 0)\lvec(80 10)\lvec(70 10)\lvec(70 0)\htext(72
6){\tiny $1$}

\move(80 0)\lvec(90 0)\lvec(90 10)\lvec(80 10)\lvec(80 0)\htext(82
6){\tiny $0$}

\move(60 0)\lvec(70 10)\lvec(70 0)\lvec(60 0)\lfill f:0.8
\htext(66 2){\tiny $1$}

\move(70 0)\lvec(80 10)\lvec(80 0)\lvec(70 0)\lfill f:0.8
\htext(76 2){\tiny $0$}

\move(80 0)\lvec(90 10)\lvec(90 0)\lvec(80 0)\lfill f:0.8
\htext(86 2){\tiny $1$}

\move(60 10)\lvec(70 10)\lvec(70 20)\lvec(60 20)\lvec(60
10)\htext(63 13){$2$}

\move(70 10)\lvec(80 10)\lvec(80 20)\lvec(70 20)\lvec(70
10)\htext(73 13){$2$}

\move(80 10)\lvec(90 10)\lvec(90 20)\lvec(80 20)\lvec(80
10)\htext(83 13){$2$}

\move(60 20)\lvec(70 20)\lvec(70 25)\lvec(60 25)\lvec(60
20)\htext(63 26){\tiny $3$}

\move(70 20)\lvec(80 20)\lvec(80 25)\lvec(70 25)\lvec(70
20)\htext(73 26){\tiny $3$}

\move(80 20)\lvec(90 20)\lvec(90 25)\lvec(80 25)\lvec(80
20)\htext(83 26){\tiny $3$}

\move(60 25)\lvec(70 25)\lvec(70 30)\lvec(60 30)\lvec(60
25)\htext(63 21){\tiny $3$}

\move(70 25)\lvec(80 25)\lvec(80 30)\lvec(70 30)\lvec(70
25)\htext(73 21){\tiny $3$}

\move(80 25)\lvec(90 25)\lvec(90 30)\lvec(80 30)\lvec(80
25)\htext(83 21){\tiny $3$}

\move(60 30)\lvec(70 30)\lvec(70 40)\lvec(60 40)\lvec(60
30)\htext(63 33){$2$}

\move(70 30)\lvec(80 30)\lvec(80 40)\lvec(70 40)\lvec(70
30)\htext(73 33){$2$}

\move(80 30)\lvec(90 30)\lvec(90 40)\lvec(80 40)\lvec(80
30)\htext(83 33){$2$}

\move(60 40)\lvec(70 40)\lvec(70 50)\lvec(60 50)\lvec(60
40)\htext(62 46){\tiny $0$}

\move(70 40)\lvec(80 40)\lvec(80 50)\lvec(70 50)\lvec(70
40)\htext(72 46){\tiny $1$}

\move(80 40)\lvec(90 40)\lvec(90 50)\lvec(80 50)\lvec(80
40)\htext(82 46){\tiny $0$}

\move(60 40)\lvec(70 50)\lvec(70 40)\lvec(60 40) \htext(66
42){\tiny $1$}

\move(70 40)\lvec(80 50)\lvec(80 40)\lvec(70 40) \htext(76
42){\tiny $0$}

\move(80 40)\lvec(90 50)\lvec(90 40)\lvec(80 40) \htext(86
42){\tiny $1$}

\move(60 50)\lvec(70 50)\lvec(70 60)\lvec(60 60)\lvec(60
50)\htext(63 53){$2$}

\move(70 50)\lvec(80 50)\lvec(80 60)\lvec(70 60)\lvec(70
50)\htext(73 53){$2$}

\move(80 50)\lvec(90 50)\lvec(90 60)\lvec(80 60)\lvec(80
50)\htext(83 53){$2$}
\move(60 60)\lvec(70 60)\lvec(70 65)\lvec(60 65)\lvec(60
60)\htext(63 61){\tiny $3$}

\move(70 60)\lvec(80 60)\lvec(80 65)\lvec(70 65)\lvec(70
60)\htext(73 61){\tiny $3$}

\move(80 60)\lvec(90 60)\lvec(90 65)\lvec(80 65)\lvec(80
60)\htext(83 61){\tiny $3$}
\move(60 60)\lvec(70 60)\lvec(70 70)\lvec(60 70)\lvec(60
60)\htext(63 66){\tiny $3$}

\move(70 60)\lvec(80 60)\lvec(80 70)\lvec(70 70)\lvec(70
60)\htext(73 66){\tiny $3$}

\move(80 60)\lvec(90 60)\lvec(90 70)\lvec(80 70)\lvec(80
60)\htext(83 66){\tiny $3$}
\move(60 70)\lvec(70 70)\lvec(70 80)\lvec(60 80)\lvec(60
70)\htext(63 73){$2$}

\move(70 70)\lvec(80 70)\lvec(80 80)\lvec(70 80)\lvec(70
70)\htext(73 73){$2$}

\move(80 70)\lvec(90 70)\lvec(90 80)\lvec(80 80)\lvec(80
70)\htext(83 73){$2$}
\move(60 80)\lvec(70 90)\lvec(70 80)\lvec(60 80)\htext(66
82){\tiny $1$}

\move(70 80)\lvec(80 90)\lvec(80 80)\lvec(70 80)\htext(76
82){\tiny $0$}

\move(80 80)\lvec(90 90)\lvec(90 80)\lvec(80 80)\htext(86
82){\tiny $1$}

\end{texdraw}} : $i$-component of type $\rm III_2^-$.\vskip 5mm
}
\end{ex}

Let $Y$ be a proper Young wall with the decomposition $Y=(Y_0,
Y_1, \cdots, Y_r)$ into $i$-components. To each $Y_k$, we
associate a $U_{(i)}$-module $V_k$ as follows. Then we will view
$Y$ as $Y_0\otimes Y_1\otimes\cdots\otimes Y_r$ inside $V_0\otimes
V_1\otimes\cdots\otimes V_r$.

If $Y_k$ is a trivial $i$-component, then we associate the trivial
representation $V_k=U=\mathbb{Q}(q)u$, and we identify $Y_k$ with
$u$. If $Y_k$ is an $i$-component of type $\rm III_{\infty}$ or
$\rm III_0$, then we associate the 2-dimensional representation
$V_k=V=\mathbb{Q}(q)v_0\oplus\mathbb{Q}(q)v_1$, where the
$U_{(i)}$-module action is given by (\ref{I+}). We identify $Y_k$
with a basis element of $V$ as follows:\vskip 5mm

\begin{center}
$v_0\leftrightarrow$\raisebox{-0.4\height}{\begin{texdraw}\drawdim
em \setunitscale 0.1 \linewd 0.5

\htext(-4 3){$\cdots$} \move(10 0)\lvec(20 10)\lvec(20 0)\lvec(10
0)\lfill f:0.8

\move(20 0)\lvec(30 10)\lvec(30 0)\lvec(20 0)\lfill f:0.8

\move(30 0)\lvec(40 10)\lvec(40 0)\lvec(30 0)\lfill f:0.8

\move(40 0)\lvec(50 10)\lvec(50 0)\lvec(40 0)\lfill f:0.8

\htext(46 2){\tiny $\cdot$}\htext(41 6){\tiny $*$}

\move(-10 0)\lvec(50 0)
\end{texdraw}}
\hskip 1cm
$v_1\leftrightarrow$\raisebox{-0.4\height}{\begin{texdraw}\drawdim
em \setunitscale 0.1 \linewd 0.5

\htext(-5 3){$\cdots$} \move(10 0)\lvec(20 10)\lvec(20 0)\lvec(10
0)\lfill f:0.8

\move(20 0)\lvec(30 10)\lvec(30 0)\lvec(20 0)\lfill f:0.8

\move(30 0)\lvec(40 10)\lvec(40 0)\lvec(30 0)\lfill f:0.8

\move(40 0)\lvec(50 10)\lvec(50 0)\lvec(40 0)\lfill f:0.8

\htext(46 2){\tiny $\cdot$}\htext(41 6){\tiny $*$}

\move(-10 0)\lvec(50 0)\lvec(50 10)\lvec(40 10)
\end{texdraw}}
\end{center}\vskip 5mm

\begin{center}
$v_0\leftrightarrow$\raisebox{-0.65\height}{
\begin{texdraw}
\drawdim em \setunitscale 0.1 \linewd 0.5

\move(0 -10)\lvec(0 10)\lvec(10 10)\lvec(10 -10)\lvec(0 -10)

\htext(1 16){\tiny $*$}
\end{texdraw}}
\hskip 5mm $v_1\leftrightarrow$\raisebox{-0.6\height}{
\begin{texdraw}
\drawdim em \setunitscale 0.1 \linewd 0.5

\move(0 -10)\lvec(0 10)\lvec(10 10)\lvec(10 -10)\lvec(0 -10)

\move(0 10)\lvec(0 20)\lvec(10 20)\lvec(0 10)

\htext(1 16){\tiny $*$}
\end{texdraw}}\hskip 10mm
$v_0\leftrightarrow$\raisebox{-0.85\height}{
\begin{texdraw}
\drawdim em \setunitscale 0.1 \linewd 0.5

\move(0 -10)\lvec(0 10)\lvec(10 10)\lvec(10 -10)\lvec(0 -10)

\htext(6 12){\tiny $*$}
\end{texdraw}}\hskip 5mm
$v_1\leftrightarrow$\raisebox{-0.6\height}{
\begin{texdraw}
\drawdim em \setunitscale 0.1 \linewd 0.5

\move(0 -10)\lvec(0 10)\lvec(10 10)\lvec(10 -10)\lvec(0 -10)

\move(0 10)\lvec(10 20)\lvec(10 10)\lvec(0 10)

\htext(5 12){\tiny $*$}
\end{texdraw}}
\end{center}\vskip 5mm

\begin{center}
$v_0\leftrightarrow$\raisebox{-0.6\height}{
\begin{texdraw}
\drawdim em \setunitscale 0.1 \linewd 0.5

\move(0 -10)\lvec(0 10)\lvec(10 10)\lvec(10 -10)\lvec(0 -10)

\move(0 10)\lvec(10 20)\lvec(10 10)\lvec(0 10)

\htext(6 12){\tiny $\cdot$}
\end{texdraw}}
\hskip 5mm $v_1\leftrightarrow$\raisebox{-0.6\height}{
\begin{texdraw}
\drawdim em \setunitscale 0.1 \linewd 0.5

\move(0 -10)\lvec(0 10)\lvec(10 10)\lvec(10 -10)\lvec(0 -10)

\move(0 10)\lvec(0 20)\lvec(10 20)\lvec(0 10)

\move(0 10)\lvec(10 20)\lvec(10 10)\lvec(0 10)

\htext(1 16){\tiny $*$}\htext(6 12){$\cdot$}
\end{texdraw}}\hskip 10mm
$v_0\leftrightarrow$\raisebox{-0.6\height}{
\begin{texdraw}
\drawdim em \setunitscale 0.1 \linewd 0.5

\move(0 -10)\lvec(0 10)\lvec(10 10)\lvec(10 -10)\lvec(0 -10)

\move(0 10)\lvec(0 20)\lvec(10 20)\lvec(0 10)

\htext(2 16){\tiny $\cdot$}
\end{texdraw}}\hskip 5mm
$v_1\leftrightarrow$\raisebox{-0.6\height}{
\begin{texdraw}
\drawdim em \setunitscale 0.1 \linewd 0.5

\move(0 -10)\lvec(0 10)\lvec(10 10)\lvec(10 -10)\lvec(0 -10)

\move(0 10)\lvec(10 20)\lvec(10 10)\lvec(0 10)

\move(0 10)\lvec(0 20)\lvec(10 20)\lvec(0 10)

\htext(2 16){\tiny $\cdot$}

\htext(5 12){\tiny $*$}
\end{texdraw}}\quad .\vskip 5mm
\end{center}

If $Y_k$ is an $i$-component of type $\rm III_1$, then we
associate the 5-dimensional representation
\begin{equation*}
V_k=W_1=\mathbb{Q}(q)w_0\oplus\mathbb{Q}(q)w_1\oplus\mathbb{Q}(q)w_2\oplus
\mathbb{Q}(q)u\oplus\mathbb{Q}(q)u',
\end{equation*}
where the $U_{(i)}$-module action is given by
\begin{equation}
\begin{split}
&K_iw_0=q_i^2w_0,\ e_i w_0=0, \ f_iw_0=w_1+q_iu, \\
&K_iw_1=w_1,\ e_iw_1=q_i^{-1}w_0,\ f_iw_1=q_i^{-1}w_2, \\
&K_iw_2=q_i^{-2}w_2,\ e_iw_2=w_1+q_iu,\ f_iw_2=0, \\
&K_iu=u,\ e_iu=w_0,\ f_iu=w_2, \\
&K_iu'=u',\ e_iu'=-q_iw_0,\ f_iu'=-q_iw_2.
\end{split}
\end{equation}

We identify the $i$-component $Y_k$ with a basis element of $W_1$
as follows:\vskip 5mm
\begin{center}
$w_0\leftrightarrow$\raisebox{-0.4\height}{
\begin{texdraw}
\fontsize{6}{6}\drawdim em \setunitscale 0.1 \linewd 0.5

\move(0 -10)\lvec(0 10)\lvec(20 10)\lvec(20 -10)\lvec(0 -10)

\move(10 -10)\lvec(10 10)

\move(10 10)\lvec(20 20)\lvec(20 10)\lvec(10 10)

\htext(16 12){$\cdot$}

\end{texdraw}}\hskip 7mm
$w_1\leftrightarrow$\raisebox{-0.4\height}{
\begin{texdraw}
\fontsize{6}{6}\drawdim em \setunitscale 0.1 \linewd 0.5

\move(10 10)\lvec(20 20)\lvec(20 10)\lvec(10 10)

\htext(16 12){$\cdot$}

\move(0 10)\lvec(10 20)\lvec(10 10)\lvec(0 10)

\htext(5 12){$*$}

\move(0 -10)\lvec(0 10)\lvec(20 10)\lvec(20 -10)\lvec(0 -10)

\move(10 -10)\lvec(10 10)

\end{texdraw}}\hskip 7mm
$w_2\leftrightarrow$\raisebox{-0.4\height}{
\begin{texdraw}
\fontsize{6}{6}\drawdim em \setunitscale 0.1 \linewd 0.5

\move(10 10)\lvec(20 20)\lvec(20 10)\lvec(10 10)

\htext(16 12){$\cdot$}

\move(0 10)\lvec(10 20)\lvec(10 10)\lvec(0 10)

\htext(5 12){$*$}

\move(10 10)\lvec(20 20)\lvec(10 20)\lvec(10 10)

\htext(11 15){$*$}

\move(0 -10)\lvec(0 10)\lvec(20 10)\lvec(20 -10)\lvec(0 -10)

\move(10 -10)\lvec(10 10)

\end{texdraw}}\vskip 5mm

$u\leftrightarrow$\raisebox{-0.4\height}{
\begin{texdraw}
\fontsize{6}{6}\drawdim em \setunitscale 0.1 \linewd 0.5

\move(10 10)\lvec(20 20)\lvec(20 10)\lvec(10 10)

\htext(16 12){$\cdot$}

\move(10 10)\lvec(20 20)\lvec(10 20)\lvec(10 10)

\htext(11 16){$*$}

\move(0 -10)\lvec(0 10)\lvec(20 10)\lvec(20 -10)\lvec(0 -10)

\move(10 -10)\lvec(10 10)

\end{texdraw}}\hskip 7mm
$u'\leftrightarrow$\raisebox{-0.4\height}{
\begin{texdraw}
\fontsize{6}{6}\drawdim em \setunitscale 0.1 \linewd 0.5

\move(0 10)\lvec(10 20)\lvec(0 20)\lvec(0 10)

\htext(2 16){$\cdot$}

\move(10 10)\lvec(20 20)\lvec(10 20)\lvec(10 10)

\htext(11 16){$*$}

\move(0 -10)\lvec(0 10)\lvec(20 10)\lvec(20 -10)\lvec(0 -10)

\move(10 -10)\lvec(10 10)

\end{texdraw}}
\end{center}\vskip 10mm

\begin{center}
$w_0\leftrightarrow$\raisebox{-0.4\height}{
\begin{texdraw}
\fontsize{6}{6}\drawdim em \setunitscale 0.1 \linewd 0.5

\move(0 -10)\lvec(0 10)\lvec(20 10)\lvec(20 -10)\lvec(0 -10)

\move(10 -10)\lvec(10 10)

\move(10 10)\lvec(20 20)\lvec(10 20)\lvec(10 10)

\htext(12 16){$\cdot$}

\end{texdraw}}\hskip 7mm
$w_1\leftrightarrow$\raisebox{-0.4\height}{
\begin{texdraw}
\fontsize{6}{6}\drawdim em \setunitscale 0.1 \linewd 0.5

\move(0 10)\lvec(10 20)\lvec(0 20)\lvec(0 10)

\htext(2 16){$*$}

\move(10 10)\lvec(20 20)\lvec(10 20)\lvec(10 10)

\htext(12 16){$\cdot$}

\move(0 -10)\lvec(0 10)\lvec(20 10)\lvec(20 -10)\lvec(0 -10)

\move(10 -10)\lvec(10 10)

\end{texdraw}}\hskip 7mm
$w_2\leftrightarrow$\raisebox{-0.4\height}{
\begin{texdraw}
\fontsize{6}{6}\drawdim em \setunitscale 0.1 \linewd 0.5

\move(10 10)\lvec(20 20)\lvec(20 10)\lvec(10 10)

\htext(15 12){$*$}

\move(0 10)\lvec(10 20)\lvec(0 20)\lvec(0 10)

\htext(2 16){$*$}

\move(10 10)\lvec(20 20)\lvec(10 20)\lvec(10 10)

\htext(12 16){$\cdot$}

\move(0 -10)\lvec(0 10)\lvec(20 10)\lvec(20 -10)\lvec(0 -10)

\move(10 -10)\lvec(10 10)

\end{texdraw}}\vskip 5mm

$u\leftrightarrow$\raisebox{-0.4\height}{
\begin{texdraw}
\fontsize{6}{6}\drawdim em \setunitscale 0.1 \linewd 0.5

\move(10 10)\lvec(20 20)\lvec(20 10)\lvec(10 10)

\htext(15 12){$*$}

\move(10 10)\lvec(20 20)\lvec(10 20)\lvec(10 10)

\htext(12 16){$\cdot$}

\move(0 -10)\lvec(0 10)\lvec(20 10)\lvec(20 -10)\lvec(0 -10)

\move(10 -10)\lvec(10 10)

\end{texdraw}}\hskip 7mm
$u'\leftrightarrow$\raisebox{-0.4\height}{
\begin{texdraw}
\fontsize{6}{6}\drawdim em \setunitscale 0.1 \linewd 0.5

\move(10 10)\lvec(20 20)\lvec(20 10)\lvec(10 10)

\htext(15 12){$*$}

\move(0 10)\lvec(10 20)\lvec(10 10)\lvec(0 10)

\htext(6 12){$\cdot$}

\move(0 -10)\lvec(0 10)\lvec(20 10)\lvec(20 -10)\lvec(0 -10)

\move(10 -10)\lvec(10 10)

\end{texdraw}}\quad.
\end{center}
\vskip 5mm

The $U_{(i)}$-module $W_1$ is decomposed as $W_1\cong V(2)\oplus
V(0)\oplus V(0)$, where
\begin{equation}
\begin{split}
&V(2)\cong \mathbb{Q}(q)w_0\oplus \mathbb{Q}(q)(w_1 +
q_iu)\oplus \mathbb{Q}(q)w_2, \\
&V(0) \cong \mathbb{Q}(q)(u-q_iw_1)\cong
\mathbb{Q}(q)(u'+q_i^2w_1).
\end{split}
\end{equation}
Hence, the crystal basis $(L,B)$ of $W_1$ is given by
\begin{equation}
\begin{split}
L=&\mathbb{A}_0w_0\oplus \mathbb{A}_0(w_1+q_iu)\oplus
\mathbb{A}_0 w_2 \\
& \oplus \mathbb{A}_0(u-q_iw_1)\oplus
\mathbb{A}_0(u'+q_i^2w_1), \\
B=&\{\,\overline{w_0},\overline{w_1},\overline{w_2},\overline{u},\overline{u'}\,\}
\end{split}
\end{equation}
with the crystal graph
$$\overline{w_0}\stackrel{i}{\longrightarrow}\overline{w_1}
\stackrel{i}{\longrightarrow}\overline{w_2}\hskip 1cm
\overline{u}\hskip 1cm\overline{u'}.$$

If $Y_k$ is an $i$-component of type ${\rm III}_{2l}^+$, then we
associate the $4$-dimensional representation
\begin{equation*}
V_k=W_{2l}^{\pm}=\mathbb{Q}(q)w_0\oplus \mathbb{Q}(q)w_1\oplus
\mathbb{Q}(q)w_0'\oplus \mathbb{Q}(q)w_1',
\end{equation*}
where the $U_{(i)}$-module action is given by
\begin{equation}
\begin{split}
&K_iw_0=q_iw_0, \ K_iw_1=q_i^{-1}w_1, \\
&K_iw_0'=q_iw_0', \ K_iw_1'=q_i^{-1}w_1', \\
&e_iw_0=0, \ e_iw_1=w_0\pm q_i^{2l}w_0', \\
&e_iw_0'=0, \ e_iw_1'=w_0', \\
&f_iw_0=w_1\mp q_i^{2l}w_1', \ f_iw_1=0, \\
&f_iw_0'=w_1', \ f_iw_1'=0.
\end{split}
\end{equation}

We identify the $i$-component $Y_k$ with a basis element of
$W_{2l}^{\pm}$ as follows :\vskip 5mm

$W_{2l}^+$
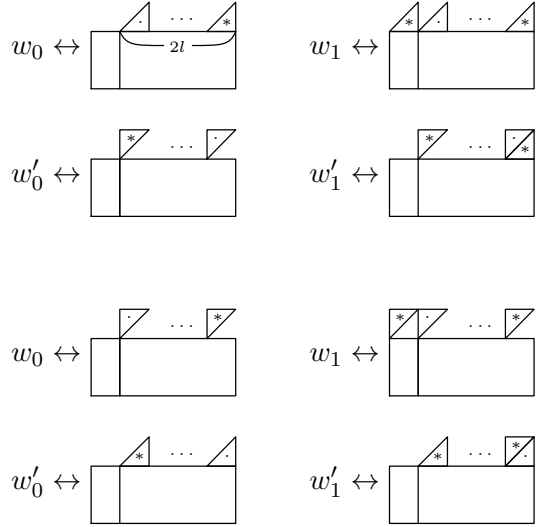
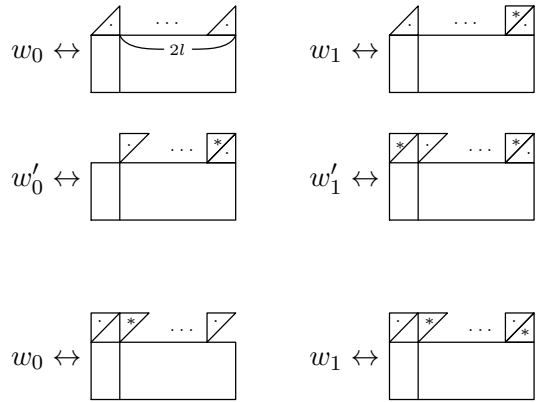
\begin{center}
$w_0\leftrightarrow$\raisebox{-0.4\height}{
\begin{texdraw}
\fontsize{6}{6}\drawdim em \setunitscale 0.1 \linewd 0.5

\move(0 -10)\lvec(0 10)\lvec(50 10)\lvec(50 -10)\lvec(0 -10)

\move(10 10)\lvec(10 -10)

\move(10 10)\lvec(20 20)\lvec(20 10)\lvec(10 10)

\htext(16 12){$\cdot$}

\move(40 10)\lvec(50 20)\lvec(50 10)\lvec(40 10)

\htext(45 12){$*$}

\htext(27 13){$\cdots$}

\move(10 10)\clvec(12 5)(15 5)(25 5)

\move(50 10)\clvec(48 5)(45 5)(35 5)

\htext(27 3){${2l}$}

\end{texdraw}}\hskip 10mm
$w_1\leftrightarrow$\raisebox{-0.4\height}{
\begin{texdraw}
\fontsize{6}{6}\drawdim em \setunitscale 0.1 \linewd 0.5

\move(0 -10)\lvec(0 10)\lvec(50 10)\lvec(50 -10)\lvec(0 -10)

\move(10 10)\lvec(10 -10)

\move(10 10)\lvec(20 20)\lvec(20 10)\lvec(10 10)

\htext(16 12){$\cdot$}

\move(40 10)\lvec(50 20)\lvec(50 10)\lvec(40 10)

\htext(45 12){$*$}

\htext(27 13){$\cdots$}

\move(0 10)\lvec(10 20)\lvec(10 10)\lvec(0 10)

\htext(5 12){$*$}

\end{texdraw}}\vskip 5mm

$w_0'\leftrightarrow$\raisebox{-0.4\height}{
\begin{texdraw}
\fontsize{6}{6}\drawdim em \setunitscale 0.1 \linewd 0.5

\move(0 -10)\lvec(0 10)\lvec(50 10)\lvec(50 -10)\lvec(0 -10)

\move(10 10)\lvec(10 -10)

\move(10 10)\lvec(10 -10)

\move(10 10)\lvec(20 20)\lvec(10 20)\lvec(10 10)

\htext(12 16){$*$}

\htext(27 13){$\cdots$}

\move(40 10)\lvec(50 20)\lvec(40 20)\lvec(40 10)

\htext(42 16){$\cdot$}

\end{texdraw}}\hskip 10mm
$w_1'\leftrightarrow$\raisebox{-0.4\height}{
\begin{texdraw}
\fontsize{6}{6}\drawdim em \setunitscale 0.1 \linewd 0.5

\move(0 -10)\lvec(0 10)\lvec(50 10)\lvec(50 -10)\lvec(0 -10)

\move(10 10)\lvec(10 -10)

\move(10 10)\lvec(20 20)\lvec(10 20)\lvec(10 10)

\htext(12 16){$*$}

\move(40 10)\lvec(50 20)\lvec(50 10)\lvec(40 10)

\htext(45 12){$*$}

\htext(27 13){$\cdots$}

\move(40 10)\lvec(50 20)\lvec(40 20)\lvec(40 10)

\htext(42 16){$\cdot$}

\end{texdraw}}
\end{center}\vskip 10mm

\begin{center}
$w_0\leftrightarrow$\raisebox{-0.4\height}{
\begin{texdraw}
\fontsize{6}{6}\drawdim em \setunitscale 0.1 \linewd 0.5

\move(0 -10)\lvec(0 10)\lvec(50 10)\lvec(50 -10)\lvec(0 -10)

\move(10 10)\lvec(10 -10)

\move(10 10)\lvec(10 -10)

\move(10 10)\lvec(20 20)\lvec(10 20)\lvec(10 10)

\htext(12 16){$\cdot$}

\htext(27 13){$\cdots$}

\move(40 10)\lvec(50 20)\lvec(40 20)\lvec(40 10)

\htext(42 16){$*$}

\end{texdraw}}\hskip 10mm
$w_1\leftrightarrow$\raisebox{-0.4\height}{
\begin{texdraw}
\fontsize{6}{6}\drawdim em \setunitscale 0.1 \linewd 0.5

\move(0 -10)\lvec(0 10)\lvec(50 10)\lvec(50 -10)\lvec(0 -10)

\move(10 10)\lvec(10 -10)

\move(10 10)\lvec(10 -10)

\move(10 10)\lvec(20 20)\lvec(10 20)\lvec(10 10)

\htext(12 16){$\cdot$}

\htext(27 13){$\cdots$}

\move(40 10)\lvec(50 20)\lvec(40 20)\lvec(40 10)

\htext(42 16){$*$}

\move(0 10)\lvec(10 20)\lvec(0 20)\lvec(0 10)

\htext(2 16){$*$}

\end{texdraw}}\vskip 5mm

$w_0'\leftrightarrow$\raisebox{-0.4\height}{
\begin{texdraw}
\fontsize{6}{6}\drawdim em \setunitscale 0.1 \linewd 0.5

\move(0 -10)\lvec(0 10)\lvec(50 10)\lvec(50 -10)\lvec(0 -10)

\move(10 10)\lvec(10 -10)

\move(10 10)\lvec(20 20)\lvec(20 10)\lvec(10 10)

\htext(15 12){$*$}

\move(40 10)\lvec(50 20)\lvec(50 10)\lvec(40 10)

\htext(46 12){$\cdot$}

\htext(27 13){$\cdots$}

\end{texdraw}}\hskip 10mm
$w_1'\leftrightarrow$\raisebox{-0.4\height}{
\begin{texdraw}
\fontsize{6}{6}\drawdim em \setunitscale 0.1 \linewd 0.5

\move(0 -10)\lvec(0 10)\lvec(50 10)\lvec(50 -10)\lvec(0 -10)

\move(10 10)\lvec(10 -10)

\move(10 10)\lvec(20 20)\lvec(20 10)\lvec(10 10)

\htext(15 12){$*$}

\move(40 10)\lvec(50 20)\lvec(50 10)\lvec(40 10)

\htext(46 12){$\cdot$}

\move(40 10)\lvec(50 20)\lvec(40 20)\lvec(40 10)

\htext(42 16){$*$}

\htext(27 13){$\cdots$}

\end{texdraw}}
\end{center}\vskip 10mm

$W_{2l}^-$
\begin{center}
$w_0\leftrightarrow$\raisebox{-0.4\height}{
\begin{texdraw}
\fontsize{6}{6}\drawdim em \setunitscale 0.1 \linewd 0.5

\move(0 -10)\lvec(0 10)\lvec(50 10)\lvec(50 -10)\lvec(0 -10)

\move(10 10)\lvec(10 -10)

\move(0 10)\lvec(10 20)\lvec(10 10)\lvec(0 10)

\htext(6 12){$\cdot$}

\move(40 10)\lvec(50 20)\lvec(50 10)\lvec(40 10)

\htext(46 12){$\cdot$}

\htext(22 13){$\cdots$}

\move(10 10)\clvec(12 5)(20 5)(25 5)

\move(50 10)\clvec(48 5)(40 5)(35 5)

\htext(27 3){$2l$}

\end{texdraw}}\hskip 10mm
$w_1\leftrightarrow$\raisebox{-0.4\height}{
\begin{texdraw}
\fontsize{6}{6}\drawdim em \setunitscale 0.1 \linewd 0.5

\move(0 -10)\lvec(0 10)\lvec(50 10)\lvec(50 -10)\lvec(0 -10)

\move(10 10)\lvec(10 -10)

\move(0 10)\lvec(10 20)\lvec(10 10)\lvec(0 10)

\htext(6 12){$\cdot$}

\move(40 10)\lvec(50 20)\lvec(50 10)\lvec(40 10)

\htext(46 12){$\cdot$}

\htext(22 13){$\cdots$}

\move(40 10)\lvec(50 20)\lvec(40 20)\lvec(40 10)

\htext(42 16){$*$}

\end{texdraw}}\vskip 5mm

$w_0'\leftrightarrow$\raisebox{-0.4\height}{
\begin{texdraw}
\fontsize{6}{6}\drawdim em \setunitscale 0.1 \linewd 0.5

\move(0 -10)\lvec(0 10)\lvec(50 10)\lvec(50 -10)\lvec(0 -10)
\move(10 10)\lvec(10 -10)

\move(10 10)\lvec(20 20)\lvec(10 20)\lvec(10 10)

\htext(12 15){$\cdot$}

\move(40 10)\lvec(50 20)\lvec(50 10)\lvec(40 10)

\htext(46 12){$\cdot$}

\htext(27 13){$\cdots$}

\move(40 10)\lvec(50 20)\lvec(40 20)\lvec(40 10)

\htext(42 16){$*$}

\end{texdraw}}\hskip 10mm
$w_1'\leftrightarrow$\raisebox{-0.4\height}{
\begin{texdraw}
\fontsize{6}{6}\drawdim em \setunitscale 0.1 \linewd 0.5

\move(0 -10)\lvec(0 10)\lvec(50 10)\lvec(50 -10)\lvec(0 -10)

\move(10 10)\lvec(10 -10)

\move(0 10)\lvec(10 20)\lvec(0 20)\lvec(0 10)

\htext(2 15){$*$}

\move(10 10)\lvec(20 20)\lvec(10 20)\lvec(10 10)

\htext(12 15){$\cdot$}

\move(40 10)\lvec(50 20)\lvec(50 10)\lvec(40 10)

\htext(47 12){$\cdot$}

\htext(27 13){$\cdots$}

\move(40 10)\lvec(50 20)\lvec(40 20)\lvec(40 10)

\htext(42 16){$*$}

\end{texdraw}}
\end{center}\vskip 10mm

\begin{center}
$w_0\leftrightarrow$\raisebox{-0.4\height}{
\begin{texdraw}
\fontsize{6}{6}\drawdim em \setunitscale 0.1 \linewd 0.5

\move(0 -10)\lvec(0 10)\lvec(50 10)\lvec(50 -10)\lvec(0 -10)

\move(10 10)\lvec(10 -10)

\move(0 10)\lvec(10 20)\lvec(0 20)\lvec(0 10)

\htext(2 16){$\cdot$}

\move(10 10)\lvec(20 20)\lvec(10 20)\lvec(10 10)

\htext(12 16){$*$}

\htext(27 13){$\cdots$}

\move(40 10)\lvec(50 20)\lvec(40 20)\lvec(40 10) \htext(42
16){$\cdot$}

\end{texdraw}}\hskip 10mm
$w_1\leftrightarrow$\raisebox{-0.4\height}{
\begin{texdraw}
\fontsize{6}{6}\drawdim em \setunitscale 0.1 \linewd 0.5

\move(0 -10)\lvec(0 10)\lvec(50 10)\lvec(50 -10)\lvec(0 -10)

\move(10 10)\lvec(10 -10)

\move(0 10)\lvec(10 20)\lvec(0 20)\lvec(0 10)

\htext(2 16){$\cdot$}

\move(10 10)\lvec(20 20)\lvec(10 20)\lvec(10 10)

\htext(12 16){$*$}

\move(40 10)\lvec(50 20)\lvec(50 10)\lvec(40 10)

\htext(45 12){$*$}

\htext(27 13){$\cdots$}

\move(40 10)\lvec(50 20)\lvec(40 20)\lvec(40 10)

\htext(42 16){$\cdot$}

\end{texdraw}}\vskip 5mm

$w_0'\leftrightarrow$\raisebox{-0.4\height}{
\begin{texdraw}
\fontsize{6}{6}\drawdim em \setunitscale 0.1 \linewd 0.5

\move(0 -10)\lvec(0 10)\lvec(50 10)\lvec(50 -10)\lvec(0 -10)

\move(10 10)\lvec(10 -10)

\move(10 10)\lvec(20 20)\lvec(20 10)\lvec(10 10)

\htext(16 12){$\cdot$}

\move(40 10)\lvec(50 20)\lvec(50 10)\lvec(40 10)

\htext(45 12){$*$}

\move(40 10)\lvec(50 20)\lvec(40 20)\lvec(40 10)

\htext(42 16){$\cdot$}

\htext(27 13){$\cdots$}

\end{texdraw}}\hskip 10mm
$w_1'\leftrightarrow$\raisebox{-0.4\height}{
\begin{texdraw}
\fontsize{6}{6}\drawdim em \setunitscale 0.1 \linewd 0.5

\move(0 -10)\lvec(0 10)\lvec(50 10)\lvec(50 -10)\lvec(0 -10)

\move(10 10)\lvec(10 -10)

\move(10 10)\lvec(20 20)\lvec(20 10)\lvec(10 10)

\htext(17 12){$\cdot$}

\move(40 10)\lvec(50 20)\lvec(50 10)\lvec(40 10)

\htext(45 12){$*$}

\htext(27 13){$\cdots$}

\move(0 10)\lvec(10 20)\lvec(10 10)\lvec(0 10)

\htext(5 12){$*$}

\move(40 10)\lvec(50 20)\lvec(40 20)\lvec(40 10)

\htext(42 16){$\cdot$}

\end{texdraw}}\quad .
\end{center}\vskip 5mm

The $U_{(i)}$-module $W_{2l}^{\pm}$ is decomposed as
$W_{2l}^{\pm}\cong V(1)\oplus V(1)$, where
\begin{equation}
\begin{split}
V(2)&\cong \mathbb{Q}(q)w_0\oplus \mathbb{Q}(q)(w_1\mp
q_i^{2l}w_1'), \\
&\cong \mathbb{Q}(q)w_0'\oplus \mathbb{Q}(q)w_1'.
\end{split}
\end{equation}
Hence, the crystal basis $(L,B)$ of $W_{2l}^{\pm}$ is given by
\begin{equation}
\begin{split}
L&=\mathbb{A}_0w_0\oplus \mathbb{A}_0(w_1\mp q_i^{2l}w_1')\oplus
\mathbb{A}_0w_0'\oplus \mathbb{A}_0w_1', \\
B&=\{\,\overline{w_0},\overline{w_1},\overline{w_0'},\overline{w_1'}\,\}
\end{split}
\end{equation}
with the crystal graph
$$\overline{w_0}\stackrel{i}{\longrightarrow}\overline{w_1}\hskip
1cm \overline{w_0'}\stackrel{i}{\longrightarrow}\overline{w_1'}.$$

If $Y_k$ is an $i$-component of type ${\rm III}_{2l+1}$, then we
associate the $6$-dimensional representation
\begin{equation*}
V_k=W_{2l+1}=\mathbb{Q}(q)w_0\oplus \mathbb{Q}(q)w_1\oplus
\mathbb{Q}(q)w_2\oplus \mathbb{Q}(q)u\oplus \mathbb{Q}(q)u'\oplus
\mathbb{Q}(q)u'',
\end{equation*}
where the $U_{(i)}$-module action is given by
\begin{equation}
\begin{split}
&K_iw_0=q_i^2w_0, \ e_iw_0=0, \ f_iw_0=w_1+q_iu+q_i^{2l+1}u'',
\\
&K_iw_1=w_1, \ e_iw_1=q_i^{-1}w_0, \ f_iw_1=q_i^{-1}w_2, \\
&K_iw_2=q_i^{-2}w_2, \ e_iw_2=w_1+q_iu+q_i^{2l+1}u'', \
f_iw_2=0, \\
&K_iu=u, \ e_iu=w_0, \ f_iu=w_2, \\
&K_iu'=u', \ e_iu'=-q_i^{2l+1}w_0, \ f_iu'=-q_i^{2l+1}w_2,
\\
&K_iu''=u'', \ e_iu''=0, \ f_iu''=0.
\end{split}
\end{equation}
We identify the $i$-component $Y_k$ with a basis element of
$W_{2l+1}$ as follows:\vskip 5mm

\begin{center}
$w_0\leftrightarrow$\raisebox{-0.4\height}{
\begin{texdraw}
\fontsize{6}{6}\drawdim em \setunitscale 0.1 \linewd 0.5

\move(0 -10)\lvec(0 10)\lvec(50 10)\lvec(50 -10)\lvec(0 -10)

\move(10 10)\lvec(10 -10)

\move(10 10)\lvec(20 20)\lvec(20 10)\lvec(10 10)

\htext(16 12){$\cdot$}

\move(40 10)\lvec(50 20)\lvec(50 10)\lvec(40 10)

\htext(46 12){$\cdot$}

\htext(27 13){$\cdots$}

\move(10 10)\clvec(12 5)(15 5)(20 5)

\move(50 10)\clvec(48 5)(45 5)(40 5)

\htext(23 3){$2l+1$}

\end{texdraw}}\hskip 10mm
$w_1\leftrightarrow$\raisebox{-0.4\height}{
\begin{texdraw}
\fontsize{6}{6}\drawdim em \setunitscale 0.1 \linewd 0.5

\move(0 -10)\lvec(0 10)\lvec(50 10)\lvec(50 -10)\lvec(0 -10)

\move(10 10)\lvec(10 -10)

\move(10 10)\lvec(20 20)\lvec(20 10)\lvec(10 10)

\htext(16 12){$\cdot$}

\move(40 10)\lvec(50 20)\lvec(50 10)\lvec(40 10)

\htext(46 12){$\cdot$}

\htext(27 13){$\cdots$}

\move(0 10)\lvec(10 20)\lvec(10 10)\lvec(0 10)

\htext(5 12){$*$}

\end{texdraw}}\vskip 5mm

$w_2\leftrightarrow$\raisebox{-0.4\height}{
\begin{texdraw}
\fontsize{6}{6}\drawdim em \setunitscale 0.1 \linewd 0.5

\move(0 -10)\lvec(0 10)\lvec(50 10)\lvec(50 -10)\lvec(0 -10)

\move(10 10)\lvec(10 -10)

\move(10 10)\lvec(20 20)\lvec(20 10)\lvec(10 10)

\htext(16 12){$\cdot$}

\move(40 10)\lvec(50 20)\lvec(50 10)\lvec(40 10)

\htext(46 12){$\cdot$}

\htext(27 13){$\cdots$}

\move(0 10)\lvec(10 20)\lvec(10 10)\lvec(0 10)

\htext(5 12){$*$}

\move(40 10)\lvec(50 20)\lvec(40 20)\lvec(40 10)

\htext(42 16){$*$}

\end{texdraw}}\hskip 10mm
$u\leftrightarrow$\raisebox{-0.4\height}{
\begin{texdraw}
\fontsize{6}{6}\drawdim em \setunitscale 0.1 \linewd 0.5

\move(0 -10)\lvec(0 10)\lvec(50 10)\lvec(50 -10)\lvec(0 -10)

\move(10 10)\lvec(10 -10)

\move(10 10)\lvec(20 20)\lvec(20 10)\lvec(10 10)

\htext(16 12){$\cdot$}

\move(40 10)\lvec(50 20)\lvec(50 10)\lvec(40 10)

\htext(46 12){$\cdot$}

\htext(27 13){$\cdots$}

\move(40 10)\lvec(50 20)\lvec(40 20)\lvec(40 10)

\htext(42 16){$*$}

\end{texdraw}}\vskip 5mm

$u'\leftrightarrow$\raisebox{-0.4\height}{
\begin{texdraw}
\fontsize{6}{6}\drawdim em \setunitscale 0.1 \linewd 0.5

\move(0 -10)\lvec(0 10)\lvec(50 10)\lvec(50 -10)\lvec(0 -10)

\move(10 10)\lvec(10 -10)

\move(0 10)\lvec(10 20)\lvec(0 20)\lvec(0 10)

\htext(2 16){$\cdot$}

\htext(20 13){$\cdots$}

\move(40 10)\lvec(50 20)\lvec(40 20)\lvec(40 10)

\htext(42 16){$*$}

\end{texdraw}}\hskip 1cm
$u''\leftrightarrow$\raisebox{-0.4\height}{
\begin{texdraw}
\fontsize{6}{6}\drawdim em \setunitscale 0.1 \linewd 0.5

\move(0 -10)\lvec(0 10)\lvec(50 10)\lvec(50 -10)\lvec(0 -10)

\move(10 10)\lvec(10 -10)

\move(10 10)\lvec(20 20)\lvec(10 20)\lvec(10 10)

\htext(12 16){$*$}

\move(40 10)\lvec(50 20)\lvec(50 10)\lvec(40 10)

\htext(46 12){$\cdot$}

\htext(27 13){$\cdots$}

\move(40 10)\lvec(50 20)\lvec(40 20)\lvec(40 10)

\htext(42 16){$*$}

\end{texdraw}}
\end{center}\vskip 10mm

\begin{center}
$w_0\leftrightarrow$\raisebox{-0.4\height}{
\begin{texdraw}
\fontsize{6}{6}\drawdim em \setunitscale 0.1 \linewd 0.5

\move(0 -10)\lvec(0 10)\lvec(50 10)\lvec(50 -10)\lvec(0 -10)

\move(10 10)\lvec(10 -10)

\move(10 10)\lvec(10 -10)

\move(10 10)\lvec(20 20)\lvec(10 20)\lvec(10 10)

\htext(12 16){$\cdot$}

\htext(27 13){$\cdots$}

\move(40 10)\lvec(50 20)\lvec(40 20)\lvec(40 10)

\htext(42 16){$\cdot$}

\end{texdraw}}\hskip 10mm
$w_1\leftrightarrow$\raisebox{-0.4\height}{
\begin{texdraw}
\fontsize{6}{6}\drawdim em \setunitscale 0.1 \linewd 0.5

\move(0 -10)\lvec(0 10)\lvec(50 10)\lvec(50 -10)\lvec(0 -10)

\move(10 10)\lvec(10 -10)

\move(0 10)\lvec(10 20)\lvec(0 20)\lvec(0 10)

\htext(2 16){$*$}

\move(10 10)\lvec(20 20)\lvec(10 20)\lvec(10 10)

\htext(12 16){$\cdot$}

\htext(27 13){$\cdots$}

\move(40 10)\lvec(50 20)\lvec(40 20)\lvec(40 10)

\htext(42 16){$\cdot$}

\end{texdraw}}\vskip 5mm

$w_2\leftrightarrow$\raisebox{-0.4\height}{
\begin{texdraw}
\fontsize{6}{6}\drawdim em \setunitscale 0.1 \linewd 0.5

\move(0 -10)\lvec(0 10)\lvec(50 10)\lvec(50 -10)\lvec(0 -10)

\move(10 10)\lvec(10 -10)

\move(0 10)\lvec(10 20)\lvec(0 20)\lvec(0 10)

\htext(2 16){$*$}

\move(10 10)\lvec(20 20)\lvec(10 20)\lvec(10 10)

\htext(12 16){$\cdot$}

\move(40 10)\lvec(50 20)\lvec(50 10)\lvec(40 10)

\htext(46 12){$*$}

\htext(27 13){$\cdots$}

\move(40 10)\lvec(50 20)\lvec(40 20)\lvec(40 10)

\htext(42 16){$\cdot$}

\end{texdraw}}\hskip 10mm
$u\leftrightarrow$\raisebox{-0.4\height}{
\begin{texdraw}
\fontsize{6}{6}\drawdim em \setunitscale 0.1 \linewd 0.5

\move(0 -10)\lvec(0 10)\lvec(50 10)\lvec(50 -10)\lvec(0 -10)

\move(10 10)\lvec(10 -10)

\move(10 10)\lvec(20 20)\lvec(10 20)\lvec(10 10)

\htext(12 16){$\cdot$}

\move(40 10)\lvec(50 20)\lvec(50 10)\lvec(40 10)

\htext(46 12){$*$}

\htext(27 13){$\cdots$}

\move(40 10)\lvec(50 20)\lvec(40 20)\lvec(40 10)

\htext(42 16){$\cdot$}

\end{texdraw}}\vskip 5mm

$u'\leftrightarrow$\raisebox{-0.4\height}{
\begin{texdraw}
\fontsize{6}{6}\drawdim em \setunitscale 0.1 \linewd 0.5

\move(0 -10)\lvec(0 10)\lvec(50 10)\lvec(50 -10)\lvec(0 -10)

\move(10 10)\lvec(10 -10)

\move(10 10)\lvec(20 20)\lvec(20 10)\lvec(10 10)

\htext(16 12){$*$}

\move(40 10)\lvec(50 20)\lvec(50 10)\lvec(40 10)

\htext(46 12){$*$}

\htext(27 13){$\cdots$}

\move(0 10)\lvec(10 20)\lvec(10 10)\lvec(0 10)

\htext(6 12){$\cdot$}

\end{texdraw}}\hskip 10mm
$u''\leftrightarrow$\raisebox{-0.4\height}{
\begin{texdraw}
\fontsize{6}{6}\drawdim em \setunitscale 0.1 \linewd 0.5

\move(0 -10)\lvec(0 10)\lvec(50 10)\lvec(50 -10)\lvec(0 -10)

\move(10 10)\lvec(10 -10)

\move(10 10)\lvec(20 20)\lvec(20 10)\lvec(10 10)

\htext(16 12){$*$}

\move(40 10)\lvec(50 20)\lvec(50 10)\lvec(40 10)

\htext(46 12){$*$}

\htext(27 13){$\cdots$}

\move(40 10)\lvec(50 20)\lvec(40 20)\lvec(40 10)

\htext(42 16){$\cdot$}

\end{texdraw}}\quad .
\end{center}\vskip 5mm

The $U_{(i)}$-module $W_{2l+1}$ is decomposed as $W_{2l+1}\cong
V(2)\oplus V(0)^{\oplus 3}$, where
\begin{equation}
\begin{split}
V(2)&\cong \mathbb{Q}(q)w_0\oplus
\mathbb{Q}(q)(w_1+q_iu+q_i^{2l+1}u'')\oplus \mathbb{Q}(q)w_2,
\\
V(0)&\cong \mathbb{Q}(q)(u-q_iw_1)\cong
\mathbb{Q}(q)(u'+q_i^{2l+2}w_1)\cong \mathbb{Q}(q)u''.
\end{split}
\end{equation}

Hence, the crystal basis $(L,B)$ of $W_{2l+1}$ is given by
\begin{equation}
\begin{split}
L=& \mathbb{A}_0w_0\oplus
\mathbb{A}_0(w_1+q_iu+q_i^{2l+1}u'')\oplus \mathbb{A}_0w_2 \\
& \oplus \mathbb{A}_0(u-q_iw_1)\oplus
\mathbb{A}_0(u'+q_i^{2l+2}w_1)\oplus \mathbb{A}_0u'', \\
B=&\{\,\overline{w_0},\overline{w_1},\overline{w_2},
\overline{u},\overline{u'},\overline{u''}\,\}
\end{split}
\end{equation}
with the crystal graph
$$\overline{w_0}\stackrel{i}{\longrightarrow}\overline{w_1}
\stackrel{i}{\longrightarrow}\overline{w_2}\hskip 1cm
\overline{u}\hskip 1cm\overline{u'}\hskip
1cm\overline{u''}.$$\vskip 5mm

\begin{ex}{\rm
In Example \ref{comp3}, we have\vskip 5mm

\begin{center}
$Y=$\raisebox{-0.3\height}{
\begin{texdraw}
\fontsize{8}{8}\drawdim em \setunitscale 0.13 \linewd 0.5

\htext(-.2 3){$\cdots$} \move(10 0)\lvec(20 10)\lvec(20 0)\lvec(10
0)\lfill f:0.8 \htext(16 2){$0$}

\move(20 0)\lvec(30 10)\lvec(30 0)\lvec(20 0)\lfill f:0.8
\htext(26 2){$1$}

\end{texdraw}}\ $\otimes$
\raisebox{-0.4\height}{\begin{texdraw} \fontsize{8}{8}\drawdim em
\setunitscale 0.13 \linewd 0.5

\move(10 0)\lvec(20 0)\lvec(20 10)\lvec(10 10)\lvec(10 0)\htext(12
6){\tiny $1$}

\move(20 0)\lvec(30 0)\lvec(30 10)\lvec(20 10)\lvec(20 0)\htext(22
6){\tiny $0$}

\move(30 0)\lvec(40 0)\lvec(40 10)\lvec(30 10)\lvec(30 0)\htext(32
6){\tiny $1$}

\move(40 0)\lvec(50 0)\lvec(50 10)\lvec(40 10)\lvec(40 0)\htext(42
6){\tiny $0$}

\move(10 0)\lvec(20 10)\lvec(20 0)\lvec(10 0)\lfill f:0.8
\htext(16 2){\tiny $0$}

\move(20 0)\lvec(30 10)\lvec(30 0)\lvec(20 0)\lfill f:0.8
\htext(26 2){\tiny $1$}

\move(30 0)\lvec(40 10)\lvec(40 0)\lvec(30 0)\lfill f:0.8
\htext(36 2){\tiny $0$}

\move(40 0)\lvec(50 10)\lvec(50 0)\lvec(40 0)\lfill f:0.8
\htext(46 2){\tiny $1$}

\move(10 10)\lvec(20 10)\lvec(20 20)\lvec(10 20)\lvec(10
10)\htext(13 13){$2$}

\move(20 10)\lvec(30 10)\lvec(30 20)\lvec(20 20)\lvec(20
10)\htext(23 13){$2$}

\move(30 10)\lvec(40 10)\lvec(40 20)\lvec(30 20)\lvec(30
10)\htext(33 13){$2$}

\move(40 10)\lvec(50 10)\lvec(50 20)\lvec(40 20)\lvec(40
10)\htext(43 13){$2$}

\move(10 20)\lvec(20 20)\lvec(20 30)\lvec(10 30)\lvec(10
20)\htext(13 26){\tiny $3$}

\move(20 20)\lvec(30 20)\lvec(30 30)\lvec(20 30)\lvec(20
20)\htext(23 26){\tiny $3$}

\move(30 20)\lvec(40 20)\lvec(40 30)\lvec(30 30)\lvec(30
20)\htext(33 26){\tiny $3$}

\move(40 20)\lvec(50 20)\lvec(50 30)\lvec(40 30)\lvec(40
20)\htext(43 26){\tiny $3$}

\move(10 20)\lvec(20 20)\lvec(20 25)\lvec(10 25)\lvec(10 20)
\htext(13 21){\tiny $3$}

\move(20 20)\lvec(30 20)\lvec(30 25)\lvec(20 25)\lvec(20 20)
\htext(23 21){\tiny $3$}

\move(30 20)\lvec(40 20)\lvec(40 25)\lvec(30 25)\lvec(30 20)
\htext(33 21){\tiny $3$}

\move(40 20)\lvec(50 20)\lvec(50 25)\lvec(40 25)\lvec(40 20)
\htext(43 21){\tiny $3$}

\move(10 30)\lvec(20 30)\lvec(20 40)\lvec(10 40)\lvec(10
30)\htext(13 33){$2$}

\move(20 30)\lvec(30 30)\lvec(30 40)\lvec(20 40)\lvec(20
30)\htext(23 33){$2$}

\move(30 30)\lvec(40 30)\lvec(40 40)\lvec(30 40)\lvec(30
30)\htext(33 33){$2$}

\move(40 30)\lvec(50 30)\lvec(50 40)\lvec(40 40)\lvec(40
30)\htext(43 33){$2$}

\move(20 40)\lvec(30 50)\lvec(20 50)\lvec(20 40) \htext(22
46){\tiny $0$}

\move(30 40)\lvec(40 50)\lvec(30 50)\lvec(30 40) \htext(32
46){\tiny $1$}

\move(40 40)\lvec(50 50)\lvec(40 50)\lvec(40 40) \htext(42
46){\tiny $0$}

\move(40 40)\lvec(50 50)\lvec(50 40)\lvec(40 40) \htext(46
42){\tiny $1$}
\end{texdraw}}\ $\otimes$
\raisebox{-0.4\height}{\begin{texdraw} \fontsize{8}{8}\drawdim em
\setunitscale 0.13 \linewd 0.5

\move(50 0)\lvec(60 0)\lvec(60 10)\lvec(50 10)\lvec(50 0)\htext(52
6){\tiny $1$}

\move(50 0)\lvec(60 10)\lvec(60 0)\lvec(50 0)\lfill f:0.8
\htext(56 2){\tiny $0$}

\move(50 10)\lvec(60 10)\lvec(60 20)\lvec(50 20)\lvec(50
10)\htext(53 13){$2$}

\move(50 20)\lvec(60 20)\lvec(60 30)\lvec(50 30)\lvec(50
20)\htext(53 26){\tiny $3$}

\move(50 20)\lvec(60 20)\lvec(60 25)\lvec(50 25)\lvec(50
20)\htext(53 21){\tiny $3$}

\move(50 30)\lvec(60 30)\lvec(60 40)\lvec(50 40)\lvec(50
30)\htext(53 33){$2$}

\move(50 40)\lvec(60 40)\lvec(60 50)\lvec(50 50) \lvec(50
40)\htext(52 46){\tiny $1$}

\move(50 40)\lvec(60 50)\lvec(60 40)\lvec(50 40) \htext(56
42){\tiny $0$}

\move(50 50)\lvec(60 50)\lvec(60 60)\lvec(50 60)\lvec(50
50)\htext(53 53){$2$}

\end{texdraw}} \ $\otimes$
\raisebox{-0.4\height}{\begin{texdraw} \fontsize{8}{8}\drawdim em
\setunitscale 0.13 \linewd 0.5

\move(60 0)\lvec(70 0)\lvec(70 10)\lvec(60 10)\lvec(60 0)\htext(62
6){\tiny $0$}

\move(70 0)\lvec(80 0)\lvec(80 10)\lvec(70 10)\lvec(70 0)\htext(72
6){\tiny $1$}

\move(80 0)\lvec(90 0)\lvec(90 10)\lvec(80 10)\lvec(80 0)\htext(82
6){\tiny $0$}

\move(60 0)\lvec(70 10)\lvec(70 0)\lvec(60 0)\lfill f:0.8
\htext(66 2){\tiny $1$}

\move(70 0)\lvec(80 10)\lvec(80 0)\lvec(70 0)\lfill f:0.8
\htext(76 2){\tiny $0$}

\move(80 0)\lvec(90 10)\lvec(90 0)\lvec(80 0)\lfill f:0.8
\htext(86 2){\tiny $1$}

\move(60 10)\lvec(70 10)\lvec(70 20)\lvec(60 20)\lvec(60
10)\htext(63 13){$2$}

\move(70 10)\lvec(80 10)\lvec(80 20)\lvec(70 20)\lvec(70
10)\htext(73 13){$2$}

\move(80 10)\lvec(90 10)\lvec(90 20)\lvec(80 20)\lvec(80
10)\htext(83 13){$2$}

\move(60 20)\lvec(70 20)\lvec(70 25)\lvec(60 25)\lvec(60
20)\htext(63 26){\tiny $3$}

\move(70 20)\lvec(80 20)\lvec(80 25)\lvec(70 25)\lvec(70
20)\htext(73 26){\tiny $3$}

\move(80 20)\lvec(90 20)\lvec(90 25)\lvec(80 25)\lvec(80
20)\htext(83 26){\tiny $3$}

\move(60 25)\lvec(70 25)\lvec(70 30)\lvec(60 30)\lvec(60
25)\htext(63 21){\tiny $3$}

\move(70 25)\lvec(80 25)\lvec(80 30)\lvec(70 30)\lvec(70
25)\htext(73 21){\tiny $3$}

\move(80 25)\lvec(90 25)\lvec(90 30)\lvec(80 30)\lvec(80
25)\htext(83 21){\tiny $3$}

\move(60 30)\lvec(70 30)\lvec(70 40)\lvec(60 40)\lvec(60
30)\htext(63 33){$2$}

\move(70 30)\lvec(80 30)\lvec(80 40)\lvec(70 40)\lvec(70
30)\htext(73 33){$2$}

\move(80 30)\lvec(90 30)\lvec(90 40)\lvec(80 40)\lvec(80
30)\htext(83 33){$2$}

\move(60 40)\lvec(70 40)\lvec(70 50)\lvec(60 50)\lvec(60
40)\htext(62 46){\tiny $0$}

\move(70 40)\lvec(80 40)\lvec(80 50)\lvec(70 50)\lvec(70
40)\htext(72 46){\tiny $1$}

\move(80 40)\lvec(90 40)\lvec(90 50)\lvec(80 50)\lvec(80
40)\htext(82 46){\tiny $0$}

\move(60 40)\lvec(70 50)\lvec(70 40)\lvec(60 40) \htext(66
42){\tiny $1$}

\move(70 40)\lvec(80 50)\lvec(80 40)\lvec(70 40) \htext(76
42){\tiny $0$}

\move(80 40)\lvec(90 50)\lvec(90 40)\lvec(80 40) \htext(86
42){\tiny $1$}

\move(60 50)\lvec(70 50)\lvec(70 60)\lvec(60 60)\lvec(60
50)\htext(63 53){$2$}

\move(70 50)\lvec(80 50)\lvec(80 60)\lvec(70 60)\lvec(70
50)\htext(73 53){$2$}

\move(80 50)\lvec(90 50)\lvec(90 60)\lvec(80 60)\lvec(80
50)\htext(83 53){$2$}
\move(60 60)\lvec(70 60)\lvec(70 65)\lvec(60 65)\lvec(60
60)\htext(63 61){\tiny $3$}

\move(70 60)\lvec(80 60)\lvec(80 65)\lvec(70 65)\lvec(70
60)\htext(73 61){\tiny $3$}

\move(80 60)\lvec(90 60)\lvec(90 65)\lvec(80 65)\lvec(80
60)\htext(83 61){\tiny $3$}
\move(60 60)\lvec(70 60)\lvec(70 70)\lvec(60 70)\lvec(60
60)\htext(63 66){\tiny $3$}

\move(70 60)\lvec(80 60)\lvec(80 70)\lvec(70 70)\lvec(70
60)\htext(73 66){\tiny $3$}

\move(80 60)\lvec(90 60)\lvec(90 70)\lvec(80 70)\lvec(80
60)\htext(83 66){\tiny $3$}
\move(60 70)\lvec(70 70)\lvec(70 80)\lvec(60 80)\lvec(60
70)\htext(63 73){$2$}

\move(70 70)\lvec(80 70)\lvec(80 80)\lvec(70 80)\lvec(70
70)\htext(73 73){$2$}

\move(80 70)\lvec(90 70)\lvec(90 80)\lvec(80 80)\lvec(80
70)\htext(83 73){$2$}
\move(60 80)\lvec(70 90)\lvec(70 80)\lvec(60 80)\htext(66
82){\tiny $1$}

\move(70 80)\lvec(80 90)\lvec(80 80)\lvec(70 80)\htext(76
82){\tiny $0$}

\move(80 80)\lvec(90 90)\lvec(90 80)\lvec(80 80)\htext(86
82){\tiny $1$}

\end{texdraw}}\ $\in V\otimes W_3\otimes U\otimes W_2^-$.
\end{center}\vskip 5mm

}
\end{ex}

\begin{rem}\label{signature}{\rm
Let $Y$ be a proper Young wall in ${\mathcal Z}(\Lambda)$ with the
decomposition $Y=(Y_0,Y_1,\cdots,Y_r)$ into $i$-components. For
each $0\leq k\leq r$, let $V_k$ be the $U_{(i)}$-module associated
with $Y_k$, whose crystal basis is $(L_k,B_k)$. We identify $Y_k$
with a basis element of $V_k$. Then $Y_k$ can also be viewed as a
crystal element in $B_k$ and hence, $\varphi_i(Y_k)$ and
$\varepsilon_i(Y_k)$ are well-defined. On the other hand, as a
part of $Y$, we defined $\varphi_i(Y_k)$ (resp.
$\varepsilon_i(Y_k)$) to be the number of $+$'s (resp. $-$'s) in
the $i$-signature of $Y_k$. It is easy to verify that these two
definitions give the same values for the $i$-component $Y_k$.}
\end{rem}

Let $Y$ be a proper Young wall in ${\mathcal Z}(\Lambda)$ with the
decomposition $Y=(Y_0,Y_1,\cdots,Y_r)$ into $i$-components, and
let $V_k$ be the $U_{(i)}$-module associated with $Y_k$ ($0\leq
k\leq r$). Recall that $Y_k$ is identified with a basis element of
$V_k$. Set
\begin{equation}
V_Y=V_0\otimes V_1\otimes\cdots\otimes V_r.
\end{equation}
We define a $\mathbb{Q}(q)$-linear map $\theta_Y : V_Y \rightarrow
{\mathcal F}(\Lambda)$ by
\begin{equation}
\theta_Y(Y_0'\otimes Y_1'\otimes\cdots\otimes
Y_r')=Y'=(Y_0',Y_1'\cdots,Y_r'),
\end{equation}
where $Y_0'\otimes Y_1'\otimes\cdots\otimes Y_r'$ runs over the
basis element of $V_Y$. Then it is easy to see that $\theta_Y$ is
injective and $Y$ is contained in ${\rm Im}\ \theta_Y$.

We are now ready to prove the relation \eqref{rel2}.

\begin{lem}\label {Vy1}
The linear map $\theta_Y$ is a $U_{(i)}$-module homomorphism. In
particular, we have
$$(e_if_i-f_ie_i)Y=\dfrac{K_i-K_i^{-1}}{q_i-q_i^{-1}}Y.$$
\end{lem}
\pf By definition of $U_{(i)}$-action on $V_Y$ and
$U_q(\frak{g})$-action on ${\mathcal F}(\Lambda)$, it is rather
tedious but straightforward to verify that
\begin{equation}\label {Weylrel1}
\theta_{Y}(x\cdot v)=x\cdot\theta_{Y}(v) \quad\text{for all $x\in
U_{(i)}$ and $v\in V_Y$}.
\end{equation}
We will prove only when $x=K_i$. Let $v=Y_0'\otimes
Y_1'\otimes\cdots\otimes Y_r'\in V_Y$ and put
$Y'=\theta_Y(v)=(Y_0',Y_1'\cdots,Y_r')$. Then we see from the
action of $K_i$ on $Y_k'\in V_k$ that $K_i v=q_i^lv$ where
$l=\sum_{k=0}^r(\varphi_i(Y_k')-\varepsilon_i(Y_k'))$. By Remark
\ref{signature}, we have $l=\varphi_i(Y')-\varepsilon_i(Y')= {\rm
wt}(Y')(h_i)$, which implies $\theta_{Y}(K_iv)=K_i\theta_{Y}(v)$.

Hence
\begin{equation}
\begin{split}
&(e_if_i-f_ie_i)Y=\theta_{Y}((e_if_i-f_ie_i)(Y_0\otimes
Y_1\otimes\cdots\otimes Y_r)) \\
&=\theta_{Y}(\frac{K_i-K_i^{-1}}{q_i-q_i^{-1}}(Y_0\otimes
Y_1\otimes\cdots\otimes Y_r))=\frac{K_i-K_i^{-1}}{q_i-q_i^{-1}}Y.
\end{split}
\end{equation}
\qed \vskip 3mm

\begin{lem}\label{Vy2} Let $Y$ be a proper Young wall in ${\mathcal Z}(\Lambda)$. For $i,j\in I$ with $i\neq
j$, we have
\begin{equation}
(e_if_j-f_je_i)Y=0.
\end{equation}
\end{lem}
\pf The proof is quite lengthy and is based on the case-by-case
check. We give a sketch of the verification and leave the details
to the reader.

For $i\in I$, we define
\begin{align*}
&\text{${\rm adm}_i(Y)$ $=$ the set of all (virtually) admissible
$i$-slots in $Y$}, \\
&\text{${\rm rmv}_i(Y)$ $=$ the set of all (virtually) removable
$i$-blocks in $Y$}.
\end{align*}

Consider
\begin{equation}
\begin{split}
e_if_jY&=\sum_{\substack{b\in{\rm adm}_j(Y) \\ c\in{\rm
rmv}_i(Y\swarrow b)}}q_j^{L_j(b;Y)}q_i^{-R_i(c;Y\swarrow
b)}((Y\swarrow b)\nearrow c), \\
f_je_iY&=\sum_{\substack{c \in{\rm rmv}_i(Y)\\
b\in{\rm adm}_j(Y\nearrow c) }}q_j^{L_j(b;Y\nearrow c
)}q_i^{-R_i(c;Y)}((Y\nearrow c)\swarrow b).
\end{split}
\end{equation}

Let $\Box$ be the type of the $i$-blocks and let $\Box'$ be the
type of the $j$-blocks. Suppose $(\Box,\Box')\neq {\rm (III,III)}$
and observe that
\begin{multline}\label {adrm}
b\in{\rm adm}_j(Y),\, c\in{\rm rmv}_i(Y\swarrow b)  \\
\text{if and only if}\ \ c \in{\rm rmv}_i(Y),\, b\in{\rm
adm}_j(Y\nearrow c).
\end{multline}

In this case, we have $((Y\swarrow b)\nearrow c)=((Y\nearrow
c)\swarrow b)$. Furthermore, if $(\alpha_i|\alpha_j)=0$
(equivalently, $i$-blocks and $j$-blocks are not adjacent in each
column of $Y$) or $b$ is located to the left of $c$, then it is
clear that
\begin{equation}
L_j(b;Y)=L_j(b;Y\nearrow c),\quad R_i(c;Y\swarrow b)=R_i(c;Y),
\end{equation}
which implies that the corresponding two summands in
$(e_if_j-f_je_i)Y$ are cancelled out. Therefore, we have only to
consider the case when $(\alpha_i|\alpha_j)\neq 0$ and $b$ is
located to the right of $c$ and show that
\begin{equation}\label {LiRj}
q_j^{L_j(b;Y)}q_i^{-R_i(c;Y\swarrow b)}=q_j^{L_j(b;Y\nearrow c
)}q_i^{-R_i(c;Y)}.
\end{equation}

On the other hand, we can verify the following :

(i) if $(\alpha_i|\alpha_i)=(\alpha_j|\alpha_j)$, i.e. $s_i=s_j$,
then we have
\begin{equation}
\begin{split}
L_j(b;Y)&=L_j(b;Y\nearrow c)+1, \\
R_i(c;Y\swarrow b)&=R_i(c;Y) +1.
\end{split}
\end{equation}

(ii) if $(\alpha_i|\alpha_i)>(\alpha_j|\alpha_j)$, i.e.
$s_i=2s_j$, then we have
\begin{equation}
\begin{split}
L_j(b;Y)&=L_j(b;Y\nearrow c)+2, \\
R_i(c;Y\swarrow b)&=R_i(c;Y) +1.
\end{split}
\end{equation}

(iii) if $(\alpha_i|\alpha_i)<(\alpha_j|\alpha_j)$, i.e.
$2s_i=s_j$, then we have
\begin{equation}
\begin{split}
L_j(b;Y)&=L_j(b;Y\nearrow c)+1, \\
R_i(c;Y\swarrow b)&=R_i(c;Y) +2.
\end{split}
\end{equation}
It follows that
\begin{equation}
s_jL_j(b;Y)-s_iR_i(c;Y\swarrow b)=s_jL_j(b;Y\nearrow c
)-s_iR_i(c;Y),
\end{equation}
which proves \eqref{LiRj}.

Now suppose that $\Box=\Box'={\rm III}$. Note that
$(\alpha_i|\alpha_i)=(\alpha_j|\alpha_j)$ and
$(\alpha_i|\alpha_j)=0$. For $b\in{\rm adm}_j(Y)$ and $c\in{\rm
rmv}_i(Y\swarrow b)$, set
\begin{equation*}
W^+_{i,j}(b,c)=q_j^{L_j(b;Y)}q_i^{-R_i(c;Y\swarrow b)}((Y\swarrow
b)\nearrow c),
\end{equation*}
and for $c \in{\rm rmv}_i(Y)$ and $b\in{\rm adm}_j(Y\nearrow c)$,
set
\begin{equation*}
W^-_{i,j}(b,c)=q_j^{L_j(b;Y\nearrow c )}q_i^{-R_i(c;Y)}((Y\nearrow
c)\swarrow b).
\end{equation*}
First, consider the following cases.

{\bf Case 1.}\vskip 5mm

\hskip 2cm \raisebox{-0.4\height}
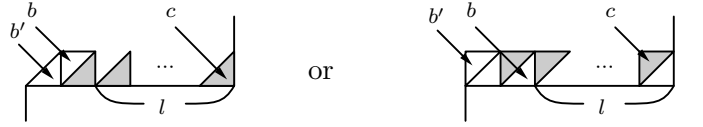
{
\begin{texdraw}
\drawdim em \setunitscale 0.12 \linewd 0.5

\move(0 0)\lvec(0 10)\lvec(60 10)\lvec(60 30)

\move(10 10)\lvec(10 20)\lvec(20 20)

\move(10 10)\lvec(20 20)\lvec(20 10)\lvec(10 10)\lfill f:0.8

\move(20 10)\lvec(30 20)\lvec(30 10)\lvec(20 10)\lfill f:0.8

\move(50 10)\lvec(60 20)\lvec(60 10)\lvec(50 10)\lfill f:0.8

\move(0 10)\lvec(10 20)

\htext(37 14){$_{\cdots}$}

\move(3 27) \arrowheadtype t:F \arrowheadsize l:4 w:2 \avec(13 17)

\htext(0 30){$_b$}

\move(-2 22) \arrowheadtype t:F \arrowheadsize l:4 w:2 \avec(8 12)

\htext(-5 24){$_{b'}$}

\move(43 27) \arrowheadtype t:F \arrowheadsize l:4 w:2 \avec(58
12)

\htext(40 30){$_c$}

\move(34 5)\clvec(26 5)(23 5)(20 10)

\move(46 5)\clvec(54 5)(57 5)(60 10)

\htext(38 2){$_{l}$}
\end{texdraw}}\hskip 1cm or \hskip 1cm
\raisebox{-0.4\height}{
\begin{texdraw}
\drawdim em \setunitscale 0.12 \linewd 0.5

\move(0 0)\lvec(0 10)\lvec(60 10)\lvec(60 30)


\move(10 10)\lvec(20 20)\lvec(10 20)\lvec(10 10)\lfill f:0.8

\move(20 10)\lvec(30 20)\lvec(20 20)\lvec(20 10)\lfill f:0.8

\move(50 10)\lvec(60 20)\lvec(50 20)\lvec(50 10)\lfill f:0.8

\move(0 10)\lvec(0 20)\lvec(10 20)\lvec(0 10)

\htext(37 14){$_{\cdots}$}

\move(3 27) \arrowheadtype t:F \arrowheadsize l:4 w:2 \avec(18 12)

\htext(0 30){$_b$}

\move(-8 26) \arrowheadtype t:F \arrowheadsize l:4 w:2 \avec(2 16)

\htext(-11 28){$_{b'}$}

\move(43 27) \arrowheadtype t:F \arrowheadsize l:4 w:2 \avec(53
17)

\htext(40 30){$_c$}

\move(34 5)\clvec(26 5)(23 5)(20 10)

\move(46 5)\clvec(54 5)(57 5)(60 10)

\htext(38 2){$_{l}$}
\end{texdraw}}\quad ,
\vskip 5mm

\noindent where $l=2m$ ($m\geq 0$), $b,b'\in{\rm adm}_j(Y)$ and
$c\in {\rm rmv}_i(Y\swarrow b)\cap{\rm rmv}_i(Y\swarrow b')$. In
this case, we have
\begin{equation*}
\begin{split}
&((Y\swarrow b')\nearrow c)=(-q_i)((Y\swarrow b)\nearrow c), \\
&L_j(b';Y)=L_j(b;Y)-1, \quad R_i(c;Y\swarrow b')=R_i(c;Y\swarrow
b),
\end{split}
\end{equation*}
which yields
\begin{equation}
W^+_{i,j}(b,c)=-W^+_{i,j}(b',c).
\end{equation}
\vskip 2mm

{\bf Case 2.} \vskip 5mm

\hskip 1.5cm\raisebox{-0.4\height}
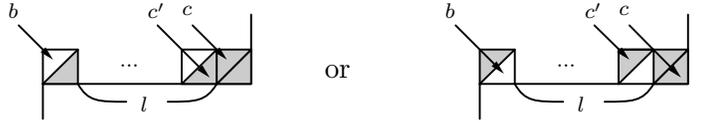
{
\begin{texdraw}
\drawdim em \setunitscale 0.12 \linewd 0.5

\move(0 0)\lvec(0 10)\lvec(60 10)\lvec(60 30)

\move(0 10)\lvec(0 20)\lvec(10 20)

\move(0 10)\lvec(10 20)\lvec(10 10)\lvec(0 10)\lfill f:0.8

\move(40 10)\lvec(50 20)\lvec(50 10)\lvec(40 10)\lfill f:0.8

\move(50 10)\lvec(60 20)\lvec(60 10)\lvec(50 10)\lfill f:0.8

\move(50 10)\lvec(60 20)\lvec(50 20)\lvec(50 10)\lfill f:0.8

\move(40 10)\lvec(40 20)\lvec(50 20)

\htext(22 14){$_{\cdots}$}

\move(-7 27) \arrowheadtype t:F \arrowheadsize l:4 w:2 \avec(3 17)

\htext(-10 29){$_b$}

\move(33 27) \arrowheadtype t:F \arrowheadsize l:4 w:2 \avec(48
12)

\htext(30 29){$_{c'}$}

\move(43 27) \arrowheadtype t:F \arrowheadsize l:4 w:2 \avec(53
17)

\htext(40 30){$_c$}

\move(24 5)\clvec(16 5)(13 5)(10 10)

\move(36 5)\clvec(44 5)(47 5)(50 10)

\htext(28 2){$_{l}$}
\end{texdraw}}\hskip 1cm or \hskip 1cm
\raisebox{-0.4\height}{
\begin{texdraw}
\drawdim em \setunitscale 0.12 \linewd 0.5

\move(0 0)\lvec(0 10)\lvec(60 10)\lvec(60 30)

\move(10 10)\lvec(10 20)

\move(0 10)\lvec(10 20)\lvec(0 20)\lvec(0 10)\lfill f:0.8

\move(40 10)\lvec(50 20)\lvec(40 20)\lvec(40 10)\lfill f:0.8

\move(50 10)\lvec(60 20)\lvec(60 10)\lvec(50 10)\lfill f:0.8

\move(50 10)\lvec(60 20)\lvec(50 20)\lvec(50 10)\lfill f:0.8

\move(40 10)\lvec(40 20)\lvec(50 20)

\htext(22 14){$_{\cdots}$}

\move(-7 27) \arrowheadtype t:F \arrowheadsize l:4 w:2 \avec(7 13)

\htext(-10 29){$_b$}

\move(33 27) \arrowheadtype t:F \arrowheadsize l:4 w:2 \avec(43
17)

\htext(30 29){$_{c'}$}

\move(43 27) \arrowheadtype t:F \arrowheadsize l:4 w:2 \avec(58
12)

\htext(40 30){$_c$}

\move(24 5)\clvec(16 5)(13 5)(10 10)

\move(36 5)\clvec(44 5)(47 5)(50 10)

\htext(28 2){$_{l}$}
\end{texdraw}}\quad ,
\vskip 5mm \noindent where $l=2m$ ($m\geq 0$), $b\in {\rm
adm}_j(Y\nearrow c)\cap{\rm adm}_j(Y\nearrow c')$ and $c,c'\in{\rm
rmv}_i(Y)$. In this case, we have
\begin{equation*}
\begin{split}
&((Y\nearrow c)\swarrow b)=(-q_i)((Y\nearrow c')\swarrow b), \\
&L_j(b;Y\nearrow c)=L_j(b;Y\nearrow c'), \quad
R_i(c;Y)=R_i(c';Y)+1,
\end{split}
\end{equation*}
which yields
\begin{equation}
W^-_{i,j}(b,c)=-W^-_{i,j}(b,c')
\end{equation}
\vskip 2mm

{\bf Case 3.} \vskip 5mm

\hskip 1.5cm\raisebox{-0.4\height}
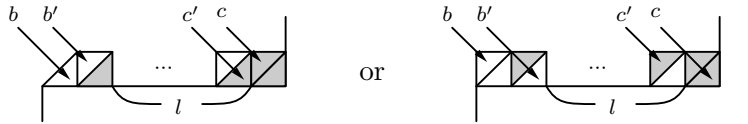
{
\begin{texdraw}
\drawdim em \setunitscale 0.12 \linewd 0.5

\move(-10 0)\lvec(-10 10)\lvec(60 10)\lvec(60 30)

\move(0 10)\lvec(0 20)\lvec(10 20)

\move(-10 10)\lvec(0 20)

\move(0 10)\lvec(10 20)\lvec(10 10)\lvec(0 10)\lfill f:0.8

\move(40 10)\lvec(50 20)\lvec(50 10)\lvec(40 10)\lfill f:0.8

\move(50 10)\lvec(60 20)\lvec(60 10)\lvec(50 10)\lfill f:0.8

\move(50 10)\lvec(60 20)\lvec(50 20)\lvec(50 10)\lfill f:0.8

\move(40 10)\lvec(40 20)\lvec(50 20)

\htext(22 14){$_{\cdots}$}

\move(-7 27) \arrowheadtype t:F \arrowheadsize l:4 w:2 \avec(3 17)

\htext(-10 29){$_{b'}$}

\move(-17 27) \arrowheadtype t:F \arrowheadsize l:4 w:2 \avec(-2
12)

\htext(-20 29){$_{b}$}

\move(33 27) \arrowheadtype t:F \arrowheadsize l:4 w:2 \avec(48
12)

\htext(30 29){$_{c'}$}

\move(43 27) \arrowheadtype t:F \arrowheadsize l:4 w:2 \avec(53
17)

\htext(40 30){$_c$}

\move(24 5)\clvec(16 5)(13 5)(10 10)

\move(36 5)\clvec(44 5)(47 5)(50 10)

\htext(28 2){$_{l}$}
\end{texdraw}}\hskip 1cm or \hskip .5cm
\raisebox{-0.4\height}{
\begin{texdraw}
\drawdim em \setunitscale 0.12 \linewd 0.5

\move(-10 0)\lvec(-10 10)\lvec(60 10)\lvec(60 30)

\move(-10 10)\lvec(-10 20)\lvec(0 20)

\move(-10 10)\lvec(0 20)

\move(0 10)\lvec(10 20)\lvec(0 20)\lvec(0 10)\lfill f:0.8

\move(40 10)\lvec(50 20)\lvec(40 20)\lvec(40 10)\lfill f:0.8

\move(50 10)\lvec(60 20)\lvec(60 10)\lvec(50 10)\lfill f:0.8

\move(50 10)\lvec(60 20)\lvec(50 20)\lvec(50 10)\lfill f:0.8

\move(10 10)\lvec(10 20)

\htext(22 14){$_{\cdots}$}

\move(-7 27) \arrowheadtype t:F \arrowheadsize l:4 w:2 \avec(8 12)

\htext(-10 29){$_{b'}$}

\move(-17 27) \arrowheadtype t:F \arrowheadsize l:4 w:2 \avec(-7
17)

\htext(-20 29){$_{b}$}

\move(33 27) \arrowheadtype t:F \arrowheadsize l:4 w:2 \avec(43
17)

\htext(30 29){$_{c'}$}

\move(43 27) \arrowheadtype t:F \arrowheadsize l:4 w:2 \avec(58
12)

\htext(40 30){$_c$}

\move(24 5)\clvec(16 5)(13 5)(10 10)

\move(36 5)\clvec(44 5)(47 5)(50 10)

\htext(28 2){$_{l}$}
\end{texdraw}}\quad ,
\vskip 5mm

\noindent where $l=2m$ ($m\geq 0$),
\begin{equation*}
\begin{split}
& b\in{\rm adm}_j(Y)\cap {\rm adm}_j(Y\nearrow c),\quad b'
\in{\rm adm}_j(Y\nearrow c), \\
\text{and}\ \ & c\in {\rm rmv}_i(Y)\cap {\rm rmv}_i(Y\swarrow
b),\quad c'\in {\rm rmv}_i(Y\swarrow b).
\end{split}
\end{equation*}
In this case, we have
\begin{align*}
W^+_{i,j}(b,c)&=W^-_{i,j}(b,c), \\
W^+_{i,j}(b,c')&=W^-_{i,j}(b',c).
\end{align*}
\vskip 2mm

{\bf Case 4.} \vskip 5mm

\hskip 1.5cm \raisebox{-0.4\height}
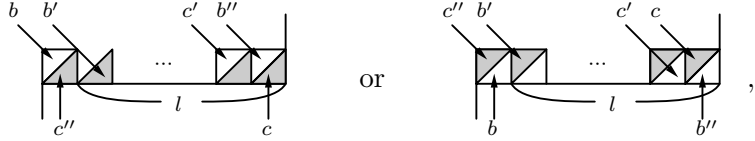
{
\begin{texdraw}
\drawdim em \setunitscale 0.12 \linewd 0.5

\move(-10 0)\lvec(-10 10)\lvec(60 10)\lvec(60 30)

\move(-10 10)\lvec(-10 20)\lvec(0 20)

\move(-10 10)\lvec(0 20)

\move(-10 10)\lvec(0 20)\lvec(0 10)\lvec(-10 10)\lfill f:0.8

\move(0 10)\lvec(10 20)\lvec(10 10)\lvec(0 10)\lfill f:0.8

\move(40 10)\lvec(50 20)\lvec(50 10)\lvec(40 10)\lfill f:0.8

\move(50 10)\lvec(60 20)\lvec(60 10)\lvec(50 10)\lfill f:0.8

\move(50 10)\lvec(60 20)\lvec(50 20)

\move(40 10)\lvec(40 20)\lvec(50 20)

\htext(22 14){$_{\cdots}$}

\move(-7 27) \arrowheadtype t:F \arrowheadsize l:4 w:2 \avec(8 12)

\htext(-10 29){$_{b'}$}

\move(-17 27) \arrowheadtype t:F \arrowheadsize l:4 w:2 \avec(-7
17)

\htext(-20 29){$_{b}$}

\move(33 27) \arrowheadtype t:F \arrowheadsize l:4 w:2 \avec(43
17)

\htext(30 29){$_{c'}$}

\move(43 27) \arrowheadtype t:F \arrowheadsize l:4 w:2 \avec(53
17)

\htext(40 29){$_{b''}$}

\move(-5 0) \arrowheadtype t:F \arrowheadsize l:4 w:2 \avec(-5 13)

\htext(-7 -5){$_{c''}$}

\move(55 0) \arrowheadtype t:F \arrowheadsize l:4 w:2 \avec(55 13)

\htext(53 -5){$_{c}$}

\move(24 5)\clvec(16 5)(3 5)(0 10)

\move(36 5)\clvec(44 5)(57 5)(60 10)

\htext(28 2){$_{l}$}
\end{texdraw}}\hskip 1cm or \hskip .5cm\raisebox{-0.4\height}{
\begin{texdraw}
\drawdim em \setunitscale 0.12 \linewd 0.5

\move(-10 0)\lvec(-10 10)\lvec(60 10)\lvec(60 30)

\move(-10 10)\lvec(-10 20)\lvec(0 20)

\move(10 10)\lvec(10 20)

\move(-10 10)\lvec(0 20)\lvec(-10 20)\lvec(-10 10)\lfill f:0.8

\move(0 10)\lvec(10 20)\lvec(0 20)\lvec(0 10)\lfill f:0.8

\move(40 10)\lvec(50 20)\lvec(40 20)\lvec(40 10)\lfill f:0.8

\move(50 10)\lvec(60 20)\lvec(50 20)\lvec(50 10)\lfill f:0.8

\move(50 10)\lvec(60 20)\lvec(50 20)

\move(40 10)\lvec(40 20)\lvec(50 20)

\htext(22 14){$_{\cdots}$}

\move(-7 27) \arrowheadtype t:F \arrowheadsize l:4 w:2 \avec(3 17)

\htext(-10 29){$_{b'}$}

\move(-17 27) \arrowheadtype t:F \arrowheadsize l:4 w:2 \avec(-7
17)

\htext(-20 29){$_{c''}$}

\move(33 27) \arrowheadtype t:F \arrowheadsize l:4 w:2 \avec(48
12)

\htext(30 29){$_{c'}$}

\move(43 27) \arrowheadtype t:F \arrowheadsize l:4 w:2 \avec(53
17)

\htext(40 29){$_{c}$}

\move(-5 0) \arrowheadtype t:F \arrowheadsize l:4 w:2 \avec(-5 13)

\htext(-7 -5){$_{b}$}

\move(55 0) \arrowheadtype t:F \arrowheadsize l:4 w:2 \avec(55 13)

\htext(53 -5){$_{b''}$}

\move(24 5)\clvec(16 5)(3 5)(0 10)

\move(36 5)\clvec(44 5)(57 5)(60 10)

\htext(28 2){$_{l}$}
\end{texdraw}}\quad ,\vskip 5mm

\noindent where $l=2m$ ($m\geq 0$),
\begin{equation*}
\begin{split}
& b\in{\rm adm}_j(Y)\cap {\rm adm}_j(Y\nearrow c),\,\, b'
\in{\rm adm}_j(Y\nearrow c), \\ & b''\in{\rm adm}_j(Y)\cap {\rm adm}_j(Y\nearrow c''), \\
 & c\in {\rm rmv}_i(Y)\cap {\rm rmv}_i(Y\swarrow
b),\,\, c'\in {\rm rmv}_i(Y\swarrow b),\\\text{and} \ \ & c''\in
{\rm rmv}_i(Y)\cap {\rm rmv}_i(Y\swarrow b'').
\end{split}
\end{equation*}
Similarly, we have
\begin{align*}
W^+_{i,j}(b,c)&=W^-_{i,j}(b,c), \\
W^+_{i,j}(b,c')&=W^-_{i,j}(b',c), \\
W^+_{i,j}(b'',c'')&=W^-_{i,j}(b'',c'').
\end{align*}
\vskip 2mm

For the other cases, it is easy to verify that \eqref{adrm} holds
and $W^+_{i,j}(b,c)=W^-_{i,j}(b,c)$. Therefore,
$(e_if_j-f_je_i)Y=0$, which completes the proof of the lemma. \qed
\vskip 3mm

{\bf Proof of Theorem \ref{Fock space}}. By Lemma \ref{Vy1}, Lemma
\ref{Vy2} and Proposition B.1 in \cite{KMPY}, the
$U_q(\frak{g})$-action on ${\mathcal F}(\Lambda)$ satisfies all
the relations in \eqref{defining rels}. Therefore, ${\mathcal
F}(\Lambda)$ becomes a $U_q(\frak{g})$-module in the category
${\mathcal O}_{int}$. \qed

\begin{rem}{\rm
If $\frak{g}=A_n^{(1)}$, the Fock space $\mathcal{F}(\Lambda)$ is
equal to the Fock space constructed by Misra and Miwa \cite{MM}
where $\mathcal{Z}(\Lambda)$ is the set of {\it Young diagrams}
and $\mathcal{Y}(\Lambda)$ is the set of {\it $n$-reduced Young
diagrams}. In \cite{KMPY}, Kashiwara, Miwa, Petersen and Yung gave
a more abstract construction of the Fock space representations of
quantum affine algebras. More precisely, for a level $l$ perfect
representation $V$ of $U'_q(\frak{g})$, they first defined the
{\it $q$-deformed wedge space $\bigwedge^r (V^{\rm aff})$} where
$V^{\rm aff}$ denotes the affinization of $V$. Then they defined
the Fock space to be the inductive limit of $q$-deformed wedge
spaces. The set of {\it normally ordered wedges} (defined by the
energy function) form a $\mathbb{Q}(q)$-basis of the Fock space.
For the level 1 case, the Fock space representation constructed in
\cite{KMPY} is isomorphic to $U_q(\frak{g})$-module ${\mathcal
F}(\Lambda)$ constructed in this paper. Moreover, there is a
bijection between the set of normally ordered wedges and the set
of proper Young walls ${\mathcal Z}(\Lambda)$ except for the case
$\frak{g}=D_{n+1}^{(2)}$. We expect that one can also construct
the higher level Fock space representations of quantum affine
algebras using combinatorics of Young walls.}
\end{rem}
\vskip 1cm

\section{Crystal basis of ${\mathcal F}(\Lambda)$.}

Let ${\mathcal L}(\Lambda)=\bigoplus_{Y\in{\mathcal
Z}(\Lambda)}\mathbb{A}_0 Y$. We will show that $({\mathcal
L}(\Lambda),{\mathcal Z}(\Lambda))$ is a crystal basis of
${\mathcal F}(\Lambda)$. In particular, the crystal of ${\mathcal
F}(\Lambda)$ is isomorphic to the $U_q(\frak{g})$-crystal
${\mathcal Z}(\Lambda)$ defined by the abstract Kashiwara
operators $\tilde{E}_i$ and $\tilde{F}_i$ ($i\in I$).

Observe that the pair $({\mathcal L}(\Lambda),{\mathcal
Z}(\Lambda))$ satisfies the first four conditions in Definition
\ref{defi:crystal basis}. For the rest three conditions, the main
step is to show that the Kashiwara operators $\tilde{e}_i$ and
$\tilde{f}_i$ ($i\in I$) on ${\mathcal Z}(\Lambda)$ induced by the
$U_q(\frak{g})$-module action on ${\mathcal F}(\Lambda)$ coincide
with the abstract Kashiwara operators $\tilde{E}_i$ and
$\tilde{F}_i$ ($i\in I$) on ${\mathcal Z}(\Lambda)$ defined in
Section 4. The proof of this step relies on the crystal basis
theory for $U_q(\frak{sl}_2)$-modules and the tensor product rule.

\begin{thm}\label {crystal basis for Fock space}
The pair $({\mathcal L}(\Lambda),{\mathcal Z}(\Lambda))$ is a
crystal basis of the Fock space representation ${\mathcal
F}(\Lambda)$. Moreover, the crystal of ${\mathcal F}(\Lambda)$ is
isomorphic to the $U_q(\frak{g})$-crystal ${\mathcal Z}(\Lambda)$
given in Section 4.
\end{thm}
\pf  We will show that the pair $({\mathcal L}(\Lambda),{\mathcal
Z}(\Lambda))$ satisfies the conditions (v), (vi) and (vii) in
Definition \ref{defi:crystal basis}.

Fix $i\in I$. Let $Y$ be a proper Young wall in ${\mathcal
Z}(\Lambda)$ with the $i$-component decomposition
$(Y_0,Y_1,\cdots,Y_r)$, and let $V_{Y}=V_0\otimes\cdots\otimes
V_{r}$ be the $U_{(i)}$-module constructed in Section 5. By Lemma
\ref{Vy1}, $\theta_{Y} : V_{Y}\rightarrow {\mathcal F}(\Lambda)$
is an injective $U_{(i)}$-module homomorphism. Since
$(L_{k},B_{k})$ is a crystal basis of $V_{k}$ ($0\leq k\leq r$),
$V_{Y}$ has a crystal basis $(L_{Y},B_{Y})$ given by
\begin{equation}
L_{Y}=L_{0}\otimes\cdots\otimes L_{r}, \ \ \
B_{Y}=B_{0}\otimes\cdots\otimes B_{r}.
\end{equation}
Furthermore, $\theta_{Y}$ satisfies
\begin{equation}\label {theta}
\theta_{Y}(L_{Y})\subset {\mathcal L}(\Lambda),\quad
\overline{\theta_{Y}}(B_{Y})\subset {\mathcal Z}(\Lambda),
\end{equation}
where $\overline{\theta_{Y}} : L_{Y}/qL_{Y}\longrightarrow
{\mathcal L}(\Lambda)/q{\mathcal L}(\Lambda)$ is the injective
$\mathbb{Q}$-linear map induced from $\theta_{Y}$.

Let $\tilde{e}_i$ and $\tilde{f}_i$  be the Kashiwara operators
induced from the $U_{(i)}$-module structure on $V_{Y}$ and
${\mathcal F}(\Lambda)$.  Since $\tilde{e}_i$ and $\tilde{f}_i$
commute with $\theta_{Y}$ and $Y$ is contained in
$\theta_{Y}(L_{Y})$, we have $\tilde{e}_iY, \tilde{f}_iY \in
{\mathcal L}(\Lambda)$. Hence, the condition (v) is satisfied.

Since $\tilde{e}_i$ and $\tilde{f}_i$ commute with
$\overline{\theta_{Y}}$, we also have $\tilde{e}_iY, \tilde{f}_iY
\in {\mathcal Z}(\Lambda)\cup \{0\} \mod{q{\mathcal L}(\Lambda)}$.
Furthermore, if $\tilde{e}_iY\neq 0$ (resp. $\tilde{f}_iY\neq
0$)$\mod{q{\mathcal L}(\Lambda)}$, then
$\tilde{e}_i\tilde{f}_iY=Y$ (resp.
$\tilde{f}_i\tilde{e}_iY=Y$)$\mod{q{\mathcal L}(\Lambda)}$, which
implies that the condition (vi) and (vii) are satisfied. Hence,
$({\mathcal L}(\Lambda),{\mathcal Z}(\Lambda))$ is a crystal basis
of ${\mathcal F}(\Lambda)$.

It remains to show that
\begin{equation}\label{eEfF}
\tilde{e}_i Y=\tilde{E}_i Y, \ \ \tilde{f}_i Y=\tilde{F}_i Y \ \
\ \mod{q{\mathcal L}(\Lambda)}.
\end{equation}
It follows from Remark \ref{signature} that the $i$-signature of
$Y_0\otimes\cdots\otimes Y_r$ (see Section 2) is equal to the
$i$-signature of $Y$ (see Section 4). Then by the tensor product
rule and the definitions of $\tilde{E}_i$ and $\tilde{F}_i$, we
have
\begin{equation}
\begin{split}
\overline{\theta_{Y}}(\tilde{e}_i(Y_0\otimes\cdots\otimes
Y_r))&=\tilde{E}_iY, \\
\overline{\theta_{Y}}(\tilde{f}_i(Y_0\otimes\cdots\otimes
Y_r))&=\tilde{F}_iY.
\end{split}
\end{equation}
Therefore, we obtain \eqref{eEfF} and conclude that the crystal of
${\mathcal F}(\Lambda)$ is isomorphic to the
$U_q(\frak{g})$-crystal ${\mathcal Z}(\Lambda)$. \qed \vskip 3mm

Using Theorem \ref{crystal basis for Fock space}, one can
decompose the Fock space ${\mathcal F}(\Lambda)$ into a direct sum
of irreducible highest weight modules over $U_q(\frak{g})$ by
locating the maximal vectors in the crystal ${\mathcal
Z}(\Lambda)$.
\begin{cor}\label {decomposition}
\begin{equation}
{\mathcal F}(\Lambda)=
\begin{cases}
\bigoplus_{m\geq 0}V(\Lambda-m\delta)^{\oplus p(m)}
& \text{if $\frak{g}\neq D^{(2)}_{n+1}$,} \\
\bigoplus_{m\geq 0}V(\Lambda-2m\delta)^{\oplus p(m)} & \text{if
$\frak{g}= D^{(2)}_{n+1}$,}
\end{cases}
\end{equation}
where $p(m)$ denotes the number of partitions of $m$.
\end{cor}
\pf We will show that the weight of each maximal vector in
${\mathcal Z}(\Lambda)$ is of the form $\Lambda-m\delta$ (resp.
$\Lambda-2m\delta$) if $\frak{g}\neq D^{(2)}_{n+1}$ (resp.
$\frak{g}=D^{(2)}_{n+1}$) for some $m\geq 0$, and that there
exists a bijection between the set of partitions of $m$ ($m\geq
0$) and the set of maximal vectors in ${\mathcal Z}(\Lambda)$ with
weight $\Lambda-m\delta$ (resp. $\Lambda-2m\delta$) if
$\frak{g}\neq D^{(2)}_{n+1}$ (resp. $\frak{g}=D^{(2)}_{n+1}$). Let
$Y=(y_k)_{k=0}^{\infty}\in{\mathcal Z}(\Lambda)$ be a maximal
vector, i.e., $\tilde{e}_iY=\tilde{E}_iY=0$ for all $i\in I$.
Suppose that $Y$ is the ground-state wall $Y_{\Lambda}$. Since
${\rm wt}(Y_{\Lambda})=\Lambda$ and ${\mathcal
Z}(\Lambda)_{\Lambda}=\{\,Y_{\Lambda}\,\}$, the multiplicity of
$V(\Lambda)$ in ${\mathcal F}(\Lambda)$ is $1$. From now on, we
assume that $Y\neq Y_{\Lambda}$. Let $l$ be the maximum such that
$|y_l|\neq 0$. Suppose that $Y'=(y_k)_{k=l+1}^{\infty}\in{\mathcal
Z}(\Lambda_j)$ for some $j\in I$. Denote by $\Box$ the type of the
$j$-block. Let $\Delta$ be the volume of the $\delta$-column.

{\bf Case 1.} $\Box={\rm I}$

This case occurs only when $\frak{g}=A^{(1)}_n$. Since $y_{l+1}$
is $j$-admissible and $Y$ is a maximal vector, there is a
removable $j$-block on top of $y_l$, and hence  $y_l$ is obtained
by adding some $\delta$-columns to the ground-state wall, or
$|y_l|=m_l\Delta$ for some $m_l\geq 1$. If $l\neq 0$, let $l'$ be
the maximum such that $l'<l$ and $|y_{l'}|>|y_l|$. Note that
$y_{l'+1}$ is $j'$-admissible for some $j'\in I$. Since $Y$ is a
maximal vector, there exists an $j'$-block on top of $y_{l'}$,
which implies that $|y_{l'}|=m_{l'}\Delta$ for some $m_{l'}\geq
1$. By repeating the above argument column by column from left to
right, we conclude that for $0\leq k\leq l$, $|y_k|=m_k\Delta$ for
some $m_k\geq 1$. Moreover, $(m_0,m_1,\cdots,m_l)$ forms a
partition and ${\rm wt}(Y)=\Lambda-(\sum_{k=0}^lm_k)\delta$.\vskip
5mm

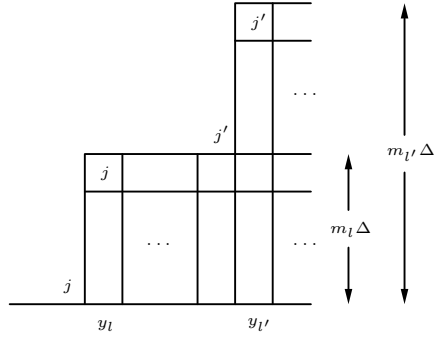
\begin{center}
\begin{texdraw}
\fontsize{6}{6}\drawdim em \setunitscale 0.13 \linewd 0.5
\arrowheadtype t:F \arrowheadsize l:4 w:2

\move(-20 0)\lvec(60 0)\move(0 0)\lvec(0 40)\lvec(60 40)\move(0
30)\lvec(60 30)\move(10 0)\lvec(10 40)

\move(40 0)\lvec(40 80)\lvec(60 80)\move(50 0)\lvec(50 80)\move(40
70)\lvec(60 70)

\move(30 0)\lvec(30 40)

\htext(-6 3){$j$}\htext(4 33){$j$}

\htext(34 43){$j'$}\htext(44 73){$j'$}

\htext(3 -7){$y_l$}\htext(43 -7){$y_{l'}$}

\htext(55 15){$\cdots$} \htext(55 55){$\cdots$}\htext(16
15){$\cdots$}

\move(70 25) \avec(70 40)

\move(70 15) \avec(70 0)

\htext(65 18){$m_l\Delta$}

\move(85 45) \avec(85 80)

\move(85 35) \avec(85 0)

\htext(80 38){$m_{l'}\Delta$}

\end{texdraw}
\end{center}\vskip 5mm

{\bf Case 2.} $\Box={\rm II}$

We see from the pattern for ${\mathcal Z}(\Lambda)$ that
$\Lambda=\Lambda_j$. By the maximality of $Y$, the $j$-signature
of $y_{l+1}$ is $+$ and the $j$-signature of $y_l$ is $-$ or $-+$.
Hence $y_l$ is obtained by adding some $\delta$-columns to the
ground-state wall. If $l\neq 0$, let $l'$ be the maximum such that
$l'<l$ and $|y_{l'}|>|y_l|$. Also by the maximality of $Y$, the
$j$-signature of $y_{l'}$ is $-$ or $-+$, which means that
$y_{l'}$ is obtained by adding some $\delta$-columns to the ground
state wall. Repeating the above argument from left to right, we
conclude that for $0\leq k\leq l$, the total volume of the blocks
added on the $k$th column is $m_k\Delta$ for some $m_k\geq 1$.
Hence, $(m_0,\cdots,m_l)$ forms a partition, and ${\rm
wt}(Y)=\Lambda-(\sum_{k=0}^lm_k)\delta$  if $\frak{g}\neq
D^{(2)}_{n+1}$, ${\rm wt}(Y)=\Lambda-2(\sum_{k=0}^lm_k)\delta$ if
$\frak{g}=D^{(2)}_{n+1}$.\vskip 5mm

\begin{center}
\begin{texdraw}
\fontsize{6}{6}\drawdim em \setunitscale 0.13 \linewd 0.5
\arrowheadtype t:F \arrowheadsize l:4 w:2

\move(60 5)\lvec(60 0)\lvec(50 0)\lvec(50 5)\move(60 5)\lvec(0
5)\lvec(0 0)\lvec(60 0)\lfill f:0.8

\move(20 0)\lvec(20 45)\move(40 50)\lvec(60 50)

\move(-10 0)\lvec(60 0)\move(10 0)\lvec(10 45)\lvec(60 45)
\move(10 40)\lvec(60 40)

\move(10 10)\lvec(60 10)

\move(40 0)\lvec(40 85)\lvec(60 85)\move(50 0)\lvec(50 85)\move(40
80)\lvec(60 80)

\htext(13 -7){$y_l$}\htext(43 -7){$y_{l'}$}

\htext(55 20){$\cdots$} \htext(55 60){$\cdots$}

\move(70 30) \avec(70 45)

\move(70 20) \avec(70 5)

\htext(65 23){$m_l\Delta$}

\move(85 55) \avec(85 85)

\move(85 40) \avec(85 5)

\htext(80 43){$m_{l'}\Delta$}

\htext(4 6){$_j$}\htext(14 41){$_j$}

\end{texdraw}
\end{center}\vskip 5mm

{\bf Case 3.} $\Box={\rm III}$

Let $j'$ be the color of the type III block, with which the
$j$-block forms a unit cube. By the maximality of $Y$, we observe
that the $y_{l+1}$ is $j$-admissible and $y_{l}$ is $j$-removable
but not $j'$-removable. Hence $y_l$ is obtained by adding some
$\delta$-columns to the ground-state wall. If $l\neq 0$, let $l'$
be the maximum such that $l'<l$ and $|y_{l'}|>|y_l|$. If
$y_{l'+1}$ is $j$-admissible, then by the maximality of $Y$,
$y_{l'}$ is $j$-removable but not $j'$-removable. On the other
hand, if $y_{l'+1}$ is $j'$-admissible, then by the maximality of
$Y$, $y_{l'}$ is $j'$-removable but not $j$-removable. As in Case
1 and 2, by repeating the above argument, we conclude that for
$0\leq k\leq l$, the total volume of the blocks added on the $k$th
column is $m_k\Delta$ for some $m_k\geq 1$. Hence,
$(m_0,m_1,\cdots,m_l)$ forms a partition and ${\rm
wt}(Y)=\Lambda-(\sum_{k=0}^lm_k)\delta$.\vskip 5mm

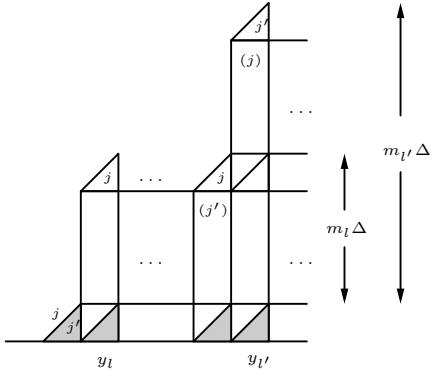
\begin{center}
\begin{texdraw}
\fontsize{6}{6}\drawdim em \setunitscale 0.13 \linewd 0.5

\move(-20 0)\lvec(60 0)\move(0 0)\lvec(0 40)\lvec(60 40)

\move(0 10)\lvec(60 10)\move(10 0)\lvec(10 40)\move(30 0)\lvec(30
40)\move(40 0)\lvec(40 80)\lvec(60 80)\move(50 0)\lvec(50 80)

\move(40 50)\lvec(60 50)

\move(-10 0)\lvec(0 10)\lvec(0 0)\lvec(-10 0)\lfill f:0.8

\move(0 0)\lvec(10 10)\lvec(10 0)\lvec(0 0)\lfill f:0.8

\move(30 0)\lvec(40 10)\lvec(40 0)\lvec(30 0)\lfill f:0.8

\move(40 0)\lvec(50 10)\lvec(50 0)\lvec(40 0)\lfill f:0.8

\move(0 40)\lvec(10 50)\lvec(10 40)\lvec(0 40)

\move(30 40)\lvec(40 50)\lvec(40 40)\lvec(30 40)

\move(40 40)\lvec(50 50)\lvec(50 40)\lvec(40 40)

\move(40 80)\lvec(50 90)\lvec(50 80)\lvec(40 80)

\htext(55 20){$\cdots$} \htext(55 60){$\cdots$}

\htext(15 20){$\cdots$} \htext(15 42){$\cdots$}

\move(70 35) \arrowheadtype t:F \arrowheadsize l:4 w:2 \avec(70
50)

\move(70 25) \arrowheadtype t:F \arrowheadsize l:4 w:2 \avec(70
10)

\htext(65 28){$m_l\Delta$}

\move(85 60) \arrowheadtype t:F \arrowheadsize l:4 w:2 \avec(85
90)

\move(85 45) \arrowheadtype t:F \arrowheadsize l:4 w:2 \avec(85
10)

\htext(80 48){$m_{l'}\Delta$}

\htext(4 -7){$y_l$}\htext(44 -7){$y_{l'}$}

\htext(-8 6){$_j$}\htext(6 42){$_j$}\htext(-4 2){$_{j'}$}

\htext(36 42){$_j$}\htext(31 33){$_{(j')}$}

\htext(46 82){$_{j'}$}\htext(42 73){$_{(j)}$}
\end{texdraw}
\end{center}\vskip 5mm

Conversely, for a given partition $(m_k)_{k=0}^{\infty}$ of a
nonnegative integer $m$, there exists a unique proper Young wall
$Y=(y_k)_{k=0}^{\infty}$ in ${\mathcal Z}(\Lambda)$ such that
$y_k$ is obtained by adding $m_k$ many $\delta$-columns to the
$k$th column of the ground state wall $Y_{\Lambda}$ (hence the
total volume of the blocks added to the $k$th column is
$m_k\Delta$). It is easy to check that $Y$ is a maximal vector
with ${\rm wt}(Y)=\Lambda-m\delta$ (resp. $\Lambda-2m\delta$) if
$\frak{g}\neq D^{(2)}_{n+1}$ (resp. $\frak{g}=D^{(2)}_{n+1}$).
\qed \vskip 1cm

\section{Generalized Lascoux-Leclerc-Thibon algorithm}
In this section, we generalize Lascoux-Leclerc-Thibon algorithm
(\cite{LLT}) to obtain an effective algorithm for constructing the
global basis $\mathcal{G}(\Lambda)$ of the basic representation
$V(\Lambda)$ of $U_q(\frak{g})$. Observe that $V(\Lambda)$ is
realized as the $U_q(\frak{g})$-submodule of ${\mathcal
F}(\Lambda)$ generated by the ground state wall $Y_{\Lambda}$.
Also recall that the crystal $B(\Lambda)$ of $V(\Lambda)$ is
isomorphic to the $U_q(\frak{g})$-crystal ${\mathcal Y}(\Lambda)$
consisting of reduced proper Young walls. Thus our goal is the
following: for each reduced proper Young wall $Y\in {\mathcal
Y}(\Lambda)$, we would like to give an algorithm of computing the
corresponding global basis element $G(Y)$ as a linear combination
of proper Young walls in ${\mathcal Z}(\Lambda)$.

For this purpose, we first investigate the action of divided
powers $f_i^{(r)}$ ($i\in I, r\geq 1$) on the proper Young walls.
Let $Y$ be a proper Young wall in ${\mathcal Y}(\Lambda)$ (not
necessarily reduced), and write
\begin{equation}
f_i^{(r)}Y=\sum_{\substack{Z\in {\mathcal Z}(\Lambda) \\ {\rm
wt}(Z)={\rm wt}(Y)-r\alpha_i }}Q_{Y,Z}(q)Z,
\end{equation}
where $Q_{Y,Z}(q)\in \mathbb{Q}(q)$. For each
$Z=(z_k)_{k=0}^{\infty}\in {\mathcal Z}(\Lambda)$ with
$Q_{Y,Z}(q)\neq 0$, there exists a unique sequence of proper Young
walls $Y=Y_{0},Y_{1},\cdots,Y_{r}=Z$ such that

\begin{itemize}
\item[(i)] $\lambda_{k+1}Y_{k+1}=Y_{k}\swarrow b_{k+1}$ for some
$\lambda_{k+1}\in\mathbb{Z}[q,q^{-1}]$
 and a (virtually) admissible $i$-slot $b_{k+1}$ of $Y_{k}$,
\item[(ii)] $b_{k+1}$ is placed on top of $b_k$ or to the right of $b_k$.
\end{itemize}

\begin{ex}\label{Y1Y2Y3}{\rm
Let $\frak{g}=A^{(2)}_4$, $\Lambda=\Lambda_0$, $i=0$, and consider
\vskip 5mm

\begin{center}
$Y=$\raisebox{-0.5\height}{
\begin{texdraw}
\drawdim em \setunitscale 0.13 \linewd 0.5


\move(0 0)\lvec(10 0)\lvec(10 10)\lvec(0 10)\lvec(0 0)\htext(3
1){\tiny $0$}

\move(10 0)\lvec(20 0)\lvec(20 10)\lvec(10 10)\lvec(10 0)\htext(13
1){\tiny $0$}

\move(20 0)\lvec(30 0)\lvec(30 10)\lvec(20 10)\lvec(20 0)\htext(23
1){\tiny $0$}

\move(30 0)\lvec(40 0)\lvec(40 10)\lvec(30 10)\lvec(30 0)\htext(33
1){\tiny $0$}

\move(40 0)\lvec(50 0)\lvec(50 10)\lvec(40 10)\lvec(40 0)\htext(43
1){\tiny $0$}

\move(0 0)\lvec(10 0)\lvec(10 5)\lvec(0 5)\lvec(0 0)\lfill f:0.8
\htext(3 6){\tiny $0$}

\move(10 0)\lvec(20 0)\lvec(20 5)\lvec(10 5)\lvec(10 0)\lfill
f:0.8 \htext(13 6){\tiny $0$}

\move(20 0)\lvec(30 0)\lvec(30 5)\lvec(20 5)\lvec(20 0)\lfill
f:0.8 \htext(23 6){\tiny $0$}

\move(30 0)\lvec(40 0)\lvec(40 5)\lvec(30 5)\lvec(30 0)\lfill
f:0.8 \htext(33 6){\tiny $0$}

\move(40 0)\lvec(50 0)\lvec(50 5)\lvec(40 5)\lvec(40 0)\lfill
f:0.8 \htext(43 6){\tiny $0$}

\move(0 10)\lvec(10 10)\lvec(10 20)\lvec(0 20)\lvec(0 10)\htext(3
13){$_1$}

\move(10 10)\lvec(20 10)\lvec(20 20)\lvec(10 20)\lvec(10
10)\htext(13 13){$_1$}

\move(20 10)\lvec(30 10)\lvec(30 20)\lvec(20 20)\lvec(20
10)\htext(23 13){$_1$}

\move(30 10)\lvec(40 10)\lvec(40 20)\lvec(30 20)\lvec(30
10)\htext(33 13){$_1$}

\move(40 10)\lvec(50 10)\lvec(50 20)\lvec(40 20)\lvec(40
10)\htext(43 13){$_1$}

\move(0 20)\lvec(10 20)\lvec(10 30)\lvec(0 30)\lvec(0 20)\htext(3
23){$_2$}

\move(10 20)\lvec(20 20)\lvec(20 30)\lvec(10 30)\lvec(10
20)\htext(13 23){$_2$}

\move(20 20)\lvec(30 20)\lvec(30 30)\lvec(20 30)\lvec(20
20)\htext(23 23){$_2$}

\move(30 20)\lvec(40 20)\lvec(40 30)\lvec(30 30)\lvec(30
20)\htext(33 23){$_2$}

\move(40 20)\lvec(50 20)\lvec(50 30)\lvec(40 30)\lvec(40
20)\htext(43 23){$_2$}

\move(0 30)\lvec(10 30)\lvec(10 40)\lvec(0 40)\lvec(0 30)\htext(3
33){$_1$}

\move(10 30)\lvec(20 30)\lvec(20 40)\lvec(10 40)\lvec(10
30)\htext(13 33){$_1$}

\move(20 30)\lvec(30 30)\lvec(30 40)\lvec(20 40)\lvec(20
30)\htext(23 33){$_1$}

\move(30 30)\lvec(40 30)\lvec(40 40)\lvec(30 40)\lvec(30
30)\htext(33 33){$_1$}

\move(40 30)\lvec(50 30)\lvec(50 40)\lvec(40 40)\lvec(40
30)\htext(43 33){$_1$}

\move(10 40)\lvec(20 40)\lvec(20 45)\lvec(10 45)\lvec(10
40)\htext(13 41){\tiny $0$}

\move(20 40)\lvec(30 40)\lvec(30 45)\lvec(20 45)\lvec(20
40)\htext(23 41){\tiny $0$}

\move(30 40)\lvec(40 40)\lvec(40 45)\lvec(30 45)\lvec(30
40)\htext(33 41){\tiny $0$}

\move(40 40)\lvec(50 40)\lvec(50 45)\lvec(40 45)\lvec(40
40)\htext(43 41){\tiny $0$}

%

\move(40 40)\lvec(50 40)\lvec(50 50)\lvec(40 50)\lvec(40
40)\htext(43 46){\tiny $0$}

\move(40 50)\lvec(50 50)\lvec(50 60)\lvec(40 60)\lvec(40
50)\htext(43 53){$_1$}

\move(40 60)\lvec(50 60)\lvec(50 70)\lvec(40 70)\lvec(40
60)\htext(43 63){$_2$}

\move(40 70)\lvec(50 70)\lvec(50 80)\lvec(40 80)\lvec(40
70)\htext(43 73){$_1$}

%
\end{texdraw}}\hskip 1cm and  \hskip 1cm
$Z=$\raisebox{-0.5\height}{
\begin{texdraw}
\drawdim em \setunitscale 0.13 \linewd 0.5

\move(0 0)\lvec(10 0)\lvec(10 10)\lvec(0 10)\lvec(0 0)\htext(3
1){\tiny $0$}

\move(10 0)\lvec(20 0)\lvec(20 10)\lvec(10 10)\lvec(10 0)\htext(13
1){\tiny $0$}

\move(20 0)\lvec(30 0)\lvec(30 10)\lvec(20 10)\lvec(20 0)\htext(23
1){\tiny $0$}

\move(30 0)\lvec(40 0)\lvec(40 10)\lvec(30 10)\lvec(30 0)\htext(33
1){\tiny $0$}

\move(40 0)\lvec(50 0)\lvec(50 10)\lvec(40 10)\lvec(40 0)\htext(43
1){\tiny $0$}
\move(0 0)\lvec(10 0)\lvec(10 5)\lvec(0 5)\lvec(0 0)\lfill f:0.8
\htext(3 6){\tiny $0$}

\move(10 0)\lvec(20 0)\lvec(20 5)\lvec(10 5)\lvec(10 0)\lfill
f:0.8 \htext(13 6){\tiny $0$}

\move(20 0)\lvec(30 0)\lvec(30 5)\lvec(20 5)\lvec(20 0)\lfill
f:0.8 \htext(23 6){\tiny $0$}

\move(30 0)\lvec(40 0)\lvec(40 5)\lvec(30 5)\lvec(30 0)\lfill
f:0.8 \htext(33 6){\tiny $0$}

\move(40 0)\lvec(50 0)\lvec(50 5)\lvec(40 5)\lvec(40 0)\lfill
f:0.8 \htext(43 6){\tiny $0$}
\move(0 10)\lvec(10 10)\lvec(10 20)\lvec(0 20)\lvec(0 10)\htext(3
13){$_1$}

\move(10 10)\lvec(20 10)\lvec(20 20)\lvec(10 20)\lvec(10
10)\htext(13 13){$_1$}

\move(20 10)\lvec(30 10)\lvec(30 20)\lvec(20 20)\lvec(20
10)\htext(23 13){$_1$}

\move(30 10)\lvec(40 10)\lvec(40 20)\lvec(30 20)\lvec(30
10)\htext(33 13){$_1$}

\move(40 10)\lvec(50 10)\lvec(50 20)\lvec(40 20)\lvec(40
10)\htext(43 13){$_1$}
\move(0 20)\lvec(10 20)\lvec(10 30)\lvec(0 30)\lvec(0 20)\htext(3
23){$_2$}

\move(10 20)\lvec(20 20)\lvec(20 30)\lvec(10 30)\lvec(10
20)\htext(13 23){$_2$}

\move(20 20)\lvec(30 20)\lvec(30 30)\lvec(20 30)\lvec(20
20)\htext(23 23){$_2$}

\move(30 20)\lvec(40 20)\lvec(40 30)\lvec(30 30)\lvec(30
20)\htext(33 23){$_2$}

\move(40 20)\lvec(50 20)\lvec(50 30)\lvec(40 30)\lvec(40
20)\htext(43 23){$_2$}
\move(0 30)\lvec(10 30)\lvec(10 40)\lvec(0 40)\lvec(0 30)\htext(3
33){$_1$}

\move(10 30)\lvec(20 30)\lvec(20 40)\lvec(10 40)\lvec(10
30)\htext(13 33){$_1$}

\move(20 30)\lvec(30 30)\lvec(30 40)\lvec(20 40)\lvec(20
30)\htext(23 33){$_1$}

\move(30 30)\lvec(40 30)\lvec(40 40)\lvec(30 40)\lvec(30
30)\htext(33 33){$_1$}

\move(40 30)\lvec(50 30)\lvec(50 40)\lvec(40 40)\lvec(40
30)\htext(43 33){$_1$}
\move(0 40)\lvec(10 40)\lvec(10 45)\lvec(0 45)\lvec(0 40)\htext(3
41){\tiny $0$}

\move(10 40)\lvec(20 40)\lvec(20 45)\lvec(10 45)\lvec(10
40)\htext(13 41){\tiny $0$}

\move(20 40)\lvec(30 40)\lvec(30 45)\lvec(20 45)\lvec(20
40)\htext(23 41){\tiny $0$}

\move(30 40)\lvec(40 40)\lvec(40 45)\lvec(30 45)\lvec(30
40)\htext(33 41){\tiny $0$}

\move(40 40)\lvec(50 40)\lvec(50 45)\lvec(40 45)\lvec(40
40)\htext(43 41){\tiny $0$}
\move(30 40)\lvec(40 40)\lvec(40 50)\lvec(30 50)\lvec(30
40)\htext(33 46){\tiny $0$}

\move(40 40)\lvec(50 40)\lvec(50 50)\lvec(40 50)\lvec(40
40)\htext(43 46){\tiny $0$}
\move(40 50)\lvec(50 50)\lvec(50 60)\lvec(40 60)\lvec(40
50)\htext(43 53){$_1$}
\move(40 60)\lvec(50 60)\lvec(50 70)\lvec(40 70)\lvec(40
60)\htext(43 63){$_2$}
\move(40 70)\lvec(50 70)\lvec(50 80)\lvec(40 80)\lvec(40
70)\htext(43 73){$_1$}
\move(40 80)\lvec(50 80)\lvec(50 85)\lvec(40 85)\lvec(40
80)\htext(43 81){\tiny $0$}

\move(40 85)\lvec(50 85)\lvec(50 90)\lvec(40 90)\lvec(40
85)\htext(43 86){\tiny $0$}
\move(5 54) \arrowheadtype t:F \arrowheadsize l:4 w:2 \avec(5 44)
\htext(2 55){$_{b_1}$}

\move(35 59) \arrowheadtype t:F \arrowheadsize l:4 w:2 \avec(35
49) \htext(32 60){$_{b_2}$}

\move(59 83) \arrowheadtype t:F \arrowheadsize l:4 w:2 \avec(49
83) \htext(60 81){$_{b_3}$}

\move(62 88) \arrowheadtype t:F \arrowheadsize l:4 w:2 \avec(49
88) \htext(63 86){$_{b_4}$}
\end{texdraw}}\quad .
\end{center}\vskip 5mm
Then, we have a sequence of proper Young walls
$Y=Y_0,Y_1,Y_2,Y_3,Y_4=Z$, where\vskip 5mm

\begin{center}
$Y_1=$\raisebox{-0.5\height}{
\begin{texdraw}
\drawdim em \setunitscale 0.13 \linewd 0.5

\move(0 0)\lvec(10 0)\lvec(10 10)\lvec(0 10)\lvec(0 0)\htext(3
1){\tiny $0$}

\move(10 0)\lvec(20 0)\lvec(20 10)\lvec(10 10)\lvec(10 0)\htext(13
1){\tiny $0$}

\move(20 0)\lvec(30 0)\lvec(30 10)\lvec(20 10)\lvec(20 0)\htext(23
1){\tiny $0$}

\move(30 0)\lvec(40 0)\lvec(40 10)\lvec(30 10)\lvec(30 0)\htext(33
1){\tiny $0$}

\move(40 0)\lvec(50 0)\lvec(50 10)\lvec(40 10)\lvec(40 0)\htext(43
1){\tiny $0$}
\move(0 0)\lvec(10 0)\lvec(10 5)\lvec(0 5)\lvec(0 0)\lfill f:0.8
\htext(3 6){\tiny $0$}

\move(10 0)\lvec(20 0)\lvec(20 5)\lvec(10 5)\lvec(10 0)\lfill
f:0.8 \htext(13 6){\tiny $0$}

\move(20 0)\lvec(30 0)\lvec(30 5)\lvec(20 5)\lvec(20 0)\lfill
f:0.8 \htext(23 6){\tiny $0$}

\move(30 0)\lvec(40 0)\lvec(40 5)\lvec(30 5)\lvec(30 0)\lfill
f:0.8 \htext(33 6){\tiny $0$}

\move(40 0)\lvec(50 0)\lvec(50 5)\lvec(40 5)\lvec(40 0)\lfill
f:0.8 \htext(43 6){\tiny $0$}
\move(0 10)\lvec(10 10)\lvec(10 20)\lvec(0 20)\lvec(0 10)\htext(3
13){$_1$}

\move(10 10)\lvec(20 10)\lvec(20 20)\lvec(10 20)\lvec(10
10)\htext(13 13){$_1$}

\move(20 10)\lvec(30 10)\lvec(30 20)\lvec(20 20)\lvec(20
10)\htext(23 13){$_1$}

\move(30 10)\lvec(40 10)\lvec(40 20)\lvec(30 20)\lvec(30
10)\htext(33 13){$_1$}

\move(40 10)\lvec(50 10)\lvec(50 20)\lvec(40 20)\lvec(40
10)\htext(43 13){$_1$}
\move(0 20)\lvec(10 20)\lvec(10 30)\lvec(0 30)\lvec(0 20)\htext(3
23){$_2$}

\move(10 20)\lvec(20 20)\lvec(20 30)\lvec(10 30)\lvec(10
20)\htext(13 23){$_2$}

\move(20 20)\lvec(30 20)\lvec(30 30)\lvec(20 30)\lvec(20
20)\htext(23 23){$_2$}

\move(30 20)\lvec(40 20)\lvec(40 30)\lvec(30 30)\lvec(30
20)\htext(33 23){$_2$}

\move(40 20)\lvec(50 20)\lvec(50 30)\lvec(40 30)\lvec(40
20)\htext(43 23){$_2$}
\move(0 30)\lvec(10 30)\lvec(10 40)\lvec(0 40)\lvec(0 30)\htext(3
33){$_1$}

\move(10 30)\lvec(20 30)\lvec(20 40)\lvec(10 40)\lvec(10
30)\htext(13 33){$_1$}

\move(20 30)\lvec(30 30)\lvec(30 40)\lvec(20 40)\lvec(20
30)\htext(23 33){$_1$}

\move(30 30)\lvec(40 30)\lvec(40 40)\lvec(30 40)\lvec(30
30)\htext(33 33){$_1$}

\move(40 30)\lvec(50 30)\lvec(50 40)\lvec(40 40)\lvec(40
30)\htext(43 33){$_1$}
\move(0 40)\lvec(10 40)\lvec(10 45)\lvec(0 45)\lvec(0 40)\htext(3
41){\tiny $0$}

\move(10 40)\lvec(20 40)\lvec(20 45)\lvec(10 45)\lvec(10
40)\htext(13 41){\tiny $0$}

\move(20 40)\lvec(30 40)\lvec(30 45)\lvec(20 45)\lvec(20
40)\htext(23 41){\tiny $0$}

\move(30 40)\lvec(40 40)\lvec(40 45)\lvec(30 45)\lvec(30
40)\htext(33 41){\tiny $0$}

\move(40 40)\lvec(50 40)\lvec(50 45)\lvec(40 45)\lvec(40
40)\htext(43 41){\tiny $0$}

\move(40 40)\lvec(50 40)\lvec(50 50)\lvec(40 50)\lvec(40
40)\htext(43 46){\tiny $0$}
\move(40 50)\lvec(50 50)\lvec(50 60)\lvec(40 60)\lvec(40
50)\htext(43 53){$_1$}
\move(40 60)\lvec(50 60)\lvec(50 70)\lvec(40 70)\lvec(40
60)\htext(43 63){$_2$}
\move(40 70)\lvec(50 70)\lvec(50 80)\lvec(40 80)\lvec(40
70)\htext(43 73){$_1$}
\end{texdraw}}\quad,\quad
$Y_2=$\raisebox{-0.5\height}{
\begin{texdraw}
\drawdim em \setunitscale 0.13 \linewd 0.5

\move(0 0)\lvec(10 0)\lvec(10 10)\lvec(0 10)\lvec(0 0)\htext(3
1){\tiny $0$}

\move(10 0)\lvec(20 0)\lvec(20 10)\lvec(10 10)\lvec(10 0)\htext(13
1){\tiny $0$}

\move(20 0)\lvec(30 0)\lvec(30 10)\lvec(20 10)\lvec(20 0)\htext(23
1){\tiny $0$}

\move(30 0)\lvec(40 0)\lvec(40 10)\lvec(30 10)\lvec(30 0)\htext(33
1){\tiny $0$}

\move(40 0)\lvec(50 0)\lvec(50 10)\lvec(40 10)\lvec(40 0)\htext(43
1){\tiny $0$}
\move(0 0)\lvec(10 0)\lvec(10 5)\lvec(0 5)\lvec(0 0)\lfill f:0.8
\htext(3 6){\tiny $0$}

\move(10 0)\lvec(20 0)\lvec(20 5)\lvec(10 5)\lvec(10 0)\lfill
f:0.8 \htext(13 6){\tiny $0$}

\move(20 0)\lvec(30 0)\lvec(30 5)\lvec(20 5)\lvec(20 0)\lfill
f:0.8 \htext(23 6){\tiny $0$}

\move(30 0)\lvec(40 0)\lvec(40 5)\lvec(30 5)\lvec(30 0)\lfill
f:0.8 \htext(33 6){\tiny $0$}

\move(40 0)\lvec(50 0)\lvec(50 5)\lvec(40 5)\lvec(40 0)\lfill
f:0.8 \htext(43 6){\tiny $0$}
\move(0 10)\lvec(10 10)\lvec(10 20)\lvec(0 20)\lvec(0 10)\htext(3
13){$_1$}

\move(10 10)\lvec(20 10)\lvec(20 20)\lvec(10 20)\lvec(10
10)\htext(13 13){$_1$}

\move(20 10)\lvec(30 10)\lvec(30 20)\lvec(20 20)\lvec(20
10)\htext(23 13){$_1$}

\move(30 10)\lvec(40 10)\lvec(40 20)\lvec(30 20)\lvec(30
10)\htext(33 13){$_1$}

\move(40 10)\lvec(50 10)\lvec(50 20)\lvec(40 20)\lvec(40
10)\htext(43 13){$_1$}
\move(0 20)\lvec(10 20)\lvec(10 30)\lvec(0 30)\lvec(0 20)\htext(3
23){$_2$}

\move(10 20)\lvec(20 20)\lvec(20 30)\lvec(10 30)\lvec(10
20)\htext(13 23){$_2$}

\move(20 20)\lvec(30 20)\lvec(30 30)\lvec(20 30)\lvec(20
20)\htext(23 23){$_2$}

\move(30 20)\lvec(40 20)\lvec(40 30)\lvec(30 30)\lvec(30
20)\htext(33 23){$_2$}

\move(40 20)\lvec(50 20)\lvec(50 30)\lvec(40 30)\lvec(40
20)\htext(43 23){$_2$}
\move(0 30)\lvec(10 30)\lvec(10 40)\lvec(0 40)\lvec(0 30)\htext(3
33){$_1$}

\move(10 30)\lvec(20 30)\lvec(20 40)\lvec(10 40)\lvec(10
30)\htext(13 33){$_1$}

\move(20 30)\lvec(30 30)\lvec(30 40)\lvec(20 40)\lvec(20
30)\htext(23 33){$_1$}

\move(30 30)\lvec(40 30)\lvec(40 40)\lvec(30 40)\lvec(30
30)\htext(33 33){$_1$}

\move(40 30)\lvec(50 30)\lvec(50 40)\lvec(40 40)\lvec(40
30)\htext(43 33){$_1$}
\move(0 40)\lvec(10 40)\lvec(10 45)\lvec(0 45)\lvec(0 40)\htext(3
41){\tiny $0$}

\move(10 40)\lvec(20 40)\lvec(20 45)\lvec(10 45)\lvec(10
40)\htext(13 41){\tiny $0$}

\move(20 40)\lvec(30 40)\lvec(30 45)\lvec(20 45)\lvec(20
40)\htext(23 41){\tiny $0$}

\move(30 40)\lvec(40 40)\lvec(40 45)\lvec(30 45)\lvec(30
40)\htext(33 41){\tiny $0$}

\move(40 40)\lvec(50 40)\lvec(50 45)\lvec(40 45)\lvec(40
40)\htext(43 41){\tiny $0$}
\move(30 40)\lvec(40 40)\lvec(40 50)\lvec(30 50)\lvec(30
40)\htext(33 46){\tiny $0$}

\move(40 40)\lvec(50 40)\lvec(50 50)\lvec(40 50)\lvec(40
40)\htext(43 46){\tiny $0$}
\move(40 50)\lvec(50 50)\lvec(50 60)\lvec(40 60)\lvec(40
50)\htext(43 53){$_1$}
\move(40 60)\lvec(50 60)\lvec(50 70)\lvec(40 70)\lvec(40
60)\htext(43 63){$_2$}
\move(40 70)\lvec(50 70)\lvec(50 80)\lvec(40 80)\lvec(40
70)\htext(43 73){$_1$}
\end{texdraw}}\quad , \quad
$Y_3=$\raisebox{-0.5\height}{
\begin{texdraw}
\drawdim em \setunitscale 0.13 \linewd 0.5

\move(0 0)\lvec(10 0)\lvec(10 10)\lvec(0 10)\lvec(0 0)\htext(3
1){\tiny $0$}

\move(10 0)\lvec(20 0)\lvec(20 10)\lvec(10 10)\lvec(10 0)\htext(13
1){\tiny $0$}

\move(20 0)\lvec(30 0)\lvec(30 10)\lvec(20 10)\lvec(20 0)\htext(23
1){\tiny $0$}

\move(30 0)\lvec(40 0)\lvec(40 10)\lvec(30 10)\lvec(30 0)\htext(33
1){\tiny $0$}

\move(40 0)\lvec(50 0)\lvec(50 10)\lvec(40 10)\lvec(40 0)\htext(43
1){\tiny $0$}
\move(0 0)\lvec(10 0)\lvec(10 5)\lvec(0 5)\lvec(0 0)\lfill f:0.8
\htext(3 6){\tiny $0$}

\move(10 0)\lvec(20 0)\lvec(20 5)\lvec(10 5)\lvec(10 0)\lfill
f:0.8 \htext(13 6){\tiny $0$}

\move(20 0)\lvec(30 0)\lvec(30 5)\lvec(20 5)\lvec(20 0)\lfill
f:0.8 \htext(23 6){\tiny $0$}

\move(30 0)\lvec(40 0)\lvec(40 5)\lvec(30 5)\lvec(30 0)\lfill
f:0.8 \htext(33 6){\tiny $0$}

\move(40 0)\lvec(50 0)\lvec(50 5)\lvec(40 5)\lvec(40 0)\lfill
f:0.8 \htext(43 6){\tiny $0$}
\move(0 10)\lvec(10 10)\lvec(10 20)\lvec(0 20)\lvec(0 10)\htext(3
13){$_1$}

\move(10 10)\lvec(20 10)\lvec(20 20)\lvec(10 20)\lvec(10
10)\htext(13 13){$_1$}

\move(20 10)\lvec(30 10)\lvec(30 20)\lvec(20 20)\lvec(20
10)\htext(23 13){$_1$}

\move(30 10)\lvec(40 10)\lvec(40 20)\lvec(30 20)\lvec(30
10)\htext(33 13){$_1$}

\move(40 10)\lvec(50 10)\lvec(50 20)\lvec(40 20)\lvec(40
10)\htext(43 13){$_1$}
\move(0 20)\lvec(10 20)\lvec(10 30)\lvec(0 30)\lvec(0 20)\htext(3
23){$_2$}

\move(10 20)\lvec(20 20)\lvec(20 30)\lvec(10 30)\lvec(10
20)\htext(13 23){$_2$}

\move(20 20)\lvec(30 20)\lvec(30 30)\lvec(20 30)\lvec(20
20)\htext(23 23){$_2$}

\move(30 20)\lvec(40 20)\lvec(40 30)\lvec(30 30)\lvec(30
20)\htext(33 23){$_2$}

\move(40 20)\lvec(50 20)\lvec(50 30)\lvec(40 30)\lvec(40
20)\htext(43 23){$_2$}
\move(0 30)\lvec(10 30)\lvec(10 40)\lvec(0 40)\lvec(0 30)\htext(3
33){$_1$}

\move(10 30)\lvec(20 30)\lvec(20 40)\lvec(10 40)\lvec(10
30)\htext(13 33){$_1$}

\move(20 30)\lvec(30 30)\lvec(30 40)\lvec(20 40)\lvec(20
30)\htext(23 33){$_1$}

\move(30 30)\lvec(40 30)\lvec(40 40)\lvec(30 40)\lvec(30
30)\htext(33 33){$_1$}

\move(40 30)\lvec(50 30)\lvec(50 40)\lvec(40 40)\lvec(40
30)\htext(43 33){$_1$}
\move(0 40)\lvec(10 40)\lvec(10 45)\lvec(0 45)\lvec(0 40)\htext(3
41){\tiny $0$}

\move(10 40)\lvec(20 40)\lvec(20 45)\lvec(10 45)\lvec(10
40)\htext(13 41){\tiny $0$}

\move(20 40)\lvec(30 40)\lvec(30 45)\lvec(20 45)\lvec(20
40)\htext(23 41){\tiny $0$}

\move(30 40)\lvec(40 40)\lvec(40 45)\lvec(30 45)\lvec(30
40)\htext(33 41){\tiny $0$}

\move(40 40)\lvec(50 40)\lvec(50 45)\lvec(40 45)\lvec(40
40)\htext(43 41){\tiny $0$}
\move(30 40)\lvec(40 40)\lvec(40 50)\lvec(30 50)\lvec(30
40)\htext(33 46){\tiny $0$}

\move(40 40)\lvec(50 40)\lvec(50 50)\lvec(40 50)\lvec(40
40)\htext(43 46){\tiny $0$}
\move(40 50)\lvec(50 50)\lvec(50 60)\lvec(40 60)\lvec(40
50)\htext(43 53){$_1$}
\move(40 60)\lvec(50 60)\lvec(50 70)\lvec(40 70)\lvec(40
60)\htext(43 63){$_2$}
\move(40 70)\lvec(50 70)\lvec(50 80)\lvec(40 80)\lvec(40
70)\htext(43 73){$_1$}
\move(40 80)\lvec(50 80)\lvec(50 85)\lvec(40 85)\lvec(40
80)\htext(43 81){\tiny $0$}
\end{texdraw}}

\end{center}\vskip 5mm

\noindent and $\lambda_1=1$, $\lambda_2=q^{-1}(1-q^8)$,
$\lambda_3=1$, $\lambda_4=[2]=(q+q^{-1})$. }
\end{ex}\vskip 3mm

Let $Q_{Y_k,Y_{k+1}}(q)$ be the coefficient of $Y_{k+1}$ in the
expression of $f_iY_{k}$.

Define
\begin{equation}
Q^{\circ}_{Y,Z}(q)=\prod_{k=0}^{r-1}Q_{Y_{k},Y_{k+1}}(q)\in\mathbb{Z}[q,q^{-1}].
\end{equation}
Note that each $b_k$ ($1\leq k\leq r$) can be viewed as an
$i$-block (not necessarily removable) in $Z$. When $b_k$'s are of
type ${\rm II}$, we define
\begin{equation}\label {JST}
\begin{split}
J_1&=\{\,k\,|\,b_{k-1} \text{ lies beneath } b_k\,\}, \\
J_2&=\{\,k\,|\,\text{there exists an $i$-block {\rm(}$\neq
b_{k-1}${\rm )} beneath $b_k$ in $Z$}\,\}, \\
J_3&=\{\,k\,|\,\text{there exists no $i$-block on top of and
beneath
$b_k$ in $Z$}\,\}, \\
S&=\{\,k\in J_2\,|\, \text{$b_k$ and $b_{k-1}$ lie in the same $i$-component of $Z$} \,\}, \\
T&=\{\,k\in J_3\,|\,k+1\in S\,\}
\end{split}
\end{equation}
Set $l=|J_1|$, $m=|J_2|$ and $n=|J_3|$. For each $k\in S$, let
$\mu_k=q\lambda_k$. Note that $2l+m+n=r$. Then we have

\begin{lem}\label {QYZ} Let $Y$ be a proper Young wall in ${\mathcal Z}(\Lambda)$,
and suppose that $Q_{Y,Z}(q)\neq 0$ for some $Z\in {\mathcal
Z}(\Lambda)$ with ${\rm wt}(Z)={\rm wt}(Y)-r\alpha_i$. Then we
have
\begin{equation}\label {qyz}
Q_{Y,Z}(q)=
\begin{cases}
Q^{\circ}_{Y,Z}(q)q_i^{\binom{r}{2}}
& \text{ if the $i$-blocks are of type {\rm I} or {\rm III}}, \\
Q^{\circ}_{Y,Z}(q) \dfrac{q^{\sigma(l,m,n)}}{[2]^{l}\prod_{k\in
S}\mu_k} & \text{ if the $i$-blocks are of type {\rm II}},
\end{cases}
\end{equation}
where
$\sigma(l,m,n)=4\binom{l}{2}+\binom{m}{2}+\binom{n}{2}+2l(m+n)+mn$.

In particular, we have $Q_{Y,Z}(q)\in\mathbb{Z}[q,q^{-1}]$.
\end{lem}

\begin{ex}{\rm In Example \ref{Y1Y2Y3}, we have
$J_1=\{4\}$, $J_2=\{2\}$ and $J_3=\{1\}$. Also, we have $S=\{2\}$,
$T=\{1\}$. Since $Q_{Y,Z}^{\circ}(q)=q^{-1}(1-q^8)[2]$,
$\mu_2=q\lambda_2=(1-q^8)$ and $\sigma(1,1,1)=5$, it follows that
\begin{equation*}
Q_{Y,Z}(q)=Q_{Y,Z}^{\circ}(q)\times \frac{q^5}{[2](1-q^8)}=q^4.
\end{equation*}
Hence, the coefficient of $Z$ in $f_0^{(4)}Y$ is $q^4$.

}
\end{ex}

{\bf Proof of Lemma \ref{QYZ}}

{\bf Case 1.} Suppose that the $i$-blocks are of type I or III. We
will use induction on $r$. For $r=1$, it is clear. Suppose that
\eqref{qyz} holds for $r-1$. Let $Z_{k}$ be the proper Young wall
in ${\mathcal Z}(\Lambda)$ ($1\leq k\leq r$) such that $b_k$ is
also a (virtually) admissible $i$-slot of $Z_k$ and $Z_{k}\swarrow
b_k=Z$ (up to scalar multiplication). By definition of
$Q_{Y,Z}(q)$, we have
\begin{equation*}
[r]_iQ_{Y,Z}(q)=\sum_{k=1}^rQ_{Y,Z_{k}}(q)Q_{Z_{k},Z}(q).
\end{equation*}
Note that
\begin{equation*}
\begin{split}
Q_{Z_{k},Z}(q)&=Q_{Z_{k},Z}^{\circ}(q), \\
Q_{Y,Z_{k}}^{\circ}(q)Q_{Z_{k},Z}^{\circ}(q)&=Q_{Y,Z}^{\circ}(q)q_i^{2(r-k)}.
\end{split}
\end{equation*}
By induction hypothesis, we have
\begin{equation*}
\begin{split}
[r]_iQ_{Y,Z}(q)&=\sum_{k=1}^rQ_{Y,Z_{k}}^{\circ}
(q)Q_{Z_{k},Z}^{\circ}(q)q_i^{\binom{r-1}{2}} \\
&=\sum_{k=1}^rQ_{Y,Z}^{\circ}(q)q_i^{2(r-k)+\binom{r-1}{2}} \\
&=\sum_{k=1}^rQ_{Y,Z}^{\circ}(q)q_i^{2(r-k)-r+1+\binom{r}{2}}
=Q_{Y,Z}^{\circ}(q)q_i^{\binom{r}{2}}[r]_i.
\end{split}
\end{equation*}
This completes the induction argument.\vskip 3mm

{\bf Case 2.} Suppose that $i$-blocks are of type II. We will also
use induction on $r$. For $r=1$, it is clear. Suppose that
\eqref{qyz} holds for $r-1$ ($r\geq 2$). Set $J=J_1\sqcup J_2
\sqcup J_3$. For $k\in J$, let $Z_{k}$ be the unique proper Young
wall in ${\mathcal Z}(\Lambda)$ such that $b_k$ is also an
admissible $i$-slot of $Z_k$ and $Z_k\swarrow b_k=\lambda_kZ$
($\lambda_k$ was already given in the definition of $b_k$). We
have
\begin{equation}
\lambda_k=
\begin{cases}
[2] & \text{if $k\in J_1$}, \\
q^{-1}(1-(-q^2)^{l_k+1}) & \text{if $k\in J_2$ and $k\geq 2$}, \\
1 & \text{otherwise}
\end{cases}
\end{equation}
for some $l_k\geq 0$, and $\mu_k=(1-(-q^2)^{l_k+1})$ for $k\in S$.
As in Case 1, we have
\begin{equation}
[r]_iQ_{Y,Z}(q)=\sum_{k\in J}Q_{Y,Z_{k}}(q)Q_{Z_k,Z}(q),
\end{equation}
and
\begin{equation*}
\begin{split}
Q_{Z_{k},Z}(q)&=Q_{Z_{k},Z}^{\circ}(q), \\
Q_{Y,Z_{k}}^{\circ}(q)Q_{Z_{k},Z}^{\circ}(q)&=
\begin{cases}
Q_{Y,Z}^{\circ}(q)q^{2(r-k)}\frac{(1-(-q^2)^{l_{k+1}})}{(1-(-q^2)^{l_{k+1}+1})}
& \text{if $k\in T$},  \\
Q_{Y,Z}^{\circ}(q)q^{2(r-k)} & \text{otherwise}.
\end{cases}
\end{split}
\end{equation*}

By induction hypothesis, we have {\allowdisplaybreaks
\begin{align*}
[r]_iQ_{Y,Z}(q)&=Q_{Y,Z}^{\circ}(q)\times\\
& \Bigg( \sum_{k\in
J_1}\frac{q^{\sigma(l-1,m,n+1)+2(r-k)}}{[2]^{l-1}\prod_{k'\in
S}\mu_{k'}} +  \sum_{k\in J_2\setminus S
}\frac{q^{\sigma(l,m-1,n)+2(r-k)}}{[2]^{l}\prod_{k'\in S}\mu_{k'}}
\\ & \quad +\sum_{k\in J_3\setminus
T}\frac{q^{\sigma(l,m,n-1)+2(r-k)}}{[2]^{l}\prod_{k'\in
S}\mu_{k'}} + \sum_{k\in
S}\frac{q^{\sigma(l,m-1,n)+2(r-k)}}{[2]^{l}\prod_{k'\in
S\setminus\{k\}}\mu_{k'}} \\ & \quad + \sum_{k\in
T}\frac{q^{\sigma(l,m,n-1)+2(r-k)}}{[2]^{l}\prod_{k'\in
S\setminus\{k+1\}}\mu_{k'}}
\frac{(1-(-q^2)^{l_{k+1}})}{(1-(-q^2)^{l_{k+1}+1})}\Bigg).
\end{align*}
On the other hand,
\begin{align*}
&\sum_{k\in
S}\frac{q^{\sigma(l,m-1,n)+2(r-k)}}{[2]^{l}\prod_{k'\in
S\setminus\{k\}}\mu_{k'}} + \sum_{k\in
T}\frac{q^{\sigma(l,m,n-1)+2(r-k)}}{[2]^{l}\prod_{k'\in
S\setminus\{k+1\}}\mu_{k'}}\frac{(1-(-q^2)^{l_{k+1}})}{(1-(-q^2)^{l_{k+1}+1})}
\\
&=\frac{1}{[2]^{l}\prod_{k'\in S}\mu_{k'}}\sum_{k\in
T}\biggl(q^{\sigma(l,m-1,n)+2(r-k-1)}
(1-(-q^2)^{l_{k+1}+1})  \\
& \hskip 5cm + q^{\sigma(l,m,n-1)+2(r-k)}
(1-(-q^2)^{l_{k+1}})\biggr)
\\
&=\frac{1}{[2]^{l}\prod_{k'\in S}\mu_{k'}}\sum_{k\in
T}\left(q^{\sigma(l,m-1,n)+2(r-k-1)}+q^{\sigma(l,m,n-1)+2(r-k)}\right).
\end{align*}
Since
\begin{equation*}
\begin{split}
\sigma(l-1,m,n+1)&=\sigma(l,m,n)-r+2\\
\sigma(l,m-1,n)&=\sigma(l,m,n-1)=\sigma(l,m,n)-r+1,
\end{split}
\end{equation*}
we have
\begin{align*}
&[r]Q_{Y,Z}(q)\\
&=Q_{Y,Z}^{\circ}(q)\frac{q^{\sigma(l,m,n)}}{[2]^{l}\prod_{k'\in
S}\mu_{k'}}\left( \sum_{k\in J_1}q^{-r+1+2(r-k+1)}+\sum_{k\in
J}q^{-r+1+2(r-k)} \right) \\
&=Q_{Y,Z}^{\circ}(q)\frac{q^{\sigma(l,m,n)}}{[2]^{l}\prod_{k'\in
S}\mu_{k'}}\left( \sum_{k=1}^rq^{-(r-1)+2(r-k)} \right)
=Q_{Y,Z}^{\circ}(q)\frac{q^{\sigma(l,m,n)}}{[2]^{l}\prod_{k'\in
S}\mu_{k'}}[r],
\end{align*}
} which completes our induction argument.

Finally, since $[2]^{l}\prod_{k\in S}\mu_{k}$ divides
$Q_{Y,Z}^{\circ}(q)$, we have $Q_{Y,Z}(q)\in\mathbb{Z}[q,q^{-1}]$.
\qed\vskip 3mm

Let $Y$ be a proper Young wall in ${\mathcal Z}(\Lambda)$ and let
$b$ be a block in $Y$. Suppose that $b$ lies in the $k$th column
of $Y$. The {\it coordinate of $b$} is defined to be $(k,l)$ where
$l$ is the number of unit cubes lying below $b$. Note that two
different blocks of type II or III have the same coordinate if
they are parts of a unit cube. Also, each coordinate corresponds
to a unit cube in a given pattern.

\begin{ex}{\rm
The block $b$ in the following figure has the coordinate
$(1,4)$.\vskip 5mm

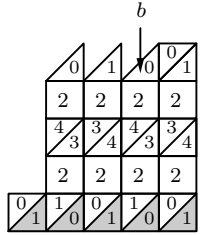
\begin{center}
\raisebox{-0.5\height}{
\begin{texdraw}
\drawdim em \setunitscale 0.13 \linewd 0.5

\move(0 0)\lvec(10 0)\lvec(10 10)\lvec(0 10)\lvec(0 0)\htext(2
6){\tiny $0$}

\move(10 0)\lvec(20 0)\lvec(20 10)\lvec(10 10)\lvec(10 0)\htext(12
6){\tiny $1$}

\move(20 0)\lvec(30 0)\lvec(30 10)\lvec(20 10)\lvec(20 0)\htext(22
6){\tiny $0$}

\move(30 0)\lvec(40 0)\lvec(40 10)\lvec(30 10)\lvec(30 0)\htext(32
6){\tiny $1$}

\move(40 0)\lvec(50 0)\lvec(50 10)\lvec(40 10)\lvec(40 0)\htext(42
6){\tiny $0$}
\move(0 0)\lvec(10 10)\lvec(10 0)\lvec(0 0)\lfill f:0.8 \htext(6
2){\tiny $1$}

\move(10 0)\lvec(20 10)\lvec(20 0)\lvec(10 0)\lfill f:0.8
\htext(16 2){\tiny $0$}

\move(20 0)\lvec(30 10)\lvec(30 0)\lvec(20 0)\lfill f:0.8
\htext(26 2){\tiny $1$}

\move(30 0)\lvec(40 10)\lvec(40 0)\lvec(30 0)\lfill f:0.8
\htext(36 2){\tiny $0$}

\move(40 0)\lvec(50 10)\lvec(50 0)\lvec(40 0)\lfill f:0.8
\htext(46 2){\tiny $1$}
\move(10 10)\lvec(20 10)\lvec(20 20)\lvec(10 20)\lvec(10
10)\htext(13 13){$_2$}

\move(20 10)\lvec(30 10)\lvec(30 20)\lvec(20 20)\lvec(20
10)\htext(23 13){$_2$}

\move(30 10)\lvec(40 10)\lvec(40 20)\lvec(30 20)\lvec(30
10)\htext(33 13){$_2$}

\move(40 10)\lvec(50 10)\lvec(50 20)\lvec(40 20)\lvec(40
10)\htext(43 13){$_2$}
\move(10 20)\lvec(20 20)\lvec(20 30)\lvec(10 30)\lvec(10
20)\htext(12 26){\tiny $4$}

\move(20 20)\lvec(30 20)\lvec(30 30)\lvec(20 30)\lvec(20
20)\htext(22 26){\tiny $3$}

\move(30 20)\lvec(40 20)\lvec(40 30)\lvec(30 30)\lvec(30
20)\htext(32 26){\tiny $4$}

\move(40 20)\lvec(50 20)\lvec(50 30)\lvec(40 30)\lvec(40
20)\htext(42 26){\tiny $3$}

\move(10 20)\lvec(20 30)\lvec(20 20)\lvec(10 20) \htext(16
22){\tiny $3$}

\move(20 20)\lvec(30 30)\lvec(30 20)\lvec(20 20) \htext(26
22){\tiny $4$}

\move(30 20)\lvec(40 30)\lvec(40 20)\lvec(30 20) \htext(36
22){\tiny $3$}

\move(40 20)\lvec(50 30)\lvec(50 20)\lvec(40 20) \htext(46
22){\tiny $4$}
\move(10 30)\lvec(20 30)\lvec(20 40)\lvec(10 40)\lvec(10
30)\htext(13 33){$_2$}

\move(20 30)\lvec(30 30)\lvec(30 40)\lvec(20 40)\lvec(20
30)\htext(23 33){$_2$}

\move(30 30)\lvec(40 30)\lvec(40 40)\lvec(30 40)\lvec(30
30)\htext(33 33){$_2$}

\move(40 30)\lvec(50 30)\lvec(50 40)\lvec(40 40)\lvec(40
30)\htext(43 33){$_2$}
\move(10 40)\lvec(20 50)\lvec(20 40)\lvec(10 40) \htext(16
42){\tiny $0$}

\move(20 40)\lvec(30 50)\lvec(30 40)\lvec(20 40) \htext(26
42){\tiny $1$}

\move(30 40)\lvec(40 50)\lvec(40 40)\lvec(30 40) \htext(36
42){\tiny $0$}

\move(40 40)\lvec(50 50)\lvec(50 40)\lvec(40 40) \htext(46
42){\tiny $1$}

\move(40 40)\lvec(50 50)\lvec(40 50)\lvec(40 40) \htext(42
46){\tiny $0$}

\move(35 54) \arrowheadtype t:F \arrowheadsize l:4 w:2 \avec(35
42)\htext(34 57){$_b$}
\end{texdraw}}
\end{center}
\vskip 5mm }
\end{ex}

Let $c=(k,l)$ ($k,l\geq 0$) be a coordinate of a block in the
pattern for ${\mathcal Z}(\Lambda)$. We define the {\it ladder at
$c$} to be the finite sequence of coordinates as follows;
\begin{equation*}
c=(k,l), (k-1, l+\Delta'), (k-2, l+2\Delta'), \cdots ,
(0,l+k\Delta'),
\end{equation*}
where $\Delta'=\Delta -1$ (resp. $\Delta$) if $\frak{g}=A_n^{(1)}$
(resp. $\frak{g}\neq A_n^{(1)}$) (This is a generalization of the
ladder for Young diagrams. See \cite{JK}).
\begin{ex}{\rm If $\frak{g}=B^{(1)}_3$ and $\Lambda=\Lambda_0$, then $\Delta=4$. There are two ladders in the following figure. The
left one is the ladder at $(7,0)$, and the right one is the ladder
at $(3,1)$.\vskip 5mm

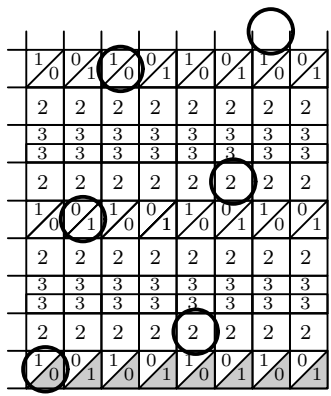
\begin{center}
\raisebox{-0.5\height}{
\begin{texdraw}
\drawdim em \setunitscale 0.13 \linewd 0.5


\move(10 0)\lvec(20 0)\lvec(20 10)\lvec(10 10)\lvec(10 0)\htext(12
6){\tiny $1$}

\move(20 0)\lvec(30 0)\lvec(30 10)\lvec(20 10)\lvec(20 0)\htext(22
6){\tiny $0$}

\move(30 0)\lvec(40 0)\lvec(40 10)\lvec(30 10)\lvec(30 0)\htext(32
6){\tiny $1$}

\move(40 0)\lvec(50 0)\lvec(50 10)\lvec(40 10)\lvec(40 0)\htext(42
6){\tiny $0$}

\move(50 0)\lvec(60 0)\lvec(60 10)\lvec(50 10)\lvec(50 0)\htext(52
6){\tiny $1$}

\move(60 0)\lvec(70 0)\lvec(70 10)\lvec(60 10)\lvec(60 0)\htext(62
6){\tiny $0$}

\move(70 0)\lvec(80 0)\lvec(80 10)\lvec(70 10)\lvec(70 0)\htext(72
6){\tiny $1$}

\move(80 0)\lvec(90 0)\lvec(90 10)\lvec(80 10)\lvec(80 0)\htext(82
6){\tiny $0$}

%

\move(10 0)\lvec(20 10)\lvec(20 0)\lvec(10 0)\lfill f:0.8
\htext(16 2){\tiny $0$}

\move(20 0)\lvec(30 10)\lvec(30 0)\lvec(20 0)\lfill f:0.8
\htext(26 2){\tiny $1$}

\move(30 0)\lvec(40 10)\lvec(40 0)\lvec(30 0)\lfill f:0.8
\htext(36 2){\tiny $0$}

\move(40 0)\lvec(50 10)\lvec(50 0)\lvec(40 0)\lfill f:0.8
\htext(46 2){\tiny $1$}

\move(50 0)\lvec(60 10)\lvec(60 0)\lvec(50 0)\lfill f:0.8
\htext(56 2){\tiny $0$}

\move(60 0)\lvec(70 10)\lvec(70 0)\lvec(60 0)\lfill f:0.8
\htext(66 2){\tiny $1$}

\move(70 0)\lvec(80 10)\lvec(80 0)\lvec(70 0)\lfill f:0.8
\htext(76 2){\tiny $0$}

\move(80 0)\lvec(90 10)\lvec(90 0)\lvec(80 0)\lfill f:0.8
\htext(86 2){\tiny $1$}

\move(10 10)\lvec(20 10)\lvec(20 20)\lvec(10 20)\lvec(10
10)\htext(13 13){$_2$}

\move(20 10)\lvec(30 10)\lvec(30 20)\lvec(20 20)\lvec(20
10)\htext(23 13){$_2$}

\move(30 10)\lvec(40 10)\lvec(40 20)\lvec(30 20)\lvec(30
10)\htext(33 13){$_2$}

\move(40 10)\lvec(50 10)\lvec(50 20)\lvec(40 20)\lvec(40
10)\htext(43 13){$_2$}

\move(50 10)\lvec(60 10)\lvec(60 20)\lvec(50 20)\lvec(50
10)\htext(53 13){$_2$}

\move(60 10)\lvec(70 10)\lvec(70 20)\lvec(60 20)\lvec(60
10)\htext(63 13){$_2$}

\move(70 10)\lvec(80 10)\lvec(80 20)\lvec(70 20)\lvec(70
10)\htext(73 13){$_2$}

\move(80 10)\lvec(90 10)\lvec(90 20)\lvec(80 20)\lvec(80
10)\htext(83 13){$_2$}
\move(10 20)\lvec(20 20)\lvec(20 30)\lvec(10 30)\lvec(10
20)\htext(13 26){\tiny $3$}

\move(20 20)\lvec(30 20)\lvec(30 30)\lvec(20 30)\lvec(20
20)\htext(23 26){\tiny $3$}

\move(30 20)\lvec(40 20)\lvec(40 30)\lvec(30 30)\lvec(30
20)\htext(33 26){\tiny $3$}

\move(40 20)\lvec(50 20)\lvec(50 30)\lvec(40 30)\lvec(40
20)\htext(43 26){\tiny $3$}

\move(50 20)\lvec(60 20)\lvec(60 30)\lvec(50 30)\lvec(50
20)\htext(53 26){\tiny $3$}

\move(60 20)\lvec(70 20)\lvec(70 30)\lvec(60 30)\lvec(60
20)\htext(63 26){\tiny $3$}

\move(70 20)\lvec(80 20)\lvec(80 30)\lvec(70 30)\lvec(70
20)\htext(73 26){\tiny $3$}

\move(80 20)\lvec(90 20)\lvec(90 30)\lvec(80 30)\lvec(80
20)\htext(83 26){\tiny $3$}

\move(10 20)\lvec(20 20)\lvec(20 25)\lvec(10 25)\lvec(10 20)
\htext(13 21){\tiny $3$}

\move(20 20)\lvec(30 20)\lvec(30 25)\lvec(20 25)\lvec(20 20)
\htext(23 21){\tiny $3$}

\move(30 20)\lvec(40 20)\lvec(40 25)\lvec(30 25)\lvec(30 20)
\htext(33 21){\tiny $3$}

\move(40 20)\lvec(50 20)\lvec(50 25)\lvec(40 25)\lvec(40 20)
\htext(43 21){\tiny $3$}

\move(50 20)\lvec(60 20)\lvec(60 25)\lvec(50 25)\lvec(50
20)\htext(53 21){\tiny $3$}

\move(60 20)\lvec(70 20)\lvec(70 25)\lvec(60 25)\lvec(60
20)\htext(63 21){\tiny $3$}

\move(70 20)\lvec(80 20)\lvec(80 25)\lvec(70 25)\lvec(70
20)\htext(73 21){\tiny $3$}

\move(80 20)\lvec(90 20)\lvec(90 25)\lvec(80 25)\lvec(80
20)\htext(83 21){\tiny $3$}
\move(10 30)\lvec(20 30)\lvec(20 40)\lvec(10 40)\lvec(10
30)\htext(13 33){$_2$}

\move(20 30)\lvec(30 30)\lvec(30 40)\lvec(20 40)\lvec(20
30)\htext(23 33){$_2$}

\move(30 30)\lvec(40 30)\lvec(40 40)\lvec(30 40)\lvec(30
30)\htext(33 33){$_2$}

\move(40 30)\lvec(50 30)\lvec(50 40)\lvec(40 40)\lvec(40
30)\htext(43 33){$_2$}

\move(50 30)\lvec(60 30)\lvec(60 40)\lvec(50 40)\lvec(50
30)\htext(53 33){$_2$}

\move(60 30)\lvec(70 30)\lvec(70 40)\lvec(60 40)\lvec(60
30)\htext(63 33){$_2$}

\move(70 30)\lvec(80 30)\lvec(80 40)\lvec(70 40)\lvec(70
30)\htext(73 33){$_2$}

\move(80 30)\lvec(90 30)\lvec(90 40)\lvec(80 40)\lvec(80
30)\htext(83 33){$_2$}

\move(10 40)\lvec(20 50)\lvec(10 50)\lvec(10 40) \htext(12
46){\tiny $1$}

\move(20 40)\lvec(30 50)\lvec(20 50)\lvec(20 40) \htext(22
46){\tiny $0$}

\move(30 40)\lvec(40 50)\lvec(30 50)\lvec(30 40) \htext(32
46){\tiny $1$}

\move(40 40)\lvec(50 50)\lvec(40 50)\lvec(40 40) \htext(42
46){\tiny $0$}

\move(40 40)\lvec(50 50)\lvec(50 40)\lvec(40 40) \htext(46
42){\tiny $1$}

\move(50 40)\lvec(60 40)\lvec(60 50)\lvec(50 50) \lvec(50
40)\htext(52 46){\tiny $1$}

\move(60 40)\lvec(70 40)\lvec(70 50)\lvec(60 50)\lvec(60
40)\htext(62 46){\tiny $0$}

\move(70 40)\lvec(80 40)\lvec(80 50)\lvec(70 50)\lvec(70
40)\htext(72 46){\tiny $1$}

\move(80 40)\lvec(90 40)\lvec(90 50)\lvec(80 50)\lvec(80
40)\htext(82 46){\tiny $0$}
\move(10 40)\lvec(20 50)\lvec(20 40)\lvec(10 40) \htext(16
42){\tiny $0$}

\move(20 40)\lvec(30 50)\lvec(30 40)\lvec(20 40) \htext(26
42){\tiny $1$}

\move(30 40)\lvec(40 50)\lvec(40 40)\lvec(30 40) \htext(36
42){\tiny $0$}

\move(40 40)\lvec(50 50)\lvec(50 40)\lvec(40 40) \htext(46
42){\tiny $1$}

\move(50 40)\lvec(60 50)\lvec(60 40)\lvec(50 40) \htext(56
42){\tiny $0$}

\move(60 40)\lvec(70 50)\lvec(70 40)\lvec(60 40) \htext(66
42){\tiny $1$}

\move(70 40)\lvec(80 50)\lvec(80 40)\lvec(70 40) \htext(76
42){\tiny $0$}

\move(80 40)\lvec(90 50)\lvec(90 40)\lvec(80 40) \htext(86
42){\tiny $1$}
\move(10 50)\lvec(20 50)\lvec(20 60)\lvec(10 60)\lvec(10
50)\htext(13 53){$_2$}

\move(20 50)\lvec(30 50)\lvec(30 60)\lvec(20 60)\lvec(20
50)\htext(23 53){$_2$}

\move(30 50)\lvec(40 50)\lvec(40 60)\lvec(30 60)\lvec(30
50)\htext(33 53){$_2$}

\move(40 50)\lvec(50 50)\lvec(50 60)\lvec(40 60)\lvec(40
50)\htext(43 53){$_2$}

\move(50 50)\lvec(60 50)\lvec(60 60)\lvec(50 60)\lvec(50
50)\htext(53 53){$_2$}

\move(60 50)\lvec(70 50)\lvec(70 60)\lvec(60 60)\lvec(60
50)\htext(63 53){$_2$}

\move(70 50)\lvec(80 50)\lvec(80 60)\lvec(70 60)\lvec(70
50)\htext(73 53){$_2$}

\move(80 50)\lvec(90 50)\lvec(90 60)\lvec(80 60)\lvec(80
50)\htext(83 53){$_2$}
\move(10 60)\lvec(20 60)\lvec(20 65)\lvec(10 65)\lvec(10 60)
\htext(13 61){\tiny $3$}

\move(20 60)\lvec(30 60)\lvec(30 65)\lvec(20 65)\lvec(20 60)
\htext(23 61){\tiny $3$}

\move(30 60)\lvec(40 60)\lvec(40 65)\lvec(30 65)\lvec(30 60)
\htext(33 61){\tiny $3$}

\move(40 60)\lvec(50 60)\lvec(50 65)\lvec(40 65)\lvec(40 60)
\htext(43 61){\tiny $3$}

\move(50 60)\lvec(60 60)\lvec(60 65)\lvec(50 65)\lvec(50
60)\htext(53 61){\tiny $3$}

\move(60 60)\lvec(70 60)\lvec(70 65)\lvec(60 65)\lvec(60
60)\htext(63 61){\tiny $3$}

\move(70 60)\lvec(80 60)\lvec(80 65)\lvec(70 65)\lvec(70
60)\htext(73 61){\tiny $3$}

\move(80 60)\lvec(90 60)\lvec(90 65)\lvec(80 65)\lvec(80
60)\htext(83 61){\tiny $3$}
\move(10 60)\lvec(20 60)\lvec(20 70)\lvec(10 70)\lvec(10
60)\htext(13 66){\tiny $3$}

\move(20 60)\lvec(30 60)\lvec(30 70)\lvec(20 70)\lvec(20
60)\htext(23 66){\tiny $3$}

\move(30 60)\lvec(40 60)\lvec(40 70)\lvec(30 70)\lvec(30
60)\htext(33 66){\tiny $3$}

\move(40 60)\lvec(50 60)\lvec(50 70)\lvec(40 70)\lvec(40
60)\htext(43 66){\tiny $3$}

\move(50 60)\lvec(60 60)\lvec(60 70)\lvec(50 70)\lvec(50
60)\htext(53 66){\tiny $3$}

\move(60 60)\lvec(70 60)\lvec(70 70)\lvec(60 70)\lvec(60
60)\htext(63 66){\tiny $3$}

\move(70 60)\lvec(80 60)\lvec(80 70)\lvec(70 70)\lvec(70
60)\htext(73 66){\tiny $3$}

\move(80 60)\lvec(90 60)\lvec(90 70)\lvec(80 70)\lvec(80
60)\htext(83 66){\tiny $3$}
\move(10 70)\lvec(20 70)\lvec(20 80)\lvec(10 80)\lvec(10
70)\htext(13 73){$_2$}

\move(20 70)\lvec(30 70)\lvec(30 80)\lvec(20 80)\lvec(20
70)\htext(23 73){$_2$}

\move(30 70)\lvec(40 70)\lvec(40 80)\lvec(30 80)\lvec(30
70)\htext(33 73){$_2$}

\move(40 70)\lvec(50 70)\lvec(50 80)\lvec(40 80)\lvec(40
70)\htext(43 73){$_2$}

\move(50 70)\lvec(60 70)\lvec(60 80)\lvec(50 80)\lvec(50
70)\htext(53 73){$_2$}

\move(60 70)\lvec(70 70)\lvec(70 80)\lvec(60 80)\lvec(60
70)\htext(63 73){$_2$}

\move(70 70)\lvec(80 70)\lvec(80 80)\lvec(70 80)\lvec(70
70)\htext(73 73){$_2$}

\move(80 70)\lvec(90 70)\lvec(90 80)\lvec(80 80)\lvec(80
70)\htext(83 73){$_2$}
\move(10 80)\lvec(20 90)\lvec(20 80)\lvec(10 80) \htext(16
82){\tiny $0$}

\move(20 80)\lvec(30 90)\lvec(30 80)\lvec(20 80) \htext(26
82){\tiny $1$}

\move(30 80)\lvec(40 90)\lvec(40 80)\lvec(30 80) \htext(36
82){\tiny $0$}

\move(40 80)\lvec(50 90)\lvec(50 80)\lvec(40 80) \htext(46
82){\tiny $1$}

\move(50 80)\lvec(60 90)\lvec(60 80)\lvec(50 80) \htext(56
82){\tiny $0$}

\move(60 80)\lvec(70 90)\lvec(70 80)\lvec(60 80)\htext(66
82){\tiny $1$}

\move(70 80)\lvec(80 90)\lvec(80 80)\lvec(70 80)\htext(76
82){\tiny $0$}

\move(80 80)\lvec(90 90)\lvec(90 80)\lvec(80 80)\htext(86
82){\tiny $1$}
\move(10 80)\lvec(20 80)\lvec(20 90)\lvec(10 90)\lvec(10
80)\htext(12 86){\tiny $1$}

\move(20 80)\lvec(30 80)\lvec(30 90)\lvec(20 90)\lvec(20
80)\htext(22 86){\tiny $0$}

\move(30 80)\lvec(40 80)\lvec(40 90)\lvec(30 90)\lvec(30
80)\htext(32 86){\tiny $1$}

\move(40 80)\lvec(50 80)\lvec(50 90)\lvec(40 90)\lvec(40
80)\htext(42 86){\tiny $0$}

\move(50 80)\lvec(60 80)\lvec(60 90)\lvec(50 90)\lvec(50
80)\htext(52 86){\tiny $1$}

\move(60 80)\lvec(70 80)\lvec(70 90)\lvec(60 90)\lvec(60
80)\htext(62 86){\tiny $0$}

\move(70 80)\lvec(80 80)\lvec(80 90)\lvec(70 90)\lvec(70
80)\htext(72 86){\tiny $1$}

\move(80 80)\lvec(90 80)\lvec(90 90)\lvec(80 90)\lvec(80
80)\htext(82 86){\tiny $0$}

\move(10 0)\lvec(5 0)\move(10 10)\lvec(5 10)\move(10 20)\lvec(5
20)\move(10 30)\lvec(5 30)\move(10 40)\lvec(5 40)\move(10
50)\lvec(5 50)\move(10 60)\lvec(5 60)\move(10 70)\lvec(5
70)\move(10 80)\lvec(5 80)\move(10 90)\lvec(5 90)

\move(10 95)\lvec(10 90)\move(20 95)\lvec(20 90)\move(30
95)\lvec(30 90)\move(40 95)\lvec(40 90)\move(50 95)\lvec(50
90)\move(60 95)\lvec(60 90)\move(70 95)\lvec(70 90)\move(80
95)\lvec(80 90)\move(90 95)\lvec(90 90)

\move(15 5)\linewd 1 \lcir r:6 \move(25 45)\linewd 1 \lcir r:6
\move(35 85)\linewd 1 \lcir r:6

\move(55 15)\linewd 1 \lcir r:6 \move(65 55)\linewd 1 \lcir r:6
\move(75 95)\linewd 1 \lcir r:6

\end{texdraw}}
\end{center}\vskip 5mm
}\end{ex}

Let $Y$ be a reduced proper Young wall in ${\mathcal Y}(\Lambda)$
and let $y_l$ be the left-most column in $Y$ with $|y_l|\neq 0$.
Take a block $b$ lying at the top of $y_l$ and let $i$ be its
color. (If the $b$ is of type III and there is another block of
type III on top of $y_l$, we take the block at the front.) Suppose
that the coordinate of $b$ is $c=(k,l)$ and let $L_c$ be the
ladder at $c$. We denote by $Y\cap L_c$ the $i$-blocks in $Y$
whose coordinates are in $L_c$. (In fact, $L_c$ is the left-most
ladder that has a nontrivial intersection with $Y$.) We define
$\overline{Y}$ to be the proper Young wall that is obtained from
$Y$ by removing all the $i$-blocks in $Y\cap L_c$. Then it is easy
to see that $\overline{Y}$ is also reduced. That is, we remove all
the $i$-blocks along the left-most ladder to obtain another
reduced proper Young wall. This process will play a crucial role
in constructing the global basis $\mathcal{G}(\Lambda)$ of
$V(\Lambda)$.
\begin{ex}{\rm \mbox{}

(a) If $\frak{g}=A_5^{(2)}$ and\vskip 5mm

\begin{center}
$Y=$\raisebox{-0.5\height}{
\begin{texdraw}
\drawdim em \setunitscale 0.13 \linewd 0.5

\move(10 0)\lvec(20 0)\lvec(20 10)\lvec(10 10)\lvec(10 0)

\move(20 0)\lvec(30 0)\lvec(30 10)\lvec(20 10)\lvec(20 0)

\move(30 0)\lvec(40 0)\lvec(40 10)\lvec(30 10)\lvec(30 0)

\move(40 0)\lvec(50 0)\lvec(50 10)\lvec(40 10)\lvec(40 0)
\move(10 10)\lvec(20 10)\lvec(20 20)\lvec(10 20)\lvec(10
10)\htext(13 13){$_2$}

\move(20 10)\lvec(30 10)\lvec(30 20)\lvec(20 20)\lvec(20
10)\htext(23 13){$_2$}

\move(30 10)\lvec(40 10)\lvec(40 20)\lvec(30 20)\lvec(30
10)\htext(33 13){$_2$}

\move(40 10)\lvec(50 10)\lvec(50 20)\lvec(40 20)\lvec(40
10)\htext(43 13){$_2$}
\move(20 20)\lvec(30 20)\lvec(30 30)\lvec(20 30)\lvec(20
20)\htext(23 23){$_3$}

\move(30 20)\lvec(40 20)\lvec(40 30)\lvec(30 30)\lvec(30
20)\htext(33 23){$_3$}

\move(40 20)\lvec(50 20)\lvec(50 30)\lvec(40 30)\lvec(40
20)\htext(43 23){$_3$}
\move(20 30)\lvec(30 30)\lvec(30 40)\lvec(20 40)\lvec(20
30)\htext(23 33){$_2$}

\move(30 30)\lvec(40 30)\lvec(40 40)\lvec(30 40)\lvec(30
30)\htext(33 33){$_2$}

\move(40 30)\lvec(50 30)\lvec(50 40)\lvec(40 40)\lvec(40
30)\htext(43 33){$_2$}
\move(20 40)\lvec(30 40)\lvec(30 50)\lvec(20 50)\lvec(20 40)

\move(30 40)\lvec(40 40)\lvec(40 50)\lvec(30 50)\lvec(30 40)

\move(40 40)\lvec(50 40)\lvec(50 50)\lvec(40 50)\lvec(40 40)
\move(20 50)\lvec(30 50)\lvec(30 60)\lvec(20 60)\lvec(20
50)\htext(23 53){$_2$}

\move(30 50)\lvec(40 50)\lvec(40 60)\lvec(30 60)\lvec(30
50)\htext(33 53){$_2$}

\move(40 50)\lvec(50 50)\lvec(50 60)\lvec(40 60)\lvec(40
50)\htext(43 53){$_2$}

\move(40 70)\lvec(50 70)\lvec(50 80)\lvec(40 80)\lvec(40
70)\htext(43 73){$_2$}
\move(30 60)\lvec(40 60)\lvec(40 70)\lvec(30 70)\lvec(30
60)\htext(33 63){$_3$}

\move(40 60)\lvec(50 60)\lvec(50 70)\lvec(40 70)\lvec(40
60)\htext(43 63){$_3$}
\move(10 0)\lvec(20 10)\lvec(20 0)\lvec(10 0)\lfill f:0.8
\htext(11 5){\tiny $1$}\htext(16 1){\tiny $0$}

\move(20 0)\lvec(30 10)\lvec(30 0)\lvec(20 0)\lfill f:0.8
\htext(22 5){\tiny $0$}\htext(26 2){\tiny $1$}

\move(30 0)\lvec(40 10)\lvec(40 0)\lvec(30 0)\lfill f:0.8
\htext(32 5){\tiny $1$}\htext(36 2){\tiny $0$}

\move(40 0)\lvec(50 10)\lvec(50 0)\lvec(40 0)\lfill f:0.8
\htext(42 5){\tiny $0$}\htext(46 2){\tiny $1$}
\move(20 40)\lvec(30 50)\lvec(30 40)\lvec(20 40)\htext(22
45){\tiny $0$}\htext(26 42){\tiny $1$}

\move(30 40)\lvec(40 50)\lvec(40 40)\lvec(30 40)\htext(32
45){\tiny $1$}\htext(36 42){\tiny $0$}

\move(40 40)\lvec(50 50)\lvec(50 40)\lvec(40 40)\htext(42
45){\tiny $0$}\htext(46 42){\tiny $1$}

\move(15 15)\linewd 1 \lcir r:6 \move(25 55)\linewd 1 \lcir r:6
\end{texdraw}}\quad ,\quad then
$\overline{Y}=$\raisebox{-0.5\height}{
\begin{texdraw}
\drawdim em \setunitscale 0.13 \linewd 0.5

\move(10 0)\lvec(20 0)\lvec(20 10)\lvec(10 10)\lvec(10 0)

\move(20 0)\lvec(30 0)\lvec(30 10)\lvec(20 10)\lvec(20 0)

\move(30 0)\lvec(40 0)\lvec(40 10)\lvec(30 10)\lvec(30 0)

\move(40 0)\lvec(50 0)\lvec(50 10)\lvec(40 10)\lvec(40 0)
%

\move(20 10)\lvec(30 10)\lvec(30 20)\lvec(20 20)\lvec(20
10)\htext(23 13){$_2$}

\move(30 10)\lvec(40 10)\lvec(40 20)\lvec(30 20)\lvec(30
10)\htext(33 13){$_2$}

\move(40 10)\lvec(50 10)\lvec(50 20)\lvec(40 20)\lvec(40
10)\htext(43 13){$_2$}
\move(20 20)\lvec(30 20)\lvec(30 30)\lvec(20 30)\lvec(20
20)\htext(23 23){$_3$}

\move(30 20)\lvec(40 20)\lvec(40 30)\lvec(30 30)\lvec(30
20)\htext(33 23){$_3$}

\move(40 20)\lvec(50 20)\lvec(50 30)\lvec(40 30)\lvec(40
20)\htext(43 23){$_3$}
\move(20 30)\lvec(30 30)\lvec(30 40)\lvec(20 40)\lvec(20
30)\htext(23 33){$_2$}

\move(30 30)\lvec(40 30)\lvec(40 40)\lvec(30 40)\lvec(30
30)\htext(33 33){$_2$}

\move(40 30)\lvec(50 30)\lvec(50 40)\lvec(40 40)\lvec(40
30)\htext(43 33){$_2$}
\move(20 40)\lvec(30 40)\lvec(30 50)\lvec(20 50)\lvec(20 40)

\move(30 40)\lvec(40 40)\lvec(40 50)\lvec(30 50)\lvec(30 40)

\move(40 40)\lvec(50 40)\lvec(50 50)\lvec(40 50)\lvec(40 40)
%

\move(30 50)\lvec(40 50)\lvec(40 60)\lvec(30 60)\lvec(30
50)\htext(33 53){$_2$}

\move(40 50)\lvec(50 50)\lvec(50 60)\lvec(40 60)\lvec(40
50)\htext(43 53){$_2$}

\move(40 70)\lvec(50 70)\lvec(50 80)\lvec(40 80)\lvec(40
70)\htext(43 73){$_2$}
\move(30 60)\lvec(40 60)\lvec(40 70)\lvec(30 70)\lvec(30
60)\htext(33 63){$_3$}

\move(40 60)\lvec(50 60)\lvec(50 70)\lvec(40 70)\lvec(40
60)\htext(43 63){$_3$}
\move(10 0)\lvec(20 10)\lvec(20 0)\lvec(10 0)\lfill f:0.8
\htext(11 5){\tiny $1$}\htext(16 1){\tiny $0$}

\move(20 0)\lvec(30 10)\lvec(30 0)\lvec(20 0)\lfill f:0.8
\htext(22 5){\tiny $0$}\htext(26 2){\tiny $1$}

\move(30 0)\lvec(40 10)\lvec(40 0)\lvec(30 0)\lfill f:0.8
\htext(32 5){\tiny $1$}\htext(36 2){\tiny $0$}

\move(40 0)\lvec(50 10)\lvec(50 0)\lvec(40 0)\lfill f:0.8
\htext(42 5){\tiny $0$}\htext(46 2){\tiny $1$}
\move(20 40)\lvec(30 50)\lvec(30 40)\lvec(20 40)\htext(22
45){\tiny $0$}\htext(26 42){\tiny $1$}

\move(30 40)\lvec(40 50)\lvec(40 40)\lvec(30 40)\htext(32
45){\tiny $1$}\htext(36 42){\tiny $0$}

\move(40 40)\lvec(50 50)\lvec(50 40)\lvec(40 40)\htext(42
45){\tiny $0$}\htext(46 42){\tiny $1$}

\end{texdraw}}\quad.
\end{center}
\vskip 5mm

(b) If $\frak{g}=A_4^{(2)}$ and\vskip 5mm

\begin{center}
$Y=$\raisebox{-0.5\height}{
\begin{texdraw}
\drawdim em \setunitscale 0.13 \linewd 0.5

\move(20 0)\lvec(30 0)\lvec(30 10)\lvec(20 10)\lvec(20 0)\htext(23
1){\tiny $0$}

\move(30 0)\lvec(40 0)\lvec(40 10)\lvec(30 10)\lvec(30 0)\htext(33
1){\tiny $0$}

\move(40 0)\lvec(50 0)\lvec(50 10)\lvec(40 10)\lvec(40 0)\htext(43
1){\tiny $0$}

\move(50 0)\lvec(60 0)\lvec(60 10)\lvec(50 10)\lvec(50 0)\htext(53
1){\tiny $0$}

\move(60 0)\lvec(70 0)\lvec(70 10)\lvec(60 10)\lvec(60 0)\htext(63
1){\tiny $0$}
\move(20 0)\lvec(30 0)\lvec(30 5)\lvec(20 5)\lvec(20 0)\lfill
f:0.8 \htext(23 6){\tiny $0$}

\move(30 0)\lvec(40 0)\lvec(40 5)\lvec(30 5)\lvec(30 0)\lfill
f:0.8 \htext(33 6){\tiny $0$}

\move(40 0)\lvec(50 0)\lvec(50 5)\lvec(40 5)\lvec(40 0)\lfill
f:0.8 \htext(43 6){\tiny $0$}

\move(50 0)\lvec(60 0)\lvec(60 5)\lvec(50 5)\lvec(50 0)\lfill
f:0.8 \htext(53 6){\tiny $0$}

\move(60 0)\lvec(70 0)\lvec(70 5)\lvec(60 5)\lvec(60 0)\lfill
f:0.8 \htext(63 6){\tiny $0$}
\move(30 10)\lvec(40 10)\lvec(40 20)\lvec(30 20)\lvec(30
10)\htext(33 13){$_1$}

\move(40 10)\lvec(50 10)\lvec(50 20)\lvec(40 20)\lvec(40
10)\htext(43 13){$_1$}

\move(50 10)\lvec(60 10)\lvec(60 20)\lvec(50 20)\lvec(50
10)\htext(53 13){$_1$}

\move(60 10)\lvec(70 10)\lvec(70 20)\lvec(60 20)\lvec(60
10)\htext(63 13){$_1$}
\move(30 20)\lvec(40 20)\lvec(40 30)\lvec(30 30)\lvec(30
20)\htext(33 23){$_2$}

\move(40 20)\lvec(50 20)\lvec(50 30)\lvec(40 30)\lvec(40
20)\htext(43 23){$_2$}

\move(50 20)\lvec(60 20)\lvec(60 30)\lvec(50 30)\lvec(50
20)\htext(53 23){$_2$}

\move(60 20)\lvec(70 20)\lvec(70 30)\lvec(60 30)\lvec(60
20)\htext(63 23){$_2$}
\move(30 30)\lvec(40 30)\lvec(40 40)\lvec(30 40)\lvec(30
30)\htext(33 33){$_1$}

\move(40 30)\lvec(50 30)\lvec(50 40)\lvec(40 40)\lvec(40
30)\htext(43 33){$_1$}

\move(50 30)\lvec(60 30)\lvec(60 40)\lvec(50 40)\lvec(50
30)\htext(53 33){$_1$}

\move(60 30)\lvec(70 30)\lvec(70 40)\lvec(60 40)\lvec(60
30)\htext(63 33){$_1$}
\move(30 40)\lvec(40 40)\lvec(40 45)\lvec(30 45)\lvec(30
40)\htext(33 41){\tiny $0$}

\move(40 40)\lvec(50 40)\lvec(50 45)\lvec(40 45)\lvec(40
40)\htext(43 41){\tiny $0$}

\move(50 40)\lvec(60 40)\lvec(60 45)\lvec(50 45)\lvec(50
40)\htext(53 41){\tiny $0$}

\move(60 40)\lvec(70 40)\lvec(70 45)\lvec(60 45)\lvec(60
40)\htext(63 41){\tiny $0$}
\move(30 40)\lvec(40 40)\lvec(40 50)\lvec(30 50)\lvec(30
40)\htext(33 46){\tiny $0$}

\move(40 40)\lvec(50 40)\lvec(50 50)\lvec(40 50)\lvec(40
40)\htext(43 46){\tiny $0$}

\move(50 40)\lvec(60 40)\lvec(60 50)\lvec(50 50)\lvec(50
40)\htext(53 46){\tiny $0$}

\move(60 40)\lvec(70 40)\lvec(70 50)\lvec(60 50)\lvec(60
40)\htext(63 46){\tiny $0$}
\move(40 50)\lvec(50 50)\lvec(50 60)\lvec(40 60)\lvec(40
50)\htext(43 53){$_1$}

\move(50 50)\lvec(60 50)\lvec(60 60)\lvec(50 60)\lvec(50
50)\htext(53 53){$_1$}

\move(60 50)\lvec(70 50)\lvec(70 60)\lvec(60 60)\lvec(60
50)\htext(63 53){$_1$}
\move(40 60)\lvec(50 60)\lvec(50 70)\lvec(40 70)\lvec(40
60)\htext(43 63){$_2$}

\move(50 60)\lvec(60 60)\lvec(60 70)\lvec(50 70)\lvec(50
60)\htext(53 63){$_2$}

\move(60 60)\lvec(70 60)\lvec(70 70)\lvec(60 70)\lvec(60
60)\htext(63 63){$_2$}
\move(40 70)\lvec(50 70)\lvec(50 80)\lvec(40 80)\lvec(40
70)\htext(43 73){$_1$}

\move(50 70)\lvec(60 70)\lvec(60 80)\lvec(50 80)\lvec(50
70)\htext(53 73){$_1$}

\move(60 70)\lvec(70 70)\lvec(70 80)\lvec(60 80)\lvec(60
70)\htext(63 73){$_1$}
\move(40 80)\lvec(50 80)\lvec(50 85)\lvec(40 85)\lvec(40
80)\htext(43 81){\tiny $0$}

\move(50 80)\lvec(60 80)\lvec(60 85)\lvec(50 85)\lvec(50
80)\htext(53 81){\tiny $0$}

\move(60 80)\lvec(70 80)\lvec(70 85)\lvec(60 85)\lvec(60
80)\htext(63 81){\tiny $0$}

\move(25 5)\linewd 1 \lcir r:6 \move(35 45)\linewd 1 \lcir r:6
\move(45 85)\linewd 1 \lcir r:6
\end{texdraw}}\quad , \quad then
$\overline{Y}=$\raisebox{-0.5\height}{
\begin{texdraw}
\drawdim em \setunitscale 0.13 \linewd 0.5

\move(30 0)\lvec(40 0)\lvec(40 10)\lvec(30 10)\lvec(30 0)\htext(33
1){\tiny $0$}

\move(40 0)\lvec(50 0)\lvec(50 10)\lvec(40 10)\lvec(40 0)\htext(43
1){\tiny $0$}

\move(50 0)\lvec(60 0)\lvec(60 10)\lvec(50 10)\lvec(50 0)\htext(53
1){\tiny $0$}

\move(60 0)\lvec(70 0)\lvec(70 10)\lvec(60 10)\lvec(60 0)\htext(63
1){\tiny $0$}
\move(30 0)\lvec(40 0)\lvec(40 5)\lvec(30 5)\lvec(30 0)\lfill
f:0.8 \htext(33 6){\tiny $0$}

\move(40 0)\lvec(50 0)\lvec(50 5)\lvec(40 5)\lvec(40 0)\lfill
f:0.8 \htext(43 6){\tiny $0$}

\move(50 0)\lvec(60 0)\lvec(60 5)\lvec(50 5)\lvec(50 0)\lfill
f:0.8 \htext(53 6){\tiny $0$}

\move(60 0)\lvec(70 0)\lvec(70 5)\lvec(60 5)\lvec(60 0)\lfill
f:0.8 \htext(63 6){\tiny $0$}

\move(30 10)\lvec(40 10)\lvec(40 20)\lvec(30 20)\lvec(30
10)\htext(33 13){$_1$}

\move(40 10)\lvec(50 10)\lvec(50 20)\lvec(40 20)\lvec(40
10)\htext(43 13){$_1$}

\move(50 10)\lvec(60 10)\lvec(60 20)\lvec(50 20)\lvec(50
10)\htext(53 13){$_1$}

\move(60 10)\lvec(70 10)\lvec(70 20)\lvec(60 20)\lvec(60
10)\htext(63 13){$_1$}
\move(30 20)\lvec(40 20)\lvec(40 30)\lvec(30 30)\lvec(30
20)\htext(33 23){$_2$}

\move(40 20)\lvec(50 20)\lvec(50 30)\lvec(40 30)\lvec(40
20)\htext(43 23){$_2$}

\move(50 20)\lvec(60 20)\lvec(60 30)\lvec(50 30)\lvec(50
20)\htext(53 23){$_2$}

\move(60 20)\lvec(70 20)\lvec(70 30)\lvec(60 30)\lvec(60
20)\htext(63 23){$_2$}
\move(30 30)\lvec(40 30)\lvec(40 40)\lvec(30 40)\lvec(30
30)\htext(33 33){$_1$}

\move(40 30)\lvec(50 30)\lvec(50 40)\lvec(40 40)\lvec(40
30)\htext(43 33){$_1$}

\move(50 30)\lvec(60 30)\lvec(60 40)\lvec(50 40)\lvec(50
30)\htext(53 33){$_1$}

\move(60 30)\lvec(70 30)\lvec(70 40)\lvec(60 40)\lvec(60
30)\htext(63 33){$_1$}
%

\move(40 40)\lvec(50 40)\lvec(50 45)\lvec(40 45)\lvec(40
40)\htext(43 41){\tiny $0$}

\move(50 40)\lvec(60 40)\lvec(60 45)\lvec(50 45)\lvec(50
40)\htext(53 41){\tiny $0$}

\move(60 40)\lvec(70 40)\lvec(70 45)\lvec(60 45)\lvec(60
40)\htext(63 41){\tiny $0$}
%

\move(40 40)\lvec(50 40)\lvec(50 50)\lvec(40 50)\lvec(40
40)\htext(43 46){\tiny $0$}

\move(50 40)\lvec(60 40)\lvec(60 50)\lvec(50 50)\lvec(50
40)\htext(53 46){\tiny $0$}

\move(60 40)\lvec(70 40)\lvec(70 50)\lvec(60 50)\lvec(60
40)\htext(63 46){\tiny $0$}
\move(40 50)\lvec(50 50)\lvec(50 60)\lvec(40 60)\lvec(40
50)\htext(43 53){$_1$}

\move(50 50)\lvec(60 50)\lvec(60 60)\lvec(50 60)\lvec(50
50)\htext(53 53){$_1$}

\move(60 50)\lvec(70 50)\lvec(70 60)\lvec(60 60)\lvec(60
50)\htext(63 53){$_1$}
\move(40 60)\lvec(50 60)\lvec(50 70)\lvec(40 70)\lvec(40
60)\htext(43 63){$_2$}

\move(50 60)\lvec(60 60)\lvec(60 70)\lvec(50 70)\lvec(50
60)\htext(53 63){$_2$}

\move(60 60)\lvec(70 60)\lvec(70 70)\lvec(60 70)\lvec(60
60)\htext(63 63){$_2$}
\move(40 70)\lvec(50 70)\lvec(50 80)\lvec(40 80)\lvec(40
70)\htext(43 73){$_1$}

\move(50 70)\lvec(60 70)\lvec(60 80)\lvec(50 80)\lvec(50
70)\htext(53 73){$_1$}

\move(60 70)\lvec(70 70)\lvec(70 80)\lvec(60 80)\lvec(60
70)\htext(63 73){$_1$}
%

\move(50 80)\lvec(60 80)\lvec(60 85)\lvec(50 85)\lvec(50
80)\htext(53 81){\tiny $0$}

\move(60 80)\lvec(70 80)\lvec(70 85)\lvec(60 85)\lvec(60
80)\htext(63 81){\tiny $0$}
\end{texdraw}}\quad .
\end{center}\vskip 5mm

(c) If $\frak{g}=B_3^{(1)}$ and\vskip 5mm

\begin{center}
$Y=$\raisebox{-0.5\height}{
\begin{texdraw}
\drawdim em \setunitscale 0.13 \linewd 0.5

\move(50 0)\lvec(60 0)\lvec(60 10)\lvec(50 10)\lvec(50 0)\htext(52
6){\tiny $1$}

\move(60 0)\lvec(70 0)\lvec(70 10)\lvec(60 10)\lvec(60 0)\htext(62
6){\tiny $0$}

\move(70 0)\lvec(80 0)\lvec(80 10)\lvec(70 10)\lvec(70 0)\htext(72
6){\tiny $1$}

\move(80 0)\lvec(90 0)\lvec(90 10)\lvec(80 10)\lvec(80 0)\htext(82
6){\tiny $0$}
\move(50 0)\lvec(60 10)\lvec(60 0)\lvec(50 0)\lfill f:0.8
\htext(56 2){\tiny $0$}

\move(60 0)\lvec(70 10)\lvec(70 0)\lvec(60 0)\lfill f:0.8
\htext(66 2){\tiny $1$}

\move(70 0)\lvec(80 10)\lvec(80 0)\lvec(70 0)\lfill f:0.8
\htext(76 2){\tiny $0$}

\move(80 0)\lvec(90 10)\lvec(90 0)\lvec(80 0)\lfill f:0.8
\htext(86 2){\tiny $1$}
\move(60 10)\lvec(70 10)\lvec(70 20)\lvec(60 20)\lvec(60
10)\htext(63 13){$_2$}

\move(70 10)\lvec(80 10)\lvec(80 20)\lvec(70 20)\lvec(70
10)\htext(73 13){$_2$}

\move(80 10)\lvec(90 10)\lvec(90 20)\lvec(80 20)\lvec(80
10)\htext(83 13){$_2$}
\move(60 20)\lvec(70 20)\lvec(70 30)\lvec(60 30)\lvec(60
20)\htext(63 26){\tiny $3$}

\move(70 20)\lvec(80 20)\lvec(80 30)\lvec(70 30)\lvec(70
20)\htext(73 26){\tiny $3$}

\move(80 20)\lvec(90 20)\lvec(90 30)\lvec(80 30)\lvec(80
20)\htext(83 26){\tiny $3$}
\move(60 20)\lvec(70 20)\lvec(70 25)\lvec(60 25)\lvec(60
20)\htext(63 21){\tiny $3$}

\move(70 20)\lvec(80 20)\lvec(80 25)\lvec(70 25)\lvec(70
20)\htext(73 21){\tiny $3$}

\move(80 20)\lvec(90 20)\lvec(90 25)\lvec(80 25)\lvec(80
20)\htext(83 21){\tiny $3$}
\move(60 30)\lvec(70 30)\lvec(70 40)\lvec(60 40)\lvec(60
30)\htext(63 33){$_2$}

\move(70 30)\lvec(80 30)\lvec(80 40)\lvec(70 40)\lvec(70
30)\htext(73 33){$_2$}

\move(80 30)\lvec(90 30)\lvec(90 40)\lvec(80 40)\lvec(80
30)\htext(83 33){$_2$}
\move(60 40)\lvec(70 40)\lvec(70 50)\lvec(60 50)\lvec(60
40)\htext(62 46){\tiny $0$}

\move(70 40)\lvec(80 40)\lvec(80 50)\lvec(70 50)\lvec(70
40)\htext(72 46){\tiny $1$}

\move(80 40)\lvec(90 40)\lvec(90 50)\lvec(80 50)\lvec(80
40)\htext(82 46){\tiny $0$}
\move(60 40)\lvec(70 50)\lvec(70 40)\lvec(60 40) \htext(66
42){\tiny $1$}

\move(70 40)\lvec(80 50)\lvec(80 40)\lvec(70 40) \htext(76
42){\tiny $0$}

\move(80 40)\lvec(90 50)\lvec(90 40)\lvec(80 40) \htext(86
42){\tiny $1$}
\move(70 50)\lvec(80 50)\lvec(80 60)\lvec(70 60)\lvec(70
50)\htext(73 53){$_2$}

\move(80 50)\lvec(90 50)\lvec(90 60)\lvec(80 60)\lvec(80
50)\htext(83 53){$_2$}
\move(70 60)\lvec(80 60)\lvec(80 65)\lvec(70 65)\lvec(70
60)\htext(73 61){\tiny $3$}

\move(80 60)\lvec(90 60)\lvec(90 65)\lvec(80 65)\lvec(80
60)\htext(83 61){\tiny $3$}
\move(70 60)\lvec(80 60)\lvec(80 70)\lvec(70 70)\lvec(70
60)\htext(73 66){\tiny $3$}

\move(80 60)\lvec(90 60)\lvec(90 70)\lvec(80 70)\lvec(80
60)\htext(83 66){\tiny $3$}
\move(70 70)\lvec(80 70)\lvec(80 80)\lvec(70 80)\lvec(70
70)\htext(73 73){$_2$}

\move(80 70)\lvec(90 70)\lvec(90 80)\lvec(80 80)\lvec(80
70)\htext(83 73){$_2$}
\move(70 80)\lvec(80 90)\lvec(80 80)\lvec(70 80)\htext(76
82){\tiny $0$}

\move(80 80)\lvec(90 90)\lvec(90 80)\lvec(80 80)\htext(86
82){\tiny $1$}

\move(55 5)\linewd 1 \lcir r:6 \move(65 45)\linewd 1 \lcir r:6
\move(75 85)\linewd 1 \lcir r:6
\end{texdraw}}\quad , \quad then
$\overline{Y}=$\raisebox{-0.5\height}{
\begin{texdraw}
\drawdim em \setunitscale 0.13 \linewd 0.5

\move(60 0)\lvec(70 0)\lvec(70 10)\lvec(60 10)\lvec(60 0)\htext(62
6){\tiny $0$}

\move(70 0)\lvec(80 0)\lvec(80 10)\lvec(70 10)\lvec(70 0)\htext(72
6){\tiny $1$}

\move(80 0)\lvec(90 0)\lvec(90 10)\lvec(80 10)\lvec(80 0)\htext(82
6){\tiny $0$}

\move(60 0)\lvec(70 10)\lvec(70 0)\lvec(60 0)\lfill f:0.8
\htext(66 2){\tiny $1$}

\move(70 0)\lvec(80 10)\lvec(80 0)\lvec(70 0)\lfill f:0.8
\htext(76 2){\tiny $0$}

\move(80 0)\lvec(90 10)\lvec(90 0)\lvec(80 0)\lfill f:0.8
\htext(86 2){\tiny $1$}

\move(60 10)\lvec(70 10)\lvec(70 20)\lvec(60 20)\lvec(60
10)\htext(63 13){$_2$}

\move(70 10)\lvec(80 10)\lvec(80 20)\lvec(70 20)\lvec(70
10)\htext(73 13){$_2$}

\move(80 10)\lvec(90 10)\lvec(90 20)\lvec(80 20)\lvec(80
10)\htext(83 13){$_2$}
\move(60 20)\lvec(70 20)\lvec(70 30)\lvec(60 30)\lvec(60
20)\htext(63 26){\tiny $3$}

\move(70 20)\lvec(80 20)\lvec(80 30)\lvec(70 30)\lvec(70
20)\htext(73 26){\tiny $3$}

\move(80 20)\lvec(90 20)\lvec(90 30)\lvec(80 30)\lvec(80
20)\htext(83 26){\tiny $3$}

\move(60 20)\lvec(70 20)\lvec(70 25)\lvec(60 25)\lvec(60
20)\htext(63 21){\tiny $3$}

\move(70 20)\lvec(80 20)\lvec(80 25)\lvec(70 25)\lvec(70
20)\htext(73 21){\tiny $3$}

\move(80 20)\lvec(90 20)\lvec(90 25)\lvec(80 25)\lvec(80
20)\htext(83 21){\tiny $3$}
\move(60 30)\lvec(70 30)\lvec(70 40)\lvec(60 40)\lvec(60
30)\htext(63 33){$_2$}

\move(70 30)\lvec(80 30)\lvec(80 40)\lvec(70 40)\lvec(70
30)\htext(73 33){$_2$}

\move(80 30)\lvec(90 30)\lvec(90 40)\lvec(80 40)\lvec(80
30)\htext(83 33){$_2$}
\move(60 40)\lvec(70 50)\lvec(60 50)\lvec(60 40)\htext(62
46){\tiny $0$}

\move(70 40)\lvec(80 40)\lvec(80 50)\lvec(70 50)\lvec(70
40)\htext(72 46){\tiny $1$}

\move(80 40)\lvec(90 40)\lvec(90 50)\lvec(80 50)\lvec(80
40)\htext(82 46){\tiny $0$}
\move(70 40)\lvec(80 50)\lvec(80 40)\lvec(70 40) \htext(76
42){\tiny $0$}

\move(80 40)\lvec(90 50)\lvec(90 40)\lvec(80 40) \htext(86
42){\tiny $1$}
\move(70 50)\lvec(80 50)\lvec(80 60)\lvec(70 60)\lvec(70
50)\htext(73 53){$_2$}

\move(80 50)\lvec(90 50)\lvec(90 60)\lvec(80 60)\lvec(80
50)\htext(83 53){$_2$}
\move(70 60)\lvec(80 60)\lvec(80 65)\lvec(70 65)\lvec(70
60)\htext(73 61){\tiny $3$}

\move(80 60)\lvec(90 60)\lvec(90 65)\lvec(80 65)\lvec(80
60)\htext(83 61){\tiny $3$}
\move(70 60)\lvec(80 60)\lvec(80 70)\lvec(70 70)\lvec(70
60)\htext(73 66){\tiny $3$}

\move(80 60)\lvec(90 60)\lvec(90 70)\lvec(80 70)\lvec(80
60)\htext(83 66){\tiny $3$}
\move(70 70)\lvec(80 70)\lvec(80 80)\lvec(70 80)\lvec(70
70)\htext(73 73){$_2$}

\move(80 70)\lvec(90 70)\lvec(90 80)\lvec(80 80)\lvec(80
70)\htext(83 73){$_2$}
\move(70 80)\lvec(80 90)\lvec(80 80)\lvec(70 80)\htext(76
82){\tiny $0$}

\move(80 80)\lvec(90 90)\lvec(90 80)\lvec(80 80)\htext(86
82){\tiny $1$}
\end{texdraw}}\quad .
\end{center}\vskip 5mm }
\end{ex}

\begin{prop}\label {QbarY Y=1} Let $Y$ be a reduced proper Young wall
in ${\mathcal Y}(\Lambda)$.
Suppose that ${\rm wt}(\overline{Y})={\rm wt}(Y)+r\alpha_i$ for
some $i\in I$ and $r\geq 1$. Then we have
\begin{equation*}
f_i^{(r)}\overline{Y}=Y+\sum_{\substack{Z\in {\mathcal Z}(\Lambda) \\
{\rm wt}(Z)={\rm wt}(Y) \\ Z\neq Y }}Q_{\overline{Y},Z}(q)Z.
\end{equation*}
That is, we have $Q_{\overline{Y},Y}(q)=1$.
\end{prop}
\pf By definition of $\overline{Y}$, there exists a unique
sequence of proper Young walls $Y_0=\overline{Y},\cdots,Y_r=Y$
such that for $1\leq k\leq r$

(i) $Y_k=Y_{k-1}\swarrow b_k$ (up to scalar multiplication) for
some admissible $i$-slot $b_k$ of $Y_{k-1}$,

(ii) there exists no admissible $i$-slot located to the left of
$b_k$.

In other words, $\{\,b_k\,|\,1\leq k\leq r\,\}$ are added to
$\overline{Y}$ from left to right and from bottom to top with no
admissible $i$-slot to the left of each $b_k$.

Suppose that $b_k$'s are of type I or III. Then it is easy to see
that
\begin{equation*}
Q_{\overline{Y},Y}^{\circ}(q)=\prod_{k=0}^{r-1}
Q_{Y_{k},Y_{k+1}}(q)=q_i^{-\binom{r}{2}},
\end{equation*}
which implies that $Q_{\overline{Y},Y}(q)=1$ by Lemma \ref{QYZ}.

Suppose that $b_k$'s are of type II. Let $J_1$, $J_2$, $J_3$ be
the sets given in \eqref{JST}. Set $l=|J_1|$, $m=|J_2|$, and
$n=|J_3|$. By definition of $\overline{Y}$, $m\leq 1$ and if
$m=1$, then $J_2=\{\,1\,\}$ and $b_1$ is placed on the column
which is part of the ground-state wall. Also, (ii) implies that
$n\leq 1$ and that if $n=1$, then $J_3=\{\,r\,\}$. Note that
$S=T=\emptyset$. Thus we have
\begin{equation*}
Q^{\circ}_{\overline{Y},Y}(q)=[2]_i^{l}
q_i^{-4\binom{l}{2}-\binom{m}{2}-\binom{n}{2}-2l(m+n)-mn},
\end{equation*}
which implies that $Q_{\overline{Y},Y}(q)=1$ by Lemma \ref{QYZ} .
\qed\vskip 3mm

Let $Y$ be a proper Young wall in ${\mathcal Z}(\Lambda)$. Let $L$
be a ladder which has a nontrivial intersection with $Y$. Suppose
that there are $r$ many $i$-blocks in $Y\cap L$ for some $r\geq 0$
and $i\in I$. Move these $i$-blocks to the first $r$ many
$i$-slots in $L$ from the bottom. Repeat this procedure ladder by
ladder until no block can be moved downward along a ladder. Then
we obtain another proper Young wall $Y^{R}$, which we call the
{\it reduced form of $Y$}. By definition, $Y^{R}$ is a reduced
proper Young wall. Moreover, we have $|Y^R|\unrhd|Y|$ and the
equality holds if and only if $Y$ is reduced.

\begin{ex}{\rm \mbox{}\vskip 5mm

\begin{center}
$Y=$\raisebox{-0.3\height}{
\begin{texdraw}
\drawdim em \setunitscale 0.12 \linewd 0.5

\move(70 0)\lvec(80 0)\lvec(80 10)\lvec(70 10)\lvec(70 0)\htext(72
6){\tiny $1$}

\move(80 0)\lvec(90 0)\lvec(90 10)\lvec(80 10)\lvec(80 0)\htext(82
6){\tiny $0$}

\move(70 0)\lvec(80 10)\lvec(80 0)\lvec(70 0)\lfill f:0.8
\htext(76 2){\tiny $0$}

\move(80 0)\lvec(90 10)\lvec(90 0)\lvec(80 0)\lfill f:0.8
\htext(86 2){\tiny $1$}

\move(80 10)\lvec(90 10)\lvec(90 20)\lvec(80 20)\lvec(80
10)\htext(83 13){$_2$}
\move(80 20)\lvec(90 20)\lvec(90 30)\lvec(80 30)\lvec(80
20)\htext(83 26){\tiny $3$}

\move(80 20)\lvec(90 20)\lvec(90 25)\lvec(80 25)\lvec(80
20)\htext(83 21){\tiny $3$}
\move(80 30)\lvec(90 30)\lvec(90 40)\lvec(80 40)\lvec(80
30)\htext(83 33){$_2$}
\move(80 40)\lvec(90 40)\lvec(90 50)\lvec(80 50)\lvec(80
40)\htext(82 46){\tiny $0$}
\move(80 40)\lvec(90 50)\lvec(90 40)\lvec(80 40) \htext(86
42){\tiny $1$}
\move(80 50)\lvec(90 50)\lvec(90 60)\lvec(80 60)\lvec(80
50)\htext(83 53){$_2$}
\move(80 60)\lvec(90 60)\lvec(90 65)\lvec(80 65)\lvec(80
60)\htext(83 61){\tiny $3$}
\move(80 60)\lvec(90 60)\lvec(90 70)\lvec(80 70)\lvec(80
60)\htext(83 66){\tiny $3$}
\move(80 70)\lvec(90 70)\lvec(90 80)\lvec(80 80)\lvec(80
70)\htext(83 73){$_2$}
\move(80 80)\lvec(90 90)\lvec(80 90)\lvec(80 80)\htext(82
86){\tiny $0$}

\move(85 85)\linewd 0.5 \lcir r:6

\move(65 5)\linewd 0.5 \lcir r:6

\linewd 0.5

\move(82 79) \arrowheadtype t:F \arrowheadsize l:4 w:2 \avec(68
11)
\end{texdraw}}\hskip 5mm $\longrightarrow$ \hskip 5mm
\raisebox{-0.3\height}{
\begin{texdraw}
\drawdim em \setunitscale 0.12 \linewd 0.5

\move(70 0)\lvec(80 0)\lvec(80 10)\lvec(70 10)\lvec(70 0)\htext(72
6){\tiny $1$}

\move(80 0)\lvec(90 0)\lvec(90 10)\lvec(80 10)\lvec(80 0)\htext(82
6){\tiny $0$}

\move(60 0)\lvec(70 10)\lvec(70 0)\lvec(60 0)\lfill f:0.8
\htext(66 2){\tiny $1$}

\move(70 0)\lvec(80 10)\lvec(80 0)\lvec(70 0)\lfill f:0.8
\htext(76 2){\tiny $0$}

\move(80 0)\lvec(90 10)\lvec(90 0)\lvec(80 0)\lfill f:0.8
\htext(86 2){\tiny $1$}

\move(80 10)\lvec(90 10)\lvec(90 20)\lvec(80 20)\lvec(80
10)\htext(83 13){$_2$}
\move(80 20)\lvec(90 20)\lvec(90 30)\lvec(80 30)\lvec(80
20)\htext(83 26){\tiny $3$}

\move(80 20)\lvec(90 20)\lvec(90 25)\lvec(80 25)\lvec(80
20)\htext(83 21){\tiny $3$}
\move(80 30)\lvec(90 30)\lvec(90 40)\lvec(80 40)\lvec(80
30)\htext(83 33){$_2$}
\move(80 40)\lvec(90 40)\lvec(90 50)\lvec(80 50)\lvec(80
40)\htext(82 46){\tiny $0$}
\move(80 40)\lvec(90 50)\lvec(90 40)\lvec(80 40) \htext(86
42){\tiny $1$}
\move(80 50)\lvec(90 50)\lvec(90 60)\lvec(80 60)\lvec(80
50)\htext(83 53){$_2$}
\move(80 60)\lvec(90 60)\lvec(90 65)\lvec(80 65)\lvec(80
60)\htext(83 61){\tiny $3$}
\move(80 60)\lvec(90 60)\lvec(90 70)\lvec(80 70)\lvec(80
60)\htext(83 66){\tiny $3$}
\move(80 70)\lvec(90 70)\lvec(90 80)\lvec(80 80)\lvec(80
70)\htext(83 73){$_2$}
\move(60 0)\lvec(70 10)\lvec(60 10)\lvec(60 0)\htext(62 6){\tiny
$0$}

\move(85 65)\linewd 0.5 \lellip rx:6 ry:16

\move(75 25)\linewd 0.5 \lellip rx:6 ry:16

\linewd 0.5

\move(80 58) \arrowheadtype t:F \arrowheadsize l:4 w:2 \avec(76
40)
\end{texdraw}}\hskip 5mm $\longrightarrow$ \hskip 5mm
\raisebox{-0.35\height}{
\begin{texdraw}
\drawdim em \setunitscale 0.12 \linewd 0.5

\move(60 0)\lvec(70 0)\lvec(70 10)\lvec(60 10)\lvec(60 0)\htext(62
6){\tiny $0$}

\move(70 0)\lvec(80 0)\lvec(80 10)\lvec(70 10)\lvec(70 0)\htext(72
6){\tiny $1$}

\move(80 0)\lvec(90 0)\lvec(90 10)\lvec(80 10)\lvec(80 0)\htext(82
6){\tiny $0$}

\move(60 0)\lvec(70 10)\lvec(70 0)\lvec(60 0)\lfill f:0.8
\htext(66 2){\tiny $1$}

\move(70 0)\lvec(80 10)\lvec(80 0)\lvec(70 0)\lfill f:0.8
\htext(76 2){\tiny $0$}

\move(80 0)\lvec(90 10)\lvec(90 0)\lvec(80 0)\lfill f:0.8
\htext(86 2){\tiny $1$}

\move(70 10)\lvec(80 10)\lvec(80 20)\lvec(70 20)\lvec(70
10)\htext(73 13){$_2$}

\move(80 10)\lvec(90 10)\lvec(90 20)\lvec(80 20)\lvec(80
10)\htext(83 13){$_2$}
\move(70 20)\lvec(80 20)\lvec(80 30)\lvec(70 30)\lvec(70
20)\htext(73 26){\tiny $3$}

\move(80 20)\lvec(90 20)\lvec(90 30)\lvec(80 30)\lvec(80
20)\htext(83 26){\tiny $3$}

\move(70 20)\lvec(80 20)\lvec(80 25)\lvec(70 25)\lvec(70
20)\htext(73 21){\tiny $3$}

\move(80 20)\lvec(90 20)\lvec(90 25)\lvec(80 25)\lvec(80
20)\htext(83 21){\tiny $3$}
\move(70 30)\lvec(80 30)\lvec(80 40)\lvec(70 40)\lvec(70
30)\htext(73 33){$_2$}

\move(80 30)\lvec(90 30)\lvec(90 40)\lvec(80 40)\lvec(80
30)\htext(83 33){$_2$}
\move(80 40)\lvec(90 40)\lvec(90 50)\lvec(80 50)\lvec(80
40)\htext(82 46){\tiny $0$}
\move(80 40)\lvec(90 50)\lvec(90 40)\lvec(80 40) \htext(86
42){\tiny $1$}
\end{texdraw}}\hskip 5mm$=Y^R$
\end{center}\vskip 5mm

}
\end{ex}

\begin{lem}\label {ord}{\rm (cf.\cite{JK})}
Let $Y$ be a reduced proper Young wall in ${\mathcal Y}(\Lambda)$
and let $Z$ be a proper Young wall in ${\mathcal Z}(\Lambda)$ such
that $|\overline{Y}|\unrhd |Z^{R}|$. Suppose that ${\rm
wt}(\overline{Y})={\rm wt}(Y)+r\alpha_i$ for some $i\in I$ and
$r\geq 1$. Then, for each $W\in{\mathcal Z}(\Lambda)$ appearing in
the expression of $f_i^{(r)}Z$, we have

{\rm (a)} $|Y|\unrhd|W^{R}|$.

{\rm (b)} If $|Y|=|W|$, then $|\overline{Y}|=|Z|$ and $Z$ is a
reduced proper Young wall.

{\rm (c)} If $|Y|=|W|$ and $\overline{Y}=Z$, then $Y=W$.
\end{lem}
\pf (a) Let $L$ be the left-most ladder which has a nontrivial
intersection with $Y$. Denote by $y_p,y_{p-1},\cdots,y_{p-s}$
($s\geq 0$) the first $s+1$ columns in $Y$ which meet $L$, and
denote by
$\overline{y}_p,\overline{y}_{p-1},\cdots,\overline{y}_{p-s}$
($s\geq 0$) the corresponding columns in $\overline{Y}$. Since
$Z^R=(z^R_k)_{k=0}^{\infty}$ is reduced, we have
\begin{equation}
|\overline{y}_{p-t}|\geq|z^R_{p-t}|\hskip 5mm\text{for $0\leq
t\leq s$}.
\end{equation}
Since $W=(w_k)_{k=0}^{\infty}$ (and hence
$W^R=(w_k^R)_{k=0}^{\infty}$) is given by adding $r$ many
$i$-blocks on $Z$, we also have
\begin{equation}
|y_{p-t}|\geq|w^R_{p-t}|\hskip 5mm\text{for $0\leq t\leq s$}.
\end{equation}
Therefore, if $s'\leq s$, then
\begin{equation}
\sum_{t=0}^{s'}|y_{p-t}|\geq\sum_{t=0}^{s'}|w^R_{p-t}|,
\end{equation}
and if $s< s'\leq p$, then
\begin{equation}
\sum_{t=0}^{s'}|y_{p-t}|=\sum_{t=0}^{s'}|\overline{y}_{p-t}|+r\geq
\sum_{t=0}^{s'}|z^R_{p-t}|+r\geq\sum_{t=0}^{s'}|w^R_{p-t}|.
\end{equation}
Hence we conclude $|Y|\unrhd|W^R|$.

(b) Suppose that $|Y|=|W|$. Since $|\overline{y}_p|\geq|z_p|$ and
$w_p$ is obtained by adding some $i$-blocks on $z_p$, we have
$\overline{y}_p=z_p$ and $y_p=w_p$. Suppose that for $0\leq u\leq
t < s$,
\begin{equation}
\overline{y}_{p-u}=z_{p-u},\quad y_{p-u}=w_{p-u}.
\end{equation}
Note that
\begin{equation}
\sum_{u=0}^{t+1}|\overline{y}_{p-u}|\geq\sum_{u=0}^{t+1}|z^R_{p-u}|
\geq\sum_{u=0}^{t+1}|z_{p-u}|.
\end{equation}
By our hypothesis, we have
$|\overline{y}_{p-t-1}|\geq|z_{p-t-1}|$. Since
$|y_{p-t-1}|=|w_{p-t-1}|$ and $w_{p-t-1}$ is obtained by adding
some $i$-blocks on $z_{p-t-1}$, we have
$\overline{y}_{p-t-1}=z_{p-t-1}$ and $y_{p-t-1}=w_{p-t-1}$. By
induction, $\overline{y}_{p-u}=z_{p-u}$ and $y_{p-u}=w_{p-u}$ for
$0\leq u\leq s$. Since all the $i$-blocks are added on
$(\overline{y}_k)_{k=p-s}^{\infty}$ and $(w_k)_{k=p-s}^{\infty}$,
we have $(\overline{y}_k)_{k=0}^{p-s-1}=(y_k)_{k=0}^{p-s-1}$ and
$(w_k)_{k=0}^{p-s-1}=(z_k)_{k=0}^{p-s-1}$, which implies
$|y_k|=|\overline{y}_k|=|w_k|=|z_k|$ for all $0\leq k\leq p-s-1$.
Hence, $|\overline{Y}|=|Z|=|Z^R|$ and $Z$ is reduced.

(c) follows directly from the proof of (b) \qed\vskip 3mm

Let $Y$ be a reduced proper Young wall in ${\mathcal Y}(\Lambda)$.
There exists a unique sequence of reduced proper Young walls
$\{Y_{k}\}_{k=0}^N$ such that $Y_{0}=Y$, $Y_{1} =
\overline{Y_{0}}$, $\cdots$, $Y_{k+1}=\overline{Y_{k}}$, $\cdots$,
$Y_{N}=\overline{Y_{N-1}} = Y_{\Lambda}$. Suppose that
$Y_{k}=\overline{Y_{k-1}}$ is obtained by removing $r_{k}$ many
$i_{k}$-blocks from $Y_{k-1}$ ($1\leq k\leq N$). We define
\begin{equation}
A(Y)=f_{i_1}^{(r_1)}\cdots f_{i_N}^{(r_N)}Y_{\Lambda}\in
V(\Lambda)_{\mathbb{A}}.
\end{equation}
\begin{ex}{\rm If $\frak{g}=A_4^{(2)}$, $\Lambda=\Lambda_0$ and \vskip 5mm

\begin{center}
$Y=$\raisebox{-0.5\height}{
\begin{texdraw}
\drawdim em \setunitscale 0.13 \linewd 0.5

\move(20 0)\lvec(30 0)\lvec(30 10)\lvec(20 10)\lvec(20 0)\htext(23
1){\tiny $0$}

\move(30 0)\lvec(40 0)\lvec(40 10)\lvec(30 10)\lvec(30 0)\htext(33
1){\tiny $0$}

\move(40 0)\lvec(50 0)\lvec(50 10)\lvec(40 10)\lvec(40 0)\htext(43
1){\tiny $0$}
\move(20 0)\lvec(30 0)\lvec(30 5)\lvec(20 5)\lvec(20 0)\lfill
f:0.8 \htext(23 6){\tiny $0$}

\move(30 0)\lvec(40 0)\lvec(40 5)\lvec(30 5)\lvec(30 0)\lfill
f:0.8 \htext(33 6){\tiny $0$}

\move(40 0)\lvec(50 0)\lvec(50 5)\lvec(40 5)\lvec(40 0)\lfill
f:0.8 \htext(43 6){\tiny $0$}
\move(30 10)\lvec(40 10)\lvec(40 20)\lvec(30 20)\lvec(30
10)\htext(33 13){$_1$}

\move(40 10)\lvec(50 10)\lvec(50 20)\lvec(40 20)\lvec(40
10)\htext(43 13){$_1$}

\move(30 20)\lvec(40 20)\lvec(40 30)\lvec(30 30)\lvec(30
20)\htext(33 23){$_2$}

\move(40 20)\lvec(50 20)\lvec(50 30)\lvec(40 30)\lvec(40
20)\htext(43 23){$_2$}

\move(30 30)\lvec(40 30)\lvec(40 40)\lvec(30 40)\lvec(30
30)\htext(33 33){$_1$}

\move(40 30)\lvec(50 30)\lvec(50 40)\lvec(40 40)\lvec(40
30)\htext(43 33){$_1$}

\move(30 40)\lvec(40 40)\lvec(40 45)\lvec(30 45)\lvec(30
40)\htext(33 41){\tiny $0$}

\move(40 40)\lvec(50 40)\lvec(50 45)\lvec(40 45)\lvec(40
40)\htext(43 41){\tiny $0$}

\move(30 40)\lvec(40 40)\lvec(40 50)\lvec(30 50)\lvec(30
40)\htext(33 46){\tiny $0$}

\move(40 40)\lvec(50 40)\lvec(50 50)\lvec(40 50)\lvec(40
40)\htext(43 46){\tiny $0$}

\move(40 50)\lvec(50 50)\lvec(50 60)\lvec(40 60)\lvec(40
50)\htext(43 53){$_1$}

\move(40 60)\lvec(50 60)\lvec(50 70)\lvec(40 70)\lvec(40
60)\htext(43 63){$_2$}

\move(40 70)\lvec(50 70)\lvec(50 80)\lvec(40 80)\lvec(40
70)\htext(43 73){$_1$}

\move(40 80)\lvec(50 80)\lvec(50 85)\lvec(40 85)\lvec(40
80)\htext(43 81){\tiny $0$}

\end{texdraw}}\quad , then we have
 $A(Y)=f_0^{(4)} f_1^{(2)} f_2^{(2)} f_1^{(2)} f_0^{(3)}
f_1 f_2 f_1 f_0 Y_{\Lambda_0}$.
\end{center}\vskip 5mm

}
\end{ex}

By definition, $\overline{A(Y)}=A(Y)$. Write
\begin{equation}
A(Y)=\sum_{Z\in{\mathcal Z}(\Lambda)} A_{Y,Z}(q)Z,
\end{equation}
where $A_{Y,Z}(q)\in \mathbb{Q}(q)$. Then, the coefficients
$A_{Y,Z}(q)$ satisfy the following properties.
\begin{prop}\label {AYZ} Let $Y$ be a reduced proper Young wall in
${\mathcal Y}(\Lambda)$. Then, for a proper Young wall
$Z\in{\mathcal Z}(\Lambda)$, we have

{\rm (a)} $A_{Y,Z}(q)\in \mathbb{Z}[q,q^{-1}]$,

{\rm (b)} $A_{Y,Z}(q)=0 \text{ unless } |Y|\unrhd |Z^R|$ and ${\rm
wt}(Y)={\rm wt}(Z)$,

{\rm (c)} if $A_{Y,Z}(q)\neq 0$ and $|Y|=|Z|$, then $Y=Z$ and
$A_{Y,Y}(q)=1$.
\end{prop}
\pf We will use induction on $l$, the number of blocks in $Y$
which have been added to $Y_{\Lambda}$. If $l=1$, it is clear.
Suppose that $l>1$, and (a)--(c) hold for $l'<l$. If
$A(Y)=f_{i_1}^{(r_1)}\cdots f_{i_N}^{(r_N)}Y_{\Lambda}$ for some
$N\geq 1$, then we have
\begin{equation}
\begin{split}
A(Y)&=f_{i_1}^{(r_1)}A(\overline{Y})
=\sum_{|\overline{Y}|\unrhd|Z^R|}A_{\overline{Y},Z}(q)f_{i_1}^{(r_1)}Z
\\
&=\sum_{|\overline{Y}|\unrhd|Z^R|}A_{\overline{Y},Z}(q) \left (
\sum_{|Y|\unrhd|W^R|}Q_{Z,W}(q)W\right ) \hskip 1cm\text{by Lemma
\ref{ord} (a)}
\\
&=\sum_{|Y|\unrhd|W^R|}\left
(\sum_{|\overline{Y}|\unrhd|Z^R|}A_{\overline{Y},Z}(q)Q_{Z,W}(q)
\right )W.
\end{split}
\end{equation}
By induction hypothesis and Lemma \ref{QYZ}, we have
\begin{equation}
A_{Y,W}(q)=\sum_{|\overline{Y}|\unrhd|Z^R|}A_{\overline{Y},Z}(q)Q_{Z,W}(q)
\in\mathbb{Z}[q,q^{-1}],
\end{equation}
and $A_{Y,W}(q)=0$ unless $|Y|\unrhd |W^R|$ and ${\rm wt}(Y)={\rm
wt}(W)$.

If $A_{Y,W}(q)\neq 0$ and $|Y|=|W|$, then Lemma \ref{ord} (b)
implies that $|Z|=|\overline{Y}|$ for $A_{\overline{Y},Z}(q)\neq
0$. Hence, $Z=\overline{Y}$ by induction hypothesis. Finally, we
have $Y=W$ by Lemma \ref{ord} (c), and  hence
$A_{Y,Y}(q)=A_{\overline{Y},\overline{Y}}(q)Q_{\overline{Y},Y}(q)=1$
by Proposition \ref{QbarY Y=1}, which completes the induction
argument. \qed\vskip 3mm

For proper Young walls $Y=(y_k)_{k=0}^{\infty}$ and
$Z=(z_k)_{k=0}^{\infty}$, we define $|Y|> |Z|$ if there exists
$k\geq 0$ such that $|y_k| > |z_k|$ and $|y_l|=|z_l|$ for all
$l>k$. Thus we have a total ordering on the set of partitions.
Note that $|Y|\unrhd|Z|$ implies $|Y|\geq |Z|$. Now we define a
total ordering $>$ on the set $\mathcal{Z}(\Lambda)$ of proper
Young walls as follows. First, we fix an arbitrary total ordering
$\succ$ on the set of proper Young walls with the same associated
partition. Then we define
\begin{multline}
Y>Z\quad \text{if and only if} \\
{\rm (i)} \,\,|Y|>|Z|\,\,\,\,\text{ or }\,\, {\rm (ii)}\,\,|Y|=|Z|
\text{ and } Y \succ Z\,\,.
\end{multline}
For example, if $|Y|\unrhd|Z|$ and $|Y|\neq|Z|$, then we have
$Y>Z$.

Let $Y$ be a reduced proper Young wall in $\mathcal{Y}(\Lambda)$.
By Proposition \ref{AYZ}, we may write
\begin{equation*}
A(Y)=Y+\sum_{Y>Z}A_{Y,Z}(q)Z.
\end{equation*}
It follows that the set ${\mathcal A}(\Lambda)=\{\,A(Y)\,|\,Y\in
{\mathcal Y}_{\circ}(\Lambda)\,\}$ is linearly independent over
$\mathbb{Q}(q)$. Since $\dim V(\Lambda)_{\lambda}= |{\mathcal
Y}(\Lambda)_{\lambda}|$ for $\lambda \leq \Lambda$, we conclude
that $\mathcal{A}(\Lambda)$ is a $\mathbb{Q}(q)$-basis of
$V(\Lambda)$.

Let $\mathcal{G}(\Lambda)=\{\,G(Y)\,|\,Y\in {\mathcal
Y}(\Lambda)\,\}$ be the global basis of $V(\Lambda)_{\mathbb{A}}$.
Then for each reduced proper Young wall $Y\in {\mathcal
Y}(\Lambda)$, we may write
\begin{equation}
G(Y)=\sum_{Z\in{\mathcal Z}(\Lambda)} G_{Y,Z}(q)Z\in
V(\Lambda)_{\mathbb{A}}\cap {\mathcal L}(\Lambda)
\end{equation}
for some $G_{Y,Z}(q)\in \mathbb{A}_0$. Since
$\mathcal{G}(\Lambda)$ is an $\mathbb{A}$-basis of
$V(\Lambda)_{\mathbb{A}}$, $G(Y)$ can be expressed as an
$\mathbb{A}$-linear combination of the vectors
$f_{i_1}^{(r_1)}\cdots f_{i_N}^{(r_N)}Y_{\Lambda}$. By Lemma
\ref{QYZ}, it is easy to see that
$G_{Y,Z}(q)\in\mathbb{Q}[q,q^{-1}]$. Moreover, since $G(Y)\equiv Y
\mod{q{\mathcal L}(\Lambda)}$, the coefficients $G_{Y,Z}(q)$
satisfy the following properties:

\hskip 2mm(i) $G_{Y,Z}(q)\in\mathbb{Q}[q]$,

\hskip 2mm(ii) $G_{Y,Z}(q)\in q\mathbb{Q}[q]$ unless $Y=Z$,

\hskip 2mm(iii) $G_{Y,Y}(q)=1$.

On the other hand, since $\mathcal{G}(\Lambda)$ and
$\mathcal{A}(\Lambda)$ are both $\mathbb{Q}(q)$-basis of
$V(\Lambda)$, there exists the transition matrix
$H=(H_{Y,W}(q))_{Y,W\in{\mathcal Y}(\Lambda)}$ such that
\begin{equation}\label {GHA}
G(Y)=\sum_{W\in {\mathcal Y}(\Lambda)} H_{Y,W}(q)A(W),
\end{equation}
where the indices are decreasing with respect to the total
ordering $>$ on $\mathcal{Y}(\Lambda)$. Since
$\overline{G(Y)}=G(Y)$ and $\overline{A(W)}=A(W)$, we have
$H_{Y,W}(q)=H_{Y,W}(q^{-1})$ for all $Y,W\in
\mathcal{Y}(\Lambda)$. The following proposition provides the last
ingredients for our algorithm.

\begin{prop}\label {HYW}\mbox{}

{\rm (a)} The coefficients $H_{Y,W}(q)$ satisfy the following
properties:

\hskip 4mm{\rm (i)} $H_{Y,W}(q)\in\mathbb{Q}[q,q^{-1}]$.

\hskip 4mm{\rm (ii)} $H_{Y,W}(q)=0$ unless $Y\geq W$ and ${\rm
wt}(Y)={\rm wt}(W)$.

\hskip 4mm{\rm (iii)} $H_{Y,Y}(q)=1$.

{\rm (b)} The set $\mathcal{A}(\Lambda)$ is an $\mathbb{A}$-basis
of $V(\Lambda)_{\mathbb{A}}$.
\end{prop}
\pf (a) Consider the following square matrices indexed by
${\mathcal Y}(\Lambda)$
\begin{equation}
G=(G_{Y,Z}(q)),\quad A=(A_{W,Z}(q)),
\end{equation}
where the indices are given by the total ordering $>$ in a
decreasing manner. Then by \eqref{GHA}, $G=HA$. Since $A$ is an
upper triangular matrix whose diagonal entries are all 1, we
conclude that $A$ is invertible and the entries of $A^{-1}$ are in
$\mathbb{Q}[q,q^{-1}]$. It follows that $H=GA^{-1}$ and
$H_{Y,W}(q)\in\mathbb{Q}[q,q^{-1}]$ for all
$Y,W\in\mathcal{Y}(\Lambda)$. This proves (i).

Next, let $W$ be the reduced proper Young wall in
$\mathcal{Y}(\Lambda)$ that is maximal with respect to the total
ordering $>$ on $\mathcal{Y}(\Lambda)$ among the ones with
$H_{Y,W}(q)\neq 0$. By the maximality of $W$ and Proposition
\ref{AYZ} (b), we have
$G_{Y,W}(q)=H_{Y,W}(q)A_{W,W}(q)=H_{Y,W}(q).$ Since
$G_{Y,W}(q)\in\mathbb{Q}[q]$ and
$G_{Y,W}(q)=H_{Y,W}(q)=H_{Y,W}(q^{-1})=G_{Y,W}(q^{-1})$,
$G_{Y,W}(q)$ must be a constant. It follows that $Y=W$ and
$H_{Y,Y}(q)=G_{Y,Y}(q)=1$. This proves (ii) and (iii).

(b) By (i) and \eqref{GHA}, every element of
$\mathcal{G}(\Lambda)$ can be expressed as an $\mathbb{A}$-linear
combination of the elements in $\mathcal{A}(\Lambda)$. Hence,
$\mathcal{A}(\Lambda)$ is an $\mathbb{A}$-basis of
$V(\Lambda)_{\mathbb{A}}$. \qed \vskip 3mm

Observe that, by Proposition \ref{HYW} (a), $H$ is invertible and
$H^{-1}$ is also an upper triangular matrix whose diagonal entries
are all $1$. Hence for each reduced proper Young wall
$Y\in{\mathcal Y}(\Lambda)_{\lambda}$ ($\lambda\leq \Lambda$),
$A(Y)$ can be expressed uniquely as
\begin{equation}\label {H'YZ}
A(Y)=G(Y)+\sum_{\substack{Z\in {\mathcal Y}(\Lambda)_{\lambda}
\\ Y>Z }}H'_{Y,Z}(q)G(Z),
\end{equation}
for some $H'_{Y,Z}(q)\in \mathbb{Q}[q,q^{-1}]$ such that
$H'_{Y,Z}(q)=H'_{Y,Z}(q^{-1})$.

Now, we are ready to give a generalized version of
Lascoux-Leclerc-Thibon algorithm for constructing the global basis
element $G(Y)$ (cf.\cite{LLT}).

Fix a weight $\lambda\leq \Lambda$ of $V(\Lambda)$, and we list
all the reduced proper Young walls in ${\mathcal
Y}(\Lambda)_{\lambda}$ using the total ordering $>$:
\begin{equation*}
Y_1>Y_2>\cdots>Y_l.
\end{equation*}
We will construct the basis element $G(Y_k)$ ($1\leq k\leq l$) in
a recursive way.

First, by \eqref{H'YZ}, we have $G(Y_l)=A(Y_l)$ because $Y_l$ is
the minimal element. Suppose that we have computed
$G(Y_{k+1}),\cdots,G(Y_l)$. Then, by \eqref{H'YZ}, there exist
uniquely determined coefficients
$\gamma_s(q)\in\mathbb{Q}[q,q^{-1}]$ ($k< s\leq l$) with
$\gamma_s(q)=\gamma_s(q^{-1})$ such that
\begin{multline}\label {GAY}
G(Y_k)=A(Y_k)-\gamma_{k+1}(q)G(Y_{k+1}) \\
-\gamma_{k+2}(q)G(Y_{k+2})-\cdots-\gamma_l(q)G(Y_l).
\end{multline}
Since $G(Y_k)\equiv Y_k \mod{q{\mathcal L}(\Lambda)}$ and
$\gamma_s(q)=\gamma_s(q^{-1})$, $\gamma_s(q)$ ($k< s\leq l$) are
determined recursively as follows:
\begin{itemize}
\item[(G.1)] if $A_{Y_k,Y_{k+1}}(q)=\sum_{i=-r}^{r'}a_iq^{i}$, then
$\gamma_{k+1}(q)=\sum_{i=1}^{r}a_{-i}(q^i +q^{-i}) + a_0$.

\item[(G.2)] if the coefficient of $Y_s$ ($s>k+1$) in $A(Y_k)-\sum_{p=k+1}^{s-1}
\gamma_p(q)G(Y_{p})$ is given by $\sum_{i=-r}^{r'}a_iq^{i}$, then
$\gamma_{s}(q)=\sum_{i=1}^{r}a_{-i}(q^i +q^{-i}) + a_0$.
\end{itemize}
Using this procedure, one can construct $G(Y_k)$ $(k=1,\cdots,l)$.

To summarize, we obtain the {\it generalized
Lascoux-Leclerc-Thibon algorithm} :
\begin{thm}\label {LLT} Let $Y$ be a reduced proper Young wall in ${\mathcal
Y}(\Lambda)$. Then the corresponding global basis element $G(Y)$
can be constructed recursively using the algorithm given in
\eqref{GAY}, {\rm ($G.1$)} and {\rm ($G.2$)}. Moreover $G(Y)$ has
the form
\begin{equation}\label {G(Y)}
G(Y)=Y+\sum_{\substack{Z\in{\mathcal Z}(\Lambda) \\
|Y|\rhd|Z|}}G_{Y,Z}(q)Z,
\end{equation}
where  $G_{Y,Z}(q)\in q\mathbb{Z}[q]$ for $Y\neq Z$.\qed
\end{thm}

By the construction of $G(Y)$ and Proposition \ref{AYZ} (b), we
have
\begin{equation}
G_{Y,Z}(q)=0 \quad \text{unless $|Y|\unrhd|Z^R|$}.
\end{equation}
Hence, we can also apply the modified algorithm introduced in
\cite{Mathas} as follows.

Let $Y\in{\mathcal Y}(\Lambda)$ be a reduced proper Young wall and
suppose that $G(\overline{Y})$  and $G(Y')$ ($Y>Y'$) have been
constructed. Set $C(Y)=f_i^{(r)}G(\overline{Y})$ where ${\rm
wt}(Y)={\rm wt}(\overline{Y})-r\alpha_i$. Note that
$C(Y)=\overline{C(Y)}$. By Lemma \ref{ord} (a) and Proposition
\ref{QbarY Y=1}, we have
\begin{equation}
C(Y)=Y+\sum_{|Y|\unrhd|Z^R|}C_{Y,Z}(q)Z
\end{equation}
for some $C_{Y,Z}(q)\in\mathbb{Q}[q,q^{-1}]$. Hence, for each
$Y'\in{\mathcal Y}(\Lambda)$, there exists uniquely determined
coefficients $\zeta_{Y,Y'}(q)\in\mathbb{Q}[q,q^{-1}]$ with
$\zeta_{Y,Y'}(q)=\zeta_{Y,Y'}(q^{-1})$ such that
\begin{equation}\label{GCY}
G(Y)=C(Y)-\sum_{Y>Y'}\zeta_{Y,Y'}(q)G(Y').
\end{equation}
Since $G(Y)\equiv Y \mod{q{\mathcal L}(\Lambda)}$ and
$\zeta_{Y,Y'}(q)=\zeta_{Y,Y'}(q^{-1})$, the coefficients
$\zeta_{Y,Y'}(q)$ are determined recursively as follows:
\begin{itemize}
\item[$(G'.1)$] if $Y'$ is the maximal one such that $Y>Y'$ and
$C_{Y,Y'}(q)=\sum_{i=-r}^ra_iq^{i}$, then
$\zeta_{Y,Y'}(q)=\sum_{i=1}^{r}a_{-i}(q^i +q^{-i}) + a_0$.

\item[$(G'.2)$] if the coefficient of $Y'$ in $C(Y)-\sum_{Y>Z>Y'}
\zeta_{Y,Z}(q)G(Z)$ is given by $\sum_{i=-r}^ra_iq^{i}$, then
$\zeta_{Y,Y'}(q)=\sum_{i=1}^{r}a_{-i}(q^i +q^{-i}) + a_0$.
\end{itemize}

To summarize, we obtain the {\it modified generalized LLT
algorithm}:

\begin{cor}
Let $Y$ be a reduced proper Young wall in ${\mathcal Y}(\Lambda)$.
Then the corresponding global basis element $G(Y)$ can be
constructed recursively using the algorithm given in \eqref{GCY},
{\rm ($G'.1$)} and {\rm ($G'.2$)}.\qed
\end{cor}

In the following, we illustrate several examples.
\begin{ex}{\rm Suppose that $\frak{g}=A_4^{(2)}$. Note that $q_0=q$, $q_1=q^2$, and
$q_2=q^4$.

(a) Let $Y$ be the one of the following reduced proper Young
walls:\vskip 5mm

\begin{center}
\raisebox{-0.4\height}{
\begin{texdraw}
\drawdim em \setunitscale 0.13 \linewd 0.5

\move(30 0)\lvec(40 0)\lvec(40 10)\lvec(30 10)\lvec(30 0)\htext(33
1){\tiny $0$}

\move(40 0)\lvec(50 0)\lvec(50 10)\lvec(40 10)\lvec(40 0)\htext(43
1){\tiny $0$}

\move(30 0)\lvec(40 0)\lvec(40 5)\lvec(30 5)\lvec(30 0)\lfill
f:0.8 \htext(33 6){\tiny $0$}

\move(40 0)\lvec(50 0)\lvec(50 5)\lvec(40 5)\lvec(40 0)\lfill
f:0.8 \htext(43 6){\tiny $0$}

\move(30 10)\lvec(40 10)\lvec(40 20)\lvec(30 20)\lvec(30
10)\htext(33 13){$_1$}

\move(40 10)\lvec(50 10)\lvec(50 20)\lvec(40 20)\lvec(40
10)\htext(43 13){$_1$}
\move(40 20)\lvec(50 20)\lvec(50 30)\lvec(40 30)\lvec(40
20)\htext(43 23){$_2$}
\move(40 30)\lvec(50 30)\lvec(50 40)\lvec(40 40)\lvec(40
30)\htext(43 33){$_1$}
\move(40 40)\lvec(50 40)\lvec(50 45)\lvec(40 45)\lvec(40
40)\htext(43 41){\tiny $0$}
\end{texdraw}} \hskip 1cm
\raisebox{-0.4\height}{
\begin{texdraw}
\drawdim em \setunitscale 0.13 \linewd 0.5

\move(30 0)\lvec(40 0)\lvec(40 10)\lvec(30 10)\lvec(30 0)\htext(33
1){\tiny $0$}

\move(40 0)\lvec(50 0)\lvec(50 10)\lvec(40 10)\lvec(40 0)\htext(43
1){\tiny $0$}

\move(30 0)\lvec(40 0)\lvec(40 5)\lvec(30 5)\lvec(30 0)\lfill
f:0.8 \htext(33 6){\tiny $0$}

\move(40 0)\lvec(50 0)\lvec(50 5)\lvec(40 5)\lvec(40 0)\lfill
f:0.8 \htext(43 6){\tiny $0$}

\move(30 10)\lvec(40 10)\lvec(40 20)\lvec(30 20)\lvec(30
10)\htext(33 13){$_1$}

\move(40 10)\lvec(50 10)\lvec(50 20)\lvec(40 20)\lvec(40
10)\htext(43 13){$_1$}
\move(30 20)\lvec(40 20)\lvec(40 30)\lvec(30 30)\lvec(30
20)\htext(33 23){$_2$}

\move(40 20)\lvec(50 20)\lvec(50 30)\lvec(40 30)\lvec(40
20)\htext(43 23){$_2$}
\move(30 30)\lvec(40 30)\lvec(40 40)\lvec(30 40)\lvec(30
30)\htext(33 33){$_1$}

\move(40 30)\lvec(50 30)\lvec(50 40)\lvec(40 40)\lvec(40
30)\htext(43 33){$_1$}
\move(40 40)\lvec(50 40)\lvec(50 45)\lvec(40 45)\lvec(40
40)\htext(43 41){\tiny $0$}
\end{texdraw}}
\hskip 1cm \raisebox{-0.4\height}{
\begin{texdraw}
\drawdim em \setunitscale 0.13 \linewd 0.5

\move(20 0)\lvec(30 0)\lvec(30 10)\lvec(20 10)\lvec(20 0)\htext(23
1){\tiny $0$}

\move(30 0)\lvec(40 0)\lvec(40 10)\lvec(30 10)\lvec(30 0)\htext(33
1){\tiny $0$}

\move(40 0)\lvec(50 0)\lvec(50 10)\lvec(40 10)\lvec(40 0)\htext(43
1){\tiny $0$}

\move(20 0)\lvec(30 0)\lvec(30 5)\lvec(20 5)\lvec(20 0)\lfill
f:0.8 \htext(23 6){\tiny $0$}

\move(30 0)\lvec(40 0)\lvec(40 5)\lvec(30 5)\lvec(30 0)\lfill
f:0.8 \htext(33 6){\tiny $0$}

\move(40 0)\lvec(50 0)\lvec(50 5)\lvec(40 5)\lvec(40 0)\lfill
f:0.8 \htext(43 6){\tiny $0$}

\move(30 10)\lvec(40 10)\lvec(40 20)\lvec(30 20)\lvec(30
10)\htext(33 13){$_1$}

\move(40 10)\lvec(50 10)\lvec(50 20)\lvec(40 20)\lvec(40
10)\htext(43 13){$_1$}
\move(40 20)\lvec(50 20)\lvec(50 30)\lvec(40 30)\lvec(40
20)\htext(43 23){$_2$}
\move(40 30)\lvec(50 30)\lvec(50 40)\lvec(40 40)\lvec(40
30)\htext(43 33){$_1$}
\move(40 40)\lvec(50 40)\lvec(50 45)\lvec(40 45)\lvec(40
40)\htext(43 41){\tiny $0$}
\move(40 40)\lvec(50 40)\lvec(50 50)\lvec(40 50)\lvec(40
40)\htext(43 46){\tiny $0$}
\end{texdraw}}\quad .
\end{center}
Then $Y$ is the smallest one among the reduced proper Young walls
with the same weight with respect to the total ordering.
Therefore, we have $A(Y)=G(Y)$ by \eqref{H'YZ} and \vskip 5mm

$G(\ \raisebox{-0.4\height}{
\begin{texdraw}
\drawdim em \setunitscale 0.13 \linewd 0.5

\move(30 0)\lvec(40 0)\lvec(40 10)\lvec(30 10)\lvec(30 0)\htext(33
1){\tiny $0$}

\move(40 0)\lvec(50 0)\lvec(50 10)\lvec(40 10)\lvec(40 0)\htext(43
1){\tiny $0$}

\move(30 0)\lvec(40 0)\lvec(40 5)\lvec(30 5)\lvec(30 0)\lfill
f:0.8 \htext(33 6){\tiny $0$}

\move(40 0)\lvec(50 0)\lvec(50 5)\lvec(40 5)\lvec(40 0)\lfill
f:0.8 \htext(43 6){\tiny $0$}

\move(30 10)\lvec(40 10)\lvec(40 20)\lvec(30 20)\lvec(30
10)\htext(33 13){$_1$}

\move(40 10)\lvec(50 10)\lvec(50 20)\lvec(40 20)\lvec(40
10)\htext(43 13){$_1$}
\move(40 20)\lvec(50 20)\lvec(50 30)\lvec(40 30)\lvec(40
20)\htext(43 23){$_2$}
\move(40 30)\lvec(50 30)\lvec(50 40)\lvec(40 40)\lvec(40
30)\htext(43 33){$_1$}
\move(40 40)\lvec(50 40)\lvec(50 45)\lvec(40 45)\lvec(40
40)\htext(43 41){\tiny $0$}
\end{texdraw}}\ )= f_1f_0^{(2)}f_1f_2f_1f_0 Y_{\Lambda}
=\raisebox{-0.4\height}{
\begin{texdraw}
\drawdim em \setunitscale 0.13 \linewd 0.5

\move(30 0)\lvec(40 0)\lvec(40 10)\lvec(30 10)\lvec(30 0)\htext(33
1){\tiny $0$}

\move(40 0)\lvec(50 0)\lvec(50 10)\lvec(40 10)\lvec(40 0)\htext(43
1){\tiny $0$}

\move(30 0)\lvec(40 0)\lvec(40 5)\lvec(30 5)\lvec(30 0)\lfill
f:0.8 \htext(33 6){\tiny $0$}

\move(40 0)\lvec(50 0)\lvec(50 5)\lvec(40 5)\lvec(40 0)\lfill
f:0.8 \htext(43 6){\tiny $0$}

\move(30 10)\lvec(40 10)\lvec(40 20)\lvec(30 20)\lvec(30
10)\htext(33 13){$_1$}

\move(40 10)\lvec(50 10)\lvec(50 20)\lvec(40 20)\lvec(40
10)\htext(43 13){$_1$}
\move(40 20)\lvec(50 20)\lvec(50 30)\lvec(40 30)\lvec(40
20)\htext(43 23){$_2$}
\move(40 30)\lvec(50 30)\lvec(50 40)\lvec(40 40)\lvec(40
30)\htext(43 33){$_1$}
\move(40 40)\lvec(50 40)\lvec(50 45)\lvec(40 45)\lvec(40
40)\htext(43 41){\tiny $0$}
\end{texdraw}}\quad + \quad q^2
\raisebox{-0.4\height}{
\begin{texdraw}
\drawdim em \setunitscale 0.13 \linewd 0.5

\move(40 0)\lvec(50 0)\lvec(50 10)\lvec(40 10)\lvec(40 0)\htext(43
1){\tiny $0$}

\move(40 45)\lvec(50 45)\lvec(50 50)\lvec(40 50)\lvec(40
45)\htext(43 46){\tiny $0$}

\move(40 0)\lvec(50 0)\lvec(50 5)\lvec(40 5)\lvec(40 0)\lfill
f:0.8 \htext(43 6){\tiny $0$}

\move(40 50)\lvec(50 50)\lvec(50 60)\lvec(40 60)\lvec(40
50)\htext(43 53){$_1$}

\move(40 10)\lvec(50 10)\lvec(50 20)\lvec(40 20)\lvec(40
10)\htext(43 13){$_1$}
\move(40 20)\lvec(50 20)\lvec(50 30)\lvec(40 30)\lvec(40
20)\htext(43 23){$_2$}
\move(40 30)\lvec(50 30)\lvec(50 40)\lvec(40 40)\lvec(40
30)\htext(43 33){$_1$}
\move(40 40)\lvec(50 40)\lvec(50 45)\lvec(40 45)\lvec(40
40)\htext(43 41){\tiny $0$}
\end{texdraw}}$\quad ,\vskip 5mm

$G(\ \raisebox{-0.4\height}{
\begin{texdraw}
\drawdim em \setunitscale 0.13 \linewd 0.5

\move(30 0)\lvec(40 0)\lvec(40 10)\lvec(30 10)\lvec(30 0)\htext(33
1){\tiny $0$}

\move(40 0)\lvec(50 0)\lvec(50 10)\lvec(40 10)\lvec(40 0)\htext(43
1){\tiny $0$}

\move(30 0)\lvec(40 0)\lvec(40 5)\lvec(30 5)\lvec(30 0)\lfill
f:0.8 \htext(33 6){\tiny $0$}

\move(40 0)\lvec(50 0)\lvec(50 5)\lvec(40 5)\lvec(40 0)\lfill
f:0.8 \htext(43 6){\tiny $0$}

\move(30 10)\lvec(40 10)\lvec(40 20)\lvec(30 20)\lvec(30
10)\htext(33 13){$_1$}

\move(40 10)\lvec(50 10)\lvec(50 20)\lvec(40 20)\lvec(40
10)\htext(43 13){$_1$}
\move(30 20)\lvec(40 20)\lvec(40 30)\lvec(30 30)\lvec(30
20)\htext(33 23){$_2$}

\move(40 20)\lvec(50 20)\lvec(50 30)\lvec(40 30)\lvec(40
20)\htext(43 23){$_2$}
\move(30 30)\lvec(40 30)\lvec(40 40)\lvec(30 40)\lvec(30
30)\htext(33 33){$_1$}

\move(40 30)\lvec(50 30)\lvec(50 40)\lvec(40 40)\lvec(40
30)\htext(43 33){$_1$}
\move(40 40)\lvec(50 40)\lvec(50 45)\lvec(40 45)\lvec(40
40)\htext(43 41){\tiny $0$}
\end{texdraw}}\ )= f_1f_2f_1f_0^{(2)}f_1f_2f_1f_0 Y_{\Lambda}
=\raisebox{-0.4\height}{
\begin{texdraw}
\drawdim em \setunitscale 0.13 \linewd 0.5

\move(30 0)\lvec(40 0)\lvec(40 10)\lvec(30 10)\lvec(30 0)\htext(33
1){\tiny $0$}

\move(40 0)\lvec(50 0)\lvec(50 10)\lvec(40 10)\lvec(40 0)\htext(43
1){\tiny $0$}

\move(30 0)\lvec(40 0)\lvec(40 5)\lvec(30 5)\lvec(30 0)\lfill
f:0.8 \htext(33 6){\tiny $0$}

\move(40 0)\lvec(50 0)\lvec(50 5)\lvec(40 5)\lvec(40 0)\lfill
f:0.8 \htext(43 6){\tiny $0$}

\move(30 10)\lvec(40 10)\lvec(40 20)\lvec(30 20)\lvec(30
10)\htext(33 13){$_1$}

\move(40 10)\lvec(50 10)\lvec(50 20)\lvec(40 20)\lvec(40
10)\htext(43 13){$_1$}
\move(30 20)\lvec(40 20)\lvec(40 30)\lvec(30 30)\lvec(30
20)\htext(33 23){$_2$}

\move(40 20)\lvec(50 20)\lvec(50 30)\lvec(40 30)\lvec(40
20)\htext(43 23){$_2$}
\move(30 30)\lvec(40 30)\lvec(40 40)\lvec(30 40)\lvec(30
30)\htext(33 33){$_1$}

\move(40 30)\lvec(50 30)\lvec(50 40)\lvec(40 40)\lvec(40
30)\htext(43 33){$_1$}
\move(40 40)\lvec(50 40)\lvec(50 45)\lvec(40 45)\lvec(40
40)\htext(43 41){\tiny $0$}
\end{texdraw}}\quad + \quad q^2
\raisebox{-0.4\height}{
\begin{texdraw}
\drawdim em \setunitscale 0.13 \linewd 0.5

\move(40 40)\lvec(50 40)\lvec(50 50)\lvec(40 50)\lvec(40
40)\htext(43 46){\tiny $0$}

\move(40 0)\lvec(50 0)\lvec(50 10)\lvec(40 10)\lvec(40 0)\htext(43
1){\tiny $0$}

\move(40 0)\lvec(50 0)\lvec(50 5)\lvec(40 5)\lvec(40 0)\lfill
f:0.8 \htext(43 6){\tiny $0$}

\move(40 50)\lvec(50 50)\lvec(50 60)\lvec(40 60)\lvec(40
50)\htext(43 53){$_1$}

\move(40 10)\lvec(50 10)\lvec(50 20)\lvec(40 20)\lvec(40
10)\htext(43 13){$_1$}
\move(40 60)\lvec(50 60)\lvec(50 70)\lvec(40 70)\lvec(40
60)\htext(43 63){$_2$}

\move(40 20)\lvec(50 20)\lvec(50 30)\lvec(40 30)\lvec(40
20)\htext(43 23){$_2$}
\move(40 70)\lvec(50 70)\lvec(50 80)\lvec(40 80)\lvec(40
70)\htext(43 73){$_1$}

\move(40 30)\lvec(50 30)\lvec(50 40)\lvec(40 40)\lvec(40
30)\htext(43 33){$_1$}
\move(40 40)\lvec(50 40)\lvec(50 45)\lvec(40 45)\lvec(40
40)\htext(43 41){\tiny $0$}
\end{texdraw}}$\quad ,\vskip 5mm

$G(\ \raisebox{-0.4\height}{
\begin{texdraw}
\drawdim em \setunitscale 0.13 \linewd 0.5

\move(20 0)\lvec(30 0)\lvec(30 10)\lvec(20 10)\lvec(20 0)\htext(23
1){\tiny $0$}

\move(30 0)\lvec(40 0)\lvec(40 10)\lvec(30 10)\lvec(30 0)\htext(33
1){\tiny $0$}

\move(40 0)\lvec(50 0)\lvec(50 10)\lvec(40 10)\lvec(40 0)\htext(43
1){\tiny $0$}

\move(20 0)\lvec(30 0)\lvec(30 5)\lvec(20 5)\lvec(20 0)\lfill
f:0.8 \htext(23 6){\tiny $0$}

\move(30 0)\lvec(40 0)\lvec(40 5)\lvec(30 5)\lvec(30 0)\lfill
f:0.8 \htext(33 6){\tiny $0$}

\move(40 0)\lvec(50 0)\lvec(50 5)\lvec(40 5)\lvec(40 0)\lfill
f:0.8 \htext(43 6){\tiny $0$}

\move(30 10)\lvec(40 10)\lvec(40 20)\lvec(30 20)\lvec(30
10)\htext(33 13){$_1$}

\move(40 10)\lvec(50 10)\lvec(50 20)\lvec(40 20)\lvec(40
10)\htext(43 13){$_1$}
\move(40 20)\lvec(50 20)\lvec(50 30)\lvec(40 30)\lvec(40
20)\htext(43 23){$_2$}
\move(40 30)\lvec(50 30)\lvec(50 40)\lvec(40 40)\lvec(40
30)\htext(43 33){$_1$}
\move(40 40)\lvec(50 40)\lvec(50 45)\lvec(40 45)\lvec(40
40)\htext(43 41){\tiny $0$}
\move(40 40)\lvec(50 40)\lvec(50 50)\lvec(40 50)\lvec(40
40)\htext(43 46){\tiny $0$}
\end{texdraw}} \ )=f_0f_1f_0^{(3)}f_1f_2f_1f_0
Y_{\Lambda}= \raisebox{-0.4\height}{
\begin{texdraw}
\drawdim em \setunitscale 0.13 \linewd 0.5

\move(20 0)\lvec(30 0)\lvec(30 10)\lvec(20 10)\lvec(20 0)\htext(23
1){\tiny $0$}

\move(30 0)\lvec(40 0)\lvec(40 10)\lvec(30 10)\lvec(30 0)\htext(33
1){\tiny $0$}

\move(40 0)\lvec(50 0)\lvec(50 10)\lvec(40 10)\lvec(40 0)\htext(43
1){\tiny $0$}

\move(20 0)\lvec(30 0)\lvec(30 5)\lvec(20 5)\lvec(20 0)\lfill
f:0.8 \htext(23 6){\tiny $0$}

\move(30 0)\lvec(40 0)\lvec(40 5)\lvec(30 5)\lvec(30 0)\lfill
f:0.8 \htext(33 6){\tiny $0$}

\move(40 0)\lvec(50 0)\lvec(50 5)\lvec(40 5)\lvec(40 0)\lfill
f:0.8 \htext(43 6){\tiny $0$}

\move(30 10)\lvec(40 10)\lvec(40 20)\lvec(30 20)\lvec(30
10)\htext(33 13){$_1$}

\move(40 10)\lvec(50 10)\lvec(50 20)\lvec(40 20)\lvec(40
10)\htext(43 13){$_1$}
\move(40 20)\lvec(50 20)\lvec(50 30)\lvec(40 30)\lvec(40
20)\htext(43 23){$_2$}
\move(40 30)\lvec(50 30)\lvec(50 40)\lvec(40 40)\lvec(40
30)\htext(43 33){$_1$}
\move(40 40)\lvec(50 40)\lvec(50 45)\lvec(40 45)\lvec(40
40)\htext(43 41){\tiny $0$}
\move(40 40)\lvec(50 40)\lvec(50 50)\lvec(40 50)\lvec(40
40)\htext(43 46){\tiny $0$}
\end{texdraw}}$\quad .\vskip 10mm

(b) Observe that  \vskip 5mm

$A(\ \raisebox{-0.4\height}{
\begin{texdraw}
\drawdim em \setunitscale 0.13 \linewd 0.5

\move(20 0)\lvec(30 0)\lvec(30 10)\lvec(20 10)\lvec(20 0)\htext(23
1){\tiny $0$}

\move(30 0)\lvec(40 0)\lvec(40 10)\lvec(30 10)\lvec(30 0)\htext(33
1){\tiny $0$}

\move(40 0)\lvec(50 0)\lvec(50 10)\lvec(40 10)\lvec(40 0)\htext(43
1){\tiny $0$}

\move(20 0)\lvec(30 0)\lvec(30 5)\lvec(20 5)\lvec(20 0)\lfill
f:0.8 \htext(23 6){\tiny $0$}

\move(30 0)\lvec(40 0)\lvec(40 5)\lvec(30 5)\lvec(30 0)\lfill
f:0.8 \htext(33 6){\tiny $0$}

\move(40 0)\lvec(50 0)\lvec(50 5)\lvec(40 5)\lvec(40 0)\lfill
f:0.8 \htext(43 6){\tiny $0$}

\move(20 10)\lvec(30 10)\lvec(30 20)\lvec(20 20)\lvec(20
10)\htext(23 13){$_1$}

\move(30 10)\lvec(40 10)\lvec(40 20)\lvec(30 20)\lvec(30
10)\htext(33 13){$_1$}

\move(40 10)\lvec(50 10)\lvec(50 20)\lvec(40 20)\lvec(40
10)\htext(43 13){$_1$}
\move(30 20)\lvec(40 20)\lvec(40 30)\lvec(30 30)\lvec(30
20)\htext(33 23){$_2$}

\move(40 20)\lvec(50 20)\lvec(50 30)\lvec(40 30)\lvec(40
20)\htext(43 23){$_2$}
\move(30 30)\lvec(40 30)\lvec(40 40)\lvec(30 40)\lvec(30
30)\htext(33 33){$_1$}

\move(40 30)\lvec(50 30)\lvec(50 40)\lvec(40 40)\lvec(40
30)\htext(43 33){$_1$}
\move(40 40)\lvec(50 40)\lvec(50 45)\lvec(40 45)\lvec(40
40)\htext(43 41){\tiny $0$}
\move(40 40)\lvec(50 40)\lvec(50 50)\lvec(40 50)\lvec(40
40)\htext(43 46){\tiny $0$}
\end{texdraw}} \ )=f_1f_0f_1f_2f_1f_0^{(3)}f_1f_2f_1f_0Y_{\Lambda}$\vskip
5mm

\hskip 1cm $=$ \raisebox{-0.4\height}{
\begin{texdraw}
\drawdim em \setunitscale 0.13 \linewd 0.5

\move(20 0)\lvec(30 0)\lvec(30 10)\lvec(20 10)\lvec(20 0)\htext(23
1){\tiny $0$}

\move(30 0)\lvec(40 0)\lvec(40 10)\lvec(30 10)\lvec(30 0)\htext(33
1){\tiny $0$}

\move(40 0)\lvec(50 0)\lvec(50 10)\lvec(40 10)\lvec(40 0)\htext(43
1){\tiny $0$}

\move(20 0)\lvec(30 0)\lvec(30 5)\lvec(20 5)\lvec(20 0)\lfill
f:0.8 \htext(23 6){\tiny $0$}

\move(30 0)\lvec(40 0)\lvec(40 5)\lvec(30 5)\lvec(30 0)\lfill
f:0.8 \htext(33 6){\tiny $0$}

\move(40 0)\lvec(50 0)\lvec(50 5)\lvec(40 5)\lvec(40 0)\lfill
f:0.8 \htext(43 6){\tiny $0$}

\move(20 10)\lvec(30 10)\lvec(30 20)\lvec(20 20)\lvec(20
10)\htext(23 13){$_1$}

\move(30 10)\lvec(40 10)\lvec(40 20)\lvec(30 20)\lvec(30
10)\htext(33 13){$_1$}

\move(40 10)\lvec(50 10)\lvec(50 20)\lvec(40 20)\lvec(40
10)\htext(43 13){$_1$}
\move(30 20)\lvec(40 20)\lvec(40 30)\lvec(30 30)\lvec(30
20)\htext(33 23){$_2$}

\move(40 20)\lvec(50 20)\lvec(50 30)\lvec(40 30)\lvec(40
20)\htext(43 23){$_2$}
\move(30 30)\lvec(40 30)\lvec(40 40)\lvec(30 40)\lvec(30
30)\htext(33 33){$_1$}

\move(40 30)\lvec(50 30)\lvec(50 40)\lvec(40 40)\lvec(40
30)\htext(43 33){$_1$}
\move(40 40)\lvec(50 40)\lvec(50 45)\lvec(40 45)\lvec(40
40)\htext(43 41){\tiny $0$}
\move(40 40)\lvec(50 40)\lvec(50 50)\lvec(40 50)\lvec(40
40)\htext(43 46){\tiny $0$}
\end{texdraw}}
$+$\quad $q^2$ \raisebox{-0.4\height}{
\begin{texdraw}
\drawdim em \setunitscale 0.13 \linewd 0.5

\move(20 0)\lvec(30 0)\lvec(30 10)\lvec(20 10)\lvec(20 0)\htext(23
1){\tiny $0$}

\move(30 0)\lvec(40 0)\lvec(40 10)\lvec(30 10)\lvec(30 0)\htext(33
1){\tiny $0$}

\move(40 0)\lvec(50 0)\lvec(50 10)\lvec(40 10)\lvec(40 0)\htext(43
1){\tiny $0$}

\move(20 0)\lvec(30 0)\lvec(30 5)\lvec(20 5)\lvec(20 0)\lfill
f:0.8 \htext(23 6){\tiny $0$}

\move(30 0)\lvec(40 0)\lvec(40 5)\lvec(30 5)\lvec(30 0)\lfill
f:0.8 \htext(33 6){\tiny $0$}

\move(40 0)\lvec(50 0)\lvec(50 5)\lvec(40 5)\lvec(40 0)\lfill
f:0.8 \htext(43 6){\tiny $0$}

\move(20 10)\lvec(30 10)\lvec(30 20)\lvec(20 20)\lvec(20
10)\htext(23 13){$_1$}

\move(30 10)\lvec(40 10)\lvec(40 20)\lvec(30 20)\lvec(30
10)\htext(33 13){$_1$}

\move(40 10)\lvec(50 10)\lvec(50 20)\lvec(40 20)\lvec(40
10)\htext(43 13){$_1$}
\move(30 20)\lvec(40 20)\lvec(40 30)\lvec(30 30)\lvec(30
20)\htext(33 23){$_2$}

\move(40 20)\lvec(50 20)\lvec(50 30)\lvec(40 30)\lvec(40
20)\htext(43 23){$_2$}
\move(40 50)\lvec(50 50)\lvec(50 60)\lvec(40 60)\lvec(40
50)\htext(43 53){$_1$}

\move(40 30)\lvec(50 30)\lvec(50 40)\lvec(40 40)\lvec(40
30)\htext(43 33){$_1$}
\move(40 40)\lvec(50 40)\lvec(50 45)\lvec(40 45)\lvec(40
40)\htext(43 41){\tiny $0$}
\move(40 40)\lvec(50 40)\lvec(50 50)\lvec(40 50)\lvec(40
40)\htext(43 46){\tiny $0$}
\end{texdraw}}\quad
$+$\quad $(1+q^4)$ \raisebox{-0.4\height}{
\begin{texdraw}
\drawdim em \setunitscale 0.13 \linewd 0.5

\move(20 0)\lvec(30 0)\lvec(30 10)\lvec(20 10)\lvec(20 0)\htext(23
1){\tiny $0$}

\move(30 0)\lvec(40 0)\lvec(40 10)\lvec(30 10)\lvec(30 0)\htext(33
1){\tiny $0$}

\move(40 0)\lvec(50 0)\lvec(50 10)\lvec(40 10)\lvec(40 0)\htext(43
1){\tiny $0$}

\move(20 0)\lvec(30 0)\lvec(30 5)\lvec(20 5)\lvec(20 0)\lfill
f:0.8 \htext(23 6){\tiny $0$}

\move(30 0)\lvec(40 0)\lvec(40 5)\lvec(30 5)\lvec(30 0)\lfill
f:0.8 \htext(33 6){\tiny $0$}

\move(40 0)\lvec(50 0)\lvec(50 5)\lvec(40 5)\lvec(40 0)\lfill
f:0.8 \htext(43 6){\tiny $0$}

\move(40 50)\lvec(50 50)\lvec(50 60)\lvec(40 60)\lvec(40
50)\htext(43 53){$_1$}

\move(30 10)\lvec(40 10)\lvec(40 20)\lvec(30 20)\lvec(30
10)\htext(33 13){$_1$}

\move(40 10)\lvec(50 10)\lvec(50 20)\lvec(40 20)\lvec(40
10)\htext(43 13){$_1$}
\move(30 20)\lvec(40 20)\lvec(40 30)\lvec(30 30)\lvec(30
20)\htext(33 23){$_2$}

\move(40 20)\lvec(50 20)\lvec(50 30)\lvec(40 30)\lvec(40
20)\htext(43 23){$_2$}
\move(30 30)\lvec(40 30)\lvec(40 40)\lvec(30 40)\lvec(30
30)\htext(33 33){$_1$}

\move(40 30)\lvec(50 30)\lvec(50 40)\lvec(40 40)\lvec(40
30)\htext(43 33){$_1$}
\move(40 40)\lvec(50 40)\lvec(50 45)\lvec(40 45)\lvec(40
40)\htext(43 41){\tiny $0$}
\move(40 40)\lvec(50 40)\lvec(50 50)\lvec(40 50)\lvec(40
40)\htext(43 46){\tiny $0$}
\end{texdraw}}\vskip 5mm

\hskip 1cm $+$\quad $q^2$ \raisebox{-0.4\height}{
\begin{texdraw}
\drawdim em \setunitscale 0.13 \linewd 0.5

\move(20 0)\lvec(30 0)\lvec(30 10)\lvec(20 10)\lvec(20 0)\htext(23
1){\tiny $0$}

\move(30 0)\lvec(40 0)\lvec(40 10)\lvec(30 10)\lvec(30 0)\htext(33
1){\tiny $0$}

\move(40 0)\lvec(50 0)\lvec(50 10)\lvec(40 10)\lvec(40 0)\htext(43
1){\tiny $0$}

\move(20 0)\lvec(30 0)\lvec(30 5)\lvec(20 5)\lvec(20 0)\lfill
f:0.8 \htext(23 6){\tiny $0$}

\move(30 0)\lvec(40 0)\lvec(40 5)\lvec(30 5)\lvec(30 0)\lfill
f:0.8 \htext(33 6){\tiny $0$}

\move(40 0)\lvec(50 0)\lvec(50 5)\lvec(40 5)\lvec(40 0)\lfill
f:0.8 \htext(43 6){\tiny $0$}

\move(40 50)\lvec(50 50)\lvec(50 60)\lvec(40 60)\lvec(40
50)\htext(43 53){$_1$}

\move(30 10)\lvec(40 10)\lvec(40 20)\lvec(30 20)\lvec(30
10)\htext(33 13){$_1$}

\move(40 10)\lvec(50 10)\lvec(50 20)\lvec(40 20)\lvec(40
10)\htext(43 13){$_1$}
\move(40 60)\lvec(50 60)\lvec(50 70)\lvec(40 70)\lvec(40
60)\htext(43 63){$_2$}

\move(40 20)\lvec(50 20)\lvec(50 30)\lvec(40 30)\lvec(40
20)\htext(43 23){$_2$}
\move(40 70)\lvec(50 70)\lvec(50 80)\lvec(40 80)\lvec(40
70)\htext(43 73){$_1$}

\move(40 30)\lvec(50 30)\lvec(50 40)\lvec(40 40)\lvec(40
30)\htext(43 33){$_1$}
\move(40 40)\lvec(50 40)\lvec(50 45)\lvec(40 45)\lvec(40
40)\htext(43 41){\tiny $0$}
\move(40 40)\lvec(50 40)\lvec(50 50)\lvec(40 50)\lvec(40
40)\htext(43 46){\tiny $0$}
\end{texdraw}}\quad $+$\quad $q$
\raisebox{-0.4\height}{
\begin{texdraw}
\drawdim em \setunitscale 0.13 \linewd 0.5

\move(30 0)\lvec(40 0)\lvec(40 10)\lvec(30 10)\lvec(30 0)\htext(33
1){\tiny $0$}

\move(40 0)\lvec(50 0)\lvec(50 10)\lvec(40 10)\lvec(40 0)\htext(43
1){\tiny $0$}

\move(30 40)\lvec(40 40)\lvec(40 45)\lvec(30 45)\lvec(30
40)\htext(33 41){\tiny $0$}

\move(30 0)\lvec(40 0)\lvec(40 5)\lvec(30 5)\lvec(30 0)\lfill
f:0.8 \htext(33 6){\tiny $0$}

\move(40 0)\lvec(50 0)\lvec(50 5)\lvec(40 5)\lvec(40 0)\lfill
f:0.8 \htext(43 6){\tiny $0$}

\move(40 50)\lvec(50 50)\lvec(50 60)\lvec(40 60)\lvec(40
50)\htext(43 53){$_1$}

\move(30 10)\lvec(40 10)\lvec(40 20)\lvec(30 20)\lvec(30
10)\htext(33 13){$_1$}

\move(40 10)\lvec(50 10)\lvec(50 20)\lvec(40 20)\lvec(40
10)\htext(43 13){$_1$}
\move(30 20)\lvec(40 20)\lvec(40 30)\lvec(30 30)\lvec(30
20)\htext(33 23){$_2$}

\move(40 20)\lvec(50 20)\lvec(50 30)\lvec(40 30)\lvec(40
20)\htext(43 23){$_2$}
\move(30 30)\lvec(40 30)\lvec(40 40)\lvec(30 40)\lvec(30
30)\htext(33 33){$_1$}

\move(40 30)\lvec(50 30)\lvec(50 40)\lvec(40 40)\lvec(40
30)\htext(43 33){$_1$}
\move(40 40)\lvec(50 40)\lvec(50 45)\lvec(40 45)\lvec(40
40)\htext(43 41){\tiny $0$}
\move(40 40)\lvec(50 40)\lvec(50 50)\lvec(40 50)\lvec(40
40)\htext(43 46){\tiny $0$}
\end{texdraw}}\quad $+$\quad $q^3$
\raisebox{-0.4\height}{
\begin{texdraw}
\drawdim em \setunitscale 0.13 \linewd 0.5

\move(40 80)\lvec(50 80)\lvec(50 85)\lvec(40 85)\lvec(40
80)\htext(43 81){\tiny $0$}

\move(30 0)\lvec(40 0)\lvec(40 10)\lvec(30 10)\lvec(30 0)\htext(33
1){\tiny $0$}

\move(40 0)\lvec(50 0)\lvec(50 10)\lvec(40 10)\lvec(40 0)\htext(43
1){\tiny $0$}

\move(30 0)\lvec(40 0)\lvec(40 5)\lvec(30 5)\lvec(30 0)\lfill
f:0.8 \htext(33 6){\tiny $0$}

\move(40 0)\lvec(50 0)\lvec(50 5)\lvec(40 5)\lvec(40 0)\lfill
f:0.8 \htext(43 6){\tiny $0$}

\move(40 50)\lvec(50 50)\lvec(50 60)\lvec(40 60)\lvec(40
50)\htext(43 53){$_1$}

\move(30 10)\lvec(40 10)\lvec(40 20)\lvec(30 20)\lvec(30
10)\htext(33 13){$_1$}

\move(40 10)\lvec(50 10)\lvec(50 20)\lvec(40 20)\lvec(40
10)\htext(43 13){$_1$}
\move(40 60)\lvec(50 60)\lvec(50 70)\lvec(40 70)\lvec(40
60)\htext(43 63){$_2$}

\move(40 20)\lvec(50 20)\lvec(50 30)\lvec(40 30)\lvec(40
20)\htext(43 23){$_2$}
\move(40 70)\lvec(50 70)\lvec(50 80)\lvec(40 80)\lvec(40
70)\htext(43 73){$_1$}

\move(40 30)\lvec(50 30)\lvec(50 40)\lvec(40 40)\lvec(40
30)\htext(43 33){$_1$}
\move(40 40)\lvec(50 40)\lvec(50 45)\lvec(40 45)\lvec(40
40)\htext(43 41){\tiny $0$}
\move(40 40)\lvec(50 40)\lvec(50 50)\lvec(40 50)\lvec(40
40)\htext(43 46){\tiny $0$}
\end{texdraw}}\quad .\vskip 5mm

On the other hand, \vskip 5mm

$G(\ \raisebox{-0.4\height}{
\begin{texdraw}
\drawdim em \setunitscale 0.13 \linewd 0.5

\move(20 0)\lvec(30 0)\lvec(30 10)\lvec(20 10)\lvec(20 0)\htext(23
1){\tiny $0$}

\move(30 0)\lvec(40 0)\lvec(40 10)\lvec(30 10)\lvec(30 0)\htext(33
1){\tiny $0$}

\move(40 0)\lvec(50 0)\lvec(50 10)\lvec(40 10)\lvec(40 0)\htext(43
1){\tiny $0$}

\move(20 0)\lvec(30 0)\lvec(30 5)\lvec(20 5)\lvec(20 0)\lfill
f:0.8 \htext(23 6){\tiny $0$}

\move(30 0)\lvec(40 0)\lvec(40 5)\lvec(30 5)\lvec(30 0)\lfill
f:0.8 \htext(33 6){\tiny $0$}

\move(40 0)\lvec(50 0)\lvec(50 5)\lvec(40 5)\lvec(40 0)\lfill
f:0.8 \htext(43 6){\tiny $0$}

\move(40 50)\lvec(50 50)\lvec(50 60)\lvec(40 60)\lvec(40
50)\htext(43 53){$_1$}

\move(30 10)\lvec(40 10)\lvec(40 20)\lvec(30 20)\lvec(30
10)\htext(33 13){$_1$}

\move(40 10)\lvec(50 10)\lvec(50 20)\lvec(40 20)\lvec(40
10)\htext(43 13){$_1$}
\move(30 20)\lvec(40 20)\lvec(40 30)\lvec(30 30)\lvec(30
20)\htext(33 23){$_2$}

\move(40 20)\lvec(50 20)\lvec(50 30)\lvec(40 30)\lvec(40
20)\htext(43 23){$_2$}
\move(30 30)\lvec(40 30)\lvec(40 40)\lvec(30 40)\lvec(30
30)\htext(33 33){$_1$}

\move(40 30)\lvec(50 30)\lvec(50 40)\lvec(40 40)\lvec(40
30)\htext(43 33){$_1$}
\move(40 40)\lvec(50 40)\lvec(50 45)\lvec(40 45)\lvec(40
40)\htext(43 41){\tiny $0$}
\move(40 40)\lvec(50 40)\lvec(50 50)\lvec(40 50)\lvec(40
40)\htext(43 46){\tiny $0$}
\end{texdraw}} \ )$
$=A(\ \raisebox{-0.4\height}{
\begin{texdraw}
\drawdim em \setunitscale 0.13 \linewd 0.5

\move(20 0)\lvec(30 0)\lvec(30 10)\lvec(20 10)\lvec(20 0)\htext(23
1){\tiny $0$}

\move(30 0)\lvec(40 0)\lvec(40 10)\lvec(30 10)\lvec(30 0)\htext(33
1){\tiny $0$}

\move(40 0)\lvec(50 0)\lvec(50 10)\lvec(40 10)\lvec(40 0)\htext(43
1){\tiny $0$}

\move(20 0)\lvec(30 0)\lvec(30 5)\lvec(20 5)\lvec(20 0)\lfill
f:0.8 \htext(23 6){\tiny $0$}

\move(30 0)\lvec(40 0)\lvec(40 5)\lvec(30 5)\lvec(30 0)\lfill
f:0.8 \htext(33 6){\tiny $0$}

\move(40 0)\lvec(50 0)\lvec(50 5)\lvec(40 5)\lvec(40 0)\lfill
f:0.8 \htext(43 6){\tiny $0$}

\move(40 50)\lvec(50 50)\lvec(50 60)\lvec(40 60)\lvec(40
50)\htext(43 53){$_1$}

\move(30 10)\lvec(40 10)\lvec(40 20)\lvec(30 20)\lvec(30
10)\htext(33 13){$_1$}

\move(40 10)\lvec(50 10)\lvec(50 20)\lvec(40 20)\lvec(40
10)\htext(43 13){$_1$}
\move(30 20)\lvec(40 20)\lvec(40 30)\lvec(30 30)\lvec(30
20)\htext(33 23){$_2$}

\move(40 20)\lvec(50 20)\lvec(50 30)\lvec(40 30)\lvec(40
20)\htext(43 23){$_2$}
\move(30 30)\lvec(40 30)\lvec(40 40)\lvec(30 40)\lvec(30
30)\htext(33 33){$_1$}

\move(40 30)\lvec(50 30)\lvec(50 40)\lvec(40 40)\lvec(40
30)\htext(43 33){$_1$}
\move(40 40)\lvec(50 40)\lvec(50 45)\lvec(40 45)\lvec(40
40)\htext(43 41){\tiny $0$}
\move(40 40)\lvec(50 40)\lvec(50 50)\lvec(40 50)\lvec(40
40)\htext(43 46){\tiny $0$}
\end{texdraw}} \ )$
$=f_0f_1f_2f_1^{(2)}f_0^{(3)}f_1f_2f_1f_0Y_{\Lambda}$ \vskip 5mm

$=$ \raisebox{-0.4\height}{
\begin{texdraw}
\drawdim em \setunitscale 0.13 \linewd 0.5

\move(20 0)\lvec(30 0)\lvec(30 10)\lvec(20 10)\lvec(20 0)\htext(23
1){\tiny $0$}

\move(30 0)\lvec(40 0)\lvec(40 10)\lvec(30 10)\lvec(30 0)\htext(33
1){\tiny $0$}

\move(40 0)\lvec(50 0)\lvec(50 10)\lvec(40 10)\lvec(40 0)\htext(43
1){\tiny $0$}

\move(20 0)\lvec(30 0)\lvec(30 5)\lvec(20 5)\lvec(20 0)\lfill
f:0.8 \htext(23 6){\tiny $0$}

\move(30 0)\lvec(40 0)\lvec(40 5)\lvec(30 5)\lvec(30 0)\lfill
f:0.8 \htext(33 6){\tiny $0$}

\move(40 0)\lvec(50 0)\lvec(50 5)\lvec(40 5)\lvec(40 0)\lfill
f:0.8 \htext(43 6){\tiny $0$}

\move(40 50)\lvec(50 50)\lvec(50 60)\lvec(40 60)\lvec(40
50)\htext(43 53){$_1$}

\move(30 10)\lvec(40 10)\lvec(40 20)\lvec(30 20)\lvec(30
10)\htext(33 13){$_1$}

\move(40 10)\lvec(50 10)\lvec(50 20)\lvec(40 20)\lvec(40
10)\htext(43 13){$_1$}
\move(30 20)\lvec(40 20)\lvec(40 30)\lvec(30 30)\lvec(30
20)\htext(33 23){$_2$}

\move(40 20)\lvec(50 20)\lvec(50 30)\lvec(40 30)\lvec(40
20)\htext(43 23){$_2$}
\move(30 30)\lvec(40 30)\lvec(40 40)\lvec(30 40)\lvec(30
30)\htext(33 33){$_1$}

\move(40 30)\lvec(50 30)\lvec(50 40)\lvec(40 40)\lvec(40
30)\htext(43 33){$_1$}
\move(40 40)\lvec(50 40)\lvec(50 45)\lvec(40 45)\lvec(40
40)\htext(43 41){\tiny $0$}
\move(40 40)\lvec(50 40)\lvec(50 50)\lvec(40 50)\lvec(40
40)\htext(43 46){\tiny $0$}
\end{texdraw}}\quad$+$\quad $q^2$
\raisebox{-0.4\height}{
\begin{texdraw}
\drawdim em \setunitscale 0.13 \linewd 0.5

\move(20 0)\lvec(30 0)\lvec(30 10)\lvec(20 10)\lvec(20 0)\htext(23
1){\tiny $0$}

\move(30 0)\lvec(40 0)\lvec(40 10)\lvec(30 10)\lvec(30 0)\htext(33
1){\tiny $0$}

\move(40 0)\lvec(50 0)\lvec(50 10)\lvec(40 10)\lvec(40 0)\htext(43
1){\tiny $0$}

\move(20 0)\lvec(30 0)\lvec(30 5)\lvec(20 5)\lvec(20 0)\lfill
f:0.8 \htext(23 6){\tiny $0$}

\move(30 0)\lvec(40 0)\lvec(40 5)\lvec(30 5)\lvec(30 0)\lfill
f:0.8 \htext(33 6){\tiny $0$}

\move(40 0)\lvec(50 0)\lvec(50 5)\lvec(40 5)\lvec(40 0)\lfill
f:0.8 \htext(43 6){\tiny $0$}

\move(40 50)\lvec(50 50)\lvec(50 60)\lvec(40 60)\lvec(40
50)\htext(43 53){$_1$}

\move(30 10)\lvec(40 10)\lvec(40 20)\lvec(30 20)\lvec(30
10)\htext(33 13){$_1$}

\move(40 10)\lvec(50 10)\lvec(50 20)\lvec(40 20)\lvec(40
10)\htext(43 13){$_1$}
\move(40 60)\lvec(50 60)\lvec(50 70)\lvec(40 70)\lvec(40
60)\htext(43 63){$_2$}

\move(40 20)\lvec(50 20)\lvec(50 30)\lvec(40 30)\lvec(40
20)\htext(43 23){$_2$}
\move(40 70)\lvec(50 70)\lvec(50 80)\lvec(40 80)\lvec(40
70)\htext(43 73){$_1$}

\move(40 30)\lvec(50 30)\lvec(50 40)\lvec(40 40)\lvec(40
30)\htext(43 33){$_1$}
\move(40 40)\lvec(50 40)\lvec(50 45)\lvec(40 45)\lvec(40
40)\htext(43 41){\tiny $0$}
\move(40 40)\lvec(50 40)\lvec(50 50)\lvec(40 50)\lvec(40
40)\htext(43 46){\tiny $0$}
\end{texdraw}}\quad $+$\quad $q$
\raisebox{-0.4\height}{
\begin{texdraw}
\drawdim em \setunitscale 0.13 \linewd 0.5

\move(30 40)\lvec(40 40)\lvec(40 45)\lvec(30 45)\lvec(30
40)\htext(33 41){\tiny $0$}

\move(30 0)\lvec(40 0)\lvec(40 10)\lvec(30 10)\lvec(30 0)\htext(33
1){\tiny $0$}

\move(40 0)\lvec(50 0)\lvec(50 10)\lvec(40 10)\lvec(40 0)\htext(43
1){\tiny $0$}

\move(30 0)\lvec(40 0)\lvec(40 5)\lvec(30 5)\lvec(30 0)\lfill
f:0.8 \htext(33 6){\tiny $0$}

\move(40 0)\lvec(50 0)\lvec(50 5)\lvec(40 5)\lvec(40 0)\lfill
f:0.8 \htext(43 6){\tiny $0$}

\move(40 50)\lvec(50 50)\lvec(50 60)\lvec(40 60)\lvec(40
50)\htext(43 53){$_1$}

\move(30 10)\lvec(40 10)\lvec(40 20)\lvec(30 20)\lvec(30
10)\htext(33 13){$_1$}

\move(40 10)\lvec(50 10)\lvec(50 20)\lvec(40 20)\lvec(40
10)\htext(43 13){$_1$}
\move(30 20)\lvec(40 20)\lvec(40 30)\lvec(30 30)\lvec(30
20)\htext(33 23){$_2$}

\move(40 20)\lvec(50 20)\lvec(50 30)\lvec(40 30)\lvec(40
20)\htext(43 23){$_2$}
\move(30 30)\lvec(40 30)\lvec(40 40)\lvec(30 40)\lvec(30
30)\htext(33 33){$_1$}

\move(40 30)\lvec(50 30)\lvec(50 40)\lvec(40 40)\lvec(40
30)\htext(43 33){$_1$}
\move(40 40)\lvec(50 40)\lvec(50 45)\lvec(40 45)\lvec(40
40)\htext(43 41){\tiny $0$}
\move(40 40)\lvec(50 40)\lvec(50 50)\lvec(40 50)\lvec(40
40)\htext(43 46){\tiny $0$}
\end{texdraw}}\quad $+$ \quad $q^3$
\raisebox{-0.4\height}{
\begin{texdraw}
\drawdim em \setunitscale 0.13 \linewd 0.5

\move(40 80)\lvec(50 80)\lvec(50 85)\lvec(40 85)\lvec(40
80)\htext(43 81){\tiny $0$}

\move(30 0)\lvec(40 0)\lvec(40 10)\lvec(30 10)\lvec(30 0)\htext(33
1){\tiny $0$}

\move(40 0)\lvec(50 0)\lvec(50 10)\lvec(40 10)\lvec(40 0)\htext(43
1){\tiny $0$}

\move(30 0)\lvec(40 0)\lvec(40 5)\lvec(30 5)\lvec(30 0)\lfill
f:0.8 \htext(33 6){\tiny $0$}

\move(40 0)\lvec(50 0)\lvec(50 5)\lvec(40 5)\lvec(40 0)\lfill
f:0.8 \htext(43 6){\tiny $0$}

\move(40 50)\lvec(50 50)\lvec(50 60)\lvec(40 60)\lvec(40
50)\htext(43 53){$_1$}

\move(30 10)\lvec(40 10)\lvec(40 20)\lvec(30 20)\lvec(30
10)\htext(33 13){$_1$}

\move(40 10)\lvec(50 10)\lvec(50 20)\lvec(40 20)\lvec(40
10)\htext(43 13){$_1$}
\move(40 60)\lvec(50 60)\lvec(50 70)\lvec(40 70)\lvec(40
60)\htext(43 63){$_2$}

\move(40 20)\lvec(50 20)\lvec(50 30)\lvec(40 30)\lvec(40
20)\htext(43 23){$_2$}
\move(40 70)\lvec(50 70)\lvec(50 80)\lvec(40 80)\lvec(40
70)\htext(43 73){$_1$}

\move(40 30)\lvec(50 30)\lvec(50 40)\lvec(40 40)\lvec(40
30)\htext(43 33){$_1$}
\move(40 40)\lvec(50 40)\lvec(50 45)\lvec(40 45)\lvec(40
40)\htext(43 41){\tiny $0$}
\move(40 40)\lvec(50 40)\lvec(50 50)\lvec(40 50)\lvec(40
40)\htext(43 46){\tiny $0$}
\end{texdraw}}\quad .\vskip 5mm

Therefore,

$G(\ \raisebox{-0.4\height}{
\begin{texdraw}
\drawdim em \setunitscale 0.13 \linewd 0.5

\move(20 0)\lvec(30 0)\lvec(30 10)\lvec(20 10)\lvec(20 0)\htext(23
1){\tiny $0$}

\move(30 0)\lvec(40 0)\lvec(40 10)\lvec(30 10)\lvec(30 0)\htext(33
1){\tiny $0$}

\move(40 0)\lvec(50 0)\lvec(50 10)\lvec(40 10)\lvec(40 0)\htext(43
1){\tiny $0$}

\move(20 0)\lvec(30 0)\lvec(30 5)\lvec(20 5)\lvec(20 0)\lfill
f:0.8 \htext(23 6){\tiny $0$}

\move(30 0)\lvec(40 0)\lvec(40 5)\lvec(30 5)\lvec(30 0)\lfill
f:0.8 \htext(33 6){\tiny $0$}

\move(40 0)\lvec(50 0)\lvec(50 5)\lvec(40 5)\lvec(40 0)\lfill
f:0.8 \htext(43 6){\tiny $0$}

\move(20 10)\lvec(30 10)\lvec(30 20)\lvec(20 20)\lvec(20
10)\htext(23 13){$_1$}

\move(30 10)\lvec(40 10)\lvec(40 20)\lvec(30 20)\lvec(30
10)\htext(33 13){$_1$}

\move(40 10)\lvec(50 10)\lvec(50 20)\lvec(40 20)\lvec(40
10)\htext(43 13){$_1$}
\move(30 20)\lvec(40 20)\lvec(40 30)\lvec(30 30)\lvec(30
20)\htext(33 23){$_2$}

\move(40 20)\lvec(50 20)\lvec(50 30)\lvec(40 30)\lvec(40
20)\htext(43 23){$_2$}
\move(30 30)\lvec(40 30)\lvec(40 40)\lvec(30 40)\lvec(30
30)\htext(33 33){$_1$}

\move(40 30)\lvec(50 30)\lvec(50 40)\lvec(40 40)\lvec(40
30)\htext(43 33){$_1$}
\move(40 40)\lvec(50 40)\lvec(50 45)\lvec(40 45)\lvec(40
40)\htext(43 41){\tiny $0$}
\move(40 40)\lvec(50 40)\lvec(50 50)\lvec(40 50)\lvec(40
40)\htext(43 46){\tiny $0$}
\end{texdraw}} \ )$
$=A(\ \raisebox{-0.4\height}{
\begin{texdraw}
\drawdim em \setunitscale 0.13 \linewd 0.5

\move(20 0)\lvec(30 0)\lvec(30 10)\lvec(20 10)\lvec(20 0)\htext(23
1){\tiny $0$}

\move(30 0)\lvec(40 0)\lvec(40 10)\lvec(30 10)\lvec(30 0)\htext(33
1){\tiny $0$}

\move(40 0)\lvec(50 0)\lvec(50 10)\lvec(40 10)\lvec(40 0)\htext(43
1){\tiny $0$}

\move(20 0)\lvec(30 0)\lvec(30 5)\lvec(20 5)\lvec(20 0)\lfill
f:0.8 \htext(23 6){\tiny $0$}

\move(30 0)\lvec(40 0)\lvec(40 5)\lvec(30 5)\lvec(30 0)\lfill
f:0.8 \htext(33 6){\tiny $0$}

\move(40 0)\lvec(50 0)\lvec(50 5)\lvec(40 5)\lvec(40 0)\lfill
f:0.8 \htext(43 6){\tiny $0$}

\move(20 10)\lvec(30 10)\lvec(30 20)\lvec(20 20)\lvec(20
10)\htext(23 13){$_1$}

\move(30 10)\lvec(40 10)\lvec(40 20)\lvec(30 20)\lvec(30
10)\htext(33 13){$_1$}

\move(40 10)\lvec(50 10)\lvec(50 20)\lvec(40 20)\lvec(40
10)\htext(43 13){$_1$}
\move(30 20)\lvec(40 20)\lvec(40 30)\lvec(30 30)\lvec(30
20)\htext(33 23){$_2$}

\move(40 20)\lvec(50 20)\lvec(50 30)\lvec(40 30)\lvec(40
20)\htext(43 23){$_2$}
\move(30 30)\lvec(40 30)\lvec(40 40)\lvec(30 40)\lvec(30
30)\htext(33 33){$_1$}

\move(40 30)\lvec(50 30)\lvec(50 40)\lvec(40 40)\lvec(40
30)\htext(43 33){$_1$}
\move(40 40)\lvec(50 40)\lvec(50 45)\lvec(40 45)\lvec(40
40)\htext(43 41){\tiny $0$}
\move(40 40)\lvec(50 40)\lvec(50 50)\lvec(40 50)\lvec(40
40)\htext(43 46){\tiny $0$}
\end{texdraw}} \ )-G(\ \raisebox{-0.4\height}{
\begin{texdraw}
\drawdim em \setunitscale 0.13 \linewd 0.5

\move(20 0)\lvec(30 0)\lvec(30 10)\lvec(20 10)\lvec(20 0)\htext(23
1){\tiny $0$}

\move(30 0)\lvec(40 0)\lvec(40 10)\lvec(30 10)\lvec(30 0)\htext(33
1){\tiny $0$}

\move(40 0)\lvec(50 0)\lvec(50 10)\lvec(40 10)\lvec(40 0)\htext(43
1){\tiny $0$}

\move(20 0)\lvec(30 0)\lvec(30 5)\lvec(20 5)\lvec(20 0)\lfill
f:0.8 \htext(23 6){\tiny $0$}

\move(30 0)\lvec(40 0)\lvec(40 5)\lvec(30 5)\lvec(30 0)\lfill
f:0.8 \htext(33 6){\tiny $0$}

\move(40 0)\lvec(50 0)\lvec(50 5)\lvec(40 5)\lvec(40 0)\lfill
f:0.8 \htext(43 6){\tiny $0$}

\move(40 50)\lvec(50 50)\lvec(50 60)\lvec(40 60)\lvec(40
50)\htext(43 53){$_1$}

\move(30 10)\lvec(40 10)\lvec(40 20)\lvec(30 20)\lvec(30
10)\htext(33 13){$_1$}

\move(40 10)\lvec(50 10)\lvec(50 20)\lvec(40 20)\lvec(40
10)\htext(43 13){$_1$}
\move(30 20)\lvec(40 20)\lvec(40 30)\lvec(30 30)\lvec(30
20)\htext(33 23){$_2$}

\move(40 20)\lvec(50 20)\lvec(50 30)\lvec(40 30)\lvec(40
20)\htext(43 23){$_2$}
\move(30 30)\lvec(40 30)\lvec(40 40)\lvec(30 40)\lvec(30
30)\htext(33 33){$_1$}

\move(40 30)\lvec(50 30)\lvec(50 40)\lvec(40 40)\lvec(40
30)\htext(43 33){$_1$}
\move(40 40)\lvec(50 40)\lvec(50 45)\lvec(40 45)\lvec(40
40)\htext(43 41){\tiny $0$}
\move(40 40)\lvec(50 40)\lvec(50 50)\lvec(40 50)\lvec(40
40)\htext(43 46){\tiny $0$}
\end{texdraw}} \ )$\vskip 5mm

\hskip 1cm $=$ \raisebox{-0.4\height}{
\begin{texdraw}
\drawdim em \setunitscale 0.13 \linewd 0.5

\move(20 0)\lvec(30 0)\lvec(30 10)\lvec(20 10)\lvec(20 0)\htext(23
1){\tiny $0$}

\move(30 0)\lvec(40 0)\lvec(40 10)\lvec(30 10)\lvec(30 0)\htext(33
1){\tiny $0$}

\move(40 0)\lvec(50 0)\lvec(50 10)\lvec(40 10)\lvec(40 0)\htext(43
1){\tiny $0$}

\move(20 0)\lvec(30 0)\lvec(30 5)\lvec(20 5)\lvec(20 0)\lfill
f:0.8 \htext(23 6){\tiny $0$}

\move(30 0)\lvec(40 0)\lvec(40 5)\lvec(30 5)\lvec(30 0)\lfill
f:0.8 \htext(33 6){\tiny $0$}

\move(40 0)\lvec(50 0)\lvec(50 5)\lvec(40 5)\lvec(40 0)\lfill
f:0.8 \htext(43 6){\tiny $0$}

\move(20 10)\lvec(30 10)\lvec(30 20)\lvec(20 20)\lvec(20
10)\htext(23 13){$_1$}

\move(30 10)\lvec(40 10)\lvec(40 20)\lvec(30 20)\lvec(30
10)\htext(33 13){$_1$}

\move(40 10)\lvec(50 10)\lvec(50 20)\lvec(40 20)\lvec(40
10)\htext(43 13){$_1$}
\move(30 20)\lvec(40 20)\lvec(40 30)\lvec(30 30)\lvec(30
20)\htext(33 23){$_2$}

\move(40 20)\lvec(50 20)\lvec(50 30)\lvec(40 30)\lvec(40
20)\htext(43 23){$_2$}
\move(30 30)\lvec(40 30)\lvec(40 40)\lvec(30 40)\lvec(30
30)\htext(33 33){$_1$}

\move(40 30)\lvec(50 30)\lvec(50 40)\lvec(40 40)\lvec(40
30)\htext(43 33){$_1$}
\move(40 40)\lvec(50 40)\lvec(50 45)\lvec(40 45)\lvec(40
40)\htext(43 41){\tiny $0$}
\move(40 40)\lvec(50 40)\lvec(50 50)\lvec(40 50)\lvec(40
40)\htext(43 46){\tiny $0$}
\end{texdraw}}
$+$\quad $q^2$ \raisebox{-0.4\height}{
\begin{texdraw}
\drawdim em \setunitscale 0.13 \linewd 0.5

\move(20 0)\lvec(30 0)\lvec(30 10)\lvec(20 10)\lvec(20 0)\htext(23
1){\tiny $0$}

\move(30 0)\lvec(40 0)\lvec(40 10)\lvec(30 10)\lvec(30 0)\htext(33
1){\tiny $0$}

\move(40 0)\lvec(50 0)\lvec(50 10)\lvec(40 10)\lvec(40 0)\htext(43
1){\tiny $0$}

\move(20 0)\lvec(30 0)\lvec(30 5)\lvec(20 5)\lvec(20 0)\lfill
f:0.8 \htext(23 6){\tiny $0$}

\move(30 0)\lvec(40 0)\lvec(40 5)\lvec(30 5)\lvec(30 0)\lfill
f:0.8 \htext(33 6){\tiny $0$}

\move(40 0)\lvec(50 0)\lvec(50 5)\lvec(40 5)\lvec(40 0)\lfill
f:0.8 \htext(43 6){\tiny $0$}

\move(20 10)\lvec(30 10)\lvec(30 20)\lvec(20 20)\lvec(20
10)\htext(23 13){$_1$}

\move(30 10)\lvec(40 10)\lvec(40 20)\lvec(30 20)\lvec(30
10)\htext(33 13){$_1$}

\move(40 10)\lvec(50 10)\lvec(50 20)\lvec(40 20)\lvec(40
10)\htext(43 13){$_1$}
\move(30 20)\lvec(40 20)\lvec(40 30)\lvec(30 30)\lvec(30
20)\htext(33 23){$_2$}

\move(40 20)\lvec(50 20)\lvec(50 30)\lvec(40 30)\lvec(40
20)\htext(43 23){$_2$}
\move(40 50)\lvec(50 50)\lvec(50 60)\lvec(40 60)\lvec(40
50)\htext(43 53){$_1$}

\move(40 30)\lvec(50 30)\lvec(50 40)\lvec(40 40)\lvec(40
30)\htext(43 33){$_1$}
\move(40 40)\lvec(50 40)\lvec(50 45)\lvec(40 45)\lvec(40
40)\htext(43 41){\tiny $0$}
\move(40 40)\lvec(50 40)\lvec(50 50)\lvec(40 50)\lvec(40
40)\htext(43 46){\tiny $0$}
\end{texdraw}}\quad
$+$\quad $q^4$ \raisebox{-0.4\height}{
\begin{texdraw}
\drawdim em \setunitscale 0.13 \linewd 0.5

\move(20 0)\lvec(30 0)\lvec(30 10)\lvec(20 10)\lvec(20 0)\htext(23
1){\tiny $0$}

\move(30 0)\lvec(40 0)\lvec(40 10)\lvec(30 10)\lvec(30 0)\htext(33
1){\tiny $0$}

\move(40 0)\lvec(50 0)\lvec(50 10)\lvec(40 10)\lvec(40 0)\htext(43
1){\tiny $0$}

\move(20 0)\lvec(30 0)\lvec(30 5)\lvec(20 5)\lvec(20 0)\lfill
f:0.8 \htext(23 6){\tiny $0$}

\move(30 0)\lvec(40 0)\lvec(40 5)\lvec(30 5)\lvec(30 0)\lfill
f:0.8 \htext(33 6){\tiny $0$}

\move(40 0)\lvec(50 0)\lvec(50 5)\lvec(40 5)\lvec(40 0)\lfill
f:0.8 \htext(43 6){\tiny $0$}

\move(40 50)\lvec(50 50)\lvec(50 60)\lvec(40 60)\lvec(40
50)\htext(43 53){$_1$}

\move(30 10)\lvec(40 10)\lvec(40 20)\lvec(30 20)\lvec(30
10)\htext(33 13){$_1$}

\move(40 10)\lvec(50 10)\lvec(50 20)\lvec(40 20)\lvec(40
10)\htext(43 13){$_1$}
\move(30 20)\lvec(40 20)\lvec(40 30)\lvec(30 30)\lvec(30
20)\htext(33 23){$_2$}

\move(40 20)\lvec(50 20)\lvec(50 30)\lvec(40 30)\lvec(40
20)\htext(43 23){$_2$}
\move(30 30)\lvec(40 30)\lvec(40 40)\lvec(30 40)\lvec(30
30)\htext(33 33){$_1$}

\move(40 30)\lvec(50 30)\lvec(50 40)\lvec(40 40)\lvec(40
30)\htext(43 33){$_1$}
\move(40 40)\lvec(50 40)\lvec(50 45)\lvec(40 45)\lvec(40
40)\htext(43 41){\tiny $0$}
\move(40 40)\lvec(50 40)\lvec(50 50)\lvec(40 50)\lvec(40
40)\htext(43 46){\tiny $0$}
\end{texdraw}}\quad .\vskip 5mm

}
\end{ex}

\begin{ex}{\rm Suppose that $\frak{g}=B_3^{(1)}$. Note that $q_0=q_1=q_2=q^2$ and
$q_3=q$.\vskip 5mm

(a) We have \vskip 5mm

$G(\ \raisebox{-0.4\height}{
\begin{texdraw}
\drawdim em \setunitscale 0.13 \linewd 0.5

\move(20 0)\lvec(30 0)\lvec(30 10)\lvec(20 10)\lvec(20 0)\htext(22
6){\tiny $0$}

\move(30 0)\lvec(40 0)\lvec(40 10)\lvec(30 10)\lvec(30 0)\htext(32
6){\tiny $1$}

\move(40 0)\lvec(50 0)\lvec(50 10)\lvec(40 10)\lvec(40 0)\htext(42
6){\tiny $0$}

\move(20 0)\lvec(30 10)\lvec(30 0)\lvec(20 0)\lfill f:0.8
\htext(26 2){\tiny $1$}

\move(30 0)\lvec(40 10)\lvec(40 0)\lvec(30 0)\lfill f:0.8
\htext(36 2){\tiny $0$}

\move(40 0)\lvec(50 10)\lvec(50 0)\lvec(40 0)\lfill f:0.8
\htext(46 2){\tiny $1$}
\move(20 10)\lvec(30 10)\lvec(30 20)\lvec(20 20)\lvec(20
10)\htext(23 13){$_2$}

\move(30 10)\lvec(40 10)\lvec(40 20)\lvec(30 20)\lvec(30
10)\htext(33 13){$_2$}

\move(40 10)\lvec(50 10)\lvec(50 20)\lvec(40 20)\lvec(40
10)\htext(43 13){$_2$}

\move(30 20)\lvec(40 20)\lvec(40 30)\lvec(30 30)\lvec(30
20)\htext(33 26){\tiny $3$}

\move(40 20)\lvec(50 20)\lvec(50 30)\lvec(40 30)\lvec(40
20)\htext(43 26){\tiny $3$}

\move(20 20)\lvec(30 20)\lvec(30 25)\lvec(20 25)\lvec(20 20)
\htext(23 21){\tiny $3$}

\move(30 20)\lvec(40 20)\lvec(40 25)\lvec(30 25)\lvec(30 20)
\htext(33 21){\tiny $3$}

\move(40 20)\lvec(50 20)\lvec(50 25)\lvec(40 25)\lvec(40 20)
\htext(43 21){\tiny $3$}

\move(30 30)\lvec(40 30)\lvec(40 40)\lvec(30 40)\lvec(30
30)\htext(33 33){$_2$}

\move(40 30)\lvec(50 30)\lvec(50 40)\lvec(40 40)\lvec(40
30)\htext(43 33){$_2$}

\move(30 40)\lvec(40 50)\lvec(40 40)\lvec(30 40) \htext(36
42){\tiny $0$}

\move(40 40)\lvec(50 50)\lvec(50 40)\lvec(40 40) \htext(46
42){\tiny $1$}

\end{texdraw}} \ )= A(\ \raisebox{-0.4\height}{
\begin{texdraw}
\drawdim em \setunitscale 0.13 \linewd 0.5

\move(20 0)\lvec(30 0)\lvec(30 10)\lvec(20 10)\lvec(20 0)\htext(22
6){\tiny $0$}

\move(30 0)\lvec(40 0)\lvec(40 10)\lvec(30 10)\lvec(30 0)\htext(32
6){\tiny $1$}

\move(40 0)\lvec(50 0)\lvec(50 10)\lvec(40 10)\lvec(40 0)\htext(42
6){\tiny $0$}

\move(20 0)\lvec(30 10)\lvec(30 0)\lvec(20 0)\lfill f:0.8
\htext(26 2){\tiny $1$}

\move(30 0)\lvec(40 10)\lvec(40 0)\lvec(30 0)\lfill f:0.8
\htext(36 2){\tiny $0$}

\move(40 0)\lvec(50 10)\lvec(50 0)\lvec(40 0)\lfill f:0.8
\htext(46 2){\tiny $1$}
\move(20 10)\lvec(30 10)\lvec(30 20)\lvec(20 20)\lvec(20
10)\htext(23 13){$_2$}

\move(30 10)\lvec(40 10)\lvec(40 20)\lvec(30 20)\lvec(30
10)\htext(33 13){$_2$}

\move(40 10)\lvec(50 10)\lvec(50 20)\lvec(40 20)\lvec(40
10)\htext(43 13){$_2$}

\move(30 20)\lvec(40 20)\lvec(40 30)\lvec(30 30)\lvec(30
20)\htext(33 26){\tiny $3$}

\move(40 20)\lvec(50 20)\lvec(50 30)\lvec(40 30)\lvec(40
20)\htext(43 26){\tiny $3$}

\move(20 20)\lvec(30 20)\lvec(30 25)\lvec(20 25)\lvec(20 20)
\htext(23 21){\tiny $3$}

\move(30 20)\lvec(40 20)\lvec(40 25)\lvec(30 25)\lvec(30 20)
\htext(33 21){\tiny $3$}

\move(40 20)\lvec(50 20)\lvec(50 25)\lvec(40 25)\lvec(40 20)
\htext(43 21){\tiny $3$}

\move(30 30)\lvec(40 30)\lvec(40 40)\lvec(30 40)\lvec(30
30)\htext(33 33){$_2$}

\move(40 30)\lvec(50 30)\lvec(50 40)\lvec(40 40)\lvec(40
30)\htext(43 33){$_2$}

\move(30 40)\lvec(40 50)\lvec(40 40)\lvec(30 40) \htext(36
42){\tiny $0$}

\move(40 40)\lvec(50 50)\lvec(50 40)\lvec(40 40) \htext(46
42){\tiny $1$}

\end{texdraw}} \ )
=f_3f_2f_0^{(2)}f_2f_3^{(2)}f_2f_1^{(2)}f_2f_3^{(2)}f_2f_0Y_{\Lambda_0}$\vskip
5mm

$=$  \raisebox{-0.4\height}{
\begin{texdraw}
\drawdim em \setunitscale 0.13 \linewd 0.5

\move(20 0)\lvec(30 0)\lvec(30 10)\lvec(20 10)\lvec(20 0)\htext(22
6){\tiny $0$}

\move(30 0)\lvec(40 0)\lvec(40 10)\lvec(30 10)\lvec(30 0)\htext(32
6){\tiny $1$}

\move(40 0)\lvec(50 0)\lvec(50 10)\lvec(40 10)\lvec(40 0)\htext(42
6){\tiny $0$}

\move(20 0)\lvec(30 10)\lvec(30 0)\lvec(20 0)\lfill f:0.8
\htext(26 2){\tiny $1$}

\move(30 0)\lvec(40 10)\lvec(40 0)\lvec(30 0)\lfill f:0.8
\htext(36 2){\tiny $0$}

\move(40 0)\lvec(50 10)\lvec(50 0)\lvec(40 0)\lfill f:0.8
\htext(46 2){\tiny $1$}
\move(20 10)\lvec(30 10)\lvec(30 20)\lvec(20 20)\lvec(20
10)\htext(23 13){$_2$}

\move(30 10)\lvec(40 10)\lvec(40 20)\lvec(30 20)\lvec(30
10)\htext(33 13){$_2$}

\move(40 10)\lvec(50 10)\lvec(50 20)\lvec(40 20)\lvec(40
10)\htext(43 13){$_2$}

\move(30 20)\lvec(40 20)\lvec(40 30)\lvec(30 30)\lvec(30
20)\htext(33 26){\tiny $3$}

\move(40 20)\lvec(50 20)\lvec(50 30)\lvec(40 30)\lvec(40
20)\htext(43 26){\tiny $3$}

\move(20 20)\lvec(30 20)\lvec(30 25)\lvec(20 25)\lvec(20 20)
\htext(23 21){\tiny $3$}

\move(30 20)\lvec(40 20)\lvec(40 25)\lvec(30 25)\lvec(30 20)
\htext(33 21){\tiny $3$}

\move(40 20)\lvec(50 20)\lvec(50 25)\lvec(40 25)\lvec(40 20)
\htext(43 21){\tiny $3$}

\move(30 30)\lvec(40 30)\lvec(40 40)\lvec(30 40)\lvec(30
30)\htext(33 33){$_2$}

\move(40 30)\lvec(50 30)\lvec(50 40)\lvec(40 40)\lvec(40
30)\htext(43 33){$_2$}

\move(30 40)\lvec(40 50)\lvec(40 40)\lvec(30 40) \htext(36
42){\tiny $0$}

\move(40 40)\lvec(50 50)\lvec(50 40)\lvec(40 40) \htext(46
42){\tiny $1$}

\end{texdraw}}\quad $+$ \quad $q^2$
 \raisebox{-0.4\height}{
\begin{texdraw}
\drawdim em \setunitscale 0.13 \linewd 0.5

\move(20 0)\lvec(30 0)\lvec(30 10)\lvec(20 10)\lvec(20 0)\htext(22
6){\tiny $0$}

\move(30 0)\lvec(40 0)\lvec(40 10)\lvec(30 10)\lvec(30 0)\htext(32
6){\tiny $1$}

\move(40 0)\lvec(50 0)\lvec(50 10)\lvec(40 10)\lvec(40 0)\htext(42
6){\tiny $0$}

\move(20 0)\lvec(30 10)\lvec(30 0)\lvec(20 0)\lfill f:0.8
\htext(26 2){\tiny $1$}

\move(30 0)\lvec(40 10)\lvec(40 0)\lvec(30 0)\lfill f:0.8
\htext(36 2){\tiny $0$}

\move(40 0)\lvec(50 10)\lvec(50 0)\lvec(40 0)\lfill f:0.8
\htext(46 2){\tiny $1$}
\move(20 10)\lvec(30 10)\lvec(30 20)\lvec(20 20)\lvec(20
10)\htext(23 13){$_2$}

\move(30 10)\lvec(40 10)\lvec(40 20)\lvec(30 20)\lvec(30
10)\htext(33 13){$_2$}

\move(40 10)\lvec(50 10)\lvec(50 20)\lvec(40 20)\lvec(40
10)\htext(43 13){$_2$}

\move(30 20)\lvec(40 20)\lvec(40 30)\lvec(30 30)\lvec(30
20)\htext(33 26){\tiny $3$}

\move(40 20)\lvec(50 20)\lvec(50 30)\lvec(40 30)\lvec(40
20)\htext(43 26){\tiny $3$}

\move(20 20)\lvec(30 20)\lvec(30 25)\lvec(20 25)\lvec(20 20)
\htext(23 21){\tiny $3$}

\move(30 20)\lvec(40 20)\lvec(40 25)\lvec(30 25)\lvec(30 20)
\htext(33 21){\tiny $3$}

\move(40 20)\lvec(50 20)\lvec(50 25)\lvec(40 25)\lvec(40 20)
\htext(43 21){\tiny $3$}

\move(30 30)\lvec(40 30)\lvec(40 40)\lvec(30 40)\lvec(30
30)\htext(33 33){$_2$}

\move(40 30)\lvec(50 30)\lvec(50 40)\lvec(40 40)\lvec(40
30)\htext(43 33){$_2$}

\move(40 40)\lvec(40 50)\lvec(50 50)\lvec(40 40) \htext(42
46){\tiny $0$}

\move(40 40)\lvec(50 50)\lvec(50 40)\lvec(40 40) \htext(46
42){\tiny $1$}

\end{texdraw}}\quad $+$ \quad $q^4$
 \raisebox{-0.4\height}{
\begin{texdraw}
\drawdim em \setunitscale 0.13 \linewd 0.5

\move(20 0)\lvec(30 0)\lvec(30 10)\lvec(20 10)\lvec(20 0)\htext(22
6){\tiny $0$}

\move(30 0)\lvec(40 0)\lvec(40 10)\lvec(30 10)\lvec(30 0)\htext(32
6){\tiny $1$}

\move(40 0)\lvec(50 0)\lvec(50 10)\lvec(40 10)\lvec(40 0)\htext(42
6){\tiny $0$}

\move(20 0)\lvec(30 10)\lvec(30 0)\lvec(20 0)\lfill f:0.8
\htext(26 2){\tiny $1$}

\move(30 0)\lvec(40 10)\lvec(40 0)\lvec(30 0)\lfill f:0.8
\htext(36 2){\tiny $0$}

\move(40 0)\lvec(50 10)\lvec(50 0)\lvec(40 0)\lfill f:0.8
\htext(46 2){\tiny $1$}
\move(40 50)\lvec(50 50)\lvec(50 60)\lvec(40 60)\lvec(40
50)\htext(43 53){$_2$}

\move(30 10)\lvec(40 10)\lvec(40 20)\lvec(30 20)\lvec(30
10)\htext(33 13){$_2$}

\move(40 10)\lvec(50 10)\lvec(50 20)\lvec(40 20)\lvec(40
10)\htext(43 13){$_2$}

\move(30 20)\lvec(40 20)\lvec(40 30)\lvec(30 30)\lvec(30
20)\htext(33 26){\tiny $3$}

\move(40 20)\lvec(50 20)\lvec(50 30)\lvec(40 30)\lvec(40
20)\htext(43 26){\tiny $3$}

\move(40 60)\lvec(50 60)\lvec(50 65)\lvec(40 65)\lvec(40 60)
\htext(43 61){\tiny $3$}

\move(30 20)\lvec(40 20)\lvec(40 25)\lvec(30 25)\lvec(30 20)
\htext(33 21){\tiny $3$}

\move(40 20)\lvec(50 20)\lvec(50 25)\lvec(40 25)\lvec(40 20)
\htext(43 21){\tiny $3$}

\move(30 30)\lvec(40 30)\lvec(40 40)\lvec(30 40)\lvec(30
30)\htext(33 33){$_2$}

\move(40 30)\lvec(50 30)\lvec(50 40)\lvec(40 40)\lvec(40
30)\htext(43 33){$_2$}

\move(40 40)\lvec(40 50)\lvec(50 50)\lvec(40 40) \htext(42
46){\tiny $0$}

\move(40 40)\lvec(50 50)\lvec(50 40)\lvec(40 40) \htext(46
42){\tiny $1$}

\end{texdraw}}\quad $+$ \quad $q^4$
 \raisebox{-0.4\height}{
\begin{texdraw}
\drawdim em \setunitscale 0.13 \linewd 0.5

\move(30 0)\lvec(40 0)\lvec(40 10)\lvec(30 10)\lvec(30 0)\htext(32
6){\tiny $1$}

\move(40 0)\lvec(50 0)\lvec(50 10)\lvec(40 10)\lvec(40 0)\htext(42
6){\tiny $0$}
\move(30 0)\lvec(40 10)\lvec(40 0)\lvec(30 0)\lfill f:0.8
\htext(36 2){\tiny $0$}

\move(40 0)\lvec(50 10)\lvec(50 0)\lvec(40 0)\lfill f:0.8
\htext(46 2){\tiny $1$}
\move(40 50)\lvec(50 50)\lvec(50 60)\lvec(40 60)\lvec(40
50)\htext(43 53){$_2$}

\move(30 10)\lvec(40 10)\lvec(40 20)\lvec(30 20)\lvec(30
10)\htext(33 13){$_2$}

\move(40 10)\lvec(50 10)\lvec(50 20)\lvec(40 20)\lvec(40
10)\htext(43 13){$_2$}

\move(30 20)\lvec(40 20)\lvec(40 30)\lvec(30 30)\lvec(30
20)\htext(33 26){\tiny $3$}

\move(40 20)\lvec(50 20)\lvec(50 30)\lvec(40 30)\lvec(40
20)\htext(43 26){\tiny $3$}

\move(40 60)\lvec(50 60)\lvec(50 65)\lvec(40 65)\lvec(40 60)
\htext(43 61){\tiny $3$}

\move(30 20)\lvec(40 20)\lvec(40 25)\lvec(30 25)\lvec(30 20)
\htext(33 21){\tiny $3$}

\move(40 20)\lvec(50 20)\lvec(50 25)\lvec(40 25)\lvec(40 20)
\htext(43 21){\tiny $3$}

\move(30 30)\lvec(40 30)\lvec(40 40)\lvec(30 40)\lvec(30
30)\htext(33 33){$_2$}

\move(40 30)\lvec(50 30)\lvec(50 40)\lvec(40 40)\lvec(40
30)\htext(43 33){$_2$}
\move(30 40)\lvec(40 50)\lvec(40 40)\lvec(30 40) \htext(36
42){\tiny $0$}

\move(40 40)\lvec(40 50)\lvec(50 50)\lvec(40 40) \htext(42
46){\tiny $0$}

\move(40 40)\lvec(50 50)\lvec(50 40)\lvec(40 40) \htext(46
42){\tiny $1$}

\end{texdraw}}\quad ,\vskip 10mm

$G(\ \raisebox{-0.5\height}{
\begin{texdraw}
\drawdim em \setunitscale 0.13 \linewd 0.5

\move(20 0)\lvec(30 0)\lvec(30 10)\lvec(20 10)\lvec(20 0)\htext(22
6){\tiny $0$}

\move(30 0)\lvec(40 0)\lvec(40 10)\lvec(30 10)\lvec(30 0)\htext(32
6){\tiny $1$}

\move(40 0)\lvec(50 0)\lvec(50 10)\lvec(40 10)\lvec(40 0)\htext(42
6){\tiny $0$}

\move(20 0)\lvec(30 10)\lvec(30 0)\lvec(20 0)\lfill f:0.8
\htext(26 2){\tiny $1$}

\move(30 0)\lvec(40 10)\lvec(40 0)\lvec(30 0)\lfill f:0.8
\htext(36 2){\tiny $0$}

\move(40 0)\lvec(50 10)\lvec(50 0)\lvec(40 0)\lfill f:0.8
\htext(46 2){\tiny $1$}
\move(20 10)\lvec(30 10)\lvec(30 20)\lvec(20 20)\lvec(20
10)\htext(23 13){$_2$}

\move(30 10)\lvec(40 10)\lvec(40 20)\lvec(30 20)\lvec(30
10)\htext(33 13){$_2$}

\move(40 10)\lvec(50 10)\lvec(50 20)\lvec(40 20)\lvec(40
10)\htext(43 13){$_2$}
\move(20 20)\lvec(30 20)\lvec(30 25)\lvec(20 25)\lvec(20 20)
\htext(23 21){\tiny $3$}

\move(30 20)\lvec(40 20)\lvec(40 25)\lvec(30 25)\lvec(30 20)
\htext(33 21){\tiny $3$}

\move(40 20)\lvec(50 20)\lvec(50 25)\lvec(40 25)\lvec(40 20)
\htext(43 21){\tiny $3$}
\end{texdraw}} \
)=A(\ \raisebox{-0.5\height}{
\begin{texdraw}
\drawdim em \setunitscale 0.13 \linewd 0.5

\move(20 0)\lvec(30 0)\lvec(30 10)\lvec(20 10)\lvec(20 0)\htext(22
6){\tiny $0$}

\move(30 0)\lvec(40 0)\lvec(40 10)\lvec(30 10)\lvec(30 0)\htext(32
6){\tiny $1$}

\move(40 0)\lvec(50 0)\lvec(50 10)\lvec(40 10)\lvec(40 0)\htext(42
6){\tiny $0$}

\move(20 0)\lvec(30 10)\lvec(30 0)\lvec(20 0)\lfill f:0.8
\htext(26 2){\tiny $1$}

\move(30 0)\lvec(40 10)\lvec(40 0)\lvec(30 0)\lfill f:0.8
\htext(36 2){\tiny $0$}

\move(40 0)\lvec(50 10)\lvec(50 0)\lvec(40 0)\lfill f:0.8
\htext(46 2){\tiny $1$}
\move(20 10)\lvec(30 10)\lvec(30 20)\lvec(20 20)\lvec(20
10)\htext(23 13){$_2$}

\move(30 10)\lvec(40 10)\lvec(40 20)\lvec(30 20)\lvec(30
10)\htext(33 13){$_2$}

\move(40 10)\lvec(50 10)\lvec(50 20)\lvec(40 20)\lvec(40
10)\htext(43 13){$_2$}
\move(20 20)\lvec(30 20)\lvec(30 25)\lvec(20 25)\lvec(20 20)
\htext(23 21){\tiny $3$}

\move(30 20)\lvec(40 20)\lvec(40 25)\lvec(30 25)\lvec(30 20)
\htext(33 21){\tiny $3$}

\move(40 20)\lvec(50 20)\lvec(50 25)\lvec(40 25)\lvec(40 20)
\htext(43 21){\tiny $3$}
\end{texdraw}} \ )=f_3f_2f_0f_3f_2f_1f_3f_2f_0Y_{\Lambda_0}$\vskip 5mm

$=$ \raisebox{-0.5\height}{
\begin{texdraw}
\drawdim em \setunitscale 0.13 \linewd 0.5

\move(20 0)\lvec(30 0)\lvec(30 10)\lvec(20 10)\lvec(20 0)\htext(22
6){\tiny $0$}

\move(30 0)\lvec(40 0)\lvec(40 10)\lvec(30 10)\lvec(30 0)\htext(32
6){\tiny $1$}

\move(40 0)\lvec(50 0)\lvec(50 10)\lvec(40 10)\lvec(40 0)\htext(42
6){\tiny $0$}

\move(20 0)\lvec(30 10)\lvec(30 0)\lvec(20 0)\lfill f:0.8
\htext(26 2){\tiny $1$}

\move(30 0)\lvec(40 10)\lvec(40 0)\lvec(30 0)\lfill f:0.8
\htext(36 2){\tiny $0$}

\move(40 0)\lvec(50 10)\lvec(50 0)\lvec(40 0)\lfill f:0.8
\htext(46 2){\tiny $1$}
\move(20 10)\lvec(30 10)\lvec(30 20)\lvec(20 20)\lvec(20
10)\htext(23 13){$_2$}

\move(30 10)\lvec(40 10)\lvec(40 20)\lvec(30 20)\lvec(30
10)\htext(33 13){$_2$}

\move(40 10)\lvec(50 10)\lvec(50 20)\lvec(40 20)\lvec(40
10)\htext(43 13){$_2$}
\move(20 20)\lvec(30 20)\lvec(30 25)\lvec(20 25)\lvec(20 20)
\htext(23 21){\tiny $3$}

\move(30 20)\lvec(40 20)\lvec(40 25)\lvec(30 25)\lvec(30 20)
\htext(33 21){\tiny $3$}

\move(40 20)\lvec(50 20)\lvec(50 25)\lvec(40 25)\lvec(40 20)
\htext(43 21){\tiny $3$}
\end{texdraw}}\quad$+$\quad $q(1-q^4)$
 \raisebox{-0.5\height}{
\begin{texdraw}
\drawdim em \setunitscale 0.13 \linewd 0.5

\move(20 0)\lvec(30 0)\lvec(30 10)\lvec(20 10)\lvec(20 0)\htext(22
6){\tiny $0$}

\move(30 0)\lvec(40 0)\lvec(40 10)\lvec(30 10)\lvec(30 0)\htext(32
6){\tiny $1$}

\move(40 0)\lvec(50 0)\lvec(50 10)\lvec(40 10)\lvec(40 0)\htext(42
6){\tiny $0$}

\move(20 0)\lvec(30 10)\lvec(30 0)\lvec(20 0)\lfill f:0.8
\htext(26 2){\tiny $1$}

\move(30 0)\lvec(40 10)\lvec(40 0)\lvec(30 0)\lfill f:0.8
\htext(36 2){\tiny $0$}

\move(40 0)\lvec(50 10)\lvec(50 0)\lvec(40 0)\lfill f:0.8
\htext(46 2){\tiny $1$}
\move(20 10)\lvec(30 10)\lvec(30 20)\lvec(20 20)\lvec(20
10)\htext(23 13){$_2$}

\move(30 10)\lvec(40 10)\lvec(40 20)\lvec(30 20)\lvec(30
10)\htext(33 13){$_2$}

\move(40 10)\lvec(50 10)\lvec(50 20)\lvec(40 20)\lvec(40
10)\htext(43 13){$_2$}
\move(40 25)\lvec(50 25)\lvec(50 30)\lvec(40 30)\lvec(40 25)
\htext(43 26){\tiny $3$}

\move(30 20)\lvec(40 20)\lvec(40 25)\lvec(30 25)\lvec(30 20)
\htext(33 21){\tiny $3$}

\move(40 20)\lvec(50 20)\lvec(50 25)\lvec(40 25)\lvec(40 20)
\htext(43 21){\tiny $3$}
\end{texdraw}}\quad$+$\quad$q(1+q^2)$
 \raisebox{-0.5\height}{
\begin{texdraw}
\drawdim em \setunitscale 0.13 \linewd 0.5

\move(20 0)\lvec(30 0)\lvec(30 10)\lvec(20 10)\lvec(20 0)\htext(22
6){\tiny $0$}

\move(30 0)\lvec(40 0)\lvec(40 10)\lvec(30 10)\lvec(30 0)\htext(32
6){\tiny $1$}

\move(40 0)\lvec(50 0)\lvec(50 10)\lvec(40 10)\lvec(40 0)\htext(42
6){\tiny $0$}

\move(20 0)\lvec(30 10)\lvec(30 0)\lvec(20 0)\lfill f:0.8
\htext(26 2){\tiny $1$}

\move(30 0)\lvec(40 10)\lvec(40 0)\lvec(30 0)\lfill f:0.8
\htext(36 2){\tiny $0$}

\move(40 0)\lvec(50 10)\lvec(50 0)\lvec(40 0)\lfill f:0.8
\htext(46 2){\tiny $1$}
\move(40 30)\lvec(50 30)\lvec(50 40)\lvec(40 40)\lvec(40
30)\htext(43 33){$_2$}

\move(30 10)\lvec(40 10)\lvec(40 20)\lvec(30 20)\lvec(30
10)\htext(33 13){$_2$}

\move(40 10)\lvec(50 10)\lvec(50 20)\lvec(40 20)\lvec(40
10)\htext(43 13){$_2$}
\move(40 25)\lvec(50 25)\lvec(50 30)\lvec(40 30)\lvec(40 25)
\htext(43 26){\tiny $3$}

\move(30 20)\lvec(40 20)\lvec(40 25)\lvec(30 25)\lvec(30 20)
\htext(33 21){\tiny $3$}

\move(40 20)\lvec(50 20)\lvec(50 25)\lvec(40 25)\lvec(40 20)
\htext(43 21){\tiny $3$}
\end{texdraw}}\quad .\vskip 5mm

(b) Observe that \vskip 5mm

$A(\ \raisebox{-0.4\height}{
\begin{texdraw}
\drawdim em \setunitscale 0.13 \linewd 0.5

\move(20 0)\lvec(30 0)\lvec(30 10)\lvec(20 10)\lvec(20 0)\htext(22
6){\tiny $0$}

\move(30 0)\lvec(40 0)\lvec(40 10)\lvec(30 10)\lvec(30 0)\htext(32
6){\tiny $1$}

\move(40 0)\lvec(50 0)\lvec(50 10)\lvec(40 10)\lvec(40 0)\htext(42
6){\tiny $0$}

\move(20 0)\lvec(30 10)\lvec(30 0)\lvec(20 0)\lfill f:0.8
\htext(26 2){\tiny $1$}

\move(30 0)\lvec(40 10)\lvec(40 0)\lvec(30 0)\lfill f:0.8
\htext(36 2){\tiny $0$}

\move(40 0)\lvec(50 10)\lvec(50 0)\lvec(40 0)\lfill f:0.8
\htext(46 2){\tiny $1$}
\move(20 10)\lvec(30 10)\lvec(30 20)\lvec(20 20)\lvec(20
10)\htext(23 13){$_2$}

\move(30 10)\lvec(40 10)\lvec(40 20)\lvec(30 20)\lvec(30
10)\htext(33 13){$_2$}

\move(40 10)\lvec(50 10)\lvec(50 20)\lvec(40 20)\lvec(40
10)\htext(43 13){$_2$}

\move(30 20)\lvec(40 20)\lvec(40 30)\lvec(30 30)\lvec(30
20)\htext(33 26){\tiny $3$}

\move(40 20)\lvec(50 20)\lvec(50 30)\lvec(40 30)\lvec(40
20)\htext(43 26){\tiny $3$}

\move(30 20)\lvec(40 20)\lvec(40 25)\lvec(30 25)\lvec(30 20)
\htext(33 21){\tiny $3$}

\move(40 20)\lvec(50 20)\lvec(50 25)\lvec(40 25)\lvec(40 20)
\htext(43 21){\tiny $3$}

\move(30 30)\lvec(40 30)\lvec(40 40)\lvec(30 40)\lvec(30
30)\htext(33 33){$_2$}

\move(40 30)\lvec(50 30)\lvec(50 40)\lvec(40 40)\lvec(40
30)\htext(43 33){$_2$}

\move(40 40)\lvec(50 50)\lvec(40 50)\lvec(40 40) \htext(41
46){\tiny $0$}

\move(40 40)\lvec(50 50)\lvec(50 40)\lvec(40 40) \htext(46
42){\tiny $1$}

\end{texdraw}} \
)=f_2f_0f_2f_3^{(2)}f_2f_1^{(2)}f_0f_2f_3^{(2)}f_2f_0Y_{\Lambda_0}$\vskip
5mm

$=$ \raisebox{-0.4\height}{
\begin{texdraw}
\drawdim em \setunitscale 0.13 \linewd 0.5

\move(20 0)\lvec(30 0)\lvec(30 10)\lvec(20 10)\lvec(20 0)\htext(22
6){\tiny $0$}

\move(30 0)\lvec(40 0)\lvec(40 10)\lvec(30 10)\lvec(30 0)\htext(32
6){\tiny $1$}

\move(40 0)\lvec(50 0)\lvec(50 10)\lvec(40 10)\lvec(40 0)\htext(42
6){\tiny $0$}

\move(20 0)\lvec(30 10)\lvec(30 0)\lvec(20 0)\lfill f:0.8
\htext(26 2){\tiny $1$}

\move(30 0)\lvec(40 10)\lvec(40 0)\lvec(30 0)\lfill f:0.8
\htext(36 2){\tiny $0$}

\move(40 0)\lvec(50 10)\lvec(50 0)\lvec(40 0)\lfill f:0.8
\htext(46 2){\tiny $1$}
\move(20 10)\lvec(30 10)\lvec(30 20)\lvec(20 20)\lvec(20
10)\htext(23 13){$_2$}

\move(30 10)\lvec(40 10)\lvec(40 20)\lvec(30 20)\lvec(30
10)\htext(33 13){$_2$}

\move(40 10)\lvec(50 10)\lvec(50 20)\lvec(40 20)\lvec(40
10)\htext(43 13){$_2$}

\move(30 20)\lvec(40 20)\lvec(40 30)\lvec(30 30)\lvec(30
20)\htext(33 26){\tiny $3$}

\move(40 20)\lvec(50 20)\lvec(50 30)\lvec(40 30)\lvec(40
20)\htext(43 26){\tiny $3$}

\move(30 20)\lvec(40 20)\lvec(40 25)\lvec(30 25)\lvec(30 20)
\htext(33 21){\tiny $3$}

\move(40 20)\lvec(50 20)\lvec(50 25)\lvec(40 25)\lvec(40 20)
\htext(43 21){\tiny $3$}

\move(30 30)\lvec(40 30)\lvec(40 40)\lvec(30 40)\lvec(30
30)\htext(33 33){$_2$}

\move(40 30)\lvec(50 30)\lvec(50 40)\lvec(40 40)\lvec(40
30)\htext(43 33){$_2$}

\move(40 40)\lvec(50 50)\lvec(40 50)\lvec(40 40) \htext(41
46){\tiny $0$}

\move(40 40)\lvec(50 50)\lvec(50 40)\lvec(40 40) \htext(46
42){\tiny $1$}

\end{texdraw}}\quad$+$\quad $q^2$
\raisebox{-0.4\height}{
\begin{texdraw}
\drawdim em \setunitscale 0.13 \linewd 0.5

\move(20 0)\lvec(30 0)\lvec(30 10)\lvec(20 10)\lvec(20 0)\htext(22
6){\tiny $0$}

\move(30 0)\lvec(40 0)\lvec(40 10)\lvec(30 10)\lvec(30 0)\htext(32
6){\tiny $1$}

\move(40 0)\lvec(50 0)\lvec(50 10)\lvec(40 10)\lvec(40 0)\htext(42
6){\tiny $0$}

\move(20 0)\lvec(30 10)\lvec(30 0)\lvec(20 0)\lfill f:0.8
\htext(26 2){\tiny $1$}

\move(30 0)\lvec(40 10)\lvec(40 0)\lvec(30 0)\lfill f:0.8
\htext(36 2){\tiny $0$}

\move(40 0)\lvec(50 10)\lvec(50 0)\lvec(40 0)\lfill f:0.8
\htext(46 2){\tiny $1$}
\move(20 10)\lvec(30 10)\lvec(30 20)\lvec(20 20)\lvec(20
10)\htext(23 13){$_2$}

\move(30 10)\lvec(40 10)\lvec(40 20)\lvec(30 20)\lvec(30
10)\htext(33 13){$_2$}

\move(40 10)\lvec(50 10)\lvec(50 20)\lvec(40 20)\lvec(40
10)\htext(43 13){$_2$}

\move(30 20)\lvec(40 20)\lvec(40 30)\lvec(30 30)\lvec(30
20)\htext(33 26){\tiny $3$}

\move(40 20)\lvec(50 20)\lvec(50 30)\lvec(40 30)\lvec(40
20)\htext(43 26){\tiny $3$}

\move(30 20)\lvec(40 20)\lvec(40 25)\lvec(30 25)\lvec(30 20)
\htext(33 21){\tiny $3$}

\move(40 20)\lvec(50 20)\lvec(50 25)\lvec(40 25)\lvec(40 20)
\htext(43 21){\tiny $3$}

\move(40 50)\lvec(50 50)\lvec(50 60)\lvec(40 60)\lvec(40
50)\htext(43 53){$_2$}

\move(40 30)\lvec(50 30)\lvec(50 40)\lvec(40 40)\lvec(40
30)\htext(43 33){$_2$}

\move(40 40)\lvec(50 50)\lvec(40 50)\lvec(40 40) \htext(41
46){\tiny $0$}

\move(40 40)\lvec(50 50)\lvec(50 40)\lvec(40 40) \htext(46
42){\tiny $1$}

\end{texdraw}}\quad$+$\quad $(1+q^4)$
\raisebox{-0.4\height}{
\begin{texdraw}
\drawdim em \setunitscale 0.13 \linewd 0.5

\move(20 0)\lvec(30 0)\lvec(30 10)\lvec(20 10)\lvec(20 0)\htext(22
6){\tiny $0$}

\move(30 0)\lvec(40 0)\lvec(40 10)\lvec(30 10)\lvec(30 0)\htext(32
6){\tiny $1$}

\move(40 0)\lvec(50 0)\lvec(50 10)\lvec(40 10)\lvec(40 0)\htext(42
6){\tiny $0$}

\move(20 0)\lvec(30 10)\lvec(30 0)\lvec(20 0)\lfill f:0.8
\htext(26 2){\tiny $1$}

\move(30 0)\lvec(40 10)\lvec(40 0)\lvec(30 0)\lfill f:0.8
\htext(36 2){\tiny $0$}

\move(40 0)\lvec(50 10)\lvec(50 0)\lvec(40 0)\lfill f:0.8
\htext(46 2){\tiny $1$}
\move(30 30)\lvec(40 30)\lvec(40 40)\lvec(30 40)\lvec(30
30)\htext(33 33){$_2$}

\move(30 10)\lvec(40 10)\lvec(40 20)\lvec(30 20)\lvec(30
10)\htext(33 13){$_2$}

\move(40 10)\lvec(50 10)\lvec(50 20)\lvec(40 20)\lvec(40
10)\htext(43 13){$_2$}

\move(30 20)\lvec(40 20)\lvec(40 30)\lvec(30 30)\lvec(30
20)\htext(33 26){\tiny $3$}

\move(40 20)\lvec(50 20)\lvec(50 30)\lvec(40 30)\lvec(40
20)\htext(43 26){\tiny $3$}

\move(30 20)\lvec(40 20)\lvec(40 25)\lvec(30 25)\lvec(30 20)
\htext(33 21){\tiny $3$}

\move(40 20)\lvec(50 20)\lvec(50 25)\lvec(40 25)\lvec(40 20)
\htext(43 21){\tiny $3$}

\move(40 50)\lvec(50 50)\lvec(50 60)\lvec(40 60)\lvec(40
50)\htext(43 53){$_2$}

\move(40 30)\lvec(50 30)\lvec(50 40)\lvec(40 40)\lvec(40
30)\htext(43 33){$_2$}

\move(40 40)\lvec(50 50)\lvec(40 50)\lvec(40 40) \htext(41
46){\tiny $0$}

\move(40 40)\lvec(50 50)\lvec(50 40)\lvec(40 40) \htext(46
42){\tiny $1$}

\end{texdraw}}\vskip 5mm

\hskip 1cm $+$\quad $q^2$ \raisebox{-0.4\height}{
\begin{texdraw}
\drawdim em \setunitscale 0.13 \linewd 0.5

\move(20 0)\lvec(30 0)\lvec(30 10)\lvec(20 10)\lvec(20 0)\htext(22
6){\tiny $0$}

\move(30 0)\lvec(40 0)\lvec(40 10)\lvec(30 10)\lvec(30 0)\htext(32
6){\tiny $1$}

\move(40 0)\lvec(50 0)\lvec(50 10)\lvec(40 10)\lvec(40 0)\htext(42
6){\tiny $0$}

\move(20 0)\lvec(30 10)\lvec(30 0)\lvec(20 0)\lfill f:0.8
\htext(26 2){\tiny $1$}

\move(30 0)\lvec(40 10)\lvec(40 0)\lvec(30 0)\lfill f:0.8
\htext(36 2){\tiny $0$}

\move(40 0)\lvec(50 10)\lvec(50 0)\lvec(40 0)\lfill f:0.8
\htext(46 2){\tiny $1$}
\move(40 70)\lvec(50 70)\lvec(50 80)\lvec(40 80)\lvec(40
70)\htext(43 73){$_2$}

\move(30 10)\lvec(40 10)\lvec(40 20)\lvec(30 20)\lvec(30
10)\htext(33 13){$_2$}

\move(40 10)\lvec(50 10)\lvec(50 20)\lvec(40 20)\lvec(40
10)\htext(43 13){$_2$}

\move(40 60)\lvec(50 60)\lvec(50 70)\lvec(40 70)\lvec(40
60)\htext(43 66){\tiny $3$}

\move(40 20)\lvec(50 20)\lvec(50 30)\lvec(40 30)\lvec(40
20)\htext(43 26){\tiny $3$}

\move(40 60)\lvec(50 60)\lvec(50 65)\lvec(40 65)\lvec(40 60)
\htext(43 61){\tiny $3$}

\move(40 20)\lvec(50 20)\lvec(50 25)\lvec(40 25)\lvec(40 20)
\htext(43 21){\tiny $3$}

\move(40 50)\lvec(50 50)\lvec(50 60)\lvec(40 60)\lvec(40
50)\htext(43 53){$_2$}

\move(40 30)\lvec(50 30)\lvec(50 40)\lvec(40 40)\lvec(40
30)\htext(43 33){$_2$}

\move(40 40)\lvec(50 50)\lvec(40 50)\lvec(40 40) \htext(41
46){\tiny $0$}

\move(40 40)\lvec(50 50)\lvec(50 40)\lvec(40 40) \htext(46
42){\tiny $1$}

\end{texdraw}}\quad$+$\quad $q^2$
\raisebox{-0.4\height}{
\begin{texdraw}
\drawdim em \setunitscale 0.13 \linewd 0.5

\move(30 0)\lvec(40 0)\lvec(40 10)\lvec(30 10)\lvec(30 0)\htext(32
6){\tiny $1$}

\move(40 0)\lvec(50 0)\lvec(50 10)\lvec(40 10)\lvec(40 0)\htext(42
6){\tiny $0$}

\move(30 0)\lvec(40 10)\lvec(40 0)\lvec(30 0)\lfill f:0.8
\htext(36 2){\tiny $0$}

\move(30 40)\lvec(40 50)\lvec(40 40)\lvec(30 40)\htext(36
42){\tiny $0$}

\move(40 0)\lvec(50 10)\lvec(50 0)\lvec(40 0)\lfill f:0.8
\htext(46 2){\tiny $1$}
\move(30 30)\lvec(40 30)\lvec(40 40)\lvec(30 40)\lvec(30
30)\htext(33 33){$_2$}

\move(30 10)\lvec(40 10)\lvec(40 20)\lvec(30 20)\lvec(30
10)\htext(33 13){$_2$}

\move(40 10)\lvec(50 10)\lvec(50 20)\lvec(40 20)\lvec(40
10)\htext(43 13){$_2$}

\move(30 20)\lvec(40 20)\lvec(40 30)\lvec(30 30)\lvec(30
20)\htext(33 26){\tiny $3$}

\move(40 20)\lvec(50 20)\lvec(50 30)\lvec(40 30)\lvec(40
20)\htext(43 26){\tiny $3$}

\move(30 20)\lvec(40 20)\lvec(40 25)\lvec(30 25)\lvec(30 20)
\htext(33 21){\tiny $3$}

\move(40 20)\lvec(50 20)\lvec(50 25)\lvec(40 25)\lvec(40 20)
\htext(43 21){\tiny $3$}

\move(40 50)\lvec(50 50)\lvec(50 60)\lvec(40 60)\lvec(40
50)\htext(43 53){$_2$}

\move(40 30)\lvec(50 30)\lvec(50 40)\lvec(40 40)\lvec(40
30)\htext(43 33){$_2$}

\move(40 40)\lvec(50 50)\lvec(40 50)\lvec(40 40) \htext(41
46){\tiny $0$}

\move(40 40)\lvec(50 50)\lvec(50 40)\lvec(40 40) \htext(46
42){\tiny $1$}

\end{texdraw}}\quad$+$\quad $q^4$
\raisebox{-0.4\height}{
\begin{texdraw}
\drawdim em \setunitscale 0.13 \linewd 0.5

\move(30 0)\lvec(40 0)\lvec(40 10)\lvec(30 10)\lvec(30 0)\htext(32
6){\tiny $1$}

\move(40 0)\lvec(50 0)\lvec(50 10)\lvec(40 10)\lvec(40 0)\htext(42
6){\tiny $0$}

\move(40 80)\lvec(50 90)\lvec(40 90)\lvec(40 80)\htext(42
86){\tiny $0$}

\move(30 0)\lvec(40 10)\lvec(40 0)\lvec(30 0)\lfill f:0.8
\htext(36 2){\tiny $0$}

\move(40 0)\lvec(50 10)\lvec(50 0)\lvec(40 0)\lfill f:0.8
\htext(46 2){\tiny $1$}
\move(40 70)\lvec(50 70)\lvec(50 80)\lvec(40 80)\lvec(40
70)\htext(43 73){$_2$}

\move(30 10)\lvec(40 10)\lvec(40 20)\lvec(30 20)\lvec(30
10)\htext(33 13){$_2$}

\move(40 10)\lvec(50 10)\lvec(50 20)\lvec(40 20)\lvec(40
10)\htext(43 13){$_2$}

\move(40 60)\lvec(50 60)\lvec(50 70)\lvec(40 70)\lvec(40
60)\htext(43 66){\tiny $3$}

\move(40 20)\lvec(50 20)\lvec(50 30)\lvec(40 30)\lvec(40
20)\htext(43 26){\tiny $3$}

\move(40 60)\lvec(50 60)\lvec(50 65)\lvec(40 65)\lvec(40 60)
\htext(43 61){\tiny $3$}

\move(40 20)\lvec(50 20)\lvec(50 25)\lvec(40 25)\lvec(40 20)
\htext(43 21){\tiny $3$}

\move(40 50)\lvec(50 50)\lvec(50 60)\lvec(40 60)\lvec(40
50)\htext(43 53){$_2$}

\move(40 30)\lvec(50 30)\lvec(50 40)\lvec(40 40)\lvec(40
30)\htext(43 33){$_2$}

\move(40 40)\lvec(50 50)\lvec(40 50)\lvec(40 40) \htext(41
46){\tiny $0$}

\move(40 40)\lvec(50 50)\lvec(50 40)\lvec(40 40) \htext(46
42){\tiny $1$}

\end{texdraw}}\quad .\vskip 5mm

On the other hand,\vskip 5mm

$G(\ \raisebox{-0.4\height}{
\begin{texdraw}
\drawdim em \setunitscale 0.13 \linewd 0.5

\move(20 0)\lvec(30 0)\lvec(30 10)\lvec(20 10)\lvec(20 0)\htext(22
6){\tiny $0$}

\move(30 0)\lvec(40 0)\lvec(40 10)\lvec(30 10)\lvec(30 0)\htext(32
6){\tiny $1$}

\move(40 0)\lvec(50 0)\lvec(50 10)\lvec(40 10)\lvec(40 0)\htext(42
6){\tiny $0$}

\move(20 0)\lvec(30 10)\lvec(30 0)\lvec(20 0)\lfill f:0.8
\htext(26 2){\tiny $1$}

\move(30 0)\lvec(40 10)\lvec(40 0)\lvec(30 0)\lfill f:0.8
\htext(36 2){\tiny $0$}

\move(40 0)\lvec(50 10)\lvec(50 0)\lvec(40 0)\lfill f:0.8
\htext(46 2){\tiny $1$}
\move(40 50)\lvec(50 50)\lvec(50 60)\lvec(40 60)\lvec(40
50)\htext(43 53){$_2$}

\move(30 10)\lvec(40 10)\lvec(40 20)\lvec(30 20)\lvec(30
10)\htext(33 13){$_2$}

\move(40 10)\lvec(50 10)\lvec(50 20)\lvec(40 20)\lvec(40
10)\htext(43 13){$_2$}

\move(30 20)\lvec(40 20)\lvec(40 30)\lvec(30 30)\lvec(30
20)\htext(33 26){\tiny $3$}

\move(40 20)\lvec(50 20)\lvec(50 30)\lvec(40 30)\lvec(40
20)\htext(43 26){\tiny $3$}

\move(30 20)\lvec(40 20)\lvec(40 25)\lvec(30 25)\lvec(30 20)
\htext(33 21){\tiny $3$}

\move(40 20)\lvec(50 20)\lvec(50 25)\lvec(40 25)\lvec(40 20)
\htext(43 21){\tiny $3$}

\move(30 30)\lvec(40 30)\lvec(40 40)\lvec(30 40)\lvec(30
30)\htext(33 33){$_2$}

\move(40 30)\lvec(50 30)\lvec(50 40)\lvec(40 40)\lvec(40
30)\htext(43 33){$_2$}

\move(40 40)\lvec(50 50)\lvec(40 50)\lvec(40 40) \htext(41
46){\tiny $0$}

\move(40 40)\lvec(50 50)\lvec(50 40)\lvec(40 40) \htext(46
42){\tiny $1$}

\end{texdraw}} \ )=A(\ \raisebox{-0.4\height}{
\begin{texdraw}
\drawdim em \setunitscale 0.13 \linewd 0.5

\move(20 0)\lvec(30 0)\lvec(30 10)\lvec(20 10)\lvec(20 0)\htext(22
6){\tiny $0$}

\move(30 0)\lvec(40 0)\lvec(40 10)\lvec(30 10)\lvec(30 0)\htext(32
6){\tiny $1$}

\move(40 0)\lvec(50 0)\lvec(50 10)\lvec(40 10)\lvec(40 0)\htext(42
6){\tiny $0$}

\move(20 0)\lvec(30 10)\lvec(30 0)\lvec(20 0)\lfill f:0.8
\htext(26 2){\tiny $1$}

\move(30 0)\lvec(40 10)\lvec(40 0)\lvec(30 0)\lfill f:0.8
\htext(36 2){\tiny $0$}

\move(40 0)\lvec(50 10)\lvec(50 0)\lvec(40 0)\lfill f:0.8
\htext(46 2){\tiny $1$}
\move(40 50)\lvec(50 50)\lvec(50 60)\lvec(40 60)\lvec(40
50)\htext(43 53){$_2$}

\move(30 10)\lvec(40 10)\lvec(40 20)\lvec(30 20)\lvec(30
10)\htext(33 13){$_2$}

\move(40 10)\lvec(50 10)\lvec(50 20)\lvec(40 20)\lvec(40
10)\htext(43 13){$_2$}

\move(30 20)\lvec(40 20)\lvec(40 30)\lvec(30 30)\lvec(30
20)\htext(33 26){\tiny $3$}

\move(40 20)\lvec(50 20)\lvec(50 30)\lvec(40 30)\lvec(40
20)\htext(43 26){\tiny $3$}

\move(30 20)\lvec(40 20)\lvec(40 25)\lvec(30 25)\lvec(30 20)
\htext(33 21){\tiny $3$}

\move(40 20)\lvec(50 20)\lvec(50 25)\lvec(40 25)\lvec(40 20)
\htext(43 21){\tiny $3$}

\move(30 30)\lvec(40 30)\lvec(40 40)\lvec(30 40)\lvec(30
30)\htext(33 33){$_2$}

\move(40 30)\lvec(50 30)\lvec(50 40)\lvec(40 40)\lvec(40
30)\htext(43 33){$_2$}

\move(40 40)\lvec(50 50)\lvec(40 50)\lvec(40 40) \htext(41
46){\tiny $0$}

\move(40 40)\lvec(50 50)\lvec(50 40)\lvec(40 40) \htext(46
42){\tiny $1$}
\end{texdraw}} \ )=f_0f_2f_3^{(2)}f_2^{(2)}f_1^{(2)}f_0f_2f_3^{(2)}f_2f_0Y_{\Lambda_0}$
\vskip 5mm

$=$ \raisebox{-0.4\height}{
\begin{texdraw}
\drawdim em \setunitscale 0.13 \linewd 0.5

\move(20 0)\lvec(30 0)\lvec(30 10)\lvec(20 10)\lvec(20 0)\htext(22
6){\tiny $0$}

\move(30 0)\lvec(40 0)\lvec(40 10)\lvec(30 10)\lvec(30 0)\htext(32
6){\tiny $1$}

\move(40 0)\lvec(50 0)\lvec(50 10)\lvec(40 10)\lvec(40 0)\htext(42
6){\tiny $0$}

\move(20 0)\lvec(30 10)\lvec(30 0)\lvec(20 0)\lfill f:0.8
\htext(26 2){\tiny $1$}

\move(30 0)\lvec(40 10)\lvec(40 0)\lvec(30 0)\lfill f:0.8
\htext(36 2){\tiny $0$}

\move(40 0)\lvec(50 10)\lvec(50 0)\lvec(40 0)\lfill f:0.8
\htext(46 2){\tiny $1$}
\move(40 50)\lvec(50 50)\lvec(50 60)\lvec(40 60)\lvec(40
50)\htext(43 53){$_2$}

\move(30 10)\lvec(40 10)\lvec(40 20)\lvec(30 20)\lvec(30
10)\htext(33 13){$_2$}

\move(40 10)\lvec(50 10)\lvec(50 20)\lvec(40 20)\lvec(40
10)\htext(43 13){$_2$}

\move(30 20)\lvec(40 20)\lvec(40 30)\lvec(30 30)\lvec(30
20)\htext(33 26){\tiny $3$}

\move(40 20)\lvec(50 20)\lvec(50 30)\lvec(40 30)\lvec(40
20)\htext(43 26){\tiny $3$}

\move(30 20)\lvec(40 20)\lvec(40 25)\lvec(30 25)\lvec(30 20)
\htext(33 21){\tiny $3$}

\move(40 20)\lvec(50 20)\lvec(50 25)\lvec(40 25)\lvec(40 20)
\htext(43 21){\tiny $3$}

\move(30 30)\lvec(40 30)\lvec(40 40)\lvec(30 40)\lvec(30
30)\htext(33 33){$_2$}

\move(40 30)\lvec(50 30)\lvec(50 40)\lvec(40 40)\lvec(40
30)\htext(43 33){$_2$}

\move(40 40)\lvec(50 50)\lvec(40 50)\lvec(40 40) \htext(41
46){\tiny $0$}

\move(40 40)\lvec(50 50)\lvec(50 40)\lvec(40 40) \htext(46
42){\tiny $1$}
\end{texdraw}}\quad $+$\quad $q^2$
\raisebox{-0.4\height}{
\begin{texdraw}
\drawdim em \setunitscale 0.13 \linewd 0.5

\move(20 0)\lvec(30 0)\lvec(30 10)\lvec(20 10)\lvec(20 0)\htext(22
6){\tiny $0$}

\move(30 0)\lvec(40 0)\lvec(40 10)\lvec(30 10)\lvec(30 0)\htext(32
6){\tiny $1$}

\move(40 0)\lvec(50 0)\lvec(50 10)\lvec(40 10)\lvec(40 0)\htext(42
6){\tiny $0$}

\move(20 0)\lvec(30 10)\lvec(30 0)\lvec(20 0)\lfill f:0.8
\htext(26 2){\tiny $1$}

\move(30 0)\lvec(40 10)\lvec(40 0)\lvec(30 0)\lfill f:0.8
\htext(36 2){\tiny $0$}

\move(40 0)\lvec(50 10)\lvec(50 0)\lvec(40 0)\lfill f:0.8
\htext(46 2){\tiny $1$}
\move(40 70)\lvec(50 70)\lvec(50 80)\lvec(40 80)\lvec(40
70)\htext(43 73){$_2$}

\move(30 10)\lvec(40 10)\lvec(40 20)\lvec(30 20)\lvec(30
10)\htext(33 13){$_2$}

\move(40 10)\lvec(50 10)\lvec(50 20)\lvec(40 20)\lvec(40
10)\htext(43 13){$_2$}

\move(40 60)\lvec(50 60)\lvec(50 70)\lvec(40 70)\lvec(40
60)\htext(43 66){\tiny $3$}

\move(40 20)\lvec(50 20)\lvec(50 30)\lvec(40 30)\lvec(40
20)\htext(43 26){\tiny $3$}

\move(40 60)\lvec(50 60)\lvec(50 65)\lvec(40 65)\lvec(40 60)
\htext(43 61){\tiny $3$}

\move(40 20)\lvec(50 20)\lvec(50 25)\lvec(40 25)\lvec(40 20)
\htext(43 21){\tiny $3$}

\move(40 50)\lvec(50 50)\lvec(50 60)\lvec(40 60)\lvec(40
50)\htext(43 53){$_2$}

\move(40 30)\lvec(50 30)\lvec(50 40)\lvec(40 40)\lvec(40
30)\htext(43 33){$_2$}

\move(40 40)\lvec(50 50)\lvec(40 50)\lvec(40 40) \htext(41
46){\tiny $0$}

\move(40 40)\lvec(50 50)\lvec(50 40)\lvec(40 40) \htext(46
42){\tiny $1$}

\end{texdraw}}\quad$+$\quad $q^2$
\raisebox{-0.4\height}{
\begin{texdraw}
\drawdim em \setunitscale 0.13 \linewd 0.5

\move(30 0)\lvec(40 0)\lvec(40 10)\lvec(30 10)\lvec(30 0)\htext(32
6){\tiny $1$}

\move(40 0)\lvec(50 0)\lvec(50 10)\lvec(40 10)\lvec(40 0)\htext(42
6){\tiny $0$}

\move(30 0)\lvec(40 10)\lvec(40 0)\lvec(30 0)\lfill f:0.8
\htext(36 2){\tiny $0$}

\move(30 40)\lvec(40 50)\lvec(40 40)\lvec(30 40)\htext(36
42){\tiny $0$}

\move(40 0)\lvec(50 10)\lvec(50 0)\lvec(40 0)\lfill f:0.8
\htext(46 2){\tiny $1$}
\move(30 30)\lvec(40 30)\lvec(40 40)\lvec(30 40)\lvec(30
30)\htext(33 33){$_2$}

\move(30 10)\lvec(40 10)\lvec(40 20)\lvec(30 20)\lvec(30
10)\htext(33 13){$_2$}

\move(40 10)\lvec(50 10)\lvec(50 20)\lvec(40 20)\lvec(40
10)\htext(43 13){$_2$}

\move(30 20)\lvec(40 20)\lvec(40 30)\lvec(30 30)\lvec(30
20)\htext(33 26){\tiny $3$}

\move(40 20)\lvec(50 20)\lvec(50 30)\lvec(40 30)\lvec(40
20)\htext(43 26){\tiny $3$}

\move(30 20)\lvec(40 20)\lvec(40 25)\lvec(30 25)\lvec(30 20)
\htext(33 21){\tiny $3$}

\move(40 20)\lvec(50 20)\lvec(50 25)\lvec(40 25)\lvec(40 20)
\htext(43 21){\tiny $3$}

\move(40 50)\lvec(50 50)\lvec(50 60)\lvec(40 60)\lvec(40
50)\htext(43 53){$_2$}

\move(40 30)\lvec(50 30)\lvec(50 40)\lvec(40 40)\lvec(40
30)\htext(43 33){$_2$}

\move(40 40)\lvec(50 50)\lvec(40 50)\lvec(40 40) \htext(41
46){\tiny $0$}

\move(40 40)\lvec(50 50)\lvec(50 40)\lvec(40 40) \htext(46
42){\tiny $1$}

\end{texdraw}}
\quad$+$\quad $q^4$ \raisebox{-0.4\height}{
\begin{texdraw}
\drawdim em \setunitscale 0.13 \linewd 0.5

\move(30 0)\lvec(40 0)\lvec(40 10)\lvec(30 10)\lvec(30 0)\htext(32
6){\tiny $1$}

\move(40 0)\lvec(50 0)\lvec(50 10)\lvec(40 10)\lvec(40 0)\htext(42
6){\tiny $0$}

\move(40 80)\lvec(50 90)\lvec(40 90)\lvec(40 80)\htext(42
86){\tiny $0$}

\move(30 0)\lvec(40 10)\lvec(40 0)\lvec(30 0)\lfill f:0.8
\htext(36 2){\tiny $0$}

\move(40 0)\lvec(50 10)\lvec(50 0)\lvec(40 0)\lfill f:0.8
\htext(46 2){\tiny $1$}
\move(40 70)\lvec(50 70)\lvec(50 80)\lvec(40 80)\lvec(40
70)\htext(43 73){$_2$}

\move(30 10)\lvec(40 10)\lvec(40 20)\lvec(30 20)\lvec(30
10)\htext(33 13){$_2$}

\move(40 10)\lvec(50 10)\lvec(50 20)\lvec(40 20)\lvec(40
10)\htext(43 13){$_2$}

\move(40 60)\lvec(50 60)\lvec(50 70)\lvec(40 70)\lvec(40
60)\htext(43 66){\tiny $3$}

\move(40 20)\lvec(50 20)\lvec(50 30)\lvec(40 30)\lvec(40
20)\htext(43 26){\tiny $3$}

\move(40 60)\lvec(50 60)\lvec(50 65)\lvec(40 65)\lvec(40 60)
\htext(43 61){\tiny $3$}

\move(40 20)\lvec(50 20)\lvec(50 25)\lvec(40 25)\lvec(40 20)
\htext(43 21){\tiny $3$}

\move(40 50)\lvec(50 50)\lvec(50 60)\lvec(40 60)\lvec(40
50)\htext(43 53){$_2$}

\move(40 30)\lvec(50 30)\lvec(50 40)\lvec(40 40)\lvec(40
30)\htext(43 33){$_2$}

\move(40 40)\lvec(50 50)\lvec(40 50)\lvec(40 40) \htext(41
46){\tiny $0$}

\move(40 40)\lvec(50 50)\lvec(50 40)\lvec(40 40) \htext(46
42){\tiny $1$}

\end{texdraw}}\quad .\vskip 5mm

Hence,\vskip 5mm

$G(\ \raisebox{-0.4\height}{
\begin{texdraw}
\drawdim em \setunitscale 0.13 \linewd 0.5

\move(20 0)\lvec(30 0)\lvec(30 10)\lvec(20 10)\lvec(20 0)\htext(22
6){\tiny $0$}

\move(30 0)\lvec(40 0)\lvec(40 10)\lvec(30 10)\lvec(30 0)\htext(32
6){\tiny $1$}

\move(40 0)\lvec(50 0)\lvec(50 10)\lvec(40 10)\lvec(40 0)\htext(42
6){\tiny $0$}

\move(20 0)\lvec(30 10)\lvec(30 0)\lvec(20 0)\lfill f:0.8
\htext(26 2){\tiny $1$}

\move(30 0)\lvec(40 10)\lvec(40 0)\lvec(30 0)\lfill f:0.8
\htext(36 2){\tiny $0$}

\move(40 0)\lvec(50 10)\lvec(50 0)\lvec(40 0)\lfill f:0.8
\htext(46 2){\tiny $1$}
\move(20 10)\lvec(30 10)\lvec(30 20)\lvec(20 20)\lvec(20
10)\htext(23 13){$_2$}

\move(30 10)\lvec(40 10)\lvec(40 20)\lvec(30 20)\lvec(30
10)\htext(33 13){$_2$}

\move(40 10)\lvec(50 10)\lvec(50 20)\lvec(40 20)\lvec(40
10)\htext(43 13){$_2$}

\move(30 20)\lvec(40 20)\lvec(40 30)\lvec(30 30)\lvec(30
20)\htext(33 26){\tiny $3$}

\move(40 20)\lvec(50 20)\lvec(50 30)\lvec(40 30)\lvec(40
20)\htext(43 26){\tiny $3$}

\move(30 20)\lvec(40 20)\lvec(40 25)\lvec(30 25)\lvec(30 20)
\htext(33 21){\tiny $3$}

\move(40 20)\lvec(50 20)\lvec(50 25)\lvec(40 25)\lvec(40 20)
\htext(43 21){\tiny $3$}

\move(30 30)\lvec(40 30)\lvec(40 40)\lvec(30 40)\lvec(30
30)\htext(33 33){$_2$}

\move(40 30)\lvec(50 30)\lvec(50 40)\lvec(40 40)\lvec(40
30)\htext(43 33){$_2$}

\move(40 40)\lvec(50 50)\lvec(40 50)\lvec(40 40) \htext(41
46){\tiny $0$}

\move(40 40)\lvec(50 50)\lvec(50 40)\lvec(40 40) \htext(46
42){\tiny $1$}

\end{texdraw}} \
 ) = A(\ \raisebox{-0.4\height}{
\begin{texdraw}
\drawdim em \setunitscale 0.13 \linewd 0.5

\move(20 0)\lvec(30 0)\lvec(30 10)\lvec(20 10)\lvec(20 0)\htext(22
6){\tiny $0$}

\move(30 0)\lvec(40 0)\lvec(40 10)\lvec(30 10)\lvec(30 0)\htext(32
6){\tiny $1$}

\move(40 0)\lvec(50 0)\lvec(50 10)\lvec(40 10)\lvec(40 0)\htext(42
6){\tiny $0$}

\move(20 0)\lvec(30 10)\lvec(30 0)\lvec(20 0)\lfill f:0.8
\htext(26 2){\tiny $1$}

\move(30 0)\lvec(40 10)\lvec(40 0)\lvec(30 0)\lfill f:0.8
\htext(36 2){\tiny $0$}

\move(40 0)\lvec(50 10)\lvec(50 0)\lvec(40 0)\lfill f:0.8
\htext(46 2){\tiny $1$}
\move(20 10)\lvec(30 10)\lvec(30 20)\lvec(20 20)\lvec(20
10)\htext(23 13){$_2$}

\move(30 10)\lvec(40 10)\lvec(40 20)\lvec(30 20)\lvec(30
10)\htext(33 13){$_2$}

\move(40 10)\lvec(50 10)\lvec(50 20)\lvec(40 20)\lvec(40
10)\htext(43 13){$_2$}

\move(30 20)\lvec(40 20)\lvec(40 30)\lvec(30 30)\lvec(30
20)\htext(33 26){\tiny $3$}

\move(40 20)\lvec(50 20)\lvec(50 30)\lvec(40 30)\lvec(40
20)\htext(43 26){\tiny $3$}

\move(30 20)\lvec(40 20)\lvec(40 25)\lvec(30 25)\lvec(30 20)
\htext(33 21){\tiny $3$}

\move(40 20)\lvec(50 20)\lvec(50 25)\lvec(40 25)\lvec(40 20)
\htext(43 21){\tiny $3$}

\move(30 30)\lvec(40 30)\lvec(40 40)\lvec(30 40)\lvec(30
30)\htext(33 33){$_2$}

\move(40 30)\lvec(50 30)\lvec(50 40)\lvec(40 40)\lvec(40
30)\htext(43 33){$_2$}

\move(40 40)\lvec(50 50)\lvec(40 50)\lvec(40 40) \htext(41
46){\tiny $0$}

\move(40 40)\lvec(50 50)\lvec(50 40)\lvec(40 40) \htext(46
42){\tiny $1$}

\end{texdraw}} \
)-G(\ \raisebox{-0.4\height}{
\begin{texdraw}
\drawdim em \setunitscale 0.13 \linewd 0.5

\move(20 0)\lvec(30 0)\lvec(30 10)\lvec(20 10)\lvec(20 0)\htext(22
6){\tiny $0$}

\move(30 0)\lvec(40 0)\lvec(40 10)\lvec(30 10)\lvec(30 0)\htext(32
6){\tiny $1$}

\move(40 0)\lvec(50 0)\lvec(50 10)\lvec(40 10)\lvec(40 0)\htext(42
6){\tiny $0$}

\move(20 0)\lvec(30 10)\lvec(30 0)\lvec(20 0)\lfill f:0.8
\htext(26 2){\tiny $1$}

\move(30 0)\lvec(40 10)\lvec(40 0)\lvec(30 0)\lfill f:0.8
\htext(36 2){\tiny $0$}

\move(40 0)\lvec(50 10)\lvec(50 0)\lvec(40 0)\lfill f:0.8
\htext(46 2){\tiny $1$}
\move(40 50)\lvec(50 50)\lvec(50 60)\lvec(40 60)\lvec(40
50)\htext(43 53){$_2$}

\move(30 10)\lvec(40 10)\lvec(40 20)\lvec(30 20)\lvec(30
10)\htext(33 13){$_2$}

\move(40 10)\lvec(50 10)\lvec(50 20)\lvec(40 20)\lvec(40
10)\htext(43 13){$_2$}

\move(30 20)\lvec(40 20)\lvec(40 30)\lvec(30 30)\lvec(30
20)\htext(33 26){\tiny $3$}

\move(40 20)\lvec(50 20)\lvec(50 30)\lvec(40 30)\lvec(40
20)\htext(43 26){\tiny $3$}

\move(30 20)\lvec(40 20)\lvec(40 25)\lvec(30 25)\lvec(30 20)
\htext(33 21){\tiny $3$}

\move(40 20)\lvec(50 20)\lvec(50 25)\lvec(40 25)\lvec(40 20)
\htext(43 21){\tiny $3$}

\move(30 30)\lvec(40 30)\lvec(40 40)\lvec(30 40)\lvec(30
30)\htext(33 33){$_2$}

\move(40 30)\lvec(50 30)\lvec(50 40)\lvec(40 40)\lvec(40
30)\htext(43 33){$_2$}

\move(40 40)\lvec(50 50)\lvec(40 50)\lvec(40 40) \htext(41
46){\tiny $0$}

\move(40 40)\lvec(50 50)\lvec(50 40)\lvec(40 40) \htext(46
42){\tiny $1$}

\end{texdraw}} \ )$\vskip 5mm

\hskip 1cm$=$ \raisebox{-0.4\height}{
\begin{texdraw}
\drawdim em \setunitscale 0.13 \linewd 0.5

\move(20 0)\lvec(30 0)\lvec(30 10)\lvec(20 10)\lvec(20 0)\htext(22
6){\tiny $0$}

\move(30 0)\lvec(40 0)\lvec(40 10)\lvec(30 10)\lvec(30 0)\htext(32
6){\tiny $1$}

\move(40 0)\lvec(50 0)\lvec(50 10)\lvec(40 10)\lvec(40 0)\htext(42
6){\tiny $0$}

\move(20 0)\lvec(30 10)\lvec(30 0)\lvec(20 0)\lfill f:0.8
\htext(26 2){\tiny $1$}

\move(30 0)\lvec(40 10)\lvec(40 0)\lvec(30 0)\lfill f:0.8
\htext(36 2){\tiny $0$}

\move(40 0)\lvec(50 10)\lvec(50 0)\lvec(40 0)\lfill f:0.8
\htext(46 2){\tiny $1$}
\move(20 10)\lvec(30 10)\lvec(30 20)\lvec(20 20)\lvec(20
10)\htext(23 13){$_2$}

\move(30 10)\lvec(40 10)\lvec(40 20)\lvec(30 20)\lvec(30
10)\htext(33 13){$_2$}

\move(40 10)\lvec(50 10)\lvec(50 20)\lvec(40 20)\lvec(40
10)\htext(43 13){$_2$}

\move(30 20)\lvec(40 20)\lvec(40 30)\lvec(30 30)\lvec(30
20)\htext(33 26){\tiny $3$}

\move(40 20)\lvec(50 20)\lvec(50 30)\lvec(40 30)\lvec(40
20)\htext(43 26){\tiny $3$}

\move(30 20)\lvec(40 20)\lvec(40 25)\lvec(30 25)\lvec(30 20)
\htext(33 21){\tiny $3$}

\move(40 20)\lvec(50 20)\lvec(50 25)\lvec(40 25)\lvec(40 20)
\htext(43 21){\tiny $3$}

\move(30 30)\lvec(40 30)\lvec(40 40)\lvec(30 40)\lvec(30
30)\htext(33 33){$_2$}

\move(40 30)\lvec(50 30)\lvec(50 40)\lvec(40 40)\lvec(40
30)\htext(43 33){$_2$}

\move(40 40)\lvec(50 50)\lvec(40 50)\lvec(40 40) \htext(41
46){\tiny $0$}

\move(40 40)\lvec(50 50)\lvec(50 40)\lvec(40 40) \htext(46
42){\tiny $1$}

\end{texdraw}}\quad$+$\quad $q^2$
\raisebox{-0.4\height}{
\begin{texdraw}
\drawdim em \setunitscale 0.13 \linewd 0.5

\move(20 0)\lvec(30 0)\lvec(30 10)\lvec(20 10)\lvec(20 0)\htext(22
6){\tiny $0$}

\move(30 0)\lvec(40 0)\lvec(40 10)\lvec(30 10)\lvec(30 0)\htext(32
6){\tiny $1$}

\move(40 0)\lvec(50 0)\lvec(50 10)\lvec(40 10)\lvec(40 0)\htext(42
6){\tiny $0$}

\move(20 0)\lvec(30 10)\lvec(30 0)\lvec(20 0)\lfill f:0.8
\htext(26 2){\tiny $1$}

\move(30 0)\lvec(40 10)\lvec(40 0)\lvec(30 0)\lfill f:0.8
\htext(36 2){\tiny $0$}

\move(40 0)\lvec(50 10)\lvec(50 0)\lvec(40 0)\lfill f:0.8
\htext(46 2){\tiny $1$}
\move(20 10)\lvec(30 10)\lvec(30 20)\lvec(20 20)\lvec(20
10)\htext(23 13){$_2$}

\move(30 10)\lvec(40 10)\lvec(40 20)\lvec(30 20)\lvec(30
10)\htext(33 13){$_2$}

\move(40 10)\lvec(50 10)\lvec(50 20)\lvec(40 20)\lvec(40
10)\htext(43 13){$_2$}

\move(30 20)\lvec(40 20)\lvec(40 30)\lvec(30 30)\lvec(30
20)\htext(33 26){\tiny $3$}

\move(40 20)\lvec(50 20)\lvec(50 30)\lvec(40 30)\lvec(40
20)\htext(43 26){\tiny $3$}

\move(30 20)\lvec(40 20)\lvec(40 25)\lvec(30 25)\lvec(30 20)
\htext(33 21){\tiny $3$}

\move(40 20)\lvec(50 20)\lvec(50 25)\lvec(40 25)\lvec(40 20)
\htext(43 21){\tiny $3$}

\move(40 50)\lvec(50 50)\lvec(50 60)\lvec(40 60)\lvec(40
50)\htext(43 53){$_2$}

\move(40 30)\lvec(50 30)\lvec(50 40)\lvec(40 40)\lvec(40
30)\htext(43 33){$_2$}

\move(40 40)\lvec(50 50)\lvec(40 50)\lvec(40 40) \htext(41
46){\tiny $0$}

\move(40 40)\lvec(50 50)\lvec(50 40)\lvec(40 40) \htext(46
42){\tiny $1$}

\end{texdraw}}\quad $+$\quad $q^4$
\raisebox{-0.4\height}{
\begin{texdraw}
\drawdim em \setunitscale 0.13 \linewd 0.5

\move(20 0)\lvec(30 0)\lvec(30 10)\lvec(20 10)\lvec(20 0)\htext(22
6){\tiny $0$}

\move(30 0)\lvec(40 0)\lvec(40 10)\lvec(30 10)\lvec(30 0)\htext(32
6){\tiny $1$}

\move(40 0)\lvec(50 0)\lvec(50 10)\lvec(40 10)\lvec(40 0)\htext(42
6){\tiny $0$}

\move(20 0)\lvec(30 10)\lvec(30 0)\lvec(20 0)\lfill f:0.8
\htext(26 2){\tiny $1$}

\move(30 0)\lvec(40 10)\lvec(40 0)\lvec(30 0)\lfill f:0.8
\htext(36 2){\tiny $0$}

\move(40 0)\lvec(50 10)\lvec(50 0)\lvec(40 0)\lfill f:0.8
\htext(46 2){\tiny $1$}
\move(30 30)\lvec(40 30)\lvec(40 40)\lvec(30 40)\lvec(30
30)\htext(33 33){$_2$}

\move(30 10)\lvec(40 10)\lvec(40 20)\lvec(30 20)\lvec(30
10)\htext(33 13){$_2$}

\move(40 10)\lvec(50 10)\lvec(50 20)\lvec(40 20)\lvec(40
10)\htext(43 13){$_2$}

\move(30 20)\lvec(40 20)\lvec(40 30)\lvec(30 30)\lvec(30
20)\htext(33 26){\tiny $3$}

\move(40 20)\lvec(50 20)\lvec(50 30)\lvec(40 30)\lvec(40
20)\htext(43 26){\tiny $3$}

\move(30 20)\lvec(40 20)\lvec(40 25)\lvec(30 25)\lvec(30 20)
\htext(33 21){\tiny $3$}

\move(40 20)\lvec(50 20)\lvec(50 25)\lvec(40 25)\lvec(40 20)
\htext(43 21){\tiny $3$}

\move(40 50)\lvec(50 50)\lvec(50 60)\lvec(40 60)\lvec(40
50)\htext(43 53){$_2$}

\move(40 30)\lvec(50 30)\lvec(50 40)\lvec(40 40)\lvec(40
30)\htext(43 33){$_2$}

\move(40 40)\lvec(50 50)\lvec(40 50)\lvec(40 40) \htext(41
46){\tiny $0$}

\move(40 40)\lvec(50 50)\lvec(50 40)\lvec(40 40) \htext(46
42){\tiny $1$}

\end{texdraw}}\quad .\vskip 5mm
}
\end{ex}
\vskip 1cm

\section{Appendix}
In this section, we list the patterns for building the walls which
are given in \cite{Ka2000}.

(a) $A_n^{(1)}$ ($n\geq 1$),

\begin{center}
On $Y_{\Lambda_i}$: \raisebox{-1\height}{\begin{texdraw} \textref
h:C v:C \fontsize{8}{8}\selectfont \drawdim em \setunitscale 1.9
\move(0 0)\rlvec(-7.7 0) \move(0 1)\rlvec(-7.7 0) \move(0
2)\rlvec(-7.7 0) \move(0 3.5)\rlvec(-7.7 0) \move(0
4.5)\rlvec(-7.7 0) \move(0 5.5)\rlvec(-7.7 0) \move(0
6.5)\rlvec(-7.7 0) \move(0 0)\rlvec(0 6.7) \move(-1 0)\rlvec(0
6.7) \move(-2 0)\rlvec(0 6.7) \move(-3.5 0)\rlvec(0 6.7)
\move(-4.5 0)\rlvec(0 6.7) \move(-5.5 0)\rlvec(0 6.7) \move(-6.5
0)\rlvec(0 6.7) \move(-7.5 0)\rlvec(0 6.7) \move(-0.5 0.5)
\bsegment \htext(0 0){$i$} \htext(0 1){$i\!\!+\!\!1$} \vtext(0
2.25){$\cdots$} \htext(0 3.5){$n$} \htext(0 4.5){$0$} \htext(0
5.5){$1$} \htext(-1 0){$i\!\!-\!\!1$} \htext(-1 1){$i$} \vtext(-1
2.25){$\cdots$} \htext(-1 3.5){$n\!\!-\!\!1$} \htext(-1 4.5){$n$}
\htext(-1 5.5){$0$} \htext(-2.25 0){$\cdots$} \htext(-2.25
1){$\cdots$} \htext(-3.5 0){$2$} \htext(-3.5 1){$3$} \htext(-4.5
0){$1$} \htext(-4.5 1){$2$} \htext(-5.5 0){$0$} \htext(-5.5
1){$1$} \htext(-6.5 0){$n$} \htext(-6.5 1){$0$} \esegment
\end{texdraw}}
\end{center}%

(b) $A_{2n-1}^{(2)}$ ($n\geq 3$),

\begin{center}
On $Y_{\Lambda_0}$ : \raisebox{-1\height}{\begin{texdraw} \textref
h:C v:C \fontsize{6}{6}\selectfont \drawdim mm \setunitscale 5
\nc{\dtri}{ \bsegment \move(-1 0)\lvec(0 1)\lvec(0 0)\lvec(-1
0)\ifill f:0.7 \esegment } \move(0 0)\dtri \move(-1 0)\dtri
\move(-2 0)\dtri \move(-3 0)\dtri \move(0 0)\rlvec(-4.3 0) \move(0
1)\rlvec(-4.3 0) \move(0 2)\rlvec(-4.3 0) \move(0 3.5)\rlvec(-4.3
0) \move(0 4.5)\rlvec(-4.3 0) \move(0 6)\rlvec(-4.3 0) \move(0
7)\rlvec(-4.3 0) \move(0 8)\rlvec(-4.3 0) \move(0 9)\rlvec(-4.3 0)
\move(0 0)\rlvec(0 9.3) \move(-1 0)\rlvec(0 9.3) \move(-2
0)\rlvec(0 9.3) \move(-3 0)\rlvec(0 9.3) \move(-4 0)\rlvec(0 9.3)
\move(-1 0)\rlvec(1 1) \move(-2 0)\rlvec(1 1) \move(-3 0)\rlvec(1
1) \move(-4 0)\rlvec(1 1) \move(-1 7)\rlvec(1 1) \move(-2
7)\rlvec(1 1) \move(-3 7)\rlvec(1 1) \move(-4 7)\rlvec(1 1)
\vtext(-0.5 2.75){$\cdots$} \vtext(-0.5 5.25){$\cdots$}
\vtext(-1.5 2.75){$\cdots$} \vtext(-1.5 5.25){$\cdots$}
\vtext(-2.5 2.75){$\cdots$} \vtext(-2.5 5.25){$\cdots$}
\vtext(-3.5 2.75){$\cdots$} \vtext(-3.5 5.25){$\cdots$}
\htext(-0.25 7.27){$1$} \htext(-0.75 7.75){$0$} \htext(-1.25
7.27){$0$} \htext(-1.75 7.75){$1$} \htext(-2.25 7.27){$1$}
\htext(-2.75 7.75){$0$} \htext(-3.25 7.27){$0$} \htext(-3.75
7.75){$1$} \htext(-0.25 0.27){$1$} \htext(-0.75 0.75){$0$}
\htext(-1.25 0.27){$0$} \htext(-1.75 0.75){$1$} \htext(-2.25
0.27){$1$} \htext(-2.75 0.75){$0$} \htext(-3.25 0.27){$0$}
\htext(-3.75 0.75){$1$} \htext(-0.5 1.5){$2$} \htext(-1.5
1.5){$2$} \htext(-2.5 1.5){$2$} \htext(-3.5 1.5){$2$} \htext(-0.5
6.5){$2$} \htext(-1.5 6.5){$2$} \htext(-2.5 6.5){$2$} \htext(-3.5
6.5){$2$} \htext(-0.5 8.5){$2$} \htext(-1.5 8.5){$2$} \htext(-2.5
8.5){$2$} \htext(-3.5 8.5){$2$} \htext(-0.5 4){$n$} \htext(-1.5
4){$n$} \htext(-2.5 4){$n$} \htext(-3.5 4){$n$}
\end{texdraw}}\hskip 2cm
On $Y_{\Lambda_1}$ : \raisebox{-1\height}{\begin{texdraw} \textref
h:C v:C \fontsize{6}{6}\selectfont \drawdim mm \setunitscale 5
\nc{\dtri}{ \bsegment \move(-1 0)\lvec(0 1)\lvec(0 0)\lvec(-1
0)\ifill f:0.7 \esegment } \move(0 0)\dtri \move(-1 0)\dtri
\move(-2 0)\dtri \move(-3 0)\dtri \move(0 0)\rlvec(-4.3 0) \move(0
1)\rlvec(-4.3 0) \move(0 2)\rlvec(-4.3 0) \move(0 3.5)\rlvec(-4.3
0) \move(0 4.5)\rlvec(-4.3 0) \move(0 6)\rlvec(-4.3 0) \move(0
7)\rlvec(-4.3 0) \move(0 8)\rlvec(-4.3 0) \move(0 9)\rlvec(-4.3 0)
\move(0 0)\rlvec(0 9.3) \move(-1 0)\rlvec(0 9.3) \move(-2
0)\rlvec(0 9.3) \move(-3 0)\rlvec(0 9.3) \move(-4 0)\rlvec(0 9.3)
\move(-1 0)\rlvec(1 1) \move(-2 0)\rlvec(1 1) \move(-3 0)\rlvec(1
1) \move(-4 0)\rlvec(1 1) \move(-1 7)\rlvec(1 1) \move(-2
7)\rlvec(1 1) \move(-3 7)\rlvec(1 1) \move(-4 7)\rlvec(1 1)
\vtext(-0.5 2.75){$\cdots$} \vtext(-0.5 5.25){$\cdots$}
\vtext(-1.5 2.75){$\cdots$} \vtext(-1.5 5.25){$\cdots$}
\vtext(-2.5 2.75){$\cdots$} \vtext(-2.5 5.25){$\cdots$}
\vtext(-3.5 2.75){$\cdots$} \vtext(-3.5 5.25){$\cdots$}
\htext(-0.25 7.27){$0$} \htext(-0.75 7.75){$1$} \htext(-1.25
7.27){$1$} \htext(-1.75 7.75){$0$} \htext(-2.25 7.27){$0$}
\htext(-2.75 7.75){$1$} \htext(-3.25 7.27){$1$} \htext(-3.75
7.75){$0$} \htext(-0.25 0.27){$0$} \htext(-0.75 0.75){$1$}
\htext(-1.25 0.27){$1$} \htext(-1.75 0.75){$0$} \htext(-2.25
0.27){$0$} \htext(-2.75 0.75){$1$} \htext(-3.25 0.27){$1$}
\htext(-3.75 0.75){$0$} \htext(-0.5 1.5){$2$} \htext(-1.5
1.5){$2$} \htext(-2.5 1.5){$2$} \htext(-3.5 1.5){$2$} \htext(-0.5
6.5){$2$} \htext(-1.5 6.5){$2$} \htext(-2.5 6.5){$2$} \htext(-3.5
6.5){$2$} \htext(-0.5 8.5){$2$} \htext(-1.5 8.5){$2$} \htext(-2.5
8.5){$2$} \htext(-3.5 8.5){$2$} \htext(-0.5 4){$n$} \htext(-1.5
4){$n$} \htext(-2.5 4){$n$} \htext(-3.5 4){$n$}
\end{texdraw}}
\end{center}%

(c) $D_n^{(1)}$ ($n\geq 4$),
\begin{center}
On $Y_{\Lambda_0}$ : \raisebox{-1\height}{\begin{texdraw}\textref
h:C v:C \fontsize{6}{6}\selectfont \drawdim mm \setunitscale 5
\nc{\dtri}{ \bsegment \move(-1 0)\lvec(0 1)\lvec(0 0)\lvec(-1
0)\ifill f:0.7 \esegment } \move(0 0)\dtri \move(-1 0)\dtri
\move(-2 0)\dtri \move(-3 0)\dtri \move(0 0)\rlvec(-4.3 0) \move(0
1)\rlvec(-4.3 0) \move(0 2)\rlvec(-4.3 0) \move(0 3.5)\rlvec(-4.3
0) \move(0 4.5)\rlvec(-4.3 0) \move(0 5.5)\rlvec(-4.3 0) \move(0
6.5)\rlvec(-4.3 0) \move(0 8)\rlvec(-4.3 0) \move(0 9)\rlvec(-4.3
0) \move(0 10)\rlvec(-4.3 0) \move(0 11)\rlvec(-4.3 0) \move(0
0)\rlvec(0 11.3) \move(-1 0)\rlvec(0 11.3) \move(-2 0)\rlvec(0
11.3) \move(-3 0)\rlvec(0 11.3) \move(-4 0)\rlvec(0 11.3) \move(-1
0)\rlvec(1 1) \move(-2 0)\rlvec(1 1) \move(-3 0)\rlvec(1 1)
\move(-4 0)\rlvec(1 1) \move(-1 9)\rlvec(1 1) \move(-2 9)\rlvec(1
1) \move(-3 9)\rlvec(1 1) \move(-4 9)\rlvec(1 1) \htext(-0.3
0.25){$1$} \htext(-0.75 0.75){$0$} \htext(-0.5 1.5){$2$}
\vtext(-0.5 2.75){$\cdots$} \htext(-0.5 4){$n\!\!-\!\!2$}
\htext(-0.5 6){$n\!\!-\!\!2$} \htext(-0.5 8.5){$2$} \htext(-0.3
9.25){$1$} \htext(-0.75 9.75){$0$} \htext(-0.5 10.5){$2$}
\htext(-2.3 0.25){$1$} \htext(-2.75 0.75){$0$} \htext(-2.5
1.5){$2$} \vtext(-2.5 2.75){$\cdots$} \htext(-2.5
4){$n\!\!-\!\!2$} \htext(-2.5 6){$n\!\!-\!\!2$} \htext(-2.5
8.5){$2$} \htext(-2.3 9.25){$1$} \htext(-2.75 9.75){$0$}
\htext(-2.5 10.5){$2$} \htext(-1.3 0.25){$0$} \htext(-1.75
0.75){$1$} \htext(-1.5 1.5){$2$} \vtext(-1.5 2.75){$\cdots$}
\htext(-1.5 4){$n\!\!-\!\!2$} \htext(-1.5 6){$n\!\!-\!\!2$}
\htext(-1.5 8.5){$2$} \htext(-1.3 9.25){$0$} \htext(-1.75
9.75){$1$} \htext(-1.5 10.5){$2$} \htext(-3.3 0.25){$0$}
\htext(-3.75 0.75){$1$} \htext(-3.5 1.5){$2$} \vtext(-3.5
2.75){$\cdots$} \htext(-3.5 4){$n\!\!-\!\!2$} \htext(-3.5
6){$n\!\!-\!\!2$} \htext(-3.5 8.5){$2$} \htext(-3.3 9.25){$0$}
\htext(-3.75 9.75){$1$} \htext(-3.5 10.5){$2$} \htext(-0.3
4.75){$n$} \htext(-2.3 4.75){$n$} \htext(-1.75 5.25){$n$}
\htext(-3.75 5.25){$n$} \move(-1 4.5)\rlvec(0.5 0.5)\rmove(0.45
0.45)\rlvec(0.05 0.05) \move(-2 4.5)\rlvec(0.05 0.05)\rmove(0.45
0.45)\rlvec(0.5 0.5) \move(-3 4.5)\rlvec(0.5 0.5)\rmove(0.45
0.45)\rlvec(0.05 0.05) \move(-4 4.5)\rlvec(0.05 0.05)\rmove(0.45
0.45)\rlvec(0.5 0.5) \htext(-1.5 4.75){$n\!\!-\!\!1$} \htext(-3.5
4.75){$n\!\!-\!\!1$} \htext(-0.5 5.25){$n\!\!-\!\!1$} \htext(-2.5
5.25){$n\!\!-\!\!1$} \vtext(-0.5 7.25){$\cdots$} \vtext(-1.5
7.25){$\cdots$} \vtext(-2.5 7.25){$\cdots$} \vtext(-3.5
7.25){$\cdots$}
\end{texdraw}}\hskip 2cm
On $Y_{\Lambda_1}$ : \raisebox{-1\height}{
\begin{texdraw} \textref h:C v:C \fontsize{6}{6}\selectfont \drawdim mm
\setunitscale 5 \nc{\dtri}{ \bsegment \move(-1 0)\lvec(0 1)\lvec(0
0)\lvec(-1 0)\ifill f:0.7 \esegment } \move(0 0)\dtri \move(-1
0)\dtri \move(-2 0)\dtri \move(-3 0)\dtri \move(0 0)\rlvec(-4.3 0)
\move(0 1)\rlvec(-4.3 0) \move(0 2)\rlvec(-4.3 0) \move(0
3.5)\rlvec(-4.3 0) \move(0 4.5)\rlvec(-4.3 0) \move(0
5.5)\rlvec(-4.3 0) \move(0 6.5)\rlvec(-4.3 0) \move(0
8)\rlvec(-4.3 0) \move(0 9)\rlvec(-4.3 0) \move(0 10)\rlvec(-4.3
0) \move(0 11)\rlvec(-4.3 0) \move(0 0)\rlvec(0 11.3) \move(-1
0)\rlvec(0 11.3) \move(-2 0)\rlvec(0 11.3) \move(-3 0)\rlvec(0
11.3) \move(-4 0)\rlvec(0 11.3) \move(-1 0)\rlvec(1 1) \move(-2
0)\rlvec(1 1) \move(-3 0)\rlvec(1 1) \move(-4 0)\rlvec(1 1)
\move(-1 9)\rlvec(1 1) \move(-2 9)\rlvec(1 1) \move(-3 9)\rlvec(1
1) \move(-4 9)\rlvec(1 1) \htext(-0.3 0.25){$0$} \htext(-0.75
0.75){$1$} \htext(-0.5 1.5){$2$} \vtext(-0.5 2.75){$\cdots$}
\htext(-0.5 4){$n\!\!-\!\!2$} \htext(-0.5 6){$n\!\!-\!\!2$}
\htext(-0.5 8.5){$2$} \htext(-0.3 9.25){$0$} \htext(-0.75
9.75){$1$} \htext(-0.5 10.5){$2$} \htext(-2.3 0.25){$0$}
\htext(-2.75 0.75){$1$} \htext(-2.5 1.5){$2$} \vtext(-2.5
2.75){$\cdots$} \htext(-2.5 4){$n\!\!-\!\!2$} \htext(-2.5
6){$n\!\!-\!\!2$} \htext(-2.5 8.5){$2$} \htext(-2.3 9.25){$0$}
\htext(-2.75 9.75){$1$} \htext(-2.5 10.5){$2$} \htext(-1.3
0.25){$1$} \htext(-1.75 0.75){$0$} \htext(-1.5 1.5){$2$}
\vtext(-1.5 2.75){$\cdots$} \htext(-1.5 4){$n\!\!-\!\!2$}
\htext(-1.5 6){$n\!\!-\!\!2$} \htext(-1.5 8.5){$2$} \htext(-1.3
9.25){$1$} \htext(-1.75 9.75){$0$} \htext(-1.5 10.5){$2$}
\htext(-3.3 0.25){$1$} \htext(-3.75 0.75){$0$} \htext(-3.5
1.5){$2$} \vtext(-3.5 2.75){$\cdots$} \htext(-3.5
4){$n\!\!-\!\!2$} \htext(-3.5 6){$n\!\!-\!\!2$} \htext(-3.5
8.5){$2$} \htext(-3.3 9.25){$1$} \htext(-3.75 9.75){$0$}
\htext(-3.5 10.5){$2$} \htext(-0.3 4.75){$n$} \htext(-2.3
4.75){$n$} \htext(-1.75 5.25){$n$} \htext(-3.75 5.25){$n$}
\move(-1 4.5)\rlvec(0.5 0.5)\rmove(0.45 0.45)\rlvec(0.05 0.05)
\move(-2 4.5)\rlvec(0.05 0.05)\rmove(0.45 0.45)\rlvec(0.5 0.5)
\move(-3 4.5)\rlvec(0.5 0.5)\rmove(0.45 0.45)\rlvec(0.05 0.05)
\move(-4 4.5)\rlvec(0.05 0.05)\rmove(0.45 0.45)\rlvec(0.5 0.5)
\htext(-1.5 4.75){$n\!\!-\!\!1$} \htext(-3.5 4.75){$n\!\!-\!\!1$}
\htext(-0.5 5.25){$n\!\!-\!\!1$} \htext(-2.5 5.25){$n\!\!-\!\!1$}
\vtext(-0.5 7.25){$\cdots$} \vtext(-1.5 7.25){$\cdots$}
\vtext(-2.5 7.25){$\cdots$} \vtext(-3.5 7.25){$\cdots$}
\end{texdraw}}
\end{center}
\begin{center}
On $Y_{\Lambda_{n-1}}$ : \raisebox{-1\height}{\begin{texdraw}
\textref h:C v:C \fontsize{6}{6}\selectfont \drawdim mm
\setunitscale 5 \nc{\dtri}{ \bsegment \move(-1 0)\lvec(0 1)\lvec(0
0)\lvec(-1 0)\ifill f:0.7 \esegment } \move(0 0)\dtri \move(-1
0)\dtri \move(-2 0)\dtri \move(-3 0)\dtri \move(0 0)\rlvec(-4.3 0)
\move(0 1)\rlvec(-4.3 0) \move(0 2)\rlvec(-4.3 0) \move(0
3.5)\rlvec(-4.3 0) \move(0 4.5)\rlvec(-4.3 0) \move(0
5.5)\rlvec(-4.3 0) \move(0 6.5)\rlvec(-4.3 0) \move(0
8)\rlvec(-4.3 0) \move(0 9)\rlvec(-4.3 0) \move(0 10)\rlvec(-4.3
0) \move(0 11)\rlvec(-4.3 0) \move(0 0)\rlvec(0 11.3) \move(-1
0)\rlvec(0 11.3) \move(-2 0)\rlvec(0 11.3) \move(-3 0)\rlvec(0
11.3) \move(-4 0)\rlvec(0 11.3) \move(-1 4.5)\rlvec(1 1) \move(-2
4.5)\rlvec(1 1) \move(-3 4.5)\rlvec(1 1) \move(-4 4.5)\rlvec(1 1)
\htext(-0.3 0.25){$n$} \htext(-0.5 1.5){$n\!\!-\!\!2$} \vtext(-0.5
2.75){$\cdots$} \htext(-0.5 4){$2$} \htext(-0.3 4.75){$1$}
\htext(-0.75 5.25){$0$} \htext(-0.5 6){$2$} \vtext(-0.5
7.25){$\cdots$} \htext(-0.5 8.5){$n\!\!-\!\!2$} \htext(-0.3
9.25){$n$} \htext(-0.5 10.5){$n\!\!-\!\!2$} \htext(-2.3 0.25){$n$}
\htext(-2.5 1.5){$n\!\!-\!\!2$} \vtext(-2.5 2.75){$\cdots$}
\htext(-2.5 4){$2$} \htext(-2.3 4.75){$1$} \htext(-2.75 5.25){$0$}
\htext(-2.5 6){$2$} \vtext(-2.5 7.25){$\cdots$} \htext(-2.5
8.5){$n\!\!-\!\!2$} \htext(-2.3 9.25){$n$} \htext(-2.5
10.5){$n\!\!-\!\!2$} \htext(-1.75 0.75){$n$} \htext(-1.5
1.5){$n\!\!-\!\!2$} \vtext(-1.5 2.75){$\cdots$} \htext(-1.5
4){$2$} \htext(-1.3 4.75){$0$} \htext(-1.75 5.25){$1$} \htext(-1.5
6){$2$} \vtext(-1.5 7.25){$\cdots$} \htext(-1.5
8.5){$n\!\!-\!\!2$} \htext(-1.75 9.75){$n$} \htext(-1.5
10.5){$n\!\!-\!\!2$} \htext(-3.75 0.75){$n$} \htext(-3.5
1.5){$n\!\!-\!\!2$} \vtext(-3.5 2.75){$\cdots$} \htext(-3.5
4){$2$} \htext(-3.3 4.75){$0$} \htext(-3.75 5.25){$1$} \htext(-3.5
6){$2$} \vtext(-3.5 7.25){$\cdots$} \htext(-3.5
8.5){$n\!\!-\!\!2$} \htext(-3.75 9.75){$n$} \htext(-3.5
10.5){$n\!\!-\!\!2$} \move(-1 0)\rlvec(0.5 0.5)\rmove(0.45
0.45)\rlvec(0.05 0.05) \move(-2 0)\rlvec(0.05 0.05)\rmove(0.45
0.45)\rlvec(0.5 0.5) \move(-3 0)\rlvec(0.5 0.5)\rmove(0.45
0.45)\rlvec(0.05 0.05) \move(-4 0)\rlvec(0.05 0.05)\rmove(0.45
0.45)\rlvec(0.5 0.5) \move(-1 9)\rlvec(0.5 0.5)\rmove(0.45
0.45)\rlvec(0.05 0.05) \move(-2 9)\rlvec(0.05 0.05)\rmove(0.45
0.45)\rlvec(0.5 0.5) \move(-3 9)\rlvec(0.5 0.5)\rmove(0.45
0.45)\rlvec(0.05 0.05) \move(-4 9)\rlvec(0.05 0.05)\rmove(0.45
0.45)\rlvec(0.5 0.5) \htext(-3.5 9.25){$n\!\!-\!\!1$} \htext(-3.5
0.25){$n\!\!-\!\!1$} \htext(-2.5 9.75){$n\!\!-\!\!1$} \htext(-2.5
0.75){$n\!\!-\!\!1$} \htext(-1.5 0.25){$n\!\!-\!\!1$} \htext(-1.5
9.25){$n\!\!-\!\!1$} \htext(-0.5 9.75){$n\!\!-\!\!1$} \htext(-0.5
0.75){$n\!\!-\!\!1$}
\end{texdraw}}\hskip 2cm
On $Y_{\Lambda_n}$ : \raisebox{-1\height}{\begin{texdraw} \textref
h:C v:C \fontsize{6}{6}\selectfont \drawdim mm \setunitscale 5
\nc{\dtri}{ \bsegment \move(-1 0)\lvec(0 1)\lvec(0 0)\lvec(-1
0)\ifill f:0.7 \esegment } \move(0 0)\dtri \move(-1 0)\dtri
\move(-2 0)\dtri \move(-3 0)\dtri \move(0 0)\rlvec(-4.3 0) \move(0
1)\rlvec(-4.3 0) \move(0 2)\rlvec(-4.3 0) \move(0 3.5)\rlvec(-4.3
0) \move(0 4.5)\rlvec(-4.3 0) \move(0 5.5)\rlvec(-4.3 0) \move(0
6.5)\rlvec(-4.3 0) \move(0 8)\rlvec(-4.3 0) \move(0 9)\rlvec(-4.3
0) \move(0 10)\rlvec(-4.3 0) \move(0 11)\rlvec(-4.3 0) \move(0
0)\rlvec(0 11.3) \move(-1 0)\rlvec(0 11.3) \move(-2 0)\rlvec(0
11.3) \move(-3 0)\rlvec(0 11.3) \move(-4 0)\rlvec(0 11.3) \move(-1
4.5)\rlvec(1 1) \move(-2 4.5)\rlvec(1 1) \move(-3 4.5)\rlvec(1 1)
\move(-4 4.5)\rlvec(1 1) \htext(-0.75 0.75){$n$} \htext(-0.5
1.5){$n\!\!-\!\!2$} \vtext(-0.5 2.75){$\cdots$} \htext(-0.5
4){$2$} \htext(-0.3 4.75){$1$} \htext(-0.75 5.25){$0$} \htext(-0.5
6){$2$} \vtext(-0.5 7.25){$\cdots$} \htext(-0.5
8.5){$n\!\!-\!\!2$} \htext(-0.75 9.75){$n$} \htext(-0.5
10.5){$n\!\!-\!\!2$} \htext(-2.75 0.75){$n$} \htext(-2.5
1.5){$n\!\!-\!\!2$} \vtext(-2.5 2.75){$\cdots$} \htext(-2.5
4){$2$} \htext(-2.3 4.75){$1$} \htext(-2.75 5.25){$0$} \htext(-2.5
6){$2$} \vtext(-2.5 7.25){$\cdots$} \htext(-2.5
8.5){$n\!\!-\!\!2$} \htext(-2.75 9.75){$n$} \htext(-2.5
10.5){$n\!\!-\!\!2$} \htext(-1.3 0.25){$n$} \htext(-1.5
1.5){$n\!\!-\!\!2$} \vtext(-1.5 2.75){$\cdots$} \htext(-1.5
4){$2$} \htext(-1.3 4.75){$0$} \htext(-1.75 5.25){$1$} \htext(-1.5
6){$2$} \vtext(-1.5 7.25){$\cdots$} \htext(-1.5
8.5){$n\!\!-\!\!2$} \htext(-1.3 9.25){$n$} \htext(-1.5
10.5){$n\!\!-\!\!2$} \htext(-3.3 0.25){$n$} \htext(-3.5
1.5){$n\!\!-\!\!2$} \vtext(-3.5 2.75){$\cdots$} \htext(-3.5
4){$2$} \htext(-3.3 4.75){$0$} \htext(-3.75 5.25){$1$} \htext(-3.5
6){$2$} \vtext(-3.5 7.25){$\cdots$} \htext(-3.5
8.5){$n\!\!-\!\!2$} \htext(-3.3 9.25){$n$} \htext(-3.5
10.5){$n\!\!-\!\!2$} \move(-2 0)\rlvec(0.5 0.5)\rmove(0.45
0.45)\rlvec(0.05 0.05) \move(-1 0)\rlvec(0.05 0.05)\rmove(0.45
0.45)\rlvec(0.5 0.5) \move(-4 0)\rlvec(0.5 0.5)\rmove(0.45
0.45)\rlvec(0.05 0.05) \move(-3 0)\rlvec(0.05 0.05)\rmove(0.45
0.45)\rlvec(0.5 0.5) \move(-2 9)\rlvec(0.5 0.5)\rmove(0.45
0.45)\rlvec(0.05 0.05) \move(-1 9)\rlvec(0.05 0.05)\rmove(0.45
0.45)\rlvec(0.5 0.5) \move(-4 9)\rlvec(0.5 0.5)\rmove(0.45
0.45)\rlvec(0.05 0.05) \move(-3 9)\rlvec(0.05 0.05)\rmove(0.45
0.45)\rlvec(0.5 0.5) \htext(-0.5 0.25){$n\!\!-\!\!1$} \htext(-0.5
9.25){$n\!\!-\!\!1$} \htext(-2.5 0.25){$n\!\!-\!\!1$} \htext(-2.5
9.25){$n\!\!-\!\!1$} \htext(-1.5 0.75){$n\!\!-\!\!1$} \htext(-1.5
9.75){$n\!\!-\!\!1$} \htext(-3.5 0.75){$n\!\!-\!\!1$} \htext(-3.5
9.75){$n\!\!-\!\!1$}
\end{texdraw}}
\end{center}

(d) $A_{2n}^{(2)}$ ($n\geq 1$)
\begin{center}
On $Y_{\Lambda_0}$ : \raisebox{-1\height}{
\begin{texdraw} \textref h:C v:C \fontsize{6}{6}\selectfont \drawdim mm
\setunitscale 5 \move(0 0)\lvec(-4 0)\lvec(-4 0.5)\lvec(0
0.5)\ifill f:0.7 \move(0 0)\rlvec(-4.3 0) \move(0 0.5)\rlvec(-4.3
0) \move(0 1)\rlvec(-4.3 0) \move(0 2)\rlvec(-4.3 0) \move(0
3.5)\rlvec(-4.3 0) \move(0 4.5)\rlvec(-4.3 0) \move(0
6)\rlvec(-4.3 0) \move(0 7)\rlvec(-4.3 0) \move(0 7.5)\rlvec(-4.3
0) \move(0 8)\rlvec(-4.3 0) \move(0 9)\rlvec(-4.3 0) \move(0
0)\rlvec(0 9.3) \move(-1 0)\rlvec(0 9.3) \move(-2 0)\rlvec(0 9.3)
\move(-3 0)\rlvec(0 9.3) \move(-4 0)\rlvec(0 9.3) \htext(-0.5
0.25){$0$} \htext(-1.5 0.25){$0$} \htext(-2.5 0.25){$0$}
\htext(-3.5 0.25){$0$} \htext(-0.5 0.75){$0$} \htext(-1.5
0.75){$0$} \htext(-2.5 0.75){$0$} \htext(-3.5 0.75){$0$}
\htext(-0.5 1.5){$1$} \htext(-1.5 1.5){$1$} \htext(-2.5 1.5){$1$}
\htext(-3.5 1.5){$1$} \htext(-0.5 4){$n$} \htext(-1.5 4){$n$}
\htext(-2.5 4){$n$} \htext(-3.5 4){$n$} \htext(-0.5 6.5){$1$}
\htext(-1.5 6.5){$1$} \htext(-2.5 6.5){$1$} \htext(-3.5 6.5){$1$}
\htext(-0.5 7.25){$0$} \htext(-1.5 7.25){$0$} \htext(-2.5
7.25){$0$} \htext(-3.5 7.25){$0$} \htext(-0.5 7.75){$0$}
\htext(-1.5 7.75){$0$} \htext(-2.5 7.75){$0$} \htext(-3.5
7.75){$0$} \htext(-0.5 8.5){$1$} \htext(-1.5 8.5){$1$} \htext(-2.5
8.5){$1$} \htext(-3.5 8.5){$1$} \vtext(-0.5 2.75){$\cdots$}
\vtext(-1.5 2.75){$\cdots$} \vtext(-2.5 2.75){$\cdots$}
\vtext(-3.5 2.75){$\cdots$} \vtext(-0.5 5.25){$\cdots$}
\vtext(-1.5 5.25){$\cdots$} \vtext(-2.5 5.25){$\cdots$}
\vtext(-3.5 5.25){$\cdots$}
\end{texdraw}}
\end{center}

(e) $D_{n+1}^{(2)}$ ($n\geq 2$),
\begin{center}
On $Y_{\Lambda_0}$ : \raisebox{-1\height}{
\begin{texdraw} \textref h:C v:C \fontsize{6}{6}\selectfont \drawdim mm
\setunitscale 5 \move(0 0)\lvec(-4 0)\lvec(-4 0.5)\lvec(0
0.5)\ifill f:0.7 \move(0 0)\rlvec(-4.3 0) \move(0 0.5)\rlvec(-4.3
0) \move(0 1)\rlvec(-4.3 0) \move(0 2)\rlvec(-4.3 0) \move(0
3.5)\rlvec(-4.3 0) \move(0 4.5)\rlvec(-4.3 0) \move(0
4.5)\rlvec(-4.3 0) \move(0 5)\rlvec(-4.3 0) \move(0
5.5)\rlvec(-4.3 0) \move(0 6.5)\rlvec(-4.3 0) \move(0
8)\rlvec(-4.3 0) \move(0 9)\rlvec(-4.3 0) \move(0 9.5)\rlvec(-4.3
0) \move(0 10)\rlvec(-4.3 0) \move(0 11)\rlvec(-4.3 0) \move(0
0)\rlvec(0 11.3) \move(-1 0)\rlvec(0 11.3) \move(-2 0)\rlvec(0
11.3) \move(-3 0)\rlvec(0 11.3) \move(-4 0)\rlvec(0 11.3)
\htext(-0.5 0.25){$0$} \htext(-1.5 0.25){$0$} \htext(-2.5
0.25){$0$} \htext(-3.5 0.25){$0$} \htext(-0.5 0.75){$0$}
\htext(-1.5 0.75){$0$} \htext(-2.5 0.75){$0$} \htext(-3.5
0.75){$0$} \htext(-0.5 1.5){$1$} \htext(-1.5 1.5){$1$} \htext(-2.5
1.5){$1$} \htext(-3.5 1.5){$1$} \vtext(-0.5 2.75){$\cdots$}
\vtext(-1.5 2.75){$\cdots$} \vtext(-2.5 2.75){$\cdots$}
\vtext(-3.5 2.75){$\cdots$} \htext(-0.5 4){$n\!\!-\!\!1$}
\htext(-1.5 4){$n\!\!-\!\!1$} \htext(-2.5 4){$n\!\!-\!\!1$}
\htext(-3.5 4){$n\!\!-\!\!1$} \htext(-0.5 4.75){$n$} \htext(-1.5
4.75){$n$} \htext(-2.5 4.75){$n$} \htext(-3.5 4.75){$n$}
\htext(-0.5 5.25){$n$} \htext(-1.5 5.25){$n$} \htext(-2.5
5.25){$n$} \htext(-3.5 5.25){$n$} \htext(-0.5 6){$n\!\!-\!\!1$}
\htext(-1.5 6){$n\!\!-\!\!1$} \htext(-2.5 6){$n\!\!-\!\!1$}
\htext(-3.5 6){$n\!\!-\!\!1$} \vtext(-0.5 7.25){$\cdots$}
\vtext(-1.5 7.25){$\cdots$} \vtext(-2.5 7.25){$\cdots$}
\vtext(-3.5 7.25){$\cdots$} \htext(-0.5 8.5){$1$} \htext(-1.5
8.5){$1$} \htext(-2.5 8.5){$1$} \htext(-3.5 8.5){$1$} \htext(-0.5
9.25){$0$} \htext(-1.5 9.25){$0$} \htext(-2.5 9.25){$0$}
\htext(-3.5 9.25){$0$} \htext(-0.5 9.75){$0$} \htext(-1.5
9.75){$0$} \htext(-2.5 9.75){$0$} \htext(-3.5 9.75){$0$}
\htext(-0.5 10.5){$1$} \htext(-1.5 10.5){$1$} \htext(-2.5
10.5){$1$} \htext(-3.5 10.5){$1$}
\end{texdraw}}\hskip 2cm
On $Y_{\Lambda_n}$ : \raisebox{-1\height}{
\begin{texdraw}\textref h:C v:C \fontsize{6}{6}\selectfont \drawdim
mm \setunitscale 5 \move(0 0)\lvec(-4 0)\lvec(-4 0.5)\lvec(0
0.5)\ifill f:0.7 \move(0 0)\rlvec(-4.3 0) \move(0 0.5)\rlvec(-4.3
0) \move(0 1)\rlvec(-4.3 0) \move(0 2)\rlvec(-4.3 0) \move(0
3.5)\rlvec(-4.3 0) \move(0 4.5)\rlvec(-4.3 0) \move(0
4.5)\rlvec(-4.3 0) \move(0 5)\rlvec(-4.3 0) \move(0
5.5)\rlvec(-4.3 0) \move(0 6.5)\rlvec(-4.3 0) \move(0
8)\rlvec(-4.3 0) \move(0 9)\rlvec(-4.3 0) \move(0 9.5)\rlvec(-4.3
0) \move(0 10)\rlvec(-4.3 0) \move(0 11)\rlvec(-4.3 0) \move(0
0)\rlvec(0 11.3) \move(-1 0)\rlvec(0 11.3) \move(-2 0)\rlvec(0
11.3) \move(-3 0)\rlvec(0 11.3) \move(-4 0)\rlvec(0 11.3)
\htext(-0.5 0.25){$n$} \htext(-1.5 0.25){$n$} \htext(-2.5
0.25){$n$} \htext(-3.5 0.25){$n$} \htext(-0.5 0.75){$n$}
\htext(-1.5 0.75){$n$} \htext(-2.5 0.75){$n$} \htext(-3.5
0.75){$n$} \htext(-0.5 1.5){$n\!\!-\!\!1$} \htext(-1.5
1.5){$n\!\!-\!\!1$} \htext(-2.5 1.5){$n\!\!-\!\!1$} \htext(-3.5
1.5){$n\!\!-\!\!1$} \vtext(-0.5 2.75){$\cdots$} \vtext(-1.5
2.75){$\cdots$} \vtext(-2.5 2.75){$\cdots$} \vtext(-3.5
2.75){$\cdots$} \htext(-0.5 4){$1$} \htext(-1.5 4){$1$}
\htext(-2.5 4){$1$} \htext(-3.5 4){$1$} \htext(-0.5 4.75){$0$}
\htext(-1.5 4.75){$0$} \htext(-2.5 4.75){$0$} \htext(-3.5
4.75){$0$} \htext(-0.5 5.25){$0$} \htext(-1.5 5.25){$0$}
\htext(-2.5 5.25){$0$} \htext(-3.5 5.25){$0$} \htext(-0.5 6){$1$}
\htext(-1.5 6){$1$} \htext(-2.5 6){$1$} \htext(-3.5 6){$1$}
\vtext(-0.5 7.25){$\cdots$} \vtext(-1.5 7.25){$\cdots$}
\vtext(-2.5 7.25){$\cdots$} \vtext(-3.5 7.25){$\cdots$}
\htext(-0.5 8.5){$n\!\!-\!\!1$} \htext(-1.5 8.5){$n\!\!-\!\!1$}
\htext(-2.5 8.5){$n\!\!-\!\!1$} \htext(-3.5 8.5){$n\!\!-\!\!1$}
\htext(-0.5 9.25){$n$} \htext(-1.5 9.25){$n$} \htext(-2.5
9.25){$n$} \htext(-3.5 9.25){$n$} \htext(-0.5 9.75){$n$}
\htext(-1.5 9.75){$n$} \htext(-2.5 9.75){$n$} \htext(-3.5
9.75){$n$} \htext(-0.5 10.5){$n\!\!-\!\!1$} \htext(-1.5
10.5){$n\!\!-\!\!1$} \htext(-2.5 10.5){$n\!\!-\!\!1$} \htext(-3.5
10.5){$n\!\!-\!\!1$}
\end{texdraw}}
\end{center}

(f) $B_n^{(1)}$ ($n\geq 3$),
\begin{center}
On $Y_{\Lambda_0}$ : \raisebox{-1\height}{
\begin{texdraw}\textref h:C v:C \fontsize{6}{6}\selectfont \drawdim mm
\setunitscale 5 \nc{\dtri}{ \bsegment \move(-1 0)\lvec(0 1)\lvec(0
0)\lvec(-1 0)\ifill f:0.7 \esegment } \move(0 0)\dtri \move(-1
0)\dtri \move(-2 0)\dtri \move(-3 0)\dtri \move(0 0)\rlvec(-4.3 0)
\move(0 1)\rlvec(-4.3 0) \move(0 2)\rlvec(-4.3 0) \move(0
3.5)\rlvec(-4.3 0) \move(0 4.5)\rlvec(-4.3 0) \move(0
5.5)\rlvec(-4.3 0) \move(0 6.5)\rlvec(-4.3 0) \move(0
8)\rlvec(-4.3 0) \move(0 9)\rlvec(-4.3 0) \move(0 10)\rlvec(-4.3
0) \move(0 11)\rlvec(-4.3 0) \move(0 0)\rlvec(0 11.3) \move(-1
0)\rlvec(0 11.3) \move(-2 0)\rlvec(0 11.3) \move(-3 0)\rlvec(0
11.3) \move(-4 0)\rlvec(0 11.3) \move(-1 0)\rlvec(1 1) \move(-2
0)\rlvec(1 1) \move(-3 0)\rlvec(1 1) \move(-4 0)\rlvec(1 1)
\move(-1 9)\rlvec(1 1) \move(-2 9)\rlvec(1 1) \move(-3 9)\rlvec(1
1) \move(-4 9)\rlvec(1 1) \move(0 5)\rlvec(-4.3 0) \htext(-0.3
0.25){$1$} \htext(-0.75 0.75){$0$} \htext(-0.5 1.5){$2$}
\vtext(-0.5 2.75){$\cdots$} \htext(-0.5 4){$n\!\!-\!\!1$}
\htext(-0.5 6){$n\!\!-\!\!1$} \htext(-0.5 8.5){$2$} \htext(-0.3
9.25){$1$} \htext(-0.75 9.75){$0$} \htext(-0.5 10.5){$2$}
\htext(-2.3 0.25){$1$} \htext(-2.75 0.75){$0$} \htext(-2.5
1.5){$2$} \vtext(-2.5 2.75){$\cdots$} \htext(-2.5
4){$n\!\!-\!\!1$} \htext(-2.5 6){$n\!\!-\!\!1$} \htext(-2.5
8.5){$2$} \htext(-2.3 9.25){$1$} \htext(-2.75 9.75){$0$}
\htext(-2.5 10.5){$2$} \htext(-1.3 0.25){$0$} \htext(-1.75
0.75){$1$} \htext(-1.5 1.5){$2$} \vtext(-1.5 2.75){$\cdots$}
\htext(-1.5 4){$n\!\!-\!\!1$} \htext(-1.5 6){$n\!\!-\!\!1$}
\htext(-1.5 8.5){$2$} \htext(-1.3 9.25){$0$} \htext(-1.75
9.75){$1$} \htext(-1.5 10.5){$2$} \htext(-3.3 0.25){$0$}
\htext(-3.75 0.75){$1$} \htext(-3.5 1.5){$2$} \vtext(-3.5
2.75){$\cdots$} \htext(-3.5 4){$n\!\!-\!\!1$} \htext(-3.5
6){$n\!\!-\!\!1$} \htext(-3.5 8.5){$2$} \htext(-3.3 9.25){$0$}
\htext(-3.75 9.75){$1$} \htext(-3.5 10.5){$2$} \htext(-0.5
4.75){$n$} \htext(-2.5 4.75){$n$} \htext(-1.5 5.25){$n$}
\htext(-3.5 5.25){$n$} \htext(-1.5 4.75){$n$} \htext(-3.5
4.75){$n$} \htext(-0.5 5.25){$n$} \htext(-2.5 5.25){$n$}
\vtext(-0.5 7.25){$\cdots$} \vtext(-1.5 7.25){$\cdots$}
\vtext(-2.5 7.25){$\cdots$} \vtext(-3.5 7.25){$\cdots$}
\end{texdraw}}\hskip 1cm On $Y_{\Lambda_1}$:\raisebox{-1\height}{
\begin{texdraw} \textref h:C v:C \fontsize{6}{6}\selectfont \drawdim mm
\setunitscale 5 \nc{\dtri}{ \bsegment \move(-1 0)\lvec(0 1)\lvec(0
0)\lvec(-1 0)\ifill f:0.7 \esegment } \move(0 0)\dtri \move(-1
0)\dtri \move(-2 0)\dtri \move(-3 0)\dtri \move(0 0)\rlvec(-4.3 0)
\move(0 1)\rlvec(-4.3 0) \move(0 2)\rlvec(-4.3 0) \move(0
3.5)\rlvec(-4.3 0) \move(0 4.5)\rlvec(-4.3 0) \move(0
5.5)\rlvec(-4.3 0) \move(0 6.5)\rlvec(-4.3 0) \move(0
8)\rlvec(-4.3 0) \move(0 9)\rlvec(-4.3 0) \move(0 10)\rlvec(-4.3
0) \move(0 11)\rlvec(-4.3 0) \move(0 0)\rlvec(0 11.3) \move(-1
0)\rlvec(0 11.3) \move(-2 0)\rlvec(0 11.3) \move(-3 0)\rlvec(0
11.3) \move(-4 0)\rlvec(0 11.3) \move(-1 0)\rlvec(1 1) \move(-2
0)\rlvec(1 1) \move(-3 0)\rlvec(1 1) \move(-4 0)\rlvec(1 1)
\move(-1 9)\rlvec(1 1) \move(-2 9)\rlvec(1 1) \move(-3 9)\rlvec(1
1) \move(-4 9)\rlvec(1 1) \move(0 5)\rlvec(-4.3 0) \htext(-0.3
0.25){$0$} \htext(-0.75 0.75){$1$} \htext(-0.5 1.5){$2$}
\vtext(-0.5 2.75){$\cdots$} \htext(-0.5 4){$n\!\!-\!\!1$}
\htext(-0.5 6){$n\!\!-\!\!1$} \htext(-0.5 8.5){$2$} \htext(-0.3
9.25){$0$} \htext(-0.75 9.75){$1$} \htext(-0.5 10.5){$2$}
\htext(-2.3 0.25){$0$} \htext(-2.75 0.75){$1$} \htext(-2.5
1.5){$2$} \vtext(-2.5 2.75){$\cdots$} \htext(-2.5
4){$n\!\!-\!\!1$} \htext(-2.5 6){$n\!\!-\!\!1$} \htext(-2.5
8.5){$2$} \htext(-2.3 9.25){$0$} \htext(-2.75 9.75){$1$}
\htext(-2.5 10.5){$2$} \htext(-1.3 0.25){$1$} \htext(-1.75
0.75){$0$} \htext(-1.5 1.5){$2$} \vtext(-1.5 2.75){$\cdots$}
\htext(-1.5 4){$n\!\!-\!\!1$} \htext(-1.5 6){$n\!\!-\!\!1$}
\htext(-1.5 8.5){$2$} \htext(-1.3 9.25){$1$} \htext(-1.75
9.75){$0$} \htext(-1.5 10.5){$2$} \htext(-3.3 0.25){$1$}
\htext(-3.75 0.75){$0$} \htext(-3.5 1.5){$2$} \vtext(-3.5
2.75){$\cdots$} \htext(-3.5 4){$n\!\!-\!\!1$} \htext(-3.5
6){$n\!\!-\!\!1$} \htext(-3.5 8.5){$2$} \htext(-3.3 9.25){$1$}
\htext(-3.75 9.75){$0$} \htext(-3.5 10.5){$2$} \htext(-0.5
4.75){$n$} \htext(-2.5 4.75){$n$} \htext(-1.5 5.25){$n$}
\htext(-3.5 5.25){$n$} \htext(-1.5 4.75){$n$} \htext(-3.5
4.75){$n$} \htext(-0.5 5.25){$n$} \htext(-2.5 5.25){$n$}
\vtext(-0.5 7.25){$\cdots$} \vtext(-1.5 7.25){$\cdots$}
\vtext(-2.5 7.25){$\cdots$} \vtext(-3.5 7.25){$\cdots$}
\end{texdraw}}\vskip 5mm

On $Y_{\Lambda_n}$ : \raisebox{-1\height}{
\begin{texdraw}\textref h:C v:C \fontsize{6}{6}\selectfont \drawdim mm
\setunitscale 5 \move(0 0)\lvec(-4 0)\lvec(-4 0.5)\lvec(0
0.5)\ifill f:0.7 \move(0 0)\rlvec(-4.3 0) \move(0 1)\rlvec(-4.3 0)
\move(0 2)\rlvec(-4.3 0) \move(0 3.5)\rlvec(-4.3 0) \move(0
4.5)\rlvec(-4.3 0) \move(0 5.5)\rlvec(-4.3 0) \move(0
6.5)\rlvec(-4.3 0) \move(0 8)\rlvec(-4.3 0) \move(0 9)\rlvec(-4.3
0) \move(0 10)\rlvec(-4.3 0) \move(0 11)\rlvec(-4.3 0) \move(0
0)\rlvec(0 11.3) \move(-1 0)\rlvec(0 11.3) \move(-2 0)\rlvec(0
11.3) \move(-3 0)\rlvec(0 11.3) \move(-4 0)\rlvec(0 11.3) \move(-1
4.5)\rlvec(1 1) \move(-2 4.5)\rlvec(1 1) \move(-3 4.5)\rlvec(1 1)
\move(-4 4.5)\rlvec(1 1) \htext(-0.5 0.25){$n$} \htext(-0.5
1.5){$n\!\!-\!\!1$} \vtext(-0.5 2.75){$\cdots$} \htext(-0.5
4){$2$} \htext(-0.3 4.75){$1$} \htext(-0.75 5.25){$0$} \htext(-0.5
6){$2$} \vtext(-0.5 7.25){$\cdots$} \htext(-0.5
8.5){$n\!\!-\!\!1$} \htext(-0.5 9.25){$n$} \htext(-0.5
10.5){$n\!\!-\!\!1$} \htext(-2.5 0.25){$n$} \htext(-2.5
1.5){$n\!\!-\!\!1$} \vtext(-2.5 2.75){$\cdots$} \htext(-2.5
4){$2$} \htext(-2.3 4.75){$1$} \htext(-2.75 5.25){$0$} \htext(-2.5
6){$2$} \vtext(-2.5 7.25){$\cdots$} \htext(-2.5
8.5){$n\!\!-\!\!1$} \htext(-2.5 9.25){$n$} \htext(-2.5
10.5){$n\!\!-\!\!1$} \htext(-1.5 0.75){$n$} \htext(-1.5
1.5){$n\!\!-\!\!1$} \vtext(-1.5 2.75){$\cdots$} \htext(-1.5
4){$2$} \htext(-1.3 4.75){$0$} \htext(-1.75 5.25){$1$} \htext(-1.5
6){$2$} \vtext(-1.5 7.25){$\cdots$} \htext(-1.5
8.5){$n\!\!-\!\!1$} \htext(-1.5 9.75){$n$} \htext(-1.5
10.5){$n\!\!-\!\!1$} \htext(-3.5 0.75){$n$} \htext(-3.5
1.5){$n\!\!-\!\!1$} \vtext(-3.5 2.75){$\cdots$} \htext(-3.5
4){$2$} \htext(-3.3 4.75){$0$} \htext(-3.75 5.25){$1$} \htext(-3.5
6){$2$} \vtext(-3.5 7.25){$\cdots$} \htext(-3.5
8.5){$n\!\!-\!\!1$} \htext(-3.5 9.75){$n$} \htext(-3.5
10.5){$n\!\!-\!\!1$} \move(0 0.5)\rlvec(-4.3 0) \move(0
9.5)\rlvec(-4.3 0) \htext(-3.5 9.25){$n$} \htext(-3.5 0.25){$n$}
\htext(-2.5 9.75){$n$} \htext(-2.5 0.75){$n$} \htext(-1.5
0.25){$n$} \htext(-1.5 9.25){$n$} \htext(-0.5 9.75){$n$}
\htext(-0.5 0.75){$n$}
\end{texdraw}}
\end{center}
%


{\small

\end{document}